\documentclass[oneside,english,reqno]{book}
\usepackage[T1]{fontenc}
\usepackage[latin9]{inputenc}
\usepackage{babel}
\usepackage{array}
\usepackage{longtable}
\usepackage{refstyle}
\usepackage{float}
\usepackage{textcomp}
\usepackage{mathrsfs}
\usepackage{mathtools}
\usepackage{url}
\usepackage{enumitem}
\usepackage{multirow}
\usepackage{amsthm}
\usepackage{amsmath}
\usepackage{amssymb}
\usepackage{stmaryrd}
\usepackage{splitidx}
\makeindex
\usepackage{graphicx}
\usepackage{esint}
\usepackage[all]{xy}
\PassOptionsToPackage{normalem}{ulem}
\usepackage{ulem}
\usepackage{nomencl}
\providecommand{\printnomenclature}{\printglossary}
\providecommand{\makenomenclature}{\makeglossary}
\makenomenclature
\newindex[Index]{idx}
\newindex[Quotes: index of credits.]{nam}
\usepackage[unicode=true,pdfusetitle,
 bookmarks=true,bookmarksnumbered=false,bookmarksopen=false,
 breaklinks=false,pdfborder={0 0 1},backref=false,colorlinks=false]
 {hyperref}
\hypersetup{
 colorlinks=true,citecolor=blue,linkcolor=blue,linktocpage=true,urlcolor=blue}

\makeatletter

\AtBeginDocument{\providecommand\figref[1]{\ref{fig:#1}}}
\AtBeginDocument{\providecommand\chapref[1]{\ref{chap:#1}}}
\AtBeginDocument{\providecommand\secref[1]{\ref{sec:#1}}}
\AtBeginDocument{\providecommand\tabref[1]{\ref{tab:#1}}}
\AtBeginDocument{\providecommand\exerref[1]{\ref{exer:#1}}}
\AtBeginDocument{\providecommand\defref[1]{\ref{def:#1}}}
\AtBeginDocument{\providecommand\lemref[1]{\ref{lem:#1}}}
\AtBeginDocument{\providecommand\remref[1]{\ref{rem:#1}}}
\AtBeginDocument{\providecommand\corref[1]{\ref{cor:#1}}}
\AtBeginDocument{\providecommand\thmref[1]{\ref{thm:#1}}}
\AtBeginDocument{\providecommand\exaref[1]{\ref{exa:#1}}}
\AtBeginDocument{\providecommand\propref[1]{\ref{prop:#1}}}
\AtBeginDocument{\providecommand\eqref[1]{\ref{eq:#1}}}
\AtBeginDocument{\providecommand\exref[1]{\ref{ex:#1}}}
\let\SF@@footnote\footnote
\def\footnote{\ifx\protect\@typeset@protect
    \expandafter\SF@@footnote
  \else
    \expandafter\SF@gobble@opt
  \fi
}
\expandafter\def\csname SF@gobble@opt \endcsname{\@ifnextchar[%]
  \SF@gobble@twobracket
  \@gobble
}
\edef\SF@gobble@opt{\noexpand\protect
  \expandafter\noexpand\csname SF@gobble@opt \endcsname}
\def\SF@gobble@twobracket[#1]#2{}
\providecommand{\tabularnewline}{\\}
\RS@ifundefined{subref}
  {\def\RSsubtxt{section~}\newref{sub}{name = \RSsubtxt}}
  {}
\RS@ifundefined{thmref}
  {\def\RSthmtxt{theorem~}\newref{thm}{name = \RSthmtxt}}
  {}
\RS@ifundefined{lemref}
  {\def\RSlemtxt{lemma~}\newref{lem}{name = \RSlemtxt}}
  {}

\usepackage{enumitem}		% customizable list environments
\theoremstyle{plain}
\ifx\thechapter\undefined
\newtheorem{thm}{\protect\theoremname}
\else
\newtheorem{thm}{\protect\theoremname}[chapter]
\fi
  \theoremstyle{definition}
  \newtheorem{defn}[thm]{\protect\definitionname}
  \theoremstyle{plain}
  \newtheorem{lem}[thm]{\protect\lemmaname}
  \theoremstyle{definition}
  \newtheorem{example}[thm]{\protect\examplename}
  \theoremstyle{remark}
  \newtheorem{rem}[thm]{\protect\remarkname}
  \theoremstyle{definition}
  \newtheorem{xca}[thm]{\protect\exercisename}
  \theoremstyle{plain}
  \newtheorem{cor}[thm]{\protect\corollaryname}
  \theoremstyle{plain}
  \newtheorem{prop}[thm]{\protect\propositionname}
  \theoremstyle{remark}
  \newtheorem{note}[thm]{\protect\notename}
 \newlist{casenv}{enumerate}{4}
 \setlist[casenv]{leftmargin=*,align=left,widest={iiii}}
 \setlist[casenv,1]{label={{\itshape\ \casename} \arabic*.},ref=\arabic*}
 \setlist[casenv,2]{label={{\itshape\ \casename} \roman*.},ref=\roman*}
 \setlist[casenv,3]{label={{\itshape\ \casename\ \alph*.}},ref=\alph*}
 \setlist[casenv,4]{label={{\itshape\ \casename} \arabic*.},ref=\arabic*}
  \theoremstyle{plain}
  \newtheorem{conjecture}[thm]{\protect\conjecturename}
  \theoremstyle{remark}
  \newtheorem{claim}[thm]{\protect\claimname}

\usepackage{needspace,bbm,appendix}

\date{}

\def\exercisename{Exercise}
\def\listexercisename{List of Exercises}

\newcommand\listofexercises{%
    \if@twocolumn
      \@restonecoltrue\onecolumn
    \else
      \@restonecolfalse
    \fi
    \chapter*{\listexercisename}%
      \@mkboth{\MakeUppercase\listexercisename}%
              {\MakeUppercase\listexercisename}%
    \@starttoc{loe}%
    \if@restonecol\twocolumn\fi
    }

\newcommand{\myexercise}[1]
{
\addcontentsline{loe}{figure}{\numberline{\thexca}{#1}}
}

\let\oldlof=\listoffigures
\renewcommand{\listoffigures}{
\cleardoublepage
\phantomsection
\oldlof
}

\setlist[description]{labelindent=0pt,style=multiline,leftmargin=3.5cm}

\setlist[description]{labelindent=0pt,style=multiline,leftmargin=3.5cm}

\renewcommand*\env@matrix[1][*\c@MaxMatrixCols c]{%
  \hskip -\arraycolsep
  \let\@ifnextchar\new@ifnextchar
  \array{#1}}

 \renewcommand*\l@figure{\@dottedtocline{1}{1em}{3.2em}}

\usepackage{ifthen}
\renewenvironment{figure}[1][]{%
 \ifthenelse{\equal{#1}{}}{%
   \@float{figure}
 }{%
   \@float{figure}[#1]%
 }%
 \centering
}{%
 \end@float
}

\renewenvironment{table}[1][]{%
 \ifthenelse{\equal{#1}{}}{%
   \@float{table}
 }{%
   \@float{table}[#1]%
 }%
 \centering
}{%
 \end@float
}

\newref{lem}{refcmd={Lemma \ref{#1}}}
\newref{thm}{refcmd={Theorem \ref{#1}}}
\newref{cor}{refcmd={Corollary \ref{#1}}}
\newref{sec}{refcmd={Section \ref{#1}}}
\newref{chap}{refcmd={Chapter \ref{#1}}}
\newref{prop}{refcmd={Proposition \ref{#1}}}
\newref{exa}{refcmd={Example \ref{#1}}}
\newref{exer}{refcmd={Exercise \ref{#1}}}
\newref{tab}{refcmd={Table \ref{#1}}}
\newref{fig}{refcmd={Figure \ref{#1}}}
\newref{sub}{refcmd={Section \ref{#1}}}
\newref{def}{refcmd={Definition \ref{#1}}}
\newref{rem}{refcmd={Remark \ref{#1}}}

\@ifundefined{showcaptionsetup}{}{%
 \PassOptionsToPackage{caption=false}{subfig}}
\usepackage{subfig}
\makeatother

  \providecommand{\claimname}{Claim}
  \providecommand{\conjecturename}{Conjecture}
  \providecommand{\corollaryname}{Corollary}
  \providecommand{\definitionname}{Definition}
  \providecommand{\examplename}{Example}
  \providecommand{\exercisename}{Exercise}
  \providecommand{\lemmaname}{Lemma}
  \providecommand{\notename}{Note}
  \providecommand{\propositionname}{Proposition}
  \providecommand{\remarkname}{Remark}
 \providecommand{\casename}{Case}
\providecommand{\theoremname}{Theorem}

\begin{document}
\frontmatter

\title{Noncommutative analysis,\\
Multivariable spectral theory for operators in Hilbert space,\\
Probability, and Unitary Representations}

\author{Palle Jorgensen, and Feng Tian}

\maketitle
\newpage{}
\begin{quote}
```The Journal of Functional Analysis' is dedicated to the broadening
of the horizons of functional analysis. Accordingly, it encourages
original research papers of high quality from all branches of science,
provided the core and flavor are of a functional analytic character
and the paper is in accordance with contemporary mathematical standards.''

--- From the cover page of `The Journal of Functional Analysis'
(written in 1967)
\end{quote}
\newpage{}
\begin{quote}
Dedicated to the memory of William B. Arveson (see \cite{MJD15})
\\
(22 November 1934 \textendash{} 15 November 2011).\vspace{2cm}

\textquotedblleft What goes around has come around, and today quantum
information theory (QIT) has led us back into a finite-dimensional
context. Completely positive maps on matrix algebras are the objects
that are dual to quantum channels; in fact, the study of quantum channels
reduces to the study of completely positive maps of matrix algebras
that preserve the unit. This is an area that is still undergoing vigorous
development in efforts to understand entanglement, entropy, and channel-capacity
in QIT.\textquotedblright{} 

--- William B. Arveson (written around 2009.)\sindex[nam]{Arveson, W., (1934-2011)}
\end{quote}
\setcounter{tocdepth}{-1}

\chapter*{Foreword{\small{}}\protect \\
{\small{}by Wayne Polyzou, Professor of Physics, University of Iowa}}

\chaptermark{Foreword}
\setcounter{tocdepth}{1}
\addcontentsline{toc}{chapter}{Foreword}

Progress in science, engineering and mathematics comes fast and it
often requires a significant effort to keep up with the advances in
other fields that impact applications. Functional analysis (especially
operators in Hilbert space, unitary representations of Lie groups,
and spectral theory) is one discipline that impacts my physics research.
Bringing my students up to speed with the subject facilitates their
ability to efficiently perform research, however the typical curriculum
in functional analysis courses is not directed to practitioners whose
primary objective is applications. This is also reflected in the many
excellent available texts on the subject, which primarily focus on
the mathematics, and are directed at students aspiring to a career
in mathematics.

I have been fortunate to have Palle Jorgensen as a colleague. He participates
in a weekly joint mathematical physics seminar that is attended by
faculty and students from both departments. It provides a forum to
address questions related to the role of mathematics in physics research.
Professor Jorgensen has a healthy appreciation of applications of
functional analysis; in these seminars he has been at the center of
discussions on a diverse range of applications involving wavelets,
reflection positivity, path integrals, entanglement, financial mathematics,
and algebraic field theory.

A number of the mathematically inclined students in my department
have benefited from taking the functional analysis course taught by
Professor Jorgensen. These students are motivated to enroll in his
class because the course material includes a significant discussion
of applications of functional analysis to subjects that interest them.

This book is based on the course that Professor Jorgensen teaches
on functional analysis. It fills in a gap that is not addressed by
the many excellent available texts on functional analysis, by using
applications to motivate basic results in functional analysis. The
way that it uses applications makes the material more accessible to
students; particularly for students who will eventually find careers
in related disciplines. The book also points to additional reference
material for students who are motivated to learn more about a specific
topic.\vspace{1em}

\begin{flushright}
W. Polyzou\vspace{2em}
\par\end{flushright}

\noindent \rule{1\textwidth}{1pt}\vspace{2em}

Over the decades, Functional Analysis, and the theory of operators
in Hilbert space, have been enriched and inspired on account of demands
from neighboring fields, within mathematics, \emph{harmonic analysis}
(wavelets and signal processing), \emph{numerical analysis} (finite
element methods, discretization), \emph{PDEs} (diffusion equations,
scattering theory), \emph{representation theory}; \emph{iterated function
systems} (fractals, Julia sets, chaotic dynamical systems), \emph{ergodic
theory}, \emph{operator algebras}, and many more. And neighboring
areas, \emph{probability}/\emph{statistics} (for example stochastic
processes, It\={o} and Malliavin calculus), \emph{physics} (representation
of Lie groups, quantum field theory), and \emph{spectral theory} for
Schrödinger operators. \sindex[nam]{Lie, Sophus} \index{quantum field theory}

\index{harmonic}

\index{ergodic}

\index{groups!Lie}

\index{signal}

\index{Lie!group}

\index{Lie!algebra}

\index{signal processing}

The book is based on a course sequence (two-semesters 313-314) taught,
over the years at the University of Iowa, by the first-named author.
The students in the course made up a mix: some advanced undergraduates,
but most of them, first or second year graduate students (from math,
as well as some from physics and stats.)

We have subsequently expanded the course notes taken by the second-named
author: we completed several topics from the original notes, and we
added a number of others, so that the book is now self-contained,
and it covers a unified theme; and yet it stresses a multitude of
applications. And it offers flexibility for users.

A glance at the table of contents makes it clear that our aim, and
choice of topics, is different from that of more traditional Functional
Analysis courses. This is deliberate. For example, in our choice of
topics, we placed emphasis on the use of Hilbert space techniques
which are then again used in our treatment of central themes of applied
functional analysis.

We have also strived for a more accessible book, and yet aimed squarely
at applications; --- we have been serious about motivation: Rather
than beginning with the four big theorems in Functional Analysis,
our point of departure is an initial choice of topics from applications.
And we have aimed for flexibility of use; acknowledging that students
and instructors will invariably have a host of diverse goals in teaching
beginning analysis courses. And students come to the course with a
varied background. Indeed, over the years we found that students have
come to the Functional Analysis sequence from other and different
areas of math, and even from other departments; and so we have presented
the material in a way that minimizes the need for prerequisites. We
also found that well motivated students are easily able to fill in
what is needed from measure theory, or from a facility with the four
big theorems of Functional Analysis. And we found that the approach
\textquotedblleft learn-by-using\textquotedblright{} has a comparative
advantage. \index{Four Big Theorems in Functional Analysis}

\paragraph{Analysis of Continuous Systems vs Discrete (Networks and Graphs)}

\begin{figure}
\subfloat[A binary tree model.]{\protect\includegraphics[width=0.35\textwidth]{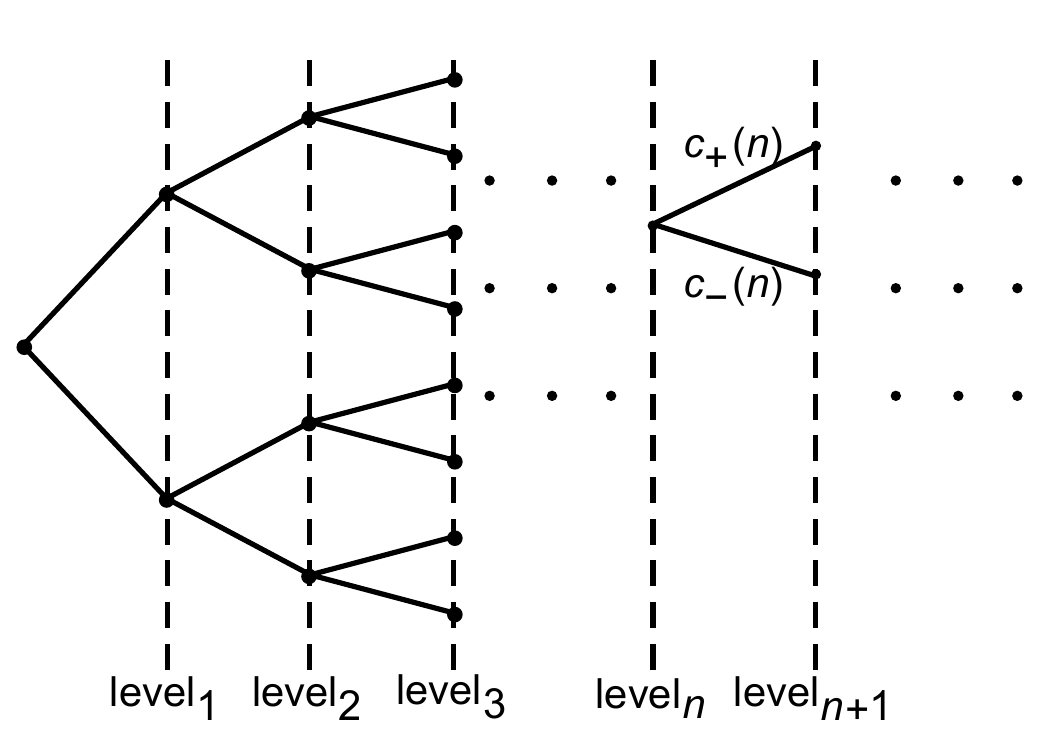}

}\hfill{}\subfloat[A Bratteli diagram. See \cite{BJO04}.]{\protect\includegraphics[width=0.6\textwidth]{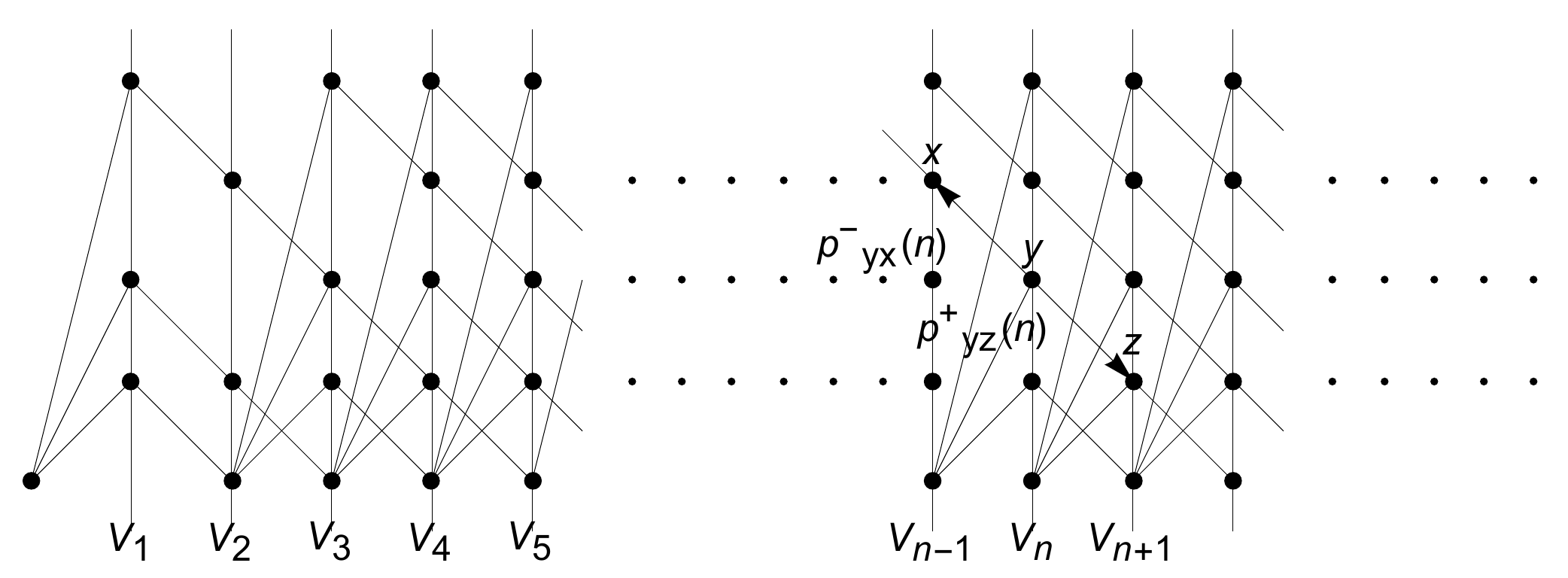}

}

\protect\caption{\label{fig:nt0}Examples of infinite weighted network; for details,
see \chapref{gLap}.}
\end{figure}

A new theme here, going beyond traditional books in the subject, is
applications of functional and harmonic analysis to \textquotedblleft \emph{large
networks},\textquotedblright{} so to discrete problems. More precisely,
we study infinite network models. Such models can often be represented
as follows: By a pair of sets, $V$ (vertices), and $E$, (edges).
In addition, one specifies a positive function $c$ defined on the
edge set $E$. (In electrical network models, $c$ represents conductance.)
There are then two associated operators $\Delta$ and $P$, each depending
on the triple $(V,E,c)$. (See \figref{nt0}.) Both operators represent
actions (i.e., operations) on appropriate spaces of functions, more
precisely functions defined on the infinite vertex set $V$. For the
networks of interest to us, the vertex set $V$ will be infinite,
reflecting statistical and stochastic properties; and it will have
additional geometric and ergodic theoretic properties. We are therefore
faced with a variety of choices of infinite-dimensional function spaces.
Many questions are of spectral theoretic flavor, and as a result,
the useful choices of function spaces will be Hilbert spaces.

But even restricting to Hilbert spaces, there are at least three natural
(and useful) candidates: (i) the plain $l^{2}$ sequence space, so
an $l^{2}$-space of functions on $V$, (ii) a suitably weighted $l^{2}$-space,
and finally (iii), an energy Hilbert space $\mathscr{H}_{E}$. (The
latter is an abstraction of more classical notions of Dirichlet spaces.)
Which one of the three to use depends on the particular operator considered,
and also on the questions asked.

In \emph{infinite network models}, both the Laplacian $\Delta$, and
the Markov operator $P$, will have infinite by infinite matrix representations.
Each of these infinite by infinite matrices will have the following
property: it will have non-zero entries localized only in finite bands
containing the infinite matrix-diagonal (i.e., they are infinite banded
matrices.) See \secref{dirac} and \figref{mp}. Thus, the standard
algebraic matrix operations will be well defined. \index{Laplacian!graph-}
\index{graph!network-}

Functional analytic and spectral theoretic tools now enter as follows:
In passing to appropriate Hilbert spaces, we arrive at various classes
of Hilbert space-operators. In the present setting, the operators
in question will be Hermitian, some unbounded, and some bounded. The
Laplacian $\Delta$ will typically be an unbounded operator, albeit
semibounded. When $\Delta$ is realized in the energy Hilbert space
$\mathscr{H}_{E}$, we must introduce boundary value considerations
in order to get selfadjoint extensions. By contrast, for the Markov
operator $P$, there is a weighted $l^{2}$-space such that $P$ is
a bounded, selfadjoint operator. Moreover, its spectrum is then contained
in the finite interval $[-1,1]$. In all of the operator realizations
by selfadjoint operators, $\Delta$ or $P$, the corresponding spectra
may be continuous, or may have a mix of spectral types, continuous
(singular or Lebesgue), and discrete.

For the operator theory, and graph Laplacians, for infinite network
models, we refer to \chapref{gLap}.

\paragraph*{An Apology}

The central themes in our book are as follows: (i) Operators in Hilbert
Space with emphasis on unbounded operators and non-commutativity,
(ii) Multivariable Spectral Theory, (iii) Noncommutative Analysis,
(iv) Probability with emphasis on Gaussian processes, and (v) Unitary
Representations. But more importantly, it is our goal to stress the
mutual interconnection between these five themes, or in fact central
areas of modern analysis. And these are interrelationships which in
the literature, at least up to now, have usually not been thought
of as especially related. But here we stress, and elaborate in detail,
on how a number of key theorems in anyone of these five areas crucially
impact advances in the others. Nonetheless, for readers expecting
a rehash of the standard list of topics from books in Functional Analysis
of past generations, therefore perhaps an apology is in order. 

The number of topics making up Functional Analysis is vast; and when
applications are added, the size and diversity are daunting. A glance
at the many books out there (see the partial list in Appendix \ref{chap:books})
will give readers an idea of the vast scope. It is by necessity that
we have made choices; and that readers will in all likelihood have
favorite topics not covered here. And there are probably surprises
too; -- things we cover here that are not typically included in standard
Functional Analysis books. We apologize to readers who had expected
a different table of contents. But we hope our choices are justified
in our discussion in Part \ref{part:intro} below.

Our glaring omissions among the big classical areas of Functional
Analysis include more technical aspects of the theory of Banach spaces.
Even in our consideration of $L^{p}$ spaces we have favored $p=1,2$,
or $\infty$. Although we have included some fundamentals from Banach
space theory, in this, we made a selection of only a few topics which
are of direct relevance to the concrete applications that we do include.
As for our bias in the choice of $L^{p}$ spaces, we can excuse this
in part by the familiar availability of interpolation theorems (the
interpolation refers to values of $p$), starting with the \emph{Riesz-Thorin
theorem} and related; see e.g., \cite{MR2381892,MR623243}. Moreover,
there is a host of books out there dealing with the exciting and deep
areas of Banach space theory, both new and classical; and we refer
readers to \cite{MR1392325,MR950981} for a sample.

Emphasis: In our applications, such as to physics, and statistics,
we have concentrated on those analysis tools that are directly needed
for the goal at hand. To fit the material into a single volume, we
have been forced to omit a number of classical areas of functional
analysis, and to concentrate on those that serve the applications
we have selected. And in particular we have omitted a number of proofs,
or reduced our discussion of some proofs to a few hints, or to exercises.
We feel this is justified as there are many great books out there
(see Appendix \ref{chap:books}) which contain complete proofs of
the big theorems in functional analysis; for example, the books by
W. Rudin, or by P. Lax.

\paragraph{Note on Presentation and Exercises}

In presenting our results, we have aimed for a reader-friendly account:
We found it helpful to include worked examples in order illustrate
abstract ideas. Main theorems hold in various degrees of generality,
but when appropriate, we have not chosen to present details in their
highest level of generality. Rather, we typically give the result
in a setting where the idea is more transparent, and easier to grasp.
But we do include comments about the more general versions; sketching
them in rough outline. The more general versions will typically be
easier for readers to follow, and to appreciate, after ideas have
been fleshed out in simpler contexts. We have made a second choice
in order to make it easier for students to grasp both ideas and the
technical details: We have included a lot of worked examples. At the
end of each of these examples, we then outline how details (from the
example in question) serve to illustrate one or more features in the
general theorems elsewhere in the book. Finally, we have made generous
use of both tables and figures. These are listed with page-references
at the end of the book.

We shall be using some terminology from neighboring areas. And in
order to help readers absorb that, we have included in Appendix \ref{chap:term}
a summary, with cited references, of some key notions from \emph{quantum
theory}, \emph{signal processing}, \emph{stochastic processes}, \emph{unitary
representations}, and from \emph{wavelet theory}.

Our selection of Exercises varies in level of difficulty, and they
vary in purpose as well. Some are really easy, aimed mainly for the
benefit of beginners; getting used a definition, or a concept. Others
are more traditional homework Exercises that can be assigned in a
standard course. And yet others are quite demanding.

We believe the material covered here is accessible to beginning graduate
students. In our choice of topics and presentation, we have aimed
for a user friendly book which hopefully will also be of interest
to both pure and applied mathematicians, as well as to students and
scientists, in anyone of a number of neighboring areas. In our presentation,
we have stressed applications, and also interconnections between disparate
areas.

\textbf{Reader\textquoteright s guide to References.} In the Reference
list, and in citations, we have included both books and research papers.
For the various themes, we have aimed at citing both original sources,
as well as timely papers; but we also cite brand new research. As
for the latter, i.e., the cited papers in the References dealing with
recent research (relating to the present topics), we mention a few,
followed by citations:

\begin{sloppypar}
\begin{itemize}
\item \textbf{Spectral theory:} \cite{AH13,CJK+12,HdSS12,Hel13,JoPe13b,JPT12,JPT12-1}.
\item \textbf{The theory of frames, including Kadison-Singer:} \cite{Wea03,Cas13,1407.4768,MR2830604,KOPT13,MSS15,MR2859705}.
(The paper \cite{MSS15} is the solution to K-S.) 
\item \textbf{Stochastic processes and applications:} \cite{MR2883397,MR2966130,MR3231624,Jor14}. 
\item \textbf{Analysis of infinite networks:} \cite{RAKK05,MR3042410,AJSV13,JoPe13,JT14c}.
\item \textbf{Representations of groups and algebras:} \cite{MR2261099,CM13,DHL09,JO00,JP12,JPS05,MR2681883,DJ08,JLH06}.
\item \textbf{Quantization and quantum information:} \cite{OdHo13,AAR13,CJK+12,MR2290770,Fan10,Maa10,OR07}.
\end{itemize}
\end{sloppypar}

\chapter*{Preface}

\addcontentsline{toc}{chapter}{Preface}

\emph{There are already many books in Functional Analysis, so why
another?}

The main reason is that we feel there is a need: in the teaching at
the beginning graduate level; more flexibility, more options for students
and instructors in pursuing new directions. And aiming for a book
which will help students with primary interests elsewhere to acquire
a facility with tools of a functional analytic flavor, say in spectral
theory for operators in Hilbert space, in commutative and non-commutative
harmonic analysis, in PDE, in numerical analysis, in stochastic processes,
or in physics.

\chapter*{Acknowledgments}

\addcontentsline{toc}{chapter}{Acknowledgment}

The first named author thanks his students in the Functional Analysis
sequence 313-314. Also he thanks postdocs, other colleagues, collaborators
over the years, and his mentors; --- the combined list includes the
following: D. Alpay, W. Arveson, S. Bezuglyi, O. Bratteli, P. Casazza,
I. Cho, D. Dutkay, B. Fuglede, S. Gelfand, T. Hida, R. Kadison, W.
Klink, K. Kornelson, D. Larson, G. Mackey, P. Muhly, Ekaterina Nathanson,
E. Nelson, R. Niedzialomski, J. Packer, E. Pearse, R.S. Phillips,
W. Polyzou, R.T. Powers, C. Radin, D. Robinson, W. Rudin, S. Sakai,
I.E. Segal, K. Shuman, Myung-Sin Song, M. H. Stone, R. Werner. 

We are especially grateful to Professor Wayne Polyzou (of The University
of Iowa Physics Department) for reading our manuscript, giving us
corrections, and making suggestions. And for graciously agreeing to
write a Foreword.

\tableofcontents{}

\nomenclature{BMO}{bounded mean oscillation}

\nomenclature{CP}{completely positive map}

\nomenclature{i.i.d}{independent identically distributed (system of random variables)}

\nomenclature{KS}{Kadison-Singer}

\nomenclature{ONB}{orthonormal basis (in Hilbert space)}

\nomenclature{PVM}{projection valued measure (the condition $P\left(A\cap B\right)=P\left(A\right)P\left(B\right)$ is part of the definition)}

\nomenclature{RKHS}{reproducing kernel Hilbert space}

\nomenclature{PDE}{partial differential equation; examples: the heat equation, diffusion equation, the wave equation, the Laplace equation, the Schrödinger equation.}

\nomenclature{PDO}{partial differential operator}

\nomenclature{ODE}{ordinary differential equation}

\nomenclature{SDE}{stochastic differential equation}

\nomenclature{Proj}{projection}

\nomenclature{irrep}{irreducible representation}

\nomenclature{ind}{induced representation}

\nomenclature{span}{linear span}

\nomenclature{conv}{convex hull}

\nomenclature{ext}{set of extreme-points}

\nomenclature{supp}{support of a function, a measure, or a distribution}

\nomenclature{$\aleph_{0}$}{aleph-sub 0, cardinality of $\mathbb{N}$}

\nomenclature{Rep($\mathfrak{A},\mathscr{H}$)}{representations of an algebra $\mathfrak{A}$ acting on a Hilbert space $\mathscr{H}$}

\nomenclature{Rep($G,\mathscr{H}$)}{representations of a group $G$ acting on a Hilbert space $\mathscr{H}$}

\nomenclature{$\Re\left\{ z\right\} $}{the real part of $z\in\mathbb{C}$}

\nomenclature{$\Im\left\{ z\right\} $}{the imaginary part of $z\in\mathbb{C}$}

\nomenclature{$\mathscr{O}_{N}$}{the Cuntz-algebra, i.e., generators $\left\{ s_{i}\right\} _{i=1}^{N}$ and relations, $s_{i}^{*}s_{j}=\delta_{i,j}\mathbbm{1}$, and $\sum_{1}^{N}s_{i}s_{i}^{*}=\mathbbm{1}$.}

\nomenclature{$NR_{T}$}{numerical range of a given operator $T$}

\printnomenclature[2cm]{}

\chapter*{Notation}

\addcontentsline{toc}{chapter}{Notation}

\paragraph{\emph{Operators in Hilbert space}}

\begin{description}
\item [{$A=A^{*}$}] selfadjoint
\item [{$A\subset A^{*}$}] symmetric (also called Hermitian)
\item [{$A\subset-A^{*}$}] skew-symmetric
\item [{$AA^{*}=A^{*}A$}] normal
\item [{$UU^{*}=U^{*}U=I$}] unitary
\item [{$I,\;\mbox{or}\;I_{\mathscr{H}}$}] identity operator in a given
Hilbert space $\mathscr{H}$, i.e., $I\left(v\right)=v$, $\forall v\in\mathscr{H}$
\item [{$P=P^{*}=P^{2}$}] projection
\item [{$\mathscr{G}\left(A\right)$}] graph of operator\vspace{3pt}
\item [{$\left\langle \cdot,\cdot\right\rangle $}] inner product of a
given Hilbert space $\mathscr{H}$, i.e., $\left\langle v,w\right\rangle $
for $v,w\in\mathscr{H}$; linear in the second variable.
\item [{$\left|v\left\rangle \right\langle w\right|$}] Dirac ket-bra vectors
for rank-one operator
\item [{$\bigwedge$}] lattice operation ``minimum'' applied to projections
\item [{$\bigvee$}] lattice operation ``maximum'' applied to projections
\item [{$\bigcap$}] set-theoretic intersection
\item [{$\bigcup$}] set-theoretic union
\item [{$\mbox{ess sup}$}] essential supremum
\item [{$l^{p}$}] sequence space, $l^{p}$-summable
\item [{$L^{p}\left(\mu\right)$}] $L^{p}$-integrable functions on a $\mu$
measure space
\item [{span}] all linear combination of a specified subset
\item [{$\overline{\mbox{span}}$}] closure of span
\item [{$E^{*}$}] the dual of a given normed space $E$ ($E^{*}$ is a
Banach space)
\item [{$E^{**}$}] double-dual
\item [{$\left(..\right)'$}] commutant of a set of operators
\item [{$\left(..\right)''$}] double-commutant
\item [{$\widehat{\cdot}$}] Fourier transform, or Gelfand transform
\item [{$\chi_{E}$}] indicator function of a set $E$
\item [{$\delta$}] Dirac delta ``function'' 
\item [{$*$}] convolution
\item [{sp}] (spec) spectrum
\item [{res}] resolvent set
\item [{$\mathcal{B}\left(X\right)$}] Borel sets, i.e., the sigma-algebra
generated by the open sets in a topological space $X$
\item [{$\mathscr{B}\left(\mathscr{H}\right)$}] all bounded linear operators
$\mathscr{H}\longrightarrow\mathscr{H}$
\item [{$\mathscr{F}R\left(\mathscr{H}\right)$}] all finite-rank operators
in $\mathscr{H}$
\item [{$\mbox{tr (trace)}$}] the trace functional
\item [{$\mathscr{T}_{1}\left(\mathscr{H}\right)$}] all trace class operators
in $\mathscr{B}\left(\mathscr{H}\right)$
\item [{$\mbox{Proj}\left(\mathscr{H}\right)$}] the lattice of all orthogonal
projections $P$ in a fixed Hilbert space, i.e., $P=P^{2}=P^{*}$.
\item [{$dom\left(A\right)$}] the domain of some linear operator $A$
\item [{$ran\left(A\right)\;\left(\mbox{or}\;R\left(A\right)\right)$}] the
range of $A$
\item [{$\mbox{Ker}\left(A\right)$}] the kernel of $A$
\item [{$M_{\varphi}$}] the operator of multiplication by some function
acting in some $L^{2}\left(\mu\right)$, or a multiplier in some RKHS. 
\item [{$T_{\varphi}=P_{+}M_{\varphi}P_{+}$}] Toeplitz-operator with symbol
$\varphi$
\item [{$\mu\circ T^{-1}$}] transformation of measure, i.e., $\left(\mu\circ T^{-1}\right)\left(\triangle\right)=\mu\left(T^{-1}\left(\triangle\right)\right)$,
$\triangle\in$ sigma-algebra, $T^{-1}\left(\triangle\right)=\left\{ x:Tx\in\triangle\right\} $.
\item [{$\int^{\oplus}$}] direct integral decomposition
\item [{$\oplus$}] orthogonal sum
\item [{$\otimes$}] tensor product
\item [{$P,\:Q$}] notation used for pairs of projections, but \uline{also}
for the momentum and position operators from quantum mechanics
\item [{$\mathbb{P}$}] some probability measure
\item [{$\mathbb{E}=\mathbb{E}_{\mathbb{P}}$}] expectation $\mathbb{E}_{\mathbb{P}}\left(X\right)=\int_{\Omega}Xd\mathbb{P}$
\item [{$G$}] Lie group
\item [{$\mathfrak{g}$}] Lie algebra
\item [{$\mathfrak{g}\xrightarrow{\;\exp\;}G$}] exponential mapping
\item [{$\mathcal{U}$}] representation of some Lie group $G$ 
\item [{$d\mathcal{U}$}] representation of the Lie algebra $\mathfrak{g}$
corresponding to $G$; the derived representation.
\item [{$C^{*}$-algebra}] an algebra $\mathfrak{A}$ with involution $\mathfrak{A}\ni a\rightarrow a^{*}\in\mathfrak{A}$,
$a^{**}=a$, $\left(ab\right)^{*}=b^{*}a^{*}$, $a,b\in\mathfrak{A}$;
and norm $\left\Vert \cdot\right\Vert $, such that $\left(\mathfrak{A},\left\Vert \cdot\right\Vert \right)$
is complete; and $\left\Vert ab\right\Vert \leq\left\Vert a\right\Vert \left\Vert b\right\Vert $,
$a,b\in\mathfrak{A}$, holds, as well as $\left\Vert a^{*}a\right\Vert =\left\Vert a\right\Vert ^{2}$,
$a\in\mathfrak{A}$.
\item [{$W^{*}$-algebra}] (also called \emph{von Neumann algebra}, or
a \emph{ring of operators}) A $W^{*}$-algebra is a $C^{*}$-algebra
$\mathfrak{A}$ which has the following additional property: There
is a Banach space $\left(X_{*},\left\Vert \cdot\right\Vert _{*}\right)$
such that $\mathfrak{A}$, with its  $C^{*}$-norm (i.e., $\left\Vert a^{*}a\right\Vert =\left\Vert a\right\Vert ^{2}$,
$a\in\mathfrak{A}$), is the dual, $\mathfrak{A}=\left(X_{*}\right)^{*}$.
When such a Banach space $X_{*}$ exists, it is called a \emph{pre-dual}.
(This characterization of $W^{*}$-algebra is due to Sakai \cite{MR0442701}.)
\end{description}

\paragraph{\emph{Operations on subspaces of Hilbert spaces $\mathscr{H}$}}
\begin{description}
\item [{$\mathscr{T}\subset\mathscr{H}$}] \begin{flushleft}
some subspace in $\mathscr{H}$
\par\end{flushleft}
\item [{$\mathscr{T}^{\perp}$}] \begin{flushleft}
ortho-complement 
\[
\mathscr{T}^{\perp}=\left\{ h\in\mathscr{H}\::\:\left\langle h,s\right\rangle =0,\:\forall s\in\mathscr{T}\right\} =\mathscr{H}\ominus\mathscr{T}
\]

\par\end{flushleft}
\end{description}
\begin{flushleft}
$\mathscr{T}^{\perp\perp}=\overline{span}\mathscr{T}$ 
\par\end{flushleft}

\paragraph{\emph{Normal or not! It depends:}}
\begin{itemize}
\item An \emph{operator} $T$ (bounded or not) is normal iff (Def.)
\[
T^{*}T=TT^{*}
\]

\item A \emph{state} $s$ on a $*$-algebra $\mathfrak{A}$ is normal if
it allows a representation $\left(\mathscr{H},\rho\right)$ where
$\mathscr{H}$ is a Hilbert space, and $\rho$ is a positive trace-class
operator in $\mathscr{H}$ such that $trace\left(\rho\right)=1$,
and 
\[
s\left(A\right)=trace\left(\rho A\right),\;\forall A\in\mathfrak{A}.
\]

\item A \emph{random variable} $X$ on a probability space $\left(\Omega,\mathcal{F},\mathbb{P}\right)$
is said to be normal iff (Def.) its distribution is normal, i.e.,
$\exists\,m\in\mathbb{R}$, $\sigma>0$ such that 
\[
\mathbb{P}\left(\left\{ \omega\in\Omega\:\big|\:a\leq X\left(\omega\right)\leq b\right\} \right)=\int_{a}^{b}\frac{1}{\sigma\sqrt{2\pi}}e^{-\frac{1}{2}\left(\frac{x-m}{\sigma}\right)^{2}}dx.
\]
\index{random variable}
\end{itemize}
\mainmatter

\part{Introduction and Motivation\label{part:intro}}

Below we outline the main areas covered inside the book. We offer
some tips for the reader, and conclude with a list of applications.

\section{Motivation}

More traditional books on Functional Analysis, and operators in Hilbert
space, tend to postpone applications till after all the big theorems
from the theory have been covered. The purpose of the present book
is to give students a tour of a selection of applications. We aim
to do this by first offering a crash course in theory topics tailor
made for the purpose (part \ref{part:topics}). In order to stress
the interplay between theory and applications (part \ref{part:app})
we have emphasized the traffic in both directions. We believe that
the multitude of new applications makes Functional Analysis both a
powerful, versatile, and timeless tool in mathematics.

A glance at existing books in Functional Analysis and related areas
(see list of reviews in the Appendix \ref{chap:books}) shows that
related books so far already display a rich variety, even if they
may have the same title \textquotedblleft Functional Analysis\textquotedblright{}
or \textquotedblleft Functional Analysis with a subtitle, or a narrowing
of the focus.\textquotedblright{}

Still \emph{the aims, and the contents} of these other books go in
a different directions than ours. One thing they have in common is
an emphasis on the Four Big Theorems in Functional Analysis, The \emph{Hahn-Banach
Theorem}, The \emph{Open Mapping Theorem}, The \emph{Uniform Boundedness
Principle}, The \emph{Closed Range Theorem}, and duality principles.

\index{Hahn-Banach theorem}

\index{open mapping theorem}

\index{uniform boundedness principle}

\index{closed range theorem}

\index{Four Big Theorems in Functional Analysis}

By contrast, we do as follows; rather we select a list of topics and
applications that acquire a degree of elegance when presented in a
functional analytic setting. There are several reasons for this different
approach, the main ones are as follows: 
\begin{enumerate}[label=(\roman{enumi}),ref=\roman{enumi}]
\item  The subject is ever changing on account of demands from neighboring
fields; 
\item Students come to our graduate functional analysis course with a diversity
of backgrounds, and a \textquotedblleft one-size fits all\textquotedblright{}
approach is not practical; 
\item Well-motivated students can pick up on their own what is needed from
the Four Big Theorems; 
\item Concentrating on the Four Big Theorems leaves too little time for
a variety of neighboring areas, both within mathematics, and in neighboring
sciences. 
\item Also the more traditional approach, beginning with the Four Big Theorems
is already in many existing books (see the Appendix \ref{chap:books}).
\end{enumerate}
A glance at the \emph{Table of Contents} will reflect our aim: beginning
with tools from Hilbert space in Chapters \ref{chap:basic} \& \ref{chap:lin},
but motivated by quantum physics; a preview of the Spectral Theorem
in \chapref{sp}; some basic tools from the theory of operator algebras
in \chapref{GNS}, with an emphasis on the Gelfand-Naimark-Segal (GNS)
construction; and stressing the many links between properties of states,
and the parallel properties of the representations, and the operator
algebras they generate.\index{Theorem!Spectral-}

In \chapref{cp}, and motivated by physics and harmonic\index{harmonic}
analysis, we discuss dilation theory. This is the general phenomenon
(pioneered by Stinespring and Arveson) of studying problems in an
initial Hilbert space by passing to an enlarged Hilbert space. \index{space!Hilbert-}

\chapref{bm} (Brownian motion), while different from the others,
still fits perfectly, and inviting application of the tools already
discussed in the first four chapters. The applications we cover in
\chapref{groups} are primarily to representations of groups and algebras.
\chapref{KS} is an application of theorems from Chapters \ref{chap:sp}-\ref{chap:GNS}
to the problem named after Kadison and Singer, now abbreviated the
KS-problem. It is 50 years old, is motivated by Dirac\textquoteright s
formulation of quantum physics (observables, states, and measurements);
and it was solved only a year ago (as of present).\index{observable}\index{measurement}

The last three chapters are, \ref{chap:ext}: s\emph{elfadjoint extensions},
\ref{chap:gLap}: \emph{graph-Laplacian}, and \ref{chap:RKHS}: \emph{reproducing
kernel Hilbert spaces} (RKHSs), and they are somewhat more technical,
but they are critical for a host of the questions inside the book.
Some readers may be familiar with this material already. If not, a
quick reading of Chapters \ref{chap:lin}, \ref{chap:ext}, and \ref{chap:gLap}
may be useful. Similarly, in the appendix, to help students orient
themselves, we give a birds-eye view of, in the order of 20 books
out there, all of which cover an approach to Functional Analysis,
and its many applications. \index{extension!-of operator}\index{extension!selfadjoint-}
\index{Laplacian!graph-}

\section{Key Themes in the Book: A bird\textquoteright s eye preview }

While each of our central themes has found book presentations, the
particular interconnection and applications that are the focus of
the present book, have not previously been explored in a textbook
form. To the extent they are in the literature at all, it will be
in the form of research papers.

\subsection{\label{sub:oph}Operators in Hilbert Space}

The notion of a \emph{Hilbert space}, is one of the most successful
axiomatic constructions in modern analysis. It was John von Neumann
who coined the term Hilbert space. While, historically, the concept
originated with problems from partial differential equations (PDE),
potential theory, quantum physics, and ergodic theory, it has since
found a host of other applications involving the part of functional
analysis dealing with infinite-dimensional function spaces; such areas
as: the study of unitary representations of groups (for example symmetry
groups from physics), complex function theory (Hardy spaces of holomorphic
functions), applications to probability, to stochastic processes,
to signal processing, to thermodynamics (heat transfer, ergodic theory.)
One reason the von Neumann-Hilbert axioms have proved especially successful
is their versatility in dealing with optimization problems arising
in the study of infinite-dimensional function spaces. This is so despite
the fact that the Hilbert space axioms themselves are formulated in
the abstract, independently of the particular context where they are
applied. Specifically, the axioms entail a given vector space $\mathscr{H}$,
equipped with an inner product (part of the axiom system), which in
turn induces a norm. The last axiom is that $\mathscr{H}$ must be
complete with respect to this norm.

With this one then proceeds to devise a host of coordinate systems,
orthonormal bases (ONB). A noteworthy family of ONBs of more recent
vintage are wavelet bases.

Among more recent areas of application, we mention machine learning;
a sub-area of artificial intelligence. In its current version, machine
learning models are formulated in the setting of reproducing kernel
Hilbert spaces (RKHS); see \chapref{RKHS} below, and \cite{SZ07}.
Indeed, in modern machine learning theory, the RKHSs play a critical
role in the construction of optimization algorithms. A second use
of RKHSs is in the solution of maximum-likelihood problems from probability
theory.

As for the study of linear transformations (operators), our present
dual emphasis will be unbounded operators, and non-commutativity.
Specifically, we study systems of densely defined linear operators.
A key motivation for this emphasis is again quantum mechanics: Indeed
quantum mechanical observables (momentum, position, energy, etc) correspond
to non-commuting selfadjoint unbounded operators in Hilbert space. 

The first two Hilbert spaces most students encounter are $l^{2}$
and $L^{2}\left(\mathbb{R}\right)$:

\needspace{3\baselineskip}
\begin{itemize}
\item[] $l^{2}$: sequences $\left(x_{n}\right)_{n=1}^{\infty}$ such that
\[
\left\Vert x\right\Vert _{l^{2}}^{2}=\sum_{n=1}^{\infty}\left|x_{n}\right|^{2}<\infty.
\]

\item[] $L^{2}\left(\mathbb{R}\right)$: measurable functions on $\mathbb{R}$
such that 
\[
\left\Vert f\right\Vert _{L^{2}}^{2}=\int_{\mathbb{R}}\left|f\left(x\right)\right|^{2}dx<\infty.
\]

\end{itemize}

These two examples serve to illustrate the axiom system for Hilbert
space which we shall study in  \secref{Hilbert} below. 

The norms for $l^{2}$ and for $L^{2}\left(\mathbb{R}\right)$ come
from associated inner products $\left\langle \cdot,\cdot\right\rangle $,
for example, $\left\langle x,y\right\rangle _{l^{2}}=\sum_{n=1}^{\infty}\overline{x}_{n}y_{n}$,
$\forall x,y\in l^{2}$. The system of vectors $\delta_{1},\delta_{2},\cdots$
in $l^{2}$, given by 
\[
\delta_{k}\left(n\right)=\delta_{k,n}=\begin{cases}
1 & \mbox{if \ensuremath{n=k}}\\
0 & \mbox{if \ensuremath{n\neq k}}
\end{cases}
\]
satisfies $\left\langle \delta_{k},\delta_{l}\right\rangle _{l^{2}}=\delta_{k,l}$
(the orthonormality property); and, for all $x=\left(x_{n}\right)_{1}^{\infty}\in l^{2}$,
we have
\[
\lim_{N\rightarrow\infty}\left\Vert x-\sum_{k=1}^{N}x_{k}\delta_{k}\right\Vert _{l^{2}}^{2}=\lim_{N\rightarrow\infty}\sum_{k=N+1}^{\infty}\left|x_{k}\right|^{2}=0;
\]
the second property (called ``\emph{total}'') that orthonormal bases
(ONBs) have.

One naturally wonders ``what are analogous ONBs for the second mentioned
Hilbert space $L^{2}\left(\mathbb{R}\right)?$'' And we shall turn
to this question also in \chapref{basic} below: There are two classes,
(i) \emph{special functions}, of which the best known are the Hermite
functions (\secref{bdd}); and (ii) \emph{wavelet-bases} (Sections
\ref{sec:Hilbert} and \ref{sec:end}).

The chapters of special relevance to these topics are: Chapters \ref{chap:basic},
\ref{chap:lin}, \ref{chap:KS}, \ref{chap:ext}, and \ref{chap:RKHS}.
Sections of special relevance include \ref{sec:proj}, \ref{sec:multi},
\ref{sec:normal}, \ref{sec:MPVM}, \ref{sec:sa}, \ref{sec:Cayley},
and \ref{sec:fried}.

\subsection{Multivariable Spectral Theory}

In this setting we are dealing with more than one operator at a time.
A host of applications naturally present themselves, again with the
same applications as mentioned in sect \ref{sub:oph} above. From
in the late nineteen thirties we have the study of selfadjoint algebras
in the works of Murray-von Neumann and of Gelfand-Naimark. Later,
this was followed up with systematic studies of non-selfadjoint algebras;
e.g., the work of Kadison-Singer. Other studies of multivariate operator
theory emphasize analogues of analyticity, both in the commutative
as well as in non-commutative settings. It has had remarkable successes,
including applications in other areas of mathematics such as complex
and algebraic geometry, and non-commutative geometry. In the multivariable
case, some researchers consider either $n$-tuples of operators, or
representations of algebras with generators and relations; while others
have adopted the language of Hilbert modules; for example, modules
over algebras of holomorphic functions, polynomials or entire functions
depending on the given number $n$ (commuting) complex variables. 

Historically, the first important multivariable problem in operator
theory was perhaps the relations of Heisenberg for a pair of linear
operators $P$ and $Q$ with dense domain $\mathscr{D}$ in a fixed
Hilbert space $\mathscr{H}$. The relations require that
\begin{equation}
PQf-QPf=-i\,f\label{eq:hr}
\end{equation}
holds for all $f\in\mathscr{D}$. 

By now (\ref{eq:hr}) is well understood, but there are many subtle
points; all of which make important connections to what we call multivariable
spectral theory; for examlpe (\ref{eq:hr}) does not have solutions
for bounded operators in $\mathscr{H}$.

We shall also consider a variety of multivariable systems of bounded
operators; -- in this case, it is usually in the setting of non-normal
operators (so in particular non-selfadjoint), for example for: (i)
finite sets of bounded commuting operators in a fixed Hilbert space
$\mathscr{H}$; (ii) algebras $\mathfrak{A}$ of operators on $\mathscr{H}$
such that the pair $\left(\mathfrak{A},\mathscr{H}\right)$ forms
a module; and (iii) finite systems of isometries in some Hilbert space
$\mathscr{H}$; and (iv) sets of isometries subject to the added condition
that the ranges forms a system of orthogonal subspaces of $\mathscr{H}$
with sum equal to $\mathscr{H}$. The relations on sets of isometries
described in (iv) are called the Cuntz relations (see sect \ref{sec:ov})
and \cite{MR0467330,BJ02,BJO04}. The Cuntz relations correspond to
representations of a $C^{*}$-algebra, called the Cuntz-algebra. It
has many applications, some of which will be studied, see e.g., sect
\ref{sec:egCalgebras}. A good reference for (i) is \cite{Arv98}.
\index{Cuntz relations}

The chapters of special relevance to these topics are: Chapters \ref{chap:cp},
\ref{chap:groups}, and \ref{chap:ext}. Sections of special relevance
include  \ref{sec:cpgns}, \ref{sec:sspring}, \ref{sec:ind}, \ref{sec:Heisenberg},
\ref{sec:dreprep}, and \ref{sec:fried}.

\subsection{\label{sec:na}Noncommutative Analysis}

The above multivariable settings are part of a wider theme: \emph{noncommutative
analysis}, a field which extends (classical commutative) Fourier analysis.
This began with the study of locally compact groups from physics,
and their unitary representations. The case of compact groups encompasses
the Peter-Weyl theorem from the 1920s, but needs from number theory
(mathematics), and from relativistic quantum physics, have dictated
extensions to non-compact (and non-commutative) groups, typically
Lie groups. \index{noncommutative!-analysis}

A more recent area of noncommutative analysis is the study of \emph{free
probability}, which we shall only touch on tangentially inside the
book. It is an exciting and new, rapidly growing, research direction;
with new advances in theory as well as in applications. Fortunately,
there are already nice and accessible book treatments, see e.g., \cite{Spe11}
and the sources cited there. In free probability, we study systems
of non-commutative random variables. As stochastic processes, they
are not Gaussian. Rather the notion of free independence dictates
the semicircle-law (not the Gaussian distribution). The rigorous study
of free probability entails such operator algebraic notions as free
products. We emphasize that the important new notion of free independence
is dictated by non-commutativity, and that it generalizes the more
familiar notion of independence which was used previously in probability.
The subject was initiated by Dan Voiculescu in the 1980ties. Its applications
up to now include: random matrix theory, representations of symmetric
groups, large deviations of stochastic processes, and quantum information
theory.

We use the term \textquotedblleft noncommutative analysis\textquotedblright{}
more broadly than the related one, \textquotedblleft noncommutative
geometry.\textquotedblright{} The latter owes much to the pioneering
work of Alain Connes, see e.g., \cite{MR2382238}. \sindex[nam]{Connes, A., (1947-)}\index{Connes, A.}In
broad outline, it covers the role von Neumann algebra theory plays
in noncommutative considerations in geometry and in quantum physics
(the Standard Model);\index{Standard Model} in noncommutative metric
theory and spaces, noncommutativity in topology, spectral triples,
differential geometry, cyclic cohomology, cyclic homology, K-theory,
and M-theory. In more detail, noncommutative geometry\index{noncommutative!-geometry}
(NCG) is concerned with a geometric approach to the construction of
spaces that are locally presented via noncommutative algebras of operators.
This is the framework of, what in physics, is referred to as \textquotedblleft local
quantum field theory.\textquotedblright{} The prime applications of
NCG are to particle physics where A. Connes has developed a noncommutative
standard model. Some of the other successes of NCG include extensions
of known topological invariants to formal duals of noncommutative
operator algebras. Via a Connes-Chern character map, this has led
to the discovery of a new homology theory of noncommutative operator
algebras; and to a new non-commutative theory of characteristic classes;
and to generalizations of the classical index theorems. \index{quantum field theory}

The Standard Model of particle physics deals with the electromagnetic,
weak, and strong nuclear interactions, and with classifications of
all the known subatomic particles, the \textquotedbl{}theory of almost
everything.\textquotedbl{} It received a boost in the mid-1970s after
an experimental confirmation of the existence of quarks; and later
of the tau neutrino, and the Higgs boson (2013).

Helpful references here are \cite{Arv76,MR1141333,BR79,DJ08,Gli61,JM84,Jor94,Jor11,MR0374334,MR548728}.

The chapters of special relevance to these topics are: Chapters \ref{chap:GNS},
\ref{chap:cp}, and \ref{chap:groups}. Sections of special relevance
include \ref{sec:sr}, \ref{sec:KMil}, \ref{sec:sspring}, \ref{sec:end},
\ref{sec:ind}, \ref{sec:dreprep}, and \ref{sec:multirep}.

\subsection{Probability}

Probability theory originated with the need for quantification of
uncertainty, as it arises for example in quantum physics, and in financial
markets. In the 1930ties, Kolmogorov's formulated precise mathematical
axioms of probability space $\Omega$, sample points, events as specified
subsets, in a prescribed sigma-algebra $\mathcal{F}$ of subsets of
$\Omega$, and a probability measure, defined on $\mathcal{F}$. 

Our present focus will be a subclass of stochastic processes, the
Gaussian processes, especially those which are derived from stochastic
integration defined relative to Brownian motion. 

Brownian motion is the simplest of the continuous-time stochastic
(meaning probabilistic) processes. It is a limit of simpler stochastic
processes going by the name random walks; a fact which reflects the
universality of the normal distribution, the Gaussians.

It is not an accident that we have focused on problems from quantum
physics and from probability. With some over simplification, it is
fair to say that Hilbert's 6th problem asked for a mathematical rigorous
treatment of these two areas. In 1900, when Hilbert formulated his
23 problems, these two areas did not yet have mathematically rigorous
foundations. 

The topic from probability that shall concern us the most is that
of Brownian motion. In a nutshell, a Brownian motion may be thought
as this way: There is a probability measure $\mathbb{P}$ on a sigma-algebra
of subsets of the continuous functions $\omega$ on $\mathbb{R}$
such that
\[
B_{t}\left(\omega\right)=\omega\left(t\right),\quad t\in\mathbb{R},\:\omega\in C\left(\mathbb{R}\right)
\]
satisfy a number of axioms of which we mention here only that for
each $t\in\mathbb{R}$, $B_{t}$ has a Gaussian distribution relative
to $\mathbb{P}$ such that
\[
\int_{C\left(\mathbb{R}\right)}\left|B_{t}\left(\omega\right)-B_{s}\left(\omega\right)\right|^{2}d\mathbb{P}\left(\omega\right)=\left|t-s\right|,\quad s,t\in\mathbb{R},
\]
and 
\[
\int_{C\left(\mathbb{R}\right)}B_{t}\left(\omega\right)d\mathbb{P}\left(\omega\right)=0,\quad\forall t\in\mathbb{R}.
\]

Helpful references here are \cite{MR2966130,MR3231624,MR0112604,MR2053326,Ito06,Nel67,parthasarathy1982probability,Sla03}.

The chapters of special relevance to these topics are: Chapters \ref{chap:bm},
and \ref{chap:RKHS}. Sections of special relevance include \ref{sec:dbm},
\ref{sec:stoch}-\ref{sec:types}, and Figures \ref{fig:geobm}, \ref{fig:geobm2},
\ref{fig:cyl}, \ref{fig:expExpSp}, and \ref{fig:expExt}.

\subsection{Unitary Representations}

An early motivation (see also \secref{na} above) is work of J. von
Neumann and I.E. Segal. They showed that, if $G$ is a locally compact
unimodular group such that the associated von Neumann group algebra
is of type I, then the regular representation of $G$, acting on the
Hilbert space $L^{2}\left(G\right)$ relative to Haar measure, as
a unitary representation, is a direct integral of irreducible unitary
representations (\textquotedblleft irreps\textquotedblright{} for
short.) This leads to a notion of a unitary dual for $G$, defined
as the set of equivalence classes (under unitary equivalence) of such
representations, the \textquotedblleft irreps.\textquotedblright{} 

But for general locally compact groups, including countable discrete
groups, the von Neumann group algebra typically is not of type I and
the regular translation-representation of $G$ cannot be expressed
in terms of building blocks of \textquotedblleft irreps.\textquotedblright{} 

The applications of our present results on unitary representations
will include those discussed above, so in particular, applications
to quantum physics, and to probability, especially to Gaussian stochastic
processes.

The chapters of special relevance to these topics are: Chapters \ref{chap:lin},
\ref{chap:GNS}, and especially \ref{chap:groups}. Sections of special
relevance include \ref{sec:dga}, \ref{sec:char}, \ref{sec:gns},
\ref{sec:NHspace}, \ref{sec:abelianc}, \ref{sec:ureprep}, \ref{sec:gaarding},
 and \ref{sec:dreprep}.

\section{Note on Cited Books and Papers }

For readers looking for references on the foundations, our suggestions
are as follows: Operators in Hilbert space: \cite{Arv76,Arv72}. Quantum
mechanics: \cite{OR07,Wei03,GG02}, and \cite{MR1939631,MR965583,MR675039}.
Non-commutative functional analysis and algebras of operators: \cite{BJKR84,BR81,BR79}.
Unitary representations of groups: \cite{Mac92,Mac85,Mac52}. 

\index{representation!unitary}

In our use of citations we adopted the following dual approach. Inside
the chapters, as the material is developed, we include citations to
key sources that we rely on; -- but this is done sparingly so as not
to interrupt the narrative too much.

To remedy sparse citations inside chapters, and, in order to help
the reader orient herself in the literature, each of the 11 chapters
concludes with a little bibliographical section, summarizing papers
and books of special relevance to the topic inside the text. Thus
there is a separate list of citations which concludes each chapter.
Readers who do not find a particular citation inside the chapter itself
will likely be able to locate it from the \emph{end-of-chapter-list}.

\section{Reader Guide}

\begin{flushleft}
Below we explain chapter by chapter how the six areas in \tabref{lsp}
are covered. 
\par\end{flushleft}

\renewcommand{\arraystretch}{2}

\begin{table}
\begin{tabular}{|c|>{\centering}p{0.4\textwidth}|>{\centering}p{0.4\textwidth}|}
\hline 
 & Subject & Example\tabularnewline
\hline 
\hline 
\enskip{}A\enskip{} & analysis & $f\left(x\right)-f\left(0\right)=\int_{0}^{x}f'\left(y\right)dy$\tabularnewline
\hline 
B & dynamical systems & functions on fractals, Cantor set, etc.\tabularnewline
\hline 
C & PDE & Sobolev spaces\tabularnewline
\hline 
D & numerical analysis & discretization\tabularnewline
\hline 
E & measures / probability theory & probability space $\left(\Omega,\mathcal{F},\mathbb{P}\right)$\tabularnewline
\hline 
F & quantum theory & Hilbert spaces of quantum states\tabularnewline
\hline 
\end{tabular}

\protect\caption{\label{tab:lsp}Examples of Linear Spaces: Banach spaces, Banach algebras,
Hilbert spaces $\mathscr{H}$, linear operators act in $\mathscr{H}$.
\index{Banach space}\index{space!Banach-} \index{Sobolev space}\index{space!Sobolev-}}
\end{table}

\index{Cantor!-set}

\renewcommand{\arraystretch}{1}

Ch 1: Areas A, E, F.

Ch 2: Areas A, C, F.

Ch 3: Areas B, C, E, F.

Ch 4: Areas A, B, E, F.

Ch 5: Areas A, E, F.

Ch 6: Area E.

Ch 7: Areas D, F.

Ch 8: Areas E, F.

Ch 9: Areas B, C, D.

Ch 10: Areas A, F.

Ch 11: Areas A, B, C, D, E.

\vspace{1em}

\begin{flushleft}
In more detail, the six areas in \tabref{lsp} may be fleshed out
as follows:
\par\end{flushleft}

Examples of subjects within \emph{area A} include measure theory,
transforms, construction of bases, Fourier series, Fourier transforms,
wavelets, and wavelet transforms, as well as a host of operations
in analysis. 

Subjects from \emph{area B} include solutions to ordinary differential
equations (ODEs), and the output of iteration schemes, such as the
Newton iteration algorithm. Also included are ergodic theory; and
the study of fractals, including harmonic analysis on fractals.

\emph{Area C} encompasses the study of the three types of linear PDEs,
elliptic, parabolic and hyperbolic. Sample questions: weak solutions,
\emph{a priori} estimates, diffusion equations, and scattering theory.

\emph{Area D} encompasses discretization, algorithms (Newton etc),
estimation of error terms, approximation (for example wavelet approximation,
and the associated algorithms.) 

\emph{Area E} encompasses probability theory, stochastic processes
(including Brownian motion), and path-space integration. 

Finally \emph{area F} includes the theory of unbounded operators in
Hilbert space, the three versions of the Spectral Theorem, as well
as representations of Lie groups, and of algebras generated by the
commutation relations coming from physics.\index{commutation relations}

\section{A Word About the Exercises}

All the chapters have exercises. The topics in the last three chapters
are more specialized, and exercises seem less natural there. The purpose
of the exercises is to improve and facilitate the use of the book
in courses; -- to help students and instructors. There is a total
of 147 exercises. To help with classroom use, we have listed them
in the back, numbering chapter-by-chapter. Each exercise is given
a name identification. Here is a sample: Exercise \exerref{laxmilgram}
(Lax-Milgram), \ref{exer:haar} (The Haar wavelet), \ref{exer:resolvent}
(the resolvent identity), \ref{exer:powers} (Powers-Størmer), \ref{exer:rp1}
(time-reflection), \ref{exer:prod} (extreme measures), \ref{exer:mult-1}
(multiplicity), \ref{exer:PW} (a formula from Peter-Weyl), and so
on; \dots{} \ref{exer:zego} (Szegö-kernel).

The degree of difficulty of the exercises varies from one to the next,
some are relatively easy; for example, serving to give the reader
a chance to practice definitions or new concepts; -- and some are
quite difficult. But of the exercises all interact naturally with
the topics developed in the various chapters. This is why we have
integrated them into the development of the topics, chapter for chapter.
And this is also why some chapters have many exercises, such as \chapref{basic}
with a total of 38 exercises; -- \chapref{sp} has 12 exercises; and
\chapref{GNS} has 44 exercises in all. In all of the chapters, we
have mixed and interspaced the placement of exercises with the central
themes: some exercises supplement examples, and some theorems, within
each chapter.

There are two lists after the Appendices, a \emph{List of Exercises},
and a \emph{list of all the figures}. The second should help readers
with cross-references; and the first with use of the exercises in
course-assignments. 

\emph{The Appendices} themselves serve to aid readers navigate the
book-literature. Appendix \ref{chap:books} includes \emph{telegraphic
reviews}, and Appendix \ref{chap:bios} is a collection of \emph{biographical
sketches} of the pioneers in the subject.

\section{\label{sec:app}List of Applications}

Part of the discussion below will make use of terminology from neighboring
areas, such as physics, engineering, and statistics. For readers who
might be encountering this for the first time, we have a terminology
section in the back. It is a section in the Appendix \ref{chap:term},
called \textquotedblleft Terminology from neighboring areas.\textquotedblright{}
The Appendix also includes other lists: a list of telegraphic reviews
of related books; biographical sketches; diagrams with lines illustrating
interconnections between disparate areas; a list of figures, and a
list of all the Exercises inside the text. Each Exercise is given
a descriptive name.

Each of the chapters is illustrated with examples and applications.
A recurrent theme is the important notion of \emph{positive definite
functions}, and their realization in Hilbert space. Applications to
\emph{Wiener measure} and path-space are included in Chapters \ref{chap:basic},
\ref{chap:bm}, and \ref{chap:RKHS}.\index{space!Hilbert-}

More applications, starting in \chapref{basic} are: (i) a variance
formula for the Haar-wavelet basis in $L^{2}(0,1)$; (ii) a formula
for perturbations of diagonal operators; and (iii) the $\infty\times\infty$
matrix representation of \emph{Heisenberg\textquoteright s commutation
relations} in the case of the canonical pair, momentum and position
operators from quantum mechanics.

\index{operators!momentum-}

\index{operators!position-}

The case of Heisenberg\textquoteright s commutation relations motivates
the need for a systematic study of \emph{unbounded operators} in Hilbert
space. This is started in  \chapref{lin}, and resumed then systematically
in Chapters \ref{chap:sp} (the \emph{Spectral Theorem}), and \ref{chap:ext}
(the theory of \emph{extensions of symmetric operators with dense
domain}; -- von Neumann indices, and All That). Both our study of
selfadjoint operators and normal operators, and their spectral theory,
throughout the book is motivated by the axioms from quantum theory:
observables, states, measurements, and the uncertainty principle.
Our systematic treatment of \emph{projection valued measures}, and
\emph{quantum states}, in \chapref{sp} is a case in point. \index{normal operator}
\index{selfadjoint operator} \index{axioms}

This goes for our theme in \chapref{GNS} as well, the \emph{Gelfand-Naimark-Segal}
(GNS) representation. \emph{Quantum states} must be realized in Hilbert
space, but what is the relevant Hilbert space, when a \emph{quantum
observable} is prepared in a state? To answer this we must realize
the observables as selfadjoint operators affiliated with a suitable
$C^{*}$-algebra, say $\mathfrak{A}$, or von Neumann algebra. States
on these algebras then become positive linear functionals. The GNS
construction is a device for constructing representation in Hilbert
space for every state, defined as a \emph{positive linear functional}
on $\mathfrak{A}$. In this construction, the pure states are matched
up with irreducible representations.

In \secref{NHspace}, our application is to the subject of \textquotedblleft \emph{reflection-positivity}\textquotedblright{}
from quantum physics. This notion came up first in a renormalization\index{renormalization}
question in physics: \textquotedblleft How to realize observables
in relativistic quantum field theory (RQFT)?\textquotedblright{} \index{affiliated with}

The material in \chapref{cp} has applications to signal processing;
-- to the construction of sub-band filters, and filter banks. These
applications are discussed in \chapref{cp}; -- included as one of
the applications of a certain family of representations of the Cuntz
relations. Other applications of these representations include wavelet
filters. \index{filter bank}\index{signal}\index{signal processing}

\section{Groups and Physics}
\begin{quote}
\textquotedblleft The importance of group theory was emphasized very
recently when some physicists using group theory predicted the existence
of a particle that had never been observed before, and described the
properties it should have. Later experiments proved that this particle
really exists and has those properties.\textquotedblright{} 

--- Irving Adler\vspace{2em}
\end{quote}
Recall that in RQFT, the symmetry group is the Poincaré group, but
its physical representations are often illusive. Starting with papers
by Osterwalder-Schrader in the 1970ties (see e.g., \cite{MR0376002,MR887102,JO00}),
it was suggested to instead begin with representations of the Euclidian
group, and then to get to the Poincaré group through the back door,
via an analytic continuation (a c-dual group construction), and a
\emph{renormalization}. This lead to a systematic study of renormalizations
for the Hilbert space of quantum states. The \textquotedblleft c-dual\textquotedblright{}
here refers to an \emph{analytic continuation} which links the two
groups. This in turn is accomplished with the use of a certain reflection,
and a corresponding change in the inner product. In a simplified summary,
the construction is as follows: Starting with the inner product in
the initial Hilbert, say $\mathscr{H}$, and a unitary representation
admitting a reflection $\mathcal{J}$, we then pass to a certain invariant
subspace of $\mathscr{H}$, and use $\mathcal{J}$ in the definition
of the new inner product. The result is a physical energy operator
(dual of time) with the correct positive spectrum for the relativistic
problem, hence \textquotedblleft \emph{reflection-positivity}.\textquotedblright{}
The invariant subspace refers to invariance only in a positive time
direction. All of this is presented in \secref{NHspace}, and illustrated
with an example. \index{Osterwalder-Schrader}\index{renormalization}

\index{representation!unitary}

Chapter \ref{chap:cp} deals with the same theme; only there the states
are \emph{operator valued}. From the theory of Stinespring and Arveson
we know that there is then a different positivity notion, \emph{complete
positivity} (CP).

Among the Hilbert spaces we encounter are $L^{2}$ spaces of random
variables on a probability space $\left(\Omega,\mathcal{F},\mathbb{P}\right)$.
The case of Brownian motion is studied in \chapref{bm}, and again
in \chapref{RKHS}.\index{space!Hilbert-}

In  \chapref{groups}, we introduce families of \emph{unitary representations}
of groups, and $*$-representations of algebras; each one motivated
by an application from physics, or from signal-processing. We are
stressing examples as opposed to general theory. \index{signal}

\chapref{KS} is devoted to the \emph{Kadison-Singer} problem (KS).
\index{Kadison-Singer problem}It is a problem from operator algebras,
but originating with Dirac\textquoteright s presentation of quantum
mechanics. By choosing a suitable orthonormal basis (ONB) we may take
for Hilbert space the sequence $l^{2}\left(\mathbb{N}\right)$ space,
square-summable sequences. Dirac was interested in the algebra $\mathscr{B}\left(l^{2}\left(\mathbb{N}\right)\right)$
of all bounded operators in $l^{2}\left(\mathbb{N}\right)$. But with
the $\infty\times\infty$ matrix representation for elements in $\mathscr{B}\left(l^{2}\left(\mathbb{N}\right)\right)$,
we can talk about the \emph{maximal abelian subalgebra} $\mathscr{D}$
of all diagonal operators in $\mathscr{B}\left(l^{2}\left(\mathbb{N}\right)\right)$.
Note $\mathscr{D}$ is just a copy of $l^{\infty}\left(\mathbb{N}\right)$.
The Dirac-KS question is this: \textquotedblleft Does every pure state
on $\mathscr{D}$ have a \emph{unique} pure-state extension to $\mathscr{B}\left(l^{2}\left(\mathbb{N}\right)\right)$?\textquotedblright{}
\index{state!pure-}

The problem was solved in the affirmative; just a year ago (see \cite{MSS15}),
and we sketch the framework for the KS problem. However the details
of the solution are far beyond the scope of our book.

The application in \chapref{gLap} is to potential theory of \emph{infinite
networks}; mathematically infinite graphs $G=\left(V,E\right)$, $V$
the specified set of vertices, and $E$ the edges. Our emphasis is
electrical networks, and the functions include energy, conductance,
resistance, voltage, and current. The main operator here is the so
called \emph{graph Laplacian}.

The new applications in  \chapref{RKHS} include \emph{scattering
theory},\emph{ learning theory} (as it is used in machine learning
and in pattern recognition.)\index{graph!network-}

\part{\label{part:topics}Topics from Functional Analysis and Operators
in Hilbert Space: a selection}

\chapter{Elementary Facts\label{chap:basic}}
\begin{quotation}
\textquotedblleft \dots{} the {[}quantum mechanical{]} observables
are operators on a Hilbert space. The algebra of operators on a Hilbert
space is noncommutative. It is this noncommutativity of operators
on a Hilbert space that provides a precise formulation of {[}Heisenberg\textquoteright s{]}
uncertainty principle: There are operator solutions to equations like
$pq-qp=1$. This equation has no commutative counterpart. In fact,
it has no solution in operators p,q acting on a finite dimensional
space. So if you're interested in the dynamics of quantum theory,
you must work with operators rather than functions and, more precisely,
operators on infinite dimensional spaces.\textquotedblright{} \sindex[nam]{Arveson, W., (1934-2011)}

--- William B. Arveson (1934-2011. The quote is from 2009.)\vspace{1em}\\
I received an early copy of Heisenberg's first work a little before
publication and I studied it for a while and within a week or two
I saw that the noncommutation was really the dominant characteristic
of Heisenberg's new theory. It was really more important than Heisenberg's
idea of building up the theory in terms of quantities closely connected
with experimental results. So I was led to concentrate on the idea
of noncommutation and to see how the ordinary dynamics which people
had been using until then should be modified to include it. 

--- P. A. M. Dirac\sindex[nam]{Dirac, P.A.M., (1902-1984)}\vspace{1em}

Problems worthy

\hspace{8pt}of attack 

prove their worth

\hspace{8pt}by hitting back. 

--- Piet Hein\sindex[nam]{Hein, P., (1905-1996)}\vspace{2em}
\end{quotation}
Below we outline some \emph{basic concepts}, \emph{ideas}, and \emph{examples}
which will be studied inside the book itself. While they represent
only a sample, and we favor the setting of Hilbert space, the details
below still tie in nicely with diverse tools and techniques not directly
related to Hilbert space.

The discussion below concentrates on topics connected to Hilbert space,
but we will also have occasion to use some other basic facts from
functional analysis; e.g., duality and Hahn-Banach. We have collected
those, in a condensed form, in an Appendix at the end of the chapter.

From linear algebra we know precisely what square matrices $M$ can
be diagonalized; the \emph{normal matrices}, i.e., $M^{*}M=MM^{*}$.
More precisely, a matrix is normal if and only if it is conjugate
to a diagonal matrix. More general square matrices don\textquoteright t
diagonalize, but they admit a Jordan form.

In the infinite dimensional case, -- while infinite matrices are useful,
the axiomatic setting of \emph{Hilbert space} and \emph{linear operators}
has proved more successful than an infinite matrix formulation; and,
following von Neumann and Stone, we will make precise the notion of
normal operators. Because of applications, the case of \emph{unbounded
operators} is essential. In separate chapters, we will prepare the
ground for this. \index{space!Hilbert-}

The\emph{ Spectral Theorem} (see \cite{Sto90, Yos95, Ne69, RS75, DS88b})
states that a linear operator $T$ (in Hilbert space) is \emph{normal},
i.e., $T^{*}T=TT^{*}$, if and only it is unitarily equivalent to
a \emph{multiplication operator} in some $L^{2}$ space, i.e., multiplication
by a measurable function, and the function may be unbounded. The implied
Hilbert space $L^{2}$ is with respect to some measure space, which
of course will depend on the normal operator $T$, given at the outset.
Hence the classification of \emph{normal operators} is equivalent
to the classification of \emph{measure spaces}; -- a technically quite
subtle problem. \index{normal operator}

There is a second (and equivalent) version of the Spectral Theorem,
one based on \emph{projection valued measures} (PVMs), and we will
present this as well. It is a powerful tool in the theory of unitary
representations of locally compact groups (see \chapref{groups} below),
and in a host of areas of pure and applied mathematics.

\index{Theorem!Spectral-}

\index{representation!unitary}\index{axioms}

It is natural to ask whether there is an analogue of the finite-dimensional
Jordan form; i.e., extending from finite to the infinite dimensional
case. The short answer is \textquotedblleft no,\textquotedblright{}
although there are partial results. They are beyond the scope of this
book.

In our first two chapters below we prepare the ground for the statement
and proof of the \emph{Spectral Theorem}, but we hasten to add that
there are several versions. In the bounded case, for \emph{compact
selfadjoint operators} (\secref{spcpt}), the analogue to the spectral
theorem from linear algebra is closest, i.e., eigenvalues and eigenvectors.
Going beyond this will entail an understanding of \emph{continuous
spectrum} (\secref{MPVM}), and of \emph{multiplicity theory} in the
measure theoretic category (\secref{mult}). \index{eigenvalue}

With a few exceptions, we will assume that all of the Hilbert spaces
considered are separable; i.e., that their orthonormal bases (ONBs)
are countable. The exceptions to this will include the $L^{2}$-space
of the Bohr completion of the reals $\mathbb{R}$. See \exerref{bohr}.

\section{\label{sec:topics}A Sample of Topics}
\begin{quote}
\textquotedblleft Too many people write papers that are very abstract
and at the end they may give some examples. It should be the other
way around. You should start with understanding the interesting examples
and build up to explain what the general phenomena are.\textquotedblright{} 

--- Sir Michael Atiyah\vspace{2em}
\end{quote}
Classical functional analysis is roughly divided into two branches,
each with a long list of subbranches:
\begin{itemize}
\item study of function spaces (Banach space, Hilbert space)
\item applications in physics, statistics, and to engineering
\end{itemize}
Within pure mathematics, it is manifested in the list below:
\begin{itemize}
\item representation theory of groups and algebras, among a long list of
diverse topics 
\end{itemize}

We will consider three classes of algebraic objects of direct functional
analytic relevance: (i) generators and relations; (ii) algebras, and
(iii) groups.

In the case of (i), we illustrate the ideas with the \emph{canonical
commutation relation} \index{Heisenberg, W.K.!commutation relation}\index{commutation relations}
\begin{equation}
PQ-QP=-i\,I,\quad i=\sqrt{-1}.\label{eq:cc}
\end{equation}
The objective is to build a Hilbert space such that the symbols $P$
and $Q$ are represented by unbounded essentially selfadjoint operators
(see \cite{RS75, Ne69, vN32a, DS88b}), each defined on a common dense
domain in some Hilbert space, and with the operators satisfying (\ref{eq:cc})
on this domain. (See technical points inside the present book, and
in the cited references.) \index{operators!adjoint-} \index{operators!essentially selfadjoint-}
\index{selfadjoint operator}

In class (ii), we consider both $C^{*}$-\emph{algebras} and \emph{von
Neumann algebras} (also called $W^{*}$-algebras); and in case (iii),
our focus is on \emph{unitary representations} of the group $G$ under
consideration. The group may be abelian or non-abelian, continuous
or discrete, locally compact or not. Our present focus will be the
case when $G$ is a \emph{Lie group}. In this case, we will study
its representations with the use of the corresponding \emph{Lie algebra}. 

\index{algebras!$C^{*}$-algebra}

\index{algebras!von Neumann algebra ($W^{*}$-algebra)}

\index{groups!abelian}

\index{groups!non-abelian}

\index{Lie!group}

\index{Lie!algebra}

\index{von Neumann, J.}

\index{essentially selfadjoint operator}
\begin{itemize}
\item $C^{*}$-algebras, von Neumann algebras
\end{itemize}

We will be considering $C^{*}$-, and $W^{*}$-algebras axiomatically.
In doing this we use the theorem by S. Sakai to the effect that the
$W^{*}$-algebras consist of the subset of the $C^{*}$-algebras that
are the dual of a Banach space. If the $W^{*}$-algebra is given,
the Banach space is called the pre-dual\index{pre-dual}. Representations
will be studied with the use of states, and we stress the theorem
of Gelfand, Naimark, and Segal (GNS) linking states with cyclic representations.

\index{Banach space!dual-}

\index{Banach space!pre-dual}

\index{space!Banach--}

\index{space!dual-}\index{axioms}
\begin{itemize}
\item wavelets theory
\end{itemize}

A wavelet is a special basis for a suitable $L^{2}$-space which is
given by generators and relations, plus self-similarity. Our approach
to wavelets will be a mix of functional analysis and harmonic\index{harmonic}
analysis, and we will stress a correspondence between a family of
representations of a particular $C^{*}$-algebra, called the Cuntz-algebra,
on one side and wavelets on the other.

\index{wavelet}

\index{algebras!Cuntz algebra}

\index{Cuntz-algebra}
\begin{itemize}
\item harmonic analysis
\end{itemize}

Our approach to harmonic analysis will be general, -- encompassing
anyone of a set of direct sum (or integral) decompositions. Further
our presentation will rely on representations.
\begin{itemize}
\item analytic number theory 
\end{itemize}

Our notions from analytic number theory will be those that connect
to groups, and representations; such as the study of automorphic forms,
and of properties of generalized zeta-functions; see e.g., \cite{MR2290770,MR2261099,MR960228,MR963118}.

Our brief bird\textquoteright s eye view of the topics above is only
preliminary, only hints; and most questions will be addressed in more
detail inside the book.

As for references, the literature on the commutation relations (\ref{eq:cc})
is extensive, and we refer to \cite{Sza04,Nel59,Fug82,Pou73}.

Some of the questions regarding the commutation relations involve
the subtle difference between (\ref{eq:cc}) itself vs its group version,
-- often referred to as the \emph{Weyl relations}, or the integrated
form. As for the other themes mentioned above, operator algebras,
math physics, wavelets and harmonic analysis, the reader will find
selected references to these themes at the end of this chapter.

A glance at the table of contents makes it clear that we venture into
a few topics at the cross roads of mathematics and physics; and a
disclaimer is in order. In the 1930s, David Hilbert encouraged the
pioneers in quantum physics to axiomatize the theory that was taking
shape then with the initial papers by Heisenberg, Schrödinger, Dirac.
Others like J. von Neumann joined into this program. These endeavors
were partly in response to \emph{Hilbert\textquoteright s Sixth Problem}
\cite{Wig76}; \index{Hilbert's six problem}\index{axioms}

\textquotedblleft Give a mathematical treatment of the axioms of physics\textquotedblright{} 

in his famous list of 23 problems announced in 1900 \cite{MR1557926}.
At this time, quantum physics barely existed. Max Planck\textquoteright s
hypothesis on discreteness of atomic energy-measurements is usually
dated a little after the turn of the Century. \index{measurement}

\textbf{Quantum mechanics} is a first quantized quantum theory that
supersedes classical mechanics at the atomic and subatomic levels.
It is a fundamental branch of physics that provides the underlying
mathematical framework for many fields of physics and chemistry. The
term \textquotedblleft quantum mechanics\textquotedblright{} is sometimes
used in a more general sense, to mean quantum physics.\index{Hilbert's six problem}

With hindsight, we know that there are considerable limitations to
the use of axioms in physics. While a number of important questions
in quantum physics have mathematical formulations, others depend on
physical intuition. Hence in our discussion of questions at the cross-roads
of mathematics and physics, we will resort to hand-waiving.
\begin{quotation}
\textquotedblleft For those who are not shocked when they first come
across quantum theory cannot possibly have understood it.\textquotedblright{} 

Niels Bohr, --- quoted in W. Heisenberg, Physics and Beyond (1971).
\sindex[nam]{Bohr, N., (1885-1962)}
\end{quotation}

\section{\label{sec:duality}Duality}

The ``functional'' in the name ``Functional Analysis'' derives
from the abstract notion of a \emph{linear functional}: Let $E$ be
a vector space over a field $\mathbb{F}$ (we shall take $\mathbb{F}=\mathbb{R}$,
or $\mathbb{C}$ below.) \index{functional}
\begin{defn}
\label{def:Ba1}A function $\varphi:E\rightarrow\mathbb{R}$ (or $\mathbb{C}$)
is said to be a \emph{linear functional}, if we require
\[
\varphi\left(u+\lambda v\right)=\varphi\left(u\right)+\lambda\varphi\left(v\right),\;\forall\lambda\in\mathbb{R},\:\forall u,v\in E.
\]
If $E$ comes with a topology (for example from a norm, or from a
system of seminorms), we will consider continuous linear functionals.
Occasionally, continuity will be implicit.
\end{defn}
\index{continuous linear functional}

\begin{defn}
\label{def:Ba2}The set of all continuous linear functionals is denoted
$E^{*}$, and it is called the \emph{dual space}. (In many examples
there is a natural identification of $E^{*}$ as illustrated in \tabref{dual}.)
\end{defn}

\begin{defn}
If $E$ is a normed vector space, and if it is complete in the given
norm, we say that $E$ is a \emph{Banach space}. \index{Banach space}\index{space!Banach-}\end{defn}
\begin{lem}
Let $E$ be a normed space with dual $E^{*}$. For $\varphi\in E^{*}$,
set 
\[
\left\Vert \varphi\right\Vert _{*}:=\sup_{\left\Vert x\right\Vert =1}\left|\varphi\left(x\right)\right|.
\]
Then $\left(E^{*},\left\Vert \cdot\right\Vert _{*}\right)$ is a Banach
space.\end{lem}
\begin{proof}
An exercise.
\end{proof}
Given a Banach space \emph{E}, there are typically three steps involved
in the discovery of an explicit form for the \emph{dual Banach space}
$E^{*}$. \tabref{dual} illustrates this in two examples, but there
are many more to follow; -- for example, the case when $E=$ the Hardy
space $\mathbb{H}_{1}$, or $E=$ the trace-class operators on a Hilbert
space. 

Moreover, the same idea based on a duality-pairing applies \emph{mutatis
mutandis}, to other topological vector spaces as well, for example,
to those from Schwartz\textquoteright{} theory of distributions. 

\index{distribution!Schwartz-}\index{space!Hilbert-}\index{space!Hardy-}

The three steps are as follows: 

\textbf{Step 1.} Given $E$, then first come up with a second Banach
space $F$ as a candidate for the dual Banach space $E^{*}$. (Note
that $E^{*}$ is so far, \emph{a priori}, only an abstraction.)\index{space!dual-}

\textbf{Step 2.} Set up a bilinear and non-degenerate pairing, say
$p$, between the two Banach spaces $E$ and $F$, and check that
$p\left(\cdot,\cdot\right)$ is continuous on $E\times F$. Rescale
such that 
\[
\left|p\left(x,y\right)\right|\leq\left\Vert x\right\Vert _{E}\left\Vert y\right\Vert _{F},\;\forall x\in E,\:y\in F.
\]
This way, via $p$, we design a linear and isometric embedding of
$F$ into $E^{*}$.

\textbf{Step 3.} Verify that the embedding from step 2 is \textquotedblleft onto\textquotedblright{}
$E^{*}$. If \textquotedblleft yes\textquotedblright{} we say that
$F$ \textquotedblleft is\textquotedblright{} the dual Banach space.
(Example, the dual of $\mathbb{H}_{1}$ is BMO \cite{MR0280994}.)

Examples of Banach spaces include: (i) $l^{p}$, $1\leq p\leq\infty$,
and (ii) $L^{p}\left(\mu\right)$ where $\mu$ is a positive measure
on some given measure space; details below.
\begin{example}
\label{exa:lp}$l^{p}$: all $p$-summable sequences. 

A sequence $x=\left(x_{k}\right)_{k\in\mathbb{N}}$ is in $l^{p}$
iff $\sum_{k\in\mathbb{N}}\left|x_{k}\right|^{p}<\infty$, and then
\[
\left\Vert x\right\Vert _{p}:=\left(\sum_{k\in\mathbb{N}}\left|x_{k}\right|^{p}\right)^{\frac{1}{p}}.
\]

\end{example}

\begin{example}
\label{exa:Lp}$L^{p}$: all $p$-integrable function with respect
to some fixed measure $\mu$.

Let $F:\mathbb{R}\rightarrow\mathbb{R}$ be monotone increasing, i.e.,
$x\leq y$ $\Longrightarrow$ $F\left(x\right)\leq F\left(y\right)$;
then there is a Borel measure $\mu$ on $\mathbb{R}$ (see \cite{Rud87})
such that $\mu\left((x,y]\right)=F\left(y\right)-F\left(x\right)$;
and $\int\varphi d\mu$ will be the limit of the Stieltjes sums: 
\[
\sum_{i}\varphi\left(x_{i}\right)\left(F\left(x_{i+1}\right)-F\left(x_{i}\right)\right),\;\mbox{where}\;x_{1}<x_{2}<\cdots<x_{n}.
\]
We say that $\varphi\in L^{p}\left(\mu\right)$ iff $\int_{\mathbb{R}}\left|\varphi\right|^{p}d\mu$
is well defined and finite; then 
\[
\left\Vert \varphi\right\Vert _{p}=\left(\int_{\mathbb{R}}\left|\varphi\left(x\right)\right|^{p}d\mu\left(x\right)\right)^{\frac{1}{p}}.
\]

By Stieltjes integral, $\int\left|\varphi\right|^{p}d\mu=\int\left|\varphi\right|^{p}dF$.\index{integral!Stieltjes-}
Here, we give the definition of $L^{p}\left(\mu\right)$ in the case
where $\mu=dF$, but it applies more generally. \index{space!$L^{p}$-}

For completeness of $l^{p}$ and of $L^{p}\left(\mu\right)$, see
\cite{Rud87}. \index{Borel measure}
\end{example}

\renewcommand{\arraystretch}{2}

\begin{table}
\begin{tabular}{|>{\centering}p{0.25\textwidth}|>{\centering}p{0.25\textwidth}|>{\centering}p{0.4\textwidth}|}
\hline 
$E$ & $E^{*}$ & how?\tabularnewline
\hline 
$l^{p}$, $1\leq p<\infty$ with $l^{p}$-norm & $l^{q}$, $\frac{1}{p}+\frac{1}{q}=1$ & $x=\left(x_{i}\right)\in l^{p}$, $y=\left(y_{i}\right)\in l^{q}$,
$\varphi_{y}\left(x\right)=\sum_{i}x_{i}y_{i}$\tabularnewline
\hline 
$C\left(I\right)$, $I=\left[0,1\right]$ with max-norm & signed Borel measures $\mu$ on $I$, of bounded variation & $\varphi_{\mu}\left(f\right)=\int_{0}^{1}f\left(x\right)d\mu\left(x\right)$,
$\forall f\in C\left(I\right).$\tabularnewline
\hline 
$C^{\infty}\left(\mathbb{R}\right)$, system of seminorms  & $\mathcal{E}'$ all Schwartz distributions $D$ on $\mathbb{R}$ of
compact support & $\varphi_{D}\left(f\right)=D$ applied to $f$, $f\in C^{\infty}\left(\mathbb{R}\right).$\tabularnewline
\hline 
\end{tabular}

\protect\caption{\label{tab:dual}Examples of dual spaces.}
\end{table}

\renewcommand{\arraystretch}{1}

\begin{figure}
\includegraphics[width=0.4\textwidth]{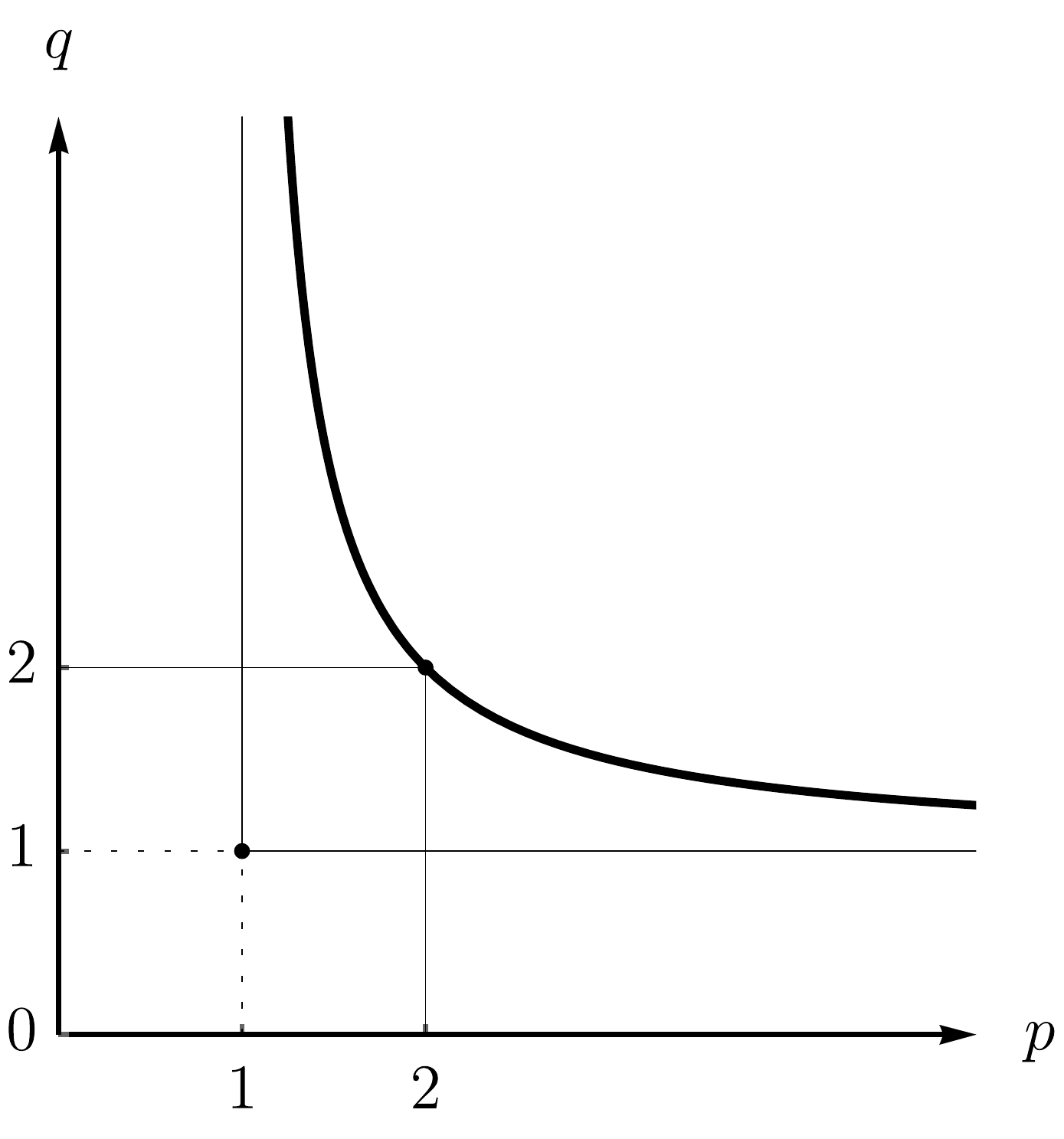}

\protect\caption{\label{fig:conj}Dual exponents for the $L^{p}$ spaces, $\frac{1}{p}+\frac{1}{q}=1$.}
\end{figure}

\index{distribution!Schwartz-}
\begin{rem}
At the foundation of analysis of $L^{p}$-spaces (including $l^{p}$
for the case of counting-measure) is Hölder's inequality; see e.g.,
\cite[ch 3]{Rud87}. Recall conjugate pairs $p,q\in[1,\infty)$, $\frac{1}{p}+\frac{1}{q}=1$,
or equivalently $p+q=pq$; see \figref{conj}. 

\index{inequality!Hölder}

We present Hölder's inequality without proof: Fix a measure space
$\left(X,\mathcal{F},\mu\right)$. If $p,q$ are conjugate, $1<p<\infty$,
then for measurable functions $f,g$ we have:
\begin{equation}
\left|\int_{X}fgd\mu\right|\leq\left(\int_{X}\left|f\right|^{p}d\mu\right)^{\frac{1}{p}}\left(\int_{X}\left|g\right|^{q}d\mu\right)^{\frac{1}{q}}.\label{eq:hd1}
\end{equation}
If $p=1$, $q=\infty$, and we have:
\begin{equation}
\left|\int_{X}fgd\mu\right|\leq\left(\int_{X}\left|f\right|d\mu\right)\mbox{ess sup}_{x\in X}\left|g\left(x\right)\right|;\label{eq:hd2}
\end{equation}
where 
\[
\left\Vert g\right\Vert _{\infty}:=\mbox{ess sup}\left|g\right|
\]
denotes essential supremum, i.e., neglecting sets of $\mu$-measure
zero.
\end{rem}

The following result is basic in the subject.
\begin{thm}[Hahn-Banach]
\label{thm:banach} Let $E$ be a Banach space, and let $x\in E\backslash\left\{ 0\right\} $,
then there is a $\varphi\in E^{*}$ such that $\varphi\left(x\right)=\left\Vert x\right\Vert $,
and $\left\Vert \varphi\right\Vert _{E^{*}}=1$. \end{thm}
\begin{rem}
In Examples \ref{exa:lp} and \ref{exa:Lp} above, i.e., $l^{p}$
and $L^{p}\left(\mu\right)$, it is possible to identify the needed
elements in $E^{*}$. But the power of Theorem \ref{thm:banach} is
that it yields existence for all Banach spaces, i.e., when $E$ is
given only by the axioms from Definitions \ref{def:Ba1}-\ref{def:Ba2}.\index{space!Banach-}\index{axioms}\end{rem}
\begin{defn}
The weak-$*$ topology on $E^{*}$ is the weakest topology which makes
all the linear functionals 
\[
E^{*}\ni l\longrightarrow l\left(x\right)\in\mathbb{C}
\]
continuous, as $x$ ranges over $E$.\end{defn}
\begin{xca}[weak-$*$ neighborhoods]
\label{exer:wnhd}\myexercise{weak-$*$ neighborhoods} Show that
the neighborhoods of $0$ in $E^{*}$ have a basis of open sets $\mathscr{N}$
indexed as follows:

Let $\epsilon\in\mathbb{R}_{+}$, $n\in\mathbb{N}$, and $x_{1},\ldots,x_{n}\in E$,
and set 
\[
\mathscr{N}_{\epsilon,x_{1},\ldots,x_{n}}:=\left\{ l\in E^{*}\::\:\left|l\left(x_{i}\right)\right|<\epsilon,\;i=1,\cdots,n\right\} .
\]

\end{xca}
Terminology. The subsets of $E^{*}$ in \exerref{wnhd} are often
called cylinder sets. They form a basis for the weak-$*$ topology.
They also generate a sigma algebra of subsets of $E^{*}$, often called
the cylinder sigma algebra. We will be using it in Sections \ref{sec:pspace}
(pg. \pageref{sec:pspace}), \ref{sec:dbm} (pg. \pageref{sec:dbm}),
and \ref{sec:stoch} (pg. \pageref{sec:stoch}) below.
\begin{xca}[weak-$*$ vs norm]
\label{exer:wnorm}\myexercise{weak-$*$  vs norm}Let $1<p\leq\infty$
be fixed. Set $l^{p}=l^{p}\left(\mathbb{N}\right)$, and show that
$\left\{ x\in l^{p}\::\:\left\Vert x\right\Vert _{l^{p}}\leqq1\right\} $
is weak-$*$ compact, but not norm-compact.

\uline{Hint}: By weak-$*$, we refer to $l^{p}=\left(l^{q}\right)^{*}$,
$\frac{1}{p}+\frac{1}{q}=1$. 
\end{xca}

\begin{xca}[Be careful with weak-$*$ limits.]
\myexercise{Be careful with weak-$*$ limits}Settings as in the previous
exercise, but now with $p=2$. Let $\left\{ e_{k}\right\} _{k\in\mathbb{N}}$
be the standard ONB in $l^{2}$, i.e., 
\begin{equation}
e_{k}\left(i\right)=\delta_{i,k},\;\forall i,k\in\mathbb{N}.\label{eq:dd1}
\end{equation}
Show that $0$ in $l^{2}$ is a weak $*$-limit of the sequence $\left\{ e_{k}\right\} _{k\in\mathbb{N}}$.
Conclude that $\left\{ x\in l^{2}\::\:\left\Vert x\right\Vert _{2}=1\right\} $
is \uline{not} weak-$*$ closed.

\uline{Hint}: By Parseval, we have, for all $x\in l^{2}$, 
\[
\left\Vert x\right\Vert _{2}^{2}=\sum_{k\in\mathbb{N}}\left|\left\langle e_{k},x\right\rangle _{2}\right|^{2},
\]
so $\lim_{k\rightarrow\infty}\left\langle e_{k},x\right\rangle _{2}=0$. 
\end{xca}

\subsection{\label{sub:dmeas}Duality and Measures}
\begin{defn}
Let $E_{i}$, $i=1,2$, be Banach spaces, and let $T:E_{1}\rightarrow E_{2}$
be a linear mapping. We say that $T$ is bounded (continuous) iff
(Def.) $\exists\,C<\infty$, such that
\begin{equation}
\left\Vert Tx\right\Vert _{2}\leq C\left\Vert x\right\Vert _{1},\;\forall x\in E_{1}.\label{eq:ba1}
\end{equation}

\end{defn}

\begin{figure}
\includegraphics[width=0.5\textwidth]{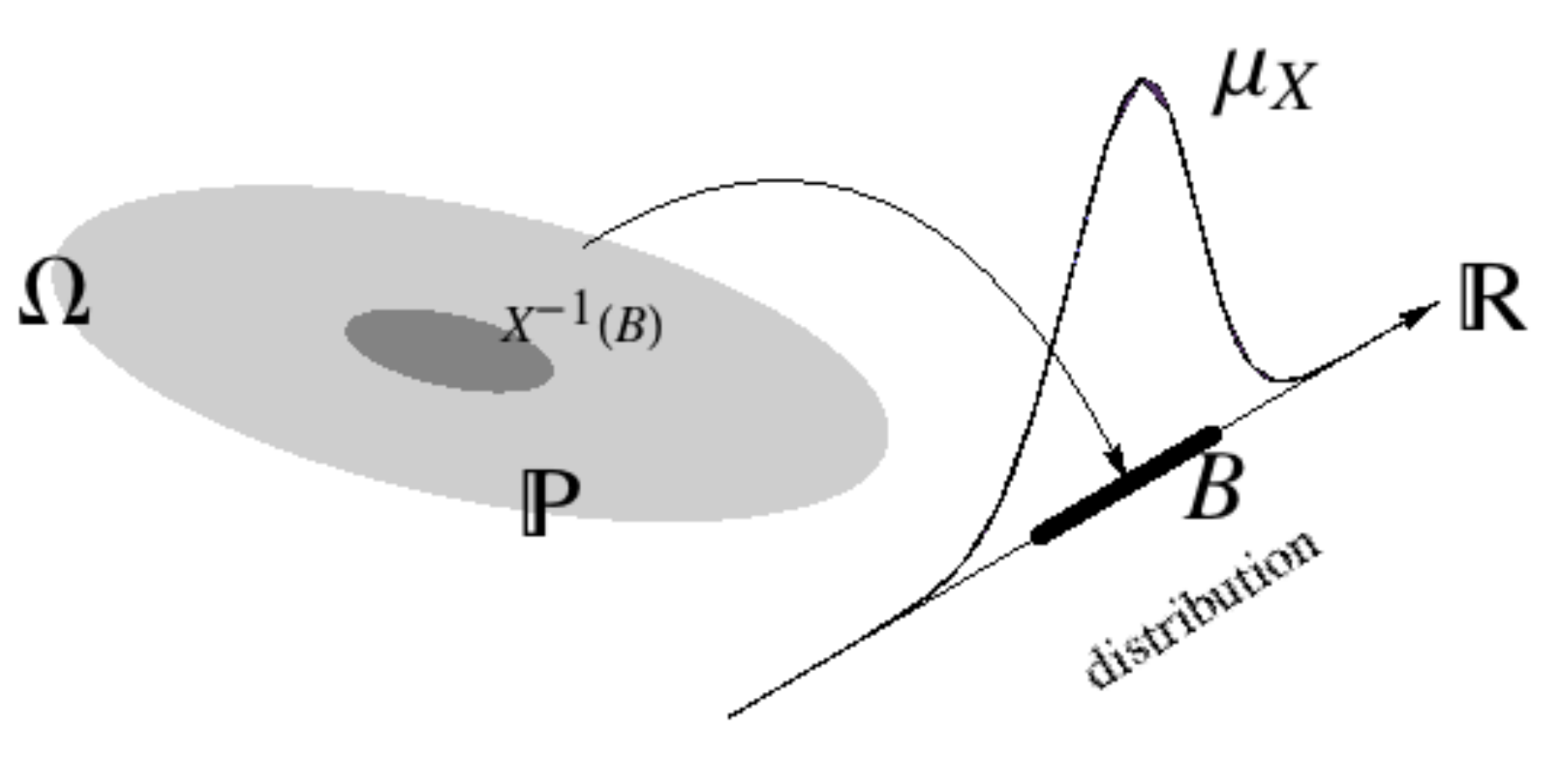}

\protect\caption{\label{fig:rv} A measurement $X$; $X^{-1}\left(B\right)=\left\{ \omega\in\Omega\::\:X\left(\omega\right)\in B\right\} $.
A random variable and its distribution. }
\end{figure}

\begin{defn}
Define $T^{*}:E_{2}^{*}\rightarrow E_{1}^{*}$ by 
\begin{equation}
\left(T^{*}\varphi_{2}\right)\left(x\right)=\varphi_{2}\left(Tx\right),\;\forall x\in E_{1},\;\forall\varphi_{2}\in E_{2}^{*}.\label{eq:ba2}
\end{equation}
We shall adopt the following equivalent notation:
\begin{equation}
\left\langle T^{*}\varphi_{2},x\right\rangle =\left\langle \varphi_{2},Tx\right\rangle ,,\;\forall x\in E_{1},\;\forall\varphi_{2}\in E_{2}^{*}.\label{eq:ba3}
\end{equation}
(Here $E^{*}$ denotes ``dual Banach space.'') It is immediate that
(\ref{eq:ba1}) implies 
\begin{equation}
\left\Vert T^{*}\varphi_{2}\right\Vert _{*}\leq C\left\Vert \varphi_{2}\right\Vert _{*},\;\forall\varphi_{2}\in E_{2}^{*}.\label{eq:ba4}
\end{equation}

\end{defn}
\textbf{Application.} Let $\Omega_{k}$, $k=1,2$, be compact spaces,
and let $\Psi:\Omega_{2}\rightarrow\Omega_{1}$, be a continuous function.
Set
\begin{equation}
Tf=f\circ\Psi,\;\forall f\in C\left(\Omega_{1}\right).\label{eq:ba5}
\end{equation}
Recall the dual Banach spaces:
\begin{eqnarray}
C\left(\Omega_{k}\right)^{*} & = & \mbox{the respective signed measures on \ensuremath{\Omega_{k}}}\label{eq:ba6}\\
 &  & \mbox{of bounded variation, \ensuremath{k=1,2};}\nonumber \\
\left\Vert \mu\right\Vert _{*} & = & \left|\mu\right|\left(\Omega\right)\left(=\mbox{variation of \ensuremath{\mu}}\right)\label{eq:ba7}\\
 & = & \sup\sum_{i}\left|\mu\left(E_{i}\right)\right|,\label{eq:ba8}
\end{eqnarray}
where $\{E_{i}\:|\:E_{i}\in\mathcal{B}\left(\Omega\right)\}$ in (\ref{eq:ba8})
runs over all partitions of $\Omega$.
\begin{xca}[Transformation of measures]
\myexercise{Transformation of measures} Apply (\ref{eq:ba3})-(\ref{eq:ba4})
to show that 
\[
\left(T^{*}\mu_{2}\right)\left(E\right)=\mu_{2}\left(\Psi^{-1}\left(E\right)\right),\;\forall E\in\mathcal{B}\left(\Omega_{1}\right),
\]
or stated equivalently
\[
\int_{\Omega_{1}}f\:d\mu\left(T^{*}\mu_{2}\right)=\int_{\Omega_{2}}\left(f\circ\Psi\right)d\mu_{2},\;\forall f\in C\left(\Omega_{1}\right),\;\forall\mu_{2}\in C\left(\Omega_{2}\right)^{*}.
\]
See Figures \ref{fig:tm} and \ref{fig:tm1} below.\end{xca}
\begin{rem}
We shall make use of the following special case of \emph{pull-back
of measures}. It underlies the notion of ``the distribution of a
\emph{random variable} (math lingo, a measurable function)'' from
statistics. See Figures \ref{fig:rv} and \ref{fig:tm}. We shall
make use of it below, both in the case of a single random variable,
or an indexed family (called a \emph{stochastic process}.)\end{rem}
\begin{defn}
Let $\left(\Omega,\mathcal{F},\mathbb{P}\right)$ be a \emph{probability
space}:
\begin{itemize}
\item[$\Omega$:]  a set, called ``the sample space''.
\item[$\mathcal{F}$:]  a sigma-algebra of subsets of $\Omega$. Elements in $\mathcal{F}$
are called events. 
\item[$\mathbb{P}$:]  a probability measure defined on $\mathcal{F}$, so $\mathbb{P}$
is positive, sigma-additive, and $\mathbb{P}\left(\Omega\right)=1$. 
\end{itemize}

We say a function $X:\Omega\rightarrow\mathbb{R}$ is a \emph{random
variable} iff (Def.) the following implication holds:
\begin{equation}
B\in\mathcal{B}\left(\mathbb{R}\right)\Longrightarrow X^{-1}\left(B\right)\in\mathcal{F};\;\mbox{see Fig}\:\ref{fig:rv}.\label{eq:rv1}
\end{equation}

\end{defn}

So if $X$ is a fixed random variable, there is an induced measure
$\mu_{X}$ on $\mathbb{R}$, a positive Borel measure. It is the pull-back
via $X$, i.e., 
\begin{equation}
\mu_{X}\left(B\right)=\mathbb{P}\left(X^{-1}\left(B\right)\right),\;\forall B\in\mathcal{B}\left(\mathbb{R}\right).\label{eq:rv2}
\end{equation}
If $\mu_{X}$ is Gaussian, see \secref{multi}, we say that $X$ is
a \emph{Gaussian random variable}. If $\mu_{X}$ is uniform, we say
that $X$ is uniformly distributed; and similarly for the other probability
distributions on $\mathbb{R}$; see \tabref{pk} in \secref{dga}
below.

\begin{figure}
\includegraphics[width=0.5\textwidth]{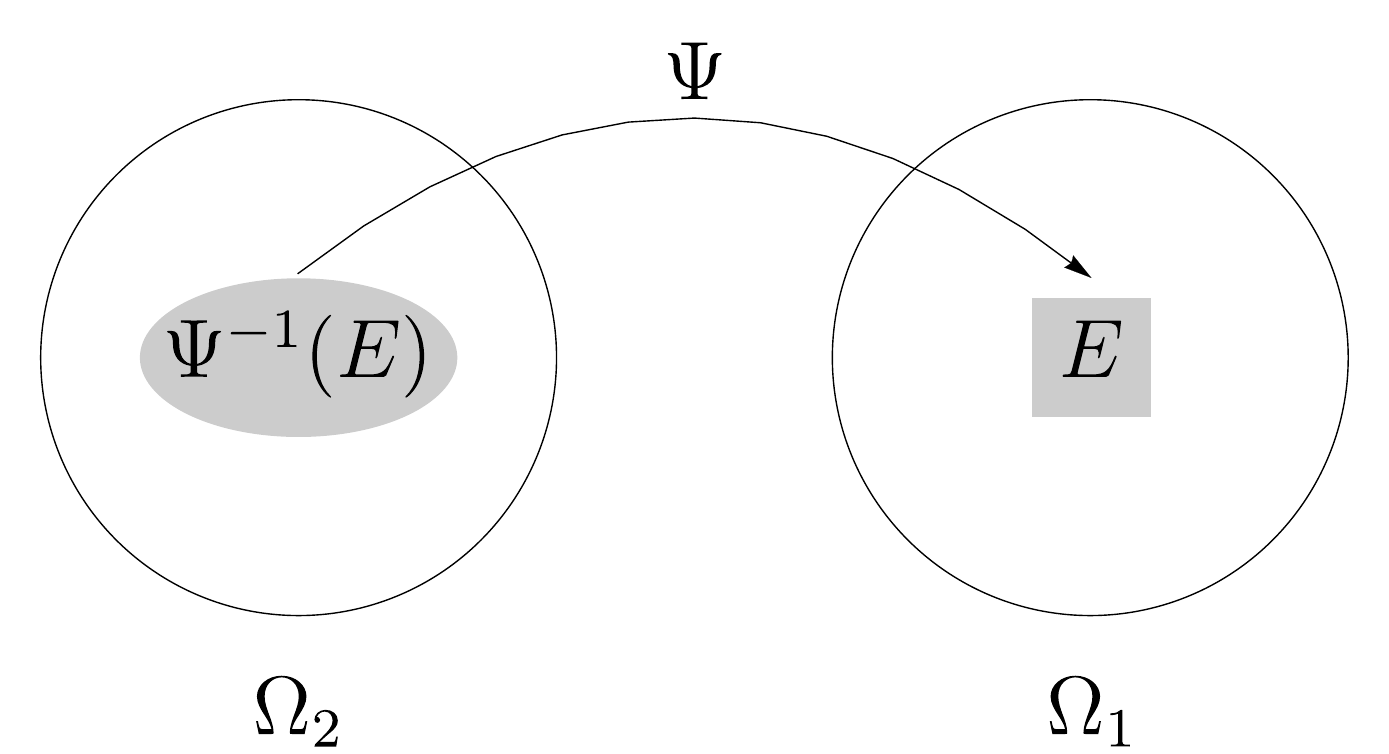}

\protect\caption{\label{fig:tm} $\Psi^{-1}\left(E\right)=\left\{ \omega\in\Omega_{2}\::\:\Psi\left(\omega\right)\in E\right\} $,
pull-back.}
\end{figure}

\renewcommand{\arraystretch}{2}

\begin{figure}
\begin{tabular}{|>{\centering}p{0.1\columnwidth}|>{\centering}p{0.3\columnwidth}|>{\centering}p{0.5\columnwidth}|}
\hline 
$\longrightarrow$ & compact spaces & $\Omega_{2}\xrightarrow{\quad\Psi\quad}\Omega_{1}$\tabularnewline
\hline 
\textbf{$\longleftarrow$} & Banach spaces & $C\left(\Omega_{2}\right)\xleftarrow{\quad T\quad}C\left(\Omega_{1}\right)$\tabularnewline
\hline 
\textbf{$\longrightarrow$} & duals: measures & $\mathcal{M}\left(\Omega_{2}\right)\xrightarrow{\quad T^{*}\quad}\mbox{\ensuremath{\mathcal{M}}}\left(\Omega_{1}\right)$\tabularnewline
\hline 
\end{tabular}

\protect\caption{\label{fig:tm1} Contra-variance (from point transformations, to transformation
of functions, to transformation of measures).}
\end{figure}

\renewcommand{\arraystretch}{1}

\subsection{Other Spaces in Duality}

Below we consider three spaces of functions on $\mathbb{R}$, and
their duals. These are basics of the L. Schwartz' theory of distributions:
\begin{itemize}[itemsep=0.5em,label=]
\item $\mathcal{D}:=C_{c}^{\infty}\left(\mathbb{R}\right)$ = all $C^{\infty}$-functions
on $\mathbb{R}$ having compact support;
\item $\mathcal{S}:=\mathcal{S}\left(\mathbb{R}\right)$ = all $C^{\infty}$-functions
on $\mathbb{R}$ such that $x^{k}f^{\left(n\right)}\in L^{2}\left(\mathbb{R}\right)$
for all $k,n\in\mathbb{N}$;
\item $\mathcal{E}:=C^{\infty}\left(\mathbb{R}\right)$ = all $C^{\infty}$-functions
on $\mathbb{R}$ (without support restriction).
\end{itemize}
Each of the three spaces of test functions $\mathcal{D}$, $\mathcal{S}$,
and $\mathcal{E}$ have countable families of seminorms, turning them
into topological vector spaces (TVS). The two, $\mathcal{S}$ and
$\mathcal{E}$ are Fréchet spaces, while $\mathcal{D}$ is an inductive
limit of Fréchet spaces (abbreviated LF.) 

For $\mathcal{S}$, the seminorms are the max of the absolute value
of the above listed functions, so indexed by $k$ and $n$. For the
other two, $\mathcal{E}$, and $\mathcal{D}$, the seminorms are indexed
by a number $n$ of derivatives, and by compact intervals, say $[-k,k]$.
For each $n$ and $k$, we $\max|f^{\left(n\right)}\left(x\right)|$
over $[-k,k]$. As TVSs, these three spaces in turn are the building
blocks of Schwartz' theory of distributions, see \cite{MR0107812}
and \cite{Tre06}. In each case, the dual space will be defined with
reference to the respective topologies. See details below. 

Clearly, 
\begin{equation}
\mathcal{D}\hookrightarrow\mathcal{S}\hookrightarrow\mathcal{E};\label{eq:sd1}
\end{equation}
but all of the three spaces come with  a natural system of seminorms
turning them into topological vector spaces, and we have the continuous
inclusions $\mathcal{D}\hookrightarrow\mathcal{S}$, and $\mathcal{S}\hookrightarrow\mathcal{E}$. 

Hence for the duals, we have 
\begin{equation}
\mathcal{E}'\hookrightarrow\mathcal{S}'\hookrightarrow\mathcal{D}';\;\mbox{where}\label{eq:sd2}
\end{equation}

\begin{itemize}[itemsep=0.5em,label=]
\item $\mathcal{E}'$ = the space of all compactly supported distributions
on $\mathbb{R}$;
\item $\mathcal{S}'$ = the space of all tempered distributions on $\mathbb{R}$;
and
\item $\mathcal{D}'$ = all distributions on $\mathbb{R}$.
\end{itemize}
\Needspace*{3\baselineskip}
\begin{xca}[Gelfand triple]
\myexercise{Gelfand triple}~
\begin{enumerate}
\item Using \tabref{dual}, show that $L^{2}\left(\mathbb{R}\right)$ is
contained in $\mathcal{S}'$ (= tempered distributions.)
\item Using self-duality of $L^{2}$, i.e., $(L^{2})^{*}\simeq L^{2}$ (by
Riesz), make precise the following double inclusions:
\begin{equation}
\mathcal{S}\hookrightarrow L^{2}\hookrightarrow\mathcal{S}'\label{eq:dl1}
\end{equation}
where each inclusion mapping in (\ref{eq:dl1}) is continuous with
respect to the respective topologies; the Fréchet topology on $\mathcal{S}$,
the norm-topology on $L^{2}$, and the weak-$*$ (dual) topology on
$\mathcal{S}'$. (The system (\ref{eq:dl1}) is an example of a \emph{Gelfand
triple}, see \secref{rhs}.)
\end{enumerate}
\end{xca}

\section{Transfinite Induction (Zorn and All That)}

Let $(X,\leq)$ be a \emph{partially ordered} set. By partial ordering,
we mean a binary relation ``$\leq$'' on the set $X$, such that
(i) $x\leq x$; (ii) $x\leq y$ and $y\leq x$ implies $x=y$; and
(iii) $x\leq y$ and $y\leq z$ implies $x\leq z$. 

A subset $C\subset X$ is said to be a \emph{chain}, or \emph{totally
ordered}, if $x,y\in C$ implies that either $x\leq y$ or $y\leq x$.
Zorn's lemma says that if every chain has a majorant then there exists
a maximal element in $X$. 

\begin{figure}
\includegraphics[width=0.5\textwidth]{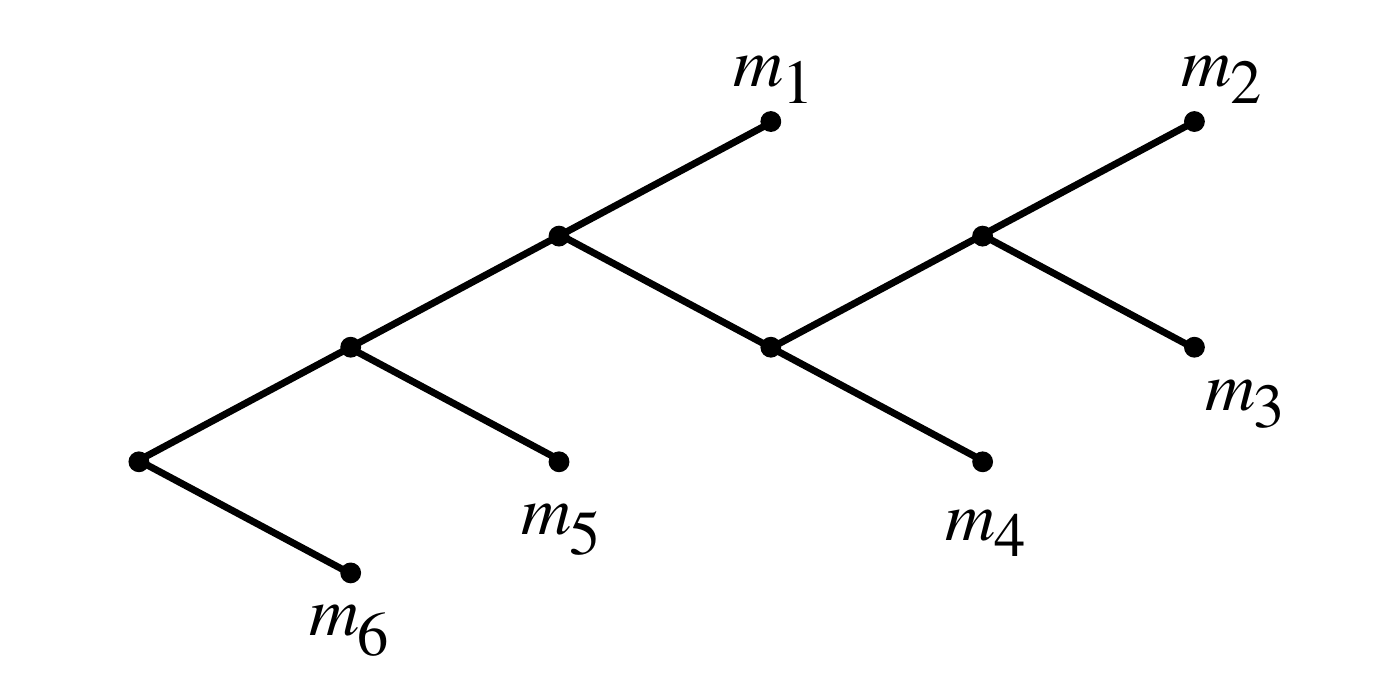}

\protect\caption{\label{fig:zorn1}Finite Tree (natural order on the set of vertices).
Examples of maximal elements: $m_{1},m_{2},\ldots$}
\end{figure}

\begin{figure}
\includegraphics[width=0.5\textwidth]{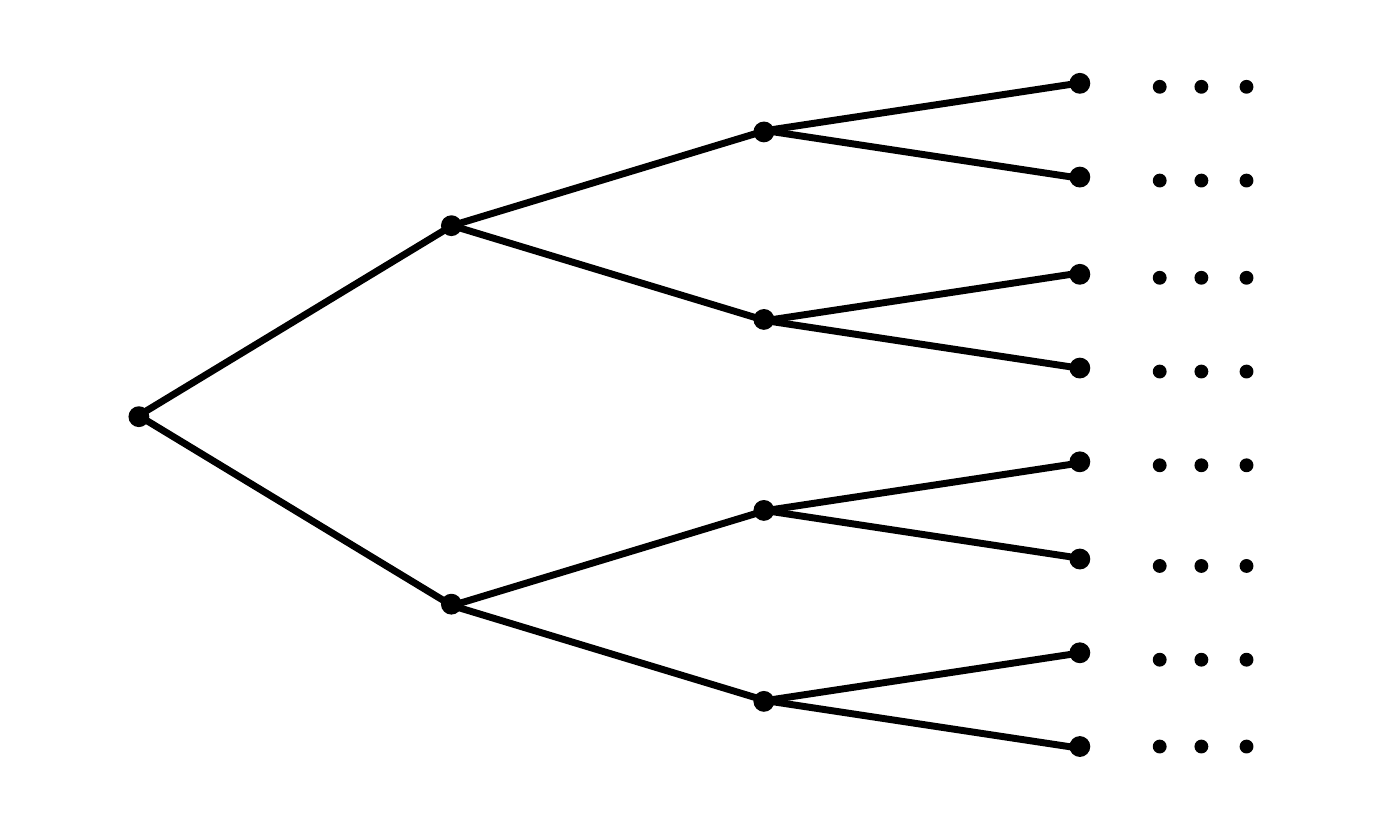}

All finite words in the alphabet $\left\{ 0,1\right\} $, continued
indefinitely.

\protect\caption{\label{fig:zorn2}Infinite Tree (no maximal element!)}
\end{figure}

\index{partially ordered set}\index{Zorn's lemma}
\begin{thm}[Zorn]
 Let $(X,\leq)$ be a partially ordered set. If every chain $C$
in $X$ has a majorant (upper bound), then there exists an element
$m$ in $X$ so that $x\ge m$ implies $x=m$. 
\end{thm}
An illuminating example of a partially ordered set is the binary tree
model (Figs \ref{fig:zorn1}-\ref{fig:zorn2}). Another example is
when $X$ is a family of subsets of a given set, partially ordered
by inclusion. 

Zorn's lemma lies at the foundation of set theory. It is in fact an
axiom and is equivalent to the axiom of choice, and to Hausdorff's
maximality principle. \index{axioms}
\begin{thm}[Hausdorff Maximality Principle]
 Let $(X,\leq)$ be a partially ordered set, then there exists a
maximal totally ordered subset $L$ in $X$. 
\end{thm}

The axiom of choice is equivalent to the following statement on infinite
products, which itself is extensively used in functional analysis. 
\begin{thm}[axiom of choice]
\label{thm:aoc} Let $A_{\alpha}$ be a family of nonempty sets indexed
by $\alpha\in I$. Then the infinite Cartesian product 
\[
\Omega=\prod_{\alpha\in I}A_{\alpha}=\left\{ \omega:I\rightarrow\cup_{\alpha\in I}A_{\alpha}\:\big|\:\omega\left(\alpha\right)\in A_{\alpha}\right\} 
\]
is nonempty.
\end{thm}
The point of using the axiom of choice is that, if the index set is
uncountable, there is no way to verify whether $(x_{\alpha})$ is
in $\Omega$, or not. It is just impossible to check for each $\alpha$
that $x_{\alpha}$ is contained in $A_{\alpha}$. 
\begin{rem}
A more dramatic consequence of the axiom of choice is the mind-boggling
\emph{Banach-Tarski paradox}; see e.g., \cite{MR3104075}. It states:
For the solid ball $B$ in 3-dimensional space, there exists a decomposition
of $B$ into a finite number of disjoint subsets, which can then in
turn be put back together again, but in a different way which will
yield two identical copies of the original ball $B$; -- stated informally
as: \textquotedbl{}A pea can be chopped up and reassembled into the
Sun.\textquotedbl{} The axiom of choice allows for the construction
of nonmeasurable sets, i.e., sets that do not have a volume, and that
for their construction would require performing an uncountably infinite
number of choices. \index{Banach-Tarski paradox}
\end{rem}
In case the set is countable, we simply apply the down to earth standard
induction. Note that the standard mathematical induction is equivalent
to the Peano's axiom: \emph{Every nonempty subset of the set of natural
number has a unique smallest element.} The power of transfinite induction
is that it applies to uncountable sets as well. \index{Peano's axiom}

In applications, the key of using the transfinite induction is to
cook up, in a clear way, a partially ordered set, so that the maximal
element turns out to be the object to be constructed. \index{partially ordered set}

Examples include \emph{Hahn-Banach extension theorem}, \emph{Krein-Milman}'s
theorem \index{Theorem!Krein-Milman's-}on compact convex set, existence
of orthonormal bases in Hilbert space, \emph{Tychnoff's theorem} on
infinite Cartesian product of compact spaces (follows immediately
from the axiom of choice.)

\index{extension!-of functional}

\index{compact}

\index{convex}

\index{product!infinite}

\index{Krein-Milman}
\begin{thm}[Tychonoff]
 Let $A_{\alpha}$ be a family of compact sets indexed by $\alpha\in I$.
Then the infinite Cartesian product $\prod_{\alpha}A_{\alpha}$ is
compact with respect to the product topology. \index{product topology}\index{Tychonoff (compactness)}
\end{thm}
We will apply transfinite induction (Zorn's lemma\index{Zorn's lemma})
to show that every infinite dimensional Hilbert space has an orthonormal
basis (ONB). 

\index{orthonormal basis (ONB)} \index{Theorem!Tychonoff-}

\section{Basics of Hilbert Space Theory\label{sec:Hilbert}}

Key to functional analysis is the idea of \emph{normed vector spaces}.
The interesting ones are infinite-dimensional. To use them effectively
in the solution of problems, we must be able to take limits, hence
the assumption of completeness. A complete normed linear space is
called a \emph{Banach space}. But for applications in physics, statistics,
and in engineering it often happens that the norm comes from an \emph{inner
product}; -- this is the case of \emph{Hilbert space}. With an inner
product, one is typically able to get much more precise results, than
in the less structured case of Banach space. (Many Banach spaces are
not Hilbert spaces.)

The more interesting Hilbert spaces typically arise in concrete applications
as infinite-dimensional spaces of function. And as such, they have
proved indispensable tools in the study of partial differential equations
(PDE), in quantum mechanics, in Fourier analysis, in signal processing,
in representations of groups, and in ergodic theory. The term Hilbert
space was originally coined by John von Neumann, who identified the
axioms that now underlie these diverse applied areas. Examples include
spaces of \emph{square-integrable functions} (e.g., the $L^{2}$ random
variables of a probability space), \emph{Sobolev spaces}, Hilbert
spaces of Schwartz distributions, and \emph{Hardy spaces} of holomorphic
functions; -- to mention just a few.

One reason for their success is that geometric intuition from finite
dimensions carries over: e.g., the Pythagorean Theorem, the parallelogram
law; and, for optimization problems, the important notion of \textquotedblleft \emph{orthogonal
projection}.\textquotedblright{} And the idea (from linear algebra)
of diagonalizing a normal matrix; -- the \emph{spectral theorem}.\index{signal processing}

Linear mappings (transformations) between Hilbert spaces are called
linear operators, or simply \textquotedblleft operators.\textquotedblright{}
They include partial differential operators (PDOs), and many others. 
\begin{defn}
Let $X$ be a vector space over $\mathbb{C}$. 

A norm on $X$ is a mapping $\left\Vert \cdot\right\Vert :X\rightarrow\mathbb{C}$
such that 
\begin{itemize}
\item $\left\Vert cx\right\Vert =\left|c\right|\left\Vert x\right\Vert $,
$c\in\mathbb{C}$, $x\in X$;
\item $\left\Vert x\right\Vert \geq0$; $\left\Vert x\right\Vert =0$ implies
$x=0$, for all $x\in X$; 
\item $\left\Vert x+y\right\Vert \leq\left\Vert x\right\Vert +\left\Vert y\right\Vert $,
for all $x,y\in X$.
\end{itemize}
\end{defn}

\begin{defn}
Let $\left(X,\left\Vert \cdot\right\Vert \right)$ be a normed space.
$X$ is called a \emph{Banach space} if it is complete with respect
to the induced metric 
\[
d\left(x,y\right):=\left\Vert x-y\right\Vert ,\;x,y\in X.
\]

\end{defn}

\begin{defn}
\label{def:ip}Let $X$ be vector space over $\mathbb{C}$. An inner
product on $X$ is a function $\left\langle \cdot,\cdot\right\rangle :X\times X\rightarrow\mathbb{C}$
so that for all $x,y\in\mathscr{H}$, and $c\in\mathbb{C}$, we have
\begin{itemize}
\item $\left\langle x,\cdot\right\rangle :X\rightarrow\mathbb{C}$ is linear
(linearity)
\item $\left\langle x,y\right\rangle =\overline{\left\langle y,x\right\rangle }$
(conjugation) 
\item $\left\langle x,x\right\rangle \geq0$; and $\left\langle x,x\right\rangle =0$
implies $x=0$ (positivity)
\end{itemize}
\end{defn}
\begin{rem}
The abstract formulation of Hilbert space was invented by von Neumann
in 1925. It fits precisely with the axioms of quantum mechanics (spectral
lines, etc.) A few years before von Neumann\textquoteright s formulation,
Max Born had translated Heisenberg\textquoteright s quantum mechanics
into modern mathematics. In 1924, in a break-through paper, Heisenberg
had invented quantum mechanics, but he had not been precise about
the mathematics. His use of \textquotedblleft matrices\textquotedblright{}
was highly intuitive. It was only in the subsequent years, with the
axiomatic language of Hilbert space, that the group of physicists
and mathematicians around Hilbert in Göttingen were able to give the
theory the form it now has in modern textbooks.
\end{rem}
\index{von Neumann, J.}

\index{Cauchy-Schwarz' inequality}

\index{inequality!Cauchy-Schwarz}\index{axioms}

\begin{lem}[Cauchy-Schwarz]
\label{lem:CS}\footnote{Hermann Amandus Schwarz (1843 - 1921), German mathematician, contemporary
of Weierstrass, and known for his work in complex analysis. He is
the one in many theorems in books on analytic functions. We will often
refer to (\ref{eq:h1}) as simply ``Schwarz''. The abbreviation
is useful because we use it a lot. 

There are other two ``Schwartz'' (with a ``t''): 

Laurent Schwartz (1915 - 2002), French mathematician, Fields Medal
in 1950 for his work of distribution theory. 

Jack Schwartz (1930 - 2009), American mathematician, author of the
famous book ``Linear Operators''. }Let $\left(X,\left\langle \cdot,\cdot\right\rangle \right)$ be an
inner product space, then
\begin{equation}
\left|\left\langle x,y\right\rangle \right|^{2}\leq\left\langle x,x\right\rangle \left\langle y,y\right\rangle ,\;\forall x,y\in X.\label{eq:h1}
\end{equation}
\end{lem}
\begin{proof}
By the positivity axiom in the definition of an inner product, we
see that 
\[
\sum_{i,j=1}^{2}\overline{c_{i}}c_{j}\left\langle x_{i},x_{j}\right\rangle =\left\langle \sum_{i=1}^{2}c_{i}x_{i},\sum_{j=1}^{2}c_{j}x_{j}\right\rangle \geq0,\;\forall c_{1},c_{2}\in\mathbb{C};
\]
i.e., the matrix 
\[
\left[\begin{array}{cc}
\left\langle x_{1},x_{1}\right\rangle  & \left\langle x_{1},x_{2}\right\rangle \\
\left\langle x_{2},x_{1}\right\rangle  & \left\langle x_{2},x_{2}\right\rangle 
\end{array}\right]
\]
is positive definite. Hence the above matrix has nonnegative determinant,
and (\ref{eq:h1}) follows. \index{positive definite!-matrix}\end{proof}
\begin{cor}
Let $\left(X,\left\langle \cdot,\cdot\right\rangle \right)$ be an
inner product space, then 
\begin{equation}
\left\Vert x\right\Vert :=\sqrt{\left\langle x,x\right\rangle },\;x\in X\label{eq:h2}
\end{equation}
defines a norm. \end{cor}
\begin{proof}
It suffices to check the triangle inequality (\defref{ip}). For all
$x,y\in X$, we have (with the use of \lemref{CS}): 
\begin{eqnarray*}
\left\Vert x+y\right\Vert ^{2} & = & \left\langle x+y,x+y\right\rangle \\
 & = & \left\Vert x\right\Vert ^{2}+\left\Vert y\right\Vert ^{2}+2\Re\left\{ \left\langle x,y\right\rangle \right\} \\
 & \leq & \left\Vert x\right\Vert ^{2}+\left\Vert y\right\Vert ^{2}+2\left\Vert x\right\Vert \left\Vert y\right\Vert \quad\left(\text{by}\:\left(\ref{eq:h1}\right)\right)\\
 & = & \left(\left\Vert x\right\Vert +\left\Vert y\right\Vert \right)^{2}
\end{eqnarray*}
and the corollary follows.\end{proof}
\begin{defn}
An inner product space $\left(X,\left\langle \cdot,\cdot\right\rangle \right)$
is called a \emph{Hilbert space} if $X$ is complete with respect
to the metric 
\[
d\left(x,y\right)=\left\Vert x-y\right\Vert ,\;x,y\in X;
\]
where the RHS is given by (\ref{eq:h2}). \end{defn}
\begin{xca}[Hilbert completion]
\myexercise{Hilbert completion}Let $\left(X,\left\langle \cdot,\cdot\right\rangle \right)$
be an inner-product space (\defref{ip}), and let $\mathscr{H}$ be
its \emph{metric completion} with respect to the norm in (\ref{eq:h2}).
Show that $\left\langle \cdot,\cdot\right\rangle $ on $X\times X$
extends by limit to a sesquilinear form $\left\langle \cdot,\cdot\right\rangle ^{\sim}$
on $\mathscr{H}\times\mathscr{H}$ ; and that $\mathscr{H}$ with
$\left\langle \cdot,\cdot\right\rangle ^{\sim}$ \uline{is} a Hilbert
space. 
\end{xca}
\index{completion!Hilbert-}\index{sesquilinear form}

\begin{figure}
\begin{minipage}[t]{1\columnwidth}%
\begin{center}
\subfloat[$\varphi\left(s,t\right)=s\wedge t$]{\protect\begin{centering}
\protect\includegraphics[width=0.3\textwidth]{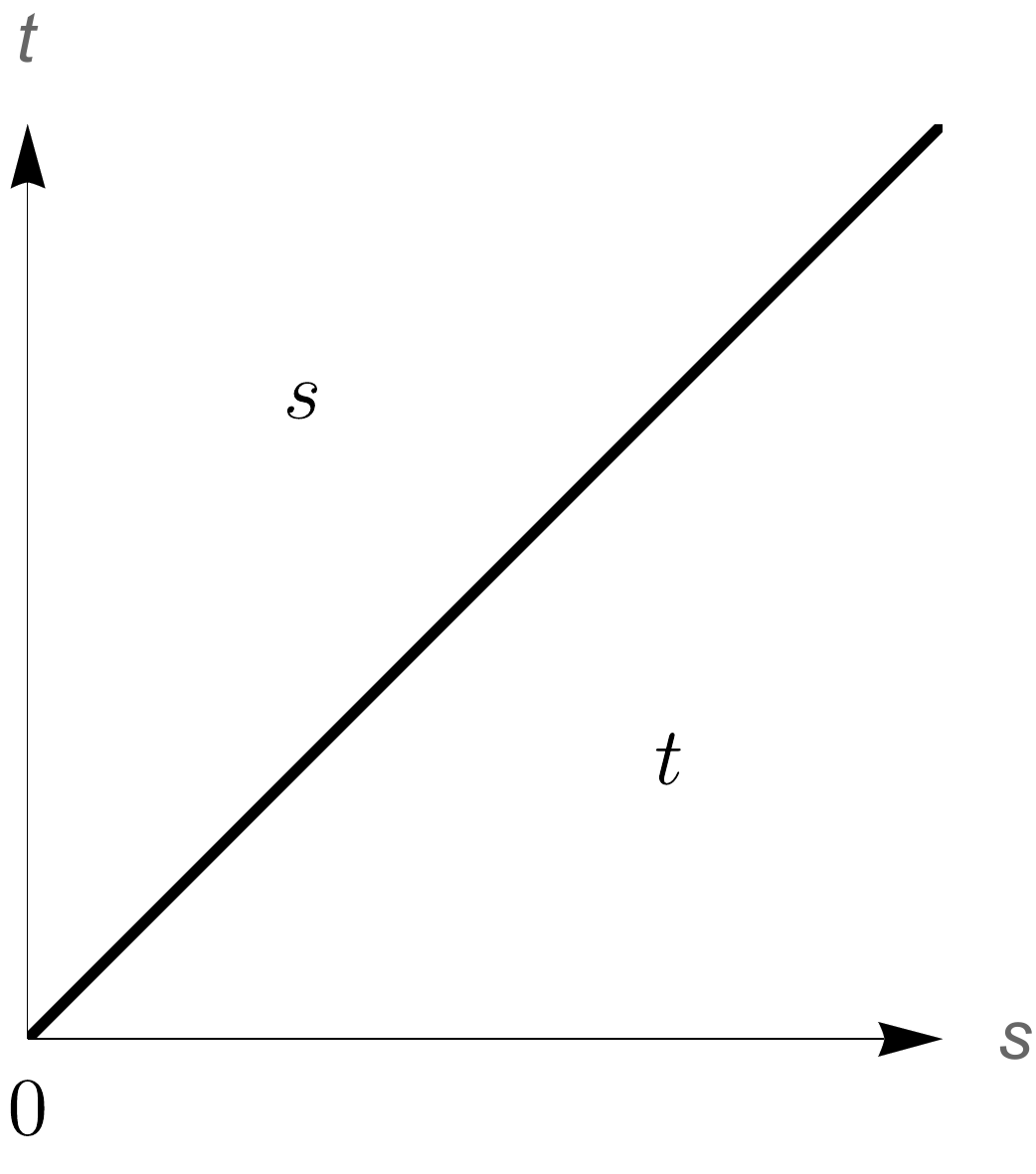}\protect
\par\end{centering}

}\hfill{}\subfloat[with $t$ fixed]{\protect\begin{centering}
\protect\includegraphics[width=0.3\textwidth]{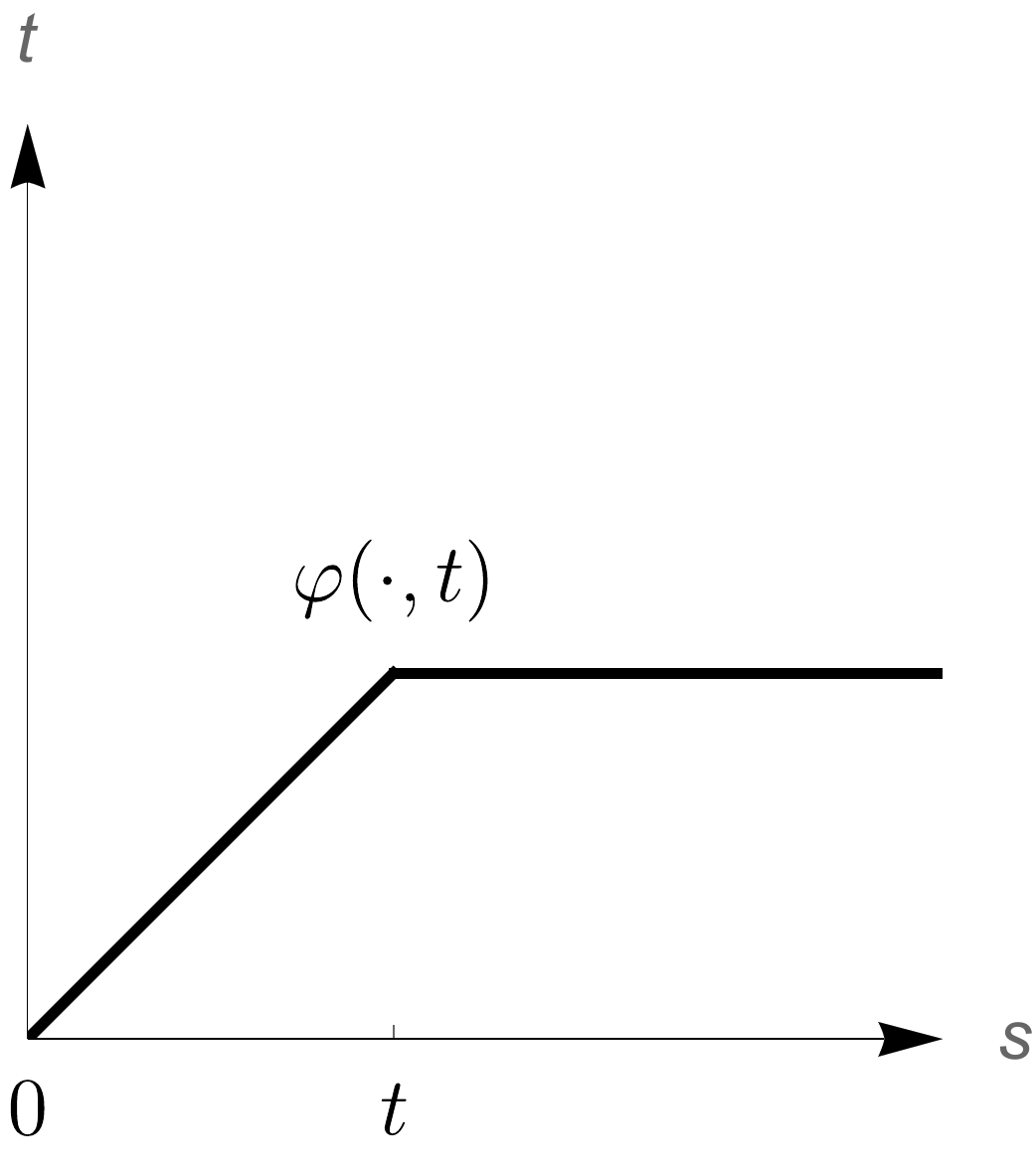}\protect
\par\end{centering}

}
\par\end{center}%
\end{minipage}

\centering{}\protect\caption{\label{fig:gns}Covariance function of Brownian motion.}
\end{figure}

\begin{figure}
\includegraphics[width=0.5\textwidth]{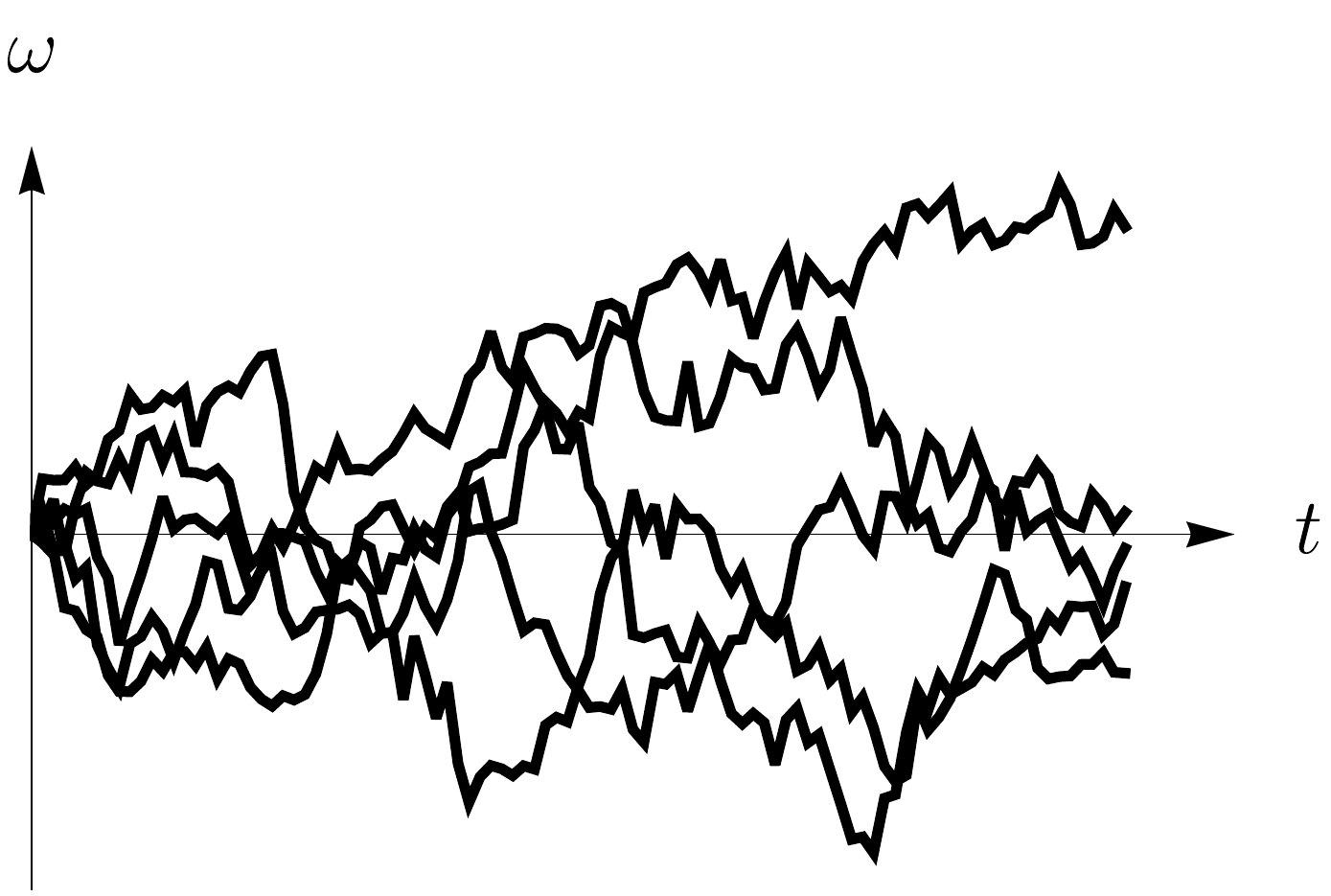}

\protect\caption{\label{fig:bmp}A set of Brownian sample-paths generated by a Monte-Carlo
computer simulation.}
\end{figure}

\Needspace*{3\baselineskip}
\begin{xca}[$L^{2}$ of a measure-space]
\label{exer:L2m}\myexercise{$L^{2}$ of a measure-space}Let $\left(M,\mathcal{B},\mu\right)$
be as follows:
\begin{enumerate}[label=]
\item $M:$ locally compact Hausdorff space;
\item $\mathcal{B}:$ the Borel sigma-algebra, i.e., generated by the open
subsets of $M$;
\item $\mu:$ a fixed positive measure defined on $\mathcal{B}$.
\end{enumerate}

Let $\mathcal{F}:=span\left\{ \chi_{E}\:|\:E\in\mathcal{B}\right\} $,
and on linear combinations, set
\begin{equation}
\left\Vert \sum_{i:\text{ finite}}c_{i}\chi_{E_{i}}\right\Vert _{\mathscr{H}}^{2}=\sum_{i}\left|c_{i}\right|^{2}\mu\left(E_{i}\right)\label{eq:hc}
\end{equation}
where $E_{i}\in\mathcal{B}$, and $E_{i}\cap E_{j}=\emptyset$ ($i\neq j$),
are assumed.

Show that the Hilbert-completion of $\mathcal{F}$ with respect to
to (\ref{eq:hc}) agrees with the standard definitions \cite{Rud87,parthasarathy1982probability}
of the $L^{2}\left(\mu\right)$-space.\end{xca}
\begin{rem}
\label{rem:gns}An extremely useful method to build Hilbert spaces
is the GNS construction. For details, see  \chapref{GNS}. 

The idea is to start with a positive definite function $\varphi:X\times X\rightarrow\mathbb{C}$,
defined on an arbitrary set $X$. We say $\varphi$ is \emph{positive
definite}, if for all $n\in\mathbb{N}$, \index{positive definite!-function}
\begin{equation}
\sum_{i,j=1}^{n}\overline{c_{i}}c_{j}\varphi\left(x_{i},x_{j}\right)\geq0\label{eq:h3}
\end{equation}
for all system of coefficients $c_{1},\ldots,c_{n}\in\mathbb{C}$,
and all $x_{1},\ldots,x_{n}\in X$. 

Given $\varphi$, set 
\[
H_{0}:=\left\{ \sum_{\mbox{finite}}c_{x}\delta_{x}:x\in X,c_{x}\in\mathbb{C}\right\} =span_{\mathbb{C}}\left\{ \delta_{x}:x\in X\right\} ,
\]
and define a sesquilinear form\index{sesquilinear form} on $H_{0}$
by 

\[
\left\langle \sum c_{x}\delta_{x},\sum d_{y}\delta_{y}\right\rangle _{\varphi}:=\sum\overline{c_{x}}d_{y}\varphi\left(x,y\right).
\]
Note that 
\[
\left\Vert \sum c_{x}\delta_{x}\right\Vert _{\varphi}^{2}:=\left\langle \sum c_{x}\delta_{x},\sum c_{x}\delta_{x}\right\rangle _{\varphi}=\sum_{x,y}\overline{c_{x}}c_{y}\varphi\left(x,y\right)\geq0
\]
by assumption. (All summations are finite.)

However, $\left\langle \cdot,\cdot\right\rangle _{\varphi}$ is in
general not an inner product since the strict positivity axiom may
not be satisfied. Hence one has to pass to a quotient space by letting
\[
N=\left\{ f\in H_{0}\:\big|\:\left\langle f,f\right\rangle _{\varphi}=0\right\} ,
\]
and set $\mathscr{H}:=$ completion of the quotient space $H_{0}/N$
with respect to $\left\Vert \cdot\right\Vert _{\varphi}$. (The fact
that $N$ is really a subspace follows from (\ref{eq:h1}).) $\mathscr{H}$
is a Hilbert space.\index{space!Hilbert-}\index{completion!Hilbert-}\end{rem}
\begin{cor}
\label{cor:gns}Let $X$ be a set, and let $\varphi:X\times X\rightarrow\mathbb{C}$
be a function. Then $\varphi$ is positive definite if and only if
there is a Hilbert space $\mathscr{H}=\mathscr{H}_{\varphi}$, and
a function $\Phi:X\rightarrow\mathscr{H}$ such that 
\begin{equation}
\varphi\left(x,y\right)=\left\langle \Phi\left(x\right),\Phi\left(y\right)\right\rangle _{\mathscr{H}}\label{eq:pd}
\end{equation}
for all $\left(x,y\right)\in X\times X$, where $\left\langle \cdot,\cdot\right\rangle _{\mathscr{H}}$
denotes the inner product in $\mathscr{H}$. \index{positive definite!-function}

Given a solution $\Phi$ satisfying (\ref{eq:pd}), then we say that
$\mathscr{H}$ is minimal if 
\begin{equation}
\mathscr{H}=\overline{span}\left\{ \Phi\left(x\right)\::\:x\in X\right\} .\label{eq:pd2}
\end{equation}
Given two minimal solutions, $\Phi_{i}:X\rightarrow\mathscr{H}_{i}$,
$i=1,2$ (both satisfying (\ref{eq:pd})); then there is a unitary
isomorphism $\mathcal{U}:\mathscr{H}_{1}\rightarrow\mathscr{H}_{2}$
such that
\begin{equation}
\mathcal{U}\Phi_{1}\left(x\right)=\Phi_{2}\left(x\right)\mathcal{U},\;\forall x\in X.\label{eq:pd3}
\end{equation}
\end{cor}
\begin{proof}
These conclusions follow from  \remref{gns}, and the definitions.
(The missing details are left as an exercise to the student.)\end{proof}
\begin{rem}
\label{rem:gns-1}It is possible to be more explicit about choice
of the pair $\left(\Phi,\mathscr{H}\right)$ in \corref{gns}, where
$\varphi:X\times X\rightarrow\mathbb{C}$ is a given positive definite
function. We may in fact choose $\mathscr{H}$ to be $L^{2}\left(\Omega,\mathcal{F},\mathbb{P}\right)$
where $\mathbb{P}=\mathbb{P}_{\varphi}$ depends on $\varphi$, and
$\left(\Omega,\mathcal{F},\mathbb{P}\right)$ is a probability space. \end{rem}
\begin{example}[Wiener-measure]
\label{exa:bm}In \remref{gns-1}, take $X=[0,\infty)=\mathbb{R}_{+}\cup\left\{ 0\right\} $,
and set 
\[
\varphi\left(s,t\right)=s\wedge t=\min\left(s,t\right),
\]
see \figref{gns}. In this case, we may then take $\Omega=C\left(\mathbb{R}\right)=$
all continuous functions on $\mathbb{R}$, and $\Phi_{t}\left(\omega\right):=\omega\left(t\right)$,
$t\in[0,\infty)$, $\omega\in C\left(\mathbb{R}\right)$. 

Further, the sigma-algebra\index{sigma-algebra} in $C\left(\mathbb{R}\right)$,
$\mathcal{F}:=Cyl$, is generated by cylinder-sets, and $\mathbb{P}$
is the \emph{Wiener-measure}; and $\Phi$ on $L^{2}\left(C\left(\mathbb{R}\right),Cyl,\mathbb{P}\right)$
is the standard \emph{Brownian motion}, i.e., $\Phi:[0,\infty)\rightarrow L^{2}\left(C\left(\mathbb{R}\right),\mathbb{P}\right)$
is a Gaussian process with\index{cylinder-set}\index{Gaussian process}
\[
\mathbb{E}_{\mathbb{P}}\left(\Phi\left(s\right)\Phi\left(t\right)\right)=\int_{C\left(\mathbb{R}\right)}\Phi_{s}\left(\omega\right)\Phi_{t}\left(\omega\right)d\mathbb{P}\left(\omega\right)=s\wedge t.
\]

The process $\left\{ \Phi_{t}\right\} $ is called the Brownian motion;
its properties include that each $\Phi_{t}$ is a Gaussian random
variable. We refer to \chapref{bm} for full details. \figref{bmp}
shows a set of sample path of the standard Brownian motion.
\end{example}
\index{Wiener, N.}

\begin{xca}[Product of two positive definite functions]
\myexercise{Product of two positive definite functions}Let $\varphi$
and $\psi$ be positive definite functions $X\times X\longrightarrow\mathbb{C}$
(see (\ref{eq:h3})) and set 
\[
\xi\left(x,y\right)=\varphi\left(x,y\right)\psi\left(x,y\right),\;\forall\left(x,y\right)\in X\times X.
\]
Show that $\xi=\varphi\cdot\psi$ is again positive definite.

\uline{Hint}: Use \remref{gns} and the fact that every positive
$n\times n$ matrix $B$ has the form $B=A^{*}A$. Fix $n$, and $x_{1},\ldots,x_{n}\in X$.
Apply this to $B_{ij}:=\varphi\left(x_{i},x_{j}\right)$.
\end{xca}
We resume the discussion of stochastic processes in \chapref{RKHS}
below.
\begin{rem}
We see from the proof of \lemref{CS} that the Cauchy-Schwarz inequality
holds for all positive definite functions.\end{rem}
\begin{defn}
Let $\mathscr{H}_{i}$, $i=1,2$ be two Hilbert spaces. 

A linear operator $J:\mathscr{H}_{1}\rightarrow\mathscr{H}_{2}$ is
said to be an \emph{isometry} iff (Def.)
\[
\left\Vert Jx\right\Vert _{\mathscr{H}_{2}}=\left\Vert x\right\Vert _{\mathscr{H}_{1}},\;\forall x\in\mathscr{H}_{1}.
\]
Note that $J$ is \uline{not} assumed ``onto.''
\end{defn}
\index{isometry!It-@It\=o-} \index{integral!It-@It\=o-}\index{integral!Stochastic-}
\begin{xca}[An It\=o-isometry]
\myexercise{An It\=o-isometry} Let $T$ be a locally compact Hausdorff
space; and let $\mu$ be a positive Borel measure, i.e., consider
the measure space $(T,\mathcal{B}\left(T\right),\mu)$, where $\mathcal{B}\left(T\right)$
is the Borel sigma-algebra.
\begin{enumerate}
\item \label{enu:isog1}Show that there is a measure space $(\Omega,\mathcal{F},\mathbb{P}^{\left(\mu\right)})$
depending on $\mu$; and a function (Gaussian process)
\begin{equation}
\Phi:\mathcal{B}\left(T\right)\longrightarrow L^{2}\left(\Omega,\mathbb{P}\right)\label{eq:isog1}
\end{equation}
such that every $\Phi_{A}$, for $A\in\mathcal{B}\left(T\right)$,
is a Gaussian random variable, such that $\mathbb{E}(\Phi_{A})=0$,
and 
\begin{equation}
\mathbb{E}\left(\Phi_{A}\Phi_{B}\right)=\mu\left(A\cap B\right)\label{eq:isog2}
\end{equation}
holds for all $A,B\in\mathcal{B}\left(T\right)$. The expectation
$\mathbb{E}$ in (\ref{eq:isog2}) is with respect to $\mathbb{P}=\mathbb{P}^{\left(\mu\right)}$. 
\item Show that there is an isometry $J=J_{\left(Ito\right)}$ from $L^{2}\left(T,\mu\right)$
into $L^{2}\left(\Omega,\mathbb{P}\right)$ such that 
\begin{equation}
\mathbb{E}\left(\left|\int_{T}f\left(t\right)d\Phi_{t}\right|^{2}\right)=\int_{T}\left|f\right|^{2}d\mu\label{eq:isog3}
\end{equation}
where the expression $\int_{T}f\left(t\right)d\Phi_{t}$ on the LHS
in (\ref{eq:isog3}) is the $L^{2}$-limit of finite sums (simple
functions):
\begin{equation}
\sum_{i}c_{i}\Phi_{A_{i}},\label{eq:isog4}
\end{equation}
$c_{i}\in\mathbb{R}$, finite indexing; and $A_{i}\in\mathcal{B}\left(T\right)$,
$A_{i}\cap A_{j}=\emptyset$, $i\neq j$ (i.e., disjointness.) \\
Because of (\ref{eq:isog3}), we can set 
\begin{equation}
Jf=\int_{T}f\left(t\right)d\Phi_{t}\in L^{2}(\Omega,\mathbb{P}^{\left(\mu\right)}).
\end{equation}

\end{enumerate}

\uline{Hint}: 
\begin{enumerate}
\item This is an application of \corref{gns} (\secref{Hilbert}), applied
to $X=\mathcal{B}\left(T\right)$. Note that 
\begin{equation}
\mathcal{B}\left(T\right)\times\mathcal{B}\left(T\right):\left(A,B\right)\longmapsto\mu\left(A\cap B\right)\label{eq:isog5}
\end{equation}
is positive definite. Hence, the existence of the Gaussian process
\[
\left(\Omega,\mathcal{F},\mathbb{P},\Phi\right)
\]
subject to the conditions in part (\ref{enu:isog1}) of the exercise,
follows from \corref{gns}.
\item Consider the simple functions in (\ref{eq:isog4}), and use (\ref{eq:isog2}).
Then derive the following:
\begin{equation}
\mathbb{E}\left(\left|\sum_{i}c_{i}\Phi_{A_{i}}\right|^{2}\right)=\sum_{i}\left|c_{i}\right|^{2}\mu\left(A_{i}\right)\label{eq:isog6}
\end{equation}

\end{enumerate}

Hence, the \emph{It\={o}-isometry}, defined initially only on simple
functions, is isometric. Justify the extension by limits to all of
$L^{2}\left(T,\mu\right)$. In your last step, taking the limit over
partitions, make use of the conclusion from \exerref{L2m}.

\end{xca}
\begin{rem}
The construction in the exercise is an example of a \emph{stochastic
process}, and a \emph{stochastic integral}. Both subjects are resumed
in Chapters \ref{chap:bm} and \ref{chap:RKHS} below.
\end{rem}
\index{stochastic process}

\subsection{\label{sub:ONB}Orthonormal Bases}
\begin{defn}
Let $\mathscr{H}$ be a Hilbert space. A family of vectors $\{u_{\alpha}\}$
in $\mathscr{H}$ is said to be an orthonormal basis if 
\begin{enumerate}
\item $\left\langle u_{\alpha},u_{\beta}\right\rangle _{\mathscr{H}}=\delta_{\alpha\beta}$
and 
\item $\overline{span}\left\{ u_{\alpha}\right\} =\mathscr{H}$. (Here ``$\overline{span}$''
means ``closure of the linear span.'')
\end{enumerate}
\end{defn}

We are now ready to prove the existence of \emph{orthonormal bases}
for any Hilbert space. The key idea is to cook up a partially ordered
set\index{partially ordered set} satisfying all the requirements
for transfinite induction, so that each maximal element turns out
to be an orthonormal basis (ONB). Notice that all we have at hands
are the abstract axioms of a Hilbert space, and nothing else. Everything
will be developed out of these axioms. A separate issue is \emph{constructive}
ONBs, for example, wavelets or orthogonal polynomials.\index{orthogonal!-vectors}\index{axioms}

There is a new theory which generalizes the notion of ONB; called
``\uline{frame}'', and it is discussed in \chapref{KS} below,
along with some more applications.
\begin{thm}
\label{thm:ONB}Every Hilbert space $\mathscr{H}$ has an orthonormal
basis.
\end{thm}
To start out, we need the following lemmas. 
\begin{lem}
\label{lem:dense}Let $\mathscr{H}$ be a Hilbert space and $S\subset\mathscr{H}$.
Then the following are equivalent:
\begin{enumerate}
\item $x\perp S$ implies $x=0$
\item $\overline{span}\{S\}=\mathscr{H}$
\end{enumerate}
\end{lem}

\begin{proof}
Now we prove \thmref{ONB}. If $\mathscr{H}=\left\{ 0\right\} $ then
the proof is done. Otherwise, let $u_{1}\in\mathscr{H}$. If $\left\Vert u_{1}\right\Vert \neq1$,
it can be normalized by $u_{1}/\left\Vert u_{1}\right\Vert $. Hence
we may assume $\left\Vert u_{1}\right\Vert =1$. If $span\{u_{1}\}=\mathscr{H}$
the proof finishes again, otherwise there exists $u_{2}\notin span\{u_{1}\}$.
By \lemref{dense}, we may assume $\left\Vert u_{2}\right\Vert =1$
and $u_{1}\perp u_{2}$. It follows that there exists a collection
$S$ of orthonormal vectors in $\mathscr{H}$. 

Let $\mathbb{P}\left(S\right)$ be the set of all orthonormal sets
partially ordered by set inclusion. Let $C\subset\mathbb{P}(S)$ be
any chain and let $M:=\bigcup_{E\in C}E$. $M$ is clearly a majorant
of $C$. 

In fact, $M$ is in the partially ordered system. For if $x,y\in M$,
there exist $E_{x}$ and $E_{y}$ in $C$ so that $x\in E_{x}$ and
$y\in E_{y}$. Since $C$ is a chain, we may assume $E_{x}\leq E_{y}$.
Hence $x,y\in E_{y}$, and so $x\perp y$. This shows that $M$ is
in the partially ordered system and a majorant.

By Zorn's lemma, there exists a maximal element $m\in\mathbb{P}\left(S\right)$.
It remains to show that the closed span of $m$ is $\mathscr{H}$.
Suppose this is false, then by \lemref{dense} there exists a vector
$x\in\mathscr{H}$ so that $x\perp\overline{span}\left\{ m\right\} $. 

Since $m\cup\{x\}\geq m$, and $m$ is assumed maximal, it follows
that $x\in m$. This implies $x\perp x$. Therefore $x=0$, by the
positivity axiom of the definition of Hilbert space.\end{proof}
\begin{cor}
\label{cor:onb}Let $\mathscr{H}$ be a Hilbert space, then $\mathscr{H}$
is isomorphic to the $l^{2}$ space of the index set of an ONB of
$\mathscr{H}$. Specifically, given an ONB $\left\{ u_{\alpha}\right\} _{\alpha\in J}$
in $\mathscr{H}$, where $J$ is some index set, then 
\begin{equation}
v=\sum_{\alpha\in J}\left\langle u_{\alpha},v\right\rangle u_{\alpha},\;\mbox{and}\label{eq:h4}
\end{equation}
\begin{equation}
\left\Vert v\right\Vert ^{2}=\sum_{\alpha\in J}\left|\left\langle u_{\alpha},v\right\rangle \right|^{2},\;\forall v\in\mathscr{H}.\label{eq:h5}
\end{equation}
Moreover, 
\begin{equation}
\left\langle u,v\right\rangle =\sum_{\alpha\in J}\left\langle u,u_{\alpha}\right\rangle \left\langle u_{\alpha},v\right\rangle ,\;\forall u,v\in\mathscr{H}.\label{eq:h6}
\end{equation}

In Dirac's notation (see \secref{dirac}), (\ref{eq:h4})-(\ref{eq:h6})
can be written in the following operator identity: 
\begin{equation}
I_{\mathscr{H}}=\sum_{\alpha\in J}\left|u_{\alpha}\left\rangle \right\langle u_{\alpha}\right|.
\end{equation}
(Note: eq. (\ref{eq:h5}) is called the Parseval identity.)
\end{cor}
\index{Parseval identity}\index{space!Hilbert-}
\begin{proof}
Set $\mathscr{H}_{0}:=span\left\{ u_{\alpha}\right\} $. Then, for
all $v\in\mathscr{H}_{0}$, we have 
\[
v=\sum_{\text{finite}}\left\langle u_{\alpha},v\right\rangle u_{\alpha}
\]
and 
\[
\left\Vert v\right\Vert ^{2}=\sum_{\text{finite}}\left|\left\langle u_{\alpha},v\right\rangle \right|^{2}.
\]
Thus, the map 
\begin{equation}
\mathscr{H}_{0}\ni v\longmapsto\widehat{v}:=\left(\left\langle u_{\alpha},v\right\rangle \right)\in C_{c}\left(J\right)\label{eq:par3}
\end{equation}
is an isometric isomorphism; where $C_{c}$ denotes all the $l^{2}$-sequences
indexed by $J$, vanishing outside some finite subset of $J$.

Since $\mathscr{H}_{0}$ is dense in $\mathscr{H}$, and $C_{c}\left(J\right)$
is dense in $l^{2}\left(J\right)$, it follows that (\ref{eq:par3})
extends to a unitary operator from $\mathscr{H}$ onto $l^{2}\left(A\right)$,
see  \exerref{Fischer}. Thus, (\ref{eq:h4})-(\ref{eq:h5}) hold. 

Using the \emph{polarization identity} (\lemref{pi}) in both $\mathscr{H}$
and $l^{2}\left(J\right)$, we conclude that 
\begin{eqnarray*}
\left\langle u,v\right\rangle _{\mathscr{H}} & = & \frac{1}{4}\sum_{k=0}^{3}i^{k}\left\Vert v+i^{k}u\right\Vert _{\mathscr{H}}^{2}\\
 & \underset{\text{\ensuremath{\left(\ref{eq:par3}\right)}}}{=} & \frac{1}{4}\sum_{k=0}^{3}i^{k}\left\Vert \widehat{v}+i^{k}\widehat{u}\right\Vert _{l^{2}\left(J\right)}^{2}\\
 & = & \left\langle \widehat{u},\widehat{v}\right\rangle _{l^{2}\left(J\right)}\\
 & = & \sum_{j\in J}\left\langle u,u_{\alpha}\right\rangle _{\mathscr{H}}\left\langle u_{\alpha},v\right\rangle _{\mathscr{H}},\;\forall u,v\in\mathscr{H},
\end{eqnarray*}
which is the assertion in (\ref{eq:h6}).\end{proof}
\begin{xca}[Fischer]
\label{exer:Fischer}\myexercise{Fischer}Fix an ONB $\left\{ u_{\alpha}\right\} _{\alpha\in J}$
as in \corref{onb}, and set 
\[
Tv=\left(\left\langle u_{\alpha},v\right\rangle _{\mathscr{H}}\right)_{\alpha\in J}.
\]
Then show that $T:\mathscr{H}\longrightarrow l^{2}\left(J\right)$
is a unitary isomorphism of $\mathscr{H}$ \uline{onto} $l^{2}\left(J\right)$.\end{xca}
\begin{rem}
The correspondence $\mathscr{H}\longleftrightarrow l^{2}$ (index
set of an ONB) is \emph{functorial}, and an \emph{isomorphism}. Hence,
there seems to be just one Hilbert space. But this is misleading,
because numerous interesting realizations of an abstract Hilbert space
come in when we make a choice of the ONB. The question as to which
Hilbert space to use is equivalent to a good choice of an ONB; in
$L^{2}\left(\mathbb{R}\right)$, for example, a wavelet ONB.
\end{rem}

\begin{defn}
A Hilbert space $\mathscr{H}$ is said to be \emph{separable} iff
(Def.) it has an ONB with cardinality of $\mathbb{N}$, (this cardinal
is denoted $\aleph_{0}$).
\end{defn}
Many theorems stated first in the separable case also carry over to
non-separable; but in the more general cases, there are both surprises,
and, in some cases, substantial technical (set-theoretic) complications.

As a result, we shall make the blanket assumption that our Hilbert
spaces \uline{are} separable, (unless stated otherwise.)
\begin{xca}[ONBs and cardinality]
\myexercise{ONBs and cardinality}\label{exer:card}Let $\mathscr{H}$
be a Hilbert space, not necessarily assumed separable, and let $\left\{ u_{\alpha}\right\} _{\alpha\in A}$
and $\left\{ v_{\beta}\right\} _{\beta\in B}$ be two ONBs for $\mathscr{H}$.

Show that $A$ and $B$ have the same cardinality, i.e., that there
is a set-theoretic bijection of $A$ onto $B$.\end{xca}
\begin{defn}
~
\begin{enumerate}
\item Let $A$ be a set, and $p:A\rightarrow\mathbb{R}_{+}$ a function
on $A$. We say that the sum $\sum_{\alpha\in A}p\left(\alpha\right)$
is well-defined and finite iff (Def.) 
\[
\sup_{F\subset A,\,F\:\text{finite}}\sum_{\alpha\in F}p\left(\alpha\right)<\infty;
\]
and we set $\sum_{\alpha\in A}p\left(\alpha\right)$ equal to this
supremum. 
\item Let $A$ be a set. By $l^{2}\left(A\right)$ we mean the set of functions
$f:A\rightarrow\mathbb{C}$, such that 
\[
\sum_{\alpha\in A}\left|f\left(\alpha\right)\right|^{2}<\infty.
\]

\end{enumerate}
\end{defn}
\begin{xca}[$l^{2}\left(A\right)$]
\myexercise{$l^2({A})$}\label{exer:L2-1}Let $A$ be a set (general;
not necessarily countable) then show that $l^{2}\left(A\right)$ is
a Hilbert space.

\uline{Hint}: For $f,g\in l^{2}\left(A\right)$, introduce the
inner product $\sum_{\alpha\in A}\overline{f\left(\alpha\right)}g\left(\alpha\right)$,
by using Cauchy-Schwarz for every finite subset of $A$.
\end{xca}
\index{inequality!Cauchy-Schwarz}

\begin{xca}[A functor from sets to Hilbert space]
\myexercise{A functor from sets to Hilbert space}\label{exer:L2-2}Let
$A$ and $B$ be sets, and let $\psi:A\rightarrow B$ be a bijective
function, then show that there is an induced unitary isomorphism of
$l^{2}\left(A\right)$ onto $l^{2}\left(B\right)$.\end{xca}
\begin{example}[Wavelets]
Suppose $\mathscr{H}$ is separable (i.e., having a countable ONB),
for instance let $\mathscr{H}=L^{2}\left(\mathbb{R}\right)$. Then
\[
\mathscr{H}\cong l^{2}\left(\mathbb{N}\right)\cong l^{2}\left(\mathbb{N}\times\mathbb{N}\right),
\]
and it follows that potentially we could choose a doubly indexed basis
\[
\left\{ \psi_{j,k}\::\:j,k\in\mathbb{N}\right\} 
\]
for $L^{2}\left(\mathbb{R}\right)$. It turns out that this is precisely
the setting of wavelet basis! What's even better is that in the $l^{2}$
space, there are all kinds of diagonalized operators, which correspond
to selfadjoint (or normal) operators in $L^{2}$. Among these operators
in $L^{2}$, we single out the following two:
\begin{align}
\mbox{scaling:} & f\left(x\right)\xrightarrow{\;U_{j}\;}2^{j/2}f\left(2^{j}x\right)\label{eq:wa1}\\
\mbox{translation:} & f\left(x\right)\xrightarrow{\;V_{k}\;}f\left(x-k\right)\label{eq:wa2}
\end{align}
for all $j,k\in\mathbb{Z}$. However, $U_{j}$ and $V_{k}$ are NOT
diagonalized \emph{simultaneously} though. See below for details!\index{Hilbert space!$l^{2}$}\index{space!Hilbert-}\end{example}
\begin{rem}
The two unitary actions $U_{j}$ and $V_{k}$, $j,k\in\mathbb{Z}$,
in (\ref{eq:wa1}) and (\ref{eq:wa2}) satisfy the following important
commutation relation:
\begin{equation}
V_{k}U_{j}=U_{j}V_{2^{j}k};\label{eq:waa1}
\end{equation}
or equivalent:
\begin{equation}
U_{j}^{-1}V_{k}U_{j}=V_{2^{j}k}.\label{eq:waa2}
\end{equation}
Verify details!\end{rem}
\begin{defn}
We say a rational number is a \emph{dyadic fraction} or \emph{dyadic
rational} if it has the form of $\frac{a}{2^{b}}$, where $a\in\mathbb{Z}$,
and $b\in\mathbb{N}$. \index{dyadic rationals}
\end{defn}

In the language of groups, the pair in (\ref{eq:waa1}) \& (\ref{eq:waa2})
forms a representation of a \emph{semidirect product}; or, equivalently,
of the discrete dyadic $ax+b$ group (see \secref{em} for more details):
The latter group consists of all $2\times2$ matrices 
\[
\begin{bmatrix}2^{j} & \frac{k}{2^{l}}\\
0 & 1
\end{bmatrix};\quad j,k\in\mathbb{Z},\;l\in\mathbb{N}.
\]
This group is often referred to as one of the \emph{Baumslag--Solitar
groups}; see e.g., \cite{MR3220588,DJ08}. 

\index{groups!Baumslag-Solitar}

\index{groups!$ax+b$}

\index{product!semidirect}

\subsection{\label{sec:bdd}Bounded Operators in Hilbert Space}
\begin{defn}
A bounded operator in a Hilbert space $\mathscr{H}$ is a linear mapping
$T:\mathscr{H}\rightarrow\mathscr{H}$ such that 
\[
\left\Vert T\right\Vert :=\sup\left\{ \left\Vert Tx\right\Vert _{\mathscr{H}}:\left\Vert x\right\Vert _{\mathscr{H}}\leq1\right\} <\infty.
\]
We denote by $\mathscr{B}\left(\mathscr{H}\right)$ the algebra of
all bounded operators in $\mathscr{H}$. 
\end{defn}
Setting $\left(ST\right)\left(v\right)=S\left(T\left(v\right)\right)$,
$v\in\mathscr{H}$, $S,T\in\mathscr{B}\left(\mathscr{H}\right)$,
we have 
\[
\left\Vert ST\right\Vert \leq\left\Vert S\right\Vert \left\Vert T\right\Vert .
\]

\begin{figure}
\includegraphics[width=0.6\textwidth]{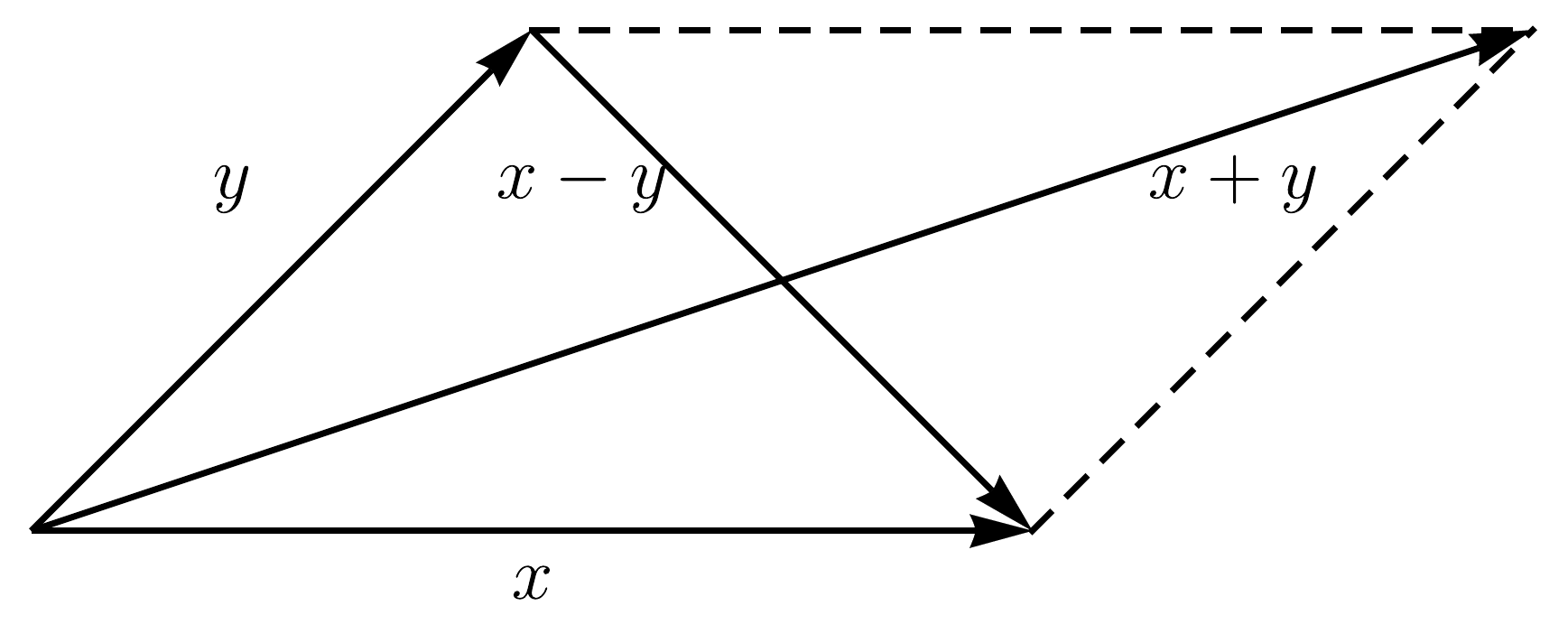}

\protect\caption{\label{fig:paral}The parallelogram law.}
\end{figure}

\begin{figure}
\includegraphics[width=0.7\textwidth]{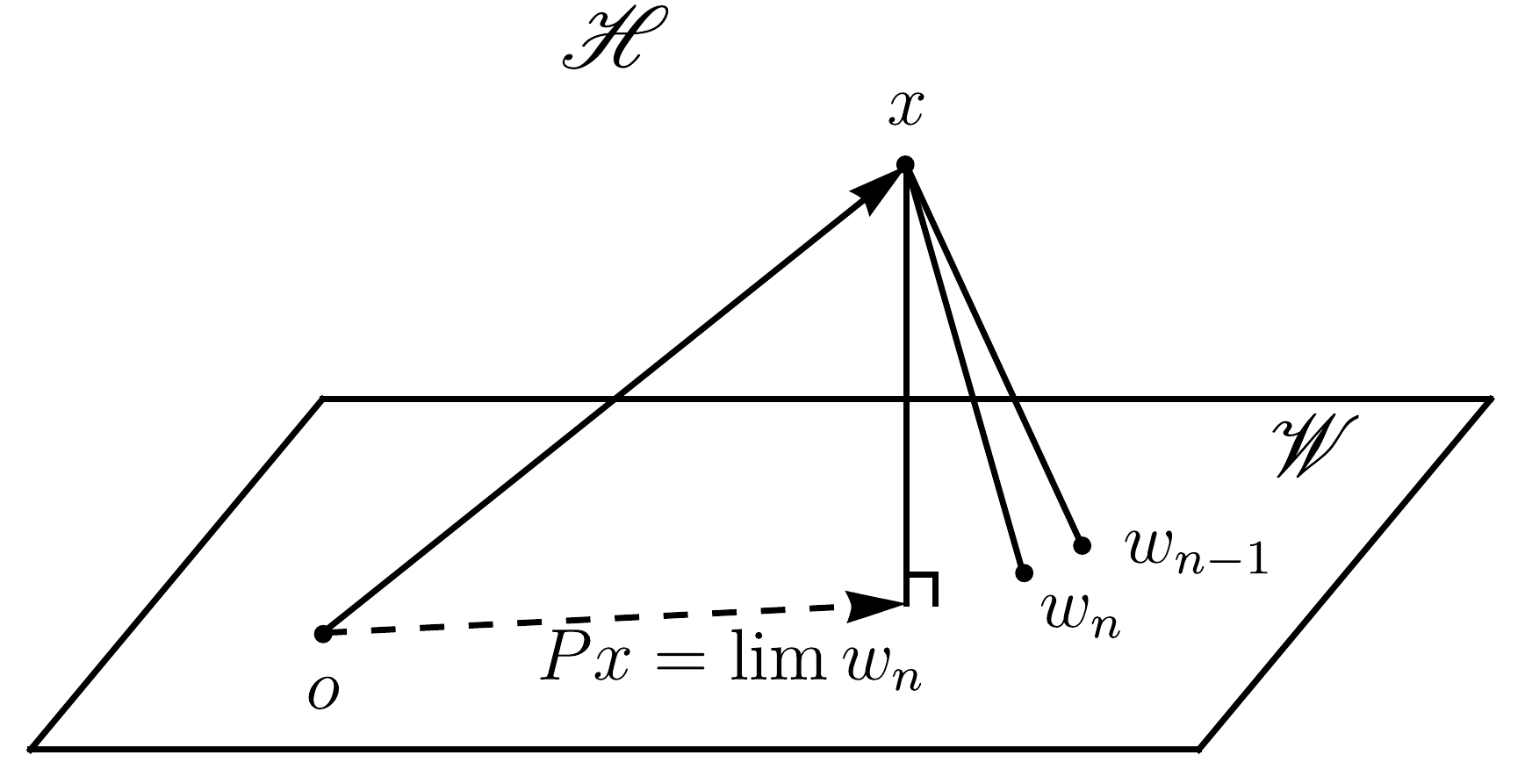}

\protect\caption{\label{fig:proj}$\left\Vert x-Px\right\Vert =\inf\left\{ \left\Vert x-w\right\Vert :w\in\mathscr{W}\right\} $.
Projection from optimization in Hilbert space. }
\end{figure}

\begin{lem}[Riesz]
There is a bijection $\mathscr{H}\ni h\longmapsto l_{h}$ between
$\mathscr{H}$ and the space of all bounded linear functionals on
$\mathscr{H}$, where 
\begin{eqnarray*}
l_{h}\left(x\right) & := & \left\langle h,x\right\rangle ,\;\forall x\in\mathscr{H},\;\mbox{and}\\
\left\Vert l_{h}\right\Vert  & := & \sup\left\{ \left|l\left(x\right)\right|:\left\Vert x\right\Vert _{\mathscr{H}}\leq1\right\} <\infty.
\end{eqnarray*}
Moreover, $\left\Vert l_{h}\right\Vert =\left\Vert h\right\Vert $.
\index{functional} \index{Riesz' theorem}\index{Theorem!Riesz-}\end{lem}
\begin{cor}
For all $T\in\mathscr{B}\left(\mathscr{H}\right)$, there exists a
unique operator $T^{*}\in\mathscr{B}\left(\mathscr{H}\right)$, called
the adjoint of $T$, such that 
\[
\left\langle x,Ty\right\rangle =\left\langle T^{*}x,y\right\rangle ,\;\forall x,y\in\mathscr{H};
\]
and $\left\Vert T^{*}\right\Vert =\left\Vert T\right\Vert $. \end{cor}
\begin{proof}
Let $T\in\mathscr{B}\left(\mathscr{H}\right)$, then it follows from
the Cauchy-Schwarz inequality that 
\[
\left|\left\langle x,Ty\right\rangle \right|\leq\left\Vert Tx\right\Vert \left\Vert y\right\Vert \leq\left\Vert T\right\Vert \left\Vert x\right\Vert \left\Vert y\right\Vert .
\]
Hence the mapping $y\longmapsto\left\langle x,Ty\right\rangle $ is
a bounded linear functional on $\mathscr{H}$. By Riesz's theorem,
there exists a unique $h_{x}\in\mathscr{H}$, such that $\left\langle x,Ty\right\rangle =\left\langle h_{x},y\right\rangle $,
for all $y\in\mathscr{H}$. Set $T^{*}x:=h_{x}$. One checks that
$T^{*}$ linear, bounded, and in fact $\left\Vert T^{*}\right\Vert =\left\Vert T\right\Vert $.
\end{proof}
\index{inequality!Cauchy-Schwarz}
\begin{xca}[The $C^{*}$ property]
\myexercise{The $C^{*} property$}\label{exer:Tnorm}Let $T\in\mathscr{B}\left(\mathscr{H}\right)$,
then prove 
\begin{equation}
\left\Vert T^{*}T\right\Vert =\left\Vert T\right\Vert ^{2}.\label{eq:cnorm}
\end{equation}
\end{xca}
\begin{defn}
\label{def:op}Let $T\in\mathscr{B}\left(\mathscr{H}\right)$. Then,
\begin{itemize}
\item $T$ is \emph{normal} if $TT^{*}=T^{*}T$
\item $T$ is \emph{selfadjoint} if $T=T^{*}$
\item $T$ is \emph{unitary} is $T^{*}T=TT^{*}=I_{\mathscr{H}}\left(=\mbox{the identity operator}\right)$ 
\item $T$ is a (selfadjoint) \emph{projection} if $T=T^{*}=T^{2}$
\end{itemize}
\end{defn}
For $T\in\mathscr{B}\left(\mathscr{H}\right)$, we may write 
\begin{eqnarray*}
R & = & \frac{1}{2}\left(T+T^{*}\right)\\
S & = & \frac{1}{2i}\left(T-T^{*}\right)
\end{eqnarray*}
then both $R$ and $S$ are selfadjoint, and\index{operators!adjoint-}
\[
T=R+iS.
\]
This is similar the to decomposition of a complex number into its
real and imaginary parts. Notice also that $T$ is normal if and only
if $R$ and $S$ commute. (Prove this!) Thus the study of a family
of commuting normal operators is equivalent to the study of a family
of commuting selfadjoint\index{operators!selfadjoint} operators. 
\begin{xca}[The group of all unitary operators ]
\myexercise{The group of all unitary operators}Let $\mathscr{H}$
be a fixed Hilbert space, and denote by $G_{\mathscr{H}}$ the unitary
operators in $\mathscr{H}$ (see \defref{op}).
\begin{enumerate}
\item Show that $G_{\mathscr{H}}$ is a group, and that $T^{-1}=T^{*}$
for all $T\in G_{\mathscr{H}}$.
\item Let $\left\{ u_{\alpha}\right\} _{\alpha\in J}$ be an ONB in $\mathscr{H}$,
and let $v_{\alpha}:=T\left(u_{\alpha}\right)$, $\alpha\in J$; then
show that $\left\{ v_{\alpha}\right\} _{\alpha\in J}$ is also an
ONB.
\item Show that, for any pair of ONBs $\left\{ u_{\alpha}\right\} _{J}$
, $\left\{ w_{\alpha}\right\} _{J}$ with the same index set $J$,
there is then a unique $T\in G_{\mathscr{H}}$ such that $w_{\alpha}=T\left(u_{\alpha}\right)$,
$\alpha\in J$. We say that $G_{\mathscr{H}}$ \emph{acts transitively}
on the set of all ONBs in $\mathscr{H}$.
\end{enumerate}
\end{xca}
\begin{thm}
\label{thm:Hproj}Let $\mathscr{H}$ be a Hilbert space. There is
a one-to-one correspondence between selfadjoint projections and closed
subspaces of $\mathscr{H}$ (\figref{proj}), 
\[
\left[\mbox{Closed subspace}\;\mathscr{M}\subset\mathscr{H}\right]\longleftrightarrow\mbox{Projections}.
\]
\end{thm}
\begin{proof}
Let $P$ be a selfadjoint projection in $\mathscr{H}$, i.e., $P^{2}=P=P^{*}$.
Then 
\[
\mathscr{M}=P\mathscr{H}=\left\{ x\in\mathscr{H}:Px=x\right\} 
\]
is a closed subspace in $\mathscr{H}$. Let $P^{\perp}:=I-P$ be the
complement of $P$, so that 
\[
P^{\perp}\mathscr{H}=\left\{ x\in\mathscr{H}:P^{\perp}x=x\right\} =\left\{ x\in\mathscr{H}:Px=0\right\} .
\]
Since $PP^{\perp}=P(1-P)=P-P^{2}=P-P=0$, we have $P\mathscr{H}\perp P^{\perp}\mathscr{H}$. 

Conversely, let $\mathscr{W}\subsetneq\mathscr{H}$ be a closed subspace.
Note the following ``parallelogram law'' holds: 
\begin{equation}
\left\Vert x+y\right\Vert ^{2}+\left\Vert x-y\right\Vert ^{2}=2(\left\Vert x\right\Vert ^{2}+\left\Vert y\right\Vert ^{2}),\;\forall x,y\in\mathscr{H};\label{eq:paral}
\end{equation}
see  \figref{paral} for an illustration.

Let $x\in\mathscr{H}\backslash\mathscr{W}$, and set 
\[
d:=\inf_{w\in\mathscr{W}}\left\Vert x-w\right\Vert .
\]
The key step in the proof is showing that the infimum is attained;
see \figref{proj}. 

By definition, there exists a sequence $\left\{ w_{n}\right\} $ in
$\mathscr{W}$ so that $\left\Vert w_{n}-x\right\Vert \rightarrow0$
as $n\rightarrow\infty$. Applying (\ref{eq:paral}) to $x-w_{n}$
and $x-w_{m}$, we get
\begin{eqnarray*}
 &  & \left\Vert \left(x-w_{n}\right)+\left(x-w_{m}\right)\right\Vert ^{2}+\left\Vert \left(x-w_{n}\right)-\left(x-w_{m}\right)\right\Vert ^{2}\\
 & = & 2\left(\left\Vert x-w_{n}\right\Vert ^{2}+\left\Vert x-w_{m}\right\Vert ^{2}\right);
\end{eqnarray*}
which simplifies to
\begin{eqnarray}
\left\Vert w_{n}-w_{m}\right\Vert ^{2} & = & 2\left(\left\Vert x-w_{n}\right\Vert ^{2}+\left\Vert x-w_{m}\right\Vert ^{2}\right)-4\left\Vert x-\frac{w_{n}+w_{m}}{2}\right\Vert ^{2}\nonumber \\
 & \leq & 2\left(\left\Vert x-w_{n}\right\Vert ^{2}+\left\Vert x-w_{m}\right\Vert ^{2}\right)-4d.\label{eq:paral1}
\end{eqnarray}
Notice here all we require is $\frac{1}{2}\left(w_{n}+w_{m}\right)\in\mathscr{W}$,
hence the argument carries over if we simply assume $\mathscr{W}$
is a closed convex subset in $\mathscr{H}$. We conclude from (\ref{eq:paral1})
that $\left\Vert w_{n}-w_{m}\right\Vert \rightarrow0$, and so $\left\{ w_{n}\right\} $
is a Cauchy sequence. Since $\mathscr{H}$ is complete, there is a
unique limit, 
\begin{equation}
Px:=\lim_{n\rightarrow\infty}w_{n}\in\mathscr{W}
\end{equation}
and
\begin{equation}
d=\left\Vert x-Px\right\Vert \left(=\inf_{w\in\mathscr{W}}\left\Vert x-w\right\Vert \right).\label{eq:min}
\end{equation}
See \figref{proj}.

Set $P^{\perp}x:=x-Px$. We proceed to verify that $P^{\perp}x\in\mathscr{W}^{\perp}$.
By the minimizing property in (\ref{eq:min}), we have 
\begin{eqnarray}
\left\Vert P^{\perp}x\right\Vert ^{2} & \leq & \left\Vert P^{\perp}x+tw\right\Vert ^{2}\nonumber \\
 & = & \left\Vert P^{\perp}x\right\Vert ^{2}+\left|t\right|^{2}\left\Vert w\right\Vert ^{2}+t\left\langle P^{\perp}x,w\right\rangle +\overline{t}\left\langle w,P^{\perp}x\right\rangle \label{eq:min2}
\end{eqnarray}
for all $t\in\mathbb{C}$, and all $w\in\mathscr{W}$. Assuming $w\neq0$
(the non-trivial case), and setting 
\[
t=-\frac{\left\langle w,P^{\perp}x\right\rangle }{\left\Vert w\right\Vert ^{2}}
\]
in (\ref{eq:min2}), it follows that
\[
0\leq-\frac{\left|\left\langle w,P^{\perp}x\right\rangle \right|^{2}}{\left\Vert w\right\Vert ^{2}}\Longrightarrow\left\langle w,P^{\perp}x\right\rangle =0,\;\forall w\in\mathscr{W}.
\]
This shows that $P^{\perp}x\in\mathscr{W}^{\perp}$, for all $x\in\mathscr{H}$.

For uniqueness, suppose $P_{1}$ and $P_{2}$ both have the stated
properties, then for all $x\in\mathscr{H}$, we have 
\[
x=P_{1}x+P_{1}^{\perp}x=P_{2}x+P_{2}^{\perp}x;\;\mbox{i.e., }
\]
\[
P_{1}x-P_{2}x=P_{2}^{\perp}x-P_{1}^{\perp}x\in\mathscr{W}\cap\mathscr{W}^{\perp}=\left\{ 0\right\} 
\]
thus, $P_{1}x=P_{2}x$, $\forall x\in\mathscr{H}$. 

We leave the rest to the reader. See, e.g., \cite{Rud73}, \cite[p.62]{Ne69}.
\end{proof}

\begin{xca}[Riesz]
\myexercise{Riesz}\label{exer:riesz}As a corollary to \thmref{Hproj},
prove the following version of Riesz' theorem. Let $\mathscr{H}$
be a fixed Hilbert space: \index{Riesz' theorem}\index{Theorem!Riesz-}

For every $l\in\mathscr{H}^{*}$, show that there is a unique $h\left(=h_{l}\right)\in\mathscr{H}$
such that
\begin{equation}
l\left(f\right)=\left\langle h,f\right\rangle ,\;\forall f\in\mathscr{H}.
\end{equation}

\end{xca}

\begin{xca}[Lax-Milgram \cite{MR1892228}]
\label{exer:laxmilgram}\myexercise{Lax-Milgram}Let $B:\mathscr{H}\times\mathscr{H}\rightarrow\mathbb{C}$
be sesquilinear, and suppose there is a finite constant $c$ such
that \sindex[nam]{Lax, P.D., (1926 --)}\index{Theorem!Lax-Milgram-}
\[
\left|B\left(h,k\right)\right|\leq c\left\Vert h\right\Vert \left\Vert k\right\Vert ;
\]
and $b>0$ such that 
\[
\left|B\left(h,h\right)\right|\geq b\left\Vert h\right\Vert ^{2},\;\forall h,k\in\mathscr{H}.
\]

Then prove that, for every $h\in\mathscr{H}$, there is a unique $k\left(=k_{h}\right)\in\mathscr{H}$
such that 
\begin{equation}
\left\langle h,f\right\rangle =B\left(k_{h},f\right),\;\forall f\in\mathscr{H}.
\end{equation}
\end{xca}
\begin{rem}
In view of Riesz, Lax-Milgram is an assertion about $\mathscr{H}^{*}$.
The Lax-Milgram lemma was proved with view to solving elliptic PDEs,
but in \chapref{KS} we give an application to frame expansion. \index{Lax-Milgram}
\end{rem}
\begin{flushleft}
\textbf{The Gram-Schmidt Process}
\par\end{flushleft}

\index{Gram-Schmidt orthogonalization}\index{orthogonal!-vectors}

\begin{figure}
\subfloat[$v_{1}={\displaystyle \frac{u_{1}}{\left\Vert u_{1}\right\Vert }}$]{%
\begin{minipage}[t]{0.45\columnwidth}%
\protect\begin{center}
\protect\includegraphics[scale=0.6]{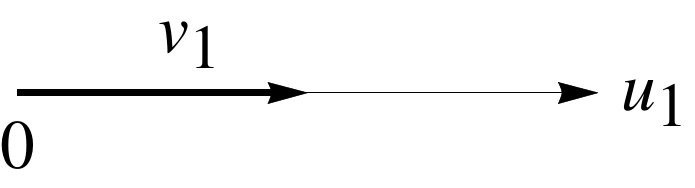}\protect
\par\end{center}%
\end{minipage}

}\hfill{}\subfloat[$v_{2}={\displaystyle \frac{u_{2}-\left\langle v_{1},u_{2}\right\rangle v_{1}}{\left\Vert u_{2}-\left\langle v_{1},u_{2}\right\rangle v_{1}\right\Vert }}$]{%
\begin{minipage}[t]{0.45\columnwidth}%
\protect\begin{center}
\protect\includegraphics[scale=0.6]{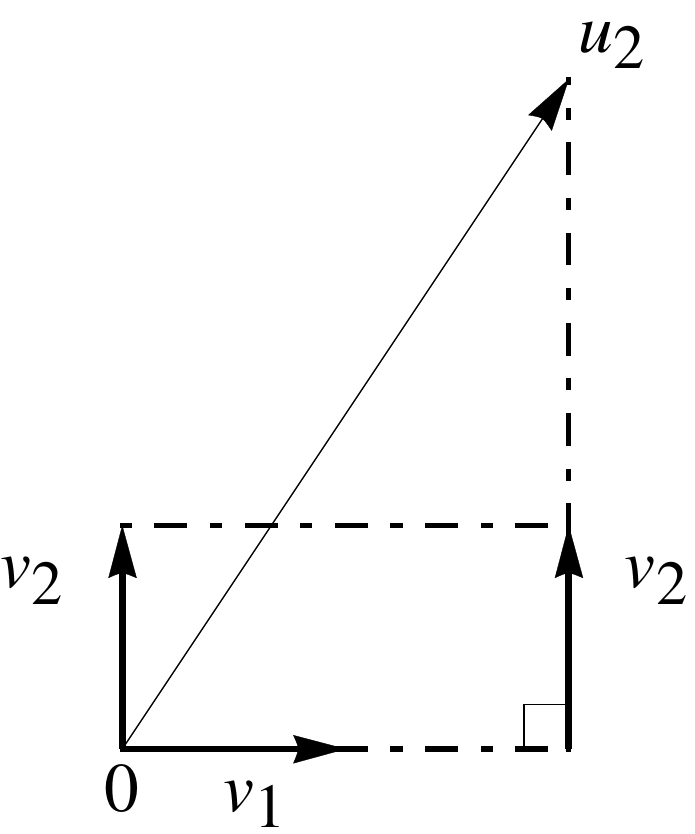}\protect
\par\end{center}%
\end{minipage}

}

\protect\caption{\label{fig:gm}The first two steps in G-S.}
\end{figure}

Every Hilbert space has an ONB, but it does not mean in practice it
is easy to select one that works well for a particular problem. The
Gram-Schmidt orthogonalization process was developed a little earlier
than von Neumann's formulation of abstract Hilbert space. It is an
important tool to get an orthonormal set out of a set of linearly
independent vectors.
\begin{lem}[Gram-Schmidt]
\label{lem:Gram-Schmidt} Let $\left\{ u_{n}\right\} $ be a sequence
of linearly independent vectors in $\mathscr{H}$, then there exists
a sequence $\{v_{n}\}$ of unit vectors so that $\left\langle v_{n},v_{k}\right\rangle =\delta_{n,k}$
and
\[
span\left\{ u_{k}\right\} _{k=1}^{n}=span\left\{ v_{k}\right\} _{k=1}^{n}
\]
for all $n;$ and therefore, 
\[
\overline{span}\left\{ u_{k}\right\} =\overline{span}\left\{ v_{k}\right\} .
\]
\end{lem}
\begin{proof}
Given $\left\{ u_{n}\right\} $ as in the statement of the lemma,
we set 
\begin{eqnarray*}
v_{1} & = & \frac{u_{1}}{\left\Vert u_{1}\right\Vert }.\\
v_{2} & = & \frac{u_{2}-\left\langle v_{1},u_{2}\right\rangle v_{1}}{\left\Vert u_{2}-\left\langle v_{1},u_{2}\right\rangle v_{1}\right\Vert },\cdots.
\end{eqnarray*}
The inductive step: Suppose we have constructed the orthonormal set
$F_{n}:=\left\{ v_{1},\ldots,v_{n}\right\} $, and let $P_{F_{n}}$
be the projection on $F_{n}$ . For the induction step, we set 

\begin{equation}
v_{n+1}:=\frac{u_{n+1}-P_{F_{n}}u_{n+1}}{\left\Vert u_{n+1}-P_{F_{n}}u_{n+1}\right\Vert },\;n=1,2,\ldots\label{eq:gs1}
\end{equation}
See \figref{gm}. Note the LHS in (\ref{eq:gs1}) a unit vector, and
orthogonal to $P_{F_{n}}\mathscr{H}$. 

The formula for $P_{F_{n}}$, the projection onto the span of $F_{n}$,
is 
\[
P_{F_{n}}=\sum_{k=1}^{n}\left|v_{k}\left\rangle \right\langle v_{k}\right|.
\]

\end{proof}

\begin{rem}
If $\mathscr{H}$ is \emph{non-separable}, the standard induction
does not work, and the transfinite induction is needed. \end{rem}
\begin{example}[Legendre (see \tabref{op})]
Let $\mathscr{H}=L^{2}\left(-1,1\right)$. The polynomials $\{1,x,x^{2},\ldots\}$
are linearly independent in $\mathscr{H}$, for if
\[
\sum_{k=1}^{n}c_{k}x^{k}=0
\]
then as an analytic function, the left-hand-side must be identically
zero. By \index{Stone-Weierstrass' theorem}Stone-Weierstrass' theorem,
$span\{1,x,x^{2},\ldots\}$ is dense in $C([-1,1])$ under the $\left\Vert \cdot\right\Vert _{\infty}$
norm. Since $\left\Vert \cdot\right\Vert _{L^{2}}\leq\left\Vert \cdot\right\Vert _{\infty}$,
it follows that $span\{1,x,x^{2},\ldots\}$ is also dense in $\mathscr{H}$. 

By the Gram-Schmidt process, we get a sequence $\{V_{n}\}_{n=1}^{\infty}$
of finite dimensional subspaces in $\mathscr{H}$, where $V_{n}$
has an orthonormal basis $\{h_{0},\ldots,h_{n}\}$, so that 
\begin{eqnarray*}
V_{n} & = & span\{1,x,\ldots,x^{n}\}\\
 & = & span\left\{ h_{0},h_{1},\ldots,h_{n}\right\} .
\end{eqnarray*}
Details: Set $h_{0}=\mathbbm{1}=$ constant function, and 
\[
h_{n+1}:=\frac{x^{n+1}-P_{n}x^{n+1}}{\left\Vert x^{n+1}-P_{n}x^{n+1}\right\Vert },\;n\in\mathbb{N}.
\]
Then the set $\left\{ h_{n}:n\in\mathbb{N}\cup\left\{ 0\right\} \right\} $
is an ONB in $\mathscr{H}$. These are the Legendre polynomials, see
\tabref{op}.
\end{example}
\index{Stone, M. H.}

The two important families of orthogonal polynomials on $\left(-1,1\right)$,
are in \tabref{op} below.

\index{orthogonal polynomial!Chebyshev}

\index{orthogonal polynomial!Legendre}

\index{orthogonal polynomial!Hermite}

\renewcommand{\arraystretch}{2}

\begin{table}
\begin{tabular}{>{\centering}p{0.15\textwidth}|>{\centering}p{0.45\textwidth}|>{\centering}p{0.4\textwidth}}
\hline 
\emph{Name} & \emph{Hilbert space} & \emph{List}\tabularnewline
\hline 
Legendre & $L^{2}\left(-1,1\right)$ & $P_{0}\left(x\right)=1$\tabularnewline
 & $\left\Vert f\right\Vert _{L}^{2}=\int_{-1}^{1}\left|f\left(x\right)\right|^{2}dx$,
$-1\leq x\leq1$ & $P_{1}\left(x\right)=x$\tabularnewline
 &  & $P_{2}\left(x\right)=\frac{1}{2}\left(3x^{2}-1\right)$\tabularnewline
 & \multirow{2}{0.45\textwidth}{orthogonal relation
\[
\mbox{}\int_{-1}^{1}P_{n}\left(x\right)P_{k}\left(x\right)dx=\delta_{n,k}\frac{2}{2n+1}
\]
} & $P_{3}\left(x\right)=\frac{1}{2}\left(5x^{3}-3x\right)$\tabularnewline
 &  & $P_{4}\left(x\right)=\frac{1}{8}\left(35x^{4}-30x^{2}+3\right)$\tabularnewline
 &  & $\vdots$\tabularnewline
\hline 
Chebyshev & $L^{2}\left(\left(-1,1\right);\dfrac{dx}{\sqrt{1-x^{2}}}\right)$ & $P_{n+1}\left(x\right)=2xP_{n}\left(x\right)-P_{n-1}\left(x\right)$\tabularnewline
 & $\left\Vert f\right\Vert _{C}^{2}=\int_{-1}^{1}\left|f\left(x\right)\right|^{2}\dfrac{dx}{\sqrt{1-x^{2}}}$ & $P_{n}\left(\cos\theta\right)=\cos\left(n\theta\right)$\tabularnewline
 & $=\int_{0}^{\pi}\left|f\left(\theta\right)\right|^{2}d\theta$, $x\in\left[-1,1\right]$ & $P_{0}\left(x\right)=1$\tabularnewline
 &  & $P_{1}\left(x\right)=x$\tabularnewline
 & \multirow{2}{0.45\textwidth}{orthogonal relation
\begin{eqnarray*}
\mbox{} &  & \int_{-1}^{1}\frac{P_{n}\left(x\right)P_{k}\left(x\right)}{\sqrt{1-x^{2}}}dx\\
 & = & \begin{cases}
0 & n\neq k\\
\frac{\pi}{2} & n=k\neq0\\
\pi & n=k=0
\end{cases}
\end{eqnarray*}
} & $P_{2}\left(x\right)=2x^{2}-1$\tabularnewline
 &  & $P_{3}\left(x\right)=4x^{3}-3x$\tabularnewline
 &  & $\vdots$\tabularnewline
 &  & \tabularnewline
\hline 
Hermite  & $L^{2}\left(\mathbb{R},e^{-x^{2}}dx\right)$ & (physics version)\tabularnewline
 & $\left\Vert f\right\Vert _{H}^{2}=\int_{-\infty}^{\infty}\left|f\left(x\right)\right|^{2}e^{-x^{2}}dx$,
$x\in\mathbb{R}$ & $P_{n}\left(x\right)=\left(-1\right)^{n}e^{x^{2}}\left(\frac{d}{dx}\right)^{n}e^{-x^{2}}$\tabularnewline
 &  & $P_{0}\left(x\right)=1$\tabularnewline
 & \multirow{2}{0.45\textwidth}{orthogonal relation 
\begin{eqnarray*}
\mbox{} &  & \ensuremath{\int_{-\infty}^{\infty}P_{n}\left(x\right)P_{m}\left(x\right)e^{-x^{2}}dx}\\
 & = & \sqrt{\pi}\;2^{n}\:n!\:\delta_{n,m}
\end{eqnarray*}
} & $P_{1}\left(x\right)=2x$\tabularnewline
 &  & $p_{2}\left(x\right)=4x^{2}-2$\tabularnewline
 & \multirow{2}{0.45\textwidth}{} & $P_{3}\left(x\right)=8x^{3}-12x$\tabularnewline
 &  & $\vdots$\tabularnewline
\end{tabular}

\protect\caption{\label{tab:op}ORTHOGONAL POLYNOMIALS. Legendre, Chebyshev, Hermite}
\end{table}

\renewcommand{\arraystretch}{1}
\begin{defn}
Let $\left\{ P_{n}\left(x\right)\right\} _{n\in\left\{ 0\right\} \cup\mathbb{N}}$
be a sequence of polynomials. We say that the expansion 
\[
G_{P}\left(x,t\right)=\sum_{n=0}^{\infty}P_{n}\left(x\right)t^{n}
\]
is the corresponding \emph{generating function}. \end{defn}
\begin{xca}[Generating functions]
\myexercise{Generating functions} see \tabref{op}. Show that the
generating functions for the three cases of orthogonal polynomials,
\emph{Legendre} (pg.\pageref{tab:op}), \emph{Chebyshev} (pg.\pageref{tab:op}),
and \emph{Hermite }(pg.\pageref{tab:op}) are as follows:
\begin{eqnarray*}
G_{L}\left(x,t\right) & = & \frac{1}{\sqrt{1-2xt+t^{2}}},\\
G_{C}\left(x,t\right) & = & \frac{1-xt}{1-2xt+t^{2}},\;\mbox{and}\\
G_{H}\left(x,t\right) & = & \sum_{\underset{(\text{modified})}{n=0}}^{\infty}P_{n}^{\left(H\right)}\left(x\right)\frac{t^{n}}{n!}=\exp\left(2xt-t^{2}\right).
\end{eqnarray*}
See \cite{MR0184042}.
\end{xca}

\begin{xca}[Recursive identities]
\myexercise{Recursive identities} Verify the following recursive
identities for the three classes of orthogonal polynomials:

\emph{Legendre}:
\[
\left(n+1\right)P_{n+1}\left(x\right)=\left(2n+1\right)xP_{n}\left(x\right)-nP_{n-1}\left(x\right)
\]

\emph{Chebyshev}:
\begin{eqnarray*}
P_{n+1}\left(x\right) & = & 2xP_{n}\left(x\right)-P_{n-1}\left(x\right),\\
2P_{m}\left(x\right)P_{n}\left(x\right) & = & P_{m+n}\left(x\right)+P_{m-n}\left(x\right)
\end{eqnarray*}

\emph{Hermite}:
\begin{eqnarray*}
P_{n+1}\left(x\right) & = & 2xP_{n}\left(x\right)-P_{n}'\left(x\right)\;(\mbox{derivative})\\
 & = & 2xP_{n}\left(x\right)-2nP_{n-1}\left(x\right),\\
P_{n}\left(x+y\right) & = & 2^{-\frac{n}{2}}\sum_{k=0}^{n}\binom{n}{k}P_{n-k}(x\sqrt{2})P_{k}(y\sqrt{2}).
\end{eqnarray*}

\end{xca}

\begin{xca}[Legendre, Chebyshev, Hermite, and Jacobi]
\myexercise{Legendre, Chebyshev, Hermite, and Jacobi} see \tabref{op}.
Find the three $\infty\times\infty$ Jacobi matrices $J$ associated
with the three systems of polynomials in \tabref{op}.

\uline{Hint}: Before writing down the three respective matrices
$J$, you must first normalize the polynomials with respect to the
respective Hilbert norms. See also \cite{Sho36}.

We shall return to the Hermite polynomials, and a corresponding system,
the Hermite functions, in Examples \ref{exa:her1} and \ref{exa:her2},
where they are used in a detailed analysis of the canonical commutation
relations, see also  \secref{topics} above; as well as the corresponding
harmonic oscillator Hamiltonian $H$. With the use of raising and
lowering operators, we show that the Hermite functions are eigenfunctions
for $H$, and we derive the spectrum for $H$ this way.
\end{xca}

\begin{xca}[The orthogonality rules]
\myexercise{The orthogonality rules} Verify the orthogonality rules
contained in \tabref{op}. 

\uline{Hint}: You can use direct computations, a clever system
of recursions, or Fourier transform (generating function).
\end{xca}

\begin{xca}[$M_{x}$ in Jacobi form]
\label{exer:mx}\myexercise{$M_{x}$ in Jacobi form}Let $J\subset\mathbb{R}$
be an interval (finite, or infinite), and let $\left\{ p_{n}\left(x\right)\right\} _{n=0}^{\infty}$
be a system of polynomial functions on $J$; then show that there
is a positive Borel measure $\mu$ on $J$, with infinite support,
but moments of all orders, such that $\left\{ p_{n}\left(x\right)\right\} _{n=0}^{\infty}$
is an ONB in $L^{2}\left(J,\mu\right)$, if and only if, the multiplication
operator $M_{x}$ has an $\infty\times\infty$ matrix representation,
with $\alpha_{n}\in\mathbb{R}$, $\beta_{n}\in\mathbb{C}$, and 
\begin{eqnarray*}
\beta_{1}p_{1}\left(x\right) & = & \left(x-\alpha_{0}\right)p_{0}\left(x\right)\\
\beta_{2}p_{2}\left(x\right) & = & \left(x-\alpha_{1}\right)p_{1}\left(x\right)-\overline{\beta_{1}}p_{0}\left(x\right)\\
 & \vdots\\
\beta_{n+1}p_{n+1}\left(x\right) & = & \left(x-\alpha_{n}\right)p_{n}\left(x\right)-\overline{\beta_{n}}p_{n-1}\left(x\right);
\end{eqnarray*}
i.e., with Jacobi matrix given in \figref{mxj}. 

\begin{figure}
\[
J:=\begin{pmatrix}\alpha_{0} & \overline{\beta_{1}} & 0\\
\beta_{1} & \alpha_{1} & \overline{\beta_{2}} & 0 &  & \smash{\scalebox{4}{0}}\\
0 & \beta_{2} & \alpha_{2} & \overline{\beta_{3}} & \ddots\\
 & \ddots & \ddots & \ddots & \ddots & \ddots\\
 &  & 0 & \beta_{n-1} & \alpha_{n-1} & \overline{\beta_{n}} & 0\\
 &  &  & 0 & \beta_{n} & \alpha_{n} & \overline{\beta_{n+1}} & \ddots\\
 & \smash{\scalebox{4}{0}} &  &  & 0 & \beta_{n+1} & \alpha_{n+1} & \ddots\\
 &  &  &  &  & \ddots & \ddots & \ddots
\end{pmatrix}
\]

\protect\caption{\label{fig:mxj}An $\infty\times\infty$ matrix representation of
$M_{x}$ in \exerref{mx}.}
\end{figure}
\end{xca}
\begin{example}[Fourier basis]
Let $\mathscr{H}=L^{2}[0,1]$. Consider the set of complex exponentials
\[
\left\{ e^{i2\pi nx}:n\in\mathbb{N}\cup\left\{ 0\right\} \right\} ,
\]
or equivalently, one may also consider 
\[
\left\{ 1,\;\cos2\pi nx,\;\sin2\pi nx\;:\;n\in\mathbb{N}\right\} .
\]
This is already an ONB in $\mathscr{H}$ and leads to Fourier series. 
\end{example}
In the next example we construct the \emph{Haar wavelet}. \index{wavelet}
\begin{defn}
A function $\psi\in L^{2}\left(\mathbb{R}\right)$ is said to generate
a wavelet if 
\begin{equation}
\psi_{j,k}\left(x\right)=2^{j/2}\psi\left(2^{j}x-k\right),\;j,k\in\mathbb{Z}\label{eq:w1}
\end{equation}
is an ONB in $L^{2}\left(\mathbb{R}\right)$.
\end{defn}
Note with the normalization in (\ref{eq:w1}) we get 
\[
\int_{\mathbb{R}}\left|\psi_{j,k}\left(x\right)\right|^{2}dx=\int_{\mathbb{R}}\left|\psi\left(x\right)\right|^{2}dx,\;\forall j,k\in\mathbb{Z}.
\]

\begin{figure}
\noindent \begin{centering}
\includegraphics[width=0.7\textwidth]{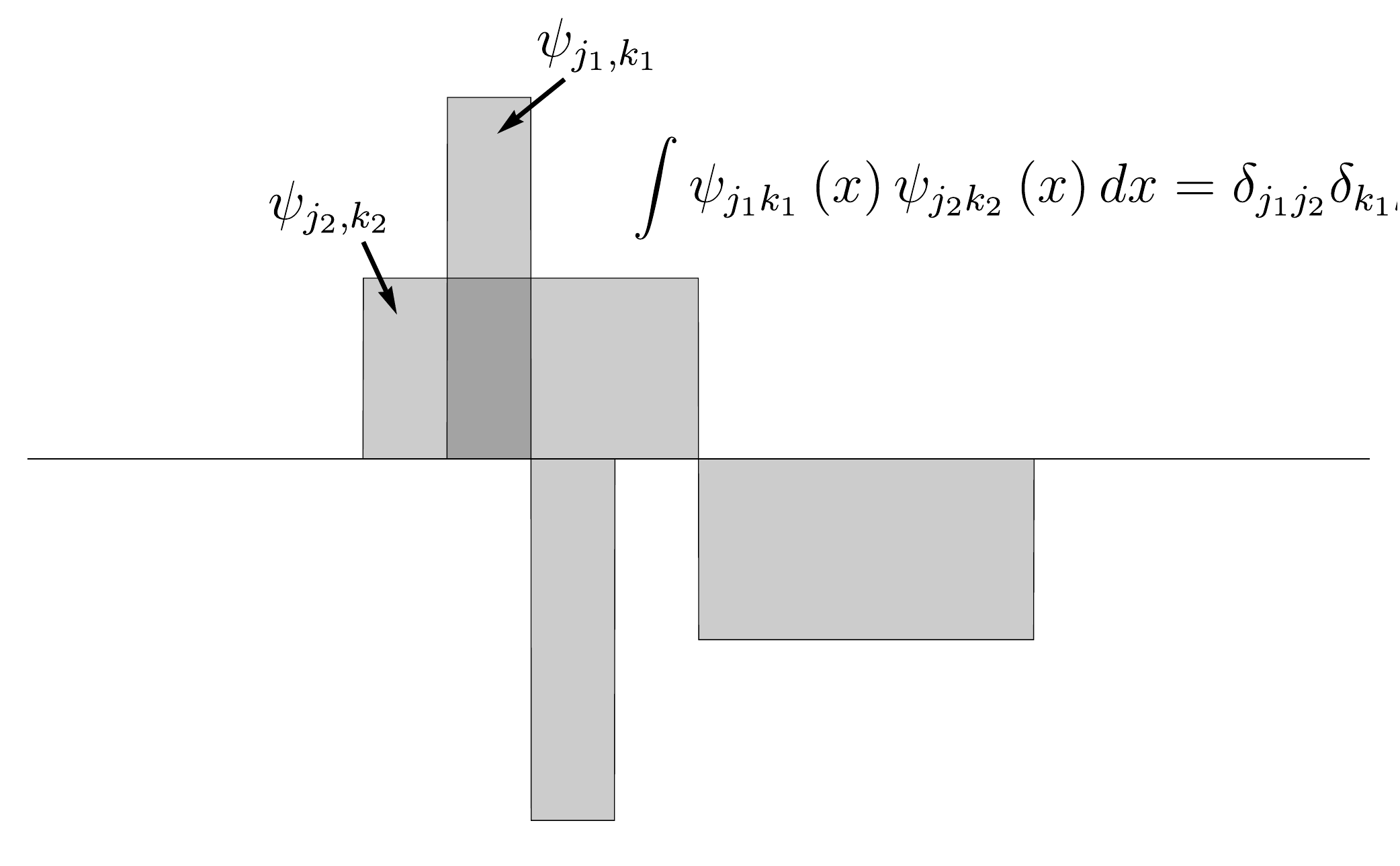}
\par\end{centering}

\protect\caption{\label{fig:haar}Haar wavelet. Scaling properties, resolutions.\index{resolution!multi-}}
\end{figure}

\begin{example}[Haar wavelet and its orthogonality relations]
\label{exa:haar}Let $\mathscr{H}=L^{2}\left(0,1\right)$, and let
$\varphi_{0}$ be the characteristic function of the unit interval
$[0,1]$. $\varphi_{0}$ is called a scaling function. Define 
\begin{eqnarray}
\varphi_{1} & := & \varphi_{0}(2x)-\varphi_{0}(2x-1),\;\mbox{and}\label{eq:hw1}\\
\psi_{j,k} & := & 2^{j/2}\varphi_{1}(2^{j}x-k),\;j,k\in\mathbb{Z}.\label{eq:hw2}
\end{eqnarray}
Claim: $\left\{ \psi_{j,k}:j,k\in\mathbb{Z}\right\} $ is an orthonormal
set (in fact, an ONB) in $L^{2}\left(0,1\right)$, when $j,k$ are
restricted as follow:
\[
k=0,1,\ldots,2^{j}-1,\;j\in\mathbb{N}\cup\left\{ 0\right\} .
\]
\end{example}
\begin{proof}
Fix $k$, if $j_{1}\neq j_{2}$, then $\psi_{j_{1},k}$ and $\psi_{j_{2},k}$
are orthogonal since their supports are nested. For fixed $j$, and
$k_{1}\neq k_{2}$, then $\psi_{j,k_{1}}$ and $\psi_{j,k_{2}}$ have
disjoint supports, and so they are also orthogonal (see \figref{haar}). 
\end{proof}

\begin{xca}[The Haar wavelet, and multiplication by $t$]
\myexercise{The Haar wavelet, and multiplication by $t$}\label{exer:haar}Let
$M:=M_{t}:L^{2}[0,1]\rightarrow L^{2}[0,1]$ be the standard multiplication
operator in $L^{2}$ of the unit-interval. Now compute the $\infty\times\infty$
matrix of $M$ relative to the orthogonal Haar wavelet basis in \exaref{haar}. \end{xca}
\begin{rem}
Even though $M_{t}$ has continuous spectrum $\left[0,1\right]$,
uniform multiplicity, it is of interest to study the diagonal part
in an $\infty\times\infty$ matrix representation of $M_{t}$. Indeed,
in the wavelet ONB in $L^{2}\left(0,1\right)$ we get the following
$\infty\times\infty$ matrix representation \index{spectrum!continuous-}
\index{multiplicity!-of spectrum}
\[
\left(M_{t}\right)_{\left(j_{1},k_{1}\right)\left(j_{2},k_{2}\right)}=\int_{0}^{1}\psi_{j_{1},k_{1}}\left(t\right)t\psi_{j_{2},k_{2}}\left(t\right)dt.
\]
The diagonal part $D$ consists of the sequence
\[
D\left(jk\right)=\int_{0}^{1}t\left(\psi_{j,k}\left(t\right)\right)^{2}dt.
\]

Anderson's theorem \cite{And79} states that $M_{t}-D\in\mathscr{K}$
(the compact operators in $L^{2}\left(0,1\right)$.) Indeed, Anderson
computed the variance
\begin{equation}
V_{jk}=\int_{0}^{1}t^{2}\psi_{j,k}^{2}\left(t\right)\,dt-\left(\int_{0}^{1}t\psi_{j,k}^{2}\left(t\right)\,dt\right)^{2}=\frac{1}{12}2^{-2j}\label{eq:hw3}
\end{equation}
for all $j\in\mathbb{N}\cup\left\{ 0\right\} $, and all $k\in\left\{ 0,1,2,\cdots,2^{j}-1\right\} $.

\uline{Generally}: It is known that if $A$ is a selfadjoint operator
acting in a separable Hilbert space, then $A=D+K$, where $D$ is
a diagonal operator, and $K$ is a compact perturbation \cite{VN35}.
(In fact, there is even a representation $A=D+K$, where $K$ is a
Hilbert-Schmidt operator. See \secref{3norm} and \chapref{sp}.)
\index{selfadjoint operator}

Note that the multiplication operator $M_{t}$ in \exerref{haar}
is bounded and selfadjoint in $L^{2}\left(0,1\right)$. Different
ONBs will yield different diagonal representations $D$. But the wavelet
basis is of special interest.

For more details on compact perturbation of linear operators in Hilbert
space, we refer to \cite{And74}.

The conclusion from Anderson is of special interests as the function
$t$ on $\left[0,1\right]$ is as \textquotedblleft nice\textquotedblright{}
as can be, while the functions from the wavelet ONB (\ref{eq:hw2})
are wiggly, and in fact get increasingly more wiggly as the scaling
degree $j$ in the wavelet ONB tends to infinity. The scaling degree
$j$ is log to the base 2 of the frequency applied to the mother wavelet
function (\ref{eq:hw1}). The conclusion from Anderson is that the
variance numbers (\ref{eq:hw3}) fall off as the inverse square of
the frequency.
\end{rem}

\begin{rem}
It is of interest to ask the analogous questions for other functions
than $t$, and for other wavelet bases, other than the Haar wavelet
basis.\end{rem}
\begin{xca}[A duality]
\myexercise{A duality}\label{exer:uz}Let $z$ be a complex number,
and $P$ be a selfadjoint projection. Show that $U(z)=zP+(I-P)$ is
unitary if and only if $\left|z\right|=1$.

\uline{Hint}: $U(z)U(z)^{*}=U(z)^{*}U(z)=\left|z\right|^{2}P+(I-P)$,
so 
\[
U\left(z\right)\:\mbox{is unitary}\Longleftrightarrow\left|z\right|=1.
\]

\end{xca}

\section{\label{sec:dirac}Dirac's Notation}
\begin{quote}
\textquotedblleft There is a great satisfaction in building good tools
for other people to use.\textquotedblright{} 

--- Freeman Dyson\vspace{2em}
\end{quote}
P.A.M. Dirac was very efficient with notation, and he introduced the
``bra-ket'' vectors \cite{Dir35,Dir47}. 

This Dirac formalism has proved extraordinarily efficient, and it
is widely used in physics. It deserves to be better known in the math
community.
\begin{defn}
Let $\mathscr{H}$ be a Hilbert space with inner product $\left\langle \cdot,\cdot\right\rangle $.
We denote by ``bra'' for vectors $\langle x|$ and ``ket'' for
vectors $|y\rangle$, for $x,y\in\mathscr{H}$. 
\end{defn}
\index{Dirac, P.A.M!bra-ket}

\index{Dirac, P.A.M!ket-bra}

\index{projection}

\index{space!Hilbert-}

With Dirac's notation, our first observation is the following lemma.
\begin{lem}
Let $v\in\mathscr{H}$ be a unit vector. The operator $\:x\mapsto\left\langle v,x\right\rangle v\:$
can be written as $P_{v}=\left|v\left\rangle \right\langle v\right|$,
i.e., a ``ket-bra'' vector. And $P_{v}$ is a rank-one selfadjoint
projection.\end{lem}
\begin{proof}
First, we see that 
\[
P_{v}^{2}=\left(\left|v\left\rangle \right\langle v\right|\right)\left(\left|v\left\rangle \right\langle v\right|\right)=|v\rangle\left\langle v,v\right\rangle \langle v|=\left|v\left\rangle \right\langle v\right|=P_{v}.
\]
Also, if $x,y\in\mathscr{H}$ then 
\[
\left\langle x,P_{v}y\right\rangle =\left\langle x,v\right\rangle \left\langle v,y\right\rangle =\left\langle \overline{\left\langle x,v\right\rangle }v,y\right\rangle =\left\langle \left\langle v,x\right\rangle v,y\right\rangle =\left\langle P_{v}x,y\right\rangle 
\]
so $P_{v}=P_{v}^{*}$. \end{proof}
\begin{cor}
Let $F=span\left\{ v_{i}\right\} $ with $\left\{ v_{i}\right\} $
a finite set of orthonormal vectors in $\mathscr{H}$, then\index{orthogonal!-vectors}
\[
P_{F}:=\sum_{v_{i}\in F}\left|v_{i}\left\rangle \right\langle v_{i}\right|
\]
is the selfadjoint projection onto $F$. \end{cor}
\begin{proof}
Indeed, we have 
\begin{eqnarray*}
P_{F}^{2} & = & \sum_{v_{i},v_{j}\in F}\left(\left|v_{i}\left\rangle \right\langle v_{i}\right|\right)\left(\left|v_{j}\left\rangle \right\langle v_{j}\right|\right)=\sum_{v_{i}\in F}\left|v_{i}\left\rangle \right\langle v_{i}\right|=P_{F}\\
P_{F}^{*} & = & \sum_{v_{i}\in F}\left(\left|v_{i}\left\rangle \right\langle v_{i}\right|\right)^{*}=\sum_{v_{i}\in F}\left|v_{i}\left\rangle \right\langle v_{i}\right|=P_{F},
\end{eqnarray*}
and we have
\begin{eqnarray}
P_{F}w=w & \Longleftrightarrow & \sum_{F}\left|\left\langle v_{i},w\right\rangle \right|^{2}=\left\Vert w\right\Vert ^{2}\label{eq:dp1}\\
 & \Longleftrightarrow & w\in F.\nonumber 
\end{eqnarray}
Since we may take the limit in (\ref{eq:dp1}), it follows that the
Corollary also holds if $F$ is infinite-dimensional, i.e., $F=\overline{span}\left\{ v_{i}\right\} $,
closure. \end{proof}
\begin{rem}
More generally, any rank-one operator can be written in Dirac notation
as 
\[
\left|u\left\rangle \right\langle v\right|:\mathscr{H}\ni x\longmapsto\left\langle v,x\right\rangle u\in\mathscr{H}.
\]
With the bra-ket notation, it is easy to verify that the set of rank-one
operators forms an algebra, which easily follows from the fact that
\[
\left(\left|v_{1}\left\rangle \right\langle v_{2}\right|\right)\left(\left|v_{3}\left\rangle \right\langle v_{4}\right|\right)=\left\langle v_{2},v_{3}\right\rangle \left|v_{1}\left\rangle \right\langle v_{4}\right|.
\]
The moment that an orthonormal basis is selected, the algebra of operators
on $\mathscr{H}$ will be translated to the algebra of matrices (infinite).
See \lemref{mat}.\end{rem}
\begin{xca}[Finite-rank reduction]
\label{exer:frank}\myexercise{Finite rank reduction} Let $\mathscr{H}$
be a Hilbert space. For all $x,y\in\mathscr{H}$, let $\left|x\left\rangle \right\langle y\right|$
denote the corresponding (Dirac) rank-1 operator. \index{operators!finite rank-}
\begin{enumerate}
\item \label{enu:fr1}Let $A,B\in\mathscr{B}\left(\mathscr{H}\right)$.
Verify that
\begin{eqnarray}
A\left|x\left\rangle \right\langle y\right| & = & \left|Ax\left\rangle \right\langle y\right|,\;\mbox{and}\label{eq:rk1}\\
\left|x\left\rangle \right\langle y\right|B & = & \left|x\left\rangle \right\langle B^{*}y\right|.\label{eq:rk2}
\end{eqnarray}
In particular, $\mathscr{F}R\left(\mathscr{H}\right)$ is a two-sided
ideal in $\mathscr{B}\left(\mathscr{H}\right)$.
\item For all $x,y\in\mathscr{H}$, set 
\[
w_{x,y}\left(A\right):=\left\langle x,Ay\right\rangle ,\;\forall A\in\mathscr{B}\left(\mathscr{H}\right).
\]
For $\left\Vert x\right\Vert =1$, set $w_{x}\left(A\right):=\left\langle x,Ax\right\rangle $,
i.e., a \emph{pure state} on $\mathscr{B}\left(\mathscr{H}\right)$.
\\
Let $\left\{ x_{i}\right\} _{i=1}^{n}\subset\mathscr{H}$, $\left\Vert x_{i}\right\Vert =1$,
and set $T$ by
\begin{equation}
T=\sum_{i}\left|x_{i}\left\rangle \right\langle x_{i}\right|.\label{eq:rk3}
\end{equation}
Show that then 
\begin{equation}
TAT=\sum_{i}\sum_{j}w_{x_{i},x_{j}}\left(A\right)\underset{\text{Dirac-rank-1}}{\underbrace{\left|x_{i}\left\rangle \right\langle x_{j}\right|}}.\label{eq:rk4}
\end{equation}
In particular, if $n=1$, and $T=\left|x_{1}\left\rangle \right\langle x_{1}\right|$,
we have: 
\begin{equation}
TAT=w_{x_{1}}\left(A\right)T.\label{eq:rk5}
\end{equation}

\item Use part (\ref{enu:fr1}), and  \thmref{3norms} below, to give a
quick proof that the compact operators form an ideal in $\mathscr{B}\left(\mathscr{H}\right)$.
\index{state!pure-} \index{ideal}
\end{enumerate}
\end{xca}

\begin{xca}[Numerical Range and Toeplitz-Hausdorff]
\label{exer:nrange} \myexercise{Numerical Range and Toeplitz-Hausdorff}The
set 
\[
\left\{ w_{x}\left(A\right)\::\:\left\Vert x\right\Vert =1\right\} \subset\mathbb{C}
\]
is called the \emph{numerical range} of $A$, $NR_{A}$. Show that
$NR_{A}$ is convex. 

\uline{Hint}: Difficult! It is called the \emph{Toeplitz-Hausdorff
theorem}; see e.g., \cite{Hal67,MR0171168}. (There are few assertions
that are true for \uline{all} bounded operators. The Toeplitz-Hausdorff
theorem is one on a short list.)

\index{Toeplitz-Hausdorff Theorem} \index{numerical range}\end{xca}
\begin{lem}
Let $\left\{ u_{\alpha}\right\} _{\alpha\in J}$ be an ONB in $\mathscr{H}$,
then we may write 
\[
I_{\mathscr{H}}=\sum_{\alpha\in J}\left|u_{\alpha}\left\rangle \right\langle u_{\alpha}\right|.
\]
\end{lem}
\begin{proof}
This is equivalent to the decomposition 
\[
v=\sum_{\alpha\in J}\left\langle u_{\alpha},v\right\rangle u_{\alpha},\;\forall v\in\mathscr{H}.
\]

\end{proof}

A selection of ONB makes a representation of the algebra of operators
acting on $\mathscr{H}$ by infinite matrices. We check that, using
Dirac's notation, the algebra of operators really becomes the algebra
of infinite matrices.

For $A\in\mathscr{B}\left(\mathscr{H}\right)$, and $\left\{ u_{i}\right\} $
an ONB, set 
\[
\left(M_{A}\right)_{i,j}=\left\langle u_{i},Au_{j}\right\rangle _{\mathscr{H}}.
\]

\begin{figure}
\includegraphics[width=0.4\textwidth]{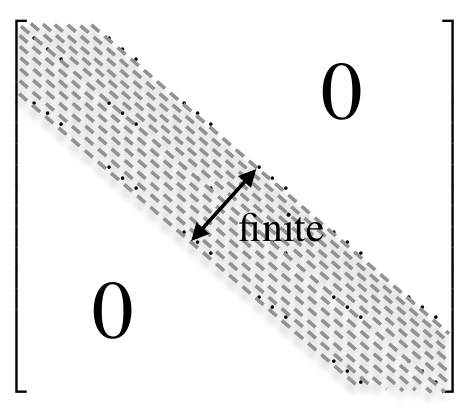}

\protect\caption{\label{fig:mp}$\infty\times\infty$ \emph{banded matrix}. Supported
on a band around the diagonal.}
\end{figure}

Most of the operators we use in the math physics problems are unbounded,
so it is a big deal that the conclusion about matrix product is valid
for unbounded operators subject to the condition that the chosen ONB
is in the domain of such operators.
\begin{lem}[matrix product]
\label{lem:mat}Assume some ONB $\left\{ u_{i}\right\} _{i\in J}$
satisfies $u_{i}\in dom\left(A^{*}\right)\cap dom\left(B\right)$;
then $M_{AB}=M_{A}M_{B}$, i.e., $\left(M_{AB}\right)_{ij}=\sum_{k}\left(M_{A}\right)_{ik}\left(M_{B}\right)_{kj}$.\end{lem}
\begin{proof}
By $AB$ we mean the operator given by 
\[
\left(AB\right)\left(u\right)=A\left(B\left(u\right)\right).
\]

Pick an ONB $\{u_{i}\}$ in $\mathscr{H}$, and the two operators
as stated. We denote by $M_{A}=A_{ij}:=\left\langle u_{i},Au_{j}\right\rangle $
the matrix of $A$ under the ONB. We compute $\left\langle u_{i},ABu_{j}\right\rangle $.
\begin{align*}
(M_{A}M_{B})_{ij}=\sum_{k}A_{ik}B_{kj} & =\sum_{k}\left\langle u_{i},Au_{k}\right\rangle \left\langle u_{k},Bu_{j}\right\rangle \\
 & =\sum_{k}\left\langle A^{*}u_{i},u_{k}\right\rangle \left\langle u_{k},Bu_{j}\right\rangle \\
 & =\left\langle A^{*}u_{i},Bu_{j}\right\rangle \quad[\mbox{by Parseval}]\\
 & =\left\langle u_{i},ABu_{j}\right\rangle \\
 & =\left(M_{AB}\right)_{ij}
\end{align*}
where we used that $I=\sum\left|u_{i}\left\rangle \right\langle u_{i}\right|$.\end{proof}
\begin{xca}[Matrix product of $\infty\times\infty$ banded matrices]
\myexercise{Matrix product of $\infty\times\infty$ banded matrices}Consider
two linear operators $A$ and $B$ both defined on a dense subspace
$\mathscr{D}$ in a fixed Hilbert space $\mathscr{H}$. Suppose $\mathscr{D}$
contains an ONB $\left\{ e_{i}\right\} _{i\in\mathbb{N}}$, and that
the corresponding matrices $M_{A}$ and $M_{B}$ with respect to $\left\{ e_{i}\right\} $
are both banded. Then show that the matrix-product
\begin{equation}
M_{AB}=M_{A}M_{B}\label{eq:mp}
\end{equation}
is well defined, and is again banded. See \figref{mp}.

\uline{Hint}: Use  \lemref{mat}, and the equation 
\[
\left(M_{AB}\right)_{i,j}=\left\langle e_{i},ABe_{j}\right\rangle =\left\langle A^{*}e_{i},Be_{j}\right\rangle ,
\]
and note that $\mathscr{D}\subset dom\left(A^{*}\right)$.
\end{xca}

\begin{rem}
Here are two\emph{ open questions} regarding \emph{banded} operators/matrices.
\begin{enumerate}
\item Fix a separable Hilbert space $\mathscr{H}$, is there an intrinsic
geometric characterization of the linear operators $A$ (with dense
domain) in $\mathscr{H}$ which admit an ONB $\left\{ e_{i}\right\} _{i\in\mathbb{N}}$
such that the matrix
\begin{equation}
\left(M_{A}\right)_{i,j}:=\left\langle e_{i},Ae_{j}\right\rangle \label{eq:oq1}
\end{equation}
is banded?
\item Given an ONB $\left\{ e_{i}\right\} _{i\in\mathbb{N}}$, what is the
$*$-algebra $\mathfrak{A}$ of unbounded operators with dense domain
\begin{equation}
\mathscr{D}=span\left\{ e_{i}\right\} \label{eq:oq2}
\end{equation}
such that every $A\in\mathfrak{A}$ is banded with respect to $\left\{ e_{i}\right\} _{i\in\mathbb{N}}$?
\end{enumerate}
\end{rem}
\begin{xca}[A transform]
\myexercise{A transform}\label{exer:onb}Let $\mathscr{H}$ be a
Hilbert space, and $A$ a set which indexes a fixed ONB $\left\{ v_{\alpha}\right\} _{\alpha\in A}$.
Now define $T:\mathscr{H}\rightarrow l^{2}\left(A\right)$, by 
\[
\left(Th\right)\left(\alpha\right):=\left(\left\langle v_{\alpha},h\right\rangle \right),\;\forall\alpha\in A,\;h\in\mathscr{H}.
\]
Show that $T$ is unitary, and onto $l^{2}\left(A\right)$, i.e.,
\begin{eqnarray*}
TT^{*} & = & I_{l^{2}\left(A\right)},\;\mbox{and}\\
T^{*}T & = & I_{\mathscr{H}}.
\end{eqnarray*}
$T$ is called the analysis transformation, and $T^{*}$ the synthesis
transformation.
\end{xca}

\subsection{\label{sec:3norm}Three Norm-Completions}

Let $\mathscr{H}$ be a fixed Hilbert space, infinite-dimensional
in the discussion below. 

Let 
\begin{eqnarray*}
\mathscr{F}R\left(\mathscr{H}\right) & = & \left\{ \text{all finite-rank operators}\;\mathscr{H}\longrightarrow\mathscr{H}\right\} \\
 & = & span\left\{ \left|v\left\rangle \right\langle w\right|\;:\;v,w\in\mathscr{H}\right\} 
\end{eqnarray*}
where $\left|v\left\rangle \right\langle w\right|$ denotes the Dirac
ket-bra-operator. \index{completion!norm-}
\begin{defn}
\label{def:3norms}On $\mathscr{F}R\left(\mathscr{H}\right)$ we introduce
the following three norms: the \emph{uniform norm} (UN), the \emph{trace-norm}
(TN), and the \emph{Hilbert-Schmidt} norm as follows: \index{trace}
\begin{itemize}
\item (UN) For all $T\in\mathscr{F}R\left(\mathscr{H}\right)$, set 
\[
\left\Vert T\right\Vert _{UN}:=\sup\left\{ \left\Vert Tv\right\Vert \;:\;\left\Vert v\right\Vert =1\right\} .
\]

\item (TN) Set 
\[
\left\Vert T\right\Vert _{TN}:=\text{trace}\left(\sqrt{T^{*}T}\right);
\]

\item (HSN) Set 
\[
\left\Vert T\right\Vert _{HSN}=\left(\text{trace}\left(T^{*}T\right)\right)^{\frac{1}{2}}.
\]

\end{itemize}
\end{defn}
\begin{thm}
\label{thm:3norms}The completion of $\mathscr{F}R\left(\mathscr{H}\right)$
with respect to $\left\Vert \cdot\right\Vert _{UN}$, $\left\Vert \cdot\right\Vert _{TN}$,
and $\left\Vert \cdot\right\Vert _{HSN}$ are respectively the compact
operators, the trace-class operators, and the Hilbert-Schmidt operators.
(See \defref{3norms} above.)
\end{thm}
We will prove, as a consequence of the Spectral Theorem (\chapref{sp}),
that the $\left\Vert \cdot\right\Vert _{UN}$- completion agrees with
the usual definition of the compact operators.

\index{completion!Hilbert-Schmidt}

\index{completion!trace-}
\begin{rem}
Note that \thmref{3norms} is for Hilbert space, and it is natural
to ask \textquotedblleft what carries over to Banach space?\textquotedblright{}
Not everything by a theorem of Per Enflo \cite{Enf73}. In detail:
The assertion in the first part of \thmref{3norms} (Hilbert space)
is that every compact operator is the norm limit of finite-rank operators;
referring to the uniform norm (UN) in \defref{3norms}. But there
are Banach spaces where this is false; -- although the easy implication
is true, i.e., that operators in the norm closure of finite-rank operators
are compact.\end{rem}
\begin{defn}
Let $T\in\mathscr{B}\left(\mathscr{H}\right)$, then we say that $T$
is compact iff (Def.) $T\left(\mathscr{H}_{1}\right)$ is relatively
compact in $\mathscr{H}$, where $\mathscr{H}_{1}:=\left\{ v\in\mathscr{H}\;\big|\;\left\Vert v\right\Vert \leq1\right\} $.

A similar remark applies to the other two Banach spaces of operators.
The Hilbert-Schmidt operators forms a Hilbert space.\end{defn}
\begin{xca}[The identity-operator in infinite dimension]
\myexercise{The identity-operator in infinite dimension}\label{exer:cpt}If
$\dim\mathscr{H}=\infty$, show that the identity operator $I_{\mathscr{H}}$
is \emph{not} compact. 

\uline{Hint}: Use an ONB.\end{xca}
\begin{defn}
The following is useful in working through the arguments above.\end{defn}
\begin{lem}
Let $v,w\in\mathscr{H}$, then 
\[
trace\left(\left|v\left\rangle \right\langle w\right|\right)=\left\langle v,w\right\rangle _{\mathscr{H}}.
\]
 \end{lem}
\begin{proof}
Introduce an ONB $\left\{ u_{\alpha}\right\} $ in $\mathscr{H}$,
and compute:
\begin{eqnarray*}
\text{trace}\left(\left|v\left\rangle \right\langle w\right|\right) & = & \sum_{\alpha}\left\langle u_{\alpha},\:\left|v\left\rangle \right\langle w\right|\:u_{\alpha}\right\rangle _{\mathscr{H}}\\
 & = & \sum_{\alpha}\left\langle u_{\alpha},v\right\rangle _{\mathscr{H}}\left\langle w,u_{\alpha}\right\rangle _{\mathscr{H}}\\
 & = & \left\langle v,w\right\rangle _{\mathscr{H}}
\end{eqnarray*}
where we used Parseval in the last step of the computation.\end{proof}
\begin{xca}[Comparing norms]
\myexercise{Comparing norms}\label{exer:3norms}Let $T\in\mathscr{F}R\left(\mathscr{H}\right)$,
and let the three norms be as in  \defref{3norms}. Then show that
\begin{equation}
\left\Vert T\right\Vert _{UN}\leq\left\Vert T\right\Vert _{HSN}\leq\left\Vert T\right\Vert _{TN};\label{eq:no1}
\end{equation}
and conclude the following contractive inclusions:
\begin{eqnarray*}
\left\{ \text{Trace-class operators}\right\}  & \subset & \left\{ \text{Hilbert-Schmidt operators}\right\} \\
 & \subset & \left\{ \text{compact operators}\right\} .
\end{eqnarray*}
\end{xca}
\begin{rem}
In the literature, the following notation is often used for the three
norms in \defref{3norms}:
\begin{equation}
\begin{cases}
\left\Vert T\right\Vert _{UN}\,\,\,\,=\left\Vert T\right\Vert _{\infty}\\
\left\Vert T\right\Vert _{TN}\,\,\,\,=\left\Vert T\right\Vert _{1}\\
\left\Vert T\right\Vert _{HSN}=\left\Vert T\right\Vert _{2}
\end{cases}\label{eq:no2}
\end{equation}
Note that if $T$ is a diagonal operator, $T=\sum_{k}x_{k}\left|u_{k}\left\rangle \right\langle u_{k}\right|$
in some ONB $\left\{ u_{k}\right\} $, then the respective norms are
$\left\Vert x\right\Vert _{\infty}$, $\left\Vert x\right\Vert _{1}$,
and $\left\Vert x\right\Vert _{2}$.

With the notation in (\ref{eq:no2}), the inequalities (\ref{eq:no1})
now take the form
\[
\left\Vert T\right\Vert _{\infty}\leq\left\Vert T\right\Vert _{2}\leq\left\Vert T\right\Vert _{1},\;T\in\mathscr{F}R\left(\mathscr{H}\right).
\]
\end{rem}
\begin{xca}[Matrix entries in infinite dimensions]
\myexercise{Matrix entries in infinite dimensions}Let $\mathscr{H}$
be a separable Hilbert space. Pick an ONB $\left\{ e_{j}\right\} _{j\in J}$,
and set $E_{ij}:=\left|e_{i}\left\rangle \right\langle e_{j}\right|$,
$\left(i,j\right)\in J^{2}$. 
\begin{enumerate}
\item Show that this is an ONB in $\mathscr{H}S\left(\mathscr{H}\right)$
(= Hilbert Schmidt operators) with respect to the inner product
\begin{equation}
\left\langle A,B\right\rangle _{\mathscr{H}S}:=trace\left(A^{*}B\right),\;A,B\in\mathscr{H}S\left(\mathscr{H}\right).\label{eq:hs2-1}
\end{equation}

\item Show that the corresponding orthogonal expansion for $A\in\mathscr{H}S\left(\mathscr{H}\right)$
is 
\begin{equation}
A=\underset{\left(i,j\right)\in J^{2}}{\sum\sum}\left\langle e_{i},Ae_{j}\right\rangle _{\mathscr{H}}E_{ij}.\label{eq:hs2-2}
\end{equation}

\end{enumerate}
\end{xca}

\begin{xca}[The three steps]
\myexercise{The three steps}Let $\mathscr{H}$ be a separable Hilbert
space, and let 
\begin{enumerate}
\item \label{enu:db1}$\mathscr{B}\left(\mathscr{H}\right)$: all bounded
operators in $\mathscr{H}$; with the uniform norm. 
\item \label{enu:db2}$\mathscr{T}_{1}\left(\mathscr{H}\right)$: all trace-class
operators, with the trace-norm; see \defref{3norms}.
\end{enumerate}

Use the three steps from \secref{duality} to show that
\[
\left(\mathscr{T}_{1}\left(\mathscr{H}\right)\right)^{*}=\mathscr{B}\left(\mathscr{H}\right);
\]
i.e., that $\mathscr{B}\left(\mathscr{H}\right)$ is the dual Banach
space where the respective norms are specified as in (\ref{enu:db1})-(\ref{enu:db2}).
(For more about this duality, see also \thmref{predual}.)

\end{xca}

\subsection{\label{sub:QM}Connection to Quantum Mechanics}

One of the powerful applications of the theory of operators in Hilbert
space, and more generally of functional analysis, is in quantum mechanics
(QM) \cite{MR1939631,MR965583,MR675039}. Even the very formulation
of the central questions in QM entails the mathematics of \emph{unbounded
selfadjoint operators}, of projection valued measures (PVM), and \emph{unitary
one-parameter groups} of operators. The latter usually abbreviated
to \textquotedblleft unitary one-parameter groups.\textquotedblright{}

By contrast to what holds for the more familiar case of bounded operators,
we stress that for unbounded selfadjoint operators, mathematical precision
necessitates making a sharp distinction between the following three
notions: \emph{selfadjoint}, \emph{essentially selfadjoint}, and \emph{formally
selfadjoint} (also called Hermitian, or symmetric). See \cite{RS75, Ne69, vN32a, DS88b}.
We define them below; see especially the appendix at the end of \chapref{lin}.
One reason for this distinction is that quantum mechanical observables,
momentum $P$, position $Q$, energy etc, become selfadjoint operators
in the axiomatic language of QM. What makes it even more subtle is
that these operators are both unbounded and non-commuting (take the
case of $P$ and $Q$ which was pair from Heisenberg\textquoteright s
pioneering paper on uncertainty.) Another subtle point entails the
relationship between selfadjoint operators, projection-valued measures,
and unitary one-parameter groups (as used in the dynamical description
of states in QM, i.e., describing the solution of the wave equation
of Schrödinger.) Unitary one-parameter groups are also used in the
study of other partial differential equations, especially hyperbolic
PDEs.

The discussion which follows below will make reference to this setting
from QM, and it serves as motivation. However the more systematic
mathematical presentation of selfadjoint operators, projection-valued
measures, and unitary one-parameter groups will be postponed to later
in the book. We first need to develop a number of technical tools.
However we have included an outline of the bigger picture in the appendix
(Stone\textquoteright s Theorem), to the present chapter. Stone\textquoteright s
theorem shows that the following three notions, (i) selfadjoint operator,
(ii) projection-valued measure, and (iii) unitary one-parameter group,
are incarnations of one and the same; i.e., when one of the three
is known, anyone of the other two can be computed from it. \index{axioms}

We emphasize that there is a host of other applications of this, for
example to harmonic analysis, to statistics, and to PDE. These will
also play an important role in later chapters.

Much of the motivation for the axiomatic approach to the theory of
linear operators in Hilbert space dates back to the early days of
quantum mechanics (Planck, Heisenberg, and Schrödinger), but in the
form suggested by J. von Neumann. (von Neumann's formulation is the
one now adopted by most books on functional analysis.) Here we will
be brief, as a systematic and historical discussion is far beyond
our present scope. Suffice it to mention here that what is known as
the \textquotedbl{}matrix-mechanics\textquotedbl{} of Heisenberg takes
the form infinite by infinite\index{matrix!inftytimesinfty
@$\infty\times\infty$} matrices with entries representing, in turn, \emph{transition probabilities}\index{transition probabilities},
where \textquotedbl{}transition\textquotedbl{} refers to \textquotedbl{}jumps\textquotedbl{}
between energy levels. See (\ref{eq:HP})-(\ref{eq:HQ}) below, and
Exercises \ref{exer:Hmatrix}-\ref{exer:PQmatrix}. By contrast to
matrices, in Schrödinger's wave mechanics, the Hilbert space represents
wave solutions to Schrödinger's equation. Now this entails the study
of one-parameter groups of unitary\index{operators!unitary} operators
in $\mathscr{H}$. In modern language, with the two settings we get
the dichotomy between the case when the Hilbert space is an $l^{2}$
space (i.e., a $l^{2}$-sequence space), vs the case of Schrödinger
when $\mathscr{H}$ is an $L^{2}$-space of functions on phase-space. 

\index{quantum mechanics!observable}

\index{quantum mechanics!measurement}

\index{Schrödinger, E. R.}

\index{von Neumann, J.}

In both cases, the observable\index{observable}s are represented
by families of selfadjoint operators in the respective Hilbert spaces.
For the purpose here, we pick the pair of selfadjoint\index{operators!selfadjoint}
operators representing momentum (denoted $P$) and position (denoted
$Q$). In one degree of freedom, we only need a single pair. The canonical
commutation relation is  
\[
PQ-QP=-i\,I;\;\mbox{or}\;PQ-QP=-i\,\hbar\,I
\]
where $\hbar=\frac{h}{2\pi}$ is Planck's constant, and $i=\sqrt{-1}$.
\index{Heisenberg, W.K.!commutation relation}\index{commutation relations}
\index{selfadjoint operator}

A few years after the pioneering work of Heisenberg and Schrödinger,
J. von Neumann and M. Stone proved that the two approaches are \emph{unitarily
equivalent}, hence they produce the same ``measurements.'' In modern
lingo, the notion of measurement take the form of projection valued
measures, which in turn are the key ingredient in the modern formulation
of the spectral theorem for selfadjoint, or normal, linear operators
in Hilbert space. (See \cite{Sto90, Yos95, Ne69, RS75, DS88b}.) Because
of dictates from physics, the ``interesting'' operators, such as
$P$ and $Q$ are \emph{unbounded}. \index{Stone, M. H.}\index{Spectral Theorem}\index{Theorem!Spectral-}\index{measurement}

The first point we will discuss about the pair of operators $P$ and
$Q$ is non-commutativity. As is typical in mathematical physics,
non-commuting operators will satisfy conditions on the resulting commutators.
In the case of $P$ and $Q$, the commutation relation is called the
canonical commutation relation; see below. For reference, see \cite{Dir47,Hei69,vN31,vN32b}.

Quantum mechanics was born during the years from 1900 to 1933. It
was created to explain phenomena in black body radiation, hydrogen
atom, where a discrete pattern occurs in the frequencies of waves
in the radiation. The radiation energy turns out to be $E=\nu\hbar$,
with $\hbar$ being the Plank's constant, and $\nu$ is frequency.
Classical mechanics runs into trouble. 

During the years of 1925 and 1926, Heisenberg found a way to represent
the energy $E$ as a matrix (spectrum = energy levels), so that the
matrix entries $\left\langle v_{j},Ev_{i}\right\rangle $ represent
transition probability for transitions from energy level $i$ to energy
level $j$. (See \figref{qmprob} below.) A fundamental relation in
quantum mechanics is the commutation relation satisfied by the momentum
operator $P$ and the position operator $Q$, where 
\begin{equation}
PQ-QP=-i\,I,\quad i=\sqrt{-1}.\label{eq:pq}
\end{equation}
Heisenberg represented the operators $P,Q$ by infinite matrices,
although his solution to (\ref{eq:pq}) is not really matrices, and
not finite matrices.

\index{quantum mechanics!momentum operator}

\index{quantum mechanics!position operator}

\index{operators!momentum-}

\index{operators!position-}
\begin{lem}
Eq. (\ref{eq:pq}) has \uline{no} solutions for finite matrices,
in fact, not even for bounded operators.\end{lem}
\begin{proof}
The reason is that for matrices, there is a trace operation where
\index{trace} 
\[
trace(AB)=trace(BA).
\]
This implies the trace on the left-hand-side is zero, while the trace
on the RHS is not. 
\end{proof}
This shows that there is no finite dimensional solution to the commutation
relation above, and one is forced to work with infinite dimensional
Hilbert space and operators on it. Notice also that $P,Q$ do not
commute, and the above commutation relation leads to the uncertainty\index{uncertainty}
principle (Hilbert, Max Born, von Neumann worked out the mathematics).
It states that the statistical variance $\triangle P$ and $\triangle Q$
satisfy $\triangle P\triangle Q\geq\hbar/2$ . We will come back to
this later in \exerref{uncertainty}.

\index{commutation relations}

We will show that non-commutativity always yields ``uncertainty.''\index{Heisenberg, W.K.!uncertainty principle}

However, Heisenberg \cite{vN31,Hei69} found his ``matrix'' solutions
by tri-diagonal $\infty\times\infty$ matrices, where \index{matrix!tri-diagonal-}\index{operators!adjoint-}
\begin{equation}
P=\frac{1}{\sqrt{2}}\left[\begin{array}{ccccc}
0 & 1\\
1 & 0 & \sqrt{2}\\
 & \sqrt{2} & 0 & \sqrt{3}\\
 &  & \sqrt{3} & 0 & \ddots\\
 &  &  & \ddots & \ddots
\end{array}\right]\label{eq:HP}
\end{equation}
 and 
\begin{equation}
Q=\frac{1}{i\sqrt{2}}\left[\begin{array}{ccccc}
0 & 1\\
-1 & 0 & \sqrt{2}\\
 & -\sqrt{2} & 0 & \sqrt{3}\\
 &  & -\sqrt{3} & 0 & \ddots\\
 &  &  & \ddots & \ddots
\end{array}\right]\label{eq:HQ}
\end{equation}
the complex $i$ in front of $Q$ is to make it selfadjoint.
\begin{xca}[The canonical commutation relation]
\myexercise{The canonical commutation relation}\label{exer:Hmatrix}Using
matrix multiplication for $\infty\times\infty$ matrices verify directly
that the two matrices $P$ and $Q$ satisfy $PQ-QP=-i\,I$, where
$I$ is the identity matrix in $l^{2}\left(\mathbb{N}_{0}\right)$,
i.e., $\left(I\right)_{ij}=\delta_{ij}$. \uline{Hint}: use the
rules in  \lemref{mat}. \index{matrix!inftytimesinfty
@$\infty\times\infty$}
\end{xca}

\begin{figure}
\includegraphics[scale=0.4]{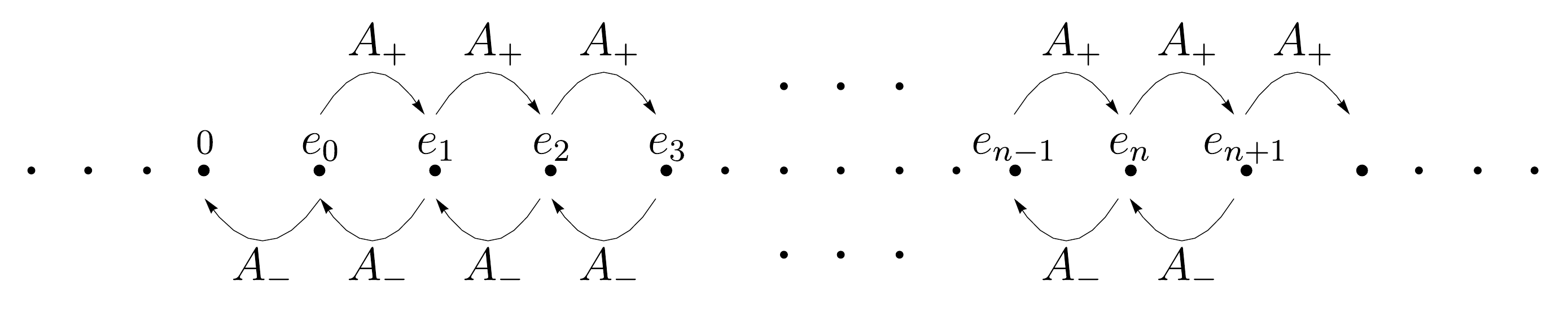}

\protect\caption{\label{fig:rl}The raising and lowering operators $A_{\pm}$. The
lowering operator $A_{-}$ kills $e_{0}$.}
\end{figure}

\begin{xca}[Raising and lowering operators (non-commutative complex variables)]
\label{exer:PQmatrix}\myexercise{Raising and lowering operators}
Set 
\begin{equation}
A_{\mp}:=P\pm iQ;\label{eq:rl}
\end{equation}
and show that the matrix representation for these operators is as
follows:
\[
A_{-}=\sqrt{2}\begin{bmatrix}0 & 1 & 0 & 0 & 0 & \cdots\\
\vdots & 0 & \sqrt{2} & 0 & 0 & \cdots\\
\vdots & \vdots & 0 & \sqrt{3} & 0 & \cdots\\
\vdots & \vdots & \vdots & 0 & \sqrt{4} & \cdots\\
\vdots & \vdots & \vdots & \vdots & \ddots & \ddots
\end{bmatrix}
\]
and 
\[
A_{+}=\sqrt{2}\begin{bmatrix}0 & 0 & \cdots & \cdots & \cdots & \cdots\\
1 & 0 & 0 & \cdots & \cdots & \cdots\\
0 & \sqrt{2} & 0 & 0 & \cdots & \cdots\\
\vdots & 0 & \sqrt{3} & 0 & 0 & \cdots\\
\vdots & \vdots & 0 & \sqrt{4} & \ddots & \ddots
\end{bmatrix}.
\]
In other words, the raising operator $A_{+}$ is a sub-banded matrix,
while the lowering operator $A_{-}$ is a supper-banded matrix. Both
$A_{+}$ and $A_{-}$ has $0$s down the diagonal. \index{matrix!raising/lowering-}
\index{matrix!banded-}

Further, show that 
\[
A_{-}A_{+}=2\begin{bmatrix}1 & 0 & \cdots & \cdots & \cdots & \cdots\\
0 & 2 & 0 & \cdots & \cdots & \cdots\\
\vdots & 0 & 3 & 0 & \cdots & \cdots\\
\vdots & \vdots & 0 & 4 & 0 & \cdots\\
\vdots & \vdots & \vdots & \vdots & \ddots & \ddots
\end{bmatrix}
\]
i.e., a diagonal matrix, with the numbers $\mathbb{N}$ down the diagonal
inside $\begin{bmatrix}\ddots\end{bmatrix}$; and 
\[
A_{+}A_{-}=2\begin{bmatrix}0 & 0 & \cdots & \cdots & \cdots & \cdots & \cdots\\
0 & 1 & 0 & \cdots & \cdots & \cdots & \cdots\\
\vdots & 0 & 2 & 0 & \cdots & \cdots & \cdots\\
\vdots &  & 0 & 3 & 0 & \cdots & \cdots\\
\vdots &  & \vdots & 0 & 4 & 0 & \cdots\\
\vdots &  & \vdots & \vdots & \vdots & \ddots & \ddots
\end{bmatrix};
\]
so that \index{matrix!diagonal-} 
\[
\frac{1}{2}\left[A_{-},A_{+}\right]=I.
\]
\end{xca}
\begin{rem}[Raising and lowering in an ONB]
In the canonical ONB $\left\{ e_{n}\right\} _{n=0}^{\infty}$ in
$l^{2}\left(\mathbb{N}_{0}\right)$, we have the following representations
of the two operators $A_{\pm}$ (see (\ref{eq:rl})):
\begin{eqnarray*}
A_{+}e_{n} & = & \sqrt{2}\sqrt{n+1}e_{n+1};\;n=0,1,2,\ldots,\;\mbox{and}\\
A_{-}e_{n} & = & \sqrt{2}\sqrt{n}e_{n-1};\:n=1,2,\ldots,\\
A_{-}e_{0} & = & 0,\;\mbox{see Fig}\:\ref{fig:rl}.
\end{eqnarray*}
The vector $e_{0}$ is called the \emph{ground state}, or \emph{vacuum
vector}.
\end{rem}

\begin{rem}
The conclusion of the discussion above is that the Heisenberg commutation
relation (\ref{eq:pq}) for pairs of selfadjoint operators has two
realizations, one in $L^{2}\left(\mathbb{R}\right)$, and the other
in $l^{2}\left(\mathbb{N}\right)$. 

In the first one we have 
\begin{eqnarray*}
\left(Pf\right)\left(x\right) & = & \frac{1}{i}\frac{d}{dx}f\\
\left(Qf\right)\left(x\right) & = & xf\left(x\right)
\end{eqnarray*}
for all $f\in\mathcal{S}\subset L^{2}\left(\mathbb{R}\right)$, where
$\mathcal{S}$ denotes the Schwartz test-function subspace in $L^{2}\left(\mathbb{R}\right)$.

The second realization is by $\infty\times\infty$ matrices, and it
is given in detail above. In \secref{Heisenberg} we shall return
to the first realization.

The \emph{Stone-von Neumann uniqueness theorem} (see \cite{vN32c,vN31})
implies the two solutions are unitarily equivalent; see  \chapref{groups}. 

\index{Stone-von Neumann Uniqueness Theorem}\index{Theorem!Stone-von Neumann uniqueness-}\end{rem}
\begin{xca}[Infinite banded matrices]
\myexercise{Infinite banded matrices}\label{exer:esa}Give an example
of two sequences 
\[
d_{1},d_{2},\ldots\in\mathbb{R},\;a_{1},a_{2},\ldots\in\mathbb{C}
\]
such that the corresponding Hermitian symmetric $\infty\times\infty$
tri-diagonal (banded) matrix $A$ in \figref{trid} satisfies $A\subset A^{*}$,
but $\overline{A}\neq A^{*}$, i.e., $A$ is \emph{not} essentially
selfadjoint when realized as a Hermitian operator in $l^{2}$.
\end{xca}
\begin{figure}
\[
A=\begin{bmatrix}d_{1} & a_{1} & 0 & \cdots & \cdots\\
\overline{a_{1}} & d_{2} & a_{2} & 0 & \cdots & \cdots\\
0 & \overline{a_{2}} & d_{3} & a_{3} & 0 & \cdots & \cdots\\
\vdots & 0 & \overline{a_{3}} & \ddots & \ddots & \ddots & \cdots & \cdots\\
\vdots & \vdots & \ddots & \ddots & \ddots & \ddots & \ddots & \cdots & \cdots\\
 & \vdots & \vdots & 0 & \overline{a_{n-2}} & d_{n-1} & a_{n-1} & 0 & \cdots\\
 &  & \vdots & \vdots & 0 & \overline{a_{n-1}} & d_{n} & a_{n} & \ddots\\
 &  &  & \vdots & \vdots & 0 & \overline{a_{n}} & d_{n+1} & \ddots\\
 &  &  &  & \vdots & \vdots & \ddots & \ddots & \ddots
\end{bmatrix}
\]

\protect\caption{\label{fig:trid}$A\subset A^{*}$ (a Hermitian Jacobi matrix).}
\end{figure}

\index{matrix!banded-}

\index{operators!essentially selfadjoint-}

\index{operators!symmetric-}
\begin{rem}[Matrices vs operators]
 Every bounded linear operator (and many unbounded operators too)
in separable Hilbert space (and, in particular, in $l^{2}$) can be
realized as a well-defined \emph{infinite \textquotedblleft square\textquotedblright{}
matrix}. In $l^{2}$ we pick the canonical ONB, but in a general Hilbert
space, a choice of ONB must be made. We saw that most rules for finite
matrices carry over the case of infinite matrices; sums, products,
and adjoints.

For instance, in order to find the matrix of the sum of two bounded
operators, just find the sum of the matrices of these operators. And
the matrix of the adjoint operator $A^{*}$ (of a bounded operator
$A$ in Hilbert space) is the adjoint matrix (conjugate transpose)
of the matrix of the operator $A$.\index{operators!adjoint-}

So while it is \textquotedblleft easy\textquotedblright{} to go from
bounded operators to infinite \textquotedblleft square\textquotedblright{}
matrices, the converse is much more subtle.\index{space!Hilbert-}\end{rem}
\begin{xca}[The Hilbert matrix]
\label{exer:HilbertMatrix}\myexercise{The Hilbert matrix}~
\begin{enumerate}
\item \label{enu:Hmt1}Show that the Hilbert matrix 
\[
H=\left(\frac{1}{1+j+k}\right)_{j,k\in\mathbb{N}}
\]
defines a bounded selfadjoint operator $T_{H}$ in $l^{2}\left(\mathbb{N}\right)$.
\index{selfadjoint operator}
\item Show that 
\[
\left\Vert T_{H}\right\Vert _{UN}=\sqrt{\pi}
\]
where $\left\Vert \cdot\right\Vert _{UN}$ denotes the uniform operator
norm 
\[
\left\Vert Tx\right\Vert =\sup\left\{ \left\Vert Tx\right\Vert ,\;x\in l^{2},\;\left\Vert x\right\Vert =1\right\} .
\]

\item Show that $T_{H}$ is positive definite.\\
\uline{Hint}: 
\[
\int_{0}^{1}x^{n}dx=\frac{1}{1+n},\;n\in\mathbb{N}.
\]
 \index{matrix!Hilbert-}\index{positive definite!-operator}
\end{enumerate}
\end{xca}
\begin{rem}
Note that the Hilbert matrix $H$ is not banded; in fact every entry
in $H$ is positive. Nonetheless, it follows from an application of
\corref{infmt} ( \secref{multi}) that the Hilbert matrix $H$ in
\exerref{HilbertMatrix} is equivalent to a banded matrix $J$; and
there is a choice of $J$ to be tri-diagonal; a Jacobi matrix; see
\figref{trid}. Since $H$ yields a bounded selfadjoint operator in
$l^{2}$, it follows from \corref{infmt} that $J$ is in fact a bounded
Jacobi-matrix with the same norm as $H$.
\end{rem}

\subsection{Probabilistic Interpretation of Parseval in Hilbert Space}

\begin{figure}
\includegraphics[width=0.5\textwidth]{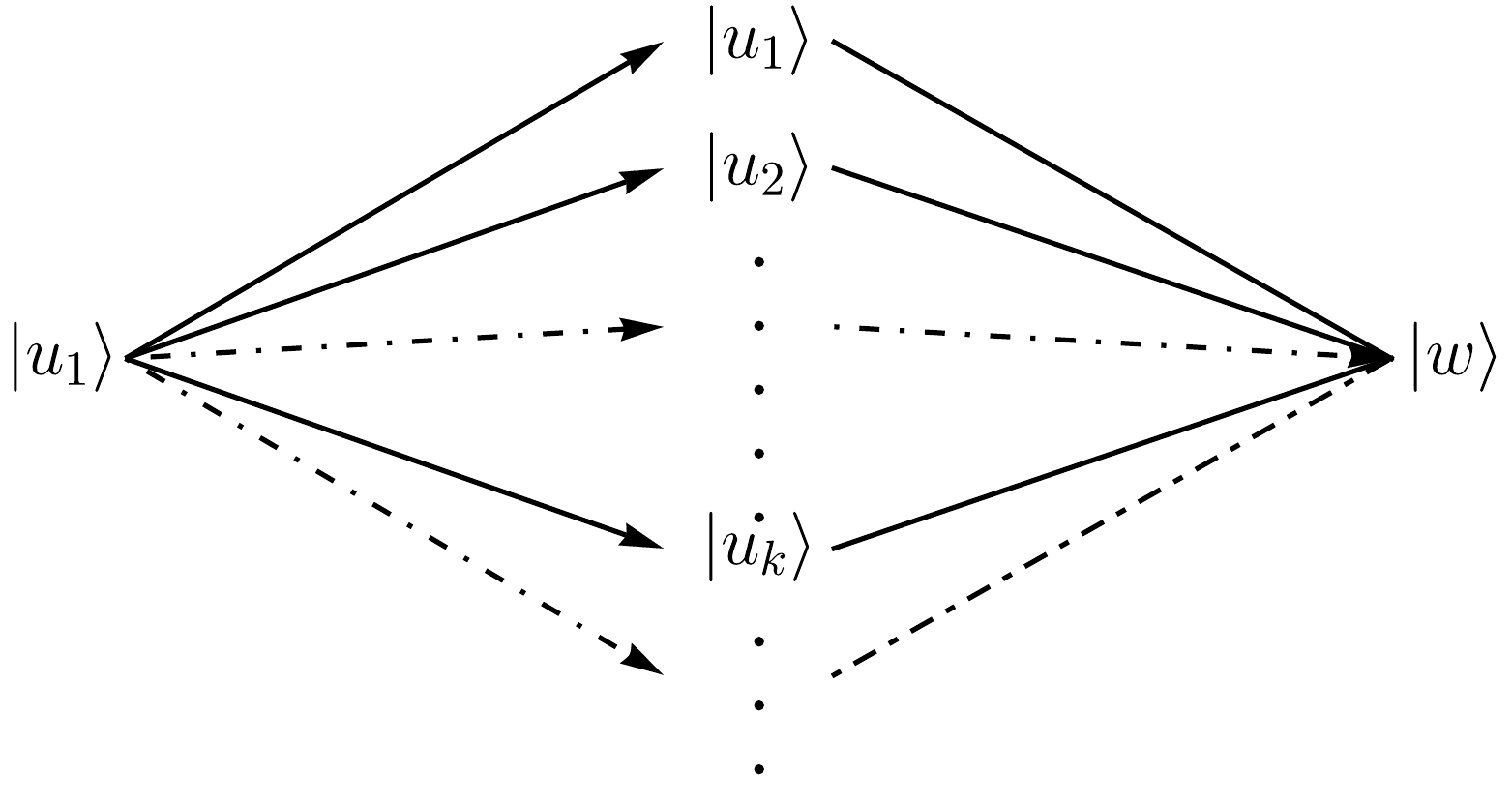}

\protect\caption{\label{fig:qmprob}Transition of quantum-states}
\end{figure}

\textbf{Case 1.} Let $\mathscr{H}$ be a complex Hilbert space, and
let $\left\{ u_{k}\right\} _{k\in\mathbb{N}}$ be an ONB, then Parseval's
formula reads:
\begin{equation}
\left\langle v,w\right\rangle _{\mathscr{H}}=\sum_{k\in\mathbb{N}}\left\langle v,u_{k}\right\rangle \left\langle u_{k},w\right\rangle ,\;\forall v,w\in\mathscr{H}.\label{eq:pr1}
\end{equation}
Translating this into a statement about ``transition probabilities''
for quantum states, $v,w\in\mathscr{H}$, with $\left\Vert v\right\Vert _{\mathscr{H}}=\left\Vert w\right\Vert _{\mathscr{H}}=1$,
we get 
\begin{equation}
\mbox{Prob}\left(v\rightarrow w\right)=\sum_{k\in\mathbb{N}}\mbox{Prob}\left(v\rightarrow u_{k}\right)\mbox{Prob}\left(u_{k}\rightarrow w\right).\label{eq:pr2}
\end{equation}
See \figref{qmprob}. The states $v$ and $w$ are said to be \emph{uncorrelated}
iff (Def.) they are orthogonal. \index{transition probabilities}

Fix a state $w\in\mathscr{H}$, then 
\[
\left\Vert w\right\Vert ^{2}=\sum_{k\in\mathbb{N}}\left|\left\langle u_{k},w\right\rangle \right|^{2}=1.
\]
The numbers $\left|\left\langle u_{i},w\right\rangle \right|^{2}$
represent a probability distribution over the index set, where $\left|\left\langle u_{k},w\right\rangle \right|^{2}$
is the probability that the quantum system is in the state $|u_{k}\rangle$.

\index{distribution!probability-}

\uline{Disclaimer:} The notation ``transition-probability'' in
(\ref{eq:pr2}) and  \figref{qmprob} is a stretch since the inner
products $\left\langle v,u_{k}\right\rangle $ are not positive. Nonetheless,
it is justified by 
\[
\sum_{k}\left|\left\langle v,u_{k}\right\rangle \right|^{2}=1
\]
when $v\in\mathscr{H}$ (is a state vector).

\textbf{Case 2. }If $P:\mathcal{B\left(\mathbb{R}\right)\rightarrow\mbox{Proj}\left(\mathscr{H}\right)}$
is a projection valued measure\index{measure!projection-valued} (Appendix
\ref{sec:stone}), we get the analogous assertions, but with integration,
as opposed to summation. In this case eq. (\ref{eq:pr1}) holds the
following form:
\begin{equation}
\left\langle v,w\right\rangle _{\mathscr{H}}=\int_{\mathbb{R}}\left\langle v,P\left(d\lambda\right)w\right\rangle _{\mathscr{H}};\label{eq:pr3}
\end{equation}
and for $v=w$, it reads:
\begin{equation}
\left\Vert v\right\Vert _{\mathscr{H}}^{2}=\int_{\mathbb{R}}\left\Vert P\left(d\lambda\right)v\right\Vert _{\mathscr{H}}^{2}.\label{eq:pr4}
\end{equation}
Recall the other axioms of $P\left(\cdot\right)$ are:\index{projection}
\begin{enumerate}
\item $P\left(A\right)=P\left(A\right)^{*}=P\left(A\right)^{2}$, $\forall A\in\mathcal{B}\left(\mathbb{R}\right)$.
\item $P\left(\cdot\right)$ is countably additive on the Borel subsets
of $\mathbb{R}$, i.e., 
\[
\sum_{j}P\left(A_{j}\right)=P\left(\cup_{j}A_{j}\right)
\]
where $A_{j}\in\mathcal{B}\left(\mathbb{R}\right)$, $A_{i}\cap A_{j}=\emptyset$,
$i\neq j$. 
\item $P\left(A\cap B\right)=P\left(A\right)P\left(B\right)$, $\forall A,B\in\mathcal{B}\left(\mathbb{R}\right)$.
\end{enumerate}

\section{\label{sec:proj}The Lattice Structure of Projections}

A lattice is a partially ordered set in which every two elements have
a supremum (also called a least upper bound or join) and an infimum
(also called a greatest lower bound or meet). \index{partially ordered set}

The purpose of the discussion below is twofold; one to identify two
cases: (i) the (easy) lattice of subsets of a fixed total set; and
(ii) the lattice of projections in a fixed Hilbert space. Secondly
we point out how non-commutativity of projections makes the comparison
of (i) and (ii) subtle; even though there are some intriguing correspondences;
see \tabref{lattice} for illustration.

\uline{Notation}: In this section we will denote projections $P$,
$Q$, etc. 

von Neumann invented the notion of abstract Hilbert space in 1928
as shown in one of the earliest papers.\footnote{Earlier authors, Schmidt and Hilbert, worked with infinite bases,
and $\infty\times\infty$ matrices.} His work was greatly motivated by quantum mechanics. In order to
express quantum mechanics logic operations, he created lattices of
projections, so that everything we do in set theory with set operation
has a counterpart in the operations of projections. See \tabref{lattice}.\index{lattice}

\renewcommand{\arraystretch}{2}

\begin{table}
\begin{longtable}{|c|c|c|c|}
\hline 
SETS & CHAR & PROJ & DEF\tabularnewline
\hline 
\hline 
$A\cap B$ & $\chi_{A}\chi_{B}$ & $P\wedge Q$ & $P\mathscr{H}\cap Q\mathscr{H}$\tabularnewline
\hline 
$A\cup B$ & $\chi_{A\cup B}$ & $P\vee Q$ & $\overline{span}\{P\mathscr{H}\cup Q\mathscr{H}\}$\tabularnewline
\hline 
$A\subset B$ & $\chi_{A}\chi_{B}=\chi_{A}$ & $P\leq Q$ & $P\mathscr{H}\subset Q\mathscr{H}$ \tabularnewline
\hline 
$A_{1}\subset A_{2}\subset\cdots$ & $\chi_{A_{i}}\chi_{A_{i+1}}=\chi_{A_{i}}$ & $P_{1}\leq P_{2}\leq\cdots$ & $P_{i}\mathscr{H}\subset P_{i+1}\mathscr{H}$\tabularnewline
\hline 
$\bigcup_{k=1}^{\infty}A_{k}$ & $\chi_{\cup_{k}A}$ & $\vee_{k=1}^{\infty}P_{k}$ & $\overline{span}\{\bigcup_{k=1}^{\infty}P_{k}\mathscr{H}\}$\tabularnewline
\hline 
$\bigcap_{k=1}^{\infty}A_{k}$ & $\chi_{\cap_{k}A_{k}}$ & $\wedge_{k=1}^{\infty}P_{k}$ & $\bigcap_{k=1}^{\infty}P_{k}\mathscr{H}$\tabularnewline
\hline 
$A\times B$ & $\left(\chi_{A\times X}\right)\left(\chi_{X\times B}\right)$ & $P\otimes Q$ & $P\otimes Q\in proj(\mathscr{H}\otimes\mathscr{K})$\tabularnewline
\hline 
\end{longtable}

\protect\caption{\label{tab:lattice}Lattice of projections in Hilbert space.}
\end{table}

\renewcommand{\arraystretch}{1}

\index{lattice}

For example, if $P$ and $Q$ are two projections in $\mathscr{B}\left(\mathscr{H}\right)$,
then 
\begin{eqnarray}
P\mathscr{H} & \subset & Q\mathscr{H}\label{eq:p1}\\
 & \Updownarrow\nonumber \\
P & = & PQ\label{eq:p2}\\
 & \Updownarrow\nonumber \\
P & \leq & Q.\label{eq:p2-1}
\end{eqnarray}
This is similar to the following equivalence relation in set theory
\begin{eqnarray}
A & \subset & B\;(\text{containment of sets})\label{eq:p3}\\
 & \Updownarrow\nonumber \\
A & = & A\cap B.\label{eq:p4}
\end{eqnarray}
In general, product and sum of projections are not projections. But
if $P\mathscr{H}\subset Q\mathscr{H}$ then the product $PQ$ is in
fact a projection. Taking adjoint in (\ref{eq:p2}) yields 
\[
P^{*}=(PQ)^{*}=Q^{*}P^{*}=QP.
\]
It follows that $PQ=QP=P$, i.e., \textit{containment of subspaces
implies the corresponding projections commute}. 

Two decades before von Neumann developed his Hilbert space theory,
Lebesgue developed his integration theory \cite{MR1504529} which
extends the classical Riemann integral. The monotone sequence of sets
$A_{1}\subset A_{2}\subset\cdots$ in Lebesgue's integration theory
also has a counterpart in the theory of Hilbert space. 
\begin{lem}
Let $P_{1}$ and $P_{2}$ be orthogonal projections acting on $\mathscr{H}$,
then 
\begin{align}
P_{1}\leq P_{2}\Longleftrightarrow\left\Vert P_{1}x\right\Vert  & \leq\left\Vert P_{2}x\right\Vert ,\;\forall x\in\mathscr{H}\label{eq:pp1}\\
 & \Updownarrow\nonumber \\
P_{1} & =P_{1}P_{2}=P_{2}P_{1}\label{eq:pp2}
\end{align}
(see  \tabref{lattice}.)
\end{lem}
\Needspace*{3\baselineskip}
\begin{proof}
Indeed, for all $x\in\mathscr{H}$, we have 
\[
\left\Vert P_{1}x\right\Vert ^{2}=\left\langle P_{1}x,P_{1}x\right\rangle =\left\langle x,P_{1}x\right\rangle =\left\langle x,P_{2}P_{1}x\right\rangle \leq\left\Vert P_{1}P_{2}x\right\Vert ^{2}\leq\left\Vert P_{2}x\right\Vert ^{2}.
\]
\end{proof}
\begin{thm}
For every monotonically increasing sequence of projections 
\[
P_{1}\leq P_{2}\leq\cdots,
\]
and setting 
\[
P:=\vee P_{k}=\lim_{k}P_{k},
\]
then $P$ defines a projection, the limit.\end{thm}
\begin{proof}
The assumption $P_{1}\leq P_{2}\leq\cdots$ implies that $\left\{ \left\Vert P_{k}x\right\Vert \right\} _{k=1}^{\infty}$,
$x\in\mathscr{H}$, is a monotone increasing sequence in $\mathbb{R}$,
and the sequence is bounded by $\left\Vert x\right\Vert $, since
$\left\Vert P_{k}x\right\Vert \leq\left\Vert x\right\Vert $, for
all $k\in\mathbb{N}$. Therefore the sequence $\left\{ P_{k}\right\} _{k=1}^{\infty}$
converges to $P\in\mathscr{B}\left(\mathscr{H}\right)$ (strongly),
and $P$ really defines a selfadjoint projection. (We use ``$\leq$''
to denote the lattice operation on projection.) Note the convergence
refers to the strong operator topology, i.e., for all $x\in\mathscr{H}$,
there exists a vector, which we denote by $Px$, so that $\lim_{k}\left\Vert P_{k}x-Px\right\Vert =0$. 
\end{proof}
The examples in  \secref{Hilbert} using Gram-Schmidt process can
now be formulated in the lattice of projections. 

\begin{figure}
\includegraphics[width=0.6\textwidth]{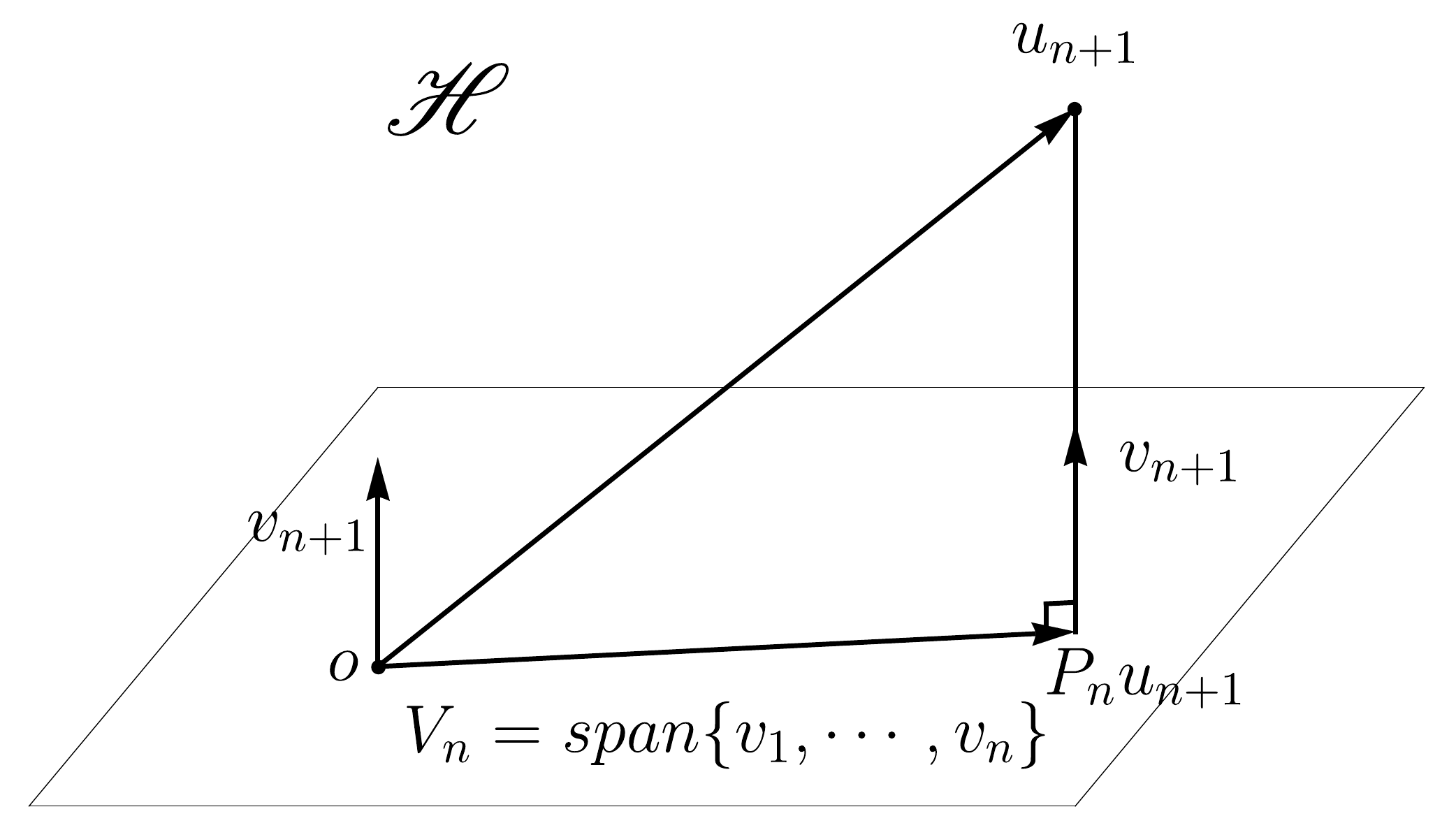}

\protect\caption{\label{fig:gmi}Gram-Schmidt: $V_{n}\protect\longrightarrow V_{n+1}$}
\end{figure}

Recall (\lemref{Gram-Schmidt}) that for a linearly independent subset
$\left\{ u_{k}\right\} \subset\mathscr{H}$, the Gram-Schmidt process
yields an orthonormal set $\left\{ v_{k}\right\} \subset\mathscr{H}$,
with $v_{1}:=u_{1}/\left\Vert u_{1}\right\Vert $, and 
\[
v_{n+1}:=\frac{u_{n+1}-P_{n}u_{n+1}}{\left\Vert u_{n+1}-P_{n}u_{n+1}\right\Vert },\;n=1,2,\ldots;
\]
where $P_{n}$ is the orthogonal projection on the $n$-dimensional
subspace 
\[
V_{n}:=span\left\{ v_{1},\ldots,v_{n}\right\} .
\]
See \figref{gmi}.

Note that 
\begin{eqnarray*}
V_{n}\subset V_{n+1}\rightarrow\bigcup_{n}V_{n} & \sim & P_{n}\leq P_{n+1}\rightarrow P\\
 &  & P_{n}^{\perp}\geq P_{n+1}^{\perp}\rightarrow P^{\perp}.
\end{eqnarray*}
Assume $\bigcup_{n}V_{n}$ is dense in $\mathscr{H}$, then $P=I$
and $P^{\perp}=0$. In lattice notations, we may write 
\begin{eqnarray*}
\vee P_{n}=\sup P_{n} & = & I\\
\wedge P_{n}^{\perp}=\inf P_{n}^{\perp} & = & 0.
\end{eqnarray*}

\begin{lem}
\label{lem:proj1}Let $P,Q\in Proj\left(\mathscr{H}\right)$, then
\[
P+Q\in Proj(\mathscr{H})\Longleftrightarrow PQ=QP=0,\;\mbox{i.e., }P\perp Q.
\]
\end{lem}
\begin{proof}
Notice that 
\begin{equation}
(P+Q)^{2}=P+Q+PQ+QP\label{eq:p4-1}
\end{equation}
and so 
\begin{eqnarray}
\left(P+Q\right)^{2} & = & P+Q\label{eq:p5}\\
 & \Updownarrow\nonumber \\
PQ+QP & = & 0.\label{eq:p6}
\end{eqnarray}

Suppose $PQ=QP=0$ then 
\[
(P+Q)^{2}=P+Q=(P+Q)^{*},\;\mbox{i.e.,}\;P+Q\in Proj\left(\mathscr{H}\right).
\]

Conversely, if $P+Q\in Proj(\mathscr{H})$, then $(P+Q)^{2}=P+Q\Longrightarrow PQ+QP=0$
by (\ref{eq:p6}). Also, $(PQ)^{*}=Q^{*}P^{*}=QP$, combining with
(\ref{eq:p6}) yields 
\begin{equation}
\left(PQ\right)^{*}=QP=-PQ.\label{eq:p7}
\end{equation}
Then, 
\[
\left(PQ\right)^{2}=P\left(QP\right)Q\underset{\left(\ref{eq:p7}\right)}{=}-P(PQ)Q=-PQ
\]
which implies $PQ\left(I+PQ\right)=0$. Hence, 
\begin{equation}
PQ=0\quad\mbox{or}\quad PQ=I.\label{eq:p8}
\end{equation}
But by (\ref{eq:p7}), $PQ$ is skew-adjoint, it follows that $PQ=0$,
and so $QP=0$.\end{proof}
\begin{rem}
Eq. (\ref{eq:p4-1}) is analogous to the following identity for characteristic
functions: 
\[
\chi_{A}+\chi_{B}=\chi_{A\cup B}-\chi_{A\cap B}
\]
Therefore, $\left(\chi_{A}+\chi_{B}\right)^{2}=\chi_{A}+\chi_{B}$
iff $\chi_{A\cap B}=0$, i.e., iff $A\cap B=\emptyset$. 
\end{rem}

The set of projections in a Hilbert space $\mathscr{H}$ is partially
ordered according to the corresponding closed subspaces partially
ordered by inclusion. Since containment implies commuting, the chain
of projections \index{partially ordered set} 
\[
P_{1}\leq P_{2}\leq\cdots
\]
is a family of commuting selfadjoint operators. By the spectral theorem
(\chapref{sp}), $\{P_{i}\}$ may be simultaneously diagonalized,
so that $P_{i}$ is unitarily equivalent to the operator of multiplication
by $\chi_{E_{i}}$ on the Hilbert space $L^{2}(X,\mu)$, where $X$
is compact and Hausdorff. Therefore the lattice structure of projections
in $\mathscr{H}$ is precisely the lattice structure of $\chi_{E}$,
or equivalently, the lattice structure of measurable sets in $X$.\index{operators!adjoint-}
\index{selfadjoint operator}
\begin{lem}
Consider $L^{2}(X,\mu)$. The following are equivalent.
\begin{enumerate}
\item $E\subset F$;
\item $\chi_{E}\chi_{F}=\chi_{F}\chi_{E}=\chi_{E}$;
\item $\left\Vert \chi_{E}f\right\Vert \leq\left\Vert \chi_{F}f\right\Vert $,
for any $f\in L^{2}$;
\item $\chi_{E}\leq\chi_{F}$, in the sense that 
\[
\left\langle f,\chi_{E}f\right\rangle \leq\left\langle f,\chi_{F}f\right\rangle ,\;\forall f\in L^{2}\left(X\right).
\]

\end{enumerate}
\end{lem}
\begin{proof}
The proof is trivial. Note that
\begin{eqnarray*}
\left\langle f,\chi_{E}f\right\rangle  & = & \int\chi_{E}\left|f\right|^{2}d\mu\\
\left\Vert \chi_{E}f\right\Vert ^{2} & = & \int\left|\chi_{E}f\right|^{2}d\mu=\int\chi_{E}\left|f\right|^{2}d\mu
\end{eqnarray*}
where we used that fact that
\[
\chi_{E}=\overline{\chi_{E}}=\chi_{E}^{2}.
\]

\end{proof}

What makes $Proj\left(\mathscr{H}\right)$ intriguing is the non-commutativity.
For example, if $P,Q\in Proj\left(\mathscr{H}\right)$ are given,
it does not follow (in general) that $P+Q\in Proj\left(\mathscr{H}\right)$;
nor that $PQP\in Proj\left(\mathscr{H}\right)$. These two conclusions
only hold if it is further assumed that $P$ and $Q$ commute; see
Lemmas \ref{lem:proj1} and \ref{lem:proj2}.
\begin{lem}
\label{lem:proj2}Let $P,Q\in Proj\left(\mathscr{H}\right)$; then
the following conditions are equivalent:
\begin{enumerate}
\item \label{enu:proj1}$PQP\in Proj\left(\mathscr{H}\right)$;
\item \label{enu:proj2}$PQ=QP$. 
\end{enumerate}
\end{lem}
\begin{proof}
First note that the operator $A=PQ-QP$ is skew-symmetric, i.e., $A^{*}=-A$,
and so its spectrum is contained in the imaginary line $i\mathbb{R}$. 

The implication (\ref{enu:proj2})$\Rightarrow$(\ref{enu:proj1})
above is immediate so assume (\ref{enu:proj1}), i.e., that 
\[
\left(PQP\right)^{2}=PQP.
\]
And using this, one checks by a direct computation that $A^{3}=0$.
But with $A^{*}=-A$, and the spectral theorem, we therefore conclude
that $A=0$, in other words, (\ref{enu:proj2}) holds.
\end{proof}

\index{Spectral Theorem}

\index{lattice}

\section{\label{sec:multi}Multiplication Operators}
\begin{xca}[Multiplication operators]
\myexercise{Multiplication operators}\label{exer:mult}Let $\left(X,\mathcal{F},\mu\right)$
be a measure space, assume $\mu$ is $\sigma$-finite. Let $L^{2}\left(\mu\right)$
be the corresponding Hilbert space. Let $\varphi$ be a locally integrable
function on $X$, and set
\[
M_{\varphi}f:=\varphi f,
\]
pointwise product, defined for 
\[
f\in dom\left(M_{\varphi}\right)=\left\{ f\in L^{2}\left(\mu\right)\::\:\varphi f\in L^{2}\left(\mu\right)\right\} .
\]
$M_{\varphi}$ is called a \emph{multiplication operator}.
\begin{enumerate}
\item Show that $M_{\varphi}$ is normal.
\item Show that $M_{\varphi}$ is selfadjoint iff $\varphi$ is $\mu$-a.e.
real-valued.
\item Show that $M_{\varphi}$ is bounded in $L^{2}\left(\mu\right)$ iff
$\varphi\in L^{\infty}\left(\mu\right)$; and, in this case, 
\begin{equation}
\left\Vert M_{\varphi}\right\Vert _{UN}=\left\Vert \varphi\right\Vert _{L^{\infty}\left(\mu\right)}.\label{eq:mul}
\end{equation}

\item Show that if $\varphi\in L^{\infty}\left(\mu\right)$, then $dom\left(M_{\varphi}\right)=L^{2}\left(\mu\right)$. 
\item Discuss the converse.
\end{enumerate}
\end{xca}

\begin{xca}[Moment theory \cite{MR0184042}]
\label{exer:moment}\myexercise{Moment theory}Let $\mu$ be a positive
Borel measure on $\mathbb{R}$ such that 
\begin{equation}
\int_{\mathbb{R}}x^{2n}d\mu\left(x\right)<\infty\label{eq:m2-1}
\end{equation}
for all $n\in\mathbb{N}$, i.e., $\mu$ has finite moments of all
orders. Let $\varphi\left(x\right)=x$, and $M=M_{\varphi}$ the corresponding
multiplication operator in $L^{2}\left(\mu\right)=$ $L^{2}\left(\mathbb{R},\mathcal{B},\mu\right)$,
i.e., 
\begin{equation}
\left(Mf\right)\left(x\right)=\left(Qf\right)\left(x\right)=xf\left(x\right),\label{eq:m2-2}
\end{equation}
for all $f\in L^{2}\left(\mu\right)$ such that $xf\in L^{2}\left(\mu\right)$.\index{matrix!banded-}
\index{moment problem}
\begin{enumerate}
\item Using Gram-Schmidt (\lemref{Gram-Schmidt}), show that $M$ has a
matrix representation by an $\infty\times\infty$ tri-diagonal (banded)
matrix as in  \figref{trid}.\vspace{1em}\\
Akhiezer calls these infinite banded matrices \emph{Jacobi matrices}.
They define formally selfadjoint (alias symmetric) operators in $l^{2}$;
unbounded of course. And these operators can only attain von Neumann
indices $(0,0)$ or $(1,1)$. Both are possible. \index{matrix!Jacobi}\index{operators!symmetric-}
\item \label{enu:m2-2}Work out a recursive formula for the two sequences
$\left(a_{n}\right)_{n\in\mathbb{N}}$ and $\left(d_{n}\right)_{n\in\mathbb{N}}$
in the expression for $M$ by the matrix of  \figref{trid} in terms
of the moments $\left(s_{n}\right)_{n\in\mathbb{N}\cup\left\{ 0\right\} }$:
\[
s_{n}:=\int_{\mathbb{R}}x^{n}d\mu\left(x\right).
\]

\item \label{enu:m2-3}Same question as in (\ref{enu:m2-2}) but for the
special case when $\mu$ is 
\[
d\mu\left(x\right)=\frac{1}{\sqrt{2\pi}}e^{-\frac{x^{2}}{2}}dx
\]
i.e., the $N\left(0,1\right)$ Gaussian measure on $\mathbb{R}$.\\
\uline{Hint}: Show first that the moments $\left\{ s_{k}\right\} _{k\in\left\{ 0\right\} \cup\mathbb{N}}$
of $\mu_{N\left(0,1\right)}$ are as follows (the Gaussian moments):
\begin{eqnarray*}
s_{2n+1} & = & 0,\;\mbox{and}\\
s_{2n} & = & \frac{\left(2n\right)!}{2^{n}\cdot n!}=\left(2n-1\right)!!\left(=\left(2n-1\right)\left(2n-3\right)\cdots5\cdot3\right).
\end{eqnarray*}

\item Give necessary and sufficient conditions for \emph{essential selfadjointness}
of the associated Jacobi matrix as an operator in $L^{2}\left(\mu\right)$,
expressed in terms $\mu$ and of the moments $\left(s_{n}\right)$
in (\ref{enu:m2-3}).
\end{enumerate}
\end{xca}
\begin{cor}
\label{cor:infmt}Let $A$ be a selfadjoint operator in a separable
Hilbert space, and suppose there is a cyclic vector $u_{0}$, $\Vert u_{0}\Vert=1$,
such that 
\begin{equation}
u_{0}\in\bigcap_{k\in\mathbb{N}}dom\left(A^{k}\right).\label{eq:cb1}
\end{equation}
Then there is an ONB $\left\{ e_{i}\right\} _{i\in\mathbb{N}}$ in
$\mathscr{H}$, contained in $dom\left(A\right)$ such that the corresponding
$\infty\times\infty$ matrix 
\begin{equation}
\left(M_{A}\right)_{i,j}=\left\langle e_{i},Ae_{j}\right\rangle ,\;i,j\in\mathbb{N}\label{eq:cb2}
\end{equation}
is banded; what is more, it is a Jacobi-matrix, see  \exerref{moment}.
\end{cor}
\index{cyclic!-vector} \index{selfadjoint operator}
\begin{proof}
Using the Spectral Theorem, we conclude that there is a measure $\mu_{0}$
on $\mathbb{R}$ and a unitary transform $W:L^{2}\left(\mathbb{R},\mu_{0}\right)\longrightarrow\mathscr{H}$
such that\end{proof}
\begin{enumerate}
\item $Wf=f\left(A\right)u_{0},\;\forall f\in L^{2}\left(\mu_{0}\right);$

\begin{enumerate}
\item $WM_{t}=AW$, where $M_{t}$ denotes multiplication by $t$ in $L^{2}\left(\mu_{0}\right)$;
and
\item $\left\langle u_{0},f\left(A\right)u_{0}\right\rangle =\int_{\mathbb{R}}f\left(t\right)d\mu_{0}\left(t\right),\;\forall f\in L^{2}\left(\mu_{0}\right)$. 
\end{enumerate}

We refer the reader to \chapref{sp} for more details.

\end{enumerate}
\begin{proof}
Now apply Gram-Schmidt to the monomials $\{t^{k}\}_{k\in\left\{ 0\right\} \cup\mathbb{N}}$,
to get orthogonal polynomials $\left\{ p_{k}\left(t\right)\right\} $
such that 
\[
span_{k\leq n}\left\{ p_{k}\left(t\right)\right\} =span_{k\leq n}\{t^{k}\}
\]
holds for all $n\in\mathbb{N}$.

Set 
\begin{equation}
e_{k}:=p_{k}\left(A\right)u_{0},\;k\in\left\{ 0\right\} \cup\mathbb{N},\label{eq:cb3}
\end{equation}
$e_{0}=u_{0}$, and this is then the desired ONB. To see this, use
the conclusion from  \exerref{moment}, together with the following:
\begin{eqnarray*}
\left\langle e_{j},e_{k}\right\rangle  & = & \left\langle p_{j}\left(A\right)u_{0},p_{k}\left(A\right)u_{0}\right\rangle \\
 & = & \left\langle u_{0},p_{j}\left(A\right)p_{k}\left(A\right)u_{0}\right\rangle \\
 & = & \left\langle u_{0},\left(p_{j}p_{k}\right)\left(A\right)u_{0}\right\rangle \\
 & \underset{\text{by (iii)}}{=} & \int_{\mathbb{R}}p_{j}\left(t\right)p_{k}\left(t\right)d\mu_{0}\left(t\right)\\
 & = & \delta_{j,k}\;\left(\mbox{by Gram-Schmidt}.\right)
\end{eqnarray*}

\end{proof}
\begin{flushleft}
Historical Note. 
\par\end{flushleft}

In \cite{MR1067743}, Lax relates an account of von Neumann and F.
Rellich speaking in Hilbert's seminar, in Göttingen (around 1930).
When they came to ``selfadjoint operator in Hilbert space,'' Erhard
Schmidt (of Gram-Schmidt) would interrupt: ``Please, young man, say
\emph{infinite matrix}.''

Ironically, von Neumann invented numerical methods for ``large''
matrices toward the end of his career \cite{MR0041539}. 

\sindex[nam]{Lax, P.D., (1926 --)}

\section*{A summary of relevant numbers from the Reference List}

\begin{flushleft}
For readers wishing to follow up sources, or to go in more depth with
topics above, we suggest: \cite{Tay86,Arv72,MR1357166,BR79,Con90,DM85,DS88b,MR1892228,RS75,MR1068530,Rud73,Rud87,MR3231624,MR3024468,MR3050315,MR2966130,AAR13,BoMc13,Hi80,Ito06,Jor14,MR3210626,BJ02,MR2254502,CW14,MR0184066,MR950981,MR0377445,JoMu80,Hel13,MR2908248,MR2382238}.
\par\end{flushleft}

\begin{subappendices}

\section{Hahn-Banach Theorems}

\emph{Version 1.} Let $S$ be a subspace of a real vector space $X$.
Let $l:S\rightarrow\mathbb{R}$ be a linear functional, and let $p:X\rightarrow\mathbb{R}$
satisfy
\begin{eqnarray}
p\left(x+y\right) & \leq & p\left(x\right)+p\left(y\right),\quad x,y\in X;\label{eq:hb1-1}\\
p\left(tx\right) & = & t\,p\left(x\right),\quad t\in\mathbb{R}_{+},\:x\in X.\label{eq:hb1-2}
\end{eqnarray}

\begin{thm}[HB1]
 Let $X,S,p$, and $l$ be as above, and assume:
\begin{equation}
l\left(x\right)\leq p\left(x\right),\quad x\in S.\label{eq:hb2}
\end{equation}
Then there is a linear functional $\widetilde{l}:X\rightarrow\mathbb{R}$,
extending $l$ on the subspace, and satisfying
\begin{equation}
\widetilde{l}\left(x\right)\leq p\left(x\right),\quad\forall x\in X.\label{eq:hb3}
\end{equation}

\end{thm}
\uline{Hint}: Introduce a partially ordered set (p.o.s.) $\left(T,m\right)$
where $S\subset T\subset X$, $T$ is a subspace, $m:T\rightarrow X$
is a linear functional extending $\left(S,l\right)$ and satisfying
\begin{equation}
m\left(x\right)\leq p\left(x\right),\quad\forall x\in T.\label{eq:hb4}
\end{equation}
Define the order $\left(T,m\right)\leq\left(T',m'\right)$ to mean
that $T\subseteq T'$ and $m'$ agrees with $m$ on $T$. (Both satisfying
(\ref{eq:hb4}).) Apply Zorn's lemma to this p.o.s., and show that
every maximal element must be a solution to (\ref{eq:hb3}).

\section{Banach-Limit}

Consider $X:=l_{\mathbb{R}}^{\infty}\left(\mathbb{N}\right)$ = all
bounded real sequences $x=\left(x_{1},x_{2},\cdots\right)$, and the
shift $\sigma:l^{\infty}\rightarrow l^{\infty}$, defined by
\[
\sigma\left(x_{1},x_{2},x_{3},\cdots\right)=\left(x_{2},x_{3},\cdots\right).
\]
Consider the subspace $S\subset X$ consisting of all convergent sequences,
and for $x\in S$, set $l\left(x\right)=\lim_{k\rightarrow\infty}x_{k}$,
i.e., it holds that, for $\forall\varepsilon\in\mathbb{R}_{+}$, $\exists\,n$
such that 
\[
\left|x_{k}-l\left(x\right)\right|<\varepsilon\quad\mbox{for}\quad\forall k\geq n.
\]

\begin{thm}[Banach]
There is a linear functional, called $LIM:X\rightarrow\mathbb{R}$,
(a Banach limit) having the following properties:
\begin{enumerate}
\item[(i)] $LIM$ is an extension of $l$ on $S$.
\item[(ii)] 
\[
\liminf_{k}x_{k}\leq LIM\left(x\right)\leq\limsup_{k}x_{k}\quad\mbox{and}
\]

\item[(iii)]  
\[
LIM\circ\sigma=LIM,
\]
i.e., $LIM$ is shift-invariant.
\end{enumerate}
\end{thm}
\begin{proof}
This is an application of (HB1) but with a modification; we set 
\[
q\left(x\right)=\limsup x_{k},\quad\mbox{and}
\]
\[
p\left(x\right)=\inf_{n}\frac{1}{n+1}\sum_{k=0}^{n}q\left(\sigma^{k}\left(x\right)\right),\quad\forall x\in X,
\]
and we note that $l\left(x\right)\leq p\left(x\right)$, $\forall x\in S$.
The rest of the proof follows that of HB1 \emph{mutatis mutandis}.\end{proof}
\begin{rem}
Note that $LIM$ is \emph{not} unique. 

It follows by the theorem above that $LIM$ is in $\left(l^{\infty}\left(\mathbb{N}\right)\right)^{*}$,
and that it is \emph{not} represented by any $y\in l^{1}\left(\mathbb{N}\right)$,
see \tabref{dual}. As a result, we have
\begin{equation}
\left(l^{\infty}\left(\mathbb{N}\right)\right)^{*}\supsetneq l^{1}\left(\mathbb{N}\right);\label{eq:bl1}
\end{equation}
i.e., the dual $\left(l^{\infty}\right)^{*}$ is (much) bigger than
$l^{1}$.\end{rem}
\begin{thm}[HB2]
 Let $X$ be a normed space, $S\subset X$ a closed subspace, $l:S\rightarrow\mathbb{R}$
a linear functional such that
\begin{equation}
\left|l\left(x\right)\right|\leq\left\Vert x\right\Vert ,\quad\forall x\in S.\label{eq:hb5}
\end{equation}
Then there is a $\widetilde{l}\in X^{*}$ such that 
\[
|\widetilde{l}\left(x\right)|\leq\left\Vert x\right\Vert ,\quad\forall x\in X,
\]
$\widetilde{l}$ extending $l$ from $S$, and 
\begin{equation}
\Vert\widetilde{l}\Vert_{X^{*}}=\left\Vert l\right\Vert _{S^{*}}.\label{eq:hb7}
\end{equation}

\end{thm}

\begin{figure}[H]
\includegraphics[width=0.5\columnwidth]{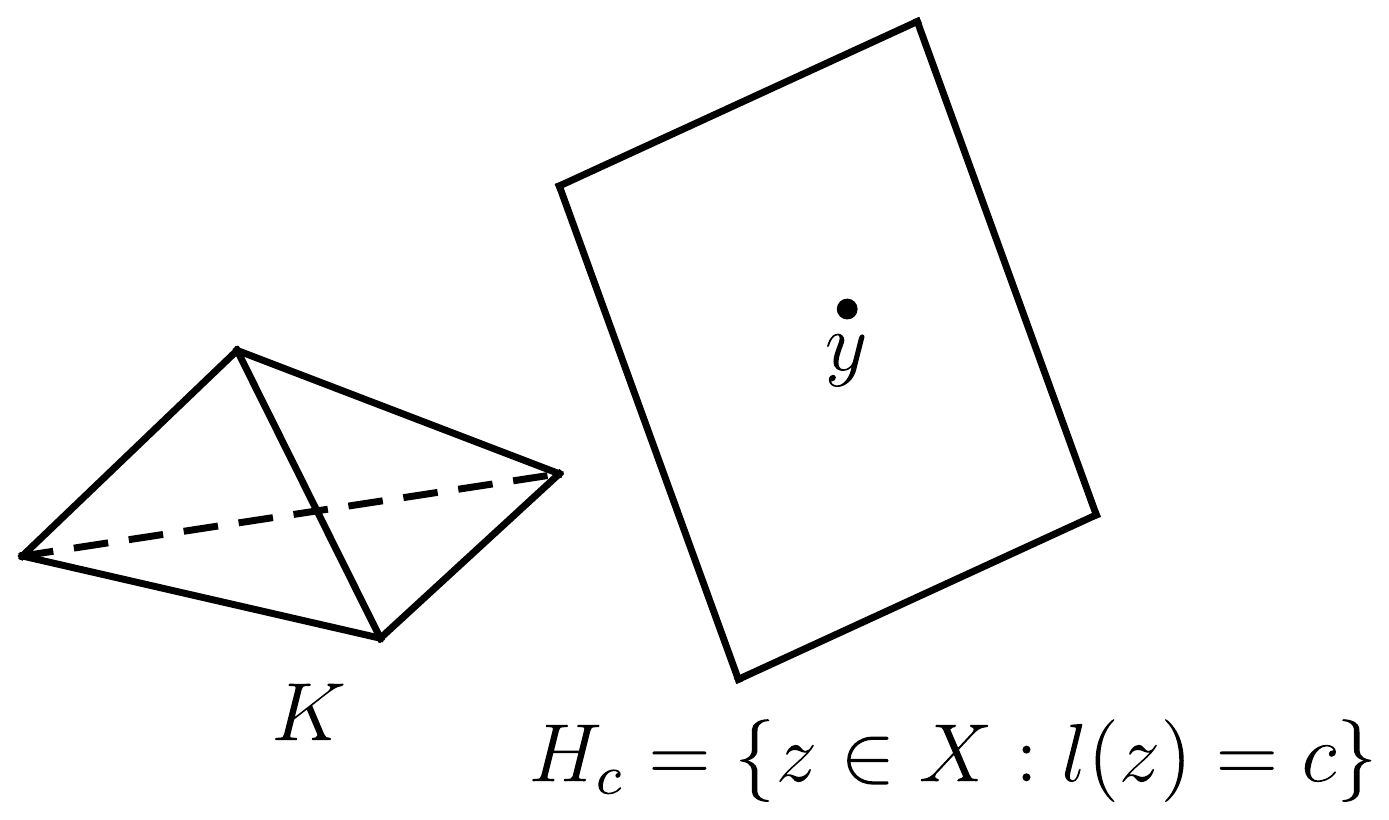}

\protect\caption{\label{fig:hb}Separation of $K$ and $y$ by the hyperplane $H_{c}$.}
\end{figure}

\begin{thm}[HB3-Separation]
 Let $X$ be a real vector space. Assume $X$ is equipped with a
topology making the two vector-operations continuous. Let $K\neq\emptyset$
be an open convex subset of $X$. Let $y\in X\backslash K$ (in the
complement). 

Then there is a linear functional $l:X\rightarrow\mathbb{R}$ such
that 
\begin{equation}
l\left(x\right)<l\left(y\right),\quad\forall x\in K.\label{eq:hb8}
\end{equation}
(In fact, there exists $c\in\mathbb{R}$ such that $l\left(x\right)<c$,
$\forall x\in K$, and $l\left(y\right)=c$.)
\end{thm}
\uline{Hint}: Assume (by translation) that $0\in K$; and set 
\begin{equation}
p_{K}\left(x\right):=\inf\left\{ a\::\:a\in\mathbb{R}_{+},\frac{x}{a}\in K\right\} ;\label{eq:hb9}
\end{equation}
and then apply version 2 to $p_{K}$, which can be shown to be sub-additive.
For the separation property, see  \figref{hb}.

\begin{thm}[HB4]
 Let $\mathfrak{A}$ be a $C^{*}$-algebra and let $\mathfrak{B}\subset\mathfrak{A}$
be a $*$-subalgebra. Let $l:\mathfrak{B}\rightarrow\mathbb{C}$ satisfy
\[
l\left(b^{*}b\right)\geq0,\;\forall b\in\mathfrak{B},\quad\mbox{and}\quad\left\Vert l\right\Vert =1\;(\mbox{positivity}),
\]
then there is a positive linear functional $\widetilde{l}:\mathfrak{A}\rightarrow\mathbb{C}$,
such that 
\begin{enumerate}
\item $\widetilde{l}$ extends $l$ on $\mathfrak{B}$; 
\item $\widetilde{l}\left(a^{*}a\right)\geq0$, $\forall a\in\mathfrak{A}$;
and
\item $\Vert\widetilde{l}\Vert=1$. 
\end{enumerate}
\end{thm}
\begin{rem}
This version is due to M. Krein, but its proof uses the same ideas
which we sketched above in versions 1-2.
\end{rem}
\end{subappendices}

\chapter{Unbounded Operators in Hilbert Space\label{chap:lin}}
\begin{quotation}
We were {[}initially{]} entirely in Heisenberg's footsteps. He had
the idea that one should take matrices, although he did not know that
his dynamical quantities were matrices.... And when one had such a
programme of formulating everything in matrix language, it takes some
effort to get rid of matrices. Though it seemed quite natural for
me to represent perturbation theory in the algebraic way, this was
not a particularly new way.

--- Max Born\vspace{1em}\sindex[nam]{Born, M., (1882-1970)}\\
...practical methods of applying quantum mechanics should be developed,
which can lead to an explanation of the main features of complex atomic
systems without too much computation. 

--- Paul Adrien Maurice Dirac\sindex[nam]{Dirac, P.A.M., (1902-1984)}\vspace{1em}\\
\textquotedblleft ...Mathematics \dots{} are not only part of a special
science, but are also closely connected with our general culture \dots ,
a bridge to the Arts and Sciences, and the seemingly so non-exact
sciences ...Our purpose is to help build such a bridge. Not for the
sake of history but for the genesis of problems, facts and proofs,
... By going back to the roots of these conceptions, back through
the dust of times past, the scars of long use would disappear, and
they would be reborn to us as creatures full of life.\textquotedblright{} 

--- Otto Toeplitz, 1926\sindex[nam]{Toeplitz, O., (1881-1940)}\vspace{2em}
\end{quotation}
Quantum physics is one of the sources of problems in Functional Analysis,
in particular the study of operators in Hilbert space. In the dictionary
translating between Quantum physics and operators in Hilbert space
we already saw that \textquotedblleft quantum observable\index{observable}s\textquotedblright{}
are \textquotedblleft selfadjoint operators.\textquotedblright{} (See
also Chapter \ref{chap:KS}, and especially  \figref{qmm}.)

As noted, even in a finite number of degrees of freedom, the relevant
operators such as momentum, position, and energy are unbounded. Most
Functional Analysis books stress the bounded case, and below we identify
questions and theorems related to key issues for unbounded linear
operators. (See, e.g., Appendix \ref{sec:stone}.) 

In this chapter, we review the basic theory of unbounded operators
in Hilbert space. For general notions, we refer to \cite{DS88a,DS88b}.\index{space!Hilbert-}

\section{\label{sec:dga}Domain, Graph, and Adjoints}

Among the classes of operators in Hilbert space, the family of selfadjoint
linear operators is crucially important for a host of applications,
e.g., to mathematical physics, and to the study of partial differential
equations (PDE). For a study of each of the three classes of linear
PDOs, elliptic, hyperbolic, and parabolic, the Spectral Theorem (see
\cite{Sto90, Yos95, Ne69, RS75, DS88b}) for associated unbounded
selfadjoint operators is a \textquotedblleft workhorse.\textquotedblright{}
\index{selfadjoint operator}

Every selfadjoint operator is densely defined, is closed, and it is
necessarily (Hermitian) symmetric. For unbounded operators, the converse
fails; although it does hold for bounded operators. It follows that
selfadjointness is a much more restricting property than the related
three properties. Moreover we will see that the distinction (between
\textquotedblleft symmetric\textquotedblright{} and selfadjoint) lies
at the heart of key issues from applications. We will further see
that symmetric operators with dense domain are automatically closable;
but they may, or may not, have selfadjoint extensions; -- again an
issue of importance in physics. \index{selfadjoint operator}

The Spectral Theorem holds for selfadjoint operators, and for normal
operators. But the case of normal operators reduces to the Spectral
Theorem for two commuting selfadjoint operators.

On account of Stone's theorem (Appendix \ref{sec:stone}) for one-parameter
unitary groups we know that the class of selfadjoint operators coincides
precisely with the infinitesimal generators of strongly continuous
one-parameter groups of unitary operators acting on Hilbert space;
-- hence applications to the Schrödinger equation, and to wave equations.
\index{normal operator} \index{strongly continuous}

\index{Schrödinger equation}

Let $\mathscr{H}$ be a complex Hilbert space. An operator $A$ is
a linear mapping whose domain $dom\left(A\right)$ and range $ran\left(A\right)$
are subspaces in $\mathscr{H}$. The kernel $ker\left(A\right)$ of
$A$ consists of all $a\in dom\left(A\right)$ such that $Aa=0$.
The operator $A$ is uniquely determined by its graph \index{graph!- of operator}
\begin{equation}
\mathscr{G}\left(A\right)=\left\{ \left(a,Aa\right):a\in dom\left(A\right)\right\} .\label{eq:operator graph}
\end{equation}
(Here the parentheses are used to denote an ordered pair, rather than
an inner product.) Thus, $\mathscr{G}\left(A\right)$ is a subspace
in $\mathscr{H}\oplus\mathscr{H}$ equipped with the inherited \emph{graph
inner product} 
\begin{align}
\left\langle a,b\right\rangle _{A} & =\left\langle a,b\right\rangle +\left\langle Aa,Ab\right\rangle ,\;\mbox{and}\label{eq:graph inner product}\\
\left\Vert a\right\Vert _{A}^{2} & =\left\langle a,a\right\rangle _{A},\;\forall a,b\in dom\left(A\right)\label{eq:graph norm}
\end{align}
In general, a subspace $K\subset\mathscr{H}\oplus\mathscr{H}$ is
the graph of an operator if and only if $\left(0,a\right)\in K$ implies
$a=0$. 

Given two operators $A$ and $B$, we say $B$ is an extension\index{extension!-of operator}
of $A$, denoted by $B\supset A$, if $\mathscr{G}\left(B\right)\supset\mathscr{G}\left(A\right)$
in $\mathscr{H}\oplus\mathscr{H}$. The operator $A$ is \emph{closable}
if $\overline{\mathscr{G}\left(A\right)}$ is the graph of an operator
$\overline{A}$, namely, the closure of $A$. We say $A$ is \emph{closed}
if $A=\bar{A}$. \index{closable operator} \index{operators!closable-}

Let $A$ be a closed operator. A dense subspace $K\subset\mathscr{H}$
is called a \emph{core} of $A$, if the closure of the restriction
$A\big|_{K}$ is equal to $A$.\index{operators!domain of}\index{operators!graph of}

Let $G$ be the group of all $3\times3$ real matrices $\begin{bmatrix}1 & a & c\\
0 & 1 & b\\
0 & 0 & 1
\end{bmatrix}$ with Lie algebra $\mathfrak{g}$ consisting of matrices $\begin{bmatrix}0 & a & c\\
0 & 0 & b\\
0 & 0 & 0
\end{bmatrix}$, $\left(a,b,c\right)\in\mathbb{R}^{3}$. Let $\mathcal{U}$ be a
strongly continuous unitary representation of $G$ acting on a Hilbert
space $\mathscr{H}$. Then for every $X\in\mathfrak{g}$, $t\longmapsto\mathcal{U}\left(e^{tX}\right)$
defines a strongly continuous one-parameter group; and hence its infinitesimal
generator, denoted $d\mathcal{U}\left(X\right)$ is a skew-adjoint
operator with dense domain in $\mathscr{H}$. \index{strongly continuous}

\index{representation!unitary}

\index{representation!strongly continuous}

For a more systematic account of the interplay between the Lie algebra
and the corresponding Lie group, as it relates to representations,
see also \chapref{groups} below, especially \exerref{lce} for the
present setting.

In the example below, we apply this to $\mathscr{H}=L^{2}\left(\mathbb{R}\right)$,
and the unitary representation $\mathcal{U}$ of $G$ is defined as
follows: For $f\in\mathscr{H}=L^{2}\left(\mathbb{R}\right)$, set
\[
\left(\mathcal{U}\left(g\right)f\right)\left(x\right)=e^{i\left(c+bx\right)}f\left(x+a\right),\;x\in\mathbb{R};
\]
called the \emph{Schrödinger representation}. Differentiating in the
three directions in the Lie algebra, we get \index{representation!Schrödinger-}
\begin{eqnarray}
\left(d\mathcal{U}\left(X_{1}\right)f\right)\left(x\right) & = & f'\left(x\right)=\frac{d}{dx}f,\label{eq:sr3-1}\\
\left(d\mathcal{U}\left(X_{2}\right)f\right)\left(x\right) & = & ix\,f\left(x\right),\;\mbox{and}\label{eq:sr3-2}\\
\left(d\mathcal{U}\left(X_{3}\right)f\right)\left(x\right) & = & i\,f\left(x\right).\label{eq:sr3-3}
\end{eqnarray}

The first two operators in (\ref{eq:sr3-1})-(\ref{eq:sr3-2}) are
often written as follows: \index{Schrödinger representation}
\[
d\mathcal{U}\left(X_{1}\right)=iP
\]
where $P$ is the \emph{momentum operator} of a single quantum mechanical
particle (wave function); and 
\[
d\mathcal{U}\left(X_{2}\right)=iQ
\]
where $Q$ is the corresponding \emph{position operator} (in a single
degree of freedom.) These two operators may be realized as follows:

\index{quantum mechanics!momentum operator}

\index{quantum mechanics!position operator}

\index{Schrödinger, E. R.}

\index{operators!momentum-}
\begin{rem}
It would appear that the function $f_{\lambda}\left(x\right)=e^{i\lambda x}$,
$\lambda$ fixed is an eigenfunction\index{eigenfunction} for $P=\frac{1}{i}\frac{d}{dx}$.
For all $\lambda$, 
\[
Pf_{\lambda}=\lambda f_{\lambda}
\]
holds pointwise, but ``spectrum'' depends on the ambient Hilbert
space $\mathscr{H}$, in this case $\mathscr{H}=L^{2}\left(\mathbb{R}\right)$;
and $f_{\lambda}\notin L^{2}\left(\mathbb{R}\right)$, so $\lambda$
is not an eigenvalue. Nonetheless, if we allow an intervals for the
$\lambda$ variable, e.g., $a<\lambda<b$, with $a$ and $b$ being
finite, then\index{eigenvalue} 
\[
F_{a,b}\left(x\right)=\int_{a}^{b}e^{i\lambda x}d\lambda=\frac{e^{ibx}-e^{iax}}{ix}
\]
is in $L^{2}\left(\mathbb{R}\right)$; and hence $P$ has \emph{continuous}
spectrum. The functions $F_{a,b}\left(\cdot\right)$ are examples
of wave-packets in quantum mechanics.

\index{spectrum!continuous-}\end{rem}
\begin{example}
The two operators $d/dx$ and $M_{x}$ in QM, are acting on $L^{2}\left(\mathbb{R}\right)$
with dense domain = the \emph{Schwartz space} $\mathcal{S}$. 
\end{example}
An alternative way to get a dense common domain, a way that works
for all representations, is to use \emph{Gårding space}, or $C^{\infty}$-\emph{vectors}.
\index{Gårding!-vector}

Let $u\in\mathscr{H}$ and define
\[
u_{\varphi}:=\int_{G}\varphi(g)\mathcal{U}_{g}u\,dg
\]
where $\varphi\in C_{c}^{\infty}$, and $\mathcal{U}\in Rep(G,\mathscr{H})$.
Let $\varphi_{\epsilon}$ be an approximation of identity. Then for
functions on $G$, $\varphi_{\epsilon}\star\psi\rightarrow\psi$ as
$\epsilon\rightarrow0$; and for $C^{\infty}$ vectors, $u_{\varphi_{\epsilon}}\rightarrow u$,
as $\epsilon\rightarrow0$ in $\mathscr{H}$, i.e., in the $\left\Vert \cdot\right\Vert _{\mathscr{H}}$-
norm. \index{Gårding!-space}

The set $\{u_{\varphi}\}$ is dense in $\mathscr{H}$. It is called
the Gårding space, or $C^{\infty}$ vectors, or Schwartz space. Notice
that not only $u_{\varphi}$ is dense in $\mathscr{H}$, their derivatives
are also dense in $\mathscr{H}$. 

Differentiating $\mathcal{U}_{g}$, we then get a Lie algebra representation
\[
\rho\left(X\right):=\frac{d}{dt}\big|_{t=0}\mathcal{U}\left(e^{tX}\right)=d\mathcal{U}\left(X\right).
\]

\begin{lem}
$\left\Vert u_{\varphi_{\epsilon}}-u\right\Vert \rightarrow0$, as
$\epsilon\rightarrow0$.\end{lem}
\begin{proof}
Since $u-u_{\varphi_{\epsilon}}=\int\varphi_{\epsilon}(g)(u-\mathcal{U}_{g}u)dg$,
we have
\begin{eqnarray*}
\left\Vert u_{\varphi_{\epsilon}}-u\right\Vert  & = & \left\Vert \int_{G}\varphi_{\epsilon}(g)(u-\mathcal{U}_{g}u)dg\right\Vert \\
 & \leq & \int_{G}\varphi_{\epsilon}(g)\left\Vert u-\mathcal{U}_{g}u\right\Vert dg
\end{eqnarray*}
where the integration on $G$ is with respect to Haar measure, and
where we used the fact that $\int_{G}\varphi_{\epsilon}=1$. Notice
that we always assume the representations are norm continuous in the
$g$ variable, otherwise it is almost impossible to get anything interesting.
i.e., assume $\mathcal{U}$ being strongly continuous. So for all
$\delta>0$, there is a neighborhood $\mathcal{O}$ of $e\in G$ so
that $\left\Vert u-\mathcal{U}_{g}u\right\Vert <\delta$ for all $g\in\mathcal{O}$.
Choose $\epsilon_{\delta}$ so that $\varphi_{\epsilon}$ is supported
in $\mathcal{O}$ for all $\epsilon<\epsilon_{\delta}$. Then the
statement is proved. \index{strongly continuous}

\index{representation!strongly continuous}\end{proof}
\begin{cor}
For all $X\in\mathfrak{g}$, $d\mathcal{U}\left(X\right)$ is essentially
skew-adjoint on the Gårding domain.\end{cor}
\begin{proof}
Let $\mathcal{U}\in Rep\left(G,\mathscr{H}\right)$ be a unitary strongly
continuous representation, where $G$ is a Lie group with Lie algebra
$\mathfrak{g}$. Let $X\in\mathfrak{g}$; we claim that $d\mathcal{U}\left(X\right)$
is essentially skew-adjoint (see \cite{RS75, Ne69, vN32a, DS88b})
on the Gårding space, i.e., the span of vectors \index{Gårding!-vector}
\index{Lie!group} 
\begin{equation}
v_{\varphi}=\int_{G}\varphi\left(g\right)\mathcal{U}_{g}v\:dg\label{eq:gd1}
\end{equation}
where $\varphi\in C_{c}^{\infty}\left(G\right)$. 

Hence we must show that 
\begin{equation}
\ker\left(d\mathcal{U}\left(X\right)^{*}\pm I\right)=0.\label{eq:gd2}
\end{equation}
The argument is the same in both cases of (\ref{eq:gd2}). Hence we
must show that, if $w\in\mathscr{H}$ satisfies 
\begin{equation}
\left\langle d\mathcal{U}\left(X\right)v_{\varphi}-v_{\varphi},w\right\rangle =0\label{eq:gd3}
\end{equation}
for all $v_{\varphi}$ in (\ref{eq:gd1}), then $w=0$. Now if we
view $X$ as an invariant vector field on $G$, then (\ref{eq:gd3})
states that the continuous function
\begin{equation}
f_{w}\left(g\right):=\mathcal{\mathcal{U}}\left(g\right)w\label{eq:gd4}
\end{equation}
is a weak solution to the ODE 
\[
Xf_{w}=f_{w};
\]
equivalently
\begin{equation}
f_{w,X}\left(t\right):=\mathcal{U}\left(\exp\left(tX\right)\right)w\label{eq:dg5}
\end{equation}
satisfies 
\begin{equation}
\frac{d}{dt}f_{w,X}\left(t\right)=f_{w,X}\left(t\right),\quad t\in\mathbb{R};\label{eq:gd6}
\end{equation}
and so 
\begin{equation}
f_{w,X}\left(t\right)=\mbox{const}\cdot e^{t},\;t\in\mathbb{R}.\label{eq:gd7}
\end{equation}
But, since $\mathcal{U}$ is unitary, $f_{w,X}$ in (\ref{eq:dg5})
is bounded; so the constant in (\ref{eq:gd7}) is zero. Hence $f_{w,X}\left(t\right)\equiv0$.
But $f_{w,X}\left(0\right)=w$, and so $w=0$.
\end{proof}
Now, let $A$ be an arbitrary linear operator in a Hilbert space $\mathscr{H}$
with $dom\left(A\right)$ dense in $\mathscr{H}$.
\begin{thm}
The following are equivalent.
\begin{enumerate}
\item $A=\bar{A}$.
\item $\mathscr{G}\left(A\right)=\overline{\mathscr{G}\left(A\right)}$.
\item $dom\left(A\right)$ is a Hilbert space with respect to the graph
inner product $\left\langle \cdot,\cdot\right\rangle _{A}$.
\item If $\left\{ \left(a_{n},Aa_{n}\right)\right\} _{n=1}^{\infty}$ is
a sequence in $\mathscr{G}\left(A\right)$, and $\left(a_{n},Aa_{n}\right)\rightarrow\left(a,b\right)$
as $n\rightarrow\infty$, then $\left(a,b\right)\in\mathscr{G}\left(A\right)$.
In particular, $b=Aa$. (The round braces $\left(\cdot,\cdot\right)$
mean ``pair-of vectors.'' )
\end{enumerate}
\end{thm}
\begin{proof}
All follow from definitions.
\end{proof}
Let $X$ be a vector space over $\mathbb{C}$. Suppose there are two
norms defined on $X$, such that 
\begin{equation}
\left\Vert \cdot\right\Vert _{1}\leq\left\Vert \cdot\right\Vert _{2}.\label{eq:n-1-1}
\end{equation}
Let $\overline{X_{i}}$ be the completion of $X$ with respect to
$\left\Vert \cdot\right\Vert _{i}$, $i=1,2$. The ordering (\ref{eq:n-1-1})
implies the identify map 
\[
\varphi:\left(X,\left\Vert \cdot\right\Vert _{2}\right)\rightarrow\left(X,\left\Vert \cdot\right\Vert _{1}\right)
\]
is continuous, hence it has a unique continuous extension\index{extension!-of functional}
$\widetilde{\varphi}$ to $\overline{X_{2}}$; and (\ref{eq:n-1-1})
passes to the closure $\overline{X_{2}}$. If $\widetilde{\varphi}$
is injective, $\overline{X_{2}}$ is embedded into $\overline{X_{1}}$
as a dense subspace. In that case, $\left\Vert \cdot\right\Vert _{1}$
and $\left\Vert \cdot\right\Vert _{2}$ are said to be\emph{ topologically
consistent}.

\index{completion!norm-}
\begin{lem}
$\left\Vert \cdot\right\Vert _{1}$ and $\left\Vert \cdot\right\Vert _{2}$
are topologically equivalent if and only if 
\[
\begin{Bmatrix}\left\{ x_{n}\right\} \subset X\mbox{ is Cauchy under }\left\Vert \cdot\right\Vert _{2}\\
\mbox{(hence Cauchy under \ensuremath{\left\Vert \cdot\right\Vert _{1}})}\\
\left\Vert x_{n}\right\Vert _{1}\rightarrow0
\end{Bmatrix}\Longrightarrow\left\Vert x_{n}\right\Vert _{2}\rightarrow0.
\]
\end{lem}
\begin{proof}
Note $\widetilde{\varphi}$ is linear, and 
\[
\ker\widetilde{\varphi}=\begin{Bmatrix}x\in\overline{X_{2}}\:\big|\:\exists\left(x_{n}\right)\subset X,\:\left\Vert x_{n}-x\right\Vert _{2}\rightarrow0,\:\widetilde{\varphi}\left(x\right)=0\\
(\mbox{note }\widetilde{\varphi}\left(x\right)=\lim_{n}\varphi\left(x_{n}\right)=\lim_{n}x_{n}\mbox{ in }\overline{X_{1}})
\end{Bmatrix}.
\]
The lemma follows from this.\end{proof}
\begin{lem}
\label{lem:gnorm}The graph norm of $A$ is topologically equivalent
to $\left\Vert \cdot\right\Vert +\left\Vert A\cdot\right\Vert $. \end{lem}
\begin{proof}
This follows from the estimate \index{graph!- of operator}
\[
\frac{1}{2}\left(\left\Vert x\right\Vert +\left\Vert Ax\right\Vert \right)^{2}\leq\left\Vert x\right\Vert ^{2}+\left\Vert Ax\right\Vert ^{2}\leq\left(\left\Vert x\right\Vert +\left\Vert Ax\right\Vert \right)^{2},\;\forall x\in dom\left(A\right).
\]
\end{proof}
\begin{thm}
\label{thm:norm-1}An operator $A$ is closable if and only if $\left\Vert \cdot\right\Vert $
and $\left\Vert \cdot\right\Vert _{A}$ are topologically equivalent.
(When they are, the completion of $dom\left(A\right)$ with respect
to $\left\Vert \cdot\right\Vert _{A}$ is identified as a subspace
of $\mathscr{H}$.)\index{closable operator}\index{completion!norm-}\end{thm}
\begin{proof}
First, assume $A$ is closable. Let $\left\{ x_{n}\right\} $ be a
sequence in $dom\left(A\right)$. Suppose $\left\{ x_{n}\right\} $
is a Cauchy sequence with respect to $\left\Vert \cdot\right\Vert _{A}$,
and $\left\Vert x_{n}\right\Vert \rightarrow0$. We need to show $\left\{ x_{n}\right\} $
converges to $0$ under the $A$-norm, i.e., $\left\Vert x_{n}\right\Vert _{A}\rightarrow0$.
Since $\left\{ \left(x_{n},Ax_{n}\right)\right\} \subset\mathscr{G}\left(A\right)$,
and $A$ is closable, it follows that $\left(x_{n},Ax_{n}\right)\rightarrow\left(0,0\right)\in\mathscr{G}\left(A\right)$.
Therefore, $\left\Vert Ax_{n}\right\Vert \rightarrow0$, and (see
\lemref{gnorm}) 
\[
\left\Vert x_{n}\right\Vert _{A}=\left\Vert x_{n}\right\Vert +\left\Vert Ax_{n}\right\Vert \rightarrow0.
\]

Conversely, assume $\left\Vert \cdot\right\Vert $ and $\left\Vert \cdot\right\Vert _{A}$
are topologically consistent. Let $\left\{ x_{n}\right\} \subset dom\left(A\right)$,
such that 
\begin{equation}
\left(x_{n},Ax_{n}\right)\rightarrow\left(0,b\right)\mbox{ in }\mathscr{H}\oplus\mathscr{H}.\label{eq:A-1-1}
\end{equation}
We proceed to show that $b=0$, which implies that $A$ is closable. 

By (\ref{eq:A-1-1}), $\left\{ x_{n}\right\} \subset dom\left(A\right)$
is a Cauchy sequence with respect to the $\left\Vert \cdot\right\Vert _{A}$-norm,
and $\left\Vert x_{n}\right\Vert \rightarrow0$. Since the two norms
are topologically consistent, then $\left\Vert x_{n}\right\Vert _{A}\rightarrow0$
and so $\left\Vert Ax_{n}\right\Vert \rightarrow0$. We conclude that
$b=0$. \end{proof}
\begin{cor}
An operator $A$ with dense domain is closable if and only if its
adjoint $A^{*}$ has dense domain.\index{operators!adjoint-}
\end{cor}
We will focus on unbounded operators. In the sequel, we will consider
densely defined Hermitian (symmetric) operators. Such operators are
necessarily closable. 

The following result is usually applied to operators whose inverses
are bounded. 
\begin{prop}
\label{prop:norm}Let $A$ be a bounded operator with domain $dom\left(A\right)$,
and acting in $\mathscr{H}$. Then $dom\left(A\right)$ is closed
in $\left\Vert \cdot\right\Vert _{A}$ if and only if it is closed
in $\left\Vert \cdot\right\Vert $. (That is, for bounded operators,
$\left\Vert \cdot\right\Vert $ and $\left\Vert \cdot\right\Vert _{A}$
are topologically equivalent. )\end{prop}
\begin{proof}
This is the result of the following estimate:
\[
\left\Vert x\right\Vert \leq\left\Vert x\right\Vert _{A}=\left\Vert x\right\Vert +\left\Vert Ax\right\Vert \leq\left(1+\left\Vert A\right\Vert \right)\left\Vert x\right\Vert ,\;\forall x\in dom\left(A\right).
\]
\end{proof}
\begin{cor}
If $A$ is a closed operator in $\mathscr{H}$ and $A^{-1}$ is bounded,
then $ran\left(A\right)$ is closed in both $\left\Vert \cdot\right\Vert $
and $\left\Vert \cdot\right\Vert _{A^{-1}}$. \end{cor}
\begin{proof}
Note the $\mathscr{G}\left(A\right)$ is closed iff $\mathscr{G}\left(A^{-1}\right)$
is closed; and $ran\left(A\right)=dom\left(A^{-1}\right)$. Now, apply
\propref{norm} to $A^{-1}$. 
\end{proof}

Let $A$ be an operator in a Hilbert space $\mathscr{H}$. The set
$\mathscr{G}\left(A\right)^{\perp}$ consists of $\left(-b^{*},b\right)$
such that $\left(-b^{*},b\right)\perp\mathscr{G}\left(A\right)$ in
$\mathscr{H}\oplus\mathscr{H}$.
\begin{prop}
\label{pro:op-adjoint operator} The following are equivalent.
\begin{enumerate}
\item $\mathscr{D}\left(A\right)$ is dense in $\mathscr{H}$.
\item $\left(b,0\right)\perp\mathscr{G}\left(A\right)\Longrightarrow b=0$.
\item If $\left(b,-b^{*}\right)\perp\mathscr{G}\left(A\right)$, the map
$b\mapsto b^{*}$ is well-defined.
\end{enumerate}
\end{prop}
\begin{proof}
Let $a\in\mathscr{D}\left(A\right)$, and $b,b^{*}\in\mathscr{H}$;
then 
\[
\left(-b^{*},b\right)\perp\left(a,Aa\right)\mbox{ in }\mathscr{H}\oplus\mathscr{H}\Longleftrightarrow\left\langle b^{*},a\right\rangle =\left\langle b,Aa\right\rangle 
\]
and the desired results follow from this.
\end{proof}
If any of the conditions is satisfied, $A^{*}:b\mapsto b^{*}$ defines
an operator, called the adjoint of $A$, such that\index{operators!adjoint-}
\begin{equation}
\left\langle b,Aa\right\rangle =\left\langle A^{*}b,a\right\rangle \label{eq:def adjoint op}
\end{equation}
for all $a\in\mathscr{D}\left(A\right)$. $\mathscr{G}\left(A\right)^{\perp}$
is the inverted graph of $A^{*}$. The adjoints are only defined for
operators with dense domains in $\mathscr{H}$.
\begin{example}
$A=d/dx$ on $L^{2}[0,1]$ with dense domain 
\[
\mathscr{D}=\left\{ f\in C^{1}\:\big|\:f\left(0\right)=f\left(1\right)=0\right\} .
\]
Integration by parts shows that $A\subset-A^{*}$.
\end{example}
For unbounded operators, $\left(AB\right){}^{*}=B^{*}A^{*}$ does
not hold in general. The situation is better if one of them is bounded. 
\begin{thm}[{\cite[Theorem 13.2]{Ru90}}]
\label{thm:ST}If $S,T,ST$ are densely defined operators then $\left(ST\right)^{*}\supset T^{*}S^{*}$.
If, in addition, $S$ is bounded then $\left(ST\right)^{*}=T^{*}S^{*}$.
\end{thm}
The next theorem follows directly from the definition of the adjoint
operators.
\begin{thm}
\label{thm:Hkerran}If $A$ is densely defined then $\mathscr{H}=\overline{\mathscr{R}\left(A\right)}\oplus\mathscr{K}\left(A^{*}\right)$.
\end{thm}
Finally, we recall some definitions.
\begin{defn}
Let $A$ be a linear operator acting in $\mathscr{H}$. $A$ is said
to be 
\begin{itemize}
\item selfadjoint if $A=A^{*}$.
\item essentially selfadjoint if $\overline{A}=A^{*}$.
\item normal if $A^{*}A=AA^{*}$.
\item regular if $\mathscr{D}\left(A\right)$ is dense in $\mathscr{H}$,
and closed in $\left\Vert \cdot\right\Vert _{A}$. 
\end{itemize}
\end{defn}
\index{operators!regular-}

\index{operators!essentially selfadjoint-}

\index{essentially selfadjoint operator}

\begin{defn}
Let $A$ be a linear operator on a Hilbert space $\mathscr{H}$. The
resolvent $R\left(A\right)$ is defined as 
\[
R\left(A\right)=\left\{ \lambda\in\mathbb{C}:\left(\lambda-A\right)^{-1}\;\mbox{exists}\right\} \;\left(\mbox{the resolvent set}\right)
\]
and the \emph{spectrum} of $A$ is the complement of $R\left(A\right)$,
and it is denoted by $sp\left(A\right)$ or $\sigma\left(A\right)$.
\end{defn}

\begin{xca}[The resolvent identity]
\myexercise{The resolvent identity}\label{exer:resolvent}Let $A$
be a linear operator in a Hilbert space $\mathscr{H}$, and, for $\lambda_{i}\in R\left(A\right)$,
$i=1,2$ consider two operators $\left(\lambda_{i}-A\right)^{-1}$.
Show that 
\[
\left(\lambda_{1}-A\right)^{-1}-\left(\lambda_{2}-A\right)^{-1}=\left(\lambda_{2}-\lambda_{1}\right)\left(\lambda_{1}-A\right)^{-1}\left(\lambda_{2}-A\right)^{-1}.
\]
This formula is called the resolvent identity. \index{resolvent identity}
\end{xca}

\section{\label{sec:char}Characteristic Matrix \index{Stone, M. H.}\index{space!Hilbert-}}

The method of characteristic matrix was developed by M.H. Stone's
\cite{Sto51}. It is extremely useful in operator theory, but has
long been overlooked in the literature. We recall some of its applications
in normal operators. \index{matrix!characteristic} \index{characteristic matrix}\index{Theorem!Stone's-}

If $\mathscr{H}$ is a fixed Hilbert space, and $A$ a given liner
operator, then its graph 
\[
\mathscr{G}\left(A\right)=\left\{ \begin{bmatrix}u\\
Au
\end{bmatrix}\::\:u\in dom\left(A\right)\right\} 
\]
is a linear subspace in $\mathscr{H}\oplus\mathscr{H}$, represented
as column vectors $\begin{bmatrix}u\\
v
\end{bmatrix}$, $u,v\in\mathscr{H}$. 

In the case where $A$ is assumed closed, we now compute the projection
onto $\overline{\mathscr{G}\left(A\right)}=$ the $\mathscr{H}\oplus\mathscr{H}$-closure
of the graph. 

Let $A$ be an operator in a Hilbert space $\mathscr{H}$. Let $P=(P_{ij})$
be the projection from $\mathscr{H}\oplus\mathscr{H}$ onto $\overline{\mathscr{G}(A)}$.
The $2\times2$ operator matrix $(P_{ij})$ of bounded operators in
$\mathscr{H}$ is called the characteristic matrix of $A$. 

Since $P^{2}=P^{*}=P$, the following identities hold 
\begin{equation}
P_{ij}^{*}=P_{ji}\label{eq:cm-entry 1}
\end{equation}
\begin{equation}
\sum_{k}P_{ik}P_{kj}=P_{ij}\label{eq:cm-entry 2}
\end{equation}
In particular, $P_{11}$ and $P_{22}$ are selfadjoint.
\begin{thm}
Let $P=(P_{ij})$ be the projection from $\mathscr{H}\oplus\mathscr{H}$
onto a closed subspace $\mathscr{K}$. The following are equivalent.
\begin{enumerate}
\item $\mathscr{K}$ is the graph of an operator. \index{graph!- of operator}
\item 
\[
\left[\begin{array}{cc}
P_{11} & P_{12}\\
P_{21} & P_{22}
\end{array}\right]\left[\begin{array}{c}
0\\
a
\end{array}\right]=\left[\begin{array}{c}
0\\
a
\end{array}\right]\Longrightarrow a=0.
\]

\item 
\[
\biggl(P_{12}a=0,P_{22}a=a\biggr)\Longrightarrow a=0.
\]

\end{enumerate}

If any of these conditions is satisfied, let $A$ be the operator
with $\mathscr{G}\left(A\right)=\mathscr{K}$, then for all $a,b\in\mathscr{H}$,
\[
\left[\begin{array}{cc}
P_{11} & P_{12}\\
P_{21} & P_{22}
\end{array}\right]\left[\begin{array}{c}
a\\
b
\end{array}\right]=\left[\begin{array}{c}
P_{11}a+P_{12}b\\
P_{21}a+P_{22}b
\end{array}\right]\in\mathscr{G}\left(A\right);
\]
i.e., 
\begin{equation}
A:(P_{11}a+P_{12}b)\mapsto P_{21}a+P_{22}b.\label{eq:cm-define A}
\end{equation}
In particular,
\begin{eqnarray}
AP_{11} & = & P_{21}\label{eq:cm-define A entry 1}\\
AP_{12} & = & P_{22}\label{eq:cm-define A entry 2}
\end{eqnarray}

\end{thm}
\begin{proof}
Let $v:=\left(a,b\right)\in\mathscr{H}\oplus\mathscr{H}$. Then $v\in\mathscr{K}$
if and only if $Pv=v$; and the theorem follows from this.
\end{proof}
The next theorem describes the adjoint operators.
\begin{thm}
Let $A$ be an operator with characteristic matrix $P=(P_{ij})$.
The following are equivalent.\index{operators!adjoint-}
\begin{enumerate}[itemsep=1em]
\item $\mathscr{D}\left(A\right)$ is dense in $\mathscr{H}$.
\item $\left[\begin{array}{c}
b\\
0
\end{array}\right]\perp\mathscr{G}\left(A\right)=0\Longrightarrow b=0.$
\item If $\begin{bmatrix}-b^{*}\\
b
\end{bmatrix}\in\mathscr{G}\left(A\right)^{\perp}$, the map $A^{*}:b\mapsto b^{*}$ is a well-defined operator.
\item $\left[\begin{array}{cc}
1-P_{11} & -P_{12}\\
-P_{21} & 1-P_{22}
\end{array}\right]\left[\begin{array}{c}
b\\
0
\end{array}\right]=\left[\begin{array}{c}
b\\
0
\end{array}\right]\Longrightarrow b=0.$
\item $\biggl(P_{11}b=0,P_{21}b=0\biggr)\Longrightarrow b=0.$ 
\end{enumerate}

If any of the above conditions is satisfied, then 
\[
\left[\begin{array}{cc}
1-P_{11} & -P_{12}\\
-P_{21} & 1-P_{22}
\end{array}\right]\left[\begin{array}{c}
a\\
b
\end{array}\right]=\left[\begin{array}{c}
(1-P_{11})a-P_{12}b\\
(1-P_{22})b-P_{21}a
\end{array}\right]\in\mathscr{G}\left(A\right)^{\perp}
\]
that is, 
\begin{equation}
A^{*}:P_{21}a-(1-P_{22})b\mapsto(1-P_{11})a-P_{12}b.\label{eq:cm-define A adjoint}
\end{equation}
In particular,
\begin{eqnarray}
A^{*}P_{21} & = & 1-P_{11}\label{eq:cm-adjont entry 1}\\
A^{*}(1-P_{22}) & = & P_{12}.\label{eq:cm-adjoint entry 2}
\end{eqnarray}

\end{thm}
\begin{proof}
$(1)\Leftrightarrow(2)\Leftrightarrow(3)$ is a restatement of Proposition
\ref{pro:op-adjoint operator}. Note the projection from $\mathscr{H}\oplus\mathscr{H}$
on $ $$\mathscr{G}\left(A\right)^{\perp}=\overline{\mathscr{G}\left(A\right)}^{\perp}$
is 
\[
1-P=\begin{bmatrix}1-P_{11} & -P_{12}\\
-P_{21} & 1-P_{22}
\end{bmatrix}
\]
and so 
\[
\left[\begin{array}{c}
b\\
0
\end{array}\right]\perp\mathscr{G}\left(A\right)\Longleftrightarrow\left(1-P\right)\left[\begin{array}{c}
b\\
0
\end{array}\right]=\left[\begin{array}{c}
b\\
0
\end{array}\right].
\]
Therefore, $(2)\Leftrightarrow(4)\Leftrightarrow(5)$. Finally, (\ref{eq:cm-define A adjoint})-(\ref{eq:cm-adjoint entry 2})
follow from the definition of $A^{*}$.\end{proof}
\begin{thm}
\label{thm:cmstone}Let $A$ be a regular operator (i.e., densely
defined, closed) with characteristic matrix $P=(P_{ij})$. 
\begin{enumerate}
\item \label{enu:cm-stone thm 1 a} The matrix entries $P_{ij}$ are given
by
\begin{equation}
\begin{array}{ccccccc}
P_{11} & = & (1+A^{*}A)^{-1} &  & P_{12} & = & A^{*}(1+AA^{*})^{-1}\\
P_{21} & = & A(1+A^{*}A)^{-1} &  & P_{22} & = & AA^{*}(1+AA^{*})^{-1}
\end{array}\label{eq:cm-regular operator, c matrix}
\end{equation}

\item \label{enu:cm-stone thm 1 b} \textup{$1-P_{22}=(1+AA^{*})^{-1}$.}
\item \label{enu:cm-stone thm 1 c} $1+A^{*}A$, $1+AA^{*}$ are selfadjoint
operators. \index{selfadjoint operator}
\item \label{enu:cm-stone thm 1 d} The following containments hold
\begin{eqnarray}
A^{*}(1+AA^{*})^{-1} & \supset & (1+A^{*}A)^{-1}A^{*}\label{eq:cm-off diag entries 1}\\
A(1+A^{*}A)^{-1} & \supset & (1+AA^{*})^{-1}A\label{eq:cm-off diag entries 2}
\end{eqnarray}

\end{enumerate}
\end{thm}
\begin{proof}
By (\ref{eq:cm-define A entry 1}) and (\ref{eq:cm-adjont entry 1}),
we have
\[
\begin{bmatrix}AP_{11}=P_{21}\\
A^{*}P_{21}=1-P_{11}
\end{bmatrix}\Longrightarrow A^{*}AP_{11}=1-P_{11},\;\mbox{i.e.,}\;(1+A^{*}A)P_{11}=1.
\]
That is, $1+A^{*}A$ is a Hermitian extension of $P_{11}^{-1}$. By
(\ref{eq:cm-entry 1}), $P_{11}$ is selfadjoint and so is $P_{11}^{-1}$.
Therefore, $1+A^{*}A=P_{11}^{-1}$, or 
\[
P_{11}=(1+A^{*}A)^{-1}.
\]
By (\ref{eq:cm-define A entry 1}), 
\[
P_{21}=AP_{11}=A(1+A^{*}A)^{-1}.
\]

Similarly, by (\ref{eq:cm-define A entry 2}) and (\ref{eq:cm-adjoint entry 2}),
we have 
\[
\begin{bmatrix}AP_{12}=P_{22}\\
A^{*}\left(1-P_{22}\right)=P_{12}
\end{bmatrix}\Longrightarrow AA^{*}\left(1-P_{22}\right)=P_{22},\;\mbox{i.e., }
\]
\[
(1+AA^{*})(1-P_{22})=1.
\]
This means $1+AA^{*}\supset\left(1-P_{22}\right)^{-1}$ is a Hermitian
extension of the selfadjoint operator $\left(1-P_{22}\right)^{-1}$
(note $P_{22}$ is selfadjoint), hence $1+AA^{*}$ is selfadjoint,
and 
\[
1-P_{22}=(1+AA^{*})^{-1}.
\]
By (\ref{eq:cm-adjoint entry 2}), 
\[
P_{12}=A^{*}(1-P_{22})=A^{*}(1+AA^{*})^{-1}.
\]
By (\ref{eq:cm-define A entry 2}),
\[
P_{22}=AP_{12}=AA^{*}(1+AA^{*})^{-1}
\]

We have proved (\ref{enu:cm-stone thm 1 a}), (\ref{enu:cm-stone thm 1 b})
and (\ref{enu:cm-stone thm 1 c}).

Finally, 
\[
P_{12}=P_{21}^{*}=(AP_{11})^{*}\supset P_{11}A^{*}
\]
yields (\ref{eq:cm-off diag entries 1}); and 
\[
P_{21}=P_{12}^{*}=(A^{*}(1-P_{22}))^{*}\supset(1-P_{22})A
\]
gives (\ref{eq:cm-off diag entries 2}). \end{proof}
\begin{xca}[$A^{**}=\overline{A}$]
\myexercise{$A^{**}=\overline{A}$} Let $A$ be a \emph{regular}
operator in a Hilbert space (i.e., we assume that $A$ has dense domain
and is closable.) Then show that
\begin{equation}
A^{**}=\overline{A}\label{eq:ga1}
\end{equation}
where $\overline{A}$ denotes the \emph{closure} of $A$; i.e., $\mathscr{G}\left(\overline{A}\right)=\overline{\mathscr{G}\left(A\right)}$.

\uline{Hint}: Establish the desired identity (\ref{eq:ga1}) by
justifying the following steps: 

Set $\chi:\mathscr{H}^{2}\longrightarrow\mathscr{H}^{2}$, 
\[
\chi\begin{pmatrix}x\\
y
\end{pmatrix}=\begin{pmatrix}-y\\
x
\end{pmatrix},\;\begin{pmatrix}x\\
y
\end{pmatrix}\in\mathscr{H}^{2}.
\]
Then 
\begin{eqnarray*}
\mathscr{G}\left(A^{**}\right) & = & \left(\chi\mathscr{G}\left(A^{*}\right)\right)^{\perp}\\
 & = & \left(\chi\left(\chi\mathscr{G}\left(A\right)\right)^{\perp}\right)^{\perp}\\
 & = & \left(\chi^{2}\mathscr{G}\left(A\right)\right)^{\perp\perp}\\
 & = & \left(\mathscr{G}\left(A\right)\right)^{\perp\perp}=\overline{\mathscr{G}\left(A\right)}=\mathscr{G}\left(\overline{A}\right).
\end{eqnarray*}

\end{xca}

\subsection{Commutants}

\index{commutant}

Let $A$, $B$ be operators in a Hilbert space $\mathscr{H}$, and
suppose $B$ is bounded. The operator $B$ is said to commute (strongly)
with $A$ if $BA\subset AB$.
\begin{lem}
\label{lem:commutant pass closure}Assume that $\overline{A}$ exists.
Then $B$ commutes with $A$ if and only if $B$ commutes with $\overline{A}$. \end{lem}
\begin{proof}
Suppose $BA\subset AB$, and we check that $B\overline{A}\subset\overline{A}B$.
The converse is trivial. For $(a,\overline{A}a)\in\mathscr{G}(\overline{A})$,
choose a sequence $(a_{n},Aa_{n})\in\mathscr{G}(A)$ such that $(a_{n},Aa_{n})\rightarrow(a,\overline{A}a)$.
By assumption, $(Ba_{n},ABa_{n})=(Ba_{n},BAa_{n})\in\mathscr{G}(A)$.
Thus, 
\[
\left(Ba_{n},ABa_{n}\right)\rightarrow\left(Ba,B\overline{A}a\right)\in\mathscr{G}(\overline{A}).
\]
That is, $Ba\in\mathscr{D}(\overline{A})$ and $\overline{A}Ba=B\overline{A}a$. \end{proof}
\begin{lem}
\label{lem:cmatrix} Let $A$ be a closed operator with characteristic
matrix $P=(P_{ij})$. Let $B$ be a bounded operator, and \index{characteristic matrix}
\[
Q_{B}:=\left[\begin{array}{cc}
B & 0\\
0 & B
\end{array}\right].
\]

\begin{enumerate}
\item $B$ commutes with $A$ $\Leftrightarrow$ $B$ leaves $\mathscr{G}(A)$
invariant $\Leftrightarrow$ $Q_{B}P=PQ_{B}P$. 
\item $B$ commutes with $P_{ij}$ $\Leftrightarrow$ $Q_{B}P=PQ_{B}$ $\Leftrightarrow$
$Q_{B^{*}}P=PQ_{B^{*}}$ $\Leftrightarrow$ $B^{*}$ commutes with
$P_{ij}$. 
\item If $B,B^{*}$ commute with $A$, then $B,B^{*}$ commute with $P_{ij}$.
\end{enumerate}
\end{lem}
\begin{proof}
Obvious.\index{operators!closed-}
\end{proof}
A closed operator is said to be \emph{affiliated} with a von Neumann
algebra\index{algebras!von Neumann algebra ($W^{*}$-algebra)} $\mathfrak{M}$
if it commutes with every unitary operator in $\mathfrak{M}'.$ By
\cite[Thm 4.1.7]{MR1468229}, every operator in $\mathfrak{M}'$ can
be written as a finite linear combination of unitary operators in
$\mathfrak{M}'$. Thus, $A$ is affiliated with $\mathfrak{M}$ if
and only if $A$ commutes with every operator in $\mathfrak{M}'$. 

\index{affiliated with}
\begin{rem}
Let $\mathfrak{M}$ be a von Neumann algebra. Let $x\in\mathfrak{M}$
such that $\left\Vert x\right\Vert \leq1$ and $x=x^{*}$. Set $y:=x+i\sqrt{1-x^{2}}$.
Then, $y^{*}y=yy^{*}=x^{2}+1-x^{2}=1$, i.e., $y$ is unitary. Also,
$x=\left(y+y^{*}\right)/2$.\end{rem}
\begin{thm}
\label{thm:cmatrix} Let $A$ be a closed operator with characteristic
matrix $P=(P_{ij})$. Let $\mathfrak{M}$ be a von Neumann algebra,
and 
\[
Q_{B}:=\left[\begin{array}{cc}
B & 0\\
0 & B
\end{array}\right],B\in\mathfrak{M}'.
\]
The following are equivalent:
\begin{enumerate}
\item $A$ is affiliated with $\mathfrak{M}$.
\item $PQ_{B}=Q_{B}P$, for all $B\in\mathfrak{M}'$.
\item $P_{ij}\in\mathfrak{M}$.
\item If $\mathscr{D}(A)$ is dense, then $A^{*}$ is affiliated with $\mathfrak{M}$. 
\end{enumerate}
\end{thm}
\begin{proof}
Notice that $\mathfrak{M}$ is selfadjoint. The equivalence of $1,2,3$
is a direct consequence of \lemref{cmatrix}. 

$P^{\perp}:=1-P$ is the projection onto the inverted graph of $A^{*}$,
should the latter exists. $PQ_{B}=Q_{B}P$ if and only if $P^{\perp}Q_{B}=Q_{B}P^{\perp}$.
Thus, $1$ is equivalent to $4$.
\end{proof}

\section{\label{sec:normal}Normal Operators}

\thmref{AAvN} below concerning operators of the form $A^{*}A$ is
an application of Stone's characteristic matrix.

\index{operators!normal}

\index{von Neumann, J.}
\begin{thm}[von Neumann]
\label{thm:AAvN} If $A$ is a regular operator in a Hilbert space
$\mathscr{H}$, then
\begin{enumerate}
\item $A^{*}A$ is selfadjoint;
\item $\mathscr{D}(A^{*}A)$ is a core of $A$, i.e., 
\[
\overline{A\big|_{\mathscr{D}\left(A^{*}A\right)}}=A;
\]

\item In particular, $\mathscr{D}(A^{*}A)$ is dense in $\mathscr{H}$.
\end{enumerate}
\end{thm}
\begin{proof}
By \thmref{cmstone}, $A^{*}A=P_{11}^{-1}-1$. Since $P_{11}$ is
selfadjoint, so is $P_{11}^{-1}$. Thus, $AA^{*}$ is selfadjoint.\index{operators!regular-}

Suppose $(a,Aa)\in\mathscr{G}(A)$ such that 
\[
(a,Aa)\perp\mathscr{G}(A\big|_{\mathscr{D}\left(A^{*}A\right)});\;\mbox{i.e.,}
\]
\[
\left\langle a,b\right\rangle +\left\langle Aa,Ab\right\rangle =\left\langle a,\left(1+A^{*}A\right)b\right\rangle =0,\;\forall b\in\mathscr{D}\left(A^{*}A\right).
\]
Since $1+A^{*}A=P_{11}^{-1}$, and $P_{11}$ is a bounded operator,
then 
\[
\mathscr{R}(1+A^{*}A)=\mathscr{D}(P_{11})=\mathscr{H}.
\]
It follows that $a\perp\mathscr{H}$, and so $a=0$. \end{proof}
\begin{thm}[von Neumann]
\label{thm:VN normal operator} Let A be a regular operator in a
Hilbert space $\mathscr{H}$. Then $A$ is normal if and only if $\mathscr{D}(A)=\mathscr{D}(A^{*})$
and $\left\Vert Aa\right\Vert =\left\Vert A^{*}a\right\Vert $, for
all $a\in\mathscr{D}(A)$.\end{thm}
\begin{proof}
Suppose $A$ is normal. Then for all $a\in\mathscr{D}\left(A^{*}A\right)\left(=\mathscr{D}\left(AA^{*}\right)\right)$,
we have 
\[
\left\Vert Aa\right\Vert ^{2}=\left\langle Aa,Aa\right\rangle =\left\langle a,A^{*}Aa\right\rangle =\left\langle a,AA^{*}a\right\rangle =\left\langle Aa,A^{*}a\right\rangle =\left\Vert A^{*}a\right\Vert ^{2};
\]
i.e., $\left\Vert Aa\right\Vert =\left\Vert A^{*}a\right\Vert $,
for all $a\in\mathscr{D}\left(A^{*}A\right)$. It follows that 
\[
\mathscr{D}\left(\overline{A\big|_{\mathscr{D}\left(A^{*}A\right)}}\right)=\mathscr{D}\left(\overline{A^{*}\big|_{\mathscr{D}\left(AA^{*}\right)}}\right).
\]
By \thmref{AAvN}, $\mathscr{D}\left(A\right)=\mathscr{D}(\overline{A\big|_{\mathscr{D}\left(A^{*}A\right)}})$
and $\mathscr{D}\left(A^{*}\right)=\mathscr{D}(\overline{A^{*}\big|_{\mathscr{D}\left(AA^{*}\right)}})$.
Therefore, $\mathscr{D}(A)=\mathscr{D}(A^{*})$ and $\left\Vert Aa\right\Vert =\left\Vert A^{*}a\right\Vert $,
for all $a\in\mathscr{D}(A)$.

Conversely, the map $Aa\mapsto A^{*}a$, $a\in\mathscr{D}(A)$, extends
uniquely to a partial isometry $V$ with initial space $\overline{\mathscr{R}(A)}$
and final space $\overline{\mathscr{R}(A^{*})}$, such that $A^{*}=VA$.
By \thmref{ST}, $A=A^{*}V^{*}$. Then $A^{*}A=A^{*}(V^{*}V)A=(A^{*}V^{*})(VA)=AA^{*}$.
Thus, $A$ is normal. 
\end{proof}
The following theorem is due to M.H. Stone. \index{isometry!partial-}
\begin{thm}
\label{thm:normal operators} Let $A$ be a regular operator in a
Hilbert space $\mathscr{H}$. Let $P=(P_{ij})$ be the characteristic
matrix of $A$. The following are equivalent.
\begin{enumerate}
\item $A$ is normal.
\item $P_{ij}$ are mutually commuting.
\item $A$ is affiliated with an abelian von Neumann algebra.
\end{enumerate}
\end{thm}
\begin{rem}
For the equivalence of $1$ and $2$, we refer to the original paper
of Stone. The most interesting part is $1\Leftrightarrow3$. The idea
of characteristic matrix gives rise to an elegant proof without reference
to the spectral theorem.\index{Theorem!Spectral-}\end{rem}
\begin{proof}[Proof of Theorem \ref{thm:normal operators}]
 Assuming $1\Leftrightarrow2$, we prove that $1\Leftrightarrow3$. 

Suppose $A$ is normal, i.e. $P_{ij}$ are mutually commuting. Then
$A$ is affiliated with the abelian von Neumann algebra $\{P_{ij}\}''$.
For if $B\in\{P_{ij}\}'$, then $B$ commutes $P_{ij}$, and so $B$
commutes with $A$ by \lemref{cmatrix}. 

Conversely, if $A$ is affiliated with an abelian von Neumann algebra
$\mathfrak{M}$, then by \thmref{cmatrix}, $P_{ij}\in\mathfrak{M}$
. This shows that $P_{ij}$ are mutually commuting, and $A$ is normal. 
\end{proof}

\section{\label{sec:polar}Polar Decomposition}

We show that the intuition behind the familiar polar decomposition
(or polar factorization) for complex numbers carries over remarkably
well to operators in Hilbert space. Indeed (\thmref{pd}) the operators
that admit a polar decomposition are precisely the regular operators,
meaning closable and with dense domain.

\index{polar decomposition}

Let $A$ be a regular operator in a Hilbert space $\mathscr{H}$.
By \thmref{AAvN}, $A^{*}A$ is a positive selfadjoint operator and
it has a unique positive square root $\left|A\right|:=\sqrt{A^{*}A}$.\index{operators!regular-}
\begin{thm}
~
\begin{enumerate}
\item $\left|A\right|:=\sqrt{A^{*}A}$ is the unique positive selfadjoint
operator $T$ satisfying $\mathscr{D}(T)=\mathscr{D}(A)$, and $\left\Vert Ta\right\Vert =\left\Vert Aa\right\Vert $
for all $a\in\mathscr{D}(A)$. 
\item $\ker\left(\left|A\right|\right)=\ker\left(A\right)$, $\overline{\mathscr{R}(\left|A\right|)}=\overline{\mathscr{R}(A^{*})}$. 
\end{enumerate}
\end{thm}
\begin{proof}
Suppose $T=\sqrt{A^{*}A}$, i.e. $T^{*}T=A^{*}A$. Let $\mathscr{D}:=\mathscr{D}(T^{*}T)=\mathscr{D}(A^{*}A)$.
By \thmref{AAvN}, $\mathscr{D}$ is a core of both $T$ and $A$.
Moreover, $\left\Vert Ta\right\Vert =\left\Vert Aa\right\Vert $,
for all $a\in\mathscr{D}$. We conclude from this norm identity that
$\mathscr{D}(T)=\mathscr{D}(A)$ and $\left\Vert Ta\right\Vert =\left\Vert Aa\right\Vert $,
for all $a\in\mathscr{D}(A)$.

Conversely, suppose $T$ has the desired properties. For all $a\in\mathscr{D}(A)=\mathscr{D}(T)$,
and $b\in\mathscr{D}(A^{*}A)$, 
\[
\left\langle Tb,Ta\right\rangle =\left\langle Ab,Aa\right\rangle =\left\langle A^{*}Ab,a\right\rangle 
\]
This implies that $Tb\in\mathscr{D}(T^{*})=\mathscr{D}(T)$, $T^{2}b=A^{*}Ab$,
for all $b\in\mathscr{D}(A^{*}A)$. That is, $T^{2}$ is a selfadjoint
extension of $A^{*}A$. Since $A^{*}A$ is selfadjoint, $T^{2}=A^{*}A$. 

The second part follows from \thmref{Hkerran}.
\end{proof}
Consequently, the map $\left|A\right|a\mapsto Aa$ extends to a unique
partial isometry $V$ with initial space $\overline{\mathscr{R}(A^{*})}$
and final space $\overline{\mathscr{R}(A)}$ (the overbar means ``norm-closure''),
such that
\begin{equation}
A=V\left|A\right|.\label{eq:polar decomposition}
\end{equation}
Equation (\ref{eq:polar decomposition}) is called the \emph{polar
decomposition} of $A$. It is clear that such decomposition is unique.
\index{polar decomposition} \index{isometry!partial-}

We have proved:
\begin{thm}
\label{thm:pd}Let $A$, $V$ and $\left|A\right|$ be as described;
then 
\[
A=V\left|A\right|.
\]

\end{thm}
Taking adjoints in (\ref{eq:polar decomposition}) yields $A^{*}=\left|A\right|V^{*}$,
so that
\begin{equation}
AA^{*}=VA^{*}AV^{*}\label{eq:AAstar and Astar A}
\end{equation}
Restrict $AA^{*}$ to $\overline{\mathscr{R}(A)}$, and restrict $A^{*}A$
restricted to $\overline{\mathscr{R}(A^{*})}$. Then the two restrictions
are unitarily equivalent. It follows that $A^{*}A$, $AA^{*}$ have
the same spectrum, aside from possibly the point $0$. 

By (\ref{eq:AAstar and Astar A}), $\left|A^{*}\right|=V\left|A\right|V^{*}=VA^{*}$,
where $\left|A^{*}\right|=\sqrt{AA^{*}}$. Apply $V^{*}$ on both
sides gives 
\begin{equation}
A^{*}=V^{*}\left|A^{*}\right|.\label{eq:polar demp Astar}
\end{equation}
By uniqueness, (\ref{eq:polar demp Astar}) is the polar decomposition
of $A^{*}$.
\begin{thm}
$A$ is affiliated with a von Neumann algebra $\mathfrak{M}$ if and
only if $\left|A\right|$ is affiliated with $\mathfrak{M}$ and $V\in\mathfrak{M}$.\end{thm}
\begin{proof}
Let $U$ be a unitary operator in $\mathfrak{M}'$. The operator $UAU^{*}$
has polar decomposition
\[
UAU^{*}=(UVU^{*})(U\left|A\right|U^{*}).
\]
By uniqueness, $A=UAU^{*}$ if and only if $V=UVU^{*}$, $\left|A\right|=U\left|A\right|U^{*}$.
Since $U$ is arbitrary, we conclude that $V\in\mathfrak{M}$, and
$A$ is affiliated with $\mathfrak{M}$. 
\end{proof}

\section*{A summary of relevant numbers from the Reference List}

\begin{flushleft}
For readers wishing to follow up sources, or to go in more depth with
topics above, we suggest: 
\par\end{flushleft}

\begin{flushleft}
Of these, refs \cite{vN32a} and \cite{DS88b} are especially central.
A more comprehensive list is: \cite{BR81,DS88b,Jor08,Kat95,MR1468230,Sto51,Sto90,Wei03,Yos95,JL01,Die75,EMC00,Jor88,Jor94,RS75,MR0184042,MR1819613,BR79,Con90,BrRo66,FL28,Fri80,MR887102,JM84,Kre46,Ne69,vN32a,Hel13}. 
\par\end{flushleft}

\begin{subappendices}

\section{\label{sec:stone}Stone's Theorem}

The gist of the result (\thmref{stone}) is as follows: Given a fixed
Hilbert space, there is then a 1-1 correspondence between any two
in pairs from the following three: (i) strongly continuous unitary
one-parameter groups $\mathcal{U}(t)$; (ii) selfadjoit operators
$H$ (generally unbounded) with dense domain; and (iii) projection
valued measures $P\left(\cdot\right)$, abbreviated PVM. Starting
with $\mathcal{U}(t)$, we say that the corresponding selfadjoint
operator $H$ is its generator, and then the PVM $P\left(\cdot\right)$
will be from the Spectral Theorem applied to $H$.
\begin{defn}[Projection valued measure (PVM)]
 Let $\mathcal{B}\left(\mathbb{R}\right)$ be the Borel sigma algebra
of subsets of $\mathbb{R}$. Let $\mathscr{H}$ be a Hilbert space.
A function $P:\mathcal{B}\left(\mathbb{R}\right)\rightarrow\mbox{Proj}\left(\mathscr{H}\right)$
is called a projection valued measure (PVM) iff (Def), $P\left(\emptyset\right)=0$;
$P\left(\mathbb{R}\right)=I_{\mathscr{H}}$; and for all $\left(E_{i}\right)_{i=1}^{\infty}$
such that $E_{i}\cap E_{j}=\emptyset$ ($i\neq j$), we have:
\begin{equation}
P\left(\bigcup\nolimits _{i}E_{i}\right)=\sum_{i}P\left(E_{i}\right)\label{eq:st1}
\end{equation}

\end{defn}

\begin{defn}
A unitary one-parameter group is a function:
\[
\mathcal{U}:\mathbb{R}\longrightarrow\big(\mbox{unitary operators in \ensuremath{\mathscr{H}}}\big)
\]
such that:
\begin{equation}
\mathcal{U}\left(s+t\right)=\mathcal{U}\left(s\right)\mathcal{U}\left(t\right),\quad\forall s,t\in\mathbb{R};\label{eq:st2}
\end{equation}
and for $\forall h\in\mathscr{H}$, 
\begin{equation}
\lim_{t\rightarrow0}\mathcal{U}\left(t\right)h=h\;\left(\mbox{strong continuity}\right).\label{eq:st3}
\end{equation}
\end{defn}
\begin{thm}[Stone's Theorem \cite{MR1892228,RS75,Rud73}]
\label{thm:stone} There is a sequence of bijective correspondences
between (\ref{enu:st1})-(\ref{enu:st3}) below, i.e., (\ref{enu:st1})$\Rightarrow$(\ref{enu:st2})$\Rightarrow$(\ref{enu:st3})$\Rightarrow$(\ref{enu:st1}): 
\begin{enumerate}
\item \label{enu:st1}PVMs $P\left(\cdot\right)$;
\item \label{enu:st2}unitary one-parameter groups $\mathcal{U}$; and
\item \label{enu:st3}selfadjoint operators $H$ with dense domain in $\mathscr{H}$.
\end{enumerate}

The correspondence is given explicitly as follows:

$\text{\ensuremath{\left(\ref{enu:st1}\right)} \ensuremath{\Rightarrow} \ensuremath{\left(\ref{enu:st2}\right)}:}$
Given $P$, a PVM, set 
\begin{equation}
\mathcal{U}\left(t\right)=\int_{\mathbb{R}}e^{i\lambda t}P\left(d\lambda\right)\label{eq:st4}
\end{equation}
where the integral on the RHS in (\ref{eq:st4}) is the limit of finite
sums of
\begin{equation}
\sum_{k}e^{i\lambda_{k}t}P\left(E_{k}\right),\;t\in\mathbb{R};\label{eq:st5}
\end{equation}
$E_{i}\cap E_{j}=\emptyset\:\left(i\neq j\right)$, $\bigcup_{k}E_{k}=\mathbb{R}$.

$\text{\ensuremath{\left(\ref{enu:st2}\right)} \ensuremath{\Rightarrow} \ensuremath{\left(\ref{enu:st3}\right)}:}$
Given $\{\mathcal{U}\left(t\right)\}_{t\in\mathbb{R}}$, set 
\[
dom\left(H\right)=\left\{ f\in\mathscr{H},\;\mbox{s.t.}\;\lim_{t\rightarrow0_{+}}\frac{1}{i\,t}\left(\mathcal{U}\left(t\right)f-f\right)\;\mbox{exists}\right\} 
\]
and
\begin{equation}
iHf=\lim_{t\rightarrow0_{+}}\frac{\mathcal{U}\left(t\right)f-f}{t},\quad f\in dom\left(H\right),\label{eq:st6}
\end{equation}
then $H^{*}=H$. 

$\text{\ensuremath{\left(\ref{enu:st3}\right)} \ensuremath{\Rightarrow} \ensuremath{\left(\ref{enu:st1}\right)}:}$
Given a selfadjoint operator $H$ with dense domain in $\mathscr{H}$;
then by the spectral theorem (\secref{normal}) there is a unique
PVM, $P\left(\cdot\right)$ such that
\begin{equation}
H=\int_{\mathbb{R}}\lambda P\left(d\lambda\right);\quad\mbox{and}\label{eq:st7}
\end{equation}
\begin{equation}
dom\left(H\right)=\left\{ f\in\mathscr{H};\;\mbox{s.t.}\;\int_{\mathbb{R}}\lambda^{2}\left\Vert P\left(d\lambda\right)f\right\Vert ^{2}<\infty\right\} .\label{eq:st8}
\end{equation}
\end{thm}
\begin{rem}
We state Stone's theorem already now even though the proof details
will require a number of technical tools to be developed systematically
only in Chapters \ref{chap:sp} and \ref{chap:GNS} below. 
\end{rem}

\begin{rem}
Note that the selfadjointness condition on $H$ in (\ref{enu:st3})
in \thmref{stone} is stronger than merely Hermitian symmetry, i.e.,
the condition
\begin{equation}
\left\langle Hu,v\right\rangle =\left\langle u,Hv\right\rangle \label{eq:stt1}
\end{equation}
for all pairs of vectors $u$ and $v$ $\in dom\left(H\right)$. We
shall discuss this important issue in much detail in Part 4 of the
book, both in connection with the theory, and its applications. The
applications are in physics, statistics, and infinite networks.

Here we limit ourselves to comments and some definitions; a full discussion
will follow in part 4 below.

\textbf{Observations.} Introducing the adjoint operator $H^{*}$,
we note that (\ref{eq:stt1}) is equivalent to 
\begin{equation}
H\subset H^{*},\;\mbox{or}\label{eq:stt2}
\end{equation}
\begin{equation}
\mathscr{G}\left(H\right)\subset\mathscr{G}\left(H^{*}\right),\label{eq:stt3}
\end{equation}
where $\mathscr{G}$ denotes the graph of the respective operators
and where (\ref{eq:stt2}) \& (\ref{eq:stt3}) mean that $dom\left(H\right)\subset dom\left(H^{*}\right)$,
and $Hu=H^{*}u$ for $\forall u\in dom\left(H\right)$. 

If (\ref{eq:stt2}) holds, then it may, or may not, have selfadjoint
extensions.

We introduce the two indices $d_{\pm}$ (deficiency-indices)
\begin{equation}
d_{\pm}=\dim\left(H^{*}\pm i\,I\right).\label{eq:stt4}
\end{equation}
The following will be proved in part 4:\end{rem}
\begin{thm}
(i) Suppose $H\subset H^{*}$, then $H$ has selfadjoint extensions
iff $d_{+}=d_{-}$. 

(ii) If $H$ has selfadjoint extensions, say $K$ (i.e., $K^{*}=K$,)
so $H\subset K$, then it follows that
\begin{equation}
H\subset K\subset H^{*}.\label{eq:stt5}
\end{equation}
So, if there are selfadjoint extensions, they lie between $H$ and
$H^{*}$.\end{thm}
\begin{defn}
If $H\subset H^{*}$, and if the closure $\overline{H}=H^{**}$ is
selfadjoint, we say that $H$ is \emph{essentially selfadjoint}.
\end{defn}
\end{subappendices}

\chapter{The Spectral Theorem\label{chap:sp}}
\begin{quotation}
As far as the laws of mathematics refer to reality, they are not certain,
and as far as they are certain, they do not refer to reality. 

--- Albert Einstein\sindex[nam]{Einstein, A., (1879-1955)}\vspace{1em}\textbf{}\\
\textbf{A large part of mathematics which becomes useful developed
with absolutely no desire to be useful, and in a situation where nobody
could possibly know in what area it would become useful; and there
were no general indications that it ever would be so.} By and large
it is uniformly true in mathematics that there is a time lapse between
a mathematical discovery and the moment when it is useful; and that
this lapse of time can be anything from 30 to 100 years, in some cases
even more; and that the whole system seems to function without any
direction, without any reference to usefulness, and without any desire
to do things which are useful. 

--- John von Neumann\sindex[nam]{von Neumann, J., (1903-1957)}\vspace{1em}\\
``The spectral theorem together with the multiplicity theory is one
of the pearls of mathematics.''

--- M. Reed and B. Simon \cite{RS75}\vspace{2em}
\end{quotation}
Most Functional Analysis books, when covering the Spectral Theorem,
stress the bounded case. Because of dictates from applications (especially
quantum physics), below we stress questions directly related to key-issues
for unbounded linear operators. These themes will be taken up again
in Chapters \ref{chap:ext} and \ref{chap:gLap}. In a number of applications,
some operator from physics may only be \textquotedblleft formally
selfadjoint\textquotedblright{} also called Hermitian; and in such
cases, one asks for selfadjoint extension\index{extension!selfadjoint-}s
(if any), \chapref{ext}. \chapref{gLap} is a particular case in
point, arising in the study of infinite graphs.\index{Spectral Theorem}\index{Theorem!Spectral-}

\section{\label{sec:qmview}An Overview}

\index{von Neumann, J.}

von Neumann's spectral theorem (see \cite{Sto90, Yos95, Ne69, RS75, DS88b})
states that an operator $A$ acting in a Hilbert space $\mathscr{H}$
is normal if and only if there exits a projection-valued measure on
$\mathbb{C}$ so that 
\begin{equation}
A=\int_{sp(A)}zP_{A}(dz)\label{eq:spA}
\end{equation}
i.e., $A$ is represented as an integral against the projection-valued
measure $P_{A}$ over its spectrum. \index{quantum mechanics!observable}
\index{integral!Borel-}

In quantum mechanics, an \emph{observable} is represented by a selfadjoint
operator. Functions of observables are again observables. This is
reflected in the spectral theorem as the functional calculus, where
we may define \index{observable} \index{functional calculus} 
\begin{equation}
\varphi(A)=\int_{sp(A)}\varphi(z)P_{A}(dz)\label{eq:spA1}
\end{equation}
using the spectral representation of $A$. \index{representation!spectral}
\index{selfadjoint operator}

When $P$ is a selfadjoint projection, $\left\langle f,Pf\right\rangle _{\mathscr{H}}=\left\Vert Pf\right\Vert _{\mathscr{H}}^{2}$
is a real number and it represents the expected value of the observable
$P$ prepared in the state $f$, unit vector in $\mathscr{H}$. Hence,
in view of (\ref{eq:spA1}), $\left\Vert P_{A}\left(\cdot\right)f\right\Vert _{\mathscr{H}}^{2}$
is a Borel probability measure on $sp\left(A\right)$, and 
\begin{equation}
\left\langle f,\varphi\left(A\right)f\right\rangle _{\mathscr{H}}=\int_{sp\left(A\right)}\varphi\left(z\right)\left\Vert P\left(dz\right)f\right\Vert _{\mathscr{H}}^{2}\label{eq:Aexp}
\end{equation}
is the expected value of the observable $\varphi\left(A\right)$.
\begin{rem}
\label{rem:expA}Let $\varphi:\mathbb{R}\rightarrow\mathbb{R}$ be
measurable and let $A=A^{*}$ be given; then, for every $f\in\mathscr{H}\backslash\left\{ 0\right\} $,
set $d\mu_{f}^{\left(A\right)}\left(\lambda\right):=\left\Vert P_{A}\left(d\lambda\right)f\right\Vert ^{2}\in\mathcal{M}_{+}\left(\mathbb{R}\right)$
(the finite positive Borel measures on $\mathbb{R}$.) Then the transformation
formula (\ref{eq:Aexp}) takes the following equivalent form:\index{Borel measure}
\begin{eqnarray}
d\mu_{f}^{\left(\varphi\left(A\right)\right)} & = & d\mu_{f}^{\left(A\right)}\circ\varphi^{-1},\;\mbox{i.e.,}\label{eq:Aexp1}\\
d\mu_{f}^{\left(\varphi\left(A\right)\right)}\left(\triangle\right) & = & d\mu_{f}^{\left(A\right)}\left(\varphi^{-1}\left(\triangle\right)\right),\;\forall\triangle\in\mathcal{B}\left(\mathbb{R}\right),\label{eq:Aexp2}
\end{eqnarray}
where $\varphi^{-1}\left(\triangle\right)=\left\{ x\::\:\varphi\left(x\right)\in\triangle\right\} $.\end{rem}
\begin{cor}
Let $A=A^{*}$, and $f\in\mathscr{H}\backslash\left\{ 0\right\} $
be given, and let $\mu_{f}$ and $\varphi$ be as in  \remref{expA},
then $\varphi\left(A\right)^{*}=\varphi\left(A\right)$, and 
\[
f\in dom\left(\varphi\left(A\right)\right)\Longleftrightarrow\varphi\in L^{2}\left(\mathbb{R},\mu_{f}\right),
\]
where ``$dom$'' is short for ``domain.''\end{cor}
\begin{proof}
This is immediate from (\ref{eq:Aexp})-(\ref{eq:Aexp2}). Indeed,
setting 
\[
d\mu_{f}\left(\lambda\right):=\left\Vert P\left(d\lambda\right)f\right\Vert _{\mathscr{H}}^{2},
\]
we get 
\[
\int_{\mathbb{R}}\left|\varphi\left(\lambda\right)\right|^{2}d\mu_{f}\left(\lambda\right)=\left\Vert \varphi\left(A\right)f\right\Vert _{\mathscr{H}}^{2}.
\]
\end{proof}
\begin{rem}
The standard diagonalization of Hermitian\index{operators!Hermitian}
matrices in linear algebra is a special case of the spectral theorem.
Recall that any Hermitian matrix $A$ can be decomposed as $A=\sum_{k}\lambda_{k}P_{k}$,
where $\lambda_{k}'s$ are the eigenvalues of $A$ and $P_{k}'s$
are the selfadjoint projections onto the eigenspaces associated with
$\lambda_{k}'s$. The projection-valued measure in this case can be
written as $P(E)=\sum_{\lambda_{k}\in E}P_{k}$, for all $E\in\mathcal{B}\left(\mathbb{R}\right)$;
i.e., the counting measure supported on $\lambda_{k}'s$.\index{Spectral Theorem}\index{diagonalization}\index{matrix!diagonal-}\index{Theorem!Spectral-}\index{eigenvalue}
\end{rem}

Quantum mechanics is stated using an abstract Hilbert space as the
state space. In practice, one has the freedom to choose exactly which
Hilbert space to use for a particular problem. Physical measurements
remain unchanged when choosing different realizations of a Hilbert
space. The concept needed here is unitary equivalence\index{unitary equivalence}.
\index{space!Hilbert-}\index{space!state-}
\begin{defn}
Let $A:\mathscr{H}_{1}\rightarrow\mathscr{H}_{1}$ and $B:\mathscr{H}_{2}\rightarrow\mathscr{H}_{2}$
be operators. $A$ is said to be unitarily equivalent to $B$ if there
exists a unitary operator $U:\mathscr{H}_{1}\rightarrow\mathscr{H}_{2}$
such that $B=UAU^{*}$.
\end{defn}
Suppose $U:\mathscr{H}_{1}\rightarrow\mathscr{H}_{2}$ is a unitary
operator, $P:\mathscr{H}_{1}\rightarrow\mathscr{H}_{1}$ is a selfadjoint
projection. Then $UPU^{*}:\mathscr{H}_{2}\rightarrow\mathscr{H}_{2}$
is a selfadjoint projection on $\mathscr{H}_{2}$. In fact, 
\[
\left(UPU^{*}\right)\left(UPU^{*}\right)=UPU^{*}
\]
where we used $UU^{*}=U^{*}U=I$, since $U$ is unitary. Let $|f_{1}\rangle$
be a state in $\mathscr{H}_{1}$ and $|f_{2}\rangle=|Uf_{1}\rangle$
be the corresponding state in $\mathscr{H}_{2}$. Then
\[
\left\langle f_{2},UPU^{*}f_{2}\right\rangle _{\mathscr{H}_{2}}=\left\langle U^{*}f_{2},PU^{*}f_{2}\right\rangle _{\mathscr{H}_{1}}=\left\langle f_{1},Pf_{1}\right\rangle _{\mathscr{H}_{1}}
\]
i.e., the observable\index{observable} $P$ has the same expectation
value. Since every selfadjoint operator is, by the spectral theorem,
decomposed into selfadjoint projections, it follows that the expectation
value of any observable remains unchanged under unitary transformations.

We will also consider family of selfadjoint operators. Heisenberg's
commutation relation $PQ-QP=-i\,I$, $i=\sqrt{-1}$, is an important
example of two non-commuting selfadjoint operators.

\begin{example}
The classical Fourier transform $\mathcal{F}:L^{2}(\mathbb{R})\rightarrow L^{2}(\mathbb{R})$
is unitary, so $\mathcal{\mathcal{F}}^{*}\mathcal{F}=\mathcal{F}\mathcal{F}^{*}=I_{L^{2}\left(\mathbb{R}\right)}$,
and in particular, the Parseval identity \index{Parseval identity}
\[
\left\Vert \mathcal{F}f\right\Vert _{L^{2}\left(\mathbb{R}\right)}^{2}=\left\Vert f\right\Vert _{L^{2}\left(\mathbb{R}\right)}^{2}
\]
holds for all $f\in L^{2}\left(\mathbb{R}\right)$.
\end{example}

\begin{example}
Let $Q$ and $P$ be the position and momentum operators in quantum
mechanics. That is, $Q=M_{x}=$ multiplication by $x$, and $P=-id/dx$
both defined on the Schwartz space $\mathcal{S}\left(\mathbb{R}\right)$--space
of rapidly decreasing functions on $\mathbb{R}$, which is dense in
the Hilbert space $L^{2}\left(\mathbb{R}\right)$. On $\mathcal{S}\left(\mathbb{R}\right)$,
the operators $P$ and $Q$ satisfy the canonical commutation relation:
$PQ-QP=-i\,I_{L^{2}\left(\mathbb{R}\right)}$. 

\index{quantum mechanics!momentum operator}

\index{quantum mechanics!position operator}

\index{Heisenberg, W.K.!commutation relation}

\index{commutation relations}
\end{example}
\index{operators!momentum-}
\begin{example}
Denote $\mathcal{F}$ the Fourier transform on $L^{2}\left(\mathbb{R}\right)$
as before. Specifically, setting 
\begin{alignat*}{1}
\left(\mathcal{F}\varphi\right)\left(x\right) & =\widehat{\varphi}\left(x\right)=\frac{1}{\sqrt{2\pi}}\int_{\mathbb{R}}\varphi\left(\xi\right)e^{-i\xi x}d\xi,\;\mbox{and}\\
\left(\mathcal{F}^{*}\psi\right)\left(\xi\right) & =\psi^{\vee}\left(\xi\right)=\frac{1}{\sqrt{2\pi}}\int_{\mathbb{R}}\psi\left(x\right)e^{i\xi x}dx,\;\xi\in\mathbb{R}.
\end{alignat*}
Note that $\mathcal{F}$ is an automorphism in $\mathcal{S}\left(\mathbb{R}\right)$,
continuous with respect to the standard l.c. topology. Moreover, 
\[
\left(\mathcal{F}^{*}Q\mathcal{F}\varphi\right)\left(\xi\right)=\mathcal{F}^{*}\left(x\widehat{\varphi}\left(x\right)\right)=\frac{1}{i}\frac{d}{d\xi}\varphi\left(\xi\right),\;\forall\varphi\in\mathcal{S}.
\]
Therefore, 
\begin{equation}
P=\mathcal{F}^{*}Q\mathcal{F}\label{eq:diagP}
\end{equation}
and so $P$ and $Q$ are unitarily equivalent. \index{transform!Fourier} 
\end{example}

A multiplication operator version of the spectral theorem is also
available. It works especially well in physics. It says that $A$
is a normal operator in $\mathscr{H}$ if and only if $A$ is unitarily
equivalent to the operator of multiplication by a measurable function
$f$ on $L^{2}(X,\mu)$, where $X$ is locally compact and Hausdorff.
The two versions are related via a measure transformation. \index{normal operator}
\begin{example}
Eq. (\ref{eq:diagP}) says that $P$ is diagonalized by Fourier transform
in the following sense. 

Let $\psi$ be any Borel function on $\mathbb{R}$, and set $M_{\psi}=$
multiplication by $\psi\left(x\right)$ in $L^{2}\left(\mathbb{R}\right)$,
with 
\begin{eqnarray}
dom\left(M_{\psi}\right) & = & \left\{ f\:\big|\:f,\psi f\in L^{2}\left(\mathbb{R}\right)\right\} \label{eq:domM}\\
 & = & \left\{ f\:\big|\:\int_{-\infty}^{\infty}\left(1+\left|\psi\left(x\right)\right|^{2}\right)\left|f\left(x\right)\right|^{2}dx<\infty\right\} ;\nonumber 
\end{eqnarray}
then we define, via eq. (\ref{eq:diagP}), 
\[
\psi\left(P\right):=\mathcal{F}^{*}\psi\left(Q\right)\mathcal{F}.
\]
In particular, given any $\triangle\in\mathcal{B}\left(\mathbb{R}\right)$,
let $\psi=\chi_{\triangle}=$ characteristic function, then 
\[
E(\triangle)=\mathcal{F}^{*}M_{\chi_{\triangle}}\mathcal{F}.
\]
One checks directly that $E\left(\triangle\right)^{2}=E\left(\triangle\right)=E\left(\triangle\right)^{*}$,
so $E\left(\triangle\right)$ is a selfadjoint projection. Indeed,
$E\left(\cdot\right)$ is a convolution operator, where 
\[
\left(E\left(\triangle\right)f\right)\left(x\right)=\int_{a}^{b}e^{i\xi x}\widehat{f}\left(\xi\right)d\xi=f\ast\left(\chi_{\left[a,b\right]}\right)^{\wedge}\left(x\right),\;\forall f\in L^{2}\left(\mathbb{R}\right).
\]
Thus, 
\[
E\left(\triangle\right)\left(L^{2}\left(\mathbb{R}\right)\right)=\left\{ f\in L^{2}\left(\mathbb{R}\right)\:\big|\:\mbox{supp}\left(\widehat{f}\right)\subset\triangle\right\} ;
\]
i.e., the space of ``band-limited'' functions, with the ``pass-band''
being $\triangle$. 
\end{example}

\index{band-limited}\index{operators!convolution-}\index{operators!convolution-}

\index{convolution}
\begin{example}
Below, it helps to denote the Fourier transformed space (or frequency
space) by $L^{2}(\widehat{\mathbb{R}})$. Fix any $f\in L^{2}\left(\mathbb{R}\right)$,
$\triangle\in\mathcal{B}\left(\mathbb{R}\right)$, then 
\begin{eqnarray*}
\mu_{f}\left(\triangle\right):=\left\Vert E\left(\triangle\right)f\right\Vert _{L^{2}\left(\mathbb{R}\right)}^{2} & = & \left\langle f,E\left(\triangle\right)f\right\rangle _{L^{2}\left(\mathbb{R}\right)}\\
 & = & \left\langle \mathcal{F}f,M_{\chi_{\triangle}}\mathcal{F}f\right\rangle _{L^{2}\left(\widehat{\mathbb{R}}\right)}\\
 & = & \int_{\triangle}\left|\widehat{f}\left(x\right)\right|^{2}dx
\end{eqnarray*}
which is a Borel measure on $\mathbb{R}$, such that 
\[
\mu_{f}\left(\mathbb{R}\right)=\int_{-\infty}^{\infty}\left|\widehat{f}\left(x\right)\right|^{2}dx=\int_{-\infty}^{\infty}\left|f\left(x\right)\right|^{2}dx=\left\Vert f\right\Vert _{L^{2}\left(\mathbb{R}\right)}^{2}
\]
by the Parseval identity.

Now, let $\psi\left(x\right)=x$ be the identity function, and note
that it is approximated pointwisely by simple functions of the form
$\sum_{\text{finite}}c_{i}\chi_{\triangle_{i}}$, where $\triangle_{i}\in\mathcal{B}\left(\mathbb{R}\right)$,
and $\triangle_{i}$'s are mutually disjoint, i.e., $x=\lim_{n\rightarrow\infty}\sum_{i=1}^{n}c_{i}\chi_{\triangle_{i}}$. 

Fix $f\in dom\left(M_{x}\right)$, see (\ref{eq:domM}), it follows
from Lebesgue dominated convergence theorem\index{Lebesgue dominated convergence theorem},
that 
\begin{eqnarray*}
\left\langle f,Pf\right\rangle _{L^{2}\left(\mathbb{R}\right)} & = & \int_{-\infty}^{\infty}\overline{\widehat{f}\left(x\right)}x\widehat{f}\left(x\right)dx\\
 & = & \lim_{n\rightarrow\infty}\int_{-\infty}^{\infty}\overline{\widehat{f}\left(x\right)}\left(\sum_{i=1}^{n}c_{i}\chi_{\triangle_{i}}\left(x\right)\right)\widehat{f}\left(x\right)dx\\
 & = & \lim_{n\rightarrow\infty}\sum_{i=1}^{n}c_{i}\left\langle f,E\left(\triangle_{i}\right)f\right\rangle _{L^{2}\left(\mathbb{R}\right)}\\
 & = & \lim_{n\rightarrow\infty}\sum_{i=1}^{n}c_{i}\left\Vert E\left(\triangle_{i}\right)f\right\Vert _{L^{2}\left(\mathbb{R}\right)}^{2}\\
 & = & \int_{-\infty}^{\infty}x\left\Vert E\left(dx\right)f\right\Vert _{L^{2}\left(\mathbb{R}\right)}^{2}.
\end{eqnarray*}
The last step above yields the projection-valued measure (PVM) version
of the spectral theorem for $P$, where we write 
\begin{equation}
P=\int_{-\infty}^{\infty}x\:dE\left(x\right).\label{eq:pvm}
\end{equation}
Consequently, we get two versions of the spectral theorem for $P=\frac{1}{i}\frac{d}{dx}\Big|_{\mathcal{S}\left(\mathbb{R}\right)}$:\index{Theorem!Spectral-}
\begin{enumerate}
\item Multiplication operator version, i.e., $P\simeq M_{x}=$ multiplication
by $x$ in $L^{2}(\widehat{\mathbb{R}})$; and
\item PVM version, as in (\ref{eq:pvm}). 
\end{enumerate}
\end{example}
This example illustrates the main ideas of the spectral theorem of
a single selfadjoint operator in Hilbert space. We will develop the
general theory in this chapter, and construct both versions of the
spectral decomposition.

\begin{example}
\label{exa:her1}Applying the Gram-Schmidt process to all polynomials
against the measure $d\mu=e^{-x^{2}/2}dx$, one gets orthogonal polynomials
$P_{n}$ in $L^{2}\left(\mu\right)$. These are called the \emph{\uline{Hermite
polynomials}}\emph{, }and the associated Hermite functions are given
by 
\[
h_{n}:=e^{-x^{2}/2}P_{n}=e^{-x^{2}}\left(\frac{d}{dx}\right)^{n}e^{x^{2}/2}.
\]
The Hermite functions (after normalization) forms an ONB in $L^{2}\left(\mathbb{R}\right)$,
which transforms $P$ and $Q$ to Heisenberg's infinite matrices in
(\ref{eq:HP})-(\ref{eq:HQ}). \index{Gram-Schmidt orthogonalization}

\index{Hermite polynomials}
\end{example}
\index{orthogonal polynomial!Hermite}

\index{Hermite function}

\begin{example}[The harmonic oscillator Hamiltonian]
\begin{flushleft}
\label{exa:her2} Let $P,Q$ be as in the previous example. We consider
the quantum Hamiltonian \index{harmonic}\index{harmonic oscillator}
\[
H:=\frac{1}{2}(Q^{2}+P^{2}-1).
\]
It can be shown that
\[
Hh_{n}=nh_{n}
\]
or equivalently,
\[
(P^{2}+Q^{2})h_{n}=(2n+1)h_{n}
\]
$n=0,1,2,\ldots$. That is, $H$ is diagonalized by the Hermite functions. 
\par\end{flushleft}

\begin{flushleft}
$H$ is called the energy operator in quantum mechanics. This explains
mathematically why the energy levels are discrete (in quanta), being
a multiple the Plank's constant $\hbar$. \index{quantum mechanics!energy operator}
\par\end{flushleft}
\end{example}

\begin{figure}
\subfloat[Harmonic oscillator $P^{2}+Q^{2}$ (bound-states). ]{%
\begin{minipage}[t]{0.4\textwidth}%
\protect\begin{center}
\protect\includegraphics[scale=0.45]{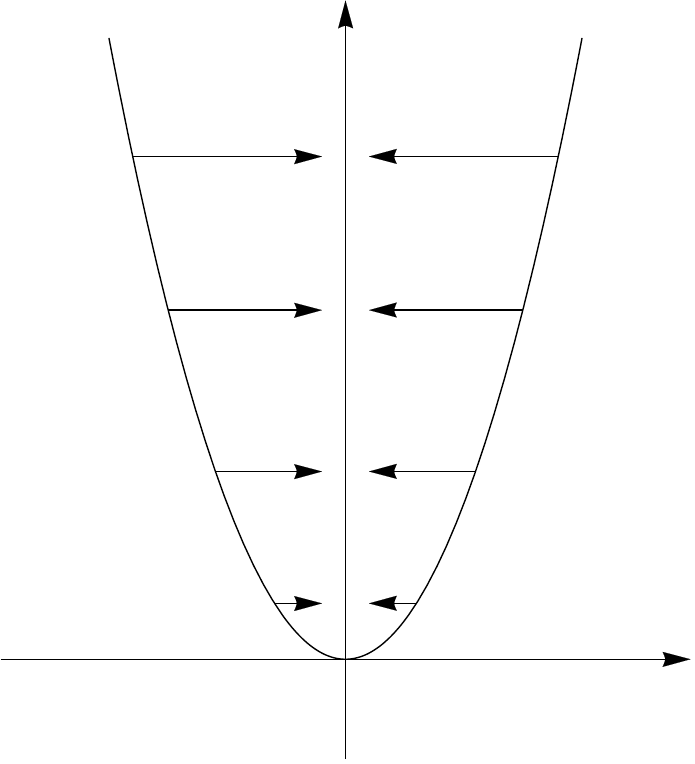}\protect
\par\end{center}%
\end{minipage}}\hfill{}\subfloat[\label{fig:rep}Repulsive potential $P^{2}-Q^{2}$. This operator
has purely continuous spectrum.]{%
\begin{minipage}[t]{0.4\textwidth}%
\protect\begin{center}
\protect\includegraphics[scale=0.45]{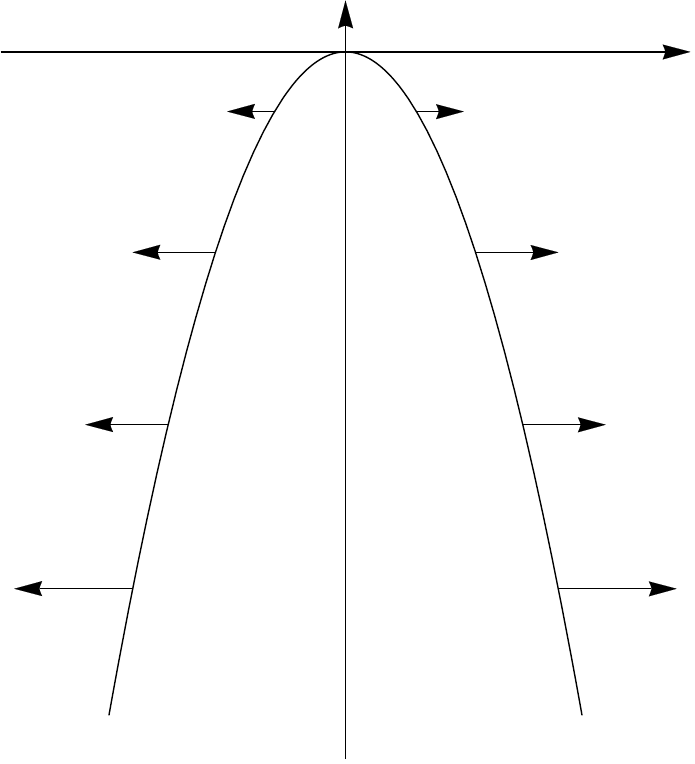}\protect
\par\end{center}%
\end{minipage}

}

\centering{}\protect\caption{\label{fig:sp}Illustration of forces: attractive vs repulsive. The
case of ``only bound states'' (a), vs continuous Lebesgue spectrum
(b).}
\end{figure}

\begin{example}[Purely discrete spectrum v.s. purely continuous spectrum]
\label{exa:PQ}The two operators $P^{2}+Q^{2}$ and $P^{2}-Q^{2}$
acting in $L^{2}\left(\mathbb{R}\right)$; see \figref{sp}.

\begin{flushleft}
\index{harmonic oscillator}
\par\end{flushleft}

\index{spectrum!continuous}
\end{example}

\begin{rem}
Note that both of the two operators $H_{\pm}:=P^{2}\pm Q^{2}$ in
 \exaref{PQ} are essentially selfadjoint as operators in $L^{2}\left(\mathbb{R}\right)$
(see \cite{RS75, Ne69, vN32a, DS88b}), and with common dense domain
equal to the Schwartz space. The potential in $H_{-}$ is repulsive,
see \figref{sp} (b). 

By comparison, the operator $H_{4}:=P^{2}-Q^{4}$ is not essentially
selfadjoint. (It can be shown that it has deficiency indices $\left(2,2\right)$.)
The following argument from physics is illuminating: For $E\in\mathbb{R}_{+}$,
consider a classical particle $x\left(t\right)$ on the energy surface
\[
S_{E}:=\left\{ x\left(t\right)\::\:\left(x'\left(t\right)\right)^{2}-\left(x\left(t\right)\right)^{4}=E\right\} .
\]
The travel time to $\pm\infty$ is finite; in fact, it is 
\[
t_{\infty}=\int_{0}^{\infty}\frac{dx}{\sqrt{E+x^{4}}}<\infty.
\]

There is a principle from quantum mechanics which implies that the
quantum mechanical particle must be assigned conditions at $\pm\infty$,
which translates into non-zero deficiency indices. (A direct computation,
which we omit, yields indices $\left(2,2\right)$.)
\end{rem}
\begin{figure}
\includegraphics[width=0.6\columnwidth]{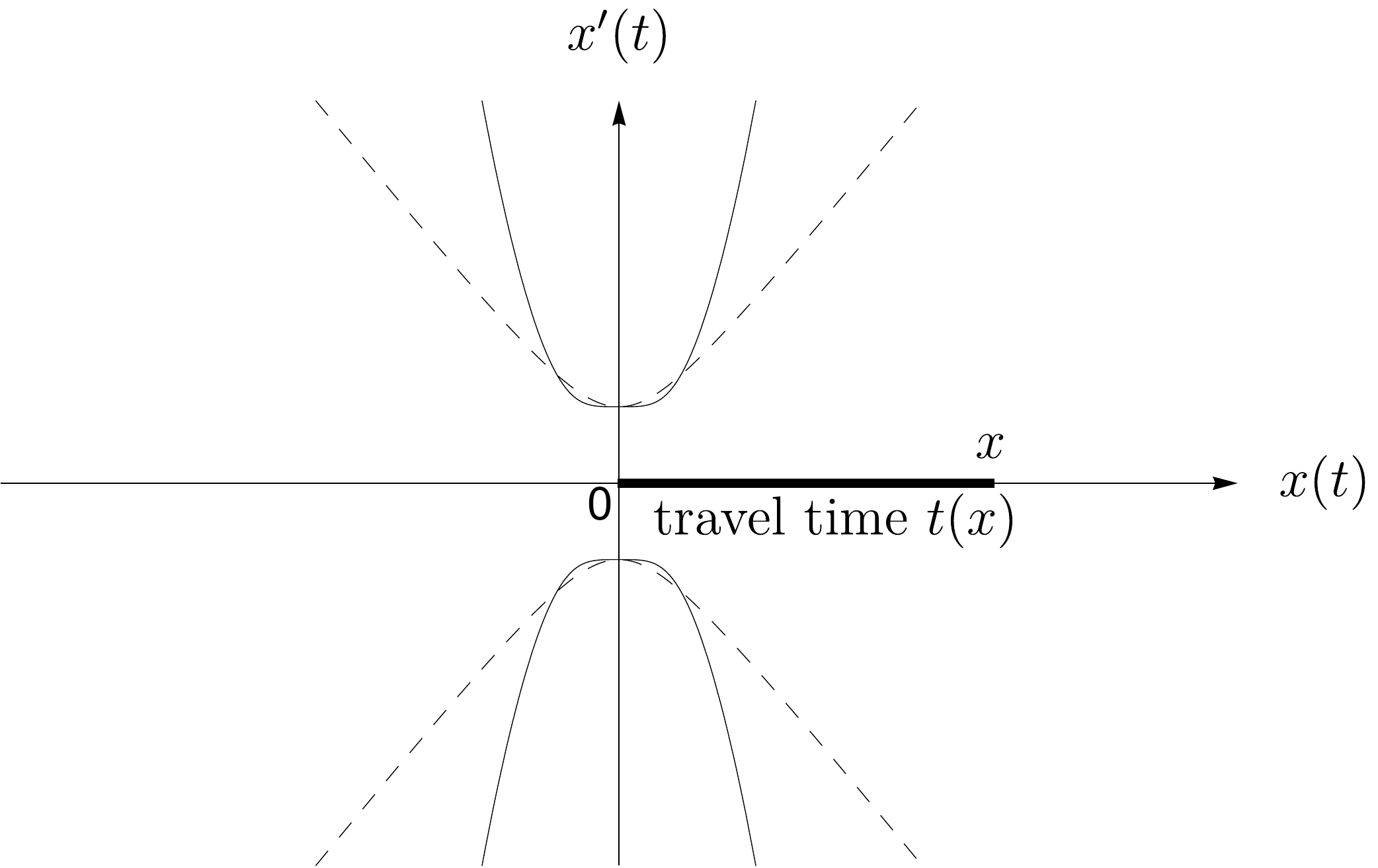}

A phase-space plot, $\left(x,x'\right)$. Travel time $t\left(x\right)={\displaystyle \int_{0}^{x}\frac{ds}{\sqrt{E+s^{4}}}}$.

\protect\caption{The energy surface $S_{E}$ for the quantum mechanical $H_{4}=P^{2}-Q^{4}$,
with $P\rightsquigarrow x'\left(t\right)$. }
\end{figure}

\begin{figure}
\includegraphics[width=0.5\columnwidth]{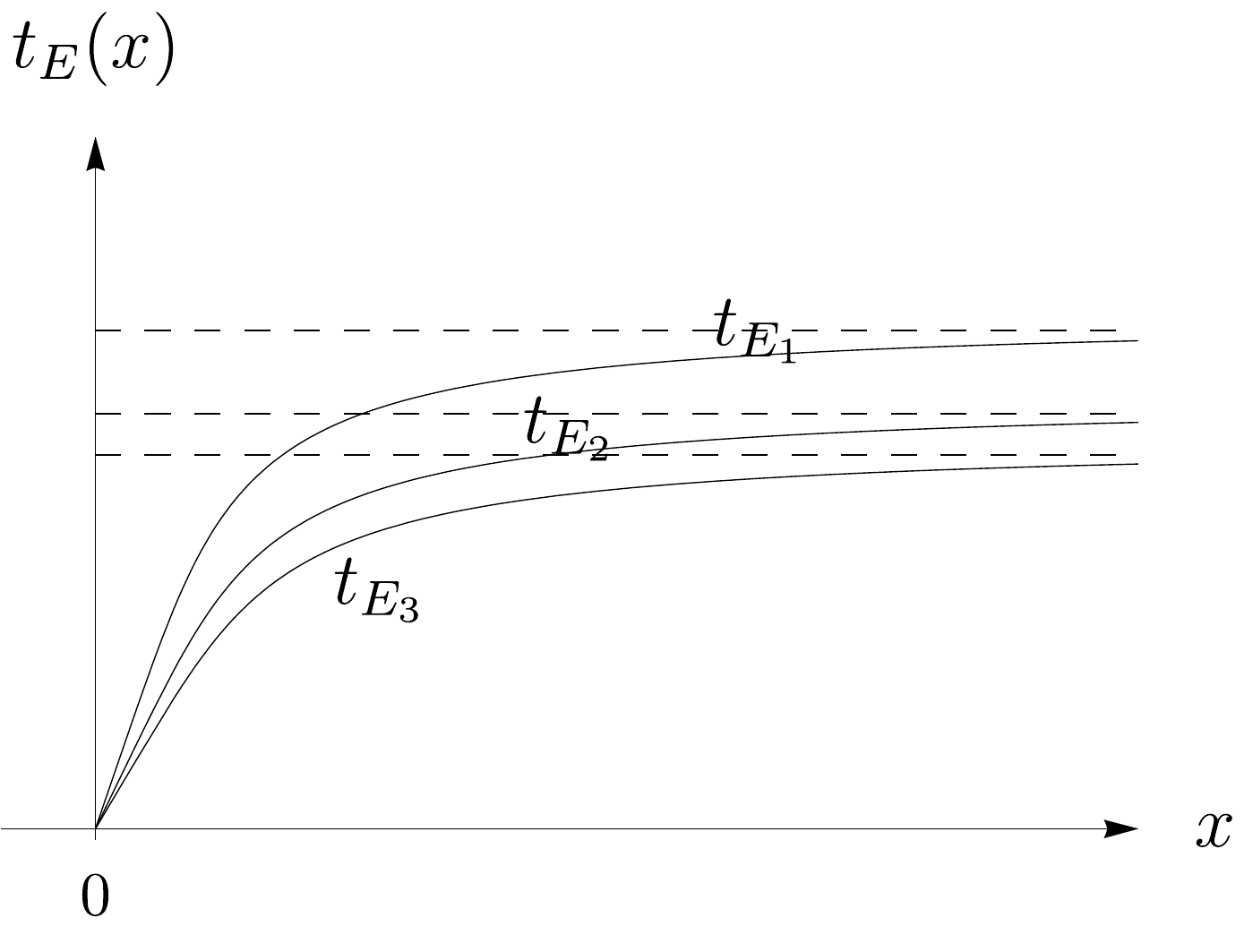}

Travel time $t_{E}\left(x\right)={\displaystyle \int_{0}^{x}\frac{ds}{\sqrt{E+s^{4}}}}$.

\protect\caption{Fix $E_{1}<E_{2}<E_{3}$, then $0<t_{E_{3}}\left(\infty\right)<t_{E_{2}}\left(\infty\right)<t_{E_{1}}\left(\infty\right)$.}
\end{figure}

In the following sections, we present some main ideas of the spectral
theorem for single normal operators acting in Hilbert space. Since
every normal operator $N$ can be written as $N=T_{1}+iT_{2}$, where
$T_{1}$ and $T_{2}$ are strongly commuting and selfadjoint, the
presentation will be focused on selfadjoint operators. \index{operators!selfadjoint}
\index{normal operator}\index{selfadjoint operator}

\section{Multiplication Operator Version}

Together, the results below serve to give a spectral representation
(by multiplication operators) for the most general case: an arbitrary
given selfadjoint (or normal) operator with dense domain in a Hilbert
space. It applies both to the bounded, and unbounded cases; and it
even applies to arbitrary families of strongly commuting selfadjoint
operators. Caution: There is a number of subtle points in such representations.
Since we aim for realizations up to unitary equivalence, care must
be exercised in treating \textquotedblleft multiplicity\textquotedblright{}
for the most general spectral types. What we present below may be
thought as a modern version of what is often called the \emph{Hahn-Hellinger
theory of spectral multiplicity}. \index{Theorem!Hahn-Hellinger-}

This version of the spectral theory states that every selfadjoint
operator $A$ is unitarily equivalent to the operator of multiplication
by a measurable function on some $L^{2}$-space. \index{operators!multiplication}\index{Spectral Theorem}
\begin{thm}
\label{thm:spA}Let $A$ be a linear operator acting in the Hilbert
space $\mathscr{H}$, then $A=A^{*}$ iff there exists a measure space
$\left(X,\mu\right)$ and a unitary operator $U:L^{2}\left(X,\mu\right)\rightarrow\mathscr{H}$
such that 
\begin{equation}
M_{\varphi}=U^{*}AU;\label{eq:mul3}
\end{equation}
where $X$ is locally compact and Hausdorff, $\varphi$ is a real-valued
$\mu$-measurable function, and
\begin{eqnarray}
M_{\varphi}f & := & \varphi f,\forall f\in dom\left(M_{\varphi}\right),\;\mbox{where}\label{eq:mul1}\\
dom\left(M_{\varphi}\right) & := & \left\{ h\in L^{2}\left(X,\mu\right)\::\:\varphi h\in L^{2}\left(X,\mu\right)\right\} .\label{eq:mul2}
\end{eqnarray}
Hence, the following diagram commutes.
\[
\xymatrix{\mathscr{H}\ar[r]^{A} & \mathscr{H}\\
L^{2}\left(X,\mu\right)\ar[u]^{U}\ar[r]^{M_{\varphi}} & L^{2}\left(X,\mu\right)\ar[u]_{U}
}
\]
If $A\in\mathscr{B}\left(\mathscr{H}\right)$, then $\varphi\in L^{\infty}\left(X,\mu\right)$
and $dom\left(M_{\varphi}\right)=\mathscr{H}$.
\end{thm}
We postpone the detailed proof till \secref{proofsp} below. 
\begin{xca}[Multiplication operators, continued]
\myexercise{Multiplication operators, continued}\label{exer:mult1}Prove
that $M_{\varphi}$ in (\ref{eq:mul1})-(\ref{eq:mul2}) is selfadjoint.
\end{xca}

\begin{xca}[Continuous spectrum]
\myexercise{Continuous spectrum}\label{exer:mult2}Let $M_{t}:L^{2}[0,1]\rightarrow L^{2}[0,1]$,
$f\left(t\right)\longmapsto tf\left(t\right)$. Show that $M_{t}$
has no eigenvalues in $L^{2}\left[0,1\right]$. \index{spectrum!continuous-}\index{eigenvalue}
\end{xca}

Before giving a proof of \thmref{spA}, we show below that one can
go one step further and get that $A$ is unitarily equivalent to the
operator of multiplication by the independent variable on some $L^{2}$-space.
This is done by a transformation of the measure $\mu$ in (\ref{eq:mul2}).

\subsection{Transformation of Measures}
\begin{defn}
Let $\varphi:X\rightarrow Y$ be a measurable function, $\mathcal{T}_{X}$
and $\mathcal{T}_{Y}$ be the respective sigma-algebra\index{sigma-algebra}s.
Fix a measure $\mu$ on $\mathcal{T}_{X}$, the measure 
\begin{equation}
\mu_{\varphi}:=\mu\circ\varphi^{-1}\label{eq:mt0}
\end{equation}
defined on $\mathcal{T}_{Y}$ is called the transformation of $\mu$
under $\varphi$. Note that 
\end{defn}
\begin{equation}
\chi_{E}\circ\varphi(x)=\chi_{\varphi^{-1}(E)}\left(x\right)\label{eq:mt}
\end{equation}
for all $E\in\mathcal{T}_{Y}$ and $x\in X$. 
\begin{lem}
\label{lem:tm}For all $\mathcal{T}$-measurable function $f$, 
\begin{equation}
\int_{X}f\circ\varphi d\mu=\int_{Y}f\,d\left(\mu\circ\varphi^{-1}\right).\label{eq:mes1}
\end{equation}
(This is a generalization of the substitution formula in calculus.)\end{lem}
\begin{proof}
For any simple function $s=\sum c_{i}\chi_{E_{i}}$, $E_{i}\in\mathcal{T}_{Y}$,
it follows from (\ref{eq:mt}) that 
\begin{eqnarray*}
\int_{X}s\circ\varphi d\mu & = & \sum c_{i}\int_{X}\chi_{E_{i}}\circ\varphi d\mu\\
 & = & \sum c_{i}\int_{X}\chi_{\varphi^{-1}\left(E_{i}\right)}d\mu\\
 & = & \sum c_{i}\mu\left(\varphi^{-1}\left(E_{i}\right)\right)\\
 & = & \int_{Y}s\:d\left(\mu\circ\varphi^{-1}\right).
\end{eqnarray*}
Note all the summations in the above calculation are finite.

Since any measurable function $f:X\rightarrow Y$ is approximated
pointwisely by simple functions, eq. (\ref{eq:mes1}) follows. \end{proof}
\begin{rem}
If $\varphi$ is nasty, even if $\mu$ is a nice measure (say the
Lebesgue measure), the transformation measure $\mu\circ\varphi^{-1}$
in (\ref{eq:mes1}) can still be nasty, e.g., it could even be singular.
\end{rem}
To simplify the discussion, we consider \uline{bounded} selfadjoint
operators below. 
\begin{cor}
Let $\varphi:X\rightarrow X$ be any measurable function. Then the
operator $Uf:=f\circ\varphi$ in $L^{2}\left(X,\mu\right)$ is isometric
iff $\mu\circ\varphi^{-1}=\mu$. Moreover, $M_{\varphi}U=UM_{t}$.
In particular, $U$ is unitary iff $\varphi$ is invertible.\end{cor}
\begin{proof}
Follows immediately from \lemref{tm}.\end{proof}
\begin{lem}
\label{lem:mt1}In \thmref{spA}, assume $A$ is bounded selfadjoint,
so that $\varphi\in L^{\infty}\left(X,\mu\right)$, and real-valued.
Let $\mu_{\varphi}:=\mu\circ\varphi^{-1}$ (eq. (\ref{eq:mt0})),
supported on the essential range of $\varphi$. Then the operator
$W:L^{2}\left(\mathbb{R},\mu_{\varphi}\right)\longrightarrow L^{2}\left(X,\mu\right)$,
by 
\begin{equation}
\left(Wf\right)\left(x\right)=f\left(\varphi\left(x\right)\right),\;\forall f\in L^{2}\left(\mu_{\varphi}\right)\label{eq:m1}
\end{equation}
is isometric, and 
\begin{equation}
WM_{t}=M_{\varphi}W,\label{eq:m2}
\end{equation}
where $M_{t}:L^{2}\left(Y,\mu_{f}\right)\longrightarrow L^{2}\left(Y,\mu_{f}\right)$,
given by 
\begin{equation}
\left(M_{t}f\right)\left(t\right)=t\,f\left(t\right);\label{eq:m3}
\end{equation}
i.e., multiplication by the identify function.\end{lem}
\begin{proof}
For all $f\in L^{2}\left(Y,\mu_{\varphi}\right)$, we have 
\[
\left\Vert f\right\Vert _{L^{2}\left(Y,\mu_{\varphi}\right)}^{2}=\int_{Y}\left|f\right|^{2}d\mu_{\varphi}=\int_{X}\left|f\circ\varphi\right|^{2}d\mu=\left\Vert Wf\right\Vert _{L^{2}\left(X,\mu\right)}^{2}
\]
so $W$ is isometric. Moreover, 
\begin{eqnarray*}
M_{\varphi}Wf & = & \varphi\left(x\right)f\left(\varphi\left(x\right)\right)\\
WM_{t}f & = & W\left(tg\left(t\right)\right)=\varphi\left(x\right)f\left(\varphi\left(x\right)\right)
\end{eqnarray*}
hence (\ref{eq:m2})-(\ref{eq:m3}) follows.\end{proof}
\begin{cor}
Let $\varphi$ be as in \lemref{mt1}. Assume $\varphi$ is invertible,
we get that $W$ in (\ref{eq:m1})-(\ref{eq:m3}) is unitary. Set
$\mathcal{F}=UW:L^{2}\left(\mathbb{R},\mu_{\varphi}\right)\longrightarrow\mathscr{H}$,
then $\mathcal{F}$ is unitary and 
\begin{equation}
M_{t}=\mathcal{F}^{*}A\mathscr{F};\label{eq:m4}
\end{equation}
i.e., the following diagram commutes.
\end{cor}
\[
\xymatrix{\mathscr{H}\ar[r]^{A} & \mathscr{H}\\
L^{2}(X,\mu)\ar[u]^{U}\ar[r]^{M_{\varphi}} & L^{2}(X,\mu)\ar[u]_{U}\\
L^{2}(\mathbb{R},\mu_{\varphi})\ar[u]^{W}\ar[r]^{M_{t}}\ar@/^{2pc}/@{-->}[luur]^{\mathcal{F}} & L^{2}(\mathbb{R},\mu_{\varphi})\ar[u]_{W}\ar@/_{2pc}/@{-->}[ruul]_{\mathcal{F}}
}
\]

\begin{rem}
\begin{flushleft}
\textbf{~}
\par\end{flushleft}
\begin{enumerate}
\item Eq. (\ref{eq:m4}) is a vast extension of \emph{diagonalizing hermitian
matrices} in linear algebra, or a generalization of Fourier transform.
\index{transform!Fourier}\index{diagonalization}
\item Given a selfadjoint operator $A$ in the Hilbert space $\mathscr{H}$,
what's involved are two algebras: the algebra of measurable functions
on $X$, treated as multiplication operators, and the algebra of operators
generated by $A$ (with identity). The two algebras are $*$-isomorphic.
The Spectral Theorem offers two useful tools:

\begin{enumerate}
\item Representing the algebra generated by $A$ by the algebra of functions.
In this direction, it helps to understand $A$.
\item Representing the algebra of functions by the algebra of operators
generated by $A$. In this direction, it reveals properties of the
function algebra and the underlying space $X$. 
\end{enumerate}
\item Let $\mathfrak{A}$ be the algebra of functions. We say that $\pi$
is a representation of $\mathfrak{A}$ on the Hilbert space $\mathscr{H}$,
denoted by $\pi\in Rep\left(\mathfrak{A},\mathscr{H}\right)$, if
$\pi:\mathfrak{A}\rightarrow\mathscr{B}\left(\mathscr{H}\right)$
is a $*$-homomorphism\index{homomorphism}, i.e., $\pi\left(gh\right)=\pi\left(g\right)\pi\left(h\right)$,
and $\pi\left(\overline{g}\right)=\pi\left(g\right)^{*}$, for all
$g,h\in\mathfrak{A}$. Given $\mathcal{F}$ as in (\ref{eq:m4}),
then 
\begin{equation}
\pi(\psi)=\mathcal{F}M_{\psi}\mathcal{F}^{*}\in Rep\left(\mathfrak{A},\mathscr{H}\right);\label{eq:srA}
\end{equation}
where the LHS in (\ref{eq:srA}) defines the operator 
\begin{equation}
\psi\left(A\right):=\pi\left(\psi\right),\;\psi\in\mathfrak{A}.\label{eq:srA1}
\end{equation}
To see that (\ref{eq:srA1}) is an algebra isomorphism, one checks
that 
\begin{eqnarray*}
\left(\psi_{1}\psi_{2}\right)\left(A\right) & = & \mathcal{F}M_{\psi_{1}\psi_{2}}\mathcal{F}^{*}\\
 & = & \mathcal{F}M_{\psi_{1}}M_{\psi_{2}}\mathcal{F}^{*}\\
 & = & \left(\mathcal{F}M_{\psi_{1}}\mathcal{F}^{*}\right)\left(\mathcal{F}M_{\psi_{2}}\mathcal{F}^{*}\right)\\
 & = & \psi_{1}\left(A\right)\psi_{2}\left(A\right)
\end{eqnarray*}
using the fact that 
\[
M_{\psi_{1}\psi_{2}}=M_{\psi_{1}}M_{\psi_{2}}
\]
i.e., multiplication operators always commute.
\item Eq. (\ref{eq:srA}) is called the \emph{spectral representation} of\emph{
$A$}. In particular, the spectral theorem of $A$ implies the following
substitution rule
\[
\sum c_{k}x^{k}\longmapsto\sum c_{k}A^{k}
\]
is well-defined, and it extends to all bounded measurable functions.
\index{representation!spectral}\index{Theorem!Spectral-}
\end{enumerate}
\end{rem}
Let $\varphi$ be as in \lemref{mt1}. For the more general case when
$\varphi$ is not necessarily invertible, so $W$ in (\ref{eq:m2})
may not be unitary, we may still diagonalize $A$, i.e., get that
$A$ is unitarily equivalent to multiplication by the independent
variable in some $L^{2}$-space; but now the corresponding $L^{2}$-space
is vector-valued, and we get a direct integral representation. This
approach is sketched in the next section.

\subsection{Direct Integral Representation}

Throughout, we assume all the Hilbert spaces are separable.

The multiplication operator version of the spectral theorem says that
$A=A^{*}$ $\Longleftrightarrow$ $A\simeq M_{\varphi}$, where 
\begin{eqnarray}
M_{\varphi}:L^{2}\left(X,\mu\right) & \longrightarrow & L^{2}\left(X,\mu\right)\label{eq:Mvphi}\\
f & \longmapsto & \varphi f.\label{eq:Mvphi1}
\end{eqnarray}
Note that $\varphi$ is real-valued. Moreover, $A$ is bounded iff
$M_{\varphi}\in L^{\infty}\left(X,\mu\right)$. When the Hilbert space
$\mathscr{H}$ is separable, we may further assume that $\mu$ is
finite, or a probability measure.

To further diagonalize $M_{\varphi}$ in the case when $\varphi$
is ``nasty'', we will need the following tool from measure theory.
\begin{defn}
\label{def:disi}Let $X$ be a locally compact and Hausdorff space,
and $\mu$ a Borel probability measure on $X$. Let $\varphi:X\rightarrow Y$
be a measurable function, and set 
\[
\nu:=\mu\circ\varphi^{-1}.
\]
A \emph{disintegration} of $\mu$ with respect to $\varphi$ is a
system of probability measures $\left\{ \mu_{y}:y\in Y\right\} $
on $X$, satisfying 
\begin{enumerate}
\item $\mu_{y}\left(X\backslash\varphi^{-1}\left(\left\{ y\right\} \right)\right)=0$,
$\nu$-a.e, i.e., $\mu_{y}$ is supported on the ``fiber'' $\varphi^{-1}\left(\left\{ y\right\} \right)$.
\item For all Borel set $E$ in $X$, the function $y\mapsto\mu_{y}\left(E\right)$
is $\nu$-measurable, and
\begin{equation}
\mu\left(E\right)=\int_{Y}\mu_{y}\left(E\right)d\nu\left(y\right).\label{eq:dism}
\end{equation}

\end{enumerate}
\end{defn}

Now, back to the Spectral Theorem. 

Let $M_{\varphi}$ be as in (\ref{eq:Mvphi})-(\ref{eq:Mvphi1}),
and let $\nu:=\mu\circ\varphi^{-1}$, i.e., a Borel probability measure
on $\mathbb{R}$. In fact, $\nu$ is supported on the essential range
of $\varphi$. 

It is well-known that, in this case, there exists a unique (up to
measure zero sets) disintegration of $\mu$ with respect to $\varphi$.
See, e.g., \cite{parthasarathy1982probability}. Therefore, we get
the direct integral decomposition \index{integral!direct-} 
\begin{equation}
L^{2}\left(\mu\right)\simeq\int^{\oplus}L^{2}(\mu_{y})d\nu\left(y\right)\;\left(\mbox{unitarily equivalent}\right)\label{eq:di}
\end{equation}
where $\left\{ \mu_{y}:y\in\mbox{essential range of }\varphi\subset\mathbb{R}\right\} $
is the system of probability measures as in \defref{disi}. 

The RHS in (\ref{eq:di}) is the Hilbert space consisting of measurable
cross-sections $f:\mathbb{R}\rightarrow\bigcup L^{2}(\mu_{y})$, where
$f\left(y\right)\in L^{2}\left(\mu_{y}\right)$, $\forall y$, and
with the inner product given by \index{direct integral} 
\begin{equation}
\left\langle f,g\right\rangle _{L^{2}\left(\nu\right)}:=\int_{Y}\left\langle f\left(y\right),g\left(y\right)\right\rangle _{L^{2}\left(\mu_{y}\right)}d\nu\left(y\right).\label{eq:di1}
\end{equation}

\begin{xca}[Direct integral Hilbert space]
\myexercise{Direct integral Hilbert space} Let the setting be as
in (\ref{eq:dism})-(\ref{eq:di}); let $\left(Y,\mathcal{F}_{Y},\nu\right)$
be a fixed measure space; and let $\left\{ \mu_{y}\right\} _{y\in Y}$
be a field of Borel measures. Show that the space of all functions
$f$ specified as follows (i)-(iii) form a Hilbert space:
\begin{enumerate}
\item[(i)]  $f:\mathbb{R}\longrightarrow\bigcup_{y\in Y}L^{2}\left(\mu_{y}\right)$;
\item[(ii)]  $y\longmapsto\left\Vert f\left(y,\cdot\right)\right\Vert _{L^{2}\left(\mu_{y}\right)}^{2}$
is measurable, and in $L^{1}\left(Y,\nu\right)$; with 
\item[(iii)]  $\int_{Y}\left\Vert f\left(y,\cdot\right)\right\Vert _{L^{2}\left(\mu_{y}\right)}^{2}d\nu\left(y\right)<\infty$. 
\end{enumerate}

Set 
\[
\left\Vert f\right\Vert _{\text{Dir. sum}}^{2}=\mbox{RHS in (iii)},
\]
and define the corresponding inner product by the RHS in (\ref{eq:di1}).

\end{xca}
\begin{thm}
Let $M_{\varphi}:L^{2}\left(X,\mu\right)\rightarrow L^{2}\left(X,\mu\right)$
be the multiplication operator in (\ref{eq:Mvphi})-(\ref{eq:Mvphi1}),
$\nu:=\mu\circ\varphi^{-1}$ as before. Then $M_{\varphi}$ is unitarily
equivalent to multiplication by the independent variable on $\int^{\oplus}L^{2}(\mu_{y})d\nu\left(y\right)$. 
\end{thm}
For details, see, e.g., \cite{Dixmier198101,Seg50}.

\subsection{\label{sec:proofsp}Proof of \thmref{spA}}

We try to get the best generalization of diagonalizing Hermitian matrices
in finite dimensional linear algebra.\index{diagonalization}\index{space!Hilbert-}

Nelson's idea \cite{Ne69} is to get from selfadjoint operators $\rightarrow$
cyclic representation of function algebra $\rightarrow$ measure $\mu$
$\rightarrow$ $L^{2}(\mu)$.\index{representation!cyclic}\index{Nelson, E}
\index{cyclic!-representation} \index{selfadjoint operator}

\emph{\uline{Sketch proof of \mbox{\thmref{spA}}:}}
\begin{enumerate}
\item Start with a single selfadjoint operator $A$ acting in an abstract
Hilbert space $\mathscr{H}$. Assume $A$ is bounded.
\item Fix $u\in\mathscr{H}$. The set $\{f(A)u\}$, as $f$ runs through
some function algebra, generates a subspace $\mathscr{H}_{u}\subset\mathscr{H}$.
$\mathscr{H}_{u}$ is called a \emph{cyclic subspace}, and $u$ the
corresponding \emph{cyclic vector}. The function algebra might be
taken as the algebra of polynomials, then later it is extended to
a much bigger algebra containing polynomials as a dense sub-algebra. 
\item Break up $\mathscr{H}$ into a direct sum of mutually disjoint cyclic
subspaces, 
\[
\mathscr{H}=\oplus_{j}\mathscr{H}_{j},
\]
with the family of cyclic vectors $u_{j}\in\mathscr{H}_{j}$. \index{cyclic!-vector}
\item Each $\mathscr{H}_{j}$ leaves $A$ invariant, and the restriction
of $A$ to each $\mathscr{H}_{j}$ is unitarily equivalent to $M_{x}$
on $L^{2}\left(sp\left(A\right),\mu_{j}\right)$, where $sp\left(A\right)$
denotes the spectrum of $A$. 
\item Piecing together all the cyclic subspace: set \index{cyclic!-subspace}
\[
X=\bigsqcup_{j}sp\left(A\right),\quad\mu=\bigsqcup_{j}\mu_{j}
\]
i.e., taking disjoint union as $u_{j}$ runs through all the cyclic
vectors. When $\mathscr{H}$ is separable, we get $\mathscr{H}=\oplus_{j\in\mathbb{N}}\mathscr{H}_{j}$,
and we may set $\mu:=\sum_{j=1}^{\infty}2^{-j}\mu_{j}$. 
\end{enumerate}
Details below.
\begin{lem}
There exists a family of cyclic vector $\left\{ u_{\alpha}\right\} $
such that $\mathscr{H}=\oplus_{\alpha}\mathscr{H}_{u_{\alpha}}$,
orthogonal sum of cyclic subspaces. \index{cyclic!-vector}\end{lem}
\begin{proof}
An application of Zorn's lemma\index{Zorn's lemma}. See \thmref{cyclic2}.\end{proof}
\begin{lem}
Set $K:=\left[-\left\Vert A\right\Vert ,\left\Vert A\right\Vert \right]$.
For each cyclic vector $u$, there exists a Borel measure $\mu_{u}$
such that $supp\left(\mu_{u}\right)\subset K$; and $\mathscr{H}_{u_{\alpha}}\simeq L^{2}\left(K,\mu_{u}\right)$. \end{lem}
\begin{proof}
The map 
\[
f\mapsto w_{u}(f):=\left\langle u,f(A)u\right\rangle _{\mathscr{H}}
\]
is a positive, bounded linear functional\index{functional} on polynomials
over $K$; the latter is dense in $C\left(K\right)$ by Stone-Weierstrass
theorem\index{Stone-Weierstrass' theorem}. Hence $w_{u}$ extends
uniquely to $C(K)$. ($w_{u}$ is a \emph{state} of the $C^{*}$-algebra
$C\left(K\right)$. See Chapter \ref{chap:GNS}.) By Riesz, there
exists a unique Borel measure $\mu_{u}$ on $K$, such that
\begin{equation}
w_{u}(f)=\left\langle u,f\left(A\right)u\right\rangle _{\mathscr{H}}=\int_{K}fd\mu_{u}.\label{eq:W0}
\end{equation}
Therefore we get $L^{2}(K,\mu_{u})$, a Hilbert space containing polynomials
as a dense subspace. Let 
\[
\mathscr{H}_{u}=\overline{span}\{f(A)u:f\in\mbox{polynomials}\}
\]
Define $W:\mathscr{H}_{u}\longrightarrow L^{2}(K,\mu_{u})$, by 
\begin{equation}
W\::\:f\left(A\right)u\longmapsto f\in L^{2}\left(\mu_{u}\right)\label{eq:W}
\end{equation}
for polynomials $f$, which then extends to $\mathscr{H}_{u}$ by
density.

\index{Riesz' theorem}

\index{extension!-of positive functional}

\index{Theorem!Riesz-}

\index{integral!Riesz-}\end{proof}
\begin{lem}
Let $W$ be the operator in (\ref{eq:W}), then 
\begin{enumerate}
\item $W$ is an isometric isomorphism; and
\item $WA=M_{t}W$, i.e., $W$ intertwines $A$ and $M_{t}$. Hence $W$
diagonalizes $A$. 
\end{enumerate}
\end{lem}
\begin{rem}
$WA=M_{t}W\Longleftrightarrow WAW^{*}=M_{t}$. In finite dimension,
it is less emphasized that the adjoint $W^{*}$ equals the inverse
$W^{-1}$. For finite dimensional case, $M_{t}=diag(\lambda_{1},\lambda_{2},\ldots,\lambda_{n})$
where the measure $\mu=\sum_{\text{finite}}\delta_{\lambda_{i}}$,
where $\delta_{\lambda_{i}}$ is the Dirac measure at the eigenvalue
$\lambda_{i}$ of $A$.\index{eigenvalue}\end{rem}
\begin{proof}
For the first part, let $f\in L^{2}\left(\mu_{u}\right)$, then 
\begin{eqnarray*}
\left\Vert f\right\Vert _{L^{2}\left(K,\mu_{u}\right)}^{2}=\int_{\mathbb{R}}\left|f\right|^{2}d\mu_{u} & = & \left\langle u,\left|f\right|^{2}(A)u\right\rangle _{\mathscr{H}}\\
 & = & \left\langle u,\bar{f}(A)f(A)u\right\rangle _{\mathscr{H}}\\
 & = & \left\langle u,f(A)^{*}f(A)u\right\rangle _{\mathscr{H}}\\
 & = & \left\langle f(A)u,f(A)u\right\rangle _{\mathscr{H}}\\
 & = & \left\Vert f(A)u\right\Vert _{\mathscr{H}}^{2}.
\end{eqnarray*}
Notice that strictly speaking, $f(A^{*})^{*}=\bar{f}(A)$. Since $A^{*}=A$,
we then get 
\[
f\left(A\right)^{*}=\overline{f}\left(A\right).
\]
Also, $f\xrightarrow{\:\pi\:}f\left(A\right)$ is a $*$-representation
of the algebra $C\left(K\right)$; i.e., $\pi$ is a homomorphism\index{homomorphism},
and $\pi(\bar{f})=f(A)^{*}$.

For the second part, let $f(t)=a_{0}+a_{1}t+\cdots a_{n}t^{n}$ be
any polynomial, then 
\begin{eqnarray*}
WAf(A) & = & WA(a_{0}+a_{1}A+a_{2}A^{2}+\cdots+a_{n}A^{n})\\
 & = & W(a_{0}A+a_{1}A^{2}+a_{2}A^{3}+\cdots+a_{n}A^{n+1})\\
 & = & a_{0}t+a_{1}t^{2}+a_{2}t^{3}+\cdots+a_{n}t^{n+1}\\
 & = & tf(t)\\
 & = & M_{t}Wf(A)
\end{eqnarray*}
thus $WA=M_{t}W$. The assertion then follows from standard approximation.

It remains to show that the isometry $\mathscr{H}_{u}\xrightarrow{\;W\;}L^{2}\left(\mu_{u}\right)$
in (\ref{eq:W}) maps \uline{onto} $L^{2}\left(\mu_{u}\right)$.
But this follows from (\ref{eq:W0}). Indeed if $f\in L^{2}\left(\mu_{u}\right)$,
then $f\left(A\right)u\in\mathscr{H}_{u}$ is well defined by the
reasoning above. As a result, for the adjoint operator $L^{2}\left(\mu_{u}\right)\xrightarrow{\:W^{*}\:}\mathscr{H}_{u}$,
we have 
\begin{equation}
W^{*}\left(f\right)=f\left(A\right)u,\;\forall f\in L^{2}\left(\mu_{u}\right).\label{eq:W1}
\end{equation}

\end{proof}
Finally we piece together all the cyclic subspaces. \index{cyclic!-subspace}
\begin{lem}
There exists a locally compact Hausdorff space $X$ and a Borel measure
$\mu$, a unitary operator $\mathcal{F}:\mathscr{H}\longrightarrow L^{2}\left(X,\mu\right)$,
such that 
\[
A=\mathcal{F}^{*}M_{\varphi}\mathcal{F}
\]
where $\varphi\in L^{\infty}\left(\mu\right)$. \end{lem}
\begin{proof}
Recall that we get a family of states $w_{j}$, with the corresponding
measures $\mu_{j}$, and Hilbert spaces $\mathscr{H}_{j}=L^{2}(\mu_{j})$.
Note that all the $L^{2}$-spaces are on $K=sp(A)$. So it's the same
underlying set, but with possibly different measures.
\end{proof}
To get a single measure space with $\mu$, Nelson \cite{Ne69} suggested
taking the disjoint union 
\[
X:=\bigcup_{j}K\times\{j\}
\]
and $\mu:=$ the disjoint union of $\mu_{j}'s$. The existence of
$\mu$ follows from Riesz. Then we get\index{Theorem!Riesz-}
\[
\mathscr{H}=\oplus\mathscr{H}_{j}\xrightarrow{\;\mathcal{F}\;}L^{2}(X,\mu).
\]
 \index{Nelson, E}
\begin{rem}
Note the representation of $L^{\infty}\left(X,\mu\right)$ onto $\mathscr{H}=\sum^{\oplus}\mathscr{H}_{j}$
is highly non unique. There we enter into the multiplicity\index{multiplicity}
theory, which starts with breaking up each $\mathscr{H}_{j}$ into
irreducible components.
\end{rem}

\begin{rem}
A representation $\pi\in Rep\left(\mathfrak{A},\mathscr{B}\left(\mathscr{H}\right)\right)$
is said to be multiplicity free if and only if $\pi\left(\mathfrak{A}\right)'$
is abelian. We say $\pi$ has multiplicity equal to $n$ if and only
if $(\pi(\mathfrak{A}))'\simeq M_{n}(\mathbb{C})$. This notation
of multiplicity generalizes the one in finite dimensional linear algebra.
See  \secref{mult}.
\end{rem}
\index{multiplicity free}
\begin{xca}[Multiplicity free]
\myexercise{Multiplicity free}\label{exer:multiplicity}Prove that
$\pi\in Rep(L^{\infty}(\mu),L^{2}(\mu))$ is multiplicity free. Conclude
that each cyclic representation is multiplicity free (i.e., it is
maximal abelian.) \index{cyclic!-representation}

\uline{Hint:} Suppose $B\in\mathscr{B}$$\left(L^{2}\left(\mu\right)\right)$
commutes with all $M_{\varphi}$, $\varphi\in L^{2}\left(\mu\right)$.
Define $g=B\mathbbm{1}$, where $\mathbbm{1}$ is the constant function.
Then, for all $\psi\in L^{2}\left(\mu\right)$, we have 
\[
B\psi=B\psi\mathbbm{1}=BM_{\psi}\mathbbm{1}=M_{\psi}B\mathbbm{1}=M_{\psi}g=\psi g=g\psi=M_{g}\psi
\]
thus $B=M_{g}$.
\end{xca}

\begin{cor}
\label{cor:unitary}$T\in\mathscr{B}\left(\mathscr{H}\right)$ is
unitary iff there exists $\mathcal{F}:\mathscr{H}\rightarrow L^{2}\left(X,d\mu\right)$,
unitary, such that 
\[
T=\mathcal{F}^{*}M_{z}\mathcal{F},
\]
where $\left|z\right|\in\mathbb{T}^{1}=\left\{ z\in\mathbb{C}:\left|z\right|=1\right\} $. \end{cor}
\begin{xca}[The Cayley transform]
\begin{flushleft}
\label{exer:ubdsa}\myexercise{The Cayley transform}Finish the proof
of \thmref{spA} for the case when $A$ is unbounded and selfadjoint. 
\par\end{flushleft}

\begin{flushleft}
\uline{Hint}: Suppose $A=A^{*}$, unbounded. The Cayley transform
\[
C_{A}:=\left(A-i\right)\left(A+i\right)^{-1}
\]
is then unitary. See \secref{Cayley}. Apply \corref{unitary} to
$C_{A}$, and convert the result back to $A$. See, e.g., \cite{Ne69,Rud73}.
\par\end{flushleft}
\end{xca}

\section{Projection-Valued Measure (PVM)\label{sec:PVM}}

A projection valued measure (PVM) $P$ satisfies the usual axioms
of measures (here Borel measures) but with the main difference: \index{Spectral Theorem}\index{Borel measure}\index{axioms}
\begin{enumerate}
\item $P\left(\triangle\right)$ is a projection for all $\triangle\in\mathcal{B}$,
i.e., $P\left(\triangle\right)=P\left(\triangle\right)^{*}=P\left(\triangle\right)^{2}$.
\item We assume that $P\left(\triangle_{1}\cap\triangle_{2}\right)=P\left(\triangle_{1}\right)P\left(\triangle_{2}\right)$
for all $\triangle_{1},\triangle_{2}\in\mathcal{B}$; this property
is called ``orthogonality''.\end{enumerate}
\begin{rem}
The notion of PVM extends the familiar notion of an ONB:

Let $\mathscr{H}$ be a separable Hilbert space, and suppose $\left\{ u_{k}\right\} _{k\in\mathbb{N}}$
is a ONB in $\mathscr{H}$; then for $\triangle\in\mathcal{B}\left(\mathbb{R}\right)$,
set 
\end{rem}
\begin{equation}
P\left(\triangle\right):=\sum_{k\in\triangle}\left|u_{k}\left\rangle \right\langle u_{k}\right|.\label{eq:Pdel}
\end{equation}

\begin{xca}[A concrete PVM]
\myexercise{A concrete PVM}\label{exer:pvm}Show that $P$ is a
PVM. 

Note that, under the summation on the RHS in (\ref{eq:Pdel}), we
used Dirac's notation $\left|u_{k}\left\rangle \right\langle u_{k}\right|$
for the rank-one projection onto $\mathbb{C}u_{k}$. Further note
that the summation is over $k$ from $\triangle$; so it varies as
$\triangle$ varies. 
\end{xca}
The PVM version of the spectral theorem says that $A=A^{*}$ iff\index{measure!projection-valued}
\[
A=\int xP_{A}(dx)
\]
where $P$ is a projection-valued measure defined on the Borel $\sigma$-algebra
of $\mathbb{R}$. 
\begin{defn}
\label{def:PVM}Let $\mathcal{B}\left(\mathbb{C}\right)$ be the Borel
$\sigma$-algebra of $\mathbb{C}$. $\mathscr{H}$ is a Hilbert space.
\[
P:\mathcal{B}(\mathbb{C})\rightarrow\mbox{Proj}\left(\mathscr{H}\right)
\]
is a \emph{projection-valued measure} (PVM), if \index{measure!projection-valued}
\begin{enumerate}
\item $P(\emptyset)=0$, $P(\mathbb{C})=I$, $P\left(A\right)$ is a projection
for all $A\in\mathcal{B}$
\item $P(A\cap B)=P(A)P(B)$
\item $P(\cup_{k}E_{k})=\sum P(E_{k})$, $E_{k}\cap E_{j}=\phi$ if $k\neq j$.
The convergence is in terms of the strong operator topology. By assumption,
the sequence of projections $\sum_{k=1}^{N}P(E_{k})$ is monotone
increasing, hence it has a limit, and 
\[
\lim_{N\rightarrow\infty}\sum_{k=1}^{N}P\left(E_{k}\right)=P\left(\cup E_{k}\right).
\]

\end{enumerate}
\end{defn}
The standard Lebesgue integration extends to PVM. 
\[
\left\langle \varphi,P(E)\varphi\right\rangle =\left\langle P(E)\varphi,P(E)\varphi\right\rangle =\left\Vert P(E)\right\Vert ^{2}\geq0
\]
since $P$ is countably additive, the map $E\mapsto P\left(E\right)$
is also countably additive. Therefore, each $\varphi\in\mathscr{H}$
induces a regular Borel measure $\mu_{\varphi}$ on the Borel $\sigma$-algebra
of $\mathbb{R}$.

For a measurable function $\psi$, 
\begin{eqnarray*}
\int\psi d\mu_{\varphi} & = & \int\psi(x)\left\langle \varphi,P(dx)\varphi\right\rangle \\
 & = & \left\langle \varphi,\left(\int\psi P(dx)\right)\varphi\right\rangle 
\end{eqnarray*}
hence we may define 
\[
\int\psi P(dx)
\]
as the operator so that for all $\varphi\in\mathscr{H}$, 
\[
\left\langle \varphi,\left(\int\psi P(dx)\right)\varphi\right\rangle .
\]

\begin{rem}
$P(E)=F\chi_{E}F^{-1}$ defines a PVM. In fact all PVMs come from
this way. In this sense, the $M_{t}$ version of the spectral theorem
is better, since it implies the PVM version. However, the PVM version
facilitates some formulations in quantum mechanics, so physicists
usually prefer this version.\index{Theorem!Spectral-}
\end{rem}

\begin{rem}
Suppose we start with the PVM version of the spectral theorem. How
to prove $(\psi_{1}\psi_{2})(A)=\psi_{1}(A)\psi_{2}(A)$? i.e. how
to check we do have an algebra isomorphism? Recall in the PVM version,
$\psi(A)$ is defined as the operator so that for all $\varphi\in\mathscr{H}$,
we have
\[
\int\psi d\mu_{\varphi}=\left\langle \varphi,\psi(A)\varphi\right\rangle .
\]
As a standard approximation technique, once starts with simple or
even step functions. Once it is worked out for simple functions, the
extension to any measurable functions is straightforward. Hence let's
suppose (WLOG) that the functions are simple.\end{rem}
\begin{lem}
We have $\psi_{1}\left(A\right)\psi_{2}\left(A\right)=\left(\psi_{1}\psi_{2}\right)\left(A\right)$.\end{lem}
\begin{proof}
Let 
\begin{eqnarray*}
\psi_{1} & = & \sum\psi_{1}(t_{i})\chi_{E_{i}}\\
\psi_{2} & = & \sum\psi_{2}(t_{j})\chi_{E_{j}}
\end{eqnarray*}
then
\begin{eqnarray*}
\int\psi_{1}P(dx)\int\psi_{2}P(dx) & = & \sum_{i,j}\psi_{1}(t_{i})\psi_{2}(t_{j})P(E_{i})P(E_{j})\\
 & = & \sum_{i}\psi_{1}(t_{i})\psi_{2}(t_{i})P(E_{i})\\
 & = & \int\psi_{1}\psi_{2}P(dx)
\end{eqnarray*}
where we used the fact that $P(A)P(B)=0$ if $A\cap B=\phi$.\end{proof}
\begin{rem}
As we delve into Nelson's lecture notes \cite{Ne69}, we notice that
on page 69, there is another unitary operator. By piecing these operators
together is precisely how we get the spectral theorem. This ``piecing''
is a vast generalization of Fourier series. \end{rem}
\begin{lem}
Pick $\varphi\in\mathscr{H}$, get the measure $\mu_{\varphi}$ where
\[
\mu_{\varphi}(\cdot)=\left\Vert P\left(\cdot\right)\varphi\right\Vert ^{2}
\]
and we have the Hilbert space $L^{2}(\mu_{\varphi})$. Take $\mathscr{H}_{\varphi}:=\overline{span}\{\psi(A)\varphi:\psi\in L^{2}(\mu_{\varphi})\}$.
Then the map
\[
\mathscr{H}\ni\psi(A)\varphi\mapsto\psi\in L^{2}\left(\mu_{\varphi}\right)
\]
is an isometry, and it extends uniquely to a unitary operator from
$\mathscr{H}_{\varphi}$ to $L^{2}(\mu_{\varphi})$. \index{isometry}\end{lem}
\begin{proof}
We have 
\begin{eqnarray*}
\left\Vert \psi(A)\varphi\right\Vert ^{2} & = & \left\langle \psi(A)\varphi,\psi(A)\varphi\right\rangle \\
 & = & \left\langle \varphi,\bar{\psi}(A)\psi(A)\varphi\right\rangle \\
 & = & \int_{\mathbb{R}}\left|\psi\left(\lambda\right)\right|^{2}\left\Vert P\left(d\lambda\right)\varphi\right\Vert ^{2}\\
 & = & \left\langle \varphi,\left|\psi\right|^{2}(A)\varphi\right\rangle \\
 & = & \int\left|\psi\right|^{2}d\mu_{\varphi}.
\end{eqnarray*}
\end{proof}
\begin{rem}
$\mathscr{H}_{\varphi}$ is called the cyclic space generated by $\varphi$.
Before we can construct $\mathscr{H}_{\varphi}$, we must make sense
of $\psi(A)\varphi$. \index{cyclic!-space}\end{rem}
\begin{lem}[{\cite[p.67]{Ne69}}]
 Let $p=a_{0}+a_{1}x+\cdots+a_{n}x^{n}$ be a polynomial. Then $\left\Vert p(A)u\right\Vert \leq\max\left|p(t)\right|$,
where $\left\Vert u\right\Vert =1$ i.e. $u$ is a state.\end{lem}
\begin{proof}
$M:=span\{u,Au,\ldots,A^{n}u\}$ is a finite dimensional subspace
in $\mathscr{H}$ (automatically closed). Let $E$ be the orthogonal
projection onto $M$. Then
\[
p(A)u=Ep(A)Eu=p(EAE)u.
\]
Since $EAE$ is a Hermitian matrix on $M$, we may apply the spectral
theorem for finite dimensional space and get
\[
EAE=\sum\lambda_{k}P_{\lambda_{k}}
\]
where $\lambda_{k}'s$ are eigenvalues associated with the projections
$P_{\lambda_{k}}$. It follows that 
\begin{eqnarray*}
p(A)u & = & p(\sum\lambda_{k}P_{\lambda_{k}})u=\left(\sum p(\lambda_{k})P_{\lambda_{k}}\right)u
\end{eqnarray*}
and 
\begin{eqnarray*}
\left\Vert p(A)u\right\Vert ^{2} & = & \sum\left|p(\lambda_{k})\right|^{2}\left\Vert P_{\lambda_{k}}u\right\Vert ^{2}\\
 & \leq & \max\left|p(t)\right|^{2}\sum\left\Vert P_{\lambda_{k}}u\right\Vert ^{2}\\
 & = & \max\left|p(t)\right|^{2}
\end{eqnarray*}
since 
\[
\sum\left\Vert P_{\lambda_{k}}u\right\Vert ^{2}=\left\Vert u\right\Vert ^{2}=1.
\]
Notice that $I=\sum P_{\lambda_{k}}$.\end{proof}
\begin{rem}
How to extend this? polynomials - continuous functions - measurable
functions. $\left[-\left\Vert A\right\Vert ,\left\Vert A\right\Vert \right]\subset\mathbb{R}$,
\[
\left\Vert EAE\right\Vert \leq\left\Vert A\right\Vert 
\]
is a uniform estimate for all truncations. Apply Stone-Weierstrass'
theorem\index{Stone-Weierstrass' theorem} to the interval $\left[-\left\Vert A\right\Vert ,\left\Vert A\right\Vert \right]$
we get that any continuous function $\psi$ is uniformly approximated
by polynomials. i.e. $\psi\sim p_{n}$. Thus
\[
\left\Vert p_{n}(A)u-p_{m}(A)u\right\Vert \leq\max\left|p_{n}-p_{m}\right|\left\Vert u\right\Vert =\left\Vert p_{n}-p_{m}\right\Vert _{\infty}\rightarrow0
\]
and $p_{n}(A)u$ is a Cauchy sequence, hence 
\[
\lim_{n}p_{n}(A)u=:\psi(A)u
\]
where we may define the operator $\psi(A)$ so that $\psi(A)u$ is
the limit of $p_{n}(A)u$.
\end{rem}

\section{\label{sec:MPVM}Convert $M_{\varphi}$ to a PVM (projection-valued
measure)}

\index{Spectral Theorem}
\begin{thm}
Let $A:\mathscr{H}\rightarrow\mathscr{H}$ be a selfadjoint operator.
Then $A$ is unitarily equivalent to the operator $M_{t}$ of multiplication
by the independent variable on the Hilbert space $L^{2}(\mu)$. There
exists a unique projection-valued measure $P$ so that 
\[
A=\int tP\left(dt\right)
\]
i.e. for all $h,k\in\mathscr{H}$, 
\[
\left\langle k,Ah\right\rangle _{\mathscr{H}}=\int t\left\langle k,P\left(dt\right)h\right\rangle _{\mathscr{H}}.
\]
\end{thm}
\begin{proof}
The uniqueness part follows from a standard argument. We will only
prove the existence of $P$. 

Let $\mathcal{F}:L^{2}(\mu)\rightarrow\mathscr{H}$ be the unitary
operator so that $A=\mathcal{F}M_{t}\mathcal{F}^{*}$. Define 
\[
P\left(E\right):=\mathcal{F}\chi_{E}\mathcal{F}^{*}
\]
for all $E$ in the Borel $\sigma$-algebra $\mathfrak{B}$ of $\mathbb{R}$.
Then $P\left(\emptyset\right)=0$, $P\left(\mathbb{R}\right)=I$;
and for all $E_{1},E_{2}\in\mathfrak{B}$, 
\begin{eqnarray*}
P\left(E_{1}\cap E_{2}\right) & = & \mathcal{F}\chi_{E_{1}\cap E_{2}}\mathcal{F}^{-1}\\
 & = & \mathcal{F}\chi_{E_{1}}\chi_{E_{2}}\mathcal{F}^{-1}\\
 & = & \left(\mathcal{F}\chi_{E_{1}}\mathcal{F}^{-1}\right)\left(\mathcal{F}\chi_{E_{2}}\mathcal{F}^{-1}\right)\\
 & = & P\left(E_{1}\right)P\left(E_{2}\right).
\end{eqnarray*}
Suppose $\{E_{k}\}$ is a sequence of mutually disjoint elements in
$\mathfrak{B}$. Let $h\in\mathscr{H}$ and write $h=\mathcal{F}\widehat{h}$
for some $\widehat{h}\in L^{2}(\mu)$. Then
\begin{eqnarray*}
\left\langle h,P\left(\cup_{k}E_{k}\right)h\right\rangle _{\mathscr{H}} & = & \left\langle \mathcal{F}\widehat{h},P(\cup E_{k})F\widehat{h}\right\rangle _{\mathscr{H}}=\left\langle \widehat{h},\mathcal{F}^{-1}P\left(\cup E_{k}\right)\mathcal{F}\widehat{h}\right\rangle _{L^{2}\left(\mu\right)}\\
 & = & \left\langle \widehat{h},\chi_{\cup E_{k}}\widehat{h}\right\rangle _{L^{2}\left(\mu\right)}=\int_{\cup E_{k}}\big|\widehat{h}\big|^{2}d\mu\\
 & = & \sum_{k}\int_{E_{k}}\big|\widehat{h}\big|^{2}d\mu=\sum_{k}\left\langle h,P(E_{k})h\right\rangle _{\mathscr{H}}.
\end{eqnarray*}
Therefore, $P$ is a projection-valued measure. 

For any $h,k\in\mathscr{H}$, write $h=\mathcal{F}\widehat{h}$ and
$k=\mathcal{F}\widehat{k}$. Then 
\begin{eqnarray*}
\left\langle k,Ah\right\rangle _{\mathscr{H}} & = & \left\langle \mathcal{F}\widehat{k},A\mathcal{F}\widehat{h}\right\rangle _{\mathscr{H}}=\left\langle \widehat{k},\mathcal{F}^{*}A\mathcal{F}\widehat{h}\right\rangle _{\mathscr{H}}\\
 & = & \left\langle \widehat{k},M_{t}\widehat{h}\right\rangle _{L^{2}\left(\mu\right)}=\int t\overline{\widehat{k}(t)}\widehat{h}(t)d\mu(t)\\
 & = & \int t\left\langle k,P\left(dt\right)h\right\rangle _{\mathscr{H}}.
\end{eqnarray*}
Thus $A=\int tP(dt)$.\end{proof}
\begin{rem}
In fact, $A$ is in the closed (under norm or strong topology) span
of $\{P(E):E\in\mathfrak{B}\}$. Equivalently, since $M_{t}=\mathcal{F}^{*}A\mathcal{F}$,
the function $f\left(t\right)=t$ is in the closed span of the set
of characteristic functions;  the latter is again a standard approximation
in measure theory. It suffices to approximate $t\chi_{[0,\infty]}$. 
\end{rem}
The wonderful idea of Lebesgue is not to partition the domain, as
was the case in Riemann integral over $\mathbb{R}^{n}$, but instead
the range. Therefore integration over an arbitrary set is made possible.
Important examples include analysis on groups. \index{operators!bounded}
\begin{prop}
\label{pro:approx_step}Let $f:[0,\infty]\rightarrow\mathbb{R}$,
$f(x)=x$, i.e. $f=x\chi_{[0,\infty]}$. Then there exists a sequence
of step functions $s_{1}\leq s_{2}\leq\cdots\leq f(x)$ such that
$\lim_{n\rightarrow\infty}s_{n}(x)=f(x)$.\end{prop}
\begin{proof}
For $n\in\mathbb{N}$, define
\[
s_{n}(x)=\begin{cases}
i2^{-n} & x\in[i2^{-n},(i+1)2^{-n})\\
n & x\in[n,\infty]
\end{cases}
\]
where $0\leq i\leq n2^{-n}-1$. Equivalently, $s_{n}$ can be written
using characteristic functions as 
\[
s_{n}=\sum_{i=0}^{n2^{n}-1}i2^{-n}\chi_{[i2^{-n},(i+1)2^{-n})}+n\chi_{[n,\infty]}.
\]
Notice that on each interval $[i2^{-n},(i+1)2^{-n})$, 
\begin{eqnarray*}
s_{n}(x) & \equiv & i2^{-n}\leq x\\
s_{n}(x)+2^{-n} & \equiv & (i+1)2^{-n}>x\\
s_{n}(x) & \leq & s_{n+1}(x).
\end{eqnarray*}
Therefore, for all $n\in\mathbb{N}$ and $x\in[0,\infty]$, 
\begin{equation}
x-2^{-n}<s_{n}(x)\leq x\label{eq:conv}
\end{equation}
and $s_{n}(x)\leq s_{n+1}(x)$. 

It follows from (\ref{eq:conv}) that 
\[
\lim_{n\rightarrow\infty}s_{n}(x)=f(x)
\]
for all $x\in[0,\infty]$.\end{proof}
\begin{cor}
Let $f(x)=x\chi_{[0,M]}(x)$. Then there exists a sequence of step
functions $s_{n}$ such that $0\leq s_{1}\leq s_{2}\leq\cdots\leq f(x)$
and $s_{n}\rightarrow f$ uniformly, as $n\rightarrow\infty$.\end{cor}
\begin{proof}
Define $s_{n}$ as in Proposition \ref{pro:approx_step}. Let $n>M$,
then by construction
\[
f(x)-2^{-n}<s_{n}(x)\leq f(x)
\]
for all $s\in[0,M]$. Hence $s_{n}\rightarrow f$ uniformly as $n\rightarrow\infty$.
\end{proof}
Proposition \ref{pro:approx_step} and its corollary immediate imply
the following.
\begin{cor}
Let $(X,S,\mu)$ be a measure space. A function (real-valued or complex-valued)
is measurable if and only if it is the point-wise limit of a sequence
of simple functions. A function is bounded measurable if and only
if it is the uniform limit of a sequence of simple functions. Let
$\{s_{n}\}$ be an approximation sequence of simple functions. Then
$s_{n}$ can be chosen such that $\left|s_{n}(x)\right|\leq\left|f(x)\right|$
for all $n=1,2,3\ldots$. \end{cor}
\begin{thm}
Let $M_{f}:L^{2}(X,S,\mu)\rightarrow L^{2}(X,S,\mu)$ be the operator
of multiplication by $f$. Then,
\begin{enumerate}
\item if $f\in L^{\infty}$, $M_{f}$ is a bounded operator, and $M_{f}$
is in the closed span of the set of selfadjoint projections under
norm topology. 
\item if $f$ is unbounded, $M_{f}$ is an unbounded operator. $M_{f}$
is in the closed span of the set of selfadjoint projections under
the strong operator topology. 
\end{enumerate}
\end{thm}
\begin{proof}
If $f\in L^{\infty}$, then there exists a sequence of simple functions
$s_{n}$ so that $s_{n}\rightarrow f$ uniformly. Hence $\left\Vert f-s_{n}\right\Vert _{\infty}\rightarrow0$,
as $n\rightarrow\infty$. 

Suppose $f$ is unbounded. By Proposition \ref{pro:approx_step} and
its corollaries, there exists a sequence of simple functions $s_{n}$
such that $\left|s_{n}(x)\right|\leq\left|f(x)\right|$ and $s_{n}\rightarrow f$
point-wisely, as $n\rightarrow\infty$. Let $h$ be any element in
the domain of $M_{f}$, i.e. 
\[
\int\left(\left|h\right|+\left|fh\right|^{2}\right)d\mu<\infty.
\]
Then
\[
\lim_{n\rightarrow\infty}\left|(f(x)-s_{n}(x))h(x)\right|^{2}=0
\]
and 
\[
\left|(f(x)-s_{n}(x))h(x)\right|^{2}\leq\mbox{const}\cdot\left|h(x)\right|^{2}.
\]
Hence by the dominated convergence theorem, \index{Lebesgue dominated convergence theorem}
\[
\lim_{n\rightarrow\infty}\int\left|(f(x)-s_{n}(x))h(x)\right|^{2}d\mu=0
\]
or equivalently,
\[
\left\Vert (f-s_{n})h\right\Vert ^{2}\rightarrow0
\]
as $n\rightarrow\infty$. i.e. $M_{s_{n}}$ converges to $M_{f}$
in the strong operator topology. 
\end{proof}

\begin{xca}[An application to numerical range]
\myexercise{An application to numerical range} Let $A$ be a bounded
normal operator in a separable Hilbert space $\mathscr{H}$; then
prove that 
\begin{equation}
NR_{A}\subseteq\overline{conv}\left(spec\left(A\right)\right);\label{eq:nr1}
\end{equation}
i.e., that the numerical range of $A$ is contained in the closed
convex hull of the spectrum of $A$. (We refer to  \exerref{nrange}
for details on ``numerical range.'')

\uline{Hint}: Since $A$ is normal, by the Spectral Theorem, it
is represented in a PVM $P_{A}\left(\cdot\right)$, i.e., taking valued
in $\mbox{Proj}\left(\mathscr{H}\right)$. For $x\in\mathscr{H}$,
$\left\Vert x\right\Vert =1$, we have 
\begin{equation}
w_{x}\left(A\right)=\left\langle x,Ax\right\rangle =\int_{spec\left(A\right)}\lambda\left\Vert P_{A}\left(d\lambda\right)x\right\Vert ^{2},\;\mbox{and}\label{eq:nr2}
\end{equation}
\begin{equation}
\int_{spec\left(A\right)}\left\Vert P_{A}\left(d\lambda\right)x\right\Vert ^{2}=\left\Vert x\right\Vert ^{2}=1,\label{eq:nr3}
\end{equation}
so $d\mu_{x}\left(\lambda\right)=\left\Vert P_{A}\left(d\lambda\right)x\right\Vert ^{2}$
is a regular Borel probability measure on $spec\left(A\right)$. Now
approximate (\ref{eq:nr2}) with simple functions on $spec\left(A\right)$.\index{Theorem!Spectral-}
\end{xca}

\paragraph{Quantum States }

Let $\mathscr{H}$ be a Hilbert space (corresponding to some quantum
system), and let $A$ be a selfadjoint operator in $\mathscr{H}$,
possibly unbounded. Vectors $f\in\mathscr{H}$ represent quantum states
if $\left\Vert f\right\Vert =1$. 

The \emph{mean} of $A$ in the state $f$ is 
\[
\left\langle f,Af\right\rangle =\int_{\mathbb{R}}\lambda\left\Vert P_{A}\left(d\lambda\right)f\right\Vert ^{2}
\]
where $P_{A}\left(\cdot\right)$ denotes the spectral resolution of
$A$.

The \emph{variance} of $A$ in the state $f$ is 
\begin{eqnarray}
v_{f}\left(A\right) & = & \left\Vert Af\right\Vert ^{2}-\left(\left\langle f,Af\right\rangle \right)^{2}\label{eq:var}\\
 & = & \int_{\mathbb{R}}\lambda^{2}\left\Vert P_{A}\left(d\lambda\right)f\right\Vert ^{2}-\left(\int_{\mathbb{R}}\lambda\left\Vert P_{A}\left(d\lambda\right)f\right\Vert ^{2}\right)^{2}\nonumber \\
 & = & \int_{\mathbb{R}}\left(\lambda-\left\langle f,Af\right\rangle \right)^{2}\left\Vert P_{A}\left(d\lambda\right)f\right\Vert ^{2}.\nonumber 
\end{eqnarray}

\begin{thm}[Uncertainty Principle]
Let $\mathscr{D}$ be a dense subspace in $\mathscr{H}$, and $A,B$
be two Hermitian operators such that $A,B:\mathscr{D}\hookrightarrow\mathscr{D}$
(i.e., $\mathscr{D}$ is assumed invariant under both $A$ and $B$.) 

Then, 
\begin{equation}
\left\Vert Ax\right\Vert \left\Vert Bx\right\Vert \geq\frac{1}{2}\left|\left\langle x,\left[A,B\right]x\right\rangle \right|,\;\forall x\in\mathscr{D};\label{eq:up2-1}
\end{equation}
where $\left[A,B\right]:=AB-BA$ is the commutator of $A$ and $B$.

In particular, setting 
\begin{eqnarray*}
A_{1} & := & A-\left\langle x,Ax\right\rangle \\
B_{1} & := & B-\left\langle x,Bx\right\rangle 
\end{eqnarray*}
then $A_{1},B_{1}$ are Hermitian, and \index{commutator} 
\[
\left[A_{1},B_{1}\right]=\left[A,B\right].
\]
Therefore, 
\[
\left\Vert A_{1}x\right\Vert \left\Vert B_{1}x\right\Vert \geq\frac{1}{2}\left|\left\langle x,\left[A,B\right]x\right\rangle \right|,\;\forall x\in\mathscr{D}.
\]
\end{thm}
\begin{proof}
By the Cauchy-Schwarz inequality (\lemref{CS}), and for $x\in\mathscr{D}$,
we have
\begin{eqnarray*}
\left\Vert Ax\right\Vert \left\Vert Bx\right\Vert  & \geq & \left|\left\langle Ax,Bx\right\rangle \right|\\
 & \geq & \left|\Im\left\{ \left\langle Ax,Bx\right\rangle \right\} \right|\\
 & = & \frac{1}{2}\left|\left\langle Ax,Bx\right\rangle -\overline{\left\langle Ax,Bx\right\rangle }\right|\\
 & = & \frac{1}{2}\left|\left\langle Ax,Bx\right\rangle -\left\langle Bx,Ax\right\rangle \right|\\
 & = & \frac{1}{2}\left|\left\langle x,\left[AB-BA\right]x\right\rangle \right|.
\end{eqnarray*}

\end{proof}
\index{inequality!Cauchy-Schwarz}
\begin{cor}
If $\left[A,B\right]=ihI$, $h\in\mathbb{R}_{+}$, and $\left\Vert x\right\Vert =1$
(i.e., is a state); then
\begin{equation}
w_{x}\left(A^{2}\right)^{\frac{1}{2}}w_{x}\left(B^{2}\right)^{\frac{1}{2}}\geq\frac{h}{2}.
\end{equation}
\end{cor}
\begin{xca}[Heisenberg's uncertainty principle]
\myexercise{Heisenberg's uncertainty}\label{exer:uncertainty}Let
$\mathscr{H}=L^{2}\left(\mathbb{R}\right)$, and $f\in L^{2}\left(\mathbb{R}\right)$
given, with $\left\Vert f\right\Vert ^{2}=\int\left|f\left(x\right)\right|^{2}dx=1$.
Suppose $f\in dom\left(P\right)\cap dom\left(Q\right)$, where $P$
and $Q$ are the momentum and position operators, respectively. (See
eq. (\ref{eq:sr3-1})-(\ref{eq:sr3-2}) in \secref{dga}.) 

Show that 
\begin{equation}
v_{f}\left(P\right)v_{f}\left(Q\right)\geq\frac{1}{4}.\label{eq:un}
\end{equation}
Inequality (\ref{eq:un}) is the mathematically precise form of Heisenberg's
uncertainty relation, often written in the form
\[
\sigma_{f}\left(P\right)\sigma_{f}\left(Q\right)\geq\frac{1}{2}
\]
where $\sigma_{f}\left(P\right)=\sqrt{v_{f}\left(P\right)}$, and
$\sigma_{f}\left(Q\right)=\sqrt{v_{f}\left(Q\right)}$. 
\end{xca}
\index{uncertainty}

\index{quantum mechanics!uncertainty principle}

\index{operators!position-}

\section{\label{sec:spcpt}The Spectral Theorem for Compact Operators}

\index{Spectral Theorem}

\subsection{Preliminaries}

The setting for the first of the two Spectral Theorems (direct integral
vs representation) we will consider is as following: (Restricting
assumptions will be relaxed in subsequent versions!)

Let $\mathscr{H}$ be a separable (typically infinite dimensional
Hilbert space assumed here!) and let $A\in\mathscr{B}\left(\mathscr{H}\right)\backslash\left\{ 0\right\} $
be compact and selfadjoint, i.e., $A=A^{*}$, and $A$ is in the $\left\Vert \cdot\right\Vert _{UN}$-closure
of $\mathscr{F}R\left(\mathscr{H}\right)$. See \secref{3norm}.\index{space!Hilbert-}
\begin{thm}
\label{thm:spcpt}With $A$ as above, there is an orthonormal set
$\left\{ u_{k}\right\} _{k\in\mathbb{N}}$, and a sequence $\left\{ \lambda_{k}\right\} _{k\in\mathbb{N}}$
in $\mathbb{R}\backslash\left\{ 0\right\} $, such that $\left|\lambda_{1}\right|\geq\left|\lambda_{2}\right|\geq\cdots\geq\left|\lambda_{k}\right|\geq\left|\lambda_{k+1}\right|\geq\cdots$,
$\lim_{k\rightarrow\infty}\lambda_{k}=0$, and
\begin{enumerate}
\item $Au_{k}=\lambda_{k}u_{k}$, $k\in\mathbb{N}$,
\item $A=\sum_{k=1}^{\infty}\lambda_{k}\left|u_{k}\left\rangle \right\langle u_{k}\right|$,
\item $\mbox{spec}\left(A\right)=\left\{ \lambda_{k}\right\} \cup\left\{ 0\right\} $, 
\item $\dim\left\{ v\in\mathscr{H}\:\big|\:Av=\lambda_{k}v\right\} <\infty$,
for all $k\in\mathbb{N}$.
\end{enumerate}

The set $\left\{ u_{k}\right\} $ extends to an ONB, containing an
ONB (possibly infinite) for the subspace
\[
\ker\left(A\right)=\left\{ v\in\mathscr{H}\:\big|\:Av=0\right\} .
\]

\end{thm}
\begin{xca}[An eigenspace]
\myexercise{An eigenspace}\label{exer:eigenspace}Let $A\in\mathscr{B}\left(\mathscr{H}\right)$
be compact, and let $\lambda\in\mathbb{C}\backslash\left\{ 0\right\} $.
Show that the eigenspace
\[
\mathscr{E}_{\lambda}:=\mbox{Ker}\left(\lambda I-A\right)
\]
is finite-dimensional.

\uline{Hint:} Assume the contrary $\dim\mathscr{E}_{\lambda}=\infty$.
Pick an ONB in $\mathscr{E}_{\lambda}$, say $\left\{ u_{i}\right\} _{i\in\mathbb{N}}$.
Since $A$ is compact, the sequence $\left\{ Au_{i}\right\} _{i\in\mathbb{N}}$
has a convergent subsequence, say $\left\{ Au_{i_{k}}\right\} $.
But then 
\[
\left\Vert Au_{i_{k}}-Au_{i_{l}}\right\Vert \xrightarrow[k,l\rightarrow\infty]{}0.
\]
On the other hand,
\[
\left\Vert Au_{i_{k}}-Au_{i_{l}}\right\Vert ^{2}=\left|\lambda\right|^{2}\left\Vert u_{i_{k}}-u_{i_{l}}\right\Vert ^{2}=2\left|\lambda\right|^{2};
\]
so a contradiction.
\end{xca}

\begin{xca}[Attaining the sup]
\myexercise{Attaining the sup}\label{exer:cptsa}Suppose $A\in\mathscr{B}\left(\mathscr{H}\right)$
is compact, and $A^{*}=A$. Suppose further that 
\[
\lambda=\sup\left\{ \left\langle x,Ax\right\rangle \::\:\left\Vert x\right\Vert =1\right\} 
\]
satisfies $\lambda>0$, strict. 
\begin{enumerate}
\item \label{enu:cptsa1}Show that $\left\Vert A\right\Vert =\lambda$.
\item \label{enu:cptsa2}Show that, if $\left\{ x_{i}\right\} _{i\in\mathbb{N}}$
satisfies:
\[
\left\Vert x_{i}\right\Vert =1,\;\left\langle x_{i},Ax_{i}\right\rangle \xrightarrow[i\rightarrow\infty]{}\lambda,
\]
then $\exists$ $x\in\mathscr{H}$, $\left\Vert x\right\Vert =1$,
and a subsequence $\left\{ x_{i_{k}}\right\} $ such that
\[
\left\Vert Ax_{i_{k}}-Ax_{i_{l}}\right\Vert \xrightarrow[k,l\rightarrow\infty]{}0,\;\mbox{and}
\]
\[
\left\langle x_{i_{k}}-x,v\right\rangle \xrightarrow[k\rightarrow\infty]{}0,\;\forall v\in\mathscr{H}.
\]

\item Conclude from (\ref{enu:cptsa1})-(\ref{enu:cptsa2}) that $Ax=\lambda x$. 
\item Make (\ref{enu:cptsa1})-(\ref{enu:cptsa2}) the first step in an
induction; thus finishing the proof of the Spectral Theorem for compact
selfadjoint operators.
\end{enumerate}
\end{xca}
We begin with some preliminaries in the preparation for the proof.

\index{sesquilinear form}
\begin{lem}[polarization identity]
\label{lem:pi}Let $X$ be a set, and $f:X\times X\rightarrow\mathbb{C}$
a sesquilinear form, conjugate linear in the first variable and linear
in the second variable. Then the following polarization identity holds:
\begin{equation}
f\left(x,y\right)=\frac{1}{4}\sum_{k=0}^{3}i^{k}f\left(y+i^{k}x,y+i^{k}x\right)\label{eq:pi}
\end{equation}
for all $x,y\in X$.\end{lem}
\begin{proof}
A direct computation shows that 
\begin{eqnarray*}
f\left(y+x,y+x\right)-f\left(y-x,y-x\right) & = & 2f\left(y,x\right)+2f\left(x,y\right),\;\mbox{and}\\
i\left(f\left(y+ix,y+ix\right)-f\left(y-ix,y-ix\right)\right) & = & -2f\left(y,x\right)+2f\left(x,y\right)
\end{eqnarray*}
Adding the above two equations yields the desired result.\index{operators!compact}\end{proof}
\begin{cor}
\label{cor:xAx}A bounded operator $A$ in $\mathscr{H}$ is selfadjoint
if and only if 
\[
\left\langle x,Ax\right\rangle \in\mathbb{R},\;\forall x\in\mathscr{H}.
\]
\end{cor}
\begin{proof}
Suppose $A$ is selfadjoint, i.e., $\left\langle x,Ay\right\rangle =\left\langle Ax,y\right\rangle $,
$\forall x,y\in\mathscr{H}$. Setting $x=y$, then 
\[
\left\langle x,Ax\right\rangle =\left\langle Ax,x\right\rangle =\overline{\left\langle x,Ax\right\rangle }\Longrightarrow\left\langle x,Ax\right\rangle \in\mathbb{R},\;\forall x\in\mathscr{H}.
\]

Conversely, suppose $\left\langle x,Ax\right\rangle \in\mathbb{R}$,
$\forall x\in\mathscr{H}$. Note that 
\[
\left(x,y\right)\longmapsto\left\langle x,Ay\right\rangle \quad\mbox{and}\quad\left(x,y\right)\longmapsto\left\langle Ax,y\right\rangle 
\]
are both sesquilinear forms defined on $\mathscr{H}$. It follows
from \lemref{pi}, that 
\begin{eqnarray}
\left\langle x,Ay\right\rangle  & = & \frac{1}{4}\sum_{k=0}^{3}i^{k}\left\langle y+i^{k}x,A\left(y+i^{k}x\right)\right\rangle ,\;\mbox{and}\label{eq:pi1}\\
\left\langle Ax,y\right\rangle  & = & \frac{1}{4}\sum_{k=0}^{3}i^{k}\left\langle A\left(y+i^{k}x\right),y+i^{k}x\right\rangle .\label{eq:pi2}
\end{eqnarray}
But by the assumption ($\left\langle x,Ax\right\rangle \in\mathbb{R}$,
$x\in\mathscr{H}$), the RHS in (\ref{eq:pi1}) and (\ref{eq:pi2})
are equal. Therefore, we conclude that $\left\langle x,Ay\right\rangle =\left\langle Ax,y\right\rangle $,
$\forall x,y\in\mathscr{H}$, i.e., $A$ is selfadjoint.
\end{proof}
\index{sesquilinear form}
\begin{thm}
\label{thm:norm}Let $\mathscr{H}$ be a Hilbert space over $\mathbb{C}$.
Let $f:\mathscr{H}\times\mathscr{H}\rightarrow\mathbb{C}$ be a sesquilinear
form, and set 
\[
M:=\sup\left\{ \left|f\left(x,y\right)\right|:\left\Vert x\right\Vert =\left\Vert y\right\Vert =1\right\} <\infty.
\]
Then there exists a unique bounded operator $A$ in $\mathscr{H}$,
satisfying 
\begin{eqnarray}
f\left(x,y\right) & = & \left\langle Ax,y\right\rangle ,\;\forall x,y\in\mathscr{H};\mbox{ and }\label{eq:sf1}\\
\left\Vert A\right\Vert  & = & M.\label{eq:sf2}
\end{eqnarray}
\end{thm}
\begin{proof}
Given $x,y\in\mathscr{H}$, nonzero, we have 
\begin{equation}
\left|f\left(\frac{x}{\left\Vert x\right\Vert },\frac{y}{\left\Vert y\right\Vert }\right)\right|\leq M\Longleftrightarrow\left|f\left(x,y\right)\right|\leq M\left\Vert x\right\Vert \left\Vert y\right\Vert .\label{eq:a-1}
\end{equation}
Thus, for each $x\in\mathscr{H}$, the map $y\mapsto f\left(x,y\right)$
is a bounded linear functional on $\mathscr{H}$. By Riesz's theorem,
there exists a unique element $\xi_{x}\in\mathscr{H}$, such that
\index{Riesz' theorem}\index{Theorem!Riesz-}
\[
f\left(x,y\right)=\left\langle \xi_{x},y\right\rangle .
\]
Set $\xi_{x}:=Ax$, $x\in\mathscr{H}$. (The uniqueness part of follows
from Riesz.)

Note the map $x\longmapsto Ax$ is linear. For if $c\in\mathbb{C}$,
then
\begin{eqnarray}
f\left(x_{1}+cx_{2},y\right) & = & \left\langle A\left(x_{1}+cx_{2}\right),y\right\rangle ,\;\mbox{and}\label{eq:a-4}\\
f\left(x_{1}+cx_{2},y\right) & = & f\left(x_{1},y\right)+\overline{c}f\left(x_{2},y\right)\nonumber \\
 & = & \left\langle Ax_{1},y\right\rangle +\overline{c}\left\langle Ax_{2},y\right\rangle \nonumber \\
 & = & \left\langle Ax_{1}+cAx_{2},y\right\rangle ;\label{eq:a-5}
\end{eqnarray}
where in (\ref{eq:a-5}), we used the fact that $f$ is conjugate
linear in the first variable. It follows that $A\left(x_{1}+cx_{2}\right)=Ax_{1}+cAx_{2}$,
i.e., $A$ is linear. 

Finally, 
\[
\left\Vert A\right\Vert =\sup_{\left\Vert x\right\Vert =1}\left\Vert Ax\right\Vert =\sup_{\left\Vert x\right\Vert =1}\left(\sup_{\left\Vert y\right\Vert =1}\left|\left\langle Ax,y\right\rangle \right|\right)
\]
and (\ref{eq:sf2}) follows.\end{proof}
\begin{cor}
Any bounded operator $A$ in $\mathscr{H}$ is uniquely determined
by the corresponding sesquilinear form $\left(x,y\right)\longmapsto\left\langle x,Ay\right\rangle $,
$\left(x,y\right)\in\mathscr{H}\times\mathscr{H}$. 
\end{cor}

\begin{cor}
\label{cor:normxy}For all $A\in\mathscr{B}\left(\mathscr{H}\right)$,
we have 
\begin{equation}
\left\Vert A\right\Vert =\sup\left\{ \left|\left\langle x,Ay\right\rangle \right|:\left\Vert x\right\Vert =\left\Vert y\right\Vert =1\right\} .\label{eq:a-2}
\end{equation}

\end{cor}

Eq. (\ref{eq:a-3}) below is the key step in the proof of \thmref{spcpt}.

\index{inequality!Cauchy-Schwarz}
\begin{cor}
\label{cor:normxx}Let $A$ be a bounded selfadjoint operator in $\mathscr{H}$,
then
\begin{equation}
\left\Vert A\right\Vert =\sup\left\{ \left|\left\langle x,Ax\right\rangle \right|:\left\Vert x\right\Vert =1\right\} .\label{eq:a-3}
\end{equation}
\end{cor}
\begin{proof}
Set $M:=\sup\left\{ \left|\left\langle x,Ax\right\rangle \right|:\left\Vert x\right\Vert =1\right\} $. 

For all unit vector $x$ in $\mathscr{H}$, we see that 
\[
\left|\left\langle x,Ax\right\rangle \right|\leq\left\Vert x\right\Vert \left\Vert Ax\right\Vert \leq\left\Vert x\right\Vert \left\Vert x\right\Vert \left\Vert A\right\Vert =\left\Vert A\right\Vert ;
\]
where the first step above uses the Cauchy-Schwarz inequality. Thus,
$M\leq\left\Vert A\right\Vert $. 

Conversely, by the polarization identity (\ref{eq:pi}), we have 
\begin{align}
4\left\langle x,Ay\right\rangle = & \left\langle A\left(x+y\right),x+y\right\rangle -\left\langle A\left(-x+y\right),-x+y\right\rangle \nonumber \\
 & +i\left\langle A\left(ix+y\right),ix+y\right\rangle -i\left\langle A\left(-ix+y\right),-ix+y\right\rangle .\label{eq:san}
\end{align}
Sine $A$ is selfadjoint, the four inner products on the RHS of (\ref{eq:san})
are all real-valued (\corref{xAx}). Therefore, 
\[
\Re\left\{ \left\langle x,Ay\right\rangle \right\} =\frac{1}{4}\left(\left\langle A\left(x+y\right),x+y\right\rangle -\left\langle A\left(-x+y\right),-x+y\right\rangle \right).
\]
Now, there exists a phase factor $e^{i\theta}$ (depending on $x,y$)
such that 
\begin{eqnarray*}
\left|\left\langle x,Ay\right\rangle \right| & = & e^{i\theta}\left\langle x,Ay\right\rangle \\
 & = & \left|\Re\left\{ \left\langle x,Ay\right\rangle \right\} \right|\\
 & = & \frac{1}{4}\left|\left\langle A\left(x+y\right),x+y\right\rangle -\left\langle A\left(-x+y\right),-x+y\right\rangle \right|\\
 & \leq & \frac{1}{4}M\left(\left\Vert x+y\right\Vert ^{2}+\left\Vert x-y\right\Vert ^{2}\right)\\
 & = & \frac{1}{4}\left\Vert M\right\Vert \left(2\left\Vert x\right\Vert ^{2}+2\left\Vert y\right\Vert ^{2}\right)=M
\end{eqnarray*}
valid for all unit vectors $x,y$ in $\mathscr{H}$, i.e., with $\left\Vert x\right\Vert =\left\Vert y\right\Vert =1$.
It follows from this and (\ref{eq:a-2}) that $\left\Vert A\right\Vert \leq M$.

Therefore, we have 
\[
M=\left\Vert A\right\Vert =\sup_{\left\Vert x\right\Vert =1}\left|\left\langle x,Ax\right\rangle \right|
\]
which is the assertion.
\end{proof}
\emph{Integral operators} with continuous kernel form an important
subclass of compact operators. 

\uline{Setting}. Let $X$ be a compact space, $\mu$ a finite positive
Borel measure on $X$, and 
\begin{equation}
K:X\times X\longrightarrow\mathbb{C}\label{eq:io1}
\end{equation}
a given function, assumed \emph{continuous} on $X\times X$. Define
\[
T_{K}:L^{2}\left(\mu\right)\longrightarrow L^{2}\left(\mu\right)\;\mbox{by}
\]
\begin{equation}
\left(T_{K}f\right)\left(x\right)=\int_{X}K\left(x,y\right)f\left(y\right)d\mu\left(y\right),\;\forall f\in L^{2}\left(\mu\right),\:\forall x\in X.\label{eq:io2}
\end{equation}

\begin{xca}[An application of Arzelà-Ascoli]
\myexercise{An application of Arzelà-Ascoli} Prove that $T_{K}$
is a compact operator in $L^{2}\left(\mu\right)$ subject to the stated
assumptions above. \index{Arzelà-Ascoli}

\uline{Hint}: 

Step 1. Show that, for $\forall x_{1},x_{2}\in X$, and $f\in L^{2}\left(\mu\right)$,
we have:
\begin{eqnarray*}
 &  & \left|T_{K}f\left(x_{1}\right)-T_{K}f\left(x_{2}\right)\right|\\
 & \leq & \sqrt{\mu\left(X\right)}\:\max_{y\in X}\left|K\left(x_{1},y\right)-K\left(x_{2},y\right)\right|\:\left\Vert f\right\Vert _{L^{2}\left(\mu\right)}.
\end{eqnarray*}

Step 2. Show that 
\[
\left|\left(T_{K}f\right)\left(x\right)\right|\leq\sqrt{\mu\left(X\right)}\:\max_{y\in X}\left|K\left(x,y\right)\right|\:\left\Vert f\right\Vert _{L^{2}\left(\mu\right)}.
\]

Step 3. Conclude from steps 1-2, and an application of Arzelà-Ascoli's
theorem that $T_{K}:L^{2}\left(\mu\right)\longrightarrow L^{2}\left(\mu\right)$
is a compact operator.
\end{xca}

\begin{xca}[Powers-Størmer \cite{MR0269230}]
\label{exer:powers}\myexercise{Powers-St\o mer}Let $\mathscr{H}$
be a Hilbert space, and let $A$ and $B$ be positive operators ($\in\mathscr{B}\left(\mathscr{H}\right)$).
Then show that 
\begin{equation}
\left\Vert A^{\frac{1}{2}}-B^{\frac{1}{2}}\right\Vert _{HS}^{2}\leq\left\Vert A-B\right\Vert _{TR}.\label{eq:ps1}
\end{equation}
(Note $A^{\frac{1}{2}}=\sqrt{A}$ is defined via the Spectral Theorem.)\index{Theorem!Spectral-}

\uline{Hint}: Difficult. (The inequality (\ref{eq:ps1}) is called
the Powers-Størmer inequality.) Set $S=A^{\frac{1}{2}}-B^{\frac{1}{2}}$,
and $T=A^{\frac{1}{2}}+B^{\frac{1}{2}}$. 

Note that (\ref{eq:ps1}) is trivial if $A-B$ is not trace-class,
so assume it has finite trace-norm. Then diagonalize $S$ in an ONB
(use the Spectral Theorem), i.e., pick an ONB $\left\{ f_{i}\right\} $
of eigenvectors with eigenvalues $\left\{ \lambda_{i}\right\} $,
$Sf_{i}=\lambda_{i}f_{i}$. Then\index{eigenvalue}
\begin{eqnarray*}
Tr\left(\left|A-B\right|\right) & = & \frac{1}{2}\sum_{i}\left\langle f_{i},\left|ST+TS\right|f_{i}\right\rangle \\
 & \geq & \sum_{i}\left|\lambda_{i}\left\langle f_{i},Tf_{i}\right\rangle \right|\\
 & \geq & \sum_{i}\left|\lambda_{i}\right|^{2}.
\end{eqnarray*}

\end{xca}

\section*{A summary of relevant numbers from the Reference List}

For readers wishing to follow up sources, or to go in more depth with
topics above, we suggest: \cite{Con90,FL28,Kat95,Kr55,MR1892228,LP89,Ne69,RS75,Rud73,Sto51,Sto90,Yos95,HJL13,DHL09,MR2367342,MR1839648,MR2354721,Fan10,MR2254502,AG93,DS88b,Hal67,MR1895530,Mac52,VN35,Hel13,MJD15}.

\part{\label{part:app}Applications}

\chapter{GNS and Representations\label{chap:GNS}}
\begin{quotation}
If one finds a difficulty in a calculation which is otherwise quite
convincing, one should not push the difficulty away; one should rather
try to make it the centre of the whole thing.

--- Werner Heisenberg\sindex[nam]{Heisenberg, W., (1901-1976)}\vspace{1em}\\
\textquotedblleft Mathematics is a way of thinking in everyday life
\dots \textquotedblright{} 

--- I.M. Gelfand\sindex[nam]{Gelfand, I.M., (1913-2009)}\vspace{1em}\\
``First quantization is a mystery; -- second quantization is a functor.'' 

--- Edward Nelson\vspace{0.5em}\\
Explanation:

A \textbf{category} is an algebraic structure that comprises \textquotedbl{}objects\textquotedbl{}
linked by morphisms, also called \textquotedbl{}arrows\textquotedbl{}.
A category has two basic properties: one allows us to compose the
arrows associatively; and the existence of an identity arrow for each
object. A functor is a type of mapping, or transformation, between
categories. Functors can be thought of as transformation between categories
that transform the rules in the first category into those of the second:
objects to objects, arrows to arrows, and diagrams to diagrams. In
small categories, functors can be thought of as morphisms.

First quantization is the replacement of classical observables, such
as energy or momentum, by operators; and of classical states by \textquotedbl{}wave
functions\textquotedbl{}. Second quantization usually refers to the
introduction of field operators when describing quantum many-body
systems. In second quantization, one passes from wave functions to
an operators; hence the non-commutativity.

Expanding upon the more traditional interpretation, Ed Nelson suggested
yet another second quantization functor; it goes from the category
of Hilbert space (Hilb) to that of probability space (Prob). In this
version, the objects in the category of Hilbert space \textquotedblleft Hilb\textquotedblright{}
are Hilbert spaces, and the morphisms are contractive linear operators.
In the category \textquotedblleft Prob,\textquotedblright{} the morphisms
are point-transformation of measures.\vspace{2em}
\end{quotation}

A state on a $C^{*}$-algebra $\mathfrak{A}$ is a positive linear
(and normalized) functional on $\mathfrak{A}$. Given a $C^{*}$-algebra
$\mathfrak{A}$, then there is a bijective correspondence between
states of $\mathfrak{A}$, on one side, and cyclic representations
of $\mathfrak{A}$ on the other; it is called the Gelfand-Naimark-Segal
construction (abbreviated GNS), and it yields an explicit correspondence
between the set of all cyclic $*$-representations of $\mathfrak{A}$,
$Rep_{cyc}\left(\mathfrak{A}\right)$; and the states of $\mathfrak{A}$,
$S(\mathfrak{A})$. It is named for Israel Gelfand, Mark Naimark,
and Irving Segal. 

The importance of the GNS construction is that it offers answers to
a host of questions about representations of algebras, and of groups;
the natural question regarding elementary building blocks.

More precisely, in the case of unitary representations of groups,
the first questions that present themselves are: \textquotedblleft What
is the \textquotedblleft right\textquotedblright{} notion of decomposition
of a given unitary representation in terms of irreducible unitary
representations of $G$?\textquotedblright{} And \textquotedblleft how
to compute decompositions?\textquotedblright{} In this generality,
there is not a precise answer. But if $G$ is assumed locally compact
and unimodular, I.E. Segal \cite{Seg50} established precise answers.
They entail direct integral theory for representations; see details
below. The power of the GNS-construction (states vs representations)
is that it allows us to answer a parallel question for states. Indeed,
there is a precise notion of \textquotedblleft building blocks\textquotedblright{}
for states, they are the pure states (extreme points).

Segal's insight was twofold: (i) Showing that the pure states correspond
to irreducible representations via the GNS correspondence. And (ii),
make the precise link between a direct integral decomposition for
states, on one side, and on the other, direct integrals for unitary
representations. Part (ii) in turn involves Krein-Milman and Choquet
theory\index{Theorem!Krein-Milman's-}; see Sections \ref{sec:gns},
\ref{sec:KMil}, and \chapref{KS}. \index{Theorem!Choquet-}

A corollary of the GNS construction is the Gelfand-Naimark theorem.
The latter characterizes $C^{*}$-algebras as precisely the norm-closed
$*$-algebras arising as $*$-subalgebras of $\mathscr{B}\left(\mathscr{H}\right)$,
the $C^{*}$-algebra of all bounded operators on a Hilbert space.
By extreme-point theory \cite{Phe01}, one shows that every $C^{*}$-algebra
has sufficiently many pure states (corresponding to irreducible representations
under GNS). As a result, the representation of $\mathfrak{A}$ arising
as a direct sum of these corresponding irreducible GNS-representations
is faithful.

Since states in quantum physics are vectors (of norm one) in Hilbert
space; two question arise: \textquotedblleft Where does the Hilbert
space come from?\textquotedblright{} And \textquotedblleft What are
the algebras of operators from which the selfadjoint observables must
be selected?\textquotedblright{} In a general framework, we offer
an answer below, it goes by the name \textquotedblleft the Gelfand-Naimark-Segal
(GNS) theorem, \textquotedblright{} which offers a direct correspondence
between states and cyclic representations.

\index{Gelfand-Naimark-Segal}

\index{observable}

Chapters \ref{chap:GNS} and \ref{chap:cp} below form a pair; --
to oversimplify, the theme in \chapref{cp} is a generalization of
that of the present chapter.

One could say that  \chapref{GNS} is about scalar valued \textquotedblleft states\textquotedblright ;
while the \textquotedblleft states\textquotedblright{} in \chapref{cp}
are operator valued. In both cases, we must specify the appropriate
notion of positivity, and this notion in the setting of \chapref{cp}
is more subtle; -- it is called \textquotedblleft complete positivity.\textquotedblright{}

But the goal in both cases is to induce in order to create representation
of some given non abelian algebra $\mathfrak{A}$ coming equipped
with a star-involution; for example a $C^{*}$-algebra. The representations,
when induced from states, will be $*$-representations; i.e., will
take the $*$-involution in $\mathfrak{A}$ to \textquotedblleft adjoint
operator\textquotedblright{} -- where \textquotedblleft adjoint\textquotedblright{}
refers to the Hilbert space of the induced representation.

In Chapters \ref{chap:GNS}-\ref{chap:cp}, this notion of induction
is developed in detail; and its counterpart for the case of unitary
representations of groups is discussed in detail in \chapref{groups}.

Historically, the two notion of \emph{induction of representations}
were used by researchers in parallel universes, for the case of operator
algebras (Chapters \ref{chap:GNS}-\ref{chap:cp}), they were pioneered
by Gelfand, Naimark, Segal, and brought to fruition by Stinespring
and Arveson. 

On the other side of the divide, in the study of unitary representations
of groups, the names are Harish-Chandra, G.W. Mackey (and more, see
cited references in \chapref{groups}); -- and this is the subject
of  \chapref{groups} below. We caution the reader that the theory
of unitary representations is a vast subject, and motivated by a number
of diverse areas, such as quantum theory, ergodic theory, harmonic
analysis, to mention just a few. And the theory of representations
of groups, and their induction, is in turn developed by different
researchers; and often with different groups $G$ in mind; -- continuous
vs discrete; Lie groups vs the more general case of locally compact
groups. The case when the group $G$ is assumed locally compact is
attractive because we then always will have left (or right-) Haar
measure at our disposal. And there is an associated left-regular representation
on the $L^{2}$ space of Haar measure. The left-invariance of Haar
measure makes this representation unitary. The analogous hold of course
for the constructions with right-invariant Haar measure; the two are
linked by the \emph{modular function} of $G$.

\index{representation!unitary}

\section{\label{sec:ov}Definitions and Facts: An Overview}

Let $\mathfrak{A}$ be an algebra over $\mathbb{C}$, with an involution
$\mathfrak{A}\ni a\mapsto a^{*}\in\mathfrak{A}$, and the unit-element
$\mathbf{1}$. Let $\mathfrak{A}_{+}$ denote the set of positive
elements in $\mathfrak{A}$; i.e.,\index{representation!GNS} 
\[
\mathfrak{A}_{+}=\left\{ b^{*}b\;\big|\;b\in\mathfrak{A}\right\} .
\]

\begin{defn}
We say $\mathfrak{A}$ is a $C^{*}$-algebra\index{algebras!$C^{*}$-algebra}
if it is complete in a norm $\left\Vert \cdot\right\Vert $, which
satisfies: 
\begin{enumerate}
\item $\left\Vert ab\right\Vert \leq\left\Vert a\right\Vert \left\Vert b\right\Vert $,
$\forall a,b\in\mathfrak{A}$;
\item $\left\Vert \mathbf{1}\right\Vert =\mathbf{1}$; 
\item \label{enu:c3}$\left\Vert b^{*}b\right\Vert =\left\Vert b\right\Vert ^{2}$,
$\forall b\in\mathfrak{A}$.\footnote{Kadison et al. in 1950's reduced the axioms of $C^{*}$-algebra from
about 6 down to just one (\ref{enu:c3}) on the $C^{*}$-norm.}
\end{enumerate}
\end{defn}
\begin{example}
Let $X$ be a compact Hausdorff space, the algebra $C\left(X\right)$
of all continuous function on $X$ is a $C^{*}$-algebra under the
sup-norm. 
\end{example}

\begin{example}
$\mathscr{B}\left(\mathscr{H}\right)$: all bounded linear operators
on a fixed Hilbert space $\mathscr{H}$ is a $C^{*}$-algebra.
\end{example}

\begin{example}
\label{exa:cuntzN}$\mathscr{O}_{N}$: the Cuntz-algebra, $N>1$;
it is the $C^{*}$-completion of $N$ generators $s_{1},s_{2},\ldots,s_{N}$
satisfying the following relations \cite{MR0467330}: \index{completion!$C^{*}$-}
\begin{enumerate}
\item \label{enu:On1}$s_{i}^{*}s_{j}=\delta_{ij}\mathbf{1}$;
\item \label{enu:On2}$\sum_{i=1}^{N}s_{i}s_{i}^{*}=\mathbf{1}$. 
\end{enumerate}

For the representations of $\mathscr{O}_{N}$, see \cite{Gli60,Gli61,BJO04}.

\end{example}
\index{algebras!Cuntz algebra}\index{Cuntz-algebra}\index{representation!- of the Cuntz algebra}

\begin{defn}
We denote $Rep\left(\mathfrak{A},\mathscr{H}\right)$ the representations
of $\mathfrak{A}$ acting on some Hilbert space $\mathscr{H}$, i.e.,
$\pi\in Rep\left(\mathfrak{A},\mathscr{H}\right)$ iff $\pi:\mathfrak{A}\rightarrow\mathscr{B}\left(\mathscr{H}\right)$
is a homomorphism\index{homomorphism} of $^{*}$-algebras, $\pi(\mathbf{1})=I_{\mathscr{H}}=$
the identity operator in $\mathscr{H}$; in particular
\begin{equation}
\left\langle \pi\left(b\right)u,v\right\rangle _{\mathscr{H}}=\left\langle u,\pi\left(b^{*}\right)v\right\rangle _{\mathscr{H}},\;\forall b\in\mathfrak{A},\forall u,v\in\mathscr{H}.\label{eq:r1}
\end{equation}
Let $S\left(\mathfrak{A}\right)$ be the states $\varphi:\mathfrak{A}\rightarrow\mathbb{C}$
on $\mathfrak{A}$; i.e., (axioms) $\varphi\in\mathfrak{A}^{*}=$
the dual of $\mathfrak{A}$, $\varphi\left(\mathbf{1}\right)=1$,
and\index{space!dual-}\index{axioms} 
\begin{equation}
\varphi\left(b^{*}b\right)\geq0,\:\forall b\in\mathfrak{A}.\label{eq:r2}
\end{equation}
\end{defn}
\begin{thm}[Gelfand-Naimark-Segal (GNS)]
\label{thm:GNS} There is a bijection:\index{Gelfand-Naimark-Segal}
\[
S\left(\mathfrak{A}\right)\longleftrightarrow\boxed{\mbox{cyclic representations, up to unitary equivalence}}
\]
as follows: 

$\longleftarrow$ (easy direction): Given $\pi\in Rep\left(\mathfrak{A},\mathscr{H}\right)$,
$u_{0}\in\mathscr{H}$, $\left\Vert u_{0}\right\Vert =1$, set 
\begin{equation}
\varphi\left(a\right)=\left\langle u_{0},\pi\left(a\right)u_{0}\right\rangle _{\mathscr{H}},\;\forall a\in\mathfrak{A}.\label{eq:r3}
\end{equation}

$\longrightarrow$ (non-trivial direction): Given $\varphi\in S\left(\mathfrak{A}\right)$,
there is a system $\left(\pi,\mathscr{H},u_{0}\right)$ such that
(\ref{eq:r3}) holds. (Notation, we set $\pi=\pi_{\varphi}$ to indicate
the state $\varphi$.) 
\end{thm}
\index{sesquilinear form}
\begin{proof}
$\left(\longrightarrow\right)$ Given $\varphi\in S\left(\mathfrak{A}\right)$,
then on $\mathfrak{A}\times\mathfrak{A}$ consider the sesquilinear
form
\begin{equation}
\left(a,b\right)\longmapsto\varphi\left(a^{*}b\right)\label{eq:s1-1}
\end{equation}
\[
\mathscr{H}_{\varphi}=\left\{ \mathfrak{A}/\left\{ b\in\mathfrak{A}\;\big|\;\varphi\left(b^{*}b\right)=0\right\} \right\} ^{\sim}
\]
where $\sim$ refers to Hilbert completion in (\ref{eq:s1-1}). Note
\[
\left|\varphi\left(a^{*}b\right)\right|^{2}\leq\varphi\left(a^{*}a\right)\varphi\left(b^{*}b\right),\!\forall a,b\in\mathfrak{A}.
\]
Set $\Omega=\mbox{class}\left(\mathbf{1}\right)$ in $\mathscr{H}_{\varphi}$,
and \index{completion!Hilbert-} 
\[
\pi_{\varphi}\left(a\right)\left(\mbox{class}\left(b\right)\right)=\mbox{class}\left(ab\right),\;\forall a,b\in\mathfrak{A};\;\left(\mbox{Schwarz.}\right)
\]
Then it is easy to show that $\left(\mathscr{H}_{\varphi},\Omega,\pi_{\varphi}\right)$
satisfies conclusion (\ref{eq:r3}), i.e., 
\[
\varphi\left(a\right)=\left\langle \Omega,\pi\left(a\right)\Omega\right\rangle _{\mathscr{H}_{\varphi}},\;\forall a\in\mathfrak{A}.
\]
\end{proof}
\begin{defn}
Let $\varphi\in S\left(\mathfrak{A}\right)$, we say it is a \emph{pure
state} iff $\varphi\in\mbox{ext}S\left(\mathfrak{A}\right)$ := the
extreme-points in $S\left(\mathfrak{A}\right)$. (See \cite{Phe01}.)
\end{defn}
\index{state!pure-}\index{pure state}
\begin{rem}[GNS-correspondence]

\begin{enumerate}
\item If $\mathfrak{A}$ is a $C^{*}$-algebra , then $S\left(\mathfrak{A}\right)\left(\subset\mathfrak{A}^{*}\right)$
is convex\index{convex} and weak {*}-compact\index{weak{*}-compact}\index{compact}. 
\item Given $\varphi\in S\left(\mathfrak{A}\right)$, and let $\pi_{\varphi}\in Rep\left(\mathfrak{A},\mathscr{H}\right)$
be the GNS-representation, see (\ref{eq:r3}); then
\begin{eqnarray*}
 & \boxed{\varphi\in\mbox{ext}S\left(\mathfrak{A}\right),\;\mbox{i.e., it is pure}}\\
 & \Updownarrow\\
 & \boxed{\pi_{\varphi}\:\mbox{is an irreducible representation}}
\end{eqnarray*}

\item If $\psi\in S\left(\mathfrak{A}\right)$, $\exists$ a measure $P_{\psi}$
on $\mbox{ext}\left(S\left(\mathfrak{A}\right)\right)$ such that
$\psi=\int w\:dP_{\psi}\left(w\right)$, and then 
\[
\pi_{\psi}=\int^{\oplus}\pi_{w}\:dP_{\psi}\left(w\right).
\]

\end{enumerate}
\end{rem}
\begin{example}[Pure states, cases where the full list is known!]
 $\varphi\in S\left(\mathfrak{A}\right)$: 
\end{example}
\renewcommand{\arraystretch}{2}

\begin{table}[H]
\begin{tabular}{|>{\centering}p{0.1\columnwidth}|>{\centering}p{0.8\columnwidth}|}
\hline 
$\mathfrak{A}$ & $\mbox{ext}S\left(\mathfrak{A}\right)$\tabularnewline
\hline 
$C\left(X\right)$ & points $x\in X$, and $\varphi=\delta_{x}$ (Dirac mass); $\varphi\left(f\right)=f\left(x\right)$,
$\forall f\in C\left(X\right)$\tabularnewline
\hline 
$\mathscr{B}\left(\mathscr{H}\right)$ & $v\in\mathscr{H}$, $\left\Vert v\right\Vert =1$, $\varphi=\varphi_{v}$;
$\varphi_{v}\left(A\right)=\left\langle v,Av\right\rangle $, $\forall A\in\mathscr{B}\left(\mathscr{H}\right)$\tabularnewline
\hline 
$\mathscr{O}_{N}$ & Partial list: $u=\left(u_{1},\ldots,u_{N}\right)\in\mathbb{C}^{N}$,
$\sum_{1}^{N}|u_{j}|^{2}=1$, $\varphi=\varphi_{u}$,

\medskip{}
specified by $\varphi\left(s_{i}s_{j}^{*}\right)=u_{i}\overline{u_{j}}$,
$\forall i,j=1,\ldots,N$; see \ref{enu:On1}-\ref{enu:On2}. \medskip{}
\tabularnewline
\hline 
\end{tabular}

\protect\caption{\label{tab:pstate}Examples of pure states.}
\end{table}

\renewcommand{\arraystretch}{1}

\index{state!pure-}
\begin{xca}[Irreducible representations]
\myexercise{Irreducible representations}\label{exer:purestate}Using
GNS, write down explicitly the irreducible representations of the
three $C^{*}$-algebras in  \tabref{pstate} corresponding to the
listed pure states. 

\uline{Hint}: In the case of $C\left(X\right)$, the representations
are one-dimensional, but in the other cases, they are infinite-dimensional,
i.e., $\dim\mathscr{H}_{\pi_{\varphi}}=\infty$. \index{pure state}
\index{state!pure-}\end{xca}
\begin{rem}
It is probably impossible to list \uline{all} pure states of $\mathscr{O}_{N}$;
see \cite{Gli60}.\end{rem}
\begin{xca}[Infinite-product measures and representations of $\mathscr{O}_{N}$]
\myexercise{Infinite-product measures and representations of $\mathscr{O}_{N}$}\label{exer:ON}
Fix $N\in\mathbb{N}$, $N>1$, and denote the cyclic group of order
$N$, 
\[
\mathbb{Z}_{N}=\mathbb{Z}/N\mathbb{Z}=\left\{ 0,1,2\ldots,N-1\right\} ,
\]
residue classes mod $N$. Let $\left(z_{i}\right)_{i=0}^{N-1}$ be
complex numbers such that $\sum_{i}\left|z_{i}\right|^{2}=1$, and
assume $z_{j}\neq0$ for all $j$; see line 3 in  \tabref{pstate}. 

Let $p$ be the probability measure on $\mathbb{Z}_{N}$ with weights
$p_{i}=\left|z_{i}\right|^{2}$, and let $\mu=\mu_{p}$ be the infinite-product
measure on $\Omega_{N}:=\vartimes_{\mathbb{N}}\mathbb{Z}_{N}=\vartimes_{\mathbb{N}}\left\{ 0,1,\cdots,N-1\right\} $,
\begin{equation}
\mu_{p}:=\vartimes_{\mathbb{N}}p=\underset{\aleph_{0}-\text{infinite}}{\underbrace{p\times p\times\cdots}}.\label{eq:d1-1}
\end{equation}
Set $\Omega_{N}\left(j\right)=\left\{ \left(x_{i}\right)\in\Omega_{N}\:;\:x_{1}=j\right\} =\left\{ j\right\} \times\Omega_{N}.$
\begin{enumerate}
\item \label{enu:d1-1}For all $x=\left(x_{1},x_{2},x_{3},\ldots\right)$,
and all $f\in L^{2}\left(\mu\right)$, set 
\[
\left(S_{j}f\right)\left(x\right)=\frac{1}{z_{j}}\chi_{\Omega_{N}\left(j\right)}\left(x\right)f\left(x_{2},x_{3},\ldots\right).
\]
Show that the adjoint operator with respect to $L^{2}\left(\mu\right)$
is
\[
\left(S_{j}^{*}f\right)\left(x\right)=z_{j}f\left(j,x_{1},x_{2},x_{3},\ldots\right);
\]
and that this system $\left\{ S_{j}\right\} _{j=0}^{N-1}$ defines
an irreducible representation of $\mathscr{O}_{N}$, i.e., is in $Rep_{irr}\left(\mathscr{O}_{N},L^{2}\left(\mu\right)\right)$.\index{operators!adjoint-}
\item Denote the representation in (\ref{enu:d1-1}) $\pi_{p}^{\left(N\right)}$,
and setting $\mathbbm{1}$ to be the constant function in $L^{2}\left(\mu\right)$,
show that we recover the pure state from line 3 in  \tabref{pstate}
corresponding to $u_{j}=z_{j}$; i.e., using the formula: 
\[
\varphi\left(s_{j}s_{k}^{*}\right)=\left\langle \mathbbm{1},S_{j}S_{k}^{*}\mathbbm{1}\right\rangle _{L^{2}\left(\mu\right)}=z_{j}\overline{z_{k}},\;\forall j,k\in\mathbb{Z}_{N}.
\]

\item Show that $\pi_{p}^{\left(N\right)}$ is \uline{not} irreducible
when restricted to the abelian subalgebra in $\mathscr{O}_{N}$ generated
by $\left\{ S_{J}S_{J}^{*}\right\} $, as $J$ ranges over all finite
words in the fixed alphabet $\mathbb{Z}_{N}$.
\end{enumerate}
\end{xca}

\begin{xca}[A representation of $\mathscr{O}_{N}$]
\myexercise{A representation of $\mathscr{O}_{N}$}\label{exer:ON1}What
can you say about the representation of $\mathscr{O}_{N}$ corresponding
to $\left(0,z_{1},\ldots,z_{N-1}\right)$, $\sum\left|z_{j}\right|^{2}=1$?

\uline{Hint}: Modify (\ref{enu:d1-1}) from  \exerref{ON}. (The
state $\varphi\left(s_{j}s_{k}^{*}\right)=z_{j}\overline{z_{k}}$
then yields $\varphi\left(s_{j}s_{0}^{*}\right)=0$, $\forall j\in\mathbb{Z}_{N}$.
)
\end{xca}

\subsection*{\noindent Groups}

Case 1. Groups contained in $\mathscr{B}\left(\mathscr{H}\right)$
where $\mathscr{H}$ is a fixed Hilbert space: 
\begin{defn}
Set 
\begin{itemize}
\item $\mathscr{B}\left(\mathscr{H}\right)^{-1}$: all bounded linear operators
in $\mathscr{H}$ with bounded inverse.
\item $\mathscr{B}\left(\mathscr{H}\right)_{uni}$: all unitary operators
$u:\mathscr{H}\rightarrow\mathscr{H}$, i.e., $u$ satisfies 
\[
uu^{*}=u^{*}u=I_{\mathscr{H}}.
\]

\end{itemize}
\end{defn}

\index{operators!unitary}
\begin{defn}
Fix a group $G$, and set: \index{homomorphism}
\begin{itemize}
\item $Rep\left(G,\mathscr{H}\right)$: all homomorphisms $\rho\in G\rightarrow\mathscr{B}\left(\mathscr{H}\right)^{-1}$
\item $Rep_{uni}\left(G,\mathscr{H}\right)$: all homomorphisms, $\rho:G\rightarrow\mathscr{B}_{uni}\left(\mathscr{H}\right)$,
i.e., 
\[
\rho(g^{-1})=\rho\left(g\right){}^{-1}=\rho\left(g\right)^{*},\:\forall g\in G
\]

\item $Rep_{cont}\left(G,\mathscr{H}\right)$: Elements $\rho\in Rep\left(G,\mathscr{H}\right)$
such that $\forall v\in\mathscr{H}$, 
\[
G\ni g\mapsto\rho\left(g\right)v
\]
is continuous from $G$ into $\mathscr{H}$; called \emph{strongly
continuous}. \index{strongly continuous}
\end{itemize}
\end{defn}
\begin{rem}
In the case of $Rep_{cont}\left(G,\mathscr{H}\right)$ it is assumed
that $G$ is a continuous group, i.e., is equipped with a topology
such that the following two operations are both continuous:
\begin{enumerate}
\item $G\times G\ni\left(g_{1},g_{2}\right)\longmapsto g_{1}g_{2}\in G$
\item $G\ni g\longmapsto g^{-1}\in G$ 
\end{enumerate}
\end{rem}
\begin{xca}[The regular representation of $G$]
\myexercise{The regular representation of $G$}\label{exer:lcg}Let
$G$ be a locally compact group with $\mu$ = a left-invariant Haar
measure. Set 
\[
\left(\rho_{L}\left(g\right)f\right)\left(x\right):=f\left(g^{-1}x\right),\;g,x\in G,\;f\in L^{2}\left(G,\mu\right).
\]
Then show that $\rho_{L}$ is a strongly continuous unitary representation
of $G$ acting in $L^{2}\left(G,\mu\right)$. 
\end{xca}
\index{representation!unitary}

\index{representation!strongly continuous}

\noindent \begin{flushleft}
\textbf{The Group Algebra}
\par\end{flushleft}

Let $G$ be a group, and set $\mathbb{C}[G]:=$ all linear combinations,
i.e., finite sums\index{algebras!group algebra} 
\begin{equation}
A=\sum_{g}A_{g}g\label{eq:r4}
\end{equation}
where $A_{g}\in\mathbb{C}$, and making $\mathfrak{A}:=\mathbb{C}[G]$
into a $*$-algebra with the following two operations on finite sums
as in (\ref{eq:r4}): $\mathfrak{A}\times\mathfrak{A}\longrightarrow\mathfrak{A}$,
given by 
\begin{equation}
\left(\sum_{g\in G}A_{g}g\right)\left(\sum_{h\in G}B_{h}h\right):=\sum_{g}\left(\sum_{hk=g}A_{h}B_{k}\right)g\label{eq:r5}
\end{equation}
and 
\begin{equation}
\left(\sum_{g\in G}A_{g}g\right)^{*}:=\sum_{g\in G}\overline{A_{g}}g^{-1}.\label{eq:r6}
\end{equation}

\begin{lem}
\label{lem:reprep}There is a bijection between $Rep_{uni}\left(G,\mathscr{H}\right)$
and $Rep\left(\mathbb{C}[G],\mathscr{H}\right)$ as follows: If $\pi\in Rep_{uni}\left(G,\mathscr{H}\right)$,
set $\widetilde{\pi}\in Rep\left(\mathbb{C}[G],\mathscr{H}\right)$:
\begin{equation}
\widetilde{\pi}\left(\sum_{g\in G}A_{g}g\right):=\sum_{g\in G}A_{g}\pi\left(g\right)\label{eq:r7}
\end{equation}
where the element $\sum_{g\in G}A_{g}g$ in (\ref{eq:r7}) is a generic
element in $\mathbb{C}[G]$, see (\ref{eq:r4}), i.e., is a finite
sum with $A_{g}\in\mathbb{C}$, for all $g\in G$.\end{lem}
\begin{xca}[Unitary representations]
\myexercise{Unitary representations}\label{exer:groupalg}Fill in
the proof details of the assertion in \lemref{reprep}.
\end{xca}
\index{representation!unitary}
\begin{example}
Let $G$ be a group, considered as a countable discrete group (the
countability is not important). Set $\mathscr{H}=l^{2}\left(G\right)$,
and 
\begin{equation}
\pi\left(g\right)\delta_{h}:=\delta_{gh},\;\forall g,h\in G.\label{eq:r8}
\end{equation}
\end{example}
\begin{xca}[A proof detail]
\myexercise{A proof detail}\label{exer:groupalg1}Show that $\pi$
in (\ref{eq:r8}) is in $Rep\left(G,l^{2}\left(G\right)\right)$. \end{xca}
\begin{defn}
Let $G$, and $\pi\in Rep\left(G,l^{2}\left(G\right)\right)$ be as
in (\ref{eq:r8}), and let $\widetilde{\pi}\in Rep\left(\mathbb{C}[G],l^{2}\left(G\right)\right)$
be the corresponding representation of $\mathbb{C}[G]$; see \lemref{reprep}.
Set 
\[
C_{red}^{*}\left(G\right):=\mbox{the norm closure of }\widetilde{\pi}\left(\mathbb{C}[G]\right)\subset\mathscr{B}\left(l^{2}\left(G\right)\right);
\]
then $C_{red}^{*}\left(G\right)$ is called the reduced $C^{*}$-algebra
of the group $G$. \end{defn}
\begin{xca}[Reduced $C^{*}$-algebra]
\myexercise{Reduced $C^{*}$-algebra}\label{exer:red}Prove that
$C_{red}^{*}\left(G\right)$ is a $C^{*}$-algebra. \end{xca}
\begin{rem}
It is known \cite{MR0374334} that $C_{red}^{*}\left(F_{2}\right)$
is simple, where $F_{2}$ is the free group on two generators. (``red''
short for reduced; it is called the reduced $C^{*}$-algebra on the
group.)
\end{rem}

\section{\label{sec:gns}The GNS Construction}

The GNS construction is a general principle for getting representations
from given data in applications, especially in quantum mechanics \cite{MR1939631,MR965583,MR675039}.
It was developed independently by I. Gelfand, M. Naimark, and I. Segal
around the 1960s, see e.g., \cite{MR0112604,Seg50}.
\begin{defn}
Let $\mathfrak{A}$ be a $*$-algebra with identity. A representation
of $\mathfrak{A}$ is a map $\pi:\mathfrak{A}\rightarrow B(\mathscr{H}_{\pi})$,
where $\mathscr{H}_{\pi}$ is a Hilbert space, such that for all $A,B\in\mathfrak{A}$,
\index{representation!of algebra}\index{space!Hilbert-}
\begin{enumerate}
\item $\pi(AB)=\pi(A)\pi(B)$
\item $\pi(A^{*})=\pi(A)^{*}$
\end{enumerate}

The $*$ operation (involution) is given on $\mathfrak{A}$ so that
$A^{**}=A$, $(AB)^{*}=B^{*}A^{*}$, $\left(\lambda A\right)^{*}=\overline{\lambda}A^{*}$,
for all $\lambda\in\mathbb{C}$.

\end{defn}
\begin{example}
The multiplication version of the spectral theorem of a single selfadjoint
operator, say $A$ acting on $\mathscr{H}$, yields a representation
of the algebra of $L^{\infty}\left(sp\left(A\right)\right)$ (or $C(sp\left(A\right))$)
as operators on $\mathscr{H}$, where 
\[
L^{\infty}\left(sp\left(A\right)\right)\ni f\xrightarrow{\;\pi\;}f\left(A\right)\in\mathscr{B}(\mathscr{H})
\]
such that $\pi(fg)=\pi(f)\pi(g)$ and $\pi(\bar{f})=\pi(f)^{*}$.
\index{selfadjoint operator}
\end{example}
The general question is given any $*$-algebra, where to get such
a representation? The answer is given by \emph{states}. One gets representations
from algebras vis states. For abelian algebras, the states are Borel\index{measure!Borel}
measures, so the measures come out as a corollary of representations. 
\begin{defn}
Let $\mathfrak{A}$ be a $*$-algebra. A state on $\mathfrak{A}$
is a linear functional\index{functional} $\varphi:\mathfrak{A}\rightarrow\mathbb{C}$
such that $\varphi(1_{\mathfrak{A}})=1$, and $\varphi(A^{*}A)\geq0$,
for all $A\in\mathfrak{A}$. \end{defn}
\begin{example}
Let $\mathfrak{A}=C(X)$, i.e., $C^{*}$-algebra of continuous functions
on a compact Hausdorff space $X$. Note that there is a natural involution
$f\mapsto f^{*}:=\overline{f}$ by complex conjugation. Let $\mu_{\varphi}$
be a Borel probability measure on $X$, then 
\[
C\left(X\right)\ni f\mapsto\varphi(f)=\int_{X}fd\mu_{\varphi}
\]
is a state. In fact, in the abelian case, all states are Borel probability\index{measure!probability}
measures.\index{space!state-}

Because of this example, we say that the GNS construction is non-commutative
measure theory.
\end{example}

\begin{example}
Let $G$ be a discrete group, and let $\mathfrak{A}=\mathbb{C}\left[G\right]$
be the group-algebra, see \secref{ov}. 

If we make the assumption (defining $\varphi$ first on points in
$G$)
\begin{equation}
\varphi\left(g\right)=\begin{cases}
1 & \quad\mbox{if}\:g=e\:\left(\mbox{the unit element in}\:G\right)\\
0 & \quad\mbox{if}\:g\in G\backslash\left\{ e\right\} ,
\end{cases}\label{eq:gg1}
\end{equation}
then the argument from above shows that $\varphi$ extends to a \emph{linear
functional} on $\mathfrak{A}$. \end{example}
\begin{xca}[{The trace state on $\mathbb{C}[G]$}]
\myexercise{The trace state on $\mathbb{C}[G]$}~
\begin{enumerate}
\item Show that $\varphi$ as defined in (\ref{eq:gg1}), extended to $\mathfrak{A}=\mathbb{C}[G]$
is a state, and if $A=\sum_{g}A_{g}g$ (finite sum), then
\begin{equation}
\varphi\left(A^{*}A\right)=\sum_{g\in G}\left|A_{g}\right|^{2};\label{eq:gg2}
\end{equation}
and moreover (the trace property):
\begin{equation}
\varphi\left(AB\right)=\varphi\left(BA\right),\;\forall A,B\in\mathfrak{A}.\label{eq:gg3}
\end{equation}

\end{enumerate}
\end{xca}
We are aiming at a proof of the GNS theorem (\thmref{GNS}), and a
way to get more general representations of $*$-algebras. Indeed,
any representation is built up by the cyclic representations (\defref{cyclic}),
and each cyclic representation is in turn given by a GNS construction.
\index{cyclic!-representation}
\begin{defn}
\label{def:cyclic}A representation $\pi\in Rep(\mathfrak{A},\mathscr{H})$
is called \emph{cyclic,} with a \emph{cyclic vector} $u\in\mathscr{H}$,
if $\mathscr{H}=\overline{span}\left\{ \pi\left(A\right)u\:\big|\:A\in\mathfrak{A}\right\} $.
\index{cyclic!-vector}\end{defn}
\begin{thm}
\label{thm:cyclic2}Given any representation $\pi\in Rep(\mathfrak{A},\mathscr{H})$,
there exists an index set $J$, and closed subspaces $\mathscr{H}_{j}\subset\mathscr{H}$
($j\in J$) such that
\begin{enumerate}
\item $\mathscr{H}_{i}\perp\mathscr{H}_{j}$, $\forall i\neq j$; 
\item $\sum_{j\in J}^{\oplus}\mathscr{H}_{j}=\mathscr{H}$; and
\item there exists cyclic vectors $v_{j}\in\mathscr{H}_{j}$ such that the
restriction of $\pi$ to $\mathscr{H}_{j}$ is cyclic. 
\end{enumerate}
\end{thm}
\begin{rem}
The proof of \ref{thm:cyclic2} is very similar to the construction
of orthonormal basis (ONB) (use Zorn's lemma!); but here we get a
family of mutually orthogonal subspaces. 

\index{orthogonal!-subspaces}

Of course, if $\mathscr{H}_{j}$'s are all one-dimensional, then it
is a decomposition into ONB. Note that not every representation is
irreducible, but every representation can be decomposed into direct
sum of cyclic\index{representation!cyclic} representations. \index{cyclic!-representation}\end{rem}
\begin{xca}[Cyclic subspaces]
\myexercise{Cyclic subspaces}\label{exer:cyclic}Prove \thmref{cyclic2}.
\uline{Hint}: pick $v_{1}\in\mathscr{H}$, and let 
\[
\mathscr{H}_{v_{1}}:=\overline{span}\left\{ \pi\left(A\right)v_{1}:A\in\mathfrak{A}\right\} ,
\]
i.e., the cyclic subspace generated by $v_{1}$. If $\mathscr{H}_{v_{1}}\neq\mathscr{H}$,
then $\exists v_{2}\in\mathscr{H}\backslash\mathscr{H}_{v_{1}}$,
and the cyclic subspace $\mathscr{H}_{v_{2}}$, so that $\mathscr{H}_{v_{1}}$
and $\mathscr{H}_{v_{2}}$ are orthogonal. If $\mathscr{H}_{v_{1}}\oplus\mathscr{H}_{v_{2}}\neq\mathscr{H}$,
we then build $\mathscr{H}_{v_{3}}$ and so on. Now use transfinite
induction or Zorn's lemma\index{Zorn's lemma} to show the family
of direct sum of mutually orthogonal cyclic subspaces is total. The
final step is exactly the same argument for the existence of an ONB
of any Hilbert space.
\end{xca}
Now we proceed to prove the \thmref{GNS} (GNS), which is restated
below. \index{cyclic!-subspace}

\index{unitary equivalence}
\begin{thm}[Gelfand-Naimark-Segal]
 There is a bijection between states $\varphi$ and cyclic representations
$\pi\in Rep\left(\mathfrak{A},\mathscr{H},u\right)$, with $\left\Vert u\right\Vert =1$;
where 
\begin{equation}
\varphi(A)=\left\langle u,\pi(A)u\right\rangle ,\;\forall A\in\mathfrak{A}.\label{eq:gns1}
\end{equation}
Moreover, fix a state $\varphi$, the corresponding cyclic representation
is unique up to unitary equivalence. Specifically, if $\left(\pi_{1},\mathscr{H}_{1},u_{1}\right)$
and $\left(\pi_{2},\mathscr{H}_{2},u_{2}\right)$ are two cyclic representations,
with cyclic vectors $u_{1},u_{2}$, respectively, satisfying 
\begin{equation}
\varphi(A)=\left\langle u_{1},\pi_{1}(A)u_{1}\right\rangle =\left\langle u_{2},\pi_{2}(A)u_{2}\right\rangle ,\;\forall A\in\mathfrak{A};\label{eq:gns4}
\end{equation}
then 
\begin{equation}
W:\pi_{1}\left(A\right)u_{1}\longmapsto\pi_{2}\left(A\right)u_{2},\;A\in\mathfrak{A}\label{eq:gns3}
\end{equation}
extends to a unitary operator from $\mathscr{H}_{1}$ onto $\mathscr{H}_{2}$,
also denoted by $W$, and such that 
\begin{equation}
\pi_{2}W=W\pi_{1},\label{eq:gns5}
\end{equation}
i.e., $W$ intertwines the two representations.\index{cyclic!-representation}\end{thm}
\begin{rem}
For the non-trivial direction, let $\varphi$ be a given state on
$\mathfrak{A}$, and we need to construct a cyclic representation
$(\pi,\mathscr{H}_{\varphi},u_{\varphi})$. Note that $\mathfrak{A}$
is an algebra, and it is also a complex vector space. Let us try to
turn $\mathfrak{A}$ into a Hilbert space and see what conditions
are needed. There is a homomorphism $\mathfrak{A}\rightarrow\mathfrak{A}$
which follows from the associative law of $\mathfrak{A}$ being an
algebra, i.e., $(AB)C=A(BC)$. To continue, $\mathfrak{A}$ should
be equipped with an inner product. Using $\varphi$, we may set $\left\langle A,B\right\rangle _{\varphi}:=\varphi\left(A^{*}B\right)$,
$\forall A,B\in\mathfrak{A}$. Then $\left\langle \cdot,\cdot\right\rangle _{\varphi}$
is linear in the second variable, and conjugate linear in the first
variable. It also satisfies $\left\langle A,A\right\rangle _{\varphi}=\varphi\left(A^{*}A\right)\geq0$.
Therefore we take $\mathscr{H}_{\varphi}:=[\mathfrak{A}/\{A:\varphi(A^{*}A)=0\}]^{cl}$. \end{rem}
\begin{proof}
Given a cyclic representation $\pi\in Rep\left(\mathfrak{A},\mathscr{H},u\right)$,
define $\varphi$ as in (\ref{eq:gns1}). Clearly $\varphi$ is linear,
and 
\begin{eqnarray*}
\varphi\left(A^{*}A\right) & = & \left\langle u,\pi(A^{*}A)u\right\rangle \\
 & = & \left\langle u,\pi(A^{*})\pi\left(A\right)u\right\rangle \\
 & = & \left\langle \pi(A)u,\pi\left(A\right)u\right\rangle \\
 & = & \left\Vert \pi\left(A\right)u\right\Vert ^{2}\geq0.
\end{eqnarray*}
Thus $\varphi$ is a state. 

Conversely, fix a state $\varphi$ on $\mathfrak{A}$. Set 
\[
\mathscr{H}_{0}:=\left\{ \sum_{i=1}^{n}c_{i}A_{i}\:\big|\:c_{i}\in\mathbb{C},\:n\in\mathbb{N}\right\} 
\]
and define the inner product 
\[
\left\langle \sum c_{i}A_{i},\sum d_{i}B_{i}\right\rangle _{\varphi}:=\sum\sum\overline{c_{i}}d_{j}\varphi\left(A_{i}^{*}B_{j}\right).
\]
Note that, by definition, 
\begin{equation}
\left\Vert \sum c_{i}A_{i}\right\Vert _{\varphi}^{2}=\left\langle \sum c_{i}A_{i},\sum c_{i}A_{i}\right\rangle _{\varphi}=\sum\sum\overline{c_{i}}c_{j}\varphi\left(A_{i}^{*}A_{j}\right)\geq0.\label{eq:gns2}
\end{equation}
The RHS of (\ref{eq:gns2}) is positive since $\varphi$ is a state.
Recall that $\varphi\left(A^{*}A\right)\geq0$, for all $A\in\mathfrak{A}$,
and this implies that for all $n\in\mathbb{N}$, the matrix $\left(\varphi\left(A_{i}^{*}A_{j}\right)\right)_{i,j=1}^{n}$
is positive definite, hence (\ref{eq:gns2}) holds. 

\emph{Proof of }(\ref{eq:gns1}): Now, let $\mathscr{H}_{\varphi}:=$
completion of $\mathscr{H}_{0}$ under $\left\langle \cdot,\cdot\right\rangle _{\varphi}$
modulo elements $s$ such that $\left\Vert s\right\Vert _{\varphi}=0$.
See \lemref{gns} below. $\mathscr{H}_{\varphi}$ is the desired cyclic
space, consisting of equivalence classes $\left[A\right]$, $\forall A\in\mathfrak{A}$.
Next, let $u_{\varphi}=[1_{\mathfrak{A}}]=$ equivalence class of
the identity element, and set 
\[
\pi\left(A\right):=\left[A\right]=[A1_{\mathfrak{A}}]=[A][1_{\mathfrak{A}}];
\]
then one checks that $\pi\in Rep\left(\mathfrak{A},\mathscr{H}_{\varphi}\right)$,
and therefore $\varphi\left(A\right)=\left\langle u_{\varphi},\pi\left(A\right)u_{\varphi}\right\rangle _{\varphi}$,
$\forall A\in\mathfrak{A}$. 

For uniqueness, let $\left(\pi_{1},\mathscr{H}_{1},u_{1}\right)$
and $\left(\pi_{2},\mathscr{H}_{2},u_{2}\right)$ be as in the statement
of the theorem, and let $W$ be as in (\ref{eq:gns3}). By (\ref{eq:gns4}),
we have 
\[
\varphi\left(A^{*}A\right)=\left\Vert \pi_{2}\left(A\right)u_{2}\right\Vert ^{2}=\left\Vert \pi_{1}\left(A\right)u_{1}\right\Vert ^{2}
\]
so that $W$ is isometric. But since $\mathscr{H}_{i}=\overline{span}\left\{ \pi_{i}\left(A\right)u_{i}:A\in\mathfrak{A}\right\} $,
$i=1,2$, then $W$ extends by density to a unitary operator from
$\mathscr{H}_{1}$ to $\mathscr{H}_{2}$. \index{cyclic!-space}

\emph{Proof of }(\ref{eq:gns5}): Finally, for all $A,B\in\mathfrak{A}$,
we have 
\begin{eqnarray*}
W\pi_{1}\left(A\right)\left(\pi_{1}\left(B\right)u_{1}\right) & = & W\pi_{1}\left(AB\right)u_{1}\\
 & = & \pi_{2}\left(AB\right)u_{2}\\
 & = & \pi_{2}\left(A\right)\left(\pi_{2}\left(B\right)u_{2}\right)\\
 & = & \pi_{2}\left(A\right)W\pi_{1}\left(B\right)u_{1};
\end{eqnarray*}
therefore, by the density argument again, we conclude that 
\[
\pi_{2}\left(A\right)=W\pi_{1}\left(A\right)\;\forall A\in\mathfrak{A}.
\]
This is the intertwining property in (\ref{eq:gns5}). \index{positive definite!-operator}\end{proof}
\begin{lem}
\textup{\label{lem:gns}$\{A\in\mathfrak{A}:\varphi(A^{*}A)=0\}$
is a closed two-sided ideal in $\mathfrak{A}$. }\index{ideal}\end{lem}
\begin{proof}
This follows from the Schwarz inequality. Note that\index{inequality!Schwarz}
\[
\begin{bmatrix}\varphi\left(A^{*}A\right) & \varphi\left(A^{*}B\right)\\
\varphi\left(B^{*}A\right) & \varphi\left(B^{*}B\right)
\end{bmatrix}
\]
is a positive definite matrix, and so its determinant is positive,
i.e., 
\begin{equation}
\left|\varphi\left(A^{*}B\right)\right|^{2}\leq\varphi\left(A^{*}A\right)\varphi\left(B^{*}B\right);\label{eq:gns6}
\end{equation}
using the fact that $\varphi\left(C^{*}\right)=\varphi\left(C\right)^{*}$,
$\forall C\in\mathfrak{A}$. The lemma follows from the estimate (\ref{eq:gns6}). \end{proof}
\begin{example}
Let $\mathfrak{A}=C[0,1]$. Set $\varphi:f\mapsto f(0)$, so that
$\varphi\left(f^{*}f\right)=\left|f\left(0\right)\right|^{2}\geq0$.
Then, 
\[
\ker\varphi=\{f\in C\left[0,1\right]\:\big|\:f(0)=0\}
\]
and $C\left[0,1\right]/\ker\varphi$ is one dimensional. The reason
is that if $f\in C[0,1]$ such that $f(0)\neq0$, then we have $f(x)\sim f(0)$
since $f(x)-f(0)\in\ker\varphi$, where $f(0)$ represents the constant
function $f(0)$ over $[0,1]$. This shows that $\varphi$ is a pure
state, since the representation has to be irreducible. \index{state!pure-}\end{example}
\begin{xca}[The GNS construction]
\myexercise{The GNS construction}\label{exer:gns}Fill in the remaining
details in the above proof of the GNS theorem.
\end{xca}
Using GNS construction we get the following structure theorem for
abstract \index{algebras!$C^{*}$-algebra}$C^{*}$-algebras. As a
result, all $C^{*}$-algebras are sub-algebras of $\mathscr{B}\left(\mathscr{H}\right)$
for some Hilbert space $\mathscr{H}$.
\begin{thm}[Gelfand-Naimark]
\label{thm:gns1} Every $C^{*}$-algebra (abelian or non-abelian)
is isometrically isomorphic to a norm-closed sub-algebra of $\mathscr{B}(\mathscr{H})$,
for some Hilbert space $\mathscr{H}$.\end{thm}
\begin{proof}
Let $\mathfrak{A}$ be any $C^{*}$-algebra, no Hilbert space $\mathscr{H}$
is given from outside. Let $S(\mathfrak{A})$ be the states on $\mathfrak{A}$,
which is a compact\index{compact} convex\index{convex} subset of
the dual space $\mathfrak{A}^{*}$. Here, compactness refers to the
weak $*$-topology\index{weak{*}-topology}.

We use Hahn-Banach theorem to show that there are plenty of states.
Specifically, $\forall a\in\mathfrak{A}$, $\exists\varphi\in\mathfrak{A}^{*}$
such that $\varphi(a)>0$. It is done first on the 1-dimensional subspace
\[
tA\mapsto t\in\mathbb{R},
\]
and then extends to $\mathfrak{A}$. (Note this is also a consequence
of Krein\textendash Milman, i.e., $S(\mathfrak{A})=cl(\mbox{pure states})$.
We will come back to this point later.)\index{Theorem!Krein-Milman's-}

For each state $\varphi$, one gets a \emph{cyclic representation}
$\left(\pi_{\varphi},\mathscr{H}_{\varphi},u_{\varphi}\right)$. Applying
transfinite induction, one concludes that $\pi:=\oplus\pi_{\varphi}$
is a representation on the Hilbert space $\mathscr{H}:=\oplus\mathscr{H}_{\varphi}$.
For details, see e.g., \cite{Rud73}.
\end{proof}
\index{Hahn-Banach theorem}
\begin{thm}
Let $\mathfrak{A}$ be an abelian $C^{*}$-algebra. Then there is
a compact Hausdorff space $X$, unique up to homeomorphism, such that
$\mathfrak{A}\cong C(X)$.
\end{thm}

\section{\label{sec:predual}States, Dual and Pre-dual}

Let $V$ be a \emph{Banach space}, i.e., (recall,  \chapref{basic}):
\index{dual}\index{pre-dual}\index{space!Banach-}
\begin{itemize}
\item $V$ is a vector space over $\mathbb{C}$; 
\item $\exists$ norm $\left\Vert \cdot\right\Vert $
\item $V$ is complete with respect to $\left\Vert \cdot\right\Vert $
\end{itemize}
The dual space $V^{*}$ consists of linear functionals\index{functional}
$l:V\rightarrow\mathbb{C}$ satisfying
\[
\left\Vert l\right\Vert :=\sup_{\left\Vert v\right\Vert =1}\left|l\left(v\right)\right|<\infty.
\]
These are the continuous linear functionals. \index{continuous linear functional}

The \emph{Hahn-Banach Theorem} implies that for all $v\in V$, $\left\Vert v\right\Vert \neq0$,
there exists $l_{v}\in V^{*}$, of norm $1$, such that $l(v)=\left\Vert v\right\Vert $.
Recall the construction is to first define $l_{v}$ on the one-dimensional
subspace spanned by the vector $v$, then use transfinite induction
to extend $l_{v}$ to all of $V$. Notice that $V^{*}$ is always
complete, even if $V$ is an incomplete normed space. In other words,
$V^{*}$ is always a Banach space. \index{Banach space!double-dual}\index{space!dual-}

Now $V$ is embedded into $V^{**}$ (as we always do this) via the
mapping
\begin{gather}
V\ni v\mapsto\psi\left(v\right)\in V^{**},\;\mbox{where}\nonumber \\
\psi\left(v\right)\left(l\right):=l\left(v\right),\;\forall l\in V^{*}.\label{eq:vv}
\end{gather}

Below we give a number of applications:
\begin{xca}[Identification by isometry]
\myexercise{Identification by isometry}Show that $V\xrightarrow{\;\psi\;}V^{**}$
in (\ref{eq:vv}) is isometric, i.e., 
\[
\left\Vert \psi\left(v\right)\right\Vert _{**}=\left\Vert v\right\Vert ,\;\forall v\in V.
\]
\end{xca}
\begin{example}
Let $X$ be a compact Hausdorff space. The algebra $C(X)$ of all
continuous functions on $X$ with the sup norm, i.e., $\left\Vert f\right\Vert _{\infty}:=\sup_{x\in X}\left|f\left(x\right)\right|$,
is a Banach space. \index{space!Banach-}
\end{example}

\begin{example}
The classical $L^{p}$ space: $(l^{p})^{*}=l^{q}$, $(L^{p})^{*}=L^{q}$,
for $1/p+1/q=1$ and $1\leq p<\infty$. If $1<p<\infty$, then $(l^{p})^{**}=l^{p}$,
i.e., these spaces are \emph{reflexive}. For $p=1$, however, we have
$(l^{1})^{*}=l^{\infty}$, but $(l^{\infty})^{*}$ is much bigger
than $l^{1}$. Also note that $(l^{p})^{*}\neq l^{q}$ except for
$p=q=2$. And $l^{p}$ is a Hilbert space iff $p=2$. \index{Hilbert space!$L^{2}$}\index{Hilbert space!$l^{2}$}
\index{Banach space!reflexive}\index{space!Hilbert-}\index{space!$L^{p}$-}
\end{example}
Let $B$ be a Banach space and denote by $B^{*}$ its dual space.
$B^{*}$ is a Banach space as well, where the norm is defined by
\[
\left\Vert f\right\Vert _{B^{*}}=\sup_{\left\Vert x\right\Vert =1}\left\{ \left|f\left(x\right)\right|\right\} .
\]
Let $B_{1}^{*}=\{f\in B^{*}:\left\Vert f\right\Vert \leq1\}$ be the
unit ball in $B^{*}$. 

\index{Banach-Alaoglu Theorem}\index{Theorem!Banach-Alaoglu-}
\begin{thm}[Banach-Alaoglu]
\label{thm:Alaoglu}$B_{1}^{*}$ is weak $*$ compact\index{weak{*}-compact}
in $B^{*}$.\end{thm}
\begin{proof}
This is proved by showing $B_{1}^{*}$ is a closed subset in $\Omega:=\prod_{\left\Vert x\right\Vert =1}\mathbb{C}_{1}$,
with $\mathbb{C}_{1}=\left\{ z\in\mathbb{C}\::\:\left|z\right|\leq1\right\} $;
and $\Omega$ is given its product topology, and is compact and Hausdorff.
\index{product topology}
\end{proof}
As an application, we have 
\begin{cor}
Let $B$ be a separable Banach space. Then every bounded sequence
in $B^{*}$ has a convergent subsequence in the weak $*$-topology.
\end{cor}

\begin{cor}
Every bounded sequence in $\mathscr{B}(\mathscr{H})$ contains a convergence
subsequence in the weak $*$-topology.
\end{cor}
We show in Theorem \ref{thm:predual} that $\mathscr{B}\left(\mathscr{H}\right)=\mathscr{T}_{1}\left(\mathscr{H}\right)^{*}$,
where $\mathscr{T}_{1}\left(\mathscr{H}\right)=$ trace-class operators.
\index{trace}

Now we turn to Hilbert space, say $\mathscr{H}$:
\begin{itemize}
\item $\mathscr{H}$ is a vector space over $\mathbb{C}$; 
\item it has an inner product $\left\langle \cdot,\cdot\right\rangle $,
and the norm $\left\Vert \cdot\right\Vert :=\sqrt{\left\langle \cdot,\cdot\right\rangle }$; 
\item $\mathscr{H}$ is complete with respect to $\left\Vert \cdot\right\Vert $;
\item $\mathscr{H}^{*}=\mathscr{H}$, i.e., $\mathscr{H}$ is reflexive; 
\item every Hilbert space has an orthonormal basis (by Zorn's lemma)
\end{itemize}
The identification $\mathscr{H}=\mathscr{H}^{*}$ is due to Riesz,
and the corresponding map is given by \index{Riesz' theorem}\index{Theorem!Riesz-}
\[
h\mapsto\left\langle h,\cdot\right\rangle \in\mathscr{H}^{*}
\]
This can also be seen by noting via an ONB that $\mathscr{H}$ is
unitarily equivalent to $l^{2}(A)$, with some index set $A$, and
$l^{2}(A)$ is reflexive.

The set of all bounded operators $\mathscr{B}(\mathscr{H})$ on $\mathscr{H}$
is a Banach space. We ask two questions: 
\begin{enumerate}
\item What is the dual $\mathscr{B}(\mathscr{H})^{*}$? 
\item Is $\mathscr{B}(\mathscr{H})$ the dual space of some Banach space?
\index{Banach space}\index{space!Banach-}
\end{enumerate}
The first question is extremely difficult and we will discuss that
later. 

For the present section, we show that 
\[
\mathscr{B}(\mathscr{H})=\mathscr{T}_{1}(\mathscr{H})^{*}
\]
where we denote by $\mathscr{T}_{1}(\mathscr{H})$ the trace-class
operators in $\mathscr{B}(\mathscr{H})$. For more details, see \thmref{predual},
and \secref{3norm}.

Let $\rho:\mathscr{H}\rightarrow\mathscr{H}$ be a compact selfadjoint
operator. Assume $\rho$ is positive, i.e., $\left\langle x,\rho x\right\rangle \geq0$
for all $x\in\mathscr{H}$. By the spectral theorem of compact operators,
we get the following decomposition
\begin{equation}
\rho=\sum\lambda_{k}P_{k}\label{eq:tr1}
\end{equation}
where $\lambda_{1}\geq\lambda_{2}\geq\cdots\rightarrow0$, and $P_{k}$
is the projection onto the finite dimensional eigenspace of $\lambda_{k}$.
\index{Spectral Theorem}

In general, we want to get rid of the assumption that $\rho\geq0$.
This is done using the polar decomposition, which we will consider
in \secref{polar} even for unbounded operators. It is much easier
for bounded operators: If $A\in\mathscr{B}\left(\mathscr{H}\right)$,
$A^{*}A$ is positive, selfadjoint, and so by the spectral theorem,
we may take $\left|A\right|:=\sqrt{A^{*}A}$. Then, one checks that
\index{polar decomposition} 
\[
\left\Vert Ax\right\Vert ^{2}=\left\langle Ax,Ax\right\rangle =\left\langle x,A^{*}Ax\right\rangle =\left\langle \sqrt{A^{*}A}x,\sqrt{A^{*}A}x\right\rangle =\left\Vert \left|A\right|x\right\Vert ^{2},
\]
thus 
\begin{equation}
\left\Vert A\right\Vert =\left\Vert \left|A\right|\right\Vert \label{eq:tr3}
\end{equation}
and there is a partial isometry $V:\mbox{range}\left(\left|A\right|\right)\rightarrow\mbox{range}\left(A\right)$,
and the following polar decomposition holds:
\begin{equation}
A=V\left|A\right|\label{eq:tr4}
\end{equation}
We will come back to this point in \secref{polar} when we consider
unbounded operators.

\index{isometry!partial-}
\begin{cor}
\label{cor:tr1}Let $A\in\mathscr{T}_{1}\left(\mathscr{H}\right)$,
then $A$ has the following decomposition
\begin{equation}
A=\sum_{n}\lambda_{n}\left|f_{n}\left\rangle \right\langle e_{n}\right|\label{eq:tr6}
\end{equation}
where $\left\{ e_{n}\right\} $ and $\left\{ f_{n}\right\} $ are
ONBs in $\mathscr{H}$. \end{cor}
\begin{proof}
Using the polar decomposition $A=V\left|A\right|$, we may first diagonalize
$\left|A\right|$ with respect to some ONB $\left\{ e_{n}\right\} $
as 
\[
\left|A\right|=\sum_{n}\lambda_{n}\left|e_{n}\left\rangle \right\langle e_{n}\right|,\;\mbox{then}
\]
\[
A=V\left|A\right|=\sum_{n}\lambda_{n}\left|Ve_{n}\left\rangle \right\langle e_{n}\right|=\sum_{n}\lambda_{n}\left|f_{n}\left\rangle \right\langle e_{n}\right|
\]
where $f_{n}:=Ve_{n}$. 
\end{proof}
With the above discussion, we may work, instead, with compact\index{compact}
operators $A:\mathscr{H}\rightarrow\mathscr{H}$ so that $A$ is a
trace class operator if $\left|A\right|$ (positive, selfadjoint)
satisfies condition (\ref{eq:tr2}). \index{trace}
\begin{defn}
Let $A$ be a compact operator with its polar decomposition $A=V\left|A\right|$,
where $\left|A\right|:=\sqrt{A^{*}A}$. Let $\left\{ \lambda_{k}\right\} _{k=1}^{\infty}$
be the eigenvalues of $\left|A\right|$, and $P_{k}$ the corresponding
spectral projections, see (\ref{eq:tr1}). We say $A$ is a trace
class operator, if\index{eigenvalue} 
\begin{equation}
\left\Vert A\right\Vert _{1}:=trace\left(\left|A\right|\right)=\sum_{n}\lambda_{n}\:\mbox{rank}\left(P_{k}\right)<\infty.\label{eq:tr2}
\end{equation}
\emph{Caution.} In our consideration of eigenvalue lists, we may of
course have multiplicity. But for compact operators, the multiplicity
is automatically finite for each non-zero eigenvalue. And if we have
sets of associated eigenvectors run through a local ONB in each of
the finite-dimensional eigenspaces, then multiplicity is counted this
way. But, alternatively, when computing a trace as a sum of eigenvalues,
then the term in such a sum must be counted with multiplicity. Or
each of the distinct numbers in an eigenvalue list can be multiplied
with the respective multiplicity. This will be clear from the context.
\end{defn}
\index{operators!trace class}

\index{polar decomposition}

We now continue the discussion from \secref{3norm} on spaces of operators. 
\begin{defn}
Let $A\in\mathscr{T}_{1}\left(\mathscr{H}\right)$, and $\left\{ e_{n}\right\} $
an ONB in $\mathscr{H}$. Set 
\begin{equation}
trace\left(A\right):=\sum_{n}\left\langle e_{n},Ae_{n}\right\rangle \label{eq:tr5}
\end{equation}
Note the RHS in (\ref{eq:tr5}) is independent of the choice of the
ONB. For if $\left\{ f_{n}\right\} $ is another ONB in $\mathscr{H}$,
using the Parseval identity repeatedly, we have \index{Parseval identity}
\begin{eqnarray*}
\sum\left\langle f_{n},Af_{n}\right\rangle  & = & \sum_{n}\sum_{m}\left\langle f_{n},e_{m}\right\rangle \left\langle e_{m},Af_{n}\right\rangle \\
 & = & \sum_{m}\sum_{n}\left\langle f_{n},e_{m}\right\rangle \left\langle A^{*}e_{m},f_{n}\right\rangle \\
 & = & \sum_{m}\left\langle A^{*}e_{m},e_{m}\right\rangle =\sum_{m}\left\langle e_{m},Ae_{m}\right\rangle .
\end{eqnarray*}
\end{defn}
\begin{cor}
\label{cor:tr2}Let $A\in\mathscr{T}_{1}\left(\mathscr{H}\right)$,
then 
\[
\left|trace\left(A\right)\right|\leq\left\Vert A\right\Vert _{1}.
\]
Therefore, the RHS in (\ref{eq:tr5}) is absolutely convergent.\end{cor}
\begin{proof}
By \corref{tr1}, there exists ONBs $\left\{ e_{n}\right\} $ and
$\left\{ f_{n}\right\} $, and $A$ has a decomposition as in (\ref{eq:tr6}).
Then, 
\begin{eqnarray*}
\left|trace\left(A\right)\right| & \leq & \sum_{n}\left|\left\langle e_{n},Ae_{n}\right\rangle \right|=\sum_{n}\lambda_{n}\left|\left\langle e_{n},f_{n}\right\rangle \right|\\
 & \leq & \sum_{n}\lambda_{n}\left\Vert e_{n}\right\Vert \left\Vert f_{n}\right\Vert =\sum_{n}\lambda_{n}=\left\Vert A\right\Vert _{1}<\infty.
\end{eqnarray*}
We have used the fact that $trace\left(A\right)$ is independent of
the choice of an ONBs.\end{proof}
\begin{lem}
\label{lem:trace}Let $\mathscr{T}_{1}(\mathscr{H})$ be the trace
class introduced above. Then, 
\begin{enumerate}
\item $\mathscr{T}_{1}(\mathscr{H})$ is a two-sided ideal in $\mathscr{B}(\mathscr{H})$. 
\item $trace(AB)=trace(BA)$
\item $\mathscr{T}_{1}(\mathscr{H})$ is a Banach space with respect to
the trace norm (\ref{eq:tr2}). 
\end{enumerate}
\end{lem}
\begin{xca}[A pre-dual]
\myexercise{A pre-dual}\label{exer:trace}Prove \lemref{trace}.\end{xca}
\begin{lem}
Let $\rho\in\mathscr{T}_{1}(\mathscr{H})$, then the map $A\mapsto trace(A\rho)$
is a state on $\mathscr{B}(\mathscr{H})$. These are called the \uline{normal}
states.\end{lem}
\begin{proof}
By \lemref{trace}, $A\rho\in\mathscr{T}_{1}(\mathscr{H})$ for all
$A\in\mathscr{B}(\mathscr{H})$. The map $A\mapsto trace(A\rho)$
is in $\mathscr{B}(\mathscr{H})^{*}$ means that the pairing $\left(A,\rho\right)\mapsto trace(A\rho)$
satisfies
\[
\left|trace(A\rho)\right|\leq\left\Vert A\right\Vert \left\Vert \rho\right\Vert _{1}.
\]
By \corref{tr2}, it suffices to verify, instead, that 
\[
\left\Vert A\rho\right\Vert _{1}\leq\left\Vert A\right\Vert \left\Vert \rho\right\Vert _{1}.
\]
Indeed, if we choose an ONB $\left\{ e_{n}\right\} $ in $\mathscr{H}$
that diagonalizes $\left|\rho\right|$, i.e., 
\[
\left|\rho\right|=\sum\lambda_{n}\left|e_{n}\left\rangle \right\langle e_{n}\right|,\;\mbox{where}\;\sum\lambda_{n}<\infty,\lambda_{n}>0,\;\forall n;
\]
then 
\begin{eqnarray*}
\left\Vert A\rho\right\Vert _{1} & = & trace\left(\sqrt{\rho^{*}A^{*}A\rho}\right)=trace\left(\sqrt{\rho^{*}\rho}\sqrt{A^{*}A}\right)\\
 & = & \sum_{n}\left\langle e_{n},\left|A\right|\left|\rho\right|e_{n}\right\rangle =\sum_{n}\lambda_{n}\left\langle e_{n},\left|A\right|e_{n}\right\rangle \\
 & \leq & \left\Vert A\right\Vert \sum_{k}\lambda_{k}=\left\Vert A\right\Vert \left\Vert \rho\right\Vert _{1}.
\end{eqnarray*}
\end{proof}
\begin{thm}
\label{thm:predual}$\mathscr{T}_{1}^{*}\left(\mathscr{H}\right)=\mathscr{B}\left(\mathscr{H}\right)$.\end{thm}
\begin{proof}
Let $l\in\mathscr{T}_{1}^{*}$. By $\mathscr{T}_{1}^{*}=\left(\mathscr{T}_{1}\right)^{*}$
we mean the dual Banach, duality with respect to the trace-norm. \index{space!dual-}

How to get an operator $A$? The operator $A$ must satisfy
\[
l(\rho)=trace(\rho A),\quad\forall\rho\in\mathscr{T}_{1}.
\]
How to pull an operator $A$ out of the hat? The idea also goes back
to Dirac. It is in fact not difficult to find $A$. Since $A$ is
determined by its matrix, it suffices to find $\left\langle f,Af\right\rangle $,
the entries in the matrix of $A$. 

For any $f_{1},f_{2}\in\mathscr{H}$, the rank-one operator $\left|f_{1}\left\rangle \right\langle f_{2}\right|$
is in $\mathscr{T}_{1}$, hence we know what $l$ does to it, i.e.,
we know the numbers $l\left(\left|f_{1}\left\rangle \right\langle f_{2}\right|\right)$.
But since $l\left(\left|f_{1}\left\rangle \right\langle f_{2}\right|\right)$
is linear in $f_{1}$, and conjugate linear in $f_{2}$, by the Riesz
theorem for Hilbert space, there exists a unique operator $A$ such
that
\[
l\left(\left|f_{1}\left\rangle \right\langle f_{2}\right|\right)=\left\langle f_{2},Af_{1}\right\rangle .
\]
Now we check that $l(\rho)=trace(\rho A)$. By \corref{tr1}, any
$\rho\in\mathscr{T}_{1}$ can be written as $\rho=\sum_{n}\lambda_{n}\left|f_{n}\left\rangle \right\langle e_{n}\right|$,
where $\left\{ e_{n}\right\} $ and $\left\{ f_{n}\right\} $ are
some ONBs in $\mathscr{H}$. Then,
\begin{eqnarray*}
trace\left(\rho A\right) & = & trace\left(\sum_{n}\lambda_{n}\left|f_{n}\left\rangle \right\langle e_{n}\right|A\right)\\
 & = & trace\left(\sum_{n}\lambda_{n}\left|Af_{n}\left\rangle \right\langle e_{n}\right|\right)\\
 & = & \sum_{m}\sum_{n}\lambda_{n}\left\langle u_{m},Af_{n}\right\rangle \left\langle e_{n},u_{m}\right\rangle \\
 & = & \sum_{n}\lambda_{n}\left(\sum_{m}\left\langle u_{m},Af_{n}\right\rangle \left\langle e_{n},u_{m}\right\rangle \right)\\
 & = & \sum_{n}\lambda_{n}\left\langle e_{n},Af_{n}\right\rangle \left(=l\left(\rho\right)\right)
\end{eqnarray*}
 where $\left\{ u_{n}\right\} $ is an ONB in $\mathscr{H}$, and
the last step follows from Parseval's identity.\end{proof}
\begin{rem}
If $B$ is the dual of a Banach space, then we say that $B$ has a
pre-dual. For example $l^{\infty}=(l^{1})^{*}$, hence $l^{1}$ is
the pre-dual of $l^{\infty}$. \index{pre-dual}\index{space!dual-}

Another example: Let $\mathbb{H}_{1}$ be hardy space of analytic
functions on the disk \cite{Rud87}. $(\mathbb{H}_{1})^{*}=$ BMO,
where BMO refers to bounded mean oscillation. It was developed by
Charles Fefferman in 1974 who won the fields medal for this theory.
See \cite{MR0280994}. (Getting hands on a specific dual space is
often a big thing.) \index{space!Hardy-}\end{rem}
\begin{defn}
Let $\mathbb{D}$ be the complex disk 
\[
\mathbb{D}=\left\{ z\in\mathbb{C}\::\:\left|z\right|<1\right\} .
\]
Consider functions $f$ analytic on $\mathbb{D}$ such that
\begin{equation}
\sup_{0<r<1}\frac{1}{2\pi}\int_{-\pi}^{\pi}\left|f\left(re^{it}\right)\right|dt<\infty.\label{eq:h1norm}
\end{equation}
This is the $\mathbb{H}_{1}$-Hardy space, and the $\mathbb{H}_{1}$-norm
is the supremum in (\ref{eq:h1norm}). (The literature on Hardy space
is extensive, and we refer to \cite{Rud87} for overview and details.)\end{defn}
\begin{thm}[C. Fefferman]
 $\mathbb{H}_{1}^{*}=\textup{BMO}$. \end{thm}
\begin{proof}
We refer to \cite{MR0280994}.\end{proof}
\begin{defn}
Let $f$ be a locally integrable function on $\mathbb{R}^{n}$, and
let $Q$ run through all $n$-cubes $\subset\mathbb{R}^{n}$. Set
\[
f_{Q}=\frac{1}{\left|Q\right|}\int_{Q}f\left(y\right)dy.
\]
We say that $f\in BMO$ iff (Def.) 
\begin{equation}
\sup_{Q}\frac{1}{\left|Q\right|}\int_{Q}\left|f\left(x\right)-f_{Q}\right|dx<\infty.\label{eq:bmo}
\end{equation}
In this case the LHS in (\ref{eq:bmo}) is the BMO-norm of $f$. Moreover,
BMO is a Banach space.\end{defn}
\begin{xca}[The Bohr compactification]
\label{exer:bohr}\myexercise{The Bohr compactification} From abstract
harmonic analysis (see \cite{Ru90}), we know that every locally Abelian
(l.c.A.) group $G$ has a Haar measure, unique up to scalar normalization.
When $G$ is a given l.c.A. group, we denote by $\widehat{G}$ its
dual group (of all continuous unitary characters.) The duality theorem
for l.c.A. groups $G$ states the following:

\begin{equation}
\mbox{\ensuremath{G} is compact}\Longleftrightarrow\mbox{\ensuremath{\widehat{G}} is discrete.}\label{eq:bo1}
\end{equation}
Moreover, in general, $G\simeq\widehat{\widehat{G}}$ where $G$ is
l.c.A.. Now consider the group $\mathbb{R}$ (the reals) with addition,
but in its \emph{discrete topology}, (usually denoted $\mathbb{R}_{d}$.)
The corresponding dual group $\mathbb{R}_{b}=\left(\mathbb{R}_{d}\right)^{\wedge}$
is therefore compact by (\ref{eq:bo1}). It is called the \emph{Bohr-compactification
}of $\mathbb{R}$. Let $d\chi$ denote its Haar measure. 
\begin{enumerate}
\item Show that $L^{2}\left(G_{b},d\chi\right)$ is a \emph{non-separable
}Hilbert space. 
\item Find an ONB in $L^{2}\left(G_{b},d\chi\right)$ indexed by $\mathbb{R}$.
\end{enumerate}
\end{xca}

\section{\label{sec:NHspace}New Hilbert Spaces From \textquotedblleft old\textquotedblright{}}

Below we consider some cases of building new Hilbert spaces from given
ones. Only sample cases are fleshed out; and they will be needed in
the sequel.

\paragraph{\uline{An Overview:}}

\subsection{GNS }

See \secref{gns}.

\subsection{Direct sum $\bigoplus_{\alpha}\mathscr{H}_{\alpha}$}
\begin{enumerate}[label=(\alph{enumi})]
\item \begin{flushleft}
 Let $\mathscr{H}_{i}$, $i=1,2$ be two given Hilbert spaces, then
the direct ``orthogonal'' sum $\mathscr{H}=\mathscr{H}_{1}\oplus\mathscr{H}_{2}$
is as follows:\index{Hilbert space!direct sum}
\begin{gather}
\mathscr{H}=\left\{ \text{symbol pairs}\;h_{1}\oplus h_{2},\;h_{i}\in\mathscr{H}_{i},\;i=1,2\right\} ,\;\mbox{and}\nonumber \\
\left\Vert h_{1}\oplus h_{2}\right\Vert _{\mathscr{H}}^{2}=\left\Vert h_{1}\right\Vert _{\mathscr{H}_{1}}^{2}+\left\Vert h_{2}\right\Vert _{\mathscr{H}_{2}}^{2}.\label{eq:hs1}
\end{gather}

\par\end{flushleft}
\item[(b)] \begin{flushleft}
Given an indexed family of Hilbert spaces $\left\{ \mathscr{H}_{\alpha}\right\} _{\alpha\in A}$
where $A$ is a set; then set $\mathscr{H}:=\oplus_{A}\mathscr{H}_{\alpha}$
to be 
\begin{equation}
\left\Vert \sum_{\alpha\in A}^{\oplus}h_{\alpha}\right\Vert _{\mathscr{H}}^{2}=\sum_{\alpha\in A}\left\Vert h_{\alpha}\right\Vert _{\mathscr{H}_{\alpha}}^{2}<\infty;\label{eq:hs2}
\end{equation}
i.e., finiteness of the sum in (\ref{eq:hs2}) is part of the definition. 
\par\end{flushleft}\end{enumerate}
\begin{xca}[Unitary operators on a direct Hilbert sum]
\myexercise{Unitary operators on a direct Hilbert sum}Let $\mathscr{H}_{i}$,
$i=1,2$, be Hilbert spaces, and set $\mathscr{H}:=\mathscr{H}_{1}\oplus\mathscr{H}_{2}$. 
\begin{enumerate}
\item[(i)] Let $G_{\mathscr{H}}$, and $G_{\mathscr{H}_{i}}$, $i=1,2$, be
the respective groups of unitary operators. Show that $G_{\mathscr{H}_{1}}\times G_{\mathscr{H}_{2}}$
is a subgroup of $G_{\mathscr{H}}$.
\item[(ii)] Let $L\in\mathscr{B}\left(\mathscr{H}\right)$, where $\mathscr{H}=\mathscr{H}_{1}\oplus\mathscr{H}_{2}$,
and suppose $L$ commutes with the group $G_{\mathscr{H}_{1}}\times G_{\mathscr{H}_{2}}$
in (i); then show that $L$ must have the following form 
\[
L=\left(\alpha I_{\mathscr{H}_{1}}\right)\times\left(\beta I_{\mathscr{H}_{2}}\right)
\]
where $\alpha,\beta\in\mathbb{C}$. (We say that the \emph{commutant}
of the group $G_{\mathscr{H}_{1}}\times G_{\mathscr{H}_{2}}$ has
this form; it is two-dimensional.)
\item[(iii)]  Let $A\in\mathscr{B}\left(\mathscr{H}_{2},\mathscr{H}_{1}\right)$,
and $B\in\mathscr{B}\left(\mathscr{H}_{1},\mathscr{H}_{2}\right)$;
show that the block-operator matrix $\begin{pmatrix}0 & A\\
B & 0
\end{pmatrix}$ defines a unitary operator in $\mathscr{H}=\mathscr{H}_{1}\oplus\mathscr{H}_{2}$
if and only if
\[
AA^{*}=I_{\mathscr{H}_{1}},\quad A^{*}A=I_{\mathscr{H}_{2}},
\]
and 
\[
BB^{*}=I_{\mathscr{H}_{2}},\quad B^{*}B=I_{\mathscr{H}_{1}},
\]
i.e., the two operators are unitary between the respective Hilbert
spaces. 
\end{enumerate}
\end{xca}

\subsection{Hilbert-Schmidt operators (continuing the discussion in 1)}

Let $\mathscr{H}$ be a fixed Hilbert space, and set 
\begin{equation}
\mathscr{H}S\left(\mathscr{H}\right):=\left\{ T\in\mathscr{B}\left(\mathscr{H}\right)\;\big|\;T^{*}T\;\mbox{is trace class}\right\} \label{eq:hs3}
\end{equation}
and set 
\begin{equation}
\left\Vert T\right\Vert _{\mathscr{H}S}^{2}:=trace\left(T^{*}T\right);\label{eq:hs4}
\end{equation}
similarly if $S,T\in\mathscr{H}S\left(\mathscr{H}\right)$, set 
\begin{equation}
\left\langle S,T\right\rangle _{\mathscr{H}S}:=trace\left(S^{*}T\right).\label{eq:hs5}
\end{equation}
Note that finiteness on the RHS in (\ref{eq:hs4}) is part of the
definition.\index{operators!Hilbert-Schmidt}

\subsection{Tensor-product $\mathscr{H}_{1}\otimes\mathscr{H}_{2}$}

Let $\mathscr{H}_{1}$ and $\mathscr{H}_{2}$ be two Hilbert spaces,
and consider finite-rank operators (rank-1 in this case):\index{Hilbert space!tensor product}\index{space!Hilbert-}
\index{product!tensor-}
\begin{gather}
\left|h_{1}\left\rangle \right\langle h_{2}\right|\;\left(\mbox{Dirac ket-bra}\right)\nonumber \\
\left|h_{1}\left\rangle \right\langle h_{2}\right|\left(u\right)=\left\langle h_{2},u\right\rangle _{\mathscr{H}_{2}}h_{1},\;\forall u\in\mathscr{H}_{2},\label{eq:hs6}
\end{gather}
so $T=\left|h_{1}\left\rangle \right\langle h_{2}\right|:\mathscr{H}_{2}\longrightarrow\mathscr{H}_{1}$
with the identification
\begin{equation}
h_{1}\otimes h_{2}\longleftrightarrow\left|h_{1}\left\rangle \right\langle h_{2}\right|.\label{eq:hs7}
\end{equation}
Set 
\begin{equation}
\left\Vert h_{1}\otimes h_{2}\right\Vert ^{2}:=trace\left(T^{*}T\right)=\left\Vert h_{1}\right\Vert _{\mathscr{H}_{1}}^{2}\left\Vert h_{2}\right\Vert _{\mathscr{H}_{2}}^{2}.\label{eq:hs8}
\end{equation}
For the Hilbert space $\mathscr{H}_{1}\otimes\mathscr{H}_{2}$ we
take the $\mathscr{H}S$-completion of the space of finite rank operators
spanned by the set in (\ref{eq:hs6}). The tensor product construction
fits with composite system in quantum mechanics. \index{completion!Hilbert-Schmidt}

\index{quantum mechanics!composite system}

\index{operators!finite rank-}

\subsection{Contractive inclusion}

Let $\mathscr{H}_{1}$ and $\mathscr{H}_{2}$ be two Hilbert spaces,
and let $T:\mathscr{H}_{1}\rightarrow\mathscr{H}_{2}$ be a contractive\index{operators!contractive}
linear operator, i.e., 
\begin{equation}
I_{\mathscr{H}_{1}}-T^{*}T\geq0.\label{eq:hs8-1}
\end{equation}
On the subspace
\begin{equation}
\mathscr{R}\left(T\right)=\left\{ Th_{1}\:\big|\:h_{1}\in\mathscr{H}_{1}\right\} \label{eq:hs9}
\end{equation}
(generally not closed in $\mathscr{H}_{2}$,) set 
\begin{equation}
\left\Vert Th_{1}\right\Vert _{\text{new}}:=\left\Vert h_{1}\right\Vert ,\;h_{1}\in\mathscr{H}_{1};\label{eq:hs10}
\end{equation}
then with $\left\Vert \cdot\right\Vert _{\text{new}}$, $\mathscr{R}\left(T\right)$
becomes a Hilbert space.

\subsection{Inflation (dilations)}

Let $T:\mathscr{H}_{1}:\longrightarrow\mathscr{H}_{2}$ be a contraction,
and set

\begin{equation}
\mathcal{U}=\begin{bmatrix}[c|c]T & \left(I_{2}-TT^{*}\right)^{\frac{1}{2}}\\
\hline \left(I_{1}-T^{*}T\right)^{\frac{1}{2}} & -T^{*}
\end{bmatrix}\label{eq:in1}
\end{equation}
The two operators in the off-diagonal slots are called the \textquotedblleft defect
operators\textquotedblright{} for the contraction $T$. Reason: the
pair of defect-operators are $\left(0,0\right)$ if and only if $T$
is a unitary isomorphism of $\mathscr{H}_{1}$ onto $\mathscr{H}_{2}$.
\begin{xca}[The Julia operator]
\myexercise{The Julia operator} Show that the matrix-block (\ref{eq:in1})
defines a \emph{unitary} operator $\mathcal{U}$ in $\mathscr{H}=\mathscr{H}_{1}\oplus\mathscr{H}_{2}$
(called the Julia operator); and that $P_{1}\mathcal{U}P_{1}=T$ where
$P_{1}$ denotes the projection of $\mathscr{H}$ onto $\mathscr{H}_{1}$. 
\end{xca}

\subsection{Reflection Positivity (or renormalization) $\left(\mathscr{H}_{+}/\mathscr{N}\right)^{\sim}$}

\paragraph{New Hilbert space from reflection positivity:\index{renormalization}}

\begin{figure}
\[ \xymatrix@R-2pc{\mathscr{H}\ar[rr]_{\mathcal{J}\mathcal{U}\mathcal{J}=\mathcal{U}^{*}}^{\mathcal{U}} &  & \mathscr{H} & \text{unitary}\\ \bigcup &  & \bigcup\\ \mathscr{H}_{+}\ar[dd] & \mathcal{U}\ar[dd] & \mathscr{H}_{+}\ar[dd] & \txt{invariant under \ensuremath{\mathcal{U}} \\ \ensuremath{\left\langle h_{+},\mathcal{J}h_{+}\right\rangle \geq0} }\\ \\ \mathscr{K}=\left(\mathscr{H}_{+}/\mathscr{N}\right)^{\sim}\ar@/_{2pc}/[rr]_{\widetilde{\mathcal{U}}} & \widetilde{\mathcal{U}} & \mathscr{K}=\left(\mathscr{H}_{+}/\mathscr{N}\right)^{\sim} & \txt{induced operator \\ \ensuremath{\mathcal{J}}-normalized\\ inner product \\ \\ \ensuremath{\widetilde{\mathcal{U}}} is contractive \\and selfadjoint} } \]

\protect\caption{\label{fig:rp}Reflection positivity. A unitary operator $\mathcal{U}$
transforms into a selfadjoint contraction $\widetilde{\mathcal{U}}$.}
\end{figure}
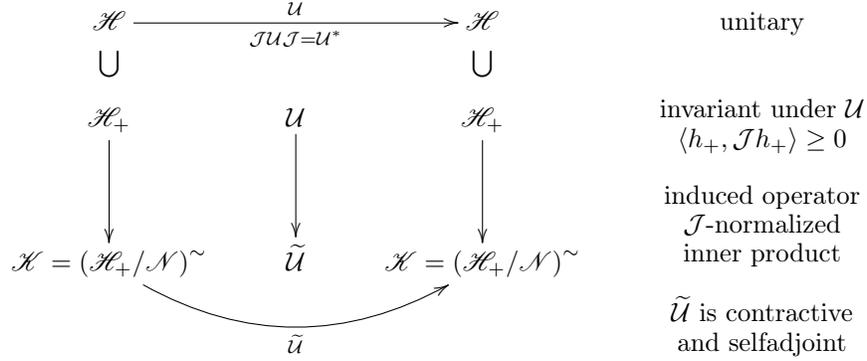

Let $\mathscr{H}$ be a given Hilbert space, $\mathscr{H}_{+}\subset\mathscr{H}$
a closed subspace, and let $\mathcal{U},\mathcal{J}:\mathscr{H}\rightarrow\mathscr{H}$
be two unitary operators, $\mathcal{J}$ satisfying the idempotency
condition
\begin{equation}
\mathcal{J}^{2}=I,\mbox{ as well as}\label{eq:rp1}
\end{equation}
\begin{equation}
\mbox{\ensuremath{\mathcal{J}\mathcal{U}\mathcal{J}}}=\mathcal{U}^{*},\:\mbox{and}\label{eq:rp2}
\end{equation}
\begin{equation}
\mathcal{U}\mathscr{H}_{+}\subset\mathscr{H}_{+};\;\mbox{and}\label{eq:rp3}
\end{equation}
finally
\begin{equation}
\left\langle h_{+},\mathcal{J}h_{+}\right\rangle \geq0,\;\forall h_{+}\in\mathscr{H}_{+}.\label{eq:rp4}
\end{equation}
Note that (\ref{eq:rp2}) states that $\mathcal{U}$ is unitarily
equivalent to its adjoint $\mathcal{U}^{*}$.\index{operators!adjoint-}
\begin{note}
Set $P_{+}:=Proj\mathscr{H}_{+}$ (= the projection onto $\mathscr{H}_{+}$),
then (\ref{eq:rp4}) is equivalent to 
\[
P_{+}\mathcal{J}P_{+}\geq0
\]
with respect to the usual ordering of operators.
\end{note}
Set 
\begin{eqnarray}
\mathscr{N} & = & \mbox{Ker}\left(P_{+}\mathcal{J}P_{+}\right)\label{eq:rp5}\\
 & = & \left\{ h_{+}\in\mathscr{H}_{+}\::\:\left\langle h_{+},\mathcal{J}h_{+}\right\rangle =0\right\} .\nonumber 
\end{eqnarray}
Set 
\begin{equation}
\mathscr{K}=\left(\mathscr{H}_{+}/\mathscr{N}\right)^{\sim}\label{eq:rp6}
\end{equation}
where ``\textasciitilde{}'' in (\ref{eq:rp6}) means Hilbert completion
with respect to the sesquilinear form: $\mathscr{H}_{+}\times\mathscr{H}_{+}\rightarrow\mathbb{C}$,
given by
\begin{equation}
\left\langle h_{+},h_{+}\right\rangle _{\mathscr{K}}:=\left\langle h_{+},\mathcal{J}h_{+}\right\rangle ,\label{eq:rp7}
\end{equation}
a renormalized inner product. \index{completion!Hilbert-}\index{sesquilinear form}
\begin{xca}[An induced operator]
\myexercise{An induced operator}\label{exer:rp}Let the setting
be as above. Show that $\widetilde{\mathcal{U}}:\mathscr{K}\rightarrow\mathscr{K}$,
given by
\begin{equation}
\widetilde{U}\left(\mbox{class}\:h_{+}\right)=\mbox{class}\left(\mathcal{U}h_{+}\right),\;h_{+}\in\mathscr{H}_{+}\label{eq:rp8}
\end{equation}
where $\mbox{class}\:h_{+}$ refers to the quotient in (\ref{eq:rp6}),
is selfadjoint and contractive (see \figref{rp}).\end{xca}
\begin{rem}
The construction outlined above is called ``reflection positivity'';
see e.g., \cite{JO00,MR965583}. It has many applications in physics
and in representation theory. \index{reflection positivity}
\end{rem}

\noindent \begin{flushleft}
\emph{Proof of the assertions in \figref{rp}}. Denote the ``new''
inner product in $\mathscr{K}$ by $\left\langle \cdot,\cdot\right\rangle _{\mathscr{K}}$,
and the initial inner product in $\mathscr{H}$ by $\left\langle \cdot,\cdot\right\rangle $. 
\par\end{flushleft}

\emph{$\widetilde{\mathcal{U}}$ is symmetric}: Let $x,y\in\mathscr{H}_{+}$,
then 
\begin{align*}
\langle x,\widetilde{\mathcal{U}}y\rangle_{\mathscr{K}} & =\left\langle x,\mathcal{J}\mathcal{U}y\right\rangle =\left\langle x,\mathcal{U}^{*}\mathcal{J}y\right\rangle \\
 & =\left\langle \mathcal{U}x,\mathcal{J}y\right\rangle =\langle\widetilde{\mathcal{U}}x,y\rangle_{\mathscr{K}}
\end{align*}
is the desired conclusion.

\emph{$\widetilde{\mathcal{U}}$ is contractive}: Let $x\in\mathscr{H}_{+}$,
then 
\begin{eqnarray*}
\left\Vert \widetilde{\mathcal{U}}x\right\Vert _{\mathscr{K}}^{2} & = & \left\langle \mathcal{U}x,\mathcal{J}\mathcal{U}x\right\rangle =\left\langle \mathcal{U}x,\mathcal{U}^{*}\mathcal{J}x\right\rangle \\
 & = & \left\langle \mathcal{U}^{2}x,\mathcal{J}x\right\rangle =\left\langle \mathcal{U}^{2}x,x\right\rangle _{\mathscr{K}}\\
 & \leq & \left\Vert \mathcal{U}^{2}x\right\Vert _{\mathscr{K}}\cdot\left\Vert x\right\Vert _{\mathscr{K}}\qquad\left(\mbox{by Schwarz in}\;\mathscr{K}\right)\\
 & \leq & \left\Vert \mathcal{U}^{4}x\right\Vert _{\mathscr{K}}^{\frac{1}{2}}\cdot\left\Vert x\right\Vert _{\mathscr{K}}^{1+\frac{1}{2}}\qquad\left(\mbox{by the first step}\right)\\
 & \leq & \left\Vert \mathcal{U}^{2^{n+1}}x\right\Vert _{\mathscr{K}}^{\frac{1}{2^{n}}}\cdot\left\Vert x\right\Vert _{\mathscr{K}}^{1+\frac{1}{2}+\cdots+\frac{1}{2^{n}}}.\qquad\left(\mbox{by iteration}\right)
\end{eqnarray*}
By the spectral-radius formula, 
\[
\lim_{n\rightarrow\infty}\left\Vert \mathcal{U}^{2^{n}}x\right\Vert _{\mathscr{K}}^{\frac{1}{2^{n}}}=1;
\]
and we get $\left\Vert \widetilde{\mathcal{U}}x\right\Vert _{\mathscr{K}}^{2}\leq\left\Vert x\right\Vert _{\mathscr{K}}^{2}$,
which is the desired contractivity. \null\hfill\qedsymbol\par
\begin{xca}[Time-reflection]
\myexercise{Time-reflection}\label{exer:rp1}Show that if $\left\{ \mathcal{U}_{t}\right\} _{t\in\mathbb{R}}$
is a unitary one-parameter group in $\mathscr{H}$ such that 
\[
\mathcal{J}\mathcal{U}_{t}\mathcal{J}=\mathcal{U}_{-t},\;t\in\mathbb{R},\;\mbox{and}
\]
\[
\mathcal{U}_{t}\mathscr{H}_{+}\subset\mathscr{H}_{+},\;t\in\mathbb{R}_{+},
\]
then 
\[
\mathcal{S}_{t}=\widetilde{\mathcal{U}}_{t}:\mathscr{K}\rightarrow\mathscr{K}
\]
is a selfadjoint contraction semigroup, $t\in\mathbb{R}_{+}$, i.e.,
there is a selfadjoint generator $L$ in $\mathscr{K}$, 
\begin{equation}
\left\langle k,Lk\right\rangle _{\mathscr{K}}\geq0,\;\forall k\in dom\left(L\right),\label{eq:rp9}
\end{equation}
where \index{semigroup} 
\begin{equation}
\mathcal{S}_{t}\left(=\widetilde{\mathcal{U}}_{t}\right)=e^{-tL},\;t\in\mathbb{R}_{+}\label{eq:rp10}
\end{equation}
and 
\begin{equation}
\mathcal{S}_{t_{1}}\mathcal{S}_{t_{2}}=\mathcal{S}_{t_{1}+t_{2}},\;t_{1},t_{2}\in\mathbb{R}_{+}.\label{eq:rp11}
\end{equation}

\end{xca}
\begin{figure}
\[
\xymatrix{
A\ar[d] & \mathscr{H}\ar[rr]^{\mathcal{U}_t=e^{-tA}} & & \mathscr{H} & A^{*}=-A\\
L & \mathscr{K}\ar[rr]^{[\mathcal{S}_{t}]_{t\in\mathbb{R}_+}}_{\mathcal{U}_{t}=e^{-tL}} & & \mathscr{K} & L^{*}=L,\,\, L\geq0
}
\]

\protect\caption{Transformation of skew-adjoint $A$ into selfadjoint semibounded $L$.
\index{operators!semibounded-}}
\end{figure}

\begin{example}[\cite{MR1895530}]
\label{exa:rp2}Fix $0<\sigma<1$, and let $\mathscr{H}\left(=\mathscr{H}_{\sigma}\right)$
be the Hilbert space of all locally integral functions on $\mathbb{R}$
satisfying
\begin{equation}
\left\Vert f\right\Vert ^{2}=\int_{\mathbb{R}}\int_{\mathbb{R}}\overline{f\left(x\right)}f\left(y\right)\left|x-y\right|^{\sigma-1}dxdy<\infty.\label{eq:rp1-1}
\end{equation}
Set 
\begin{eqnarray}
\left(\mathcal{U}\left(t\right)f\right)\left(x\right) & = & e^{\left(\sigma+1\right)t}f\left(e^{2t}x\right),\;\mbox{and}\label{eq:rp1-2}\\
\left(\mathcal{J}f\right)\left(x\right) & = & \left|x\right|^{-\sigma-1}f\left(\frac{1}{x}\right).\label{eq:rp1-3}
\end{eqnarray}
Then $\left\{ \mathcal{U}\left(t\right)\right\} _{t\in\mathbb{R}}$
and $\mathcal{J}$ satisfy the reflection property, i.e., 
\begin{equation}
\mathcal{J}\mathcal{U}\left(t\right)\mathcal{J}=\mathcal{U}\left(-t\right),\;t\in\mathbb{R}\label{eq:rp1-4}
\end{equation}
as operators in $\mathscr{H}$; and $\left\{ \mathcal{U}\left(t\right)\right\} _{t\in\mathbb{R}}$
is a unitary one-parameter group.

We now turn to the ``reflected'' version of the Hilbert norm (\ref{eq:rp1-1}):

The reflection Hilbert space $\mathscr{K}$ will be generated by the
completion of the space of functions $f$ supported in $\left(-1,1\right)$,
that satisfy
\begin{equation}
\left\Vert f\right\Vert _{\mathscr{K}}^{2}=\int_{-1}^{1}\int_{-1}^{1}\overline{f\left(x\right)}f\left(y\right)\left|1-xy\right|^{\sigma-1}dxdy<\infty.\label{eq:rp1-5}
\end{equation}
We show below that this is a Hilbert space of \emph{distributions}.\index{space!Hilbert-}

The selfadjoint contractive semigroup $\left\{ \widetilde{\mathcal{U}}\left(t\right)\right\} _{t\in\mathbb{R}_{+}}$
acting in $\mathscr{K}$ is given by the same formula as in (\ref{eq:rp1-2}),
but now acting in the Hilbert space $\mathscr{K}$ defined by (\ref{eq:rp1-5}).
Note $\widetilde{\mathcal{U}}\left(t\right)$ is only defined for
$t\in\mathbb{R}_{+}\cup\left\{ 0\right\} $. \index{semigroup}\end{example}
\begin{xca}[Renormalization]
\myexercise{Renormalization} ~
\begin{enumerate}
\item Show that the distributions $\left\{ \delta_{0}^{\left(n\right)}\right\} _{n\in\left\{ 0\right\} \cup\mathbb{N}}$
forms an orthogonal and total system in $\mathscr{K}_{\sigma}$ from
(\ref{eq:rp1-5}), for all fixed $0<\sigma<1$.
\item Show that 
\begin{equation}
\left\Vert \delta_{0}^{\left(n\right)}\right\Vert _{\mathscr{K}_{\sigma}}^{2}=n!\left(1-\sigma\right)\left(2-\sigma\right)\cdots\left(n-\sigma\right).\label{eq:rp3-1}
\end{equation}

\end{enumerate}
\end{xca}
\index{distribution!Schwartz-}\index{renormalization}

The idea of reflection positivity originated in physics. Now, when
it is carried out in concrete cases, the initial function spaces change;
but, more importantly, the inner product which produces the respective
Hilbert spaces of quantum states changes as well.

What is especially intriguing is that before reflection we may have
a Hilbert space of functions, but after the time-reflection is turned
on, then, in the new inner product, the corresponding completion magically
becomes a \emph{Hilbert space of distributions}.

Now this is illustrated already in the simple examples above, \exerref{rp1},
and \exaref{rp2}. We include details below to stress the distinction
between an abstract Hilbert-norm completion on the one hand, and a
concretely realized Hilbert space on the other.\index{space!Hilbert-}\index{completion!Hilbert-}

Constructing physical Hilbert spaces entail completions, often a completion
of a suitable space of functions. What can happen is that the completion
may fail to be a Hilbert space of \emph{functions}, but rather a suitable
Hilbert space of distributions.

Recall that a completion, say $\mathscr{H}$ is defined axiomatically,
and the ``real'' secret is revealed only when the elements in $\mathscr{H}$
are identified.

To make the idea more clear we illustrate the point by considering
functions on the interval $-1<x<1$.

Let $C_{c}^{\infty}\left(-1,1\right)$ be the $C^{\infty}$-functions
with compact supports contained in $\left(-1,1\right)$. 

A linear functional $\varphi$ on $C_{c}^{\infty}\left(-1,1\right)$
is said to be a \emph{distribution} if for $\forall$ $K\subset\left(-1,1\right)$
compact, $\forall n\in\mathbb{N}$, $\exists C=C_{K,n}$ such that
\begin{equation}
\left|\varphi\left(f\right)\right|\leq C\:\sup_{x\in K}\:\max_{0\leq j\leq n}\left|\left(\frac{d}{dx}\right)^{j}f\left(x\right)\right|,\;\forall f\in C_{c}^{\infty}\left(-1,1\right).\label{eq:rp2-1}
\end{equation}

Examples of distributions are Dirac \textquotedbl{}functions\textquotedbl{}
$\delta_{x_{0}}$, and the derivatives $\left(\frac{d}{dx}\right)^{n}\delta_{x_{0}}$,
$x_{0}\in\left(-1,1\right)$, are defined by:
\begin{equation}
\left(\left(\frac{d}{dx}\right)^{n}\delta_{x_{0}}\right)\left(f\right)=\left(-1\right)^{n}f^{\left(n\right)}\left(x_{0}\right),\;f\in C_{c}^{\infty}\left(-1,1\right).\label{eq:rp2-2}
\end{equation}
(Note: Distributions are \uline{not} functions, but in Gelfand's
rendition of the theory \cite{MR0435832} they are called \textquotedbl{}generalized
functions.\textquotedbl{})

Now equip $C_{c}^{\infty}\left(-1,1\right)$ with the sesquilinear
form from (\ref{eq:rp1-5}) in \exaref{rp2}, i.e., 
\[
\left\langle f,g\right\rangle _{\mathscr{K}_{\sigma}}:=\int_{-1}^{1}\int_{-1}^{1}\overline{f\left(x\right)}g\left(y\right)\left|1-xy\right|^{\sigma-1}dxdy.
\]

\index{sesquilinear form}
\begin{xca}[A Hilbert space of distributions]
\myexercise{A Hilbert space of distributions}~
\begin{enumerate}
\item \label{enu:rp2-1}Show that each of the distributions $\left(\frac{d}{dx}\right)^{n}\delta_{x_{0}}$,
$n\in\left\{ 0\right\} \cup\mathbb{N}$, $x_{0}\in\left(-1,1\right)$
is in the completion $\mathscr{K}_{\sigma}$ with respect to (\ref{eq:rp1-5}).
\item \label{enu:rp2-2}Compute the Hilbert norm of $\left(\frac{d}{dx}\right)^{n}\delta_{x_{0}}$
in $\mathscr{K}_{\sigma}$, i.e., find 
\begin{equation}
\left\Vert \left(\frac{d}{dx}\right)^{n}\delta_{x_{0}}\right\Vert _{\mathscr{K}_{\sigma}}\label{eq:rp2-3}
\end{equation}
for all $n\in\left\{ 0\right\} \cup\mathbb{N}$, and $x_{0}\in\left(-1,1\right)$.
\end{enumerate}

\uline{Hint}: The answer to (\ref{enu:rp2-2}) (i.e., (\ref{eq:rp2-3}))
is as follows:
\begin{itemize}
\item $n=0$: 
\[
\left\Vert \delta_{x_{0}}\right\Vert _{\mathscr{K}_{\sigma}}^{2}=\left(1-x_{0}^{2}\right)^{\sigma-1};
\]

\item $n=1$ (one derivative):
\[
\left\Vert \delta'_{x_{0}}\right\Vert _{\mathscr{K}_{\sigma}}^{2}=\left(1-\sigma\right)\left(1-x_{0}^{2}\right)^{\sigma-3}\left(1+\left(1-\sigma\right)x_{0}^{2}\right).
\]

\end{itemize}
\end{xca}
\index{distribution!Schwartz-}

\index{completion!Hilbert-}
\begin{xca}[Taylor for distributions]
\myexercise{Taylor for distributions}Fix $0<\sigma<1$, and let
$\mathscr{K}_{\sigma}$ be the corresponding Hilbert space of distributions.
As an identity in $\mathscr{K}_{\sigma}$, establish:
\begin{equation}
\delta_{x}=\sum_{n=0}^{\infty}\frac{\left(-x\right)^{n}}{n!}\delta_{0}^{\left(n\right)},\label{eq:rp3-2}
\end{equation}
valid for all $x$, $\left|x\right|<1$. 
\end{xca}
\begin{flushleft}
Historical note.
\par\end{flushleft}

Laurent Schwartz has developed a systematic study of Hilbert spaces
of distributions; see \cite{MR0179587}.\index{space!Hilbert-}

\section{A second duality principle: A metric on the set of probability measures}

Let $\left(X,d\right)$ be a separable metric space, and denote by
$\mathcal{M}_{1}\left(X\right)$ and $\mathcal{M}_{1}\left(X\times X\right)$
the corresponding sets of regular probability measures. Let $\pi_{i}$,
$i=1,2$, denote the projections: $\pi_{1}\left(x_{1},x_{2}\right)=x_{1}$,
$\pi_{2}\left(x_{1},x_{2}\right)=x_{2}$, for all $\left(x_{1},x_{2}\right)\in X\times X$.
\index{space!metric-}

For $\mu\in\mathcal{M}_{1}\left(X\times X\right)$, set 
\[
\mu^{\pi_{i}}:=\mu\circ\pi_{i}^{-1}.
\]
For $P_{i}\in\mathcal{M}_{1}\left(X\right)$, $i=1,2$, set 
\[
\mathcal{M}\left(P_{1},P_{2}\right)=\left\{ \mu\in\mathcal{M}_{1}\left(X\times X\right)\::\:\mu^{\pi_{i}}=P_{i},\:i=1,2\right\} .
\]
Finally, let $\text{Lip}_{1}=$ the Lipchitz functions on $\left(X,d\right)$,
i.e., $f\in\text{Lip}_{1}$ iff (Def.) 
\[
\left|f\left(x\right)-f\left(y\right)\right|\leq d\left(x,y\right),\;\forall x,y\in X.
\]

\index{Kantorovich-Rubinstein Theorem}
\begin{thm}[Kantorovich-Rubinstein]
 Setting 
\[
dist_{W}\left(P_{1},P_{2}\right)=\inf\left\{ \int_{X\times X}d\left(x,y\right)d\mu\left(x,y\right)\::\:\mu\in\mathcal{M}\left(P_{1},P_{2}\right)\right\} 
\]
and 
\[
dist_{K}\left(P_{1},P_{2}\right)=\sup\left\{ \int_{X}f\,d\left(P_{1}-P_{2}\right)\::\:f\in Lip_{1}\right\} 
\]
then 
\[
dist_{W}\left(P_{1},P_{2}\right)=dist_{K}\left(P_{1},P_{2}\right)
\]
for all $P_{1},P_{2}\in\mathcal{M}_{1}\left(X\right)$.\end{thm}
\begin{proof}
We omit the proof here, but refer to \cite{Kan58,KR57,Rus07}.\end{proof}
\begin{xca}[A complete metric space]
\myexercise{A complete metric space}\label{exer:meas}Show that
$\mathcal{M}_{1}\left(X\right)$ is a complete metric space when equipped
with the metric $dist_{K}$.\index{space!metric-}
\end{xca}

\begin{xca}[A distance formula]
\myexercise{A distance formula}\label{exer:meas1}Let $\left(X,d\right)$
be $\mathbb{R}$ with the usual distance $d\left(x,y\right)=\left|x-y\right|$.
For $P\in\mathcal{M}_{1}\left(\mathbb{R}\right)$ set $F_{P}\left(x\right)=P\left((-\infty,x]\right)$.
Show that then 
\[
dist\left(P_{1},P_{2}\right)=\int_{\mathbb{R}}\left|F_{P_{1}}\left(x\right)-F_{P_{2}}\left(x\right)\right|dx.
\]

\end{xca}

Let $\left(X,d\right)$ and $\mathcal{M}_{1}\left(X\right)$ be as
above. Now apply Banach's Fixed point theorem to the complete metric
space $\left(\mathcal{M}_{1}\left(X\right),dist_{K}\right)$ to get
a solution to the following:\index{Theorem!Banach's fixed-point-}
\begin{xca}[Iterated function systems]
\myexercise{Iterated function systems} Let $N\in\mathbb{N}$, and
let $\varphi_{i}:X\longrightarrow X$, $i=1,\cdots,N$ be a system
of strict contractions in $\left(X,d\right)$. On $\mathcal{M}_{1}\left(X\right)$,
set 
\begin{equation}
T_{\mu}:=\frac{1}{N}\sum_{i=1}^{N}d\mu\circ\varphi_{i}^{-1}.\label{eq:md1}
\end{equation}
Recall $\left(\mu\circ\varphi_{i}^{-1}\right)\left(\triangle\right)=\mu\left(\varphi_{i}^{-1}\left(\triangle\right)\right)$.
\begin{enumerate}
\item \label{enu:md1}Show that, if $c$ is the smallest of the contractivity
constant for $\left\{ \varphi_{i}\right\} _{i=1}^{N}$, then 
\begin{equation}
dist\left(T_{\mu},T_{\nu}\right)\leq c\:dist\left(\mu,\nu\right),\;\forall\mu,\nu\in\mathcal{M}_{1}\left(X\right).\label{eq:md2}
\end{equation}

\item \label{enu:md2}Show that there is a unique solution $\mu_{L}\in\mathcal{M}_{1}\left(X\right)$
to 
\begin{equation}
T\mu_{L}=\mu_{L},\;\mbox{i.e.,}\label{eq:md3}
\end{equation}
\begin{equation}
\int_{X}f\left(x\right)d\mu_{L}\left(x\right)=\frac{1}{N}\sum_{i=1}^{N}\int_{X}f\left(\varphi_{i}\left(x\right)\right)d\mu_{L}\left(x\right)\label{eq:md4}
\end{equation}
holds for $\forall f\in C_{b}\left(X\right)$ (= bounded continuous.)
\end{enumerate}

\uline{Hint}: The desired conclusion in (\ref{enu:md2}), i.e.,
both existence and uniqueness of $\mu_{L}$, follows from Banach's
fixed point theorem: Every strict contraction in a complete metric
space has a unique fixed-point.\index{space!metric-}\index{Theorem!Banach's fixed-point-}

\end{xca}

\begin{figure}
\noindent \begin{centering}
\includegraphics[scale=0.4]{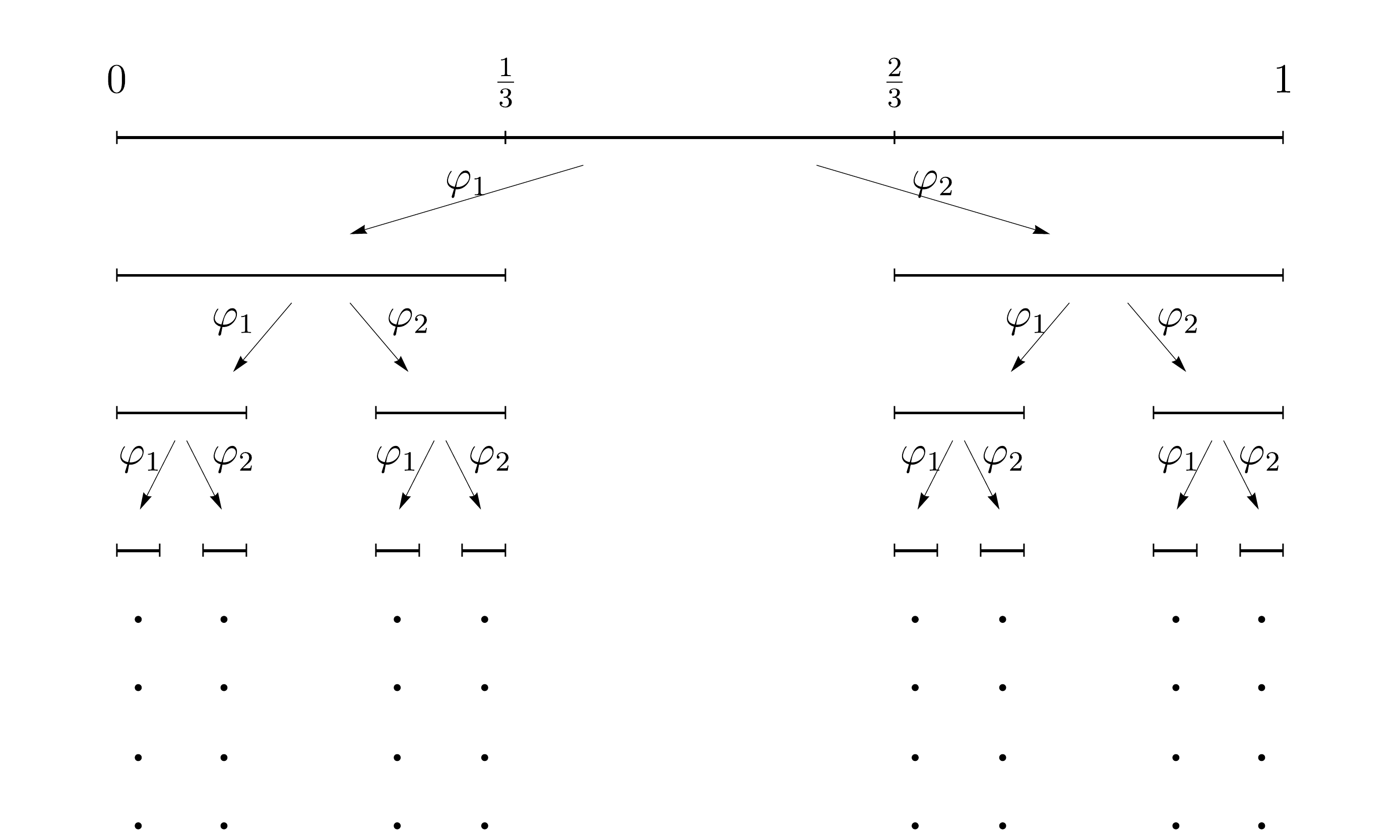}
\par\end{centering}

\protect\caption{\label{fig:can}The Middle-Third Cantor set as a limit.}
\end{figure}

\begin{figure}
\includegraphics[scale=0.4]{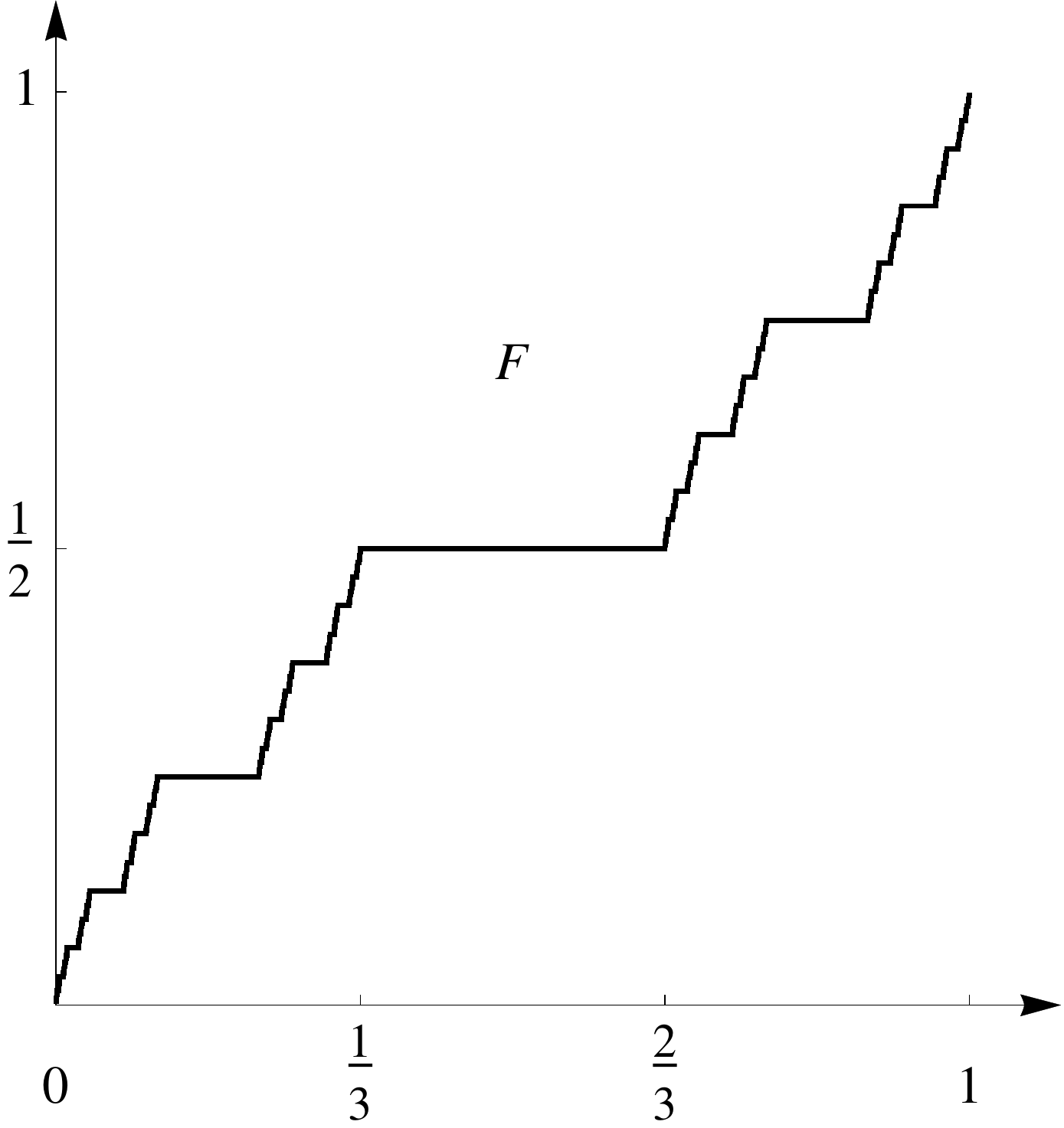}

\protect\caption{\label{fig:ds}The Devil's Staircase.}
\end{figure}

\index{Cantor!-measure}
\begin{xca}[The Middle-Third Cantor-measure]
\label{exer:cantor}\myexercise{The Middle-Third Cantor-measure}Set
$X=\left[0,1\right]=$ the unit interval with the usual metric, set
$N=2$, and 
\begin{equation}
\varphi_{1}\left(x\right)=\frac{x}{3},\quad\varphi_{2}\left(x\right)=\frac{x+2}{3},\label{eq:md5}
\end{equation}
and let $\mu_{L}$ be the corresponding measure, i.e., 
\begin{equation}
\int_{0}^{1}f\left(x\right)d\mu_{L}\left(x\right)=\frac{1}{2}\left(\int_{0}^{1}f\left(\frac{x}{3}\right)d\mu_{L}\left(x\right)+\int_{0}^{1}f\left(\frac{x+2}{3}\right)d\mu_{L}\left(x\right)\right).\label{eq:md6}
\end{equation}
Show that $\mu_{L}$ is supported on the Middle-Third Cantor set (see
\figref{can}). \index{Cantor!-middle third}\end{xca}
\begin{rem}
Starting with the Middle-Third Cantor measure $\mu_{L}$; see (\ref{eq:md6}),
we get the cumulative distribution function $F$ defined on the unit
interval $\left[0,1\right]$, 
\begin{equation}
F\left(x\right)=\mu_{L}\left(\left[0,x\right]\right).\label{eq:acan}
\end{equation}
It follows form \exerref{cantor} that the graph of $F$ is the Devil's
Staircase; see \figref{ds}. Endpoints: $F\left(0\right)=0$, and
$F\left(1\right)=1$. The union $\mathscr{O}$ of all the omitted
open intervals has total length: \index{measure!cantor}
\[
\frac{1}{3}+\frac{2}{3^{2}}+\cdots=\frac{1}{3}\sum_{n=0}^{\infty}\left(\frac{2}{3}\right)^{n}=1,
\]
and $F'\left(x\right)=0$ for all $x\in\mathscr{O}$. 

Using the argument from \exerref{wnorm} above, we get 
\[
\int_{0}^{1}dF=\int_{0}^{1}F'\left(x\right)dx=0.
\]
Since $F\left(1\right)-F\left(0\right)=1$, it would seem that the
conclusion in the Fundamental Theorem of Calculus fails. (Explain!
See e.g., \cite[ch 7]{Rud87}.)
\end{rem}

\begin{xca}[Straightening out the Devil's staircase]
\myexercise{Straightening out the Devil's staircase} Repeat the
construction from the previous exercise, but now with the two functions
$\varphi_{1},\varphi_{2}$ modified as follows:
\begin{equation}
\varphi_{1}\left(x\right)=\frac{x}{2},\quad\varphi_{2}\left(x\right)=\frac{x+1}{2};\label{eq:ds1-1}
\end{equation}
compare with (\ref{eq:md5}) above.
\end{xca}
Then rewrite formula (\ref{eq:md6}), and show that the cumulative
distribution $F$ from (\ref{eq:acan}) becomes (\figref{ds1})
\[
F\left(x\right)=\begin{cases}
0 & \;x<0\\
x & \;0\leq x\leq1\\
0 & \;x>1.
\end{cases}
\]
Explain this!

\noindent \begin{center}
\begin{figure}
\noindent \begin{centering}
\includegraphics[scale=0.4]{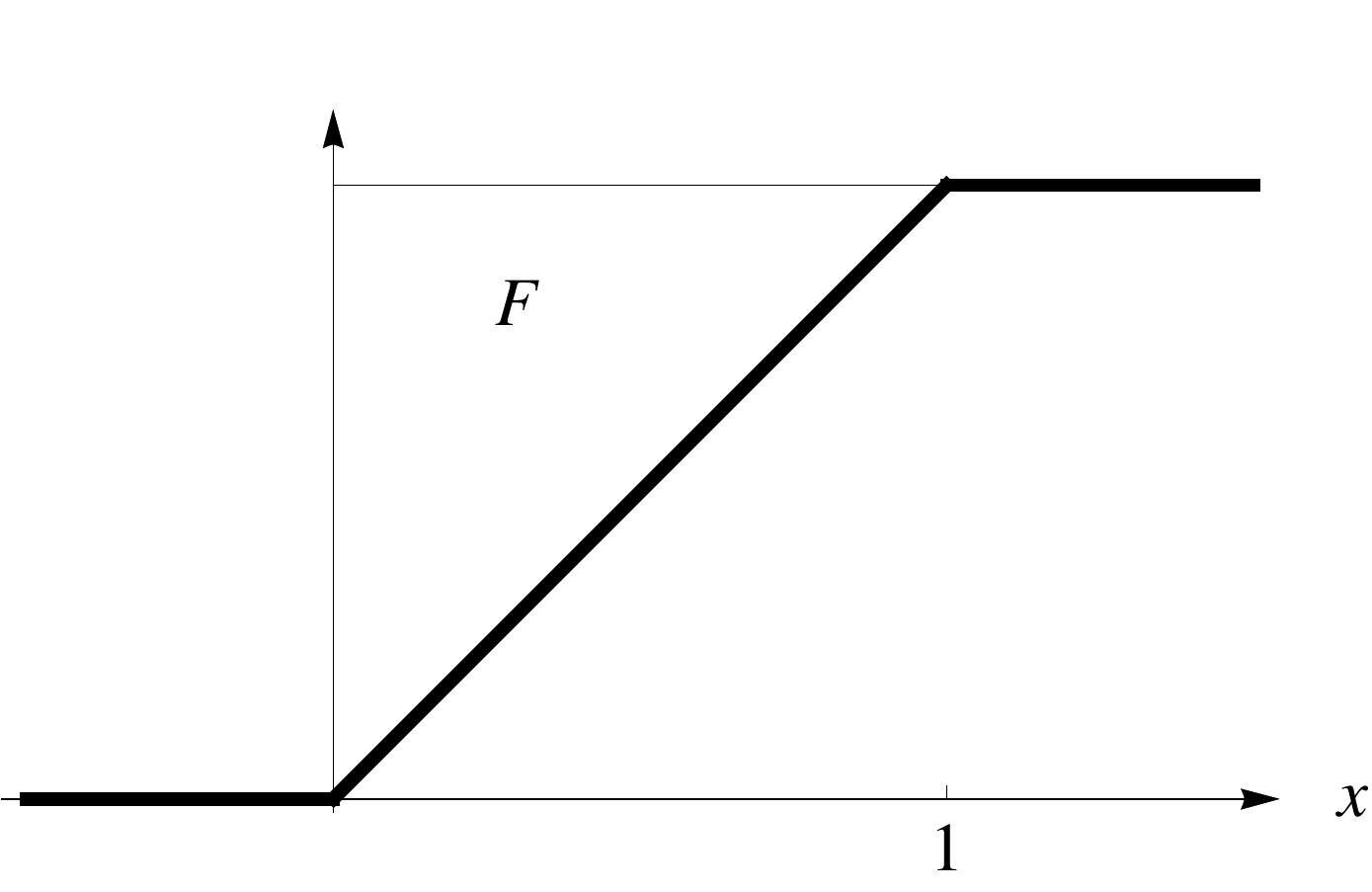}
\par\end{centering}

\protect\caption{\label{fig:ds1}Straightening out the Devil's staircase.}
\end{figure}

\par\end{center}

\section{\label{sec:abelianc}Abelian $C^{*}$-algebras}

Diagonalizing a commuting family of bounded selfadjoint operators
may be formulated in the setting of abelian $C^{*}$-algebras. By
the structure theorem of Gelfand and Naimark, every abelian $C^{*}$-algebra
containing the identity element is isomorphic to the algebra $C\left(X\right)$
of continuous functions on some compact\index{compact} Hausdorff
space $X$, which is unique up to homeomorphism. The classification
of all the representations abelian $C^{*}$-algebras, therefore, amounts
to that of $C(X)$. This problem can be understood using the idea
of $\sigma$-measures (square densities). It also leads to the multiplicity\index{multiplicity}
theory of selfadjoint operators. The best treatment on this subject
can be found in \cite{Ne69}.\index{algebras!$C^{*}$-algebra}\index{diagonalization}
\index{multiplicity!-of spectrum}

Here we discuss \emph{Gelfand's theory on abelian $C^{*}$-algebras}.
Throughout, we assume all the algebras contain unit element. 
\begin{defn}
$\mathfrak{A}$ is \emph{Banach algebra} \index{algebras!Banach algebra}
if it is a complex algebra and a Banach space such that the norm satisfies
$\left\Vert ab\right\Vert \leq\left\Vert a\right\Vert \left\Vert b\right\Vert $,
for all $a,b\in\mathfrak{A}$. 
\end{defn}
Let $\mathfrak{A}$ be an abelian Banach. Consider the closed ideals
in $\mathfrak{A}$ (since $\mathfrak{A}$ is normed, so consider closed
ideals) ordered by inclusion. By Zorn's lemma\index{Zorn's lemma},
there exists maximal ideals $M$, which are closed by maximality.
Then $\mathfrak{A}/M$ is 1-dimensional, i.e., $\mathfrak{A}/M=\{tv\}$
for some $v\in\mathfrak{A}$, and $t\in\mathbb{R}$. Therefore the
combined map
\[
\varphi:\mathfrak{A}\rightarrow\mathfrak{A}/M\rightarrow\mathbb{C},\;a\mapsto a/M\mapsto t_{a}
\]
is a (complex) homomorphism. In particular, $\mathfrak{A}\ni1_{\mathfrak{A}}\mapsto v:=1_{\mathfrak{A}}/M\in\mathfrak{A}/M\simeq\mathbb{C}$,
and $\varphi(1_{\mathfrak{A}})=1$. 

Conversely, the kernel of any homomorphism is a maximal ideal in $\mathfrak{A}$
(since the co-dimension = 1.) Therefore there is a bijection between
maximal ideas and homomorphisms. 

\index{ideal}

\index{kernel!-of operator}
\begin{lem}
\label{lem:Ba1}Let $\mathfrak{A}$ be an abelian Banach algebra.
If $a\in\mathfrak{A}$, and $\left\Vert a\right\Vert <1$, then $1_{\mathfrak{A}}-a$
is invertible. \end{lem}
\begin{proof}
It is easy to verify that $\left(1-a\right)^{-1}=1+a+a^{2}+\cdots$,
and the RHS is norm convergent. \end{proof}
\begin{cor}
Any homomorphism\index{homomorphism} $\varphi:\mathfrak{A}\rightarrow\mathbb{C}$
is a contraction. \end{cor}
\begin{proof}
Let $a\in\mathfrak{A}$, $a\neq0$. Suppose $\lambda:=\varphi\left(a\right)$
such that $\left|\lambda\right|>\left\Vert a\right\Vert $. Then $\left\Vert a/\lambda\right\Vert <1$
and so $1_{\mathfrak{A}}-a/\lambda$ is invertible by \lemref{Ba1}.
Since $\varphi$ is a homomorphism, it must map invertible element
to invertible element, hence $\varphi\left(1_{\mathfrak{A}}-a/\lambda\right)\neq0$,
i.e., $\varphi\left(a\right)\neq\lambda$, which is a contradiction.
\end{proof}
Let $X$ be the set of all \emph{maximal ideals}, identified with
all homomorphisms in $\mathfrak{A}_{1}^{*}$, where $\mathfrak{A}_{1}^{*}$
is the unit ball in $\mathfrak{A}^{*}$. Since $\mathfrak{A}_{1}^{*}$
is compact (see Banach-Alaoglu, \thmref{Alaoglu}), and $X$ is closed
in it, therefore $X$ is also compact. Here, compactness refers to
the weak{*}-topology\index{weak{*}-topology}.
\begin{defn}
The \emph{Gelfand transform} $\mathcal{F}:\mathfrak{A}\rightarrow C(X)$
is given by
\begin{equation}
\mathcal{F}(a)(\varphi)=\varphi(a),\;a\in\mathfrak{A},\varphi\in C\left(X\right).\label{eq:gel1}
\end{equation}
Hence $\mathfrak{A}/\ker\mathcal{F}$ is homomorphic to a closed subalgebra
of $C(X)$. Note $\ker\mathcal{F}=\left\{ a\in\mathfrak{A}:\varphi\left(a\right)=0,\;\forall\varphi\in X\right\} $.
It is called the \emph{radical} of $\mathfrak{A}$.
\end{defn}
The theory is takes a more pleasant form when $\mathfrak{A}$ is a
$C^{*}$-algebra. So there is an involution, and the norm satisfies
the $C^{*}$ axiom: $\left\Vert aa^{*}\right\Vert =\left\Vert a\right\Vert ^{2}$,
for all $a\in\mathfrak{A}$. \index{axioms}
\begin{thm}[Gelfand]
 If $\mathfrak{A}$ is an abelian $C^{*}$-algebra then the Gelfand
transform (\ref{eq:gel1}) is an isometric $*$-isomorphism from $\mathfrak{A}$
onto $C\left(X\right)$, where $X$ is the maximal ideal space of
$\mathfrak{A}$. \end{thm}
\begin{example}
Consider $l^{1}(\mathbb{Z})$, the convolution algebra: 
\begin{eqnarray}
\left(ab\right)_{n} & = & \sum_{k}a_{k}b_{n-k}\label{eq:ca}\\
a_{n}^{*} & = & \overline{a_{-n}}\nonumber \\
\left\Vert a\right\Vert  & = & \sum_{n}\left|a_{n}\right|\nonumber \\
1_{\mathfrak{A}} & = & \delta_{0}\;\left(\mbox{Dirac mass at }0\right);\nonumber 
\end{eqnarray}
the unit-element for the product (\ref{eq:ca}) in $l^{1}\left(\mathbb{Z}\right)$.
\index{convolution}

To identity $X$ in practice, we always start with a guess, and usually
it turns out to be correct. Since Fourier\index{transform!Fourier}
transform converts convolution to multiplication, 
\[
l^{1}\left(\mathbb{Z}\right)\ni a\xrightarrow{\quad\varphi_{z}\quad}\sum a_{n}z^{n}
\]
is a complex homomorphism. To see $\varphi_{z}$ is multiplicative,
we have 
\begin{eqnarray*}
\varphi_{z}(ab) & = & \sum(ab)_{n}z^{n}\\
 & = & \sum_{n,k}a_{k}b_{n-k}z^{n}\\
 & = & \sum_{k}a_{k}z^{k}\sum_{n}b_{n-k}z^{n-k}\\
 & = & \left(\sum_{k}a_{k}z^{k}\right)\left(\sum_{k}b_{k}z^{k}\right)\\
 & = & \varphi_{z}\left(a\right)\varphi_{z}\left(b\right).
\end{eqnarray*}
Thus $\{z:\left|z\right|=1\}$ is a subspace in the Gelfand space\index{Gelfand space}
$X$. Note that we cannot use $\left|z\right|<1$ since we are dealing
with two-sided $l^{1}$ sequence. (If the sequences were truncated,
so that $a_{n}=0$ for $n<0$ then we allow $\left|z\right|<1$. )

$\varphi_{z}$ is contractive: $\left|\varphi_{z}(a)\right|=\left|\sum a_{n}z^{n}\right|\leq\sum_{n}\left|a_{n}\right|=\left\Vert a\right\Vert $. \end{example}
\begin{xca}[The homomorphism of $l^{1}$]
\myexercise{The homomorphism of $l^{1}$}\label{exer:homo}Prove
that every homomorphism of $l^{1}\left(\mathbb{Z}\right)$ is obtained
as $\varphi_{z}$ for some $\left|z\right|=1$. Hence $X=\mathbb{T}^{1}\left(=\left\{ z\in\mathbb{C}:\left|z\right|=1\right\} \right)$.\end{xca}
\begin{example}
$l^{\infty}(\mathbb{Z})$, with $\left\Vert a\right\Vert =\sup_{n}\left|a_{n}\right|$.
The Gelfand space\index{Gelfand space} in this case is $X=\beta\mathbb{Z}$,
the Stone-\v{C}ech compactification\index{Stone-v{C}ech compactification@Stone-\v{C}ech compactification}
of $\mathbb{Z}$, which are the ultra-filters on $\mathbb{Z}$. $\beta\mathbb{Z}$
is much bigger then $p$-adic numbers. Pure states on diagonal operators
correspond to $\beta\mathbb{Z}$. See \chapref{KS} for details. \index{Stone, M. H.}
\index{state!pure-}
\end{example}

\section{\label{sec:sr}States and Representations}

Let $\mathfrak{A}$ be a $*$-algebra, a representation $\pi:\mathfrak{A}\rightarrow\mathscr{B}(\mathscr{H})$
generates a $*$-subalgebra $\pi(\mathfrak{A})$ in $\mathscr{B}(\mathscr{H})$.
By taking norm closure, one gets a $C^{*}$-algebra, i.e., a Banach
$*$-algebra with the axiom $\left\Vert a^{*}a\right\Vert =\left\Vert a\right\Vert ^{2}$.
On the other hand, by Gelfand and Naimark's theorem, all abstract
$C^{*}$-algebras are isometrically isomorphic to closed subalgebras
of $\mathscr{B}(\mathscr{H})$, for some Hilbert space $\mathscr{H}$
(\thmref{gns1}). The construction of $\mathscr{H}$ comes down to
states $S(\mathfrak{A})$ on $\mathfrak{A}$ and the GNS construction.
Therefore, the GNS construction gives rise to a bijection between
states and representations. \index{representation!GNS}

Let $\mathfrak{A}_{+}$ be the positive elements in $\mathfrak{A}$.
$s\in S(\mathfrak{A})$, $s:\mathfrak{A}\rightarrow\mathbb{C}$ and
$s(\mathfrak{A}_{+})\subset[0,\infty)$. For $C^{*}$-algebra, positive
elements can be written $f=(\sqrt{f})^{2}$ by the spectral theorem.
In general, positive elements have the form $a^{*}a$. There is a
bijection between states and GNS representations $Rep(\mathfrak{A},\mathscr{H})$,
where $s(A)=\left\langle \Omega,\pi(A)\Omega\right\rangle $.
\begin{example}
$\mathfrak{A}=C(X)$ where $X$ is a compact Hausdorff space. $s_{\mu}$
given by $s_{\mu}(a)=\int ad\mu$ is a state. The GNS construction
gives $\mathscr{H}=L^{2}(\mu)$, $\pi(f)$ is the operator of multiplication
by $f$ on $L^{2}(\mu)$. $\{\varphi1:\varphi\in C(X)\}$ is dense
in $L^{2}$, where $1$ is the cyclic vector. $s_{\mu}(f)=\left\langle \Omega,\pi(f)\Omega\right\rangle =\int1f1d\mu=\int fd\mu$,
which is also seen as the expectation of $f$ in case $\mu$ is a
probability measure.
\end{example}
We consider decomposition of representations or equivalently states,
i.e., breaking up representations corresponds to breaking up states. 

The thing that we want to do with representations comes down to the
smallest ones, i.e., the irreducible representations. Irreducible\index{representation!irreducible}
representations correspond to pure states which are extreme points
in the states (see \cite{Phe01}). \index{extreme-point} \index{pure state}
\index{state!pure-}\index{Theorem!Spectral-}

A representation $\pi:\mathfrak{A}\rightarrow\mathscr{B}(\mathscr{H})$
is irreducible, if whenever $\mathscr{H}$ breaks up into two pieces
$\mathscr{H}=\mathscr{H}_{1}\oplus\mathscr{H}_{2}$, where $\mathscr{H}_{i}$
is invariant under $\pi(\mathfrak{A})$, one of them is zero (the
other is $\mathscr{H}$). Equivalently, if $\pi=\pi_{1}\oplus\pi_{2}$,
where $\pi_{i}=\pi\big|_{\mathscr{H}_{i}}$, then one of them is zero.
This is similar to the decomposition of natural numbers into product
of primes. For example, $6=2\times3$, but $2$ and $3$ are primes
and they do not decompose further. 

Hilbert spaces are defined up to unitary equivalence. A state $\varphi$
may have equivalent representations on different Hilbert spaces (but
unitarily equivalent), however $\varphi$ does not see the distinction,
and it can only detect equivalent classes of representations. 
\begin{example}
Let $\mathfrak{A}$ be a $*$-algebra. Given two states $s_{1}$ and
$s_{2}$, by the GNS construction, we get cyclic vectors $\xi_{i}$,
and representations $\pi_{i}:\mathfrak{A}\rightarrow\mathscr{B}(\mathscr{H}_{i})$,
so that $s_{i}(A)=\left\langle \xi_{i},\pi_{i}(A)\xi_{i}\right\rangle $,
$i=1,2$. Suppose there is a unitary operator $W:\mathscr{H}_{1}\rightarrow\mathscr{H}_{2}$,
such that for all $A\in\mathfrak{A}$, 
\[
\pi_{1}(A)=W^{*}\pi_{2}(A)W.
\]
Then
\begin{eqnarray*}
\left\langle \xi_{1},\pi_{1}\left(A\right)\xi_{1}\right\rangle _{1} & = & \left\langle \xi_{1},W^{*}\pi_{2}\left(A\right)W\xi_{1}\right\rangle _{1}\\
 & = & \left\langle W\xi_{1},\pi_{2}\left(A\right)W\xi_{1}\right\rangle _{2}\\
 & = & \left\langle \xi_{2},\pi_{2}\left(A\right)\xi_{2}\right\rangle _{2},\;\forall A\in\mathfrak{A};
\end{eqnarray*}
i.e., $s_{2}(A)=s_{1}(A)$. Therefore the same state $s=s_{1}=s_{2}$
has two distinct (unitarily equivalent) representations. \end{example}
\begin{rem}
A special case of states are measures when the algebra is abelian.
Recall that all abelian $C^{*}$-algebras with identity are $C\left(X\right)$,
where $X$ is the corresponding Gelfand space\index{Gelfand space}.
Two representations are mutually singular $\pi_{1}\perp\pi_{2}$,
if and only if the two measures are mutually singular, $\mu_{1}\perp\mu_{2}$.
\index{spectrum!singular-}
\end{rem}
The theorem below is fundamental in representation theory. Recall
that if $M$ is a subset of $\mathscr{B}\left(\mathscr{H}\right)$,
the \emph{commutant} $M'$ consists of $A\in\mathscr{B}\left(\mathscr{H}\right)$
that commutes with all elements in $M$. \index{commutant}
\begin{thm}[Schur]
\label{thm:GNSReprep_Schur} Let $\pi:\mathfrak{A}\rightarrow\mathscr{B}(\mathscr{H})$
be a representation. The following are equivalent.
\begin{enumerate}
\item \label{enu:schur-1-1} $\pi$ is irreducible.
\item \label{enu:schur-1-2}The commutant $(\pi(\mathfrak{A}))'$ is one-dimensional,
i.e., $(\pi(\mathfrak{A}))'=cI_{\mathfrak{A}}$, $c\in\mathbb{C}$. 
\end{enumerate}
\end{thm}
\begin{proof}
Suppose $(\pi(\mathfrak{A}))'$ has more than one dimension. Let $X\in(\pi(\mathfrak{A}))'$,
then by taking adjoint, $X^{*}\in(\pi(\mathfrak{A}))'$. $X+X^{*}$
is selfadjoint, and $X+X^{*}\neq cI$ since by hypothesis $(\pi(\mathfrak{A}))'$
has more than one dimension. Therefore $X+X^{*}$ has a non trivial
spectral projection $P(E)$, i.e., $P(E)\notin\{0,I\}$. Let $\mathscr{H}_{1}=P(E)\mathscr{H}$
and $\mathscr{H}_{2}=(I-P(E))\mathscr{H}$. $\mathscr{H}_{1}$ and
$\mathscr{H}_{2}$ are both nonzero proper subspaces of $\mathscr{H}$.
Since $P(E)$ commutes with $\pi(A)$, for all $A\in\mathfrak{A}$,
it follows that $\mathscr{H}_{1}$ and $\mathscr{H}_{2}$ are both
invariant under $\pi$.

Conversely, suppose $(\pi(\mathfrak{A}))'$ is one-dimensional. If
$\pi$ is not irreducible, i.e., $\pi=\pi_{1}\oplus\pi_{2}$, then
for 
\[
P_{\mathscr{H}_{1}}=\left[\begin{array}{cc}
I_{\mathscr{H}_{1}} & 0\\
0 & 0
\end{array}\right],\quad P_{\mathscr{H}_{2}}=1-P_{\mathscr{H}_{1}}=\left[\begin{array}{cc}
0 & 0\\
0 & I_{\mathscr{H}_{2}}
\end{array}\right]
\]
we have
\[
P_{\mathscr{H}_{i}}\pi(A)=\pi(A)P_{\mathscr{H}_{i}},\;i=1,2
\]
for all $A\in\mathfrak{A}$. Hence $(\pi(\mathfrak{A}))'$ has more
than one dimension. \end{proof}
\begin{cor}
$\pi$ is irreducible if and only if the only projections in $(\pi(\mathfrak{A}))'$
are $0$ or $I$.
\end{cor}
Thus to test invariant subspaces, one only needs to look at projections
in the commutant. 
\begin{cor}
If $\mathfrak{A}$ is abelian, then $\pi$ is irreducible if and only
if $\mathscr{H}$ is one-dimensional.\end{cor}
\begin{proof}
Obviously, if $\dim\mathscr{H}=1$, $\pi$ is irreducible. Conversely,
by \thmref{GNSReprep_Schur}, $(\pi(\mathfrak{A}))'=cI$. Since $\pi(\mathfrak{A})$
is abelian, $\pi(\mathfrak{A})\subset\pi(\mathfrak{A})'$. Thus for
all $A\in\mathfrak{A}$, $\pi(A)=c_{A}I$, for some constant $c_{A}$. 
\end{proof}
If instead of taking the norm closure, but using the strong operator
topology, ones gets a von Neumann algebra\index{algebras!von Neumann algebra ($W^{*}$-algebra)}.
von  Neumann showed that the weak closure of $\mathfrak{A}$ is equal
to $\mathfrak{A}''$. 
\begin{cor}
$\pi$ is irreducible $\Longleftrightarrow$ $(\pi(\mathfrak{A}))'$
is 1-dimensional $\Longleftrightarrow$ $(\pi(\mathfrak{A}))''=\mathscr{B}(\mathscr{H})$.\end{cor}
\begin{rem}
In matrix notation, we write $\pi=\pi_{1}\oplus\pi_{2}$ as \index{matrix!block-}
\[
\pi(A)=\left[\begin{array}{cc}
\pi_{1}(A) & 0\\
0 & \pi_{2}(A)
\end{array}\right].
\]
If 
\[
\left[\begin{array}{cc}
X & Y\\
U & V
\end{array}\right]\in(\pi(\mathfrak{A}))'
\]
then 
\begin{eqnarray*}
\left[\begin{array}{cc}
X & Y\\
U & V
\end{array}\right]\left[\begin{array}{cc}
\pi_{1}(A) & 0\\
0 & \pi_{2}(A)
\end{array}\right] & = & \left[\begin{array}{cc}
X\pi_{1}(A) & Y\pi_{2}(A)\\
U\pi_{1}(A) & V\pi_{2}(A)
\end{array}\right]\\
\left[\begin{array}{cc}
\pi_{1}(A) & 0\\
0 & \pi_{2}(A)
\end{array}\right]\left[\begin{array}{cc}
X & Y\\
U & V
\end{array}\right] & = & \left[\begin{array}{cc}
\pi_{1}(A)X & \pi_{1}(A)Y\\
\pi_{2}(A)U & \pi_{2}(A)V
\end{array}\right].
\end{eqnarray*}
Hence 
\begin{eqnarray*}
X\pi_{1}(A) & = & \pi_{1}(A)X\\
V\pi_{2}(A) & = & \pi_{2}(A)V\\
U\pi_{1}(A) & = & \pi_{2}(A)U\\
Y\pi_{2}(A) & = & \pi_{1}(A)Y.
\end{eqnarray*}
Therefore, 
\begin{eqnarray*}
X & \in & (\pi_{1}(\mathfrak{A}))',\;V\in(\pi_{2}(\mathfrak{A}))',\;\mbox{and}\\
U,Y & \in & int\left(\pi_{1},\pi_{2}\right)=\mbox{intertwining operators of }\pi_{1},\pi_{2}.
\end{eqnarray*}
This is illustrated by the diagram below. 

\[
\xymatrix{
&\mathscr{H}_{1}\ar[r]^{\pi_{1}(A)}\ar[d]^{U} & \mathscr{H}_{1}\ar[d]_{U}&\\
&\mathscr{H}_{2}\ar[r]^{\pi_{2}(A)}\ar@/^1pc/@{.>}[lur]^{Y} & \mathscr{H}_{2}\ar@/_1pc/@{.>}[rul]_{Y} &\\
}
\]

We say $\pi_{1}$ and $\pi_{2}$ are \emph{inequivalent} if and only
if $int(\pi_{1},\pi_{2})=0$. For $\pi_{1}=\pi_{2}$, $\pi$ has \emph{multiplicity}
2. Multiplicity > 1 is equivalent to the commutant being non-abelian.
In the case $\pi=\pi_{1}\oplus\pi_{2}$ where $\pi_{1}=\pi_{2}$,
$(\pi(\mathfrak{A}))'\simeq M_{2}\left(\mathbb{C}\right)$.\index{commutant}
\index{multiplicity!-of spectrum}
\end{rem}

Schur's lemma addresses all representations. It says that a representation
$\pi:\mathfrak{A}\rightarrow\mathscr{B}(\mathscr{H})$ is irreducible
if and only if $(\pi(\mathfrak{A}))'$ is 1-dimensional. When specialize
to the GNS representation of a given state $s$, this is also equivalent
to saying that for all positive linear functional\index{functional}
$t$, $t\leq s\Rightarrow t=\lambda s$ for some $\lambda\geq0$.
This latter equivalence is obtained by using a more general result,
which relates $t$ and selfadjoint operators in the commutant $(\pi(\mathfrak{A}))'$.
\index{operators!adjoint-}

We now turn to characterize the relation between state and its GNS
representation, i.e., specialize to the GNS representation. Given
a $*$ algebra $\mathfrak{A}$, the states $S(\mathfrak{A})$ forms
a compact convex subset in the unit ball of the dual\index{dual}
$\mathfrak{A}^{*}$. 

Let $\mathfrak{A}_{+}$ be the set of positive elements in $\mathfrak{A}$.
Given $s\in S(\mathfrak{A})$, let $t$ be a positive linear functional.
By $t\leq s$, we mean $t(A)\leq s(A)$ for all $A\in\mathfrak{A}_{+}$.
We look for relation between $t$ and the commutant $(\pi(\mathfrak{A}))'$. 
\begin{lem}[Schur-Sakai-Nicodym]
 Let $t$ be a positive linear functional, and let $s$ be a state.
There is a bijection between $t$ such that $0\leq t\leq s$, and
selfadjoint operator $A$ in the commutant with $0\leq A\leq I$.
The relation is given by \index{Schur-Sakai-Nicodym} 
\[
t(\cdot)=\left\langle \Omega,\pi(\cdot)A\Omega\right\rangle 
\]
\end{lem}
\begin{rem}
This is an extension of the classical Radon-Nikodym derivative\index{Radon-Nikodym derivative}
theorem to the non-commutative setting. We may write $A=dt/ds$. The
notation $0\leq A\leq I$ refers to the partial order of selfadjoint
operators. It means that for all $\xi\in\mathscr{H}$, $0\leq$. See
\cite{MR0442701,MR1468230}. \index{selfadjoint operator} \index{Theorem!Radon-Nikodym-}\end{rem}
\begin{proof}
Easy direction, suppose $A\in(\pi(\mathfrak{A}))'$ and $0\leq A\leq I$.
As in many applications, the favorite functions one usually applies
to selfadjoint operators is the square root function $\sqrt{\cdot}$.
So let's take $\sqrt{A}$. Since $A\in(\pi(\mathfrak{A}))'$, so is
$\sqrt{A}$. We need to show $t(a)=\left\langle \Omega,\pi(a)A\Omega\right\rangle \leq s(a)$,
for all $a\geq0$ in $\mathfrak{A}$. Let $a=b^{2}$, then
\begin{eqnarray*}
t(a) & = & \left\langle \Omega,\pi(a)A\Omega\right\rangle \\
 & = & \left\langle \Omega,\pi(b^{2})A\Omega\right\rangle \\
 & = & \left\langle \Omega,\pi(b)^{*}\pi(b)A\Omega\right\rangle \\
 & = & \left\langle \pi(b)\Omega,A\pi(b)\Omega\right\rangle \\
 & \leq & \left\langle \pi(b)\Omega,\pi(b)\Omega\right\rangle \\
 & = & \left\langle \Omega,\pi(a)\Omega\right\rangle \\
 & = & s(a).
\end{eqnarray*}
Conversely, suppose $t\leq s$. Then for all $a\geq0$, $t(a)\leq s(a)=\left\langle \Omega,\pi(a)\Omega\right\rangle $.
Again write $a=b^{2}$. It follows that 
\[
t(b^{2})\leq s(b^{2})=\left\langle \Omega,\pi\left(a\right)\Omega\right\rangle =\left\Vert \pi\left(b\right)\Omega\right\Vert ^{2}.
\]
By Riesz's theorem, there is a unique $\eta$, so that 
\[
t(a)=\left\langle \pi(b)\Omega,\eta\right\rangle .
\]

Conversely, let $a=b^{2}$, then 
\[
t(b^{2})\leq s(b^{2})=\left\langle \Omega,\pi\left(a\right)\Omega\right\rangle =\left\Vert \pi\left(b\right)\Omega\right\Vert ^{2}
\]
i.e., $\pi(b)\Omega\mapsto t(b^{2})$ is a bounded quadratic form.
Therefore, there exists a unique $A\geq0$ such that \index{quadratic form}
\[
t(b^{2})=\left\langle \pi(b)\Omega,A\pi(b)\Omega\right\rangle .
\]
It is easy to see that $0\leq A\leq I$. Also, $A\in(\pi(\mathfrak{A}))'$,
the commutant of $\pi\left(\mathfrak{A}\right)$. \end{proof}
\begin{cor}
Let $s$ be a state. $(\pi,\Omega,\mathscr{H})$ is the corresponding
GNS construction. The following are equivalent.
\begin{enumerate}
\item For all positive linear functional\index{functional} $t$, $t\leq s\Rightarrow t=\lambda s$
for some $\lambda\geq0$.
\item $\pi$ is irreducible.
\end{enumerate}
\end{cor}
\begin{proof}
By Sakai-Nicodym derivative, $t\leq s$ if and only if there is a
selfadjoint operator $A\in(\pi(\mathfrak{A}))'$ so that 
\[
t(\cdot)=\left\langle \Omega,\pi(\cdot)A\Omega\right\rangle 
\]
Therefore $t=\lambda s$ if and only if $A=\lambda I$.

Suppose $t\leq s\Rightarrow t=\lambda s$ for some $\lambda\geq0$.
Then $\pi$ must be irreducible, since otherwise there exists $A\in(\pi(\mathfrak{A}))'$
with $A\neq cI$, hence $\mathfrak{A}\ni a\mapsto t(a):=\left\langle \Omega,\pi(a)A\Omega\right\rangle $
defines a positive linear functional, and $t\leq s$, however $t\neq\lambda s$.
Thus a contradiction to the hypothesis.

Conversely, suppose $\pi$ is irreducible. Then by Schur's lemma,
$(\pi(\mathfrak{A}))'$ is 1-dimensional. i.e. for all $A\in(\pi(\mathfrak{A}))'$,
$A=\lambda I$ for some $\lambda$. Therefore if $t\leq s$, by Sakai's
theorem, $t(\cdot)=\left\langle \Omega,\pi(\cdot)A\Omega\right\rangle $.
Thus $t=\lambda s$ for some $\lambda\geq0$.\end{proof}
\begin{defn}
A state $s$ is pure if it cannot be broken up into a convex combination
of two distinct states. i.e. for all states $s_{1}$ and $s_{2}$,
$s=\lambda s_{1}+(1-\lambda)s_{2}\Rightarrow s=s_{1}\mbox{ or }s=s_{2}$.
\index{space!state-}
\end{defn}
The main theorem in this section is a corollary to Sakai's theorem.
\begin{cor}
Let $s$ be a state. $(\pi,\Omega,\mathscr{H})$ is the corresponding
GNS construction. The following are equivalent.
\begin{enumerate}
\item \label{enu:ir1}$t\leq s\Rightarrow t=\lambda s$ for some $\lambda\geq0$.
\item \label{enu:ir2}$\pi$ is irreducible.
\item \label{enu:ir3}$s$ is a pure state.
\end{enumerate}
\end{cor}
\begin{proof}
By Sakai-Nicodym derivative, $t\leq s$ if and only if there is a
selfadjoint operator $A\in(\pi(\mathfrak{A}))'$ so that \index{pure state}
\index{state!pure-} 
\[
t(a)=\left\langle \Omega,\pi\left(a\right)A\Omega\right\rangle ,\;\forall a\in\mathfrak{A}.
\]
Therefore $t=\lambda s$ if and only if $A=\lambda I$.

We show that (\ref{enu:ir1})$\Leftrightarrow$(\ref{enu:ir2}) and
(\ref{enu:ir1})$\Rightarrow$(\ref{enu:ir3})$\Rightarrow$(\ref{enu:ir2}).

(\ref{enu:ir1})$\Leftrightarrow$(\ref{enu:ir2}) Suppose $t\leq s\Rightarrow t=\lambda s$,
then $\pi$ must be irreducible, since otherwise there exists $A\in(\pi(\mathfrak{A}))'$
with $A\neq cI$, hence $t(\cdot):=\left\langle \Omega,\pi(\cdot)A\Omega\right\rangle $
defines a positive linear functional with $t\leq s$, however $t\neq\lambda s$.
Conversely, suppose $\pi$ is irreducible. If $t\leq s$, then $t(\cdot)=\left\langle \Omega,\pi(\cdot)A\Omega\right\rangle $
with $A\in(\pi(\mathfrak{A}))'$. By Schur's lemma, $(\pi(\mathfrak{A}))'=\{0,\lambda I\}$.
Therefore, $A=\lambda I$ and $t=\lambda s$. 

(\ref{enu:ir1})$\Rightarrow$(\ref{enu:ir3}) Suppose $t\leq s\Rightarrow t=\lambda s$
for some $\lambda\geq0$. If $s$ is not pure, then $s=cs_{1}+(1-c)s_{2}$
where $s_{1},s_{2}$ are states and $c\in(0,1)$. By hypothesis, $s_{1}\leq s$
implies that $s_{1}=\lambda s$. It follows that $s=s_{1}=s_{2}$. 

(\ref{enu:ir3})$\Rightarrow$(\ref{enu:ir2}) Suppose $\pi$ is not
irreducible, i.e. there is a non trivial projection $P\in(\pi(\mathfrak{A}))'$.
Let $\Omega=\Omega_{1}\oplus\Omega_{2}$ where $\Omega_{1}=P\Omega$
and $\Omega_{2}=(I-P)\Omega$. Then
\begin{eqnarray*}
s(a) & = & \left\langle \Omega,\pi\left(a\right)\Omega\right\rangle \\
 & = & \left\langle \Omega_{1}\oplus\Omega_{2},\pi\left(a\right)\Omega_{1}\oplus\Omega_{2}\right\rangle \\
 & = & \left\langle \Omega_{1},\pi\left(a\right)\Omega_{1}\right\rangle +\left\langle \Omega_{2},\pi\left(a\right)\Omega_{2}\right\rangle \\
 & = & \left\Vert \Omega_{1}\right\Vert ^{2}\left\langle \frac{\Omega_{1}}{\left\Vert \Omega_{1}\right\Vert },\pi\left(a\right)\frac{\Omega_{1}}{\left\Vert \Omega_{1}\right\Vert }\right\rangle +\left\Vert \Omega_{2}\right\Vert ^{2}\left\langle \frac{\Omega_{2}}{\left\Vert \Omega_{2}\right\Vert },\pi\left(a\right)\frac{\Omega_{2}}{\left\Vert \Omega_{2}\right\Vert }\right\rangle \\
 & = & \left\Vert \Omega_{1}\right\Vert ^{2}\left\langle \frac{\Omega_{1}}{\left\Vert \Omega_{1}\right\Vert },\pi\left(a\right)\frac{\Omega_{1}}{\left\Vert \Omega_{1}\right\Vert }\right\rangle +\left(1-\left\Vert \Omega_{1}\right\Vert ^{2}\right)\left\langle \frac{\Omega_{2}}{\left\Vert \Omega_{2}\right\Vert },\pi\left(a\right)\frac{\Omega_{2}}{\left\Vert \Omega_{2}\right\Vert }\right\rangle \\
 & = & \lambda s_{1}\left(a\right)+\left(1-\lambda\right)s_{2}\left(a\right).
\end{eqnarray*}
Hence $s$ is not a pure state.
\end{proof}

\subsection{Normal States}

More general states in physics come from the mixture of particle states,
which correspond to composite system. These are called normal states
in mathematics. \index{normal state}\index{state!normal-}

Let $\rho\in\mathscr{T}_{1}\left(\mathscr{H}\right)=$ trace class
operator, such that $\rho>0$ and $tr(\rho)=1$. Define state $s_{\rho}(A):=tr(A\rho)$
, $A\in\mathscr{B}\left(\mathscr{H}\right)$. Since $\rho$ is compact,
by spectral theorem of compact operators, 
\[
\rho=\sum_{k}\lambda_{k}P_{k}
\]
such that $\lambda_{1}>\lambda_{2}>\cdots\rightarrow0$; $\sum\lambda_{k}=1$
and $P_{k}=\left|\xi_{k}\left\rangle \right\langle \xi_{k}\right|$,
i.e., the rank-1 projections. (See \secref{predual}.) We have 
\begin{itemize}
\item $s_{\rho}(I)=tr(\rho)=1$; and 
\item for all $A\in\mathscr{B}\left(\mathscr{H}\right)$, 
\begin{eqnarray*}
s_{\rho}(A)=tr(A\rho) & = & \sum_{n}\left\langle u_{n},A\rho u_{n}\right\rangle \\
 & = & \sum_{n}\left\langle A^{*}u_{n},\rho u_{n}\right\rangle \\
 & = & \sum_{n}\sum_{k}\lambda_{k}\left\langle A^{*}u_{n},\xi_{k}\right\rangle \left\langle \xi_{k},u_{n}\right\rangle \\
 & = & \sum_{k}\lambda_{k}\left(\sum_{n}\left\langle u_{n},A\xi_{k}\right\rangle \left\langle \xi_{k},u_{n}\right\rangle \right)\\
 & = & \sum_{k}\lambda_{k}\left\langle \xi_{k},A\xi_{k}\right\rangle ;
\end{eqnarray*}
where $\left\{ u_{k}\right\} $ is any ONB in $\mathscr{H}$. Hence,
\[
s_{\rho}=\sum_{k}\lambda_{k}s_{\xi_{k}}=\sum_{k}\lambda_{k}\left|\xi_{k}\left\rangle \right\langle \xi_{k}\right|
\]
i.e., $s_{\rho}$ is a convex combination of pure states $s_{\xi_{k}}:=\left|\xi_{k}\left\rangle \right\langle \xi_{k}\right|$.\end{itemize}
\begin{rem}
Notice that $tr\left(\left|\xi\left\rangle \right\langle \eta\right|\right)=\left\langle \eta,\xi\right\rangle $.
In fact, take any ONB $\{e_{n}\}$ in $\mathscr{H}$, then 
\[
tr\left(\left|\xi\left\rangle \right\langle \eta\right|\right)=\sum_{n}\left\langle e_{n}\xi\right\rangle \left\langle \eta,e_{n}\right\rangle =\left\langle \eta,\xi\right\rangle 
\]
where the last step follows from Parseval identity. (If we drop the
condition $\rho\geq0$ then we get the duality $(\mathscr{T}_{1}\mathscr{H})^{*}=\mathscr{B}(\mathscr{H})$.
See \thmref{predual}.) \index{Parseval identity}\index{state!pure-}
\end{rem}

\subsection{A Dictionary of OT and QM\protect\footnote{The abbreviation OT is for operator theory, and QM for quantum mechanics.}}
\begin{itemize}
\item states - unit vectors $\xi\in\mathscr{H}$. These are all the pure
(normal) states on $\mathscr{B}(\mathscr{H})$. \index{pure state}\index{state!pure-}
\item observable\index{observable} - selfadjoint\index{operators!selfadjoint}
operators $A=A^{*}$
\item measurement - spectrum
\end{itemize}
The spectral theorem was developed by J. von Neumann and later improved
by Dirac and others. (See \cite{Sto90, Yos95, Ne69, RS75, DS88b}.)
A selfadjoint operator $A$ corresponds to a quantum observable, and
result of a quantum measurement can be represented by the spectrum
of $A$. \index{Spectral Theorem}
\begin{itemize}
\item simple eigenvalue: $A=\lambda\left|\xi_{\lambda}\left\rangle \right\langle \xi_{\lambda}\right|$,\index{eigenvalue}
\[
s_{\xi_{\lambda}}(A)=\left\langle \xi_{\lambda},A\xi_{\lambda}\right\rangle 
\]

\item compact operator: $A=\sum_{\lambda}\lambda\left|\xi_{\lambda}\left\rangle \right\langle \xi_{\lambda}\right|$,
such that $\left\{ \xi_{\lambda}\right\} $ is an ONB of $\mathscr{H}$.
If $\xi=\sum c_{\lambda}\xi_{\lambda}$ is a unit vector, then
\[
s_{\xi}(A)=\sum_{\lambda}\lambda\left\langle \xi_{\lambda},A\xi_{\lambda}\right\rangle 
\]
where $\{|c_{\lambda}|^{2}\}_{\lambda}$ is a probability distribution
over the spectrum of $A$, and $s_{\xi}$ is the expectation value
of $A$. \index{distribution!probability-}
\item more general, allowing continuous spectrum: 
\begin{eqnarray*}
A & = & \int\lambda E(d\lambda)\\
A\xi & = & \int\lambda E(d\lambda)\xi.
\end{eqnarray*}
We may write the unit vector $\xi$ as 
\[
\xi=\int\overset{\xi_{\lambda}}{\overbrace{E(d\lambda)\xi}}
\]
so that
\[
\left\Vert \xi\right\Vert ^{2}=\int\left\Vert E(d\lambda)\xi\right\Vert ^{2}=1
\]
It is clear that $\left\Vert E(\cdot)\xi\right\Vert ^{2}$ is a probability
distribution on spectrum of $A$. $s_{\xi}(A)$ is again seen as the
expectation value of $A$ with respect to $\left\Vert E(\cdot)\xi\right\Vert ^{2}$,
since
\[
s_{\xi}(A)=\left\langle \xi,A\xi\right\rangle =\int\lambda\left\Vert E(d\lambda)\xi\right\Vert ^{2}.
\]

\end{itemize}
\index{spectrum!continuous}

\section{\label{sec:KMil}Krein-Milman, Choquet, Decomposition of States}

\index{Krein-Milman}\index{Choquet}\index{decomposition of state}

We study some examples of compact\index{compact} convex\index{convex}
sets in locally convex topological space\index{locally convex topological space}s.\footnote{Almost all spaces one works with are locally convex.}
Typical examples include the set of positive semi-definite functions,
taking values in $\mathbb{C}$ or $\mathscr{B}(\mathscr{H})$.\index{Theorem!Krein-Milman's-}
\begin{defn}
A vector space is locally convex if it has a topology which makes
the vector space operators continuous, and if the neighborhoods $\left\{ x+\mbox{Nbh}_{0}\right\} $
have a basis consisting of convex sets.
\end{defn}
The context for Krein-Milman is locally convex topological spaces.
It is in all functional analysis books. Choquet's theorem however
comes later, and it's not found in most books. A good reference is
the book by R. Phelps \cite{Phe01}. The proof of Choquet's theorem
is not specially illuminating. It uses standard integration theory.
\begin{thm}[Krein-Milman]
 Let $K$ be a compact convex set in a locally convex topological
space. Then $K$ is the closed convex hull of its extreme points $E\left(X\right)$,
i.e., 
\[
K=\overline{conv}(E(K)).
\]
\end{thm}
\begin{proof}
(sketch) If $K\supsetneqq\overline{conv}(E(K))$, we get a linear
functional $w$, such that $w$ is zero on $\overline{conv}(E(K))$
and not zero on $w\in K\backslash\overline{conv}(E(K))$. Extend $w$
by Hahn-Banach theorem to a linear functional to the whole space,
and get a contradiction.\end{proof}
\begin{note}
The dual of a normed vector space is always a Banach space, so the
theorem applies. The convex hull in an infinite dimensional space
is not always closed, so close it. A good reference to locally convex
topological space is the lovely book by F. Trèves \cite{MR2296978}.
\end{note}
A convex combination of points $(\xi_{i})$ in $K$ takes the form
$v=\sum c_{i}\xi_{i}$, where $c_{i}>0$ and $\sum c_{i}=1$. Closure
refers to taking limit, so we allow all limits of such convex combinations.
Such a $v$ is obviously in $K$, since $K$ was assumed to be convex.
The point of the Krein-Milman's theorem is the converse. 

The decomposition of states into pure states was developed by Choquet
et al; see \cite{Phe01}. The idea goes back to Krein and Choquet. 

\index{Choquet}\index{decomposition of state}\index{extreme-point}\index{pure state}\index{state!pure-}\index{Theorem!Choquet-}
\begin{thm}[Choquet]
\label{thm:ck} $K=S(\mathfrak{A})$ is a compact convex set in a
locally convex topological space. Let $E(K)$ be the set of extreme
points on $K$. Then for all $p\in K$, there exists a Borel probability
measure $\mu_{p}$, supported on a Borel set $bE(K)\supset E(K)$,
such that for all affine functions $f$, we have
\begin{equation}
f\left(p\right)=\int_{bE(X)}f\left(\xi\right)\:d\mu_{p}(\xi).\label{eq:choquet}
\end{equation}

\end{thm}
The expression in Choquet's theorem is a generalization of convex
combination. In stead of summation, it is an integral against a measure.
Since there are some bizarre cases where the extreme points $E(K)$
do not form a Borel set, the measure $\mu_{p}$ is actually supported
on $bE(K)$, such that $\mu_{p}(bE(K)-E(K))=0$. 

\index{locally convex topological space}

Applications of Choquet theory and of \thmref{ck} are manifold, and
we shall discuss some of them in \chapref{groups} below. Among them
are applications to representations of $C^{*}$-algebras; e.g., the
problem of finding \textquotedbl{}Borel-cross sections\textquotedbl{}
for the set of equivalence classes of representations of a particular
$C^{*}$-algebra. Equivalence here means \textquotedbl{}unitary equivalence.\textquotedbl{}
By a theorem of Glimm \cite{Gli60,Gli61}, we know that there are
infinite simple $C^{*}$-algebras which do not admit such Borel parameterizations.
Examples of this case include the Cuntz algebras $\mathscr{O}_{N}$,
$N>1$. Nonetheless we shall study subclasses of representations of
$\mathscr{O}_{N}$ which correspond to sub-band filters in signal
processing, and to pyramid algorithms for wavelet constructions. In
\chapref{groups} we shall also study representations of the $C^{*}$-algebra
of the free group on 2 generators, as well as the $C^{*}$-algebra
on two generators $u$ and $v$, subject to the relation $uvu^{-1}=u^{2}$.
It is called the Baumslag-Solitar algebra ($BS_{2}$), after Gilbert
Baumslag and Donald Solitar; and it is of great importance in a more
a systematic analysis of families of wavelets. It is the algebra of
a Baumslag-Solitar group. In fact, there are indexed families of Baumslag-Solitar
groups, given by their respective group presentation. They are examples
of two-generator one-relator groups, and they play an important role
in combinatorial group theory, and in geometric group theory as (counter)
examples and test-cases. 

Other examples of uses of Choquet theory in harmonic analysis and
representation theory include such decompositions from classical analysis
as Fourier transform, Laplace transform, as well as direct integral
theory for representations \cite{Sti59,Seg50}.

\index{algebras!Cuntz algebra}\index{representation!- of the Cuntz algebra}

\index{groups!free} \index{signal processing}

\begin{figure}
\includegraphics[scale=0.4]{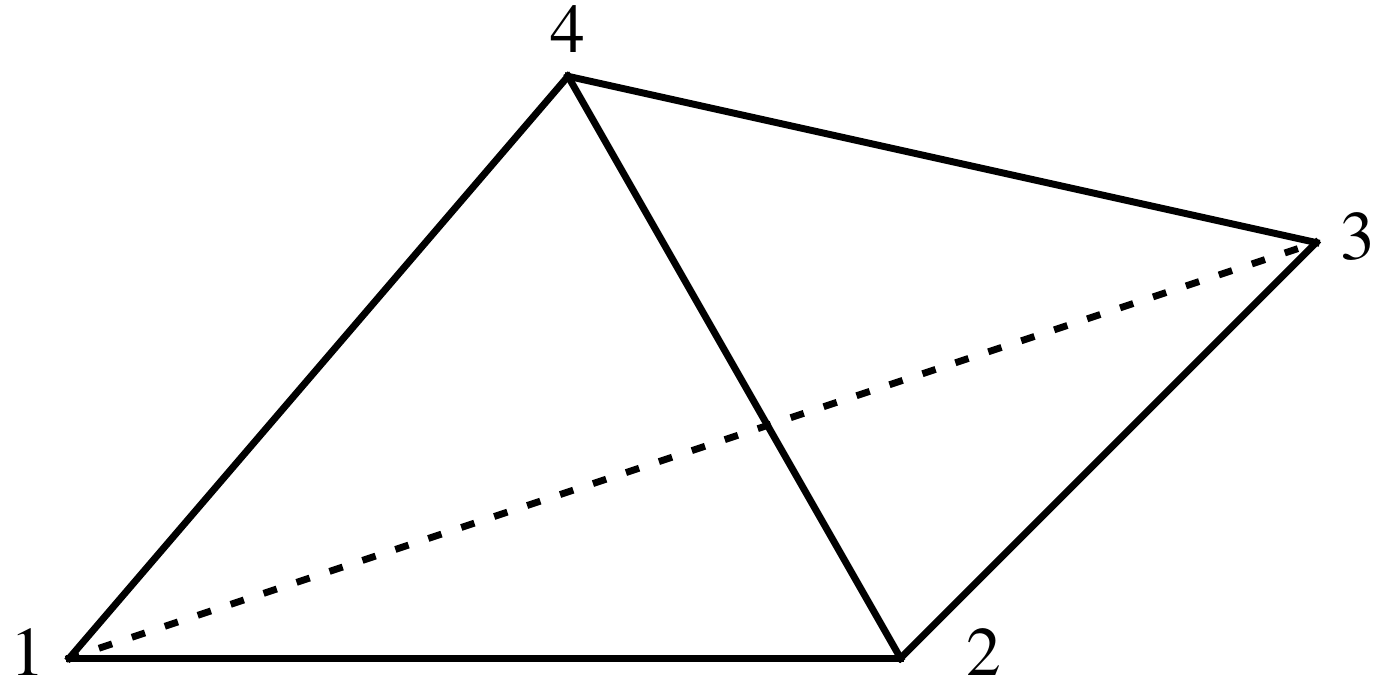}

\protect\caption{\label{fig:sim}A simplex. The four extreme points are marked. }
\end{figure}

\index{measure!probability}
\begin{note}
$\mu_{p}$ in (\ref{eq:choquet}) may not be unique. If it is unique,
$K$ is called a simplex\index{simplex}. The unit disk has its boundary
as extreme points. But representation of points in the interior using
points on the boundary is not unique. Therefore the unit disk is not
a simplex. A tetrahedron is (\figref{sim}).
\end{note}

\begin{example}
Let $(X,\mathfrak{M},\mu)$ be a measure space, where $X$ is compact
and Hausdorff. The set of all probability measures $\mathcal{P}(X)$
is a convex set. To see this, let $\mu_{1},\mu_{2}\in\mathcal{P}(X)$
and $0\leq t\leq1$, then $t\mu_{1}+(1-t)\mu_{2}$ is a measure on
$X$, moreover $(t\mu_{1}+(1-t)\mu_{2})(X)=t+1-t=1$, hence $t\mu_{1}+(1-t)\mu_{2}\in\mathcal{P}(X)$.
Usually we don't want all probability measures, but a closed subset. 
\end{example}

\begin{example}
\label{exa:cnhd}We compute extreme points in the previous example.
$K=\mathcal{P}(X)$ is compact convex in $C(X)^{*}$, which is identified
as the set of all measures due to Riesz. $C(X)^{*}$ is a Banach space
hence is always convex. The importance of being the dual of some Banach
space is that the unit ball is always weak{*}-compact\index{weak{*}-compact}
(Banach-Alaoglu, \thmref{Alaoglu}). Note the weak{*}-topology is
just the cylinder/product topology. The unit ball $B_{1}^{*}$ sits
inside the infinite product space (compact, Hausdorff) $\prod_{v\in B,\left\Vert v\right\Vert =1}D_{1}$,
where $D_{1}=\{z\in\mathbb{C}:\left|z\right|=1\}$. The weak $*$
topology on $B_{1}^{*}$ is just the restriction of the product topology
onto $B_{1}^{*}$. \index{product topology}\index{Theorem!Banach-Alaoglu-}
\end{example}
\index{Banach-Alaoglu Theorem}\index{space!Banach-}\index{space!dual-}
\begin{example}
\emph{Claim}: $E(K)=\{\delta_{x}:x\in X\}$, where $\delta_{x}$ is
the Dirac\index{measure!Dirac} measure supported at $x\in X$. By
Riesz, to know the measure is to know the linear functional. $\int fd\delta_{x}=f(x)$.
Hence we get a family of measures indexed by $X$. If $X=[0,1]$,
we get a continuous family of measures. To see these really are extreme
points, we do the GNS construction on the algebra $\mathfrak{A}=C(X)$,
with the state $\mu\in\mathcal{P}(X)$. The Hilbert space so constructed
is simply $L^{2}(\mu)$. It's clear that $L^{2}(\delta_{x})$ is 1-dimensional,
hence the representation is irreducible. We conclude that $\delta_{x}$
is a pure state, for all $x\in X$.\index{Hilbert space!$L^{2}$}

There is a bijection between state $\varphi$ and Borel measure $\mu:=\mu_{\varphi}$,\index{Borel measure}
\[
\varphi(a)=\int_{X}a\:d\mu_{\varphi}.
\]
In $C\left(X\right)$, $1_{\mathfrak{A}}=\mathbbm{1}=$ constant function.
We check that 
\[
\varphi(1_{\mathfrak{A}})=\varphi(1)=\int1d\mu=\mu(X)=1
\]
since $\mu\in\mathcal{P}\left(X\right)$ is a probability measure.
Also, if $f\geq0$ then $f=g^{2}$, with $g:=\sqrt{f}$; and 
\[
\varphi(f)=\int g^{2}d\mu\geq0.
\]
\end{example}
\begin{note}
$\nu$ is an extreme point in $\mathcal{P}(X)$ if and only if 
\[
\left(\nu\in\left[\mu_{1},\mu_{2}\right]=\mbox{convex hull of }\left\{ \mu_{1},\mu_{2}\right\} \right)\Longrightarrow\left(\nu=\mu_{1}\:\mbox{or}\:\nu=\mu_{2}\right).
\]
\end{note}
\begin{example}
Let $\mathfrak{A}=\mathscr{B}(\mathscr{H})$, and $S\left(\mathfrak{A}\right)=$
states of $\mathfrak{A}$. For each $\xi\in\mathscr{H}$, the map
$A\mapsto w_{\xi}(A):=\left\langle \xi,A\xi\right\rangle $ is a state,
called vector state. 

\emph{Claim}: $E(S)=$ vector states. 

To show this, suppose $W$ is a subspace of $\mathscr{H}$ such that
$0\varsubsetneq W\varsubsetneq\mathscr{H}$, and suppose $W$ is invariant
under the action of $\mathscr{B}(\mathscr{H})$. Then $\exists h\in\mathscr{H}$,
$h\perp W$. Choose $\xi\in W$. The wonderful rank-1 operator (due
to Dirac) $T:\xi\mapsto h$ given by $T:=\left|h\left\rangle \right\langle \xi\right|$,
shows that $h\in W$ (since $TW\subset W$ by assumption.) Hence $h\perp h$
and $h=0$. Therefore $W=\mathscr{H}$. We say $\mathscr{B}(\mathscr{H})$
acts transitively on $\mathscr{H}$.\end{example}
\begin{note}
In general, any $C^{*}$-algebra is a closed subalgebra of $\mathscr{B}(\mathscr{H})$
for some $\mathscr{H}$ (\thmref{gns1}). All the pure states on $\mathscr{B}(\mathscr{H})$
are vector states. \index{pure state}\index{state!pure-}\end{note}
\begin{example}
Let $\mathfrak{A}$ be a $*$-algebra, $S(\mathfrak{A})$ be the set
of states on $\mathfrak{A}$. $w:\mathfrak{A}\rightarrow\mathbb{C}$
is a state on $\mathfrak{A}$ if $w(1_{\mathfrak{A}})=1$ and $w(A)\geq0$,
whenever $A\geq0$. The set of completely positive (CP) maps is a
compact convex set. CP maps are generalizations of states (\chapref{cp}). \end{example}
\begin{xca}[Extreme measures]
\myexercise{Extreme measures}\label{exer:prod}Take the two state
sample space $\Omega=\prod_{1}^{\infty}\{0,1\}$ with product topology.
Assign probability measure, so that we might favor one outcome than
the other. For example, let $s=x_{1}+\cdots x_{n}$, $P_{\theta}(C_{x})=\theta^{s}(1-\theta)^{n-1}$,
i.e. $s$ heads, $(n-s)$ tails. Notice that $P_{\theta}$ is invariant
under permutation of coordinates. $x_{1},x_{2},\ldots,x_{n}\mapsto x_{\sigma(1)}x_{\sigma(2)}\ldots x_{\sigma(n)}$.
$P_{\theta}$ is a member of the set of all such invariant measures
(invariant under permutation) $P_{inv}(\Omega)$. Prove that 
\[
E(P_{inv}(\Omega))=[0,1]
\]
i.e., $P_{\theta}$ are all the possible extreme points.\end{xca}
\begin{rem}
Let $\sigma:X\rightarrow X$ be a measurable transformation. A (probability)
measure $\mu$ is \emph{ergodic}\index{ergodic} if
\[
[E\in\mathfrak{M},\sigma E=E]\Rightarrow\mu(E)\in\{0,1\}.
\]
Intuitively, it says that the whole space $X$ can't be divided non-trivially
into parts where $\mu$ is invariant. The set $X$ will be mixed up
by the transformation $\sigma$.\index{measure!ergodic}\end{rem}
\begin{xca}[Irrational rotation]
\label{exer:irot}\myexercise{Irrational rotation}Let $\theta>0$
be a fixed irrational number, and set \index{irrational rotation}
\begin{equation}
\sigma_{\theta}\left(x\right)=\theta x\;\mbox{mod}\;1\label{eq:irr}
\end{equation}
i.e., multiplication by $\theta$ modulo $1$. Show that $\sigma_{\theta}$
in (\ref{eq:irr}) is ergodic in the measure space $\mathbb{R}/\mathbb{Z}\simeq[0,1)$
with Lebesgue measure (see \figref{irr}).

\begin{figure}[H]
\centering{}\subfloat[$(x,\sqrt{2}x)$ mod $1$, $1\leq x\leq5$]{%
\begin{minipage}[t]{0.45\columnwidth}%
\protect\begin{center}
\protect\includegraphics[scale=0.35]{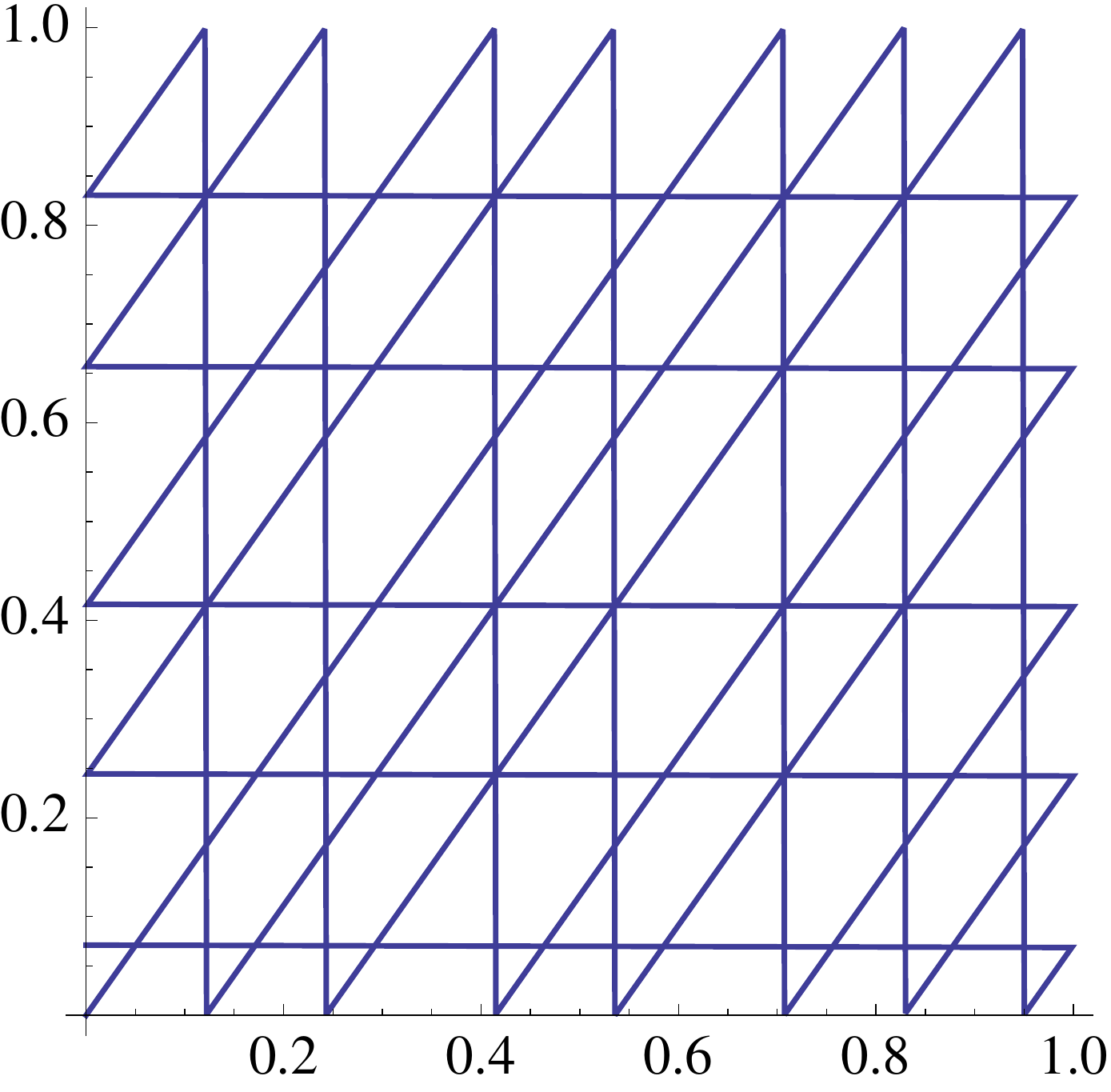}\protect
\par\end{center}%
\end{minipage}

}\hfill{}\subfloat[$(x,\sqrt{2}x)$ mod $1$, $0\leq x\leq15$]{%
\begin{minipage}[t]{0.45\columnwidth}%
\protect\begin{center}
\protect\includegraphics[scale=0.35]{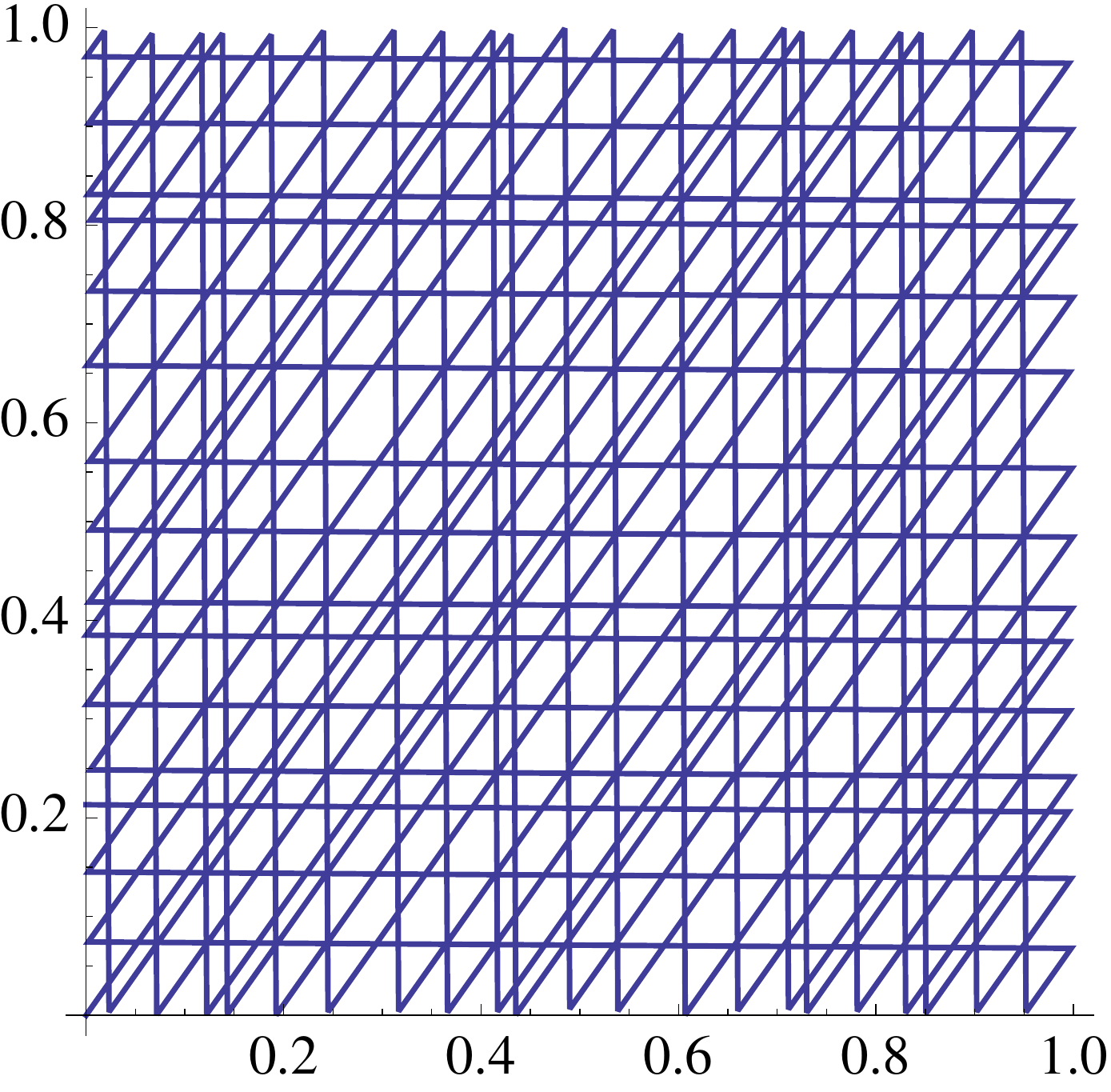}\protect
\par\end{center}%
\end{minipage}}\protect\caption{\label{fig:irr}Irrational rotation. }
\end{figure}

\end{xca}

\subsection{Noncommutative Radon-Nikodym Derivative}

\index{Radon-Nikodym derivative}

Let $w$ be a state on a $C^{*}$-algebra $\mathfrak{A}$, and let
$K$ be an operator in $\mathfrak{A}_{+}$. Set\index{Theorem!Radon-Nikodym-}
\[
w_{K}(A)=\frac{w(\sqrt{K}A\sqrt{K})}{w(K)}.
\]
Then $w_{K}$ is a state, and $w_{K}\ll w$, i.e., $w(A)=0\Rightarrow w_{K}(A)=0$.
We say that $K=\frac{dw}{dw_{K}}$ is a noncommutative Radon-Nikodym
derivative.

Check:
\begin{eqnarray*}
w_{K}(1) & = & 1\\
w_{K}(A^{*}A) & = & \frac{w(\sqrt{K}A^{*}A\sqrt{K})}{w(K)}\\
 & = & \frac{w((A\sqrt{K})^{*}(A\sqrt{K}))}{w(K)}\geq0
\end{eqnarray*}
The converse holds too \cite{MR0442701} and is called the noncommutative
Radon-Nikodym theorem. \index{noncommutative!-Radon-Nikodym derivative}

\subsection{Examples of Disintegration}
\begin{example}
$L^{2}(I)$ with Lebesgue\index{measure!Lebesgue} measure. Let 
\[
F_{x}(t)=\begin{cases}
1 & t\geq x\\
0 & t<x
\end{cases}
\]
$F_{x}$ is a monotone increasing function on $\mathbb{R}$, hence
by Riesz, we get the corresponding Riemann-Stieltjes measure $dF_{x}$.
\index{Riesz' theorem}\index{Theorem!Riesz-}\index{integral!Riesz-}

\[
d\mu=\int^{\oplus}dF_{x}(t)dx.
\]
i.e.
\[
\int fd\mu=\int dF_{x}(f)dx=\int f(x)dx.
\]
Equivalently, 
\[
d\mu=\int\delta_{x}dx
\]
i.e.
\[
\int fd\mu=\int\delta_{x}(f)dx=\int f(x)dx.
\]
$\mu$ is a state, $\delta_{x}=dF_{x}(t)$ is a pure state, $\forall x\in I$.
This is a decomposition of state into direct integral of pure states.
See \cite{Sti59,Seg50}.\index{state!pure-} \index{integral!direct-}
\end{example}

\begin{example}
$\Omega=\prod_{t\geq0}\bar{\mathbb{R}}$, $\Omega_{x}=\{w\in\Omega:w(0)=x\}$.
Kolmogorov gives rise to $P_{x}$ by conditioning $P$ with respect
to ``starting at $x$''.
\[
P=\int^{\oplus}P_{x}dx
\]
i.e.
\[
P()=\int P(\cdot|\text{start at }x)dx.
\]

\end{example}

\begin{example}
Harmonic function on $D$ \index{harmonic} 
\[
h\mapsto h(z)=\int_{\partial\mathbb{D}}\widehat{f}d\mu_{z}
\]
Poisson integration.
\end{example}

\section{\label{sec:egCalgebras}Examples of $C^{*}$-algebras}

Let $\mathscr{H}$ be an infinite-dimensional separable Hilbert space,
and let $S:\mathscr{H}\rightarrow\mathscr{H}$ be an isometry; i.e.,
we have 
\begin{equation}
S^{*}S=I_{\mathscr{H}}.\label{eq:eh1}
\end{equation}
We shall be interested in the case when $S$ is non-unitary, so the
projection 
\[
P_{S}:=SS^{*}
\]
is not $I_{\mathscr{H}}$, i.e., $P_{S}\nleq I_{\mathscr{H}}$. \index{isometry}\index{isometry!Cuntz-}
\begin{thm}[Wold, see \cite{MR0042667,Con90}]
Let $S:\mathscr{H}\rightarrow\mathscr{H}$ be an isometry. Set
\begin{equation}
\mathscr{H}_{0}:=\left\{ x\in\mathscr{H}\::\:\lim_{n\rightarrow\infty}\Vert S^{*^{n}}x\Vert=0\right\} ,\;\mbox{and}\label{eq:eh2}
\end{equation}
\begin{equation}
\mathscr{H}_{1}:=\left\{ x\in\mathscr{H}\::\:\Vert S^{*^{n}}x\Vert=\left\Vert x\right\Vert ,\;\forall n\in\mathbb{N}\right\} .\label{eq:eh3}
\end{equation}

\begin{enumerate}
\item Then 
\begin{equation}
\mathscr{H}=\mathscr{H}_{0}\oplus\mathscr{H}_{1},\label{eq:eh4}
\end{equation}
where ``$\oplus$'' in (\ref{eq:eh4}) refers to orthogonal sum,
i.e., $\mathscr{H}_{0}\perp\mathscr{H}_{1}$.
\item $S\big|_{\mathscr{H}_{0}}:\mathscr{H}_{0}\longrightarrow\mathscr{H}_{0}$
is a shift-operator;
\item $S\big|_{\mathscr{H}_{1}}:\mathscr{H}_{1}\longrightarrow\mathscr{H}_{1}$
is a unitary operator in $\mathscr{H}_{1}$. \index{shift}
\end{enumerate}
\end{thm}
\begin{xca}[Wold\textquoteright s decomposition]
\myexercise{Wold's decomposition}Carry out the details in the proof
of Wold's theorem. \index{Wold decomposition} \index{Theorem!Wold decomposition-}
\end{xca}
In summary, associate to every isometry $\mathscr{H}\xrightarrow{\;S\;}\mathscr{H}$,
there are three subspaces \index{isometry}
\begin{eqnarray*}
\mathscr{H}_{\text{shift}} & = & \mathscr{H}_{0}\;\mbox{in}\;\left(\ref{eq:eh2}\right)\\
\mathscr{H}_{\text{unit}} & = & \mathscr{H}_{1}\;\mbox{in}\;\left(\ref{eq:eh3}\right),\;\mbox{and}\\
\mathfrak{h} & = & \ker\left(S^{*}\right)\;\mbox{in}\;\left(\ref{eq:eh5}\right),\;\mbox{the multiplicity space.}
\end{eqnarray*}
For the closed subspace $\mathscr{H}_{0}$ in the shift-part of the
decomposition, it holds that $\mathscr{H}_{0}$ is the countable direct
sum of $\mathfrak{h}$ with itself.
\begin{xca}[Substitution by $z^{N}$]
\myexercise{Substitution by $z^{N}$}Let $\mathscr{H}=\mathbb{H}_{2}=\mathbb{H}_{2}\left(\mathbb{D}\right)$
be the Hardy space of the disk, and let $N\in\mathbb{N}$, $N>1$.
Set 
\begin{equation}
\left(Sf\right)\left(z\right)=f\left(z^{N}\right),\;f\in\mathbb{H}_{2},\:z\in\mathbb{D}.\label{eq:wd1}
\end{equation}
Show that the three closed subspaces for this isometry are as follows:
\begin{eqnarray*}
\mathscr{H}_{\text{unit}} & = & \mbox{the constant functions on}\;\mathbb{D}\\
 & = & \mathbb{C}e_{0},\;e_{0}\left(z\right)=z^{0}=1.\\
\mathscr{H}_{\text{shift}} & = & \mathscr{H}\ominus\mathbb{C}e_{0}\\
 & = & \left\{ f\in\mathbb{H}_{2}\::\:f\left(0\right)=0\right\} \\
\ker\left(S^{*}\right) & = & \overline{span}\left\{ z^{k}\::\:N\nmid k\;\left(\mbox{not divisible by}\;N\right)\right\} ,\;\mbox{i.e.},\\
 &  & \mbox{powers of}\;z^{k},\;k\in\left(\left\{ 0\right\} \cup\mathbb{N}\right)\backslash N\mathbb{Z},\;\mbox{so}\;k\;\mbox{not}\\
 &  & \mbox{divisible by}\;N.
\end{eqnarray*}
The isometry $S$ in (\ref{eq:wd1}) is an example of an \emph{isometry
of infinite multiplicity}. \index{isometry} \index{multiplicity!-of an isometry}
\end{xca}

\paragraph{$C^{*}$-algebras generated by isometries.}

An important family of non-abelian $C^{*}$-algebras includes those
generated by one, or more, isometries:

\paragraph{Case 1. One Isometry}

Because of Wold's decomposition, if a $C^{*}$-algebra $\mathfrak{A}$
is generated by one isometry, we may ``split off'' the one generated
by the unitary part; and then reduce the study to the case where $\mathfrak{A}$
is generated by a shift $S$.

Introduce $\mathfrak{h}:=\ker\left(S^{*}\right)$, and \index{shift}\index{isometry}
\begin{equation}
\mathfrak{h^{\infty}=\oplus_{\mathbb{N}}\mathfrak{h}}=\left\{ \left(x,x_{2},\cdots\right)\::\:x_{i}\in\mathfrak{h}\right\} \label{eq:eh5}
\end{equation}
\begin{equation}
\left\Vert \left(x_{1},x_{2},\cdots\right)\right\Vert _{\mathfrak{h}^{\infty}}^{2}:=\sum_{i=1}^{\infty}\left\Vert x_{i}\right\Vert _{\mathfrak{h}}^{2};\label{eq:eh6}
\end{equation}
and set
\begin{equation}
S_{\infty}\left(x_{1},x_{2},x_{3},\cdots\right)=\left(0,x_{1},x_{2},x_{3},\cdots\right).\label{eq:eh7}
\end{equation}

\begin{xca}[The backwards shift]
\label{exer:bshift}\myexercise{The backwards shift}Show that 
\[
S_{\infty}^{*}\left(x_{1},x_{2},x_{3},\cdots\right)=\left(x_{2},x_{3},x_{4},\cdots\right),\;\mbox{and that}
\]
\[
\left\Vert \left(S_{\infty}^{*}\right)^{n}x\right\Vert \xrightarrow[n\rightarrow\infty]{}0.
\]
\end{xca}
\begin{rem}
It would seem like the backwards shift is an overly specialized example.
Nonetheless it plays a big role in operator theory, see for example
\cite{MR1997376,MR3214905,MR2496397}, and it is an example of a wider
class of operators going by the name ``the Cowen-Douglass class,\textquotedblright{}
playing an important role in complex geometry, see \cite{MR501368}.\end{rem}
\begin{xca}[Infinite multiplicity]
\label{exer:infm}\myexercise{Infinite multiplicity}For the isometry
$\left(Sf\right)\left(z\right)=f\left(z^{N}\right)$, $f\in\mathbb{H}_{2}$,
$z\in\mathbb{D}$, write out the representation (\ref{eq:eh5})-(\ref{eq:eh7})
above. \index{isometry!Cuntz-}\index{multiplicity!-of an isometry}\index{representation!- of the Cuntz algebra}
\end{xca}

\begin{xca}[A shift is really a shift]
\myexercise{A shift is really a shift}Show that $S$ and $S_{\infty}$
are unitarily equivalent if and only if $S$ is a shift.
\end{xca}

\begin{xca}[Multiplication by $z$ is a shift in $\mathbb{H}_{2}$]
\myexercise{Multiplication by $z$ is a shift in $\mathbb{H}_{2}$}If
$\dim\mathfrak{h}=1$ (multiplicity one), show that $S$ in (\ref{eq:eh7})
is unitarily equivalent to 
\begin{equation}
\left(\widetilde{S}f\right)\left(z\right)=zf\left(z\right),\;z\in\mathbb{D},\:f\in\mathbb{H}_{2}=\mbox{the Hardy space.}\label{eq:eh8}
\end{equation}

\uline{Hint}: By $\mathbb{H}_{2}$, we mean the Hilbert space of
all analytic functions $f$ on the disk $\mathbb{D}=\left\{ z\in\mathbb{C}\::\:\left|z\right|<1\right\} $
such that 
\begin{equation}
f\left(z\right)=\sum_{k=1}^{\infty}a_{k}z^{k},\;\mbox{and}\;\left(a_{k}\right)\in l^{2}.\label{eq:eh9}
\end{equation}
We set $\left\Vert f\right\Vert _{\mathbb{H}_{2}}=\left\Vert \left(a_{k}\right)\right\Vert _{l^{2}}$.\index{shift}
\end{xca}

\begin{xca}[The two shifts in $\mathbb{H}_{2}$]
\myexercise{The two shifts in $\mathbb{H}_{2}$}Show that the adjoint
to the generator $\widetilde{S}$ (from (\ref{eq:eh8})) is
\begin{equation}
\left(\widetilde{S}^{*}f\right)\left(z\right)=\frac{f\left(z\right)-f\left(0\right)}{z},\;\forall f\in\mathbb{H}_{2},\forall z\in\mathbb{D}\backslash\left\{ 0\right\} ,\label{eq:ef1}
\end{equation}
and 
\[
\left(\widetilde{S}^{*}f\right)\left(0\right)=f'\left(0\right),\;f\in\mathbb{H}_{2}.
\]

\uline{Hint}: Show that, if $f,g\in\mathbb{H}_{2}$, then the following
holds:
\begin{equation}
\left\langle \widetilde{S}f,g\right\rangle _{\mathbb{H}_{2}}=\left\langle f,\widetilde{S}^{*}g\right\rangle _{\mathbb{H}_{2}}\label{eq:ef2}
\end{equation}
where we use formula (\ref{eq:ef1}) in computing the $\mathbb{H}_{2}$-inner
product on the RHS in (\ref{eq:ef2}).

Compare this with the result from \exerref{bshift}.
\end{xca}

\begin{xca}[A numerical range]
\myexercise{A numerical range} Let $T:=S_{\infty}^{*}$ be the backward
shift (expressed in coordinates) in \exerref{bshift}. Since $TT^{*}-T^{*}T$
is the rank-one projection $\left|e_{1}\left\rangle \right\langle e_{1}\right|$,
of course $T$ is not normal.
\begin{enumerate}
\item \label{enu:bs1-1}Show that $x_{\lambda}=\left(1,\lambda,\lambda^{2},\lambda^{3},\cdots\right)$
satisfies
\begin{equation}
Tx_{\lambda}=\lambda x_{\lambda},\;\forall\lambda\in\mathbb{C}.\label{eq:bs1-1}
\end{equation}

\item \label{enu:bs1-2}Since 
\begin{equation}
x_{\lambda}\in l^{2}\Longleftrightarrow\left|\lambda\right|<1,\label{eq:bs1-2}
\end{equation}
conclude that the point-spectrum of $T$ is $\mathbb{D}=\left\{ \lambda\in\mathbb{C}\::\:\left|\lambda\right|<1\right\} $. 
\item Combine (\ref{enu:bs1-1}) \& (\ref{enu:bs1-2}) in order to conclude
that 
\[
NR_{T}=\mathbb{D}.
\]

\end{enumerate}
\end{xca}

\begin{xca}[The finite shift]
\myexercise{The finite shift} Compare the infinite case above with
the analogous matrix case
\[
T_{3}=\begin{bmatrix}0 & 1 & 0\\
0 & 0 & 1\\
0 & 0 & 0
\end{bmatrix}\quad\mbox{and}\quad T_{3}'=\begin{bmatrix}0 & 1 & 0\\
0 & 0 & 1\\
1 & 0 & 0
\end{bmatrix}
\]
A sketch of $NR_{T_{3}}$ and $NR_{T_{3}'}$ are in \figref{T3} below.
See also \figref{bs}. 
\end{xca}
\begin{figure}
\subfloat[$NR_{T_{3}}$]{%
\begin{minipage}[t]{0.45\columnwidth}%
\noindent \protect\begin{center}
\protect\includegraphics[width=150pt,height=150pt,keepaspectratio]{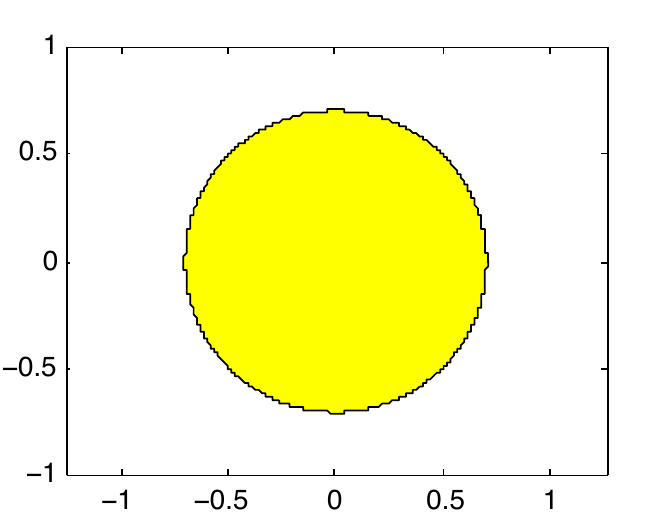}\protect
\par\end{center}%
\end{minipage}

}\hfill{}\subfloat[$NR_{T_{3}'}$]{%
\begin{minipage}[t]{0.45\columnwidth}%
\noindent \protect\begin{center}
\protect\includegraphics[width=150pt,height=150pt,keepaspectratio]{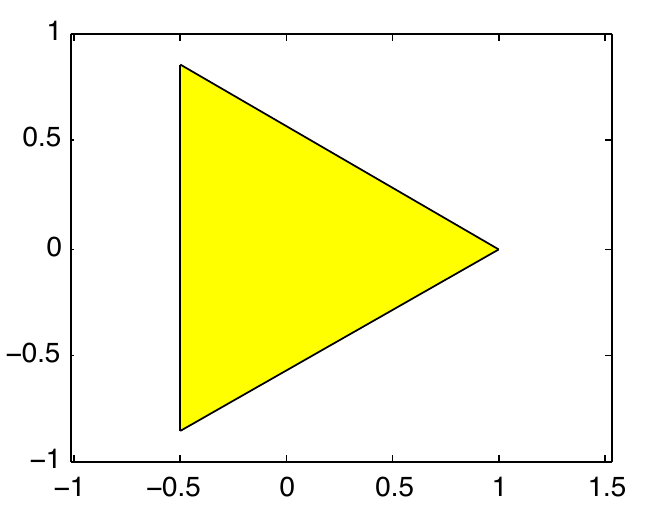}\protect
\par\end{center}%
\end{minipage}

}

\protect\caption{\label{fig:T3}The numerical range of $T_{3}$ vs $T_{3}'$}
\end{figure}

\begin{figure}
\includegraphics[width=200pt,height=200pt,keepaspectratio]{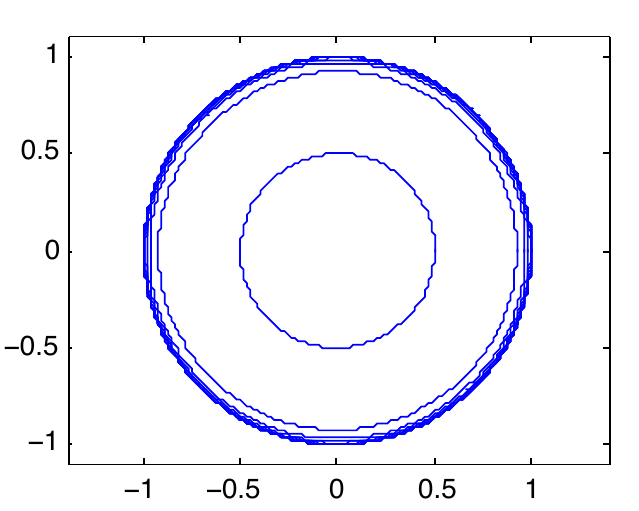}

\protect\caption{\label{fig:bs}The numerical range (NR) of the truncated finite matrices:
Expanding truncations of the infinite matrix $T$ corresponding to
the backward shift, and letting the size $\protect\longrightarrow\infty$:
$T_{3},T_{4},\cdots,T_{n},T_{n+1},\cdots$; the limit-NR fills the
open disk of radius $1$.}
\end{figure}

\begin{xca}[The Hardy space $\mathbb{H}_{2}$; a transform]
\myexercise{The Hardy space $\mathbb{H}_{2}$; a transform}\label{exer:hardy}Let
$f\left(z\right)=\sum_{k=0}^{\infty}a_{k}z^{k}$, and set 
\begin{equation}
\widetilde{f}\left(t\right)=\sum_{k=0}^{\infty}a_{k}e^{i2\pi kt},\;t\in\mathbb{R}.\label{eq:eh10}
\end{equation}
Show that 
\[
f\in\mathbb{H}_{2}\Longleftrightarrow\widetilde{f}\in L^{2}\left(\mathbb{T}\right),\;\mathbb{T}=\partial\mathbb{D};
\]
and that 
\begin{equation}
\left\Vert \widetilde{f}\right\Vert _{L^{2}\left(\mathbb{T}\right)}=\left\Vert f\right\Vert _{\mathbb{H}_{2}}\label{eq:eh11}
\end{equation}
holds.
\end{xca}
Because of \exerref{hardy}, we may identify $\mathbb{H}_{2}$ with
a closed subspace in $L^{2}\left(\mathbb{T}\right)$. Let $P_{+}$
denote the projection of $L^{2}\left(\mathbb{T}\right)$ onto $\mathbb{H}_{2}$.\index{space!Hardy-}
\begin{defn}
For $\varphi\in L^{\infty}\left(\mathbb{T}\right)$, set 
\begin{equation}
T_{\varphi}f=P_{+}\left(\varphi f\right),\;\forall f\in\mathbb{H}_{2},\label{eq:eh12}
\end{equation}
equivalently, $T_{\varphi}=P_{+}M_{\varphi}P_{+}$. 

The operator $T_{\varphi}$ in (\ref{eq:eh12}) is called a \emph{Toeplitz-operator};
and 
\begin{equation}
\mathscr{T}:=C^{*}\left(\left\{ T_{\varphi}\::\:\varphi\in L^{\infty}\left(\mathbb{T}\right)\right\} \right)\label{eq:eh12-1}
\end{equation}
is called the \emph{Toeplitz-algebra}. 
\end{defn}
\index{algebra!Toeplitz-}

\index{Toeplitz!-operator}

\index{Toeplitz!-algebra}
\begin{xca}[Multiplicity-one and $\mathbb{H}_{2}$]
\myexercise{Multiplicity-one and $\mathbb{H}_{2}$}Show that there
is a short exact sequence (in the category of $C^{*}$-algebras):
\begin{equation}
0\longrightarrow\mathscr{K}\longrightarrow\mathscr{T}\xrightarrow{\;\pi\;}\mathscr{T}/\mathscr{K}\longrightarrow0\label{eq:eh13}
\end{equation}
where $\mathscr{K}=$ the $C^{*}$-algebra of compact operators; and
$\mathscr{T}/\mathscr{K}$ is the quotient; finally 
\[
\mathscr{T}/\mathscr{K}\simeq L^{\infty}\left(\mathbb{T}\right),
\]
realized via the mapping $T_{\varphi}\xrightarrow{\;\pi\;}\varphi$
(the symbol mapping) in (\ref{eq:eh13}), i.e., $\pi\left(T_{\varphi}\right):=\varphi$,
is assigning the symbol $\varphi$ to the Toeplitz operator $T_{\varphi}$. 
\end{xca}
\index{algebra!quotient-}

\index{short exact sequence}

\begin{xca}[The Toeplitz matrices]
\myexercise{The Toeplitz matrices}Suppose $\varphi\in L^{\infty}\left(\mathbb{T}\right)$
has Fourier expansion 
\begin{equation}
\varphi\left(t\right)=\sum_{n\in\mathbb{Z}}b_{n}e^{i2\pi nt},\;t\in\mathbb{R}.\label{eq:ee1}
\end{equation}
Then show that the $\infty\times\infty$ matrix of the corresponding
Toeplitz operator $T_{\varphi}$ is as follows w.r.t the standard
ONB in $\mathbb{H}_{2}$, $\left\{ z^{n}\::\:n\in\left\{ 0\right\} \cup\mathbb{N}\right\} $.
\begin{equation}
Mat\left(T_{\varphi}\right)=\begin{bmatrix}b_{0} & b_{-1} & b_{-2} & b_{-3} & b_{-4} & \cdots & \cdots & \cdots & \cdots & \cdots\\
b_{1} & b_{0} & b_{-1} & b_{-2} & b_{-3} & \cdots & \cdots & \cdots & \cdots & \cdots\\
b_{2} & b_{1} & b_{0} & b_{-1} & b_{-2} & \cdots & \cdots & \cdots & \cdots & \cdots\\
b_{3} & b_{2} & b_{1} & b_{0} & b_{-1} & \cdots & \cdots & \cdots & \cdots & \cdots\\
b_{4} & b_{3} & b_{2} & b_{1} & b_{0} & \cdots & \cdots & \cdots & \cdots & \cdots\\
\ddots & \ddots & \ddots & \ddots & \ddots & \ddots & \ddots & \ddots & \cdots & \cdots\\
\ddots & \ddots & \ddots & \ddots & \ddots & \ddots & b_{0} & b_{-1} & b_{-2} & \cdots\\
\ddots & \ddots & \ddots & \ddots & \ddots & \ddots & b_{1} & b_{0} & b_{-1} & \cdots\\
\ddots & \ddots & \ddots & \ddots & \ddots & \ddots & b_{2} & b_{1} & b_{0} & \ddots\\
\ddots & \ddots & \ddots & \ddots & \ddots & \ddots & \ddots & \ddots & \ddots & \ddots
\end{bmatrix}\label{eq:ee2}
\end{equation}
\end{xca}
\begin{note}
Matrices of the form given in (\ref{eq:ee2}) are called \emph{Toeplitz
matrices}, i.e., with the banded pattern, constant numbers down the
diagonal lines, with $b_{0}$ in the main diagonal. 
\end{note}
\index{Toeplitz!-operator}\index{Toeplitz!-matrix}
\begin{rem}
\label{rem:ta}Note that the mapping $\varphi\longrightarrow T_{\varphi}$
(Toeplitz), $L^{\infty}\left(\mathbb{T}\right)\rightarrow\mathscr{T}$
is \uline{not} a homomorphism of the algebra $L^{\infty}\left(\mathbb{T}\right)$
into $\mathscr{T}$ = $C^{*}\left(\left\{ T_{\varphi}\right\} \right)$.
Here we view $L^{\infty}\left(\mathbb{T}\right)$ as an abelian $C^{*}$-algebra
under pointwise product, i.e., 
\[
\left(\varphi_{1}\varphi_{2}\right)\left(t\right):=\varphi_{1}\left(t\right)\varphi_{2}\left(t\right),\;\forall t\in\mathbb{R}/\mathbb{Z}.
\]

The point of mapping of the \emph{short exact sequence} (lingo from
homological algebra):
\begin{equation}
0\longrightarrow\mathscr{K}\longrightarrow\mathscr{T}\longrightarrow L^{\infty}\left(\mathbb{T}\right)\longrightarrow0\label{eq:ee2-1}
\end{equation}
is that $\varphi\longrightarrow T_{\varphi}$ is only a ``homomorphism
mod $\mathscr{K}$ (= the compact operators)'', i.e., that we have
\begin{equation}
T_{\varphi_{1}}T_{\varphi_{2}}-T_{\varphi_{1}\varphi_{2}}\in\mathscr{K}\label{eq:ee3}
\end{equation}
valid for all $\varphi_{1},\varphi_{2}\in L^{\infty}\left(\mathbb{T}\right)$. 

There is an extensive literature on (\ref{eq:ee2-1}) and (\ref{eq:ee3}),
see especially \cite{MR571362}.\end{rem}
\begin{xca}[homomorphism mod $\mathscr{K}$]
\myexercise{homomorphism mod $\mathscr{K}$} Give a direct proof
that the operator on the LHS in (\ref{eq:ee3}) is a compact operator
in $\mathbb{H}_{2}$.\end{xca}
\begin{rem}
The subject of Toeplitz operators, and Toeplitz algebras is vast (see
e.g., \cite{MR2291432}). The more restricted case where the symbol
$\varphi$ of $T_{\varphi}=P_{+}M_{\varphi}P_{+}$ is continuous (i.e.,
$\varphi\in C(S^{1})$, $S^{1}$ = the circle) is especially rich;
starting with \emph{Szegö's Index Theorem} : \index{Szegö's Index Theorem}\index{index!Szegö-}\end{rem}
\begin{defn}
Let $X$ and $Y$ be Banach spaces. A \emph{Fredholm operator} is
a bounded linear operator $T:X\rightarrow Y$, such that $\ker\left(T\right)$
and $\ker\left(T^{*}\right)$ are finite-dimensional, and $\mbox{ran}\left(T\right)$
is closed. The \emph{index} of $T$ is given as 
\[
ind\left(T\right):=\dim\left(\ker\left(T\right)\right)-\dim\left(\ker\left(T^{*}\right)\right).
\]
(The assumption on the range of $T$ in the definition is redundant
\cite{MR1921782}.)\end{defn}
\begin{thm}[Szegö \cite{MR1300213}]
 If $\varphi\in C(S^{1})$ and $\varphi$ does not vanish on $S^{1}$,
then $T_{\varphi}$ is Fredholm, and the index of $T_{\varphi}$ computes
as follows: \index{index!Fredholm-}
\begin{equation}
ind(T_{\varphi})=\dim(\ker(T_{\varphi}))-\dim(\ker(T_{\varphi}^{*}))=-\#w\left(\varphi\right)\label{eq:sz1}
\end{equation}
where $\#w\left(\varphi\right)$ in (\ref{eq:sz1}) is the winding
number 
\begin{equation}
\#w\left(\varphi\right)=\frac{1}{2\pi i}\int_{0}^{2\pi}\frac{\varphi'\left(e^{i\theta}\right)}{\varphi\left(e^{i\theta}\right)}d\theta.\label{eq:sz2}
\end{equation}
Note: $\#w(e^{in\theta})=n$, for $n\in\mathbb{Z}$, and $\ker(T_{\varphi}^{*})=(\mbox{ran}(T_{\varphi}))^{\perp}$. 
\end{thm}

\paragraph{Case 2. Multiple Isometries}

Here we refer to the Cuntz-algebra $\mathscr{O}_{N}$ (see \cite{MR0467330}),
the unique $C^{*}$-algebra $\mathscr{O}_{N}$, $N>1$, generated
by $\left\{ S_{i}\right\} _{i=1}^{N}$ and the relations 
\begin{equation}
S_{i}^{*}S_{j}=\delta_{ij},\;\mbox{and}\label{eq:eh14}
\end{equation}
\begin{equation}
\sum_{i=1}^{N}S_{i}S_{i}^{*}=\mathbf{1}.\label{eq:eh15}
\end{equation}

Cuntz showed (\cite{MR0467330}) that this is a simple $C^{*}$-algebra
(i..e, no non-trivial closed two-sided ideals), purely infinite. \index{Cuntz-algebra}\index{representation!- of the Cuntz algebra}

We shall return to the study of its \emph{representation} in \chapref{groups}.
\index{ideal}
\begin{xca}[An element in $Rep\left(\mathscr{O}_{N},\mathbb{H}_{2}\right)$]
\myexercise{An element in $Rep\left(\mathscr{O}_{N},\mathbb{H}_{2}\right)$}
Fix $N\in\mathbb{N}$, $N>1$, and consider the following operators
$\left\{ S_{k}\right\} _{k=0}^{N-1}$ acting in the Hardy space $\mathbb{H}_{2}=\mathbb{H}_{2}\left(\mathbb{D}\right)$:
\begin{equation}
\left(S_{k}f\right)\left(z\right)=z^{k}f\left(z^{N}\right),\;\forall f\in\mathbb{H}_{2},\:\forall z\in\mathbb{D},\:k=0,1,\ldots,N-1.\label{eq:mi1}
\end{equation}
Show that the operators $\left(S_{k}\right)$ in (\ref{eq:mi1}) satisfy
the $\mathscr{O}_{N}$-relations (\ref{eq:eh14})-(\ref{eq:eh15}),
i.e., that
\begin{gather*}
S_{j}^{*}S_{k}=\delta_{jk}I_{\mathbb{H}_{2}},\;\mbox{and}\\
\sum_{j=0}^{N-1}S_{j}S_{j}^{*}=I_{\mathbb{H}_{2}};
\end{gather*}
hence a \emph{representation} of $\mathscr{O}_{N}$ in $\mathbb{H}_{2}$.
\end{xca}

\begin{xca}[The multivariable Toeplitz algebra]
\myexercise{The multivariable Toeplitz algebra}For $k\in\mathbb{N}$,
set $\mathscr{H}_{k}=\mathbb{C}^{k}$ = the $k$-dimensional complex
Hilbert space with the usual inner product: \index{Toeplitz!multivariable-algebra}
\begin{equation}
\left\langle v,w\right\rangle =\sum_{j=1}^{k}\overline{v_{j}}w_{j}.\label{eq:mt1}
\end{equation}
For $k=1$, pick a normalized basis vector $\Omega$. For $N>1$,
set 
\begin{equation}
\mathscr{F}\left(\mathscr{H}_{N}\right)=\mathscr{H}_{1}\oplus\sum_{n=1}^{\infty}\oplus\mathscr{H}_{N}^{\otimes n.}\label{eq:mt2}
\end{equation}
(The letter $\mathscr{F}$ is for Fock-space.) For $f\in\mathscr{H}_{N}$,
set:
\begin{eqnarray}
T_{f}\left(\otimes_{1}^{n}h_{j}\right) & = & f\otimes\left(\otimes_{1}^{n}h_{j}\right),\;\mbox{and}\label{eq:mt3}\\
T_{f}^{*}\left(\otimes_{1}^{n}h_{j}\right) & = & \left\langle f,h_{1}\right\rangle \otimes_{2}^{n}h_{j},\;n\in\mathbb{N}.\label{eq:mt4}
\end{eqnarray}
And finally, the vacuum rule: \index{space!Fock-} 
\begin{equation}
T_{f}^{*}\Omega=0.\label{eq:mt5}
\end{equation}

\begin{enumerate}
\item Show that the following hold:
\begin{equation}
T_{f}^{*}T_{g}=\left\langle f,g\right\rangle _{N}I_{\mathscr{F}\left(\mathscr{H}_{N}\right)},\;\forall f,g\in\mathscr{H}_{N}.\label{eq:mt6}
\end{equation}
Define $T_{i}$ and $T_{i}^{*}$ from and ONB in $\mathscr{H}_{N}$,
we get 
\begin{equation}
\sum_{i=1}^{N}T_{i}T_{i}^{*}=I_{\mathscr{F}\left(\mathscr{H}_{N}\right)}-\left|\Omega\left\rangle \right\langle \Omega\right|.\label{eq:mt7}
\end{equation}
The $C^{*}$-algebra generated by $\left\{ T_{f}\::\:f\in\mathscr{H}_{N}\right\} $
is called the (multivariable) \emph{Toeplitz algebra}, and is denoted
$\mathscr{T}_{N}$. 
\item Show, with the use of (\ref{eq:mt6})-(\ref{eq:mt7}), that there
is a natural short exact sequence of $C^{*}$-algebras: 
\[
0\longrightarrow\mathscr{K}\longrightarrow\mathscr{T}_{N}\longrightarrow\mathscr{O}_{N}\longrightarrow0.
\]
Compare with (\ref{eq:ee2-1}) in \remref{ta}.
\end{enumerate}
\end{xca}

\section{Examples of Representations}

We consider the Fourier algebra. 
\begin{enumerate}
\item Discrete case: $l^{1}\left(\mathbb{Z}\right)$ and the Gelfand transform
\begin{eqnarray*}
(a*b)_{n} & = & \sum_{k}a_{k}b_{n-k}\\
(a^{*})_{n} & = & \overline{a_{-n}}\\
1_{\mathfrak{A}} & = & \delta_{0}
\end{eqnarray*}
\[
a\xrightarrow[\text{Gelfand}]{\mathcal{F}}F\left(z\right):=\sum_{n}a_{n}z^{n}
\]
We may specialize to $z=e^{it}$, $t\in\mathbb{R}\text{ mod }2\pi$.
$\{F(z)\}$ is an abelian algebra of functions, with multiplication
is given by 
\[
F(z)G(z)=\sum_{n}(a*b)_{n}z^{n}
\]
In fact, most abelian algebras can be thought of as function algebras.
\\
\\
Homomorphism:
\begin{eqnarray*}
(l^{1},*) & \xrightarrow{\mathcal{F}} & C(\mathbb{T}^{1})\\
(a_{n}) & \mapsto & F(z).
\end{eqnarray*}
If we want to write $F(z)$ as power series, then we need to drop
$a_{n}$ for $n<0$. Then $F(z)$ extends to an analytic function
over the unit disk. The representation by the sequence space
\[
\{a_{0},a_{1},\ldots\}
\]
was suggested by Hardy. We set
\[
\left\Vert F\right\Vert _{\mathbb{H}_{2}}^{2}=\sum_{k=0}^{\infty}\left|a_{k}\right|^{2};
\]
the natural isometric isomorphism. Rudin has two nice chapters on
$H^{2}$, as a Hilbert space, a RKHS. See \cite[ch16]{Rud87}.\index{reproducing kernel Hilbert space (RKHS)}\index{space!Hilbert-}\index{kernel!reproducing-}\index{space!Hardy-}
\item Continuous case: $L^{1}\left(\mathbb{R}\right)$
\begin{eqnarray*}
(f*g)(x) & = & \int_{-\infty}^{\infty}f(s)g(x-s)ds\\
f^{*}\left(x\right) & = & \overline{f\left(-x\right)}
\end{eqnarray*}
The algebra $L^{1}$ has no identity, but we may always insert one
by adding $\delta_{0}$. So $\delta_{0}$ is the homomorphism $f\longmapsto f\left(0\right)$;
and $L^{1}\left(\mathbb{R}\right)\cup\left\{ \delta_{0}\right\} $
is again a Banach $*$-algebra.\\
\\
The \emph{Gelfand map} is the classical Fourier transform, i.e., 
\[
f\xrightarrow[\text{Gelfand}]{\mathcal{F}}\hat{f}\left(\xi\right)=\int_{-\infty}^{\infty}f\left(x\right)e^{-i\xi x}dx
\]
where $\widehat{f*g}=\hat{f}\hat{g}$.\end{enumerate}
\begin{rem}
$C(\mathbb{T}^{1})$ is called the $C^{*}$-algebra completion of
$l^{1}$. $L^{\infty}(X,B,\mu)=L^{1}(\mu)^{*}$ is also a $C^{*}$-algebra.
It is a $W^{*}$-algebra, or von Neumann algebra (see, e.g., \cite{MR0442701}).
The $W^{*}$ refers to the fact that its topology comes from the weak
$*$-topology. Recall that $\mathscr{B}(\mathscr{H})$, for any Hilbert
space, is a von Neumann algebra. \index{completion!$C^{*}$-}\index{algebras!von Neumann algebra ($W^{*}$-algebra)}\end{rem}
\begin{example}
Fix $\varphi$ and set $uf=e^{i\theta}f(\theta)$, $vf=f(\theta-\varphi)$,
restrict to $[0,2\pi]$, i.e., $2\pi$ periodic functions.
\begin{eqnarray*}
vuv^{-1} & = & e^{i\varphi}u\\
vu & = & e^{i\varphi}uv
\end{eqnarray*}
$u,v$ generate a noncommutative $C^{*}$-algebra. See \cite{MR2908740,MR2681883}.
\end{example}

\begin{example}[Quantum Mechanics]
Consider the canonical commutation relation
\[
[p,q]=-i\,I,\quad i=\sqrt{-1},
\]
where $\left[x,y\right]:=xy-yx$ denotes the commutator of $x$ and
$y$ . \index{commutator}

The two symbols $p,q$ generate an algebra, but they \uline{can
not} be represented by bounded operators. But we may apply bounded
functions to them and get a $C^{*}$-algebra.
\end{example}
\index{Heisenberg, W.K.!commutation relation}\index{commutation relations}
\begin{xca}[No bounded solutions to the canonical commutation relations]
\myexercise{No bounded solutions to the canonical commutation relations}\label{exer:uncertainty1}Show
that $p,q$ and not be represented by bounded operators. \uline{Hint}:
take the trace.
\end{xca}

\begin{example}
Let $\mathscr{H}$ be an infinite dimensional Hilbert space, then
$\mathscr{H}$ is isometrically isomorphic to a proper subspace of
itself. For example, let $\{e_{n}\}$ be an ONB. $\mathscr{H}_{1}=\overline{span}\{e_{2n}\}$,
$\mathscr{H}_{2}=\overline{span}\{e_{2n+1}\}$. Let 
\begin{eqnarray*}
V_{1}(e_{n}) & = & e_{2n}\\
V_{2}(e_{n}) & = & e_{2n+1}
\end{eqnarray*}
then we get two isometries. Also, 
\begin{eqnarray*}
V_{1}V_{1}^{*}+V_{2}V_{2}^{*} & = & I\\
V_{i}^{*}V_{i} & = & I\\
V_{i}V_{i}^{*} & = & P_{i}
\end{eqnarray*}
where $P_{i}$ is a selfadjoint \emph{projection}, $i=1,2$ onto the
respective $\mathscr{H}_{i}$. This is the Cuntz algebra\index{algebras!Cuntz algebra}
$\mathcal{O}_{2}$. More general $\mathcal{O}_{N}$, $N>2$. \index{representation!- of the Cuntz algebra}

Cuntz (in 1977) showed that this is a \emph{simple} $C^{*}$-algebra,
i.e., it does not have non-trivial closed two-sided ideals. For studies
of its representations, see, e.g., \cite{Gli60,Gli61,BJO04}. \index{ideal}
\end{example}

\section{Beginning of Multiplicity Theory\label{sec:mult}}

The main question here is how to break up a representation into smaller
ones. The smallest are the irreducible representations, and the next
would be the multiplicity free representations. \index{multiplicity}

Let $\mathfrak{A}$ be an algebra. \index{multiplicity free}
\begin{itemize}
\item \emph{commutative}: e.g., function algebras
\item \emph{non-commutative}: e.g., matrix algebra, algebras generated by
representation of non-abelian groups
\end{itemize}
Smallest representation:
\begin{itemize}
\item \emph{irreducible}: $\pi\in Rep_{irr}\left(\mathfrak{A},\mathscr{H}\right)$,
where the commutant $\pi\left(\mathfrak{A}\right)'$ is 1-dimensional.
This is the starting point of further analysis.
\item \emph{multiplicity free}: Let $\pi\in Rep\left(\mathfrak{A},\mathscr{H}\right)$.
We may assume $\pi$ is cyclic, since otherwise $\pi$ can be decomposed
into a direct sum of cyclic representations, i.e., $\pi=\oplus\pi_{cyc}$;
see \thmref{cyclic2}. Then, 
\[
\pi\;\mbox{is multiplicity free}\Longleftrightarrow\pi\left(\mathfrak{A}\right)'\;\mbox{is abelian}.
\]
\index{representation!multiplicity free}
\end{itemize}
Fix a Hilbert space $\mathscr{H}$, and let $\mathfrak{C}$ be a $*$-algebra
in $\mathscr{B}\left(\mathscr{H}\right)$. The commutant $\mathfrak{C}'$
is given by 
\[
\mathfrak{C}'=\left\{ X\in\mathscr{B}\left(\mathscr{H}\right)\::\:XC=CX,\:\forall C\in\mathfrak{C}\right\} .
\]
The commutant $\mathfrak{C}'$ is also a $*$-algebra, and 
\[
\mathfrak{C\;\mbox{is abelian}\Longleftrightarrow\mathfrak{C}\subset\mathfrak{C}'}.
\]
Note that $\mathfrak{C}\subset\mathfrak{C}''$ (double-commutant.)
\begin{thm}[von Neumann]
 If $M$ is a von Neumann algebra, then $M=M''$. \end{thm}
\begin{proof}
See, e.g., \cite{BR79,MR1468229}. \end{proof}
\begin{defn}
Let $\pi\in Rep(\mathfrak{A},\mathscr{H})$. We say that $\pi$ has
\emph{multiplicity} $n$, $n\in\left\{ 0\right\} \cup\mathbb{N}$,
if $\pi\left(\frak{A}\right)'\simeq M_{n}\left(\mathbb{C}\right)$,
i.e., the commutant $\pi\left(\mathfrak{A}\right)'$ is $*$-isomorphic
to the algebra of all $n\times n$ complex matrices. $\pi$ is said
to be \emph{multiplicity-free} if $\pi\left(\mathfrak{A}\right)'\simeq\mathbb{C}I_{\mathscr{H}}$.\end{defn}
\begin{example}
Let 
\[
A=\left[\begin{array}{ccc}
1\\
 & 1\\
 &  & 2
\end{array}\right]=\left[\begin{array}{cc}
I_{2} & 0\\
0 & 2
\end{array}\right].
\]
Let $C\in M_{3}\left(\mathbb{C}\right)$, then $AC=CA$ if and only
if $C$ has the form 
\[
C=\left[\begin{array}{ccc}
a & b\\
c & d\\
 &  & 1
\end{array}\right]=\left[\begin{array}{cc}
B & 0\\
0 & 1
\end{array}\right]
\]
where $B\in M_{2}(\mathbb{C})$.
\end{example}

Let $A$ be a linear operator (not necessarily bounded) acting in
the Hilbert space $\mathscr{H}$. By the Spectral Theorem (\chapref{sp}),
we have $A=A^{*}$ if and only if 
\[
A=\int_{sp(A)}\lambda P_{A}\left(d\lambda\right);
\]
where $P_{A}$ is the corresponding projection-valued measure (PVM). 
\begin{example}
The simplest example of a PVM is when $\mathscr{H}=L^{2}\left(X,\mu\right)$,
for some compact Hausdorff space $X$, and $P\left(\omega\right):=\chi_{\omega}$,
for all Borel subsets $\omega$ in $X$. Indeed, the Spectral Theorem
states that \uline{all} PVMs come this way.
\end{example}

\begin{example}
Let $A$ be compact and selfadjoint. We may further assume that $A$
is positive, $A\geq0$, in the usual order of Hermitian operators
(i.e., $\left\langle x,Ax\right\rangle \geq0$, $\forall x\in\mathscr{H}$.)
Then by \thmref{spcpt}, $A$ has the decomposition 
\begin{equation}
A=\sum_{n=1}^{\infty}\lambda_{n}P_{n}\label{eq:cptsa}
\end{equation}
where $\lambda_{n}'s$ are the eigenvalues of $A$, such that $\lambda_{1}\geq\lambda_{2}\geq\cdots\lambda_{n}\rightarrow0$;
and $P_{n}'s$ are the selfadjoint projections onto the (finite dimensional)
eigenspace of $\lambda_{n}$. In this case, the projection-valued
measure $P_{A}$ is supported on $\mathbb{N}$, and $P_{A}\left(\left\{ n\right\} \right)=P_{n}$,
$\forall n\in\mathbb{N}$. 
\end{example}
In (\ref{eq:cptsa}), we may arrange the eigenvalues as follows:\index{eigenvalue}
\begin{equation}
\overset{s_{1}}{\overbrace{\lambda_{1}=\cdots=\lambda_{1}}}>\overset{s_{2}}{\overbrace{\lambda_{2}=\cdots=\lambda_{2}}}>\cdots>\overset{s_{n}}{\overbrace{\lambda_{n}=\cdots=\lambda_{n}}}>\cdots\rightarrow0.\label{eq:multi}
\end{equation}
We say that $\lambda_{i}$ has multiplicity $s_{i}$, i.e., the dimension
of the eigenspace of $\lambda_{i}$. Note that 
\[
\dim\mathscr{H}=\sum_{i=1}^{\infty}s_{i}.
\]

Question: What does $A$ look like if it is represented as the operator
of multiplication by the independent variable? 
\begin{example}
Let $s_{1},s_{2},\ldots$ be a sequence in $\mathbb{N}$, set 
\[
E_{k}=\left\{ x_{1}^{\left(k\right)},\ldots,x_{s_{k}}^{\left(k\right)}\right\} \subset\mathbb{C},\;\mbox{and}\;E=\bigcup_{k=1}^{\infty}E_{k}.
\]
Let $\mathscr{H}=l^{2}\left(E\right)$, and 
\[
f:=\sum_{k=1}^{\infty}\lambda_{k}\chi_{E_{k}},\;\mbox{s.t}.\:\lambda_{1}>\lambda_{2}>\cdots>\lambda_{n}\rightarrow0.
\]
Let We represent $A$ as the operator $M_{f}$ of multiplication by
$f$ on $L^{2}(X,\mu)$. Let $E_{k}=\{x_{k,1},\ldots,x_{k,s_{k}}\}\subset X$,
and let $\mathscr{H}_{k}=span\{\chi_{\{x_{k,j}\}}:j\in\{1,2,\ldots,s_{k}\}\}$.
Let Notice that $\chi_{E_{k}}$ is a rank $s_{1}$ projection. $M_{f}$
is compact if and only if it is of the given form.
\end{example}

\begin{example}
Follow the previous example, we represent $A$ as the operator $M_{t}$
of multiplication by the independent variable on some Hilbert space
$L^{2}(\mu_{f})$. For simplicity, let $\lambda>0$ and 
\[
f=\lambda\chi_{\{x_{1},x_{2}\}}=\lambda\chi_{\{x_{1}\}}+\lambda\chi_{\{x_{2}\}}
\]
i.e. $f$ is compact since it is $\lambda$ times a rank-2 projection;
$f$ is positive since $\lambda>0$. The eigenspace of $\lambda$
has two dimension, 
\[
M_{f}\chi_{\{x_{i}\}}=\lambda\chi_{\{x_{i}\}},\quad i=1,2.
\]
Define $\mu_{f}(\cdot)=\mu\circ f^{-1}(\cdot)$, then 
\[
\mu_{f}=\mu(\{x_{1}\})\delta_{\lambda}\oplus\mu(\{x_{2}\})\delta_{\lambda}\oplus\text{cont. sp }\delta_{0}
\]
and 
\[
L^{2}(\mu_{f})=L^{2}(\mu(\{x_{1}\})\delta_{\lambda})\oplus L^{2}(\mu(\{x_{2}\})\delta_{\lambda})\oplus L^{2}(\text{cont. sp }\delta_{0}).
\]
Define $U:L^{2}(\mu)\rightarrow L^{2}(\mu_{f})$ by 
\[
(Ug)=g\circ f^{-1}.
\]
$U$ is unitary, and the following diagram commute: 
\[
\xymatrix{L^{2}(X,\mu)\ar[d]^{U}\ar[r]^{M_{f}} & L^{2}(X,\mu)\ar[d]^{U}\\
L^{2}(\mathbb{R},\mu_{f})\ar[r]^{M_{t}} & L^{2}(\mathbb{R},\mu_{f})
}
\]
To check $U$ preserves the $L^{2}$-norm,
\begin{eqnarray*}
\left\Vert Ug\right\Vert ^{2} & = & \int\left\Vert g\circ f^{-1}(\{x\})\right\Vert ^{2}d\mu_{f}\\
 & = & \left\Vert g\circ f^{-1}(\{\lambda\})\right\Vert ^{2}+\left\Vert g\circ f^{-1}(\{0\})\right\Vert ^{2}\\
 & = & \left|g(x_{1})\right|^{2}\mu(\{x_{1}\})+\left|g(x_{2})\right|^{2}\mu(\{x_{2}\})+\int_{X\backslash\{x_{1},x_{2}\}}\left|g(x)\right|^{2}d\mu\\
 & = & \int_{X}\left|g(x)\right|^{2}d\mu
\end{eqnarray*}
To see $U$ diagonalizes $M_{f}$, 
\begin{eqnarray*}
M_{t}Ug & = & \lambda g(x_{1})\oplus\lambda g(x_{2})\oplus0g(t)\chi_{X\backslash\{x_{1},x_{2}\}}\\
 & = & \lambda g(x_{1})\oplus\lambda g(x_{2})\oplus0\\
UM_{f}g & = & U(\lambda g(x)\chi_{\{x_{1},x_{2}\}})\\
 & = & \lambda g(x_{1})\oplus\lambda g(x_{2})\oplus0
\end{eqnarray*}
Thus 
\[
M_{t}U=UM_{f}.
\]
\end{example}
\begin{rem}
Notice that $f$ should really be written as
\[
f=\lambda\chi_{\{x_{1},x_{2}\}}=\lambda\chi_{\{x_{1}\}}+\lambda\chi_{\{x_{2}\}}+0\chi_{X\backslash\{x_{1},x_{2}\}}
\]
since $0$ is also an eigenvalue of $M_{f}$, and the corresponding
eigenspace is the kernel of $M_{f}$. \index{kernel!-of operator}\end{rem}
\begin{example}
diagonalize $M_{f}$ on $L^{2}(\mu)$ where $f=\mathbf{\chi_{[0,1]}}$
and $\mu$ is the Lebesgue measure on $\mathbb{R}$.
\end{example}

\begin{example}
diagonalize $M_{f}$ on $L^{2}(\mu)$ where 
\[
f(x)=\begin{cases}
2x & x\in[0,1/2]\\
2-2x & x\in[1/2,1]
\end{cases}
\]
and $\mu$ is the Lebesgue measure\index{measure!Lebesgue} on $[0,1]$.\end{example}
\begin{rem}
see direct integral and disintegration of measures.
\end{rem}
In general, let $A$ be a selfadjoint operator acting on $\mathscr{H}$.
Then there exists a second Hilbert space $K$, a measure $\nu$ on
$\mathbb{R}$, and unitary transformation $F:\mathscr{H}\rightarrow L_{K}^{2}(\mathbb{R},\nu)$
such that
\[
M_{t}F=FA
\]
for measurable function $\varphi:\mathbb{R}\rightarrow K$, 
\[
\left\Vert \varphi\right\Vert _{L_{K}^{2}(\nu)}=\int\left\Vert \varphi(t)\right\Vert _{K}^{2}d\nu(t)<\infty.
\]

Examples that do have multiplicities in finite dimensional linear
algebra:
\begin{example}
2-d, $\lambda I$, $\{\lambda I\}'=M_{2}(\mathbb{C})$ which is not
abelian. Hence $mult(\lambda)=2$.
\end{example}

\begin{example}
3-d, 
\[
\left[\begin{array}{ccc}
\lambda_{1}\\
 & \lambda_{1}\\
 &  & \lambda_{2}
\end{array}\right]=\left[\begin{array}{cc}
\lambda_{1}I\\
 & \lambda_{2}
\end{array}\right]
\]
where $\lambda_{1}\neq\lambda_{2}$. The commutant is
\[
\left[\begin{array}{cc}
B\\
 & b
\end{array}\right]
\]
where $B\in M_{2}(\mathbb{C})$, and $b\in\mathbb{C}$. Therefore
the commutant is isomorphic to $M_{2}(\mathbb{C})$, and multiplicity
is equal to 2.
\end{example}

\begin{example}
The example of $M_{\varphi}$ with repetition. 
\[
M_{\varphi}\oplus M_{\varphi}:L^{2}(\mu)\oplus L^{2}(\mu)\rightarrow L^{2}(\mu)\oplus L^{2}(\mu)
\]
\[
\left[\begin{array}{cc}
M_{\varphi}\\
 & M_{\varphi}
\end{array}\right]\left[\begin{array}{c}
f_{1}\\
f_{2}
\end{array}\right]=\left[\begin{array}{c}
\varphi f_{1}\\
\varphi f_{2}
\end{array}\right]
\]
the commutant is this case is isomorphic to $M_{2}(\mathbb{C})$.
If we introduces tensor product, then representation space is also
written as$L^{2}(\mu)\otimes V_{2}$, the multiplication operator
is amplified to $M_{\varphi}\otimes I$, whose commutant is represented
as $I\otimes V_{2}$. Hence it's clear that the commutant is isomorphic
to $M_{2}(\mathbb{C})$. To check \index{product!tensor-}
\begin{eqnarray*}
(\varphi\otimes I)(I\otimes B) & = & \varphi\otimes B\\
(I\otimes B)(\varphi\otimes I) & = & \varphi\otimes B.
\end{eqnarray*}

\end{example}

\section*{A summary of relevant numbers from the Reference List}

For readers wishing to follow up sources, or to go in more depth with
topics above, we suggest: \cite{Arv76,BR79,Mac85,MR1141333,MR571362,MR0213906,MR1465320,MR1444086,MR887102,MR1997376,MR1839648,BJ02,MR0467330,Con90,Dixmier198101,Gli61,MR1468229,MR1468230,Seg50,MR0442701,Tay86,MJD15,Hal13,Hal15}.

\chapter{Completely Positive Maps\label{chap:cp}}
\begin{quotation}
\textquotedblleft Completely positive maps on von Neumann algebras
or between $C^{*}$-algebras have fascinated me since my days as a
graduate student.\textquotedblright{} 

--- William B. Arveson\sindex[nam]{Arveson, W., (1934-2011)}\vspace{1em}\\
\textquotedblleft \dots{} the development of mathematics is not something
one can predict, and it would be foolish to try. One reason we love
doing mathematics is that we don\textquoteright t know what lies ahead
that future research will uncover.\textquotedblright{}

--- Alain Connes\sindex[nam]{Connes, A., (1947-)}\vspace{2em}\index{Connes, A.}
\end{quotation}
The study of completely positive maps dates back five decades, but
because of a recent observation of Arveson (see e.g., \cite{MR2490227,MR2471599,MR2470932}),
they have acquired a brand new set of applications; applications to
quantum information theory (QIT). In this framework, one studies completely
positive maps on matrix algebras. They turn out to be the objects
that are dual to quantum channels. Even more: Arveson proved the converse:
that the study of quantum channels reduces to the study of unital
completely positive maps of matrix algebras. This work is part of
QIT, and it is still ongoing, with view to the study of entanglement,
entropy and channel-capacity.

In the last chapter we studied two question from the use of algebras
of operators in quantum physics: \textquotedblleft Where does the
Hilbert space come from?\textquotedblright{} And \textquotedblleft What
are the algebras of operators from which the selfadjoint observable\index{observable}s
must be selected?\textquotedblright{} An answer is given in \textquotedblleft the
Gelfand-Naimark-Segal (GNS) theorem;\textquotedblright{} a direct
correspondence between states and cyclic representations. But states
are scalar valued positive definite functions on $*$-algebras. For
a host of applications, one must instead consider operator valued
\textquotedblleft states.\textquotedblright{} For this a different
notion of positivity is needed, \textquotedblleft \emph{complete positivity}.\textquotedblright{}

The GNS construction gives a bijection between states and cyclic representations.
An extension to the GNS construction is Stinespring's completely positive
maps. It appeared in an early paper by Stinespring in 1955 \cite{Sti55}.
Arveson in 1970's greatly extended Stinespring's result using tensor
product \cite{Arv72}. He showed that completely positive maps are
the key in multivariable operator theory, and in noncommutative dynamics.
\index{extension!-of state} 

\index{Stinespring, W. F.}

\index{Gelfand-Naimark-Segal}

\index{completely positive map}

\index{cyclic!-subspace}

\index{cyclic!-representation}

\index{product!tensor-}

\index{Theorem!Stinespring's-}

\section{Motivation}

Let $\mathfrak{A}$ be a $*$-algebra with identity. Recall that a
functional $w:\mathfrak{A}\rightarrow\mathbb{C}$ is a \emph{state}
if $w(1_{\mathfrak{A}})=1$, $w(A^{*}A)\geq0$. If $\mathfrak{A}$
was a $C^{*}$-algebra, $A\geq0\Leftrightarrow sp(A)\geq0$, hence
we may take $B=\sqrt{A}$ and $A=B^{*}B$. 

Given a state $w$, the GNS construction gives a Hilbert space $\mathscr{K}$,
a cyclic vector $\Omega\in\mathscr{K}$, and a representation $\pi:\mathfrak{A}\rightarrow\mathscr{B}(\mathscr{K})$,
such that 
\[
w(A)=\left\langle \Omega,\pi(A)\Omega\right\rangle 
\]
\[
\mathscr{K}=\overline{span}\{\pi(A)\Omega:A\in\mathfrak{A}\}.
\]
Moreover, the Hilbert space is unique up to unitary equivalence.\index{unitary equivalence}\index{space!Hilbert-}

Stinespring modified the GNS construction as follows: Instead of a
state $w:\mathfrak{A}\rightarrow\mathbb{C}$, he considered a positive
map $\varphi:\mathfrak{A}\rightarrow\mathscr{B}(\mathscr{H})$, i.e.,
$\varphi$ maps positive elements in $\mathfrak{A}$ to positive operators
in $\mathscr{B}(\mathscr{H})$. $\varphi$ is a natural extension
of $w$, since $\mathbb{C}$ can be seen as a 1-dimensional Hilbert
space, and $w$ is a positive map $w:\mathfrak{A}\rightarrow\mathscr{B}(\mathbb{C})$.
He further realized that $\varphi$ being a positive map is not enough
to produce a Hilbert space and a representation. It turns out that
the condition to put on $\varphi$ is \emph{complete positivity}:
\begin{defn}
Let $\mathfrak{A}$ be a $*$-algebra. A map $\varphi:\mathfrak{A}\rightarrow\mathscr{B}\left(\mathscr{H}\right)$
is \emph{completely positive}, if for all $n\in\mathbb{N}$, 
\begin{equation}
\varphi\otimes I_{M_{n}}:\mathfrak{A}\otimes M_{n}\rightarrow\mathscr{B}(\mathscr{H}\otimes\mathbb{C}^{n})\label{eq:cp1}
\end{equation}
maps positive elements in $\mathfrak{A}\otimes M_{n}$ to positive
operators in $\mathscr{B}(\mathscr{H}\otimes\mathbb{C}^{n})$. $\varphi$
is called a completely positive map, or a CP map. (CP maps are developed
primarily for nonabelian algebras.)
\end{defn}
The algebra $M_{n}$ of $n\times n$ matrices can be seen as an $n^{2}$-dimensional
Hilbert space with an ONB given by the matrix units $\{e_{ij}\}_{i,j=1}^{n}$.
It is also a $*$-algebra generated by $\{e_{ij}\}_{i,j=1}^{n}$ such
that 
\[
e_{ij}e_{kl}=\begin{cases}
e_{il} & j=k\\
0 & j\neq k
\end{cases}
\]
Members of $\mathfrak{A}\otimes M_{n}$ are of the form
\[
\sum_{i,j}A_{ij}\otimes e_{ij}.
\]
In other words, $\mathfrak{A}\otimes M_{n}$ consists of precisely
the $\mathfrak{A}$-valued $n\times n$ matrices. Similarly, members
of $\mathscr{H}\otimes\mathbb{C}^{n}$ are the $n$-tuple column vectors
with $\mathscr{H}$-valued entries.

Let $I_{M_{n}}:M_{n}\rightarrow\mathscr{B}(\mathbb{C}^{n})$ be the
identity representation of $M_{n}$ onto $\mathscr{B}(\mathbb{C}^{n})$.
Then, 

\begin{equation}
\varphi\otimes I_{M_{n}}:\mathfrak{A}\otimes M_{n}\rightarrow\mathscr{B}(\mathscr{H})\otimes\mathscr{B}(\mathbb{C}^{n})\left(=\mathscr{B}(\mathscr{H}\otimes\mathbb{C}^{n})\right)\label{eq:cp2}
\end{equation}
\begin{equation}
\varphi\otimes I_{M_{n}}\left(\sum_{i,j}A_{ij}\otimes e_{ij}\right)=\sum_{i,j}\varphi(A_{ij})\otimes e_{ij}.\label{eq:cp3}
\end{equation}
Note the RHS in (\ref{eq:cp3}) is an $n\times n$ matrix with $\mathscr{B}(\mathscr{H})$-valued
entries.
\begin{rem}
The algebra $\mathscr{B}(\mathbb{C}^{n})$ of all bounded operators
on $\mathbb{C}^{n}$ is generated by the rank-one operators, i.e.,
\begin{equation}
I_{M_{n}}(e_{ij})=\left|e_{i}\left\rangle \right\langle e_{j}\right|.\label{eq:cp4}
\end{equation}
Hence the $e_{ij}$ on the LHS of (\ref{eq:cp3}) is seen as an element
in the algebra $M_{n}\left(\mathbb{C}\right)$, i.e., $n\times n$
complex matrices; while on the RHS of (\ref{eq:cp3}), $e_{ij}$ is
treated as the rank one operator $\left|e_{i}\left\rangle \right\langle e_{j}\right|\in\mathscr{B}(\mathbb{C}^{n})$.
Using Dirac's notation, when we look at $e_{ij}$ as operators, we
may write 
\[
e_{i,j}(e_{k})=|e_{i}\left\rangle \right\langle e_{j}|\:|e_{k}\rangle=\begin{cases}
\left|e_{i}\right\rangle  & j=k\\
0 & j\neq k
\end{cases}
\]
\[
e_{i,j}e_{kl}=|e_{i}\rangle\langle e_{j}|\;|e_{k}\rangle\langle e_{l}|=\begin{cases}
|e_{i}\rangle\langle e_{l}| & j=k\\
0 & j\neq k
\end{cases}
\]
This also shows that $I_{M_{n}}$ is in fact an algebra isomorphism.
\end{rem}
The CP condition in (\ref{eq:cp1}) is illustrated in the following
diagram. 
\[
\otimes\begin{cases}
\mathfrak{A}\rightarrow\mathscr{B}(\mathscr{H}): & A\mapsto\varphi(A)\\
M_{n}\rightarrow M_{n}: & x\mapsto I_{M_{n}}(X)=X\mbox{\,\ (identity representation of }M_{n})
\end{cases}
\]

It is saying that if $\sum_{i,j}A_{ij}\otimes e_{ij}$ is a positive
element in the algebra $\mathfrak{A}\otimes M_{n}$, then the $n\times n$
$\mathscr{B}(\mathscr{H})$-valued matrix $\sum_{i,j}\varphi(A_{ij})\otimes e_{ij}$
is a positive operator acting on the Hilbert space $\mathscr{H}\otimes\mathbb{C}^{n}$. 

Specifically, take any $v=\sum_{k=1}^{n}v_{k}\otimes e_{k}$ in $\mathscr{H}\otimes\mathbb{C}^{n}$,
we must have
\begin{eqnarray}
 &  & \left\langle \sum_{l}v_{l}\otimes e_{l},(\sum_{i,j}\varphi(A_{ij})\otimes e_{ij})(\sum_{k}v_{k}\otimes e_{k})\right\rangle \nonumber \\
 & = & \left\langle \sum_{l}v_{l}\otimes e_{l},\sum_{i,j,k}\varphi(A_{ij})v_{k}\otimes e_{ij}(e_{k})\right\rangle \nonumber \\
 & = & \left\langle \sum_{l}v_{l}\otimes e_{l},\sum_{i,j}\varphi(A_{ij})v_{j}\otimes e_{i}\right\rangle \nonumber \\
 & = & \sum_{i,j,l}\left\langle v_{l},\varphi\left(A_{ij}\right)v_{j}\right\rangle \left\langle e_{l},e_{i}\right\rangle \nonumber \\
 & = & \sum_{i,j}\left\langle v_{i},\varphi\left(A_{ij}\right)v_{j}\right\rangle \geq0.\label{eq:cp5}
\end{eqnarray}
Using matrix notation, the CP condition is formulated as:

For all $n\in\mathbb{N}$, and all $v\in\mathscr{H}\otimes\mathbb{C}^{n}$,
i.e., 
\[
v=\sum_{k=1}^{n}v_{k}\otimes e_{k}=\begin{bmatrix}v_{1}\\
\vdots\\
v_{n}
\end{bmatrix}
\]
we have
\begin{equation}
\left[\begin{array}{cccc}
v_{1} & v_{2} & \cdots & v_{n}\end{array}\right]\left[\begin{array}{cccc}
\varphi(A_{11}) & \varphi(A_{12}) & \cdots & \varphi(A_{1n})\\
\varphi(A_{21}) & \varphi(A_{22}) & \cdots & \varphi(A_{2n})\\
\vdots & \vdots & \ddots & \vdots\\
\varphi(A_{n1}) & \varphi(A_{n2}) & \cdots & \varphi(A_{nn})
\end{array}\right]\left[\begin{array}{c}
v_{1}\\
v_{2}\\
\vdots\\
v_{n}
\end{array}\right]\geq0.\label{eq:cp6}
\end{equation}

\section{\label{sec:cpgns}CP v.s. GNS}

The GNS construction can be reformulated as a special case of the
Stinespring's theorem \cite{Sti55}. \index{representation!GNS} \index{Theorem!Stinespring's-}

Let $\mathfrak{A}$ be a $*$-algebra, given a state $\varphi:\mathfrak{A}\rightarrow\mathbb{C}$,
there exists a triple $(\mathscr{K},\Omega,\pi)$, all depending on
$\varphi$, such that
\[
\varphi(A)=\left\langle \Omega,\pi\left(A\right)\Omega\right\rangle _{\mathscr{K}}
\]
where 
\begin{eqnarray*}
\Omega & = & \pi\left(1_{\mathfrak{A}}\right)\in\mathscr{K}\\
\mathscr{K} & = & \overline{span}\{\pi(A)\Omega:A\in\mathfrak{A}\}.
\end{eqnarray*}

The 1-dimensional Hilbert space $\mathbb{C}$ is thought of being
embedded into $\mathscr{K}$ (possibly infinite dimensional) via 
\begin{equation}
\mathbb{C}\ni t\xrightarrow{\;V\;}t\Omega\in\mathbb{C}\Omega\label{eq:cp1-1}
\end{equation}
where $\mathbb{C}\Omega=$ the one-dimensional subspace in $\mathscr{K}$
generated by the unit cyclic vector cyclic $\Omega$. 
\begin{lem}
The map $V$ in (\ref{eq:cp1-1}) is an isometry, such that $V^{*}V=I_{\mathbb{C}}:\mathbb{C}\rightarrow\mathbb{C}$,
and \index{isometry!Cuntz-} 
\begin{equation}
VV^{*}:\mathscr{K}\rightarrow\mathbb{C}\Omega\label{eq:cp7}
\end{equation}
is the projection from $\mathscr{K}$ onto the 1-d subspace $\mathbb{C}\Omega$
in $\mathscr{K}$. 

Moreover, 
\begin{equation}
\varphi\left(A\right)=V^{*}\pi\left(A\right)V,\;\forall A\in\mathfrak{A}.\label{eq:cp8}
\end{equation}
\end{lem}
\begin{proof}
Let $t\in\mathbb{C}$, then $\left\Vert Vt\right\Vert _{\mathscr{K}}=\left\Vert t\Omega\right\Vert _{\mathscr{K}}=\left|t\right|$,
and so $V$ is an isometry. 

For all $\xi\in\mathscr{K}$, we have 
\[
\left\langle \xi,Vt\right\rangle _{\mathscr{K}}=\left\langle V^{*}\xi,t\right\rangle _{\mathbb{C}}=t\overline{V^{*}\xi}.
\]
By setting $t=1$, we get
\[
V^{*}\xi=\overline{\left\langle \xi,V1\right\rangle _{\mathscr{K}}}=\overline{\left\langle \xi,\Omega\right\rangle _{\mathscr{K}}}=\left\langle \Omega,\xi\right\rangle _{\mathscr{K}}\Longleftrightarrow V^{*}=\left\langle \Omega,\cdot\right\rangle _{\mathscr{K}}.
\]
Therefore,
\[
V^{*}Vt=V^{*}\left(t\Omega\right)=\left\langle \Omega,t\Omega\right\rangle _{\mathscr{K}}=t,\;\forall t\in\mathbb{C}\Longleftrightarrow V^{*}V=I_{\mathbb{C}}
\]
\[
VV^{*}\xi=V\left(\left\langle \Omega,\xi\right\rangle _{\mathscr{K}}\right)=\left\langle \Omega,\xi\right\rangle _{\mathscr{K}}\Omega,\;\forall\xi\in\mathscr{K}\Longleftrightarrow VV^{*}=\left|\Omega\left\rangle \right\langle \Omega\right|.
\]
It follows that 
\begin{eqnarray*}
\varphi\left(A\right) & = & \left\langle \Omega,\pi\left(A\right)\Omega\right\rangle _{\mathscr{K}}\\
 & = & \left\langle V1,\pi\left(A\right)V1\right\rangle _{\mathscr{K}}\\
 & = & \left\langle 1,V^{*}\pi\left(A\right)V1\right\rangle _{\mathbb{C}}\\
 & = & V^{*}\pi\left(A\right)V,\;\forall A\in\mathfrak{A}
\end{eqnarray*}
which is the assertion in (\ref{eq:cp8}).
\end{proof}
In other words, $\Omega\longmapsto\pi(A)\Omega$ sends the unit vector
$\Omega$ from the 1-dimensional subspace $\mathbb{C}\Omega$ to the
vector $\pi(A)\Omega\in\mathscr{K}$, and $\left\langle \Omega,\pi\left(A\right)\Omega\right\rangle _{\mathscr{K}}$
cuts off the resulting vector $\pi(A)\Omega$ and only preserves the
component corresponding to the 1-d subspace $\mathbb{C}\Omega$. Notice
that the unit vector $\Omega$ is obtained from embedding the constant
$1\in\mathbb{C}$ via the map $V$, i.e., $\Omega=V1$. In matrix
notation, if we identify $\mathbb{C}$ with its image $\mathbb{C}\Omega$
in $\mathscr{K}$, then $\varphi(A)$ is put into a matrix corner:
\[
\pi(A)=\left[\begin{array}{cc}
\varphi(A) & *\\*
* & *
\end{array}\right]
\]
so that when acting on vectors,
\[
\varphi(A)=\left[\begin{array}{cc}
\Omega & 0\end{array}\right]\left[\begin{array}{cc}
\varphi(A) & *\\*
* & *
\end{array}\right]\left[\begin{array}{c}
\Omega\\
0
\end{array}\right].
\]
Equivalently
\[
\varphi(A)=P_{1}\pi(A):P_{1}\mathscr{K}\rightarrow\mathbb{C};
\]
where $P_{1}:=VV^{*}=\left|\Omega\left\rangle \right\langle \Omega\right|=$
rank-1 projection on $\mathbb{C}\Omega$.

Stinespring's construction is a generalization of the above formulation:
Let $\mathfrak{A}$ be a $*$-algebra, given a CP map $\varphi:\mathfrak{A}\rightarrow\mathscr{B}(\mathscr{H})$,
there exists a Hilbert space $\mathscr{K}\left(=\mathscr{K}_{\varphi}\right)$,
an isometry $V:\mathscr{H}\rightarrow\mathscr{K}$, and a representation
$\pi\left(=\pi_{\varphi}\right):\mathfrak{A}\rightarrow\mathscr{K}$,
such that \index{isometry} 
\[
\varphi(A)=V^{*}\pi(A)V,\;\forall A\in\mathfrak{A}.
\]
Notice that this construction starts with a possibly infinite dimensional
Hilbert space $\mathscr{H}$ (instead of the 1-dimensional Hilbert
space $\mathbb{C}$), the map $V$ embeds $\mathscr{H}$ into a bigger
Hilbert space $\mathscr{K}$. If $\mathscr{H}$ is identified with
its image in $\mathscr{K}$, then $\pi(A)$ is put into a matrix corner,
\[
\left[\begin{array}{cc}
\pi(A) & *\\*
* & *
\end{array}\right]
\]
so that when acting on vectors,
\[
\varphi(A)\xi=\left[\begin{array}{cc}
V\xi & 0\end{array}\right]\left[\begin{array}{cc}
\pi(A) & *\\*
* & *
\end{array}\right]\left[\begin{array}{c}
V\xi\\
0
\end{array}\right].
\]

This can be formulated alternatively: 

For every CP map $\varphi:\mathfrak{A}\rightarrow\mathscr{B}(\mathscr{H})$,
there is a dilated Hilbert space $\mathscr{K}\left(=\mathscr{K}_{\varphi}\right)\supset\mathscr{H}$,
a representation $\pi\left(=\pi_{\varphi}\right):\mathfrak{A}\rightarrow\mathscr{B}(\mathscr{K})$,
such that 
\[
\varphi(A)=P_{\mathscr{H}}\pi(A)
\]
\index{space!Hilbert-}i.e., $\pi(A)$ can be put into a matrix corner.
$\mathscr{K}$ is chosen as minimal in the sense that 
\[
\mathscr{K}=\overline{span}\{\pi(A)(Vh):A\in\mathfrak{A},h\in\mathscr{H}\}.
\]
\[
\xymatrix{\mathscr{H}\ar[r]^{V}\ar[d]_{\varphi\left(A\right)} & \mathscr{K}\ar[d]^{\pi\left(A\right)}\\
\mathscr{H}\ar[r]_{V} & \mathscr{K}
}
\]

\begin{note}
The containment $\mathscr{H}\subset\mathscr{K}$ comes after the identification
of $\mathscr{H}$ with its image in $\mathscr{K}$ under the isometric
embedding $V$. We write $\varphi(A)=P_{H}\pi(A)$, as opposed to
$\varphi(A)=P_{H}\pi(A)P_{H}$, since $\varphi(A)$ only acts on the
subspace $\mathscr{H}$.
\end{note}

\section{\label{sec:sspring}Stinespring's Theorem}

\index{Stinespring, W. F.}\index{completely positive map}\index{Theorem!Stinespring's-}
\begin{thm}[Stinespring \cite{Sti55}]
\label{thm:ss} Let $\mathfrak{A}$ be a $*$-algebra. The following
are equivalent:
\begin{enumerate}
\item \label{enu:ss1}$\varphi:\mathfrak{A}\rightarrow\mathscr{B}(\mathscr{H})$
is a completely positive map, and $\varphi(1_{\mathfrak{A}})=I_{\mathscr{H}}$. 
\item \label{enu:ss2}There exists a Hilbert space $\mathscr{K}$, an isometry
$V:\mathscr{H}\rightarrow\mathscr{K}$, and a representation $\pi:\mathfrak{A}\rightarrow\mathscr{B}(\mathscr{K})$
such that
\begin{equation}
\varphi(A)=V^{*}\pi(A)V,\;\forall A\in\mathfrak{A}.\label{eq:ss1}
\end{equation}

\item \label{enu:ss3}If the dilated Hilbert space $\mathscr{K}$ is taken
to be minimum, then it is unique up to unitary equivalence. Specifically,
if there are two systems $\left(V_{i},\mathscr{K}_{i},\pi_{i}\right)$,
$i=1,2$, satisfying 
\begin{eqnarray}
\varphi(A) & = & V_{i}^{*}\pi_{i}(A)V_{i}\label{eq:ss4}\\
\mathscr{K}_{i} & = & \overline{span}\{\pi_{i}(A)Vh:A\in\mathfrak{A},h\in\mathscr{H}\}\label{eq:ss5}
\end{eqnarray}
then there exists a unitary operator $W:\mathscr{K}_{1}\rightarrow\mathscr{K}_{2}$
so that
\begin{equation}
W\pi_{1}=\pi_{2}W\label{eq:ss2}
\end{equation}

\end{enumerate}
\end{thm}
\begin{proof}
~\index{isometry}

(Part (\ref{enu:ss3}), uniqueness) Let $\left(V_{i},\mathscr{K}_{i},\pi_{i}\right)$,
$i=1,2$, be as in the statement satisfying (\ref{eq:ss4})-(\ref{eq:ss5}).
Define
\[
W\pi_{1}(A)Vh=\pi_{2}(A)Vh
\]
then $W$ is an isometry, since 
\begin{eqnarray*}
\left\Vert \pi_{i}(A)Vh\right\Vert _{\mathscr{K}}^{2} & = & \left\langle \pi_{i}(A)Vh,\pi_{i}(A)Vh\right\rangle _{\mathscr{K}}\\
 & = & \left\langle h,V^{*}\pi_{i}(A^{*}A)Vh\right\rangle _{\mathscr{H}}\\
 & = & \left\langle h,\varphi(A^{*}A)h\right\rangle _{\mathscr{H}}.
\end{eqnarray*}
Hence $W$ extends uniquely to a unitary operator $W:\mathscr{K}_{1}\rightarrow\mathscr{K}_{2}$.
To see that $W$ intertwines $\pi_{1},\pi_{2}$, notice that a typical
vector in $\mathscr{K}_{i}$ is $\pi_{i}(A)Vh$, and 
\begin{eqnarray*}
W\pi_{1}(B)\pi_{1}(A)Vh & = & W\pi_{1}(BA)Vh\\
 & = & \pi_{2}(BA)Vh\\
 & = & \pi_{2}(B)\pi_{2}(A)Vh\\
 & = & \pi_{2}(B)W\pi_{1}(A)Vh.
\end{eqnarray*}
Since such vectors are dense in the respective dilated space, we conclude
that $W\pi_{1}=\pi_{2}W$, so (\ref{eq:ss2}) holds.\end{proof}
\begin{note}
$\left\Vert \pi_{1}(B)\pi_{1}(A)Vh\right\Vert ^{2}=\left\langle h,V^{*}\pi_{1}(A^{*}B^{*}BA)Vh\right\rangle $.
Fix $A\in\mathscr{B}(\mathscr{H})$, the map $B\mapsto A^{*}BA$ is
an automorphism on $\mathscr{B}(\mathscr{H})$.\end{note}
\begin{proof}
(\ref{enu:ss2}) $\Longrightarrow$ (\ref{enu:ss1})

Now suppose $\varphi(A)=V^{*}\pi(A)V$, and we verify it is completely
positive. 

Since positive elements in $\mathfrak{A}\otimes M_{n}$ are sums of
the operator matrix
\[
\sum_{i,j}A_{i}^{*}A_{j}\otimes e_{ij}=\left[\begin{array}{c}
A_{1}^{*}\\
A_{2}^{*}\\
\vdots\\
A_{n}^{*}
\end{array}\right]\left[\begin{array}{cccc}
A_{1} & A_{2} & \cdots & A_{n}\end{array}\right]
\]
it suffices to show that
\[
\varphi\otimes I_{M_{n}}\left(\sum_{i,j}A_{i}^{*}A_{j}\otimes e_{ij}\right)=\sum_{i,j}\varphi\left(A_{i}^{*}A_{j}\right)\otimes e_{ij}
\]
is a positive operator in $\mathscr{B}(\mathscr{H}\otimes\mathbb{C}^{n})$,
i.e., need to show that for all $v\in\mathscr{H}\otimes\mathbb{C}^{n}$
\begin{equation}
\left[\begin{array}{cccc}
v_{1} & v_{2} & \cdots & v_{n}\end{array}\right]\left[\begin{array}{cccc}
\varphi\left(A_{1}^{*}A_{1}\right) & \varphi\left(A_{1}^{*}A_{2}\right) & \cdots & \varphi\left(A_{1}^{*}A_{n}\right)\\
\varphi\left(A_{2}^{*}A_{1}\right) & \varphi\left(A_{2}^{*}A_{2}\right) & \cdots & \varphi\left(A_{2}^{*}A_{n}\right)\\
\vdots & \vdots & \ddots & \vdots\\
\varphi\left(A_{n}^{*}A_{1}\right) & \varphi\left(A_{n}^{*}A_{2}\right) & \cdots & \varphi\left(A_{n}^{*}A_{n}\right)
\end{array}\right]\left[\begin{array}{c}
v_{1}\\
v_{2}\\
\vdots\\
v_{n}
\end{array}\right]\geq0.\label{eq:ss3}
\end{equation}
This is true, since
\begin{eqnarray*}
\mbox{RHS}_{\left(\ref{eq:ss3}\right)} & = & \sum_{i,j}\left\langle v_{i},\varphi\left(A_{i}^{*}A_{j}\right)v_{j}\right\rangle _{\mathscr{H}}\\
 & = & \sum_{i,j}\left\langle v_{i},V^{*}\pi\left(A_{i}^{*}A_{j}\right)Vv_{j}\right\rangle _{\mathscr{H}}\\
 & = & \sum_{i,j}\left\langle \pi\left(A_{i}\right)Vv_{i},\pi\left(A_{j}\right)Vv_{j}\right\rangle _{\mathscr{K}}\\
 & = & \left\Vert \sum_{i}\pi\left(A_{i}\right)Vv_{i}\right\Vert _{\mathscr{K}}^{2}\geq0.
\end{eqnarray*}

(\ref{enu:ss1}) $\Longrightarrow$ (\ref{enu:ss2})

Given a completely positive map $\varphi$, we construct $\mathscr{K}\left(=\mathscr{K}_{\varphi}\right)$,
$V\left(=V_{\varphi}\right)$ and $\pi\left(=\pi_{\varphi}\right)$.
Recall that $\varphi:\mathfrak{A}\rightarrow\mathscr{B}(\mathscr{H})$
is a CP map means that for all $n\in\mathbb{N}$, 
\[
\varphi\otimes I_{M_{n}}:\mathfrak{A}\otimes M_{n}\rightarrow\mathscr{B}(\mathscr{H}\otimes M_{n})
\]
is positive, and 
\[
\varphi\otimes I_{M_{n}}(1_{\mathfrak{A}}\otimes I_{M_{n}})=I_{\mathscr{H}}\otimes I_{M_{n}}.
\]
The condition on the identity element can be stated using matrix notation
as
\[
\left[\begin{array}{cccc}
\varphi & 0 & \cdots & 0\\
0 & \varphi & \cdots & 0\\
\vdots & \vdots & \ddots & \vdots\\
0 & 0 & \cdots & \varphi
\end{array}\right]\left[\begin{array}{cccc}
1_{\mathfrak{A}} & 0 & \cdots & 0\\
0 & 1_{\mathfrak{A}} & \cdots & 0\\
\vdots & \vdots & \ddots & \vdots\\
0 & 0 & \cdots & 1_{\mathfrak{A}}
\end{array}\right]=\left[\begin{array}{cccc}
I_{\mathscr{H}} & 0 & \cdots & 0\\
0 & I_{\mathscr{H}} & \cdots & 0\\
\vdots & \vdots & \ddots & \vdots\\
0 & 0 & \cdots & I_{\mathscr{H}}
\end{array}\right].
\]

Let $K_{0}$ be the algebraic tensor product $\mathfrak{A}\otimes\mathscr{H}$,
i.e., 
\[
K_{0}=span\left\{ \sum_{\text{finite}}A_{i}\otimes\xi_{i}:A\in\mathfrak{A},\xi\in\mathscr{H}\right\} .
\]
Define a sesquilinear form $\left\langle \cdot,\cdot\right\rangle _{\varphi}:K_{0}\times K_{0}\rightarrow\mathbb{C}$,
by \index{sesquilinear form} 
\begin{equation}
\left\langle \sum_{i=1}^{n}A_{i}\otimes\xi_{i},\sum_{j=1}^{n}B_{j}\otimes\eta_{j}\right\rangle _{\varphi}:=\sum_{i,j}\left\langle \xi_{i},\varphi\left(A_{i}^{*}B_{j}\right)\eta_{j}\right\rangle _{\mathscr{H}}.\label{eq:sesquilinear}
\end{equation}
 By the CP condition (\ref{eq:cp1}), we have
\[
\left\langle \sum_{i=1}^{n}A_{i}\xi_{i},\sum_{j=1}^{n}A_{j}\xi_{j}\right\rangle _{\varphi}=\sum_{i,j}\left\langle \xi_{i},\varphi(A_{i}^{*}A_{j})\xi_{j}\right\rangle _{\mathscr{H}}\geq0.
\]
Let $N:=\left\{ v\in K_{0}:\left\langle v,v\right\rangle _{\varphi}=0\right\} $.
Since the Schwarz inequality holds for any sesquilinear form, it follows
that\index{inequality!Schwarz}\index{product!tensor-} 
\[
N=\left\{ v\in K_{0}:\left\langle s,v\right\rangle _{\varphi}=0,\;\forall s\in K_{0}\right\} .
\]
Thus $N$ is a closed subspace in $K_{0}$. Let $\mathscr{K}\left(=\mathscr{K}_{\varphi}\right)$
be the Hilbert space by completing $K_{0}/N$ with respect to 
\[
\left\Vert \cdot\right\Vert _{\mathscr{K}}:=\left\langle \cdot,\cdot\right\rangle _{\varphi}^{1/2}.
\]

Let $V:\mathscr{H}\rightarrow K_{0}$, by 
\[
V\xi:=1_{\mathfrak{A}}\otimes\xi,\;\forall\xi\in\mathscr{H}.
\]
Then, 
\begin{eqnarray*}
\left\Vert V\xi\right\Vert _{\varphi}^{2} & = & \left\langle 1_{\mathfrak{A}}\otimes\xi,1_{\mathfrak{A}}\otimes\xi\right\rangle _{\varphi}\\
 & = & \left\langle \xi,\varphi(1_{\mathfrak{A}}^{*}1_{\mathfrak{A}})\xi\right\rangle _{\mathscr{H}}\\
 & = & \left\langle \xi,\xi\right\rangle _{\mathscr{H}}=\left\Vert \xi\right\Vert _{\mathscr{H}}^{2}
\end{eqnarray*}
i.e., $V$ is isometric, and so $\mathscr{H}\xrightarrow{\;V\;}K_{0}$
is an isometric embedding. \index{operators!isometric}

\noindent \textbf{Claim.} 

(i) $V^{*}V=I_{\mathscr{H}}$; 

(ii) $VV^{*}=$ projection from $K_{0}$ on the subspace $1_{\mathfrak{A}}\otimes\mathscr{H}$.

Indeed, for any $A\otimes\eta\in K_{0}$, we have
\begin{eqnarray*}
\left\langle A\otimes\eta,V\xi\right\rangle _{\varphi} & = & \left\langle A\otimes\eta,1_{\mathfrak{A}}\otimes\xi\right\rangle _{\varphi}\\
 & = & \left\langle \eta,\varphi(A^{*})\xi\right\rangle _{\mathscr{H}}\\
 & = & \left\langle \varphi(A^{*})^{*}\eta,\xi\right\rangle _{\mathscr{H}}
\end{eqnarray*}
which implies that
\[
V^{*}(A\otimes\eta)=\varphi(A^{*})^{*}\eta.
\]
It follows that 
\[
V^{*}V\xi=V^{*}(1_{\mathfrak{A}}\otimes\xi)=\varphi(1_{\mathfrak{A}}^{*})^{*}\xi=\xi,\;\forall\xi\in\mathscr{H}
\]
i.e., $V^{*}V=I_{\mathscr{H}}$. Moreover, for any $A\otimes\eta\in K_{0}$,
\[
VV^{*}(A\otimes\eta)=V(\varphi(A^{*})^{*}\eta)=1_{\mathfrak{A}}\otimes\varphi(A^{*})^{*}\eta.
\]
This proves the claim. It is clear that the properties of $V$ pass
to the dilated space $\mathscr{K}\left(=\mathscr{K}_{\varphi}\right)=cl_{\varphi}\left(K_{0}/N\right)$.
\\

To finish the proof of the theorem, define $\pi\left(=\pi_{\varphi}\right)$
as follows: Set 

\[
\pi(A)\left(\sum_{j}B_{j}\otimes\eta_{j}\right):=\sum_{j}AB_{j}\otimes\eta_{j},\;\forall A\in\mathfrak{A}
\]
and extend it to $\mathscr{K}$. 

For all $\xi,\eta\in\mathscr{H}$, then, 
\begin{eqnarray*}
\left\langle \xi,V^{*}\pi(A)V\eta\right\rangle _{\mathscr{H}} & = & \left\langle V\xi,\pi(A)V\eta\right\rangle _{\mathscr{K}}\\
 & = & \left\langle 1_{\mathfrak{A}}\otimes\xi,\pi(A)1_{\mathfrak{A}}\otimes\eta\right\rangle _{\mathscr{K}}\\
 & = & \left\langle 1_{\mathfrak{A}}\otimes\xi,A\otimes\eta\right\rangle _{\mathscr{K}}\\
 & = & \left\langle \xi,\varphi(1_{\mathfrak{A}}^{*}A)\eta\right\rangle _{\mathscr{H}}\\
 & = & \left\langle \xi,\varphi(A)\eta\right\rangle _{\mathscr{H}}.
\end{eqnarray*}
We conclude that $\varphi(A)=V^{*}\pi(A)V$, for all $A\in\mathfrak{A}$.
\end{proof}
\begin{flushleft}
\textbf{Application of Stinespring's Theorem to Representations of
$\mathscr{O}_{N}$}
\par\end{flushleft}
\begin{cor}
Let $N\in\mathbb{N}$, $N>1$, and let $A_{i}\in\mathscr{B}\left(\mathscr{H}\right)$,
$1\leq i\leq N$, be a system of operators in a Hilbert space $\mathscr{H}$
such that 
\begin{equation}
\sum_{i=1}^{N}A_{i}^{*}A_{i}=I_{\mathscr{H}};\label{eq:sti1}
\end{equation}
then there is a second Hilbert space $\mathscr{K}$, and an isometry
$V:\mathscr{H}\rightarrow\mathscr{K}$, and a representation $\pi\in Rep\left(\mathscr{O}_{N},\mathscr{K}\right)$
such that 
\begin{equation}
V^{*}\pi(s_{i})V=A_{i}^{*},\;1\leq i\leq N,\label{eq:sti2}
\end{equation}
where $\left\{ s_{i}\right\} _{i=1}^{N}$ are generators for $\mathscr{O}_{N}$. \end{cor}
\begin{proof}
Given $\mathscr{O}_{N}$ with generators $\left\{ s_{i}\right\} _{i=1}^{N}$,
then set 
\[
\varphi(s_{i}s_{j}^{*})=A_{i}^{*}A_{j}
\]
using (\ref{eq:sti1}), it is easy to see that $\varphi$ is completely
positive. 

Now let $\left(\pi,\mathscr{K}\right)$ be the pair obtained from
\thmref{ss} (Stinespring); then as a block-matrix of operators, we
have as follows \index{matrix!block-}
\begin{equation}
\pi\left(s_{i}\right)^{*}=\begin{bmatrix}A_{i} & *\\
\mathbf{0} & *
\end{bmatrix}\label{eq:sti3}
\end{equation}
relative to the splitting
\begin{equation}
\mathscr{K}=V\mathscr{H}\oplus\left(\mathscr{K}\ominus V\mathscr{H}\right),\label{eq:sti4}
\end{equation}
and so $V^{*}\pi\left(s_{i}\right)^{*}V=A_{i}$, which is equivalent
to (\ref{eq:sti2}).
\end{proof}

\section{\label{sec:Comments}Comments}

In Stinespring's theorem, the dilated space comes from a general principle
(using positive definite functions) when building Hilbert spaces out
of the given data. We illustrate this point with a few familiar examples.
\index{positive definite!-function}
\begin{example}
In linear algebra, there is a bijection between inner product structures
on $\mathbb{C}^{n}$ and positive-definite $n\times n$ matrices.
Specifically, $\left\langle \cdot,\cdot\right\rangle :\mathbb{C}^{n}\times\mathbb{C}^{n}\rightarrow\mathbb{C}$
is an inner product if and only if there exists a positive definite
matrix $A$ such that
\[
\left\langle v,w\right\rangle _{A}=v^{*}Aw
\]
for all $v,w\in\mathbb{C}^{n}$. We think of $\mathbb{C}^{n}$ as
$\mathbb{C}$-valued functions on $\{1,2,\ldots,n\}$, then $\left\langle \cdot,\cdot\right\rangle _{A}$
is an inner product built on the function space. 
\end{example}
This is then extended to infinite dimensional space. 
\begin{example}
If $F$ is a positive definite function on $\mathbb{R}$, then on
$K_{0}=span\{\delta_{x}:x\in\mathbb{R}\}$, $F$ defines a sesquilinear
form $\left\langle \cdot,\cdot\right\rangle _{F}:\mathbb{R}\times\mathbb{R}\rightarrow\mathbb{C}$,
where
\[
\left\langle \sum_{i}c_{i}\delta_{x_{i}},\sum_{j}d_{j}\delta_{x_{j}}\right\rangle _{F}:=\sum_{i,j}\overline{c_{i}}d_{j}F(x_{i},x_{j}),\;\mbox{and}
\]
\[
\left\Vert \sum_{i}c_{i}\delta_{x_{i}}\right\Vert _{F}^{2}:=\left\langle \sum_{i}c_{i}\delta_{x_{i}},\sum_{j}c_{j}\delta_{x_{j}}\right\rangle _{F}=\sum_{i,j}\overline{c_{i}}c_{j}F\left(x_{i},x_{j}\right)\geq0.
\]
Let $N=\left\{ v\in K_{0}:\left\langle v,v\right\rangle =0\right\} $,
then $N$ is a closed subspace in $K_{0}$. We get a Hilbert space:
$\mathscr{K}:=cl_{F}\left(K_{0}/N\right)=$ the completion of $K_{0}/N$
with respect to $\left\Vert \cdot\right\Vert _{F}$. \index{sesquilinear form}
\end{example}
What if the index set is not $\{1,2,\ldots,n\}$ or $\mathbb{R}$,
but a $*$-algebra?
\begin{example}
$C\left(X\right)$, $X$ compact Hausdorff. It is a $C^{*}$-algebra,
where $\left\Vert f\right\Vert :=\sup_{x}\left|f\left(x\right)\right|$.
By Riesz's theorem, there is a bijection between positive states (linear
functionals) on $C(X)$ and Borel probability measures on $X$. \index{functional}
\index{Riesz' theorem}\index{Theorem!Riesz-}

Let $\mathfrak{B}\left(X\right)$ be the Borel sigma-algebra\index{sigma-algebra}
on $X$, which is also an abelian algebra: The associative multiplication
is defined as $AB:=A\cap B$. The identity element is just $X$.

Let $\mu$ be a probability measure, then $\mu(A\cap B)\geq0$, for
all $A,B\in\mathfrak{B}\left(X\right)$, and $\mu(X)=1$. Hence $\mu$
is a state. As before, we apply the GNS construction. Set 
\[
K_{0}=span\{\delta_{A}:A\in\mathfrak{M}\}=span\{\chi_{A}:A\in\mathfrak{M}\}
\]
Note the index set here is $\mathfrak{B}\left(X\right)$, and $\sum_{i}c_{i}\delta_{A_{i}}=\sum_{i}c_{i}\chi_{A_{i}}$,
i.e., these are precisely the simple functions . Define
\[
\left\langle \sum_{i}c_{i}\chi_{A_{i}},\sum_{j}d_{j}\chi_{B_{j}}\right\rangle :=\sum_{i,j}\overline{c_{i}}d_{j}\mu(A_{i}\cap B_{j})
\]
which is positive definite, since 
\[
\left\langle \sum_{i}c_{i}\chi_{A_{i}},\sum_{i}c_{i}\chi_{A_{i}}\right\rangle =\sum_{i,j}\overline{c_{i}}c_{j}\mu(A_{i}\cap A_{j})=\sum_{i}\left|c_{i}\right|^{2}\mu\left(A_{i}\right)\geq0.
\]
Here, $N=\left\{ v\in K_{0}:\left\langle v,v\right\rangle =0\right\} =\mu$-measure
zero sets, and 
\[
\mathscr{H}=cl_{\mu}\left(K_{0}/N\right)=L^{2}\left(\mu\right).
\]
 
\end{example}

\begin{example}[GNS]
Let $\mathfrak{A}$ be a $*$-algebra. The set of $\mathbb{C}$-valued
functions on $\mathfrak{A}$ is $\mathfrak{A}\otimes\mathbb{C}$,
i.e., functions of the form 
\begin{equation}
\left\{ \sum_{i}A_{i}\otimes c_{i}=\sum_{i}c_{i}\delta_{A_{i}}\right\} \label{eq:co1}
\end{equation}
with finite summation over $i$. Note that $\mathbb{C}$ is naturally
embedded into $\mathfrak{A}\otimes\mathbb{C}$ as $1_{\mathfrak{A}}\otimes\mathbb{C}$
(i.e., $c\mapsto c\delta_{1_{\mathfrak{A}}}$), and the latter is
a 1-dimensional subspace. In order to build a Hilbert space out of
(\ref{eq:co1}), one needs a positive definite function. A state $\varphi$
on $\mathfrak{A}$ does exactly the job. The sesquilinear form is
given by 
\[
\left\langle \sum_{i}c_{i}\delta_{A_{i}},\sum_{i}d_{j}\delta_{B_{j}}\right\rangle _{\varphi}:=\sum_{i,j}\overline{c_{i}}d_{j}\varphi\left(A_{i}^{*}B_{j}\right)
\]
so that 
\[
\left\Vert \sum_{i}c_{i}\delta_{A_{i}}\right\Vert _{\varphi}^{2}=\sum_{i,j}\overline{c_{i}}c_{j}\varphi\left(A_{i}^{*}A_{j}\right)\geq0.
\]
Finally, let $\mathscr{K}_{\varphi}=$ Hilbert completion of $\mathfrak{A}\otimes\mathbb{C}/\ker\varphi$.
Define $\pi(A)\delta_{B}:=\delta_{BA}$, so a ``shift'' in the index
variable, and extend to $\mathscr{K}_{\varphi}$. \index{shift}
\end{example}

In Stinespring's construction, $\mathfrak{A}\otimes\mathbb{C}$ is
replaced by $\mathfrak{A}\otimes\mathscr{H}$, i.e., in stead of working
with $\mathbb{C}$-valued functions on $\mathfrak{A}$, one looks
at $\mathscr{H}$-valued functions. Hence we are looking at functions
of the form
\[
\left\{ \sum_{i}A_{i}\otimes\xi_{i}=\sum_{i}\xi_{i}\delta_{A_{i}}\right\} 
\]
with finite summation over $i$. $\mathscr{H}$ is embedded into $\mathfrak{A}\otimes\mathscr{H}$
as $1_{\mathfrak{A}}\otimes\mathscr{H}$, by $\mathscr{H}\ni1_{\mathfrak{A}}\otimes\xi=\xi\delta_{1_{\mathfrak{A}}}$.
$1_{\mathfrak{A}}\otimes\mathscr{H}$ is in general infinite dimensional,
or we say that the function $\xi\delta_{1_{\mathfrak{A}}}$ at $1_{\mathfrak{A}}$
has infinite multiplicity\index{multiplicity}. If $\mathscr{H}$
is separable, we are actually attaching an $l^{2}$ sequence at every
point $A\in\mathfrak{A}$. 

How to build a Hilbert space out of these $\mathscr{H}$-valued functions?
The question depends on the choice of a quadratic form. If $\varphi:\mathfrak{A}\rightarrow\mathscr{B}(\mathscr{H})$
is positive, i.e., $\varphi$ maps positive elements in $\mathfrak{A}$
to positive operators on $\mathscr{H}$, then quadratic form \index{quadratic form}
\[
\left\langle A\otimes\xi,B\otimes\eta\right\rangle _{\varphi}:=\left\langle \xi,\varphi(A^{*}B)\eta\right\rangle _{\mathscr{H}}
\]
is indeed positive definite. But when extend linearly, one is in trouble.
For \index{positive definite!-function} 
\begin{eqnarray*}
 &  & \left\langle \sum_{i}A_{i}\otimes\xi_{i},\sum_{j}B_{j}\otimes\eta_{j}\right\rangle _{\varphi}\\
 & = & \sum_{i,j}\left\langle \xi_{i},\varphi(A_{i}^{*}B_{j})\eta_{j}\right\rangle _{\mathscr{H}}\\
 & = & \left[\begin{array}{cccc}
\xi_{1} & \xi_{2} & \cdots & \xi_{n}\end{array}\right]\left[\begin{array}{cccc}
\varphi(A_{1}^{*}B_{1}) & \varphi(A_{1}^{*}B_{2}) & \cdots & \varphi(A_{1}^{*}B_{n})\\
\varphi(A_{2}^{*}B_{1}) & \varphi(A_{2}^{*}B_{1}) & \cdots & \varphi(A_{2}^{*}B_{1})\\
\vdots & \vdots & \vdots & \vdots\\
\varphi(A_{n}^{*}B_{1}) & \varphi(A_{n}^{*}B_{2}) & \cdots & \varphi(A_{n}^{*}B_{n})
\end{array}\right]\left[\begin{array}{c}
\xi_{1}\\
\xi_{2}\\
\vdots\\
\xi_{n}
\end{array}\right]
\end{eqnarray*}
and it is not clear why the matrix $(\varphi(A_{i}^{*}B_{j}))_{i,j=1}^{n}$
should be a positive operator acting in $H\otimes\mathbb{C}^{n}$.
But we could very well put this extra requirement into an axiom, so
the CP condition (\ref{eq:cp1}). 
\begin{enumerate}
\item We only assume $\mathfrak{A}$ is a $*$-algebra, not necessarily
a $C^{*}$-algebra. $\varphi:\mathfrak{A}\rightarrow\mathscr{B}(\mathscr{H})$
is positive does not necessarily imply $\varphi$ is completely positive.
A counterexample for $\mathfrak{A}=M_{2}(\mathbb{C})$, and $\varphi:\mathfrak{A}\rightarrow B(\mathbb{C}^{2})\simeq M_{2}(\mathbb{C})$
given by taking transpose, i.e., $A\mapsto\varphi(A)=A^{tr}$. Then
$\varphi$ is positive, but $\varphi\otimes I_{M_{2}}$ is not.
\item The operator matrix $(A_{i}^{*}A_{j})$, which is also written as
$\sum_{i,j}A_{i}^{*}A_{j}\otimes e_{ij}$ is a positive element in
$\mathfrak{A}\otimes M_{n}$. All positive elements in $\mathfrak{A}\otimes M_{n}$
are in such form. This notation goes back again to Dirac, for the
rank-1 projections $\left|v\left\rangle \right\langle v\right|$ are
positive, and all positive operators are sums of these rank-1 operators. 
\item Given a CP map $\varphi:\mathfrak{A}\rightarrow\mathscr{B}(\mathscr{H})$,
we get a Hilbert space $\mathscr{K}_{\varphi}$, a representation
$\pi:\mathfrak{A}\rightarrow B(\mathscr{K}_{\varphi})$ and an isometry
$V:\mathscr{H}\rightarrow\mathscr{K}_{\varphi}$, such that 
\[
\varphi(A)=V^{*}\pi(A)V
\]
for all $A\in\mathfrak{A}$. $P=VV^{*}$ is a selfadjoint projection
from $\mathscr{K}_{\varphi}$ to the image of $\mathscr{H}$ under
the embedding. To see $P$ is a projection, note that 
\[
P^{2}=VV^{*}VV^{*}=V(V^{*}V)V^{*}=VV^{*}.
\]

\end{enumerate}
\begin{flushleft}
\textbf{Summary}
\par\end{flushleft}

Positive maps have been a recursive theme in functional analysis.
A classical example is $\mathfrak{A}=C_{c}(X)$ with a positive linear
functional $\Lambda:\mathfrak{A}\rightarrow\mathbb{C}$, mapping $\mathfrak{A}$
into a 1-d Hilbert space $\mathbb{C}$.\index{space!Hilbert-}

In Stinespring's formulation, $\varphi:\mathfrak{A}\rightarrow\mathscr{H}$
is a CP map, then we may write $\varphi(A)=V^{*}\pi(A)V$ where $\pi:\mathfrak{A}\rightarrow\mathscr{K}$
is a representation on a bigger Hilbert space $\mathscr{K}$ containing
$\mathscr{H}$. The containment is in the sense that $V:\mathscr{H}\hookrightarrow\mathscr{K}$
embeds $\mathscr{H}$ into $\mathscr{K}$. Notice that 
\[
V\varphi(A)=\pi(A)V\Longrightarrow\varphi(A)=V^{*}\pi(A)V
\]
but not the other way around. In Nelson's notes \cite{Ne69}, we use
the notation $\varphi\subset\pi$ for one representation being the
subrepresentation of another. To imitate the situation in linear algebra,
we may want to split an operator $T$ acting on $\mathscr{K}$ into
operators action on $\mathscr{H}$ and its complement in $\mathscr{K}$.
Let $P:\mathscr{K}\rightarrow\mathscr{H}$ be the orthogonal projection.
In matrix language,\index{Nelson, E}\index{Theorem!Stinespring's-}
\[
\left[\begin{array}{cc}
PTP & PTP^{\perp}\\
P^{\perp}TP & P^{\perp}TP^{\perp}
\end{array}\right].
\]
A better looking would be 
\[
\left[\begin{array}{cc}
PTP & 0\\
0 & P^{\perp}TP^{\perp}
\end{array}\right]=\left[\begin{array}{cc}
\varphi_{1} & 0\\
0 & \varphi_{2}
\end{array}\right]
\]
hence 
\[
\pi=\varphi_{1}\oplus\varphi_{2}.
\]
Stinespring's theorem is more general, where the off-diagonal entries
may not be zero.
\begin{xca}[Tensor with $M_{n}$]
\myexercise{Tensor with $M_{n}$}\label{exer:cp1}Let $\mathfrak{A}=\mathscr{B}(\mathscr{H})$,
$\xi_{1},\ldots,\xi_{n}\in\mathscr{H}$. The map 
\[
(\xi_{1},\ldots,\xi_{n})\mapsto(A\xi_{1},\ldots,A\xi_{n})\in\oplus^{n}\mathscr{H}
\]
is a representation of $\mathfrak{A}$ if and only if 
\[
\overset{n\mbox{ times}}{\overbrace{id_{\mathscr{H}}\oplus\cdots\oplus id_{\mathscr{H}}}}\in Rep(\mathfrak{A},\overset{n\mbox{ times}}{\overbrace{\mathscr{H}\oplus\cdots\oplus\mathscr{H}}})
\]
where in matrix notation, we have
\begin{eqnarray*}
 &  & \left[\begin{array}{cccc}
id_{\mathfrak{A}}(A) & 0 & \cdots & 0\\
0 & id_{\mathfrak{A}}(A) & \cdots & 0\\
\vdots & \vdots & \vdots & \vdots\\
0 & 0 & \cdots & id_{\mathfrak{A}}(A)
\end{array}\right]\left[\begin{array}{c}
\xi_{1}\\
\xi_{2}\\
\vdots\\
\xi_{n}
\end{array}\right]\\
 & = & \left[\begin{array}{cccc}
A & 0 & \cdots & 0\\
0 & A & \cdots & 0\\
\vdots & \vdots & \vdots & \vdots\\
0 & 0 & \cdots & A
\end{array}\right]\left[\begin{array}{c}
\xi_{1}\\
\xi_{2}\\
\vdots\\
\xi_{n}
\end{array}\right]=\left[\begin{array}{c}
A\xi_{1}\\
A\xi_{2}\\
\vdots\\
A\xi_{n}
\end{array}\right].
\end{eqnarray*}
In this case, the identity representation $id_{\mathfrak{A}}:\mathfrak{A}\rightarrow\mathscr{H}$
has multiplicity $n$.
\end{xca}

\begin{xca}[Column operators]
\myexercise{Column operators}\label{exer:cp2}Let $V_{i}:\mathscr{H}\rightarrow\mathscr{H}$,
and 
\begin{equation}
V:=\left[\begin{array}{c}
V_{1}\\
V_{2}\\
\vdots\\
V_{n}
\end{array}\right]:\mathscr{H}\rightarrow\oplus_{1}^{n}\mathscr{H}.\label{eq:co2}
\end{equation}
Let $V^{*}:\oplus_{1}^{n}\mathscr{H}\rightarrow\mathscr{H}$ be the
adjoint of $V$. Prove that $V^{*}=\left[\begin{array}{cccc}
V_{1}^{*} & V_{2}^{*} & \cdots & V_{n}^{*}\end{array}\right]$. \index{operators!adjoint-}\end{xca}
\begin{proof}
Let $\xi\in\mathscr{H}$, then $V\xi=\left[\begin{array}{c}
V_{1}\xi\\
V_{2}\xi\\
\vdots\\
V_{n}\xi
\end{array}\right]$, and 
\begin{eqnarray*}
\left\langle \left[\begin{array}{c}
\eta_{1}\\
\eta_{2}\\
\vdots\\
\eta_{n}
\end{array}\right],\left[\begin{array}{c}
V_{1}\xi\\
V_{2}\xi\\
\vdots\\
V_{n}\xi
\end{array}\right]\right\rangle  & = & \sum_{i}\left\langle \eta_{i},V_{i}\xi\right\rangle \\
 & = & \sum_{i}\left\langle V_{i}^{*}\eta_{i},\xi\right\rangle =\left\langle \left[\begin{array}{cccc}
V_{1}^{*} & V_{2}^{*} & \cdots & V_{n}\end{array}^{*}\right]\left[\begin{array}{c}
\eta_{1}\\
\eta_{2}\\
\vdots\\
\eta_{n}
\end{array}\right],\xi\right\rangle .
\end{eqnarray*}
This shows that $V^{*}=\left[\begin{array}{cccc}
V_{1}^{*} & V_{2}^{*} & \cdots & V_{n}^{*}\end{array}\right]$.\end{proof}
\begin{xca}[Row-isometry]
\myexercise{Row-isometry}\label{exer:cp3}Let $V$ be as in (\ref{eq:co2}).
The following are equivalent: \index{row-isometry}
\begin{enumerate}
\item $V$ is an isometry, i.e., $\left\Vert V\xi\right\Vert ^{2}=\left\Vert \xi\right\Vert ^{2}$,
for all $\xi\in\mathscr{H}$;
\item $\sum V_{i}^{*}V_{i}=I_{\mathscr{H}}$;
\item $V^{*}V=I_{\mathscr{H}}$.
\end{enumerate}
\end{xca}
\begin{proof}
Notice that \index{isometry} 
\[
\left\Vert V\xi\right\Vert ^{2}=\sum_{i}\left\Vert V_{i}\xi\right\Vert ^{2}=\sum_{i}\left\langle \xi,V_{i}^{*}V_{i}\xi\right\rangle =\left\langle \xi,\sum_{i}V_{i}^{*}V_{i}\xi\right\rangle .
\]
Hence $\left\Vert V\xi\right\Vert ^{2}=\left\Vert \xi\right\Vert ^{2}$
if and only if 
\[
\left\langle \xi,\sum_{i}V_{i}^{*}V_{i}\xi\right\rangle =\left\langle \xi,\xi\right\rangle 
\]
for all $\xi\in\mathscr{H}$. Equivalently, $\sum_{i}V_{i}^{*}V_{i}=I_{H}=V^{*}V$.\end{proof}
\begin{cor}[Krauss]
 Let $dim\mathscr{H}=n$. Then all the CP maps are of the form
\[
\varphi(A)=\sum_{i}V_{i}^{*}AV_{i}.
\]
(The essential part here is that for any CP mapping $\varphi$, we
get a system $\left\{ V_{i}\right\} $.)
\end{cor}
This was discovered in the physics literature by Kraus. The original
proof was very intricate, but it is a corollary of Stinespring's theorem.
When $dim\mathscr{H}=n$, let $\left\{ e_{1},\ldots e_{n}\right\} $
be an ONB. Fix a CP map $\varphi$, and get $\left(V,\mathscr{K},\pi\right)$.
Set 
\[
V_{i}:e_{i}\mapsto Ve_{i}\in\mathscr{K},\;i=1,\ldots n;
\]
then $V_{i}$ is an isometry. So we get a system of isometries, and
\[
\varphi(A)=\left[\begin{array}{cccc}
V_{1}^{*} & V_{2}^{*} & \cdots & V_{n}^{*}\end{array}\right]\left[\begin{array}{cccc}
A\\
 & A\\
 &  & \ddots\\
 &  &  & A
\end{array}\right]\left[\begin{array}{c}
V_{1}\\
V_{2}\\
\vdots\\
V_{n}
\end{array}\right].
\]
Notice that $\varphi(1)=1$ if and only if $\sum_{i}V_{i}^{*}V_{i}=1$.

\index{product!tensor-}
\begin{xca}[Tensor products]
\label{exer:cp4}\myexercise{Tensor products}Prove the following.
\begin{enumerate}
\item $\oplus_{1}^{n}\mathscr{H}\simeq\mathscr{H}\otimes\mathbb{C}^{n}$
\item $\sum_{1}^{\oplus\infty}\mathscr{H}\simeq\mathscr{H}\otimes l^{2}$
\item Given $L^{2}(X,\mathfrak{M},\mu)$, then $L^{2}(X,\mathscr{H})\simeq\mathscr{H}\otimes L^{2}(\mu)$;
where $L^{2}(X,\mathscr{H})$ consists of all measurable functions
$f:X\rightarrow\mathscr{H}$ such that
\[
\int_{X}\left\Vert f(x)\right\Vert _{\mathscr{H}}^{2}d\mu(x)<\infty
\]
and
\[
\left\langle f,g\right\rangle =\int_{X}\left\langle f\left(x\right),g\left(x\right)\right\rangle _{\mathscr{H}}d\mu\left(x\right).
\]

\item All the spaces above are Hilbert spaces.
\end{enumerate}
\end{xca}

\begin{xca}[Using tensor product in representations]
\label{exer:cp5}\myexercise{Using tensor product in representations}
Let $(X_{i},\mathfrak{M}_{i},\mu_{i})$, $i=1,2$, be measure spaces.
Let $\pi_{i}:L^{\infty}(\mu_{i})\rightarrow L^{2}(\mu_{i})$ be the
representation such that $\pi_{i}(f)$ is the operator of multiplication
by $f$ on $L^{2}(\mu_{i})$. Hence $\pi_{i}\in Rep(L^{\infty}(X_{i}),L^{2}(\mu_{i}))$,
and 
\begin{gather*}
\pi_{1}\otimes\pi_{2}\in Rep(L^{\infty}(X_{1}\times X_{2}),L^{2}(\mu_{1}\times\mu_{2})),\;\\
\pi_{1}\otimes\pi_{2}(\tilde{\varphi})\tilde{f}=\tilde{\varphi}\tilde{f}
\end{gather*}
for all $\tilde{\varphi}\in L^{\infty}(X_{1}\times X_{2})$, and all
$\tilde{f}\in L^{2}(\mu_{1}\times\mu_{2})$.

Elementary tensors: Special form
\begin{eqnarray*}
\widetilde{\varphi}\left(x_{1},x_{2}\right) & = & \varphi_{1}\left(x_{1}\right)\varphi_{2}\left(x_{2}\right),\\
\widetilde{f}\left(x_{1},x_{2}\right) & = & f_{1}\left(x_{1}\right)f_{2}\left(x_{2}\right),\\
\left(\pi_{1}\otimes\pi_{2}\right)\left(\widetilde{\varphi}\right)f & = & \pi_{1}\left(\varphi_{1}\right)f_{1}\otimes\pi_{2}\left(\varphi_{2}\right)f_{2}.
\end{eqnarray*}

\end{xca}

\begin{xca}[Transpose is not completely positive]
\myexercise{Transpose is not completely positive}~
\begin{enumerate}
\item Let $\mathfrak{A}$ be an abelian $C^{*}$-algebra; and let $\varphi:\mathfrak{A}\rightarrow\mathscr{B}\left(\mathscr{H}\right)$
be a positive mapping; then show that $\varphi$ is in fact automatically
completely positive.
\item Show that there are positive mappings which are \uline{not} completely
positive. \\
\uline{Hint}: Let $M_{n}$ be the $n\times n$ complex matrices,
and set 
\[
\varphi\left(A\right)=A^{T},\;A\in M_{n}
\]
where $A^{T}$ is the transpose matrix. If $n>1$, show that $M_{n}\xrightarrow{\varphi}M_{n}$
is positive but \uline{not} completely positive.
\end{enumerate}
\end{xca}

\section{\label{sec:end}Endomorphisms, Representations of $\mathscr{O}_{N}$,
and Numerical Range}

Let $\mathscr{H}$ be a Hilbert space, and consider endomorphisms
in $\mathscr{B}\left(\mathscr{H}\right)$, i.e., $\sigma:\mathscr{B}\left(\mathscr{H}\right)\longrightarrow\mathscr{B}\left(\mathscr{H}\right)$,
linear, and and satisfy \index{numerical range}\index{endomorphism}
\begin{eqnarray*}
\sigma\left(AB\right) & = & \sigma\left(A\right)\sigma\left(B\right)\\
\sigma\left(A^{*}\right) & = & \sigma\left(A\right)^{*},\;\forall A,B\in\mathscr{B}\left(\mathscr{H}\right),\;\mbox{and}\\
\sigma\left(I\right) & = & I.
\end{eqnarray*}

\begin{figure}
\[
\xymatrix{ &  &  & \:\ar[rr] &  & \,\ar[rd]^{S_{1}}\\
\ar[rr]^{\text{signal in}} &  & \:\ar[ru]^{S_{1}^{*}}\ar[r]_{\vdots}^{S_{2}^{*}}\ar[rd]_{S_{N}^{*}} & \:\ar[rr] &  & \,\ar[r]_{\vdots}^{S_{2}} & \bigoplus\ar[rr]^{\text{signal out}} &  & \,\\
 &  &  & \:\ar[rr] &  & \,\ar[ru]_{S_{N}}
}
\]

\protect\caption{\label{fig:cuntz}Orthogonal bands in filter bank, ``in'' = ``out''.
An application of representations of the Cuntz relations (\ref{eq:cuntz}).
\index{filter bank} \index{frequency band}\index{Cuntz relations}}
\end{figure}
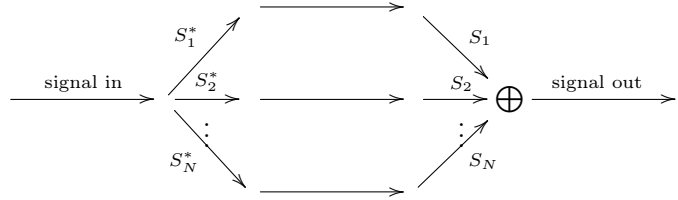

\begin{defn}
\label{def:cuntz}By a \emph{representation} $\pi$ of $\mathscr{O}_{N}$
in $\mathscr{H}$, $\pi\in Rep\left(\mathscr{O}_{N},\mathscr{H}\right)$,
we mean a system of isometries $\left(S_{i}\right)_{i=1}^{N}$ in
$\mathscr{H}$ such that
\begin{equation}
\begin{Bmatrix}S_{i}^{*}S_{j}=\delta_{ij}I\\
\sum_{i}S_{i}S_{i}^{*}=I
\end{Bmatrix}\;\text{(Cuntz relations)}\label{eq:cuntz}
\end{equation}
See  \figref{cuntz}.
\end{defn}

\begin{figure}
\includegraphics[width=0.4\textwidth]{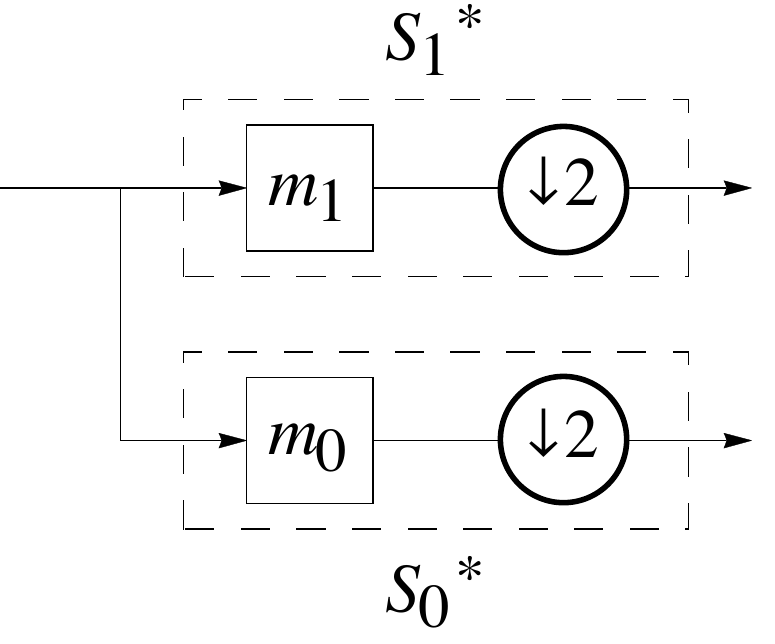}

\protect\caption{\label{fig:hbank}The Cuntz operators $S_{0}$ and $S_{1}$ in the
Haar wavelet.}
\end{figure}

\begin{figure}
\includegraphics[scale=0.4]{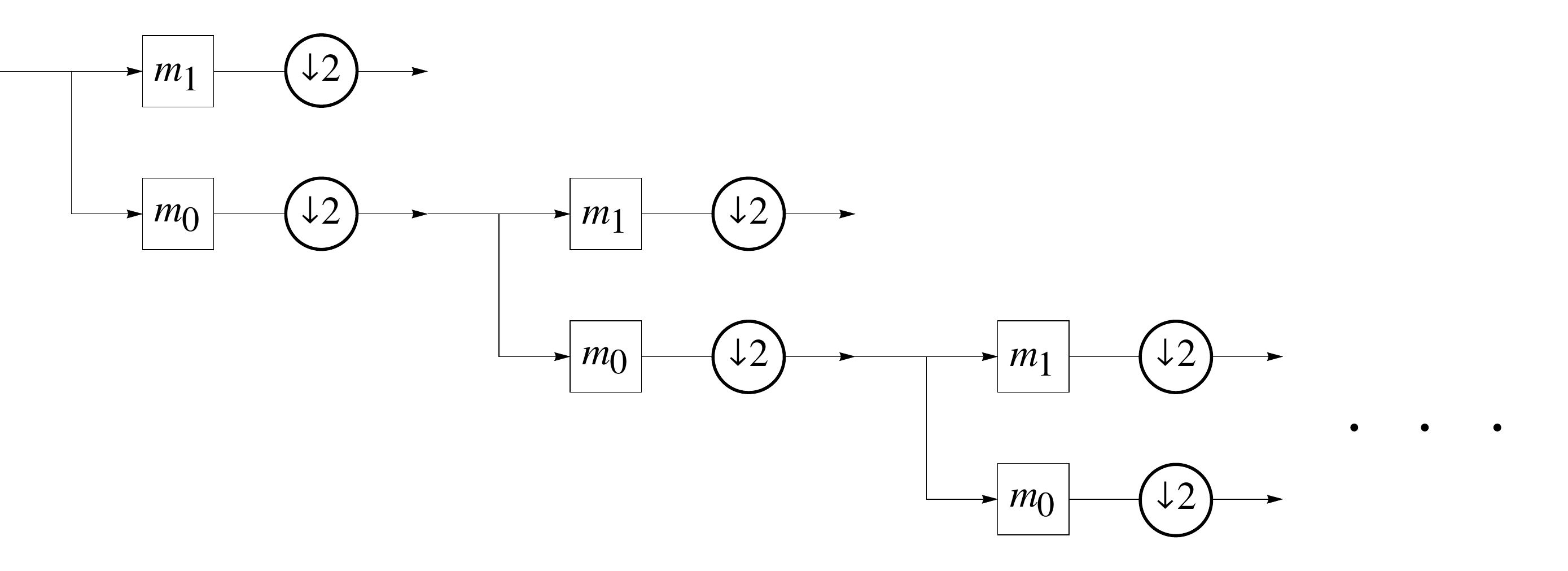}

\protect\caption{\label{fig:fba}The two-channel analysis filter bank.}
\end{figure}

\begin{figure}
\includegraphics[scale=0.4]{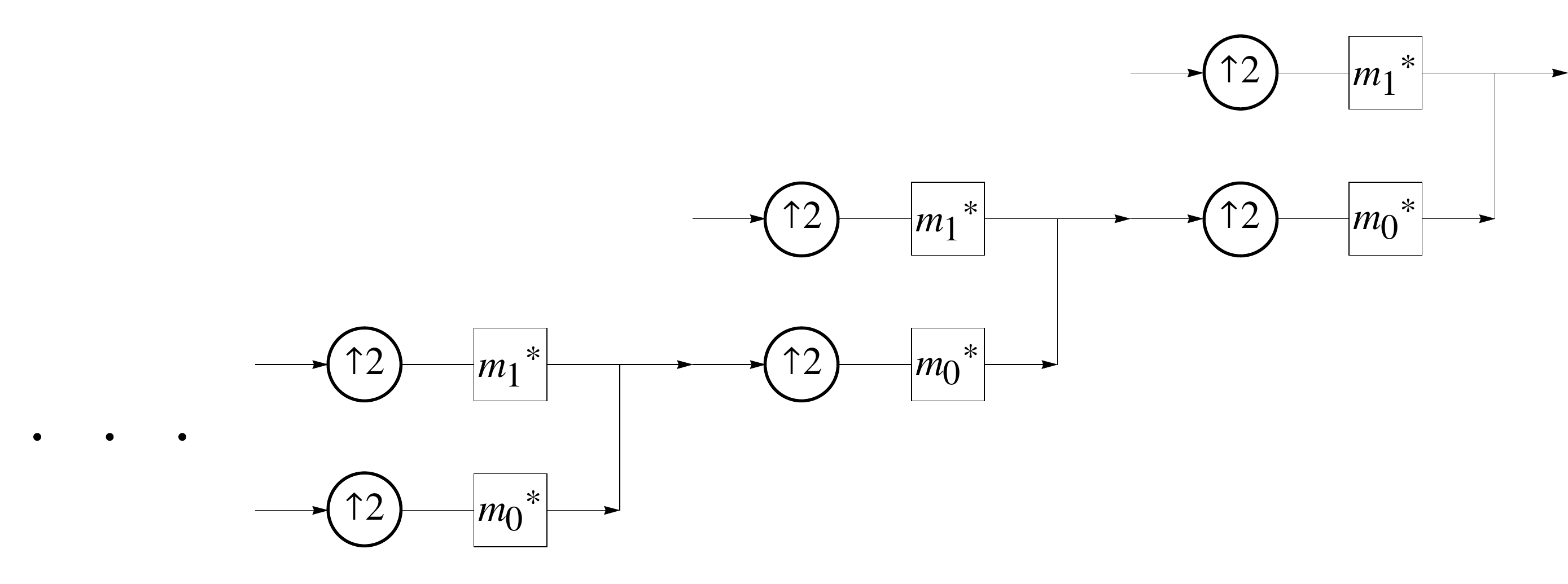}

\protect\caption{\label{fig:fbs}The two-channel synthesis filter bank.}
\end{figure}

\begin{figure}
\begin{tabular}{cc}
\includegraphics[width=0.42\textwidth]{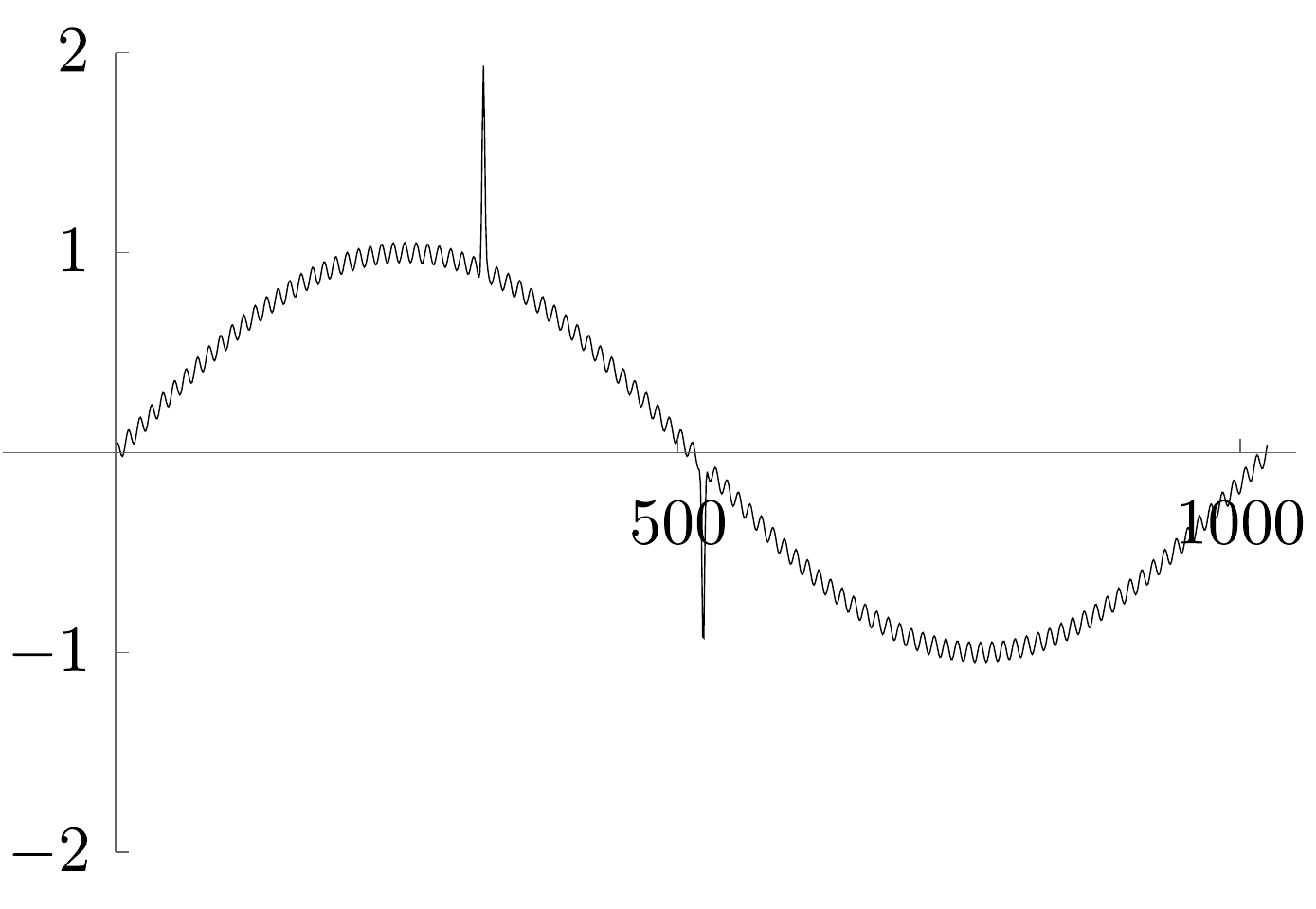} & \includegraphics[width=0.42\textwidth]{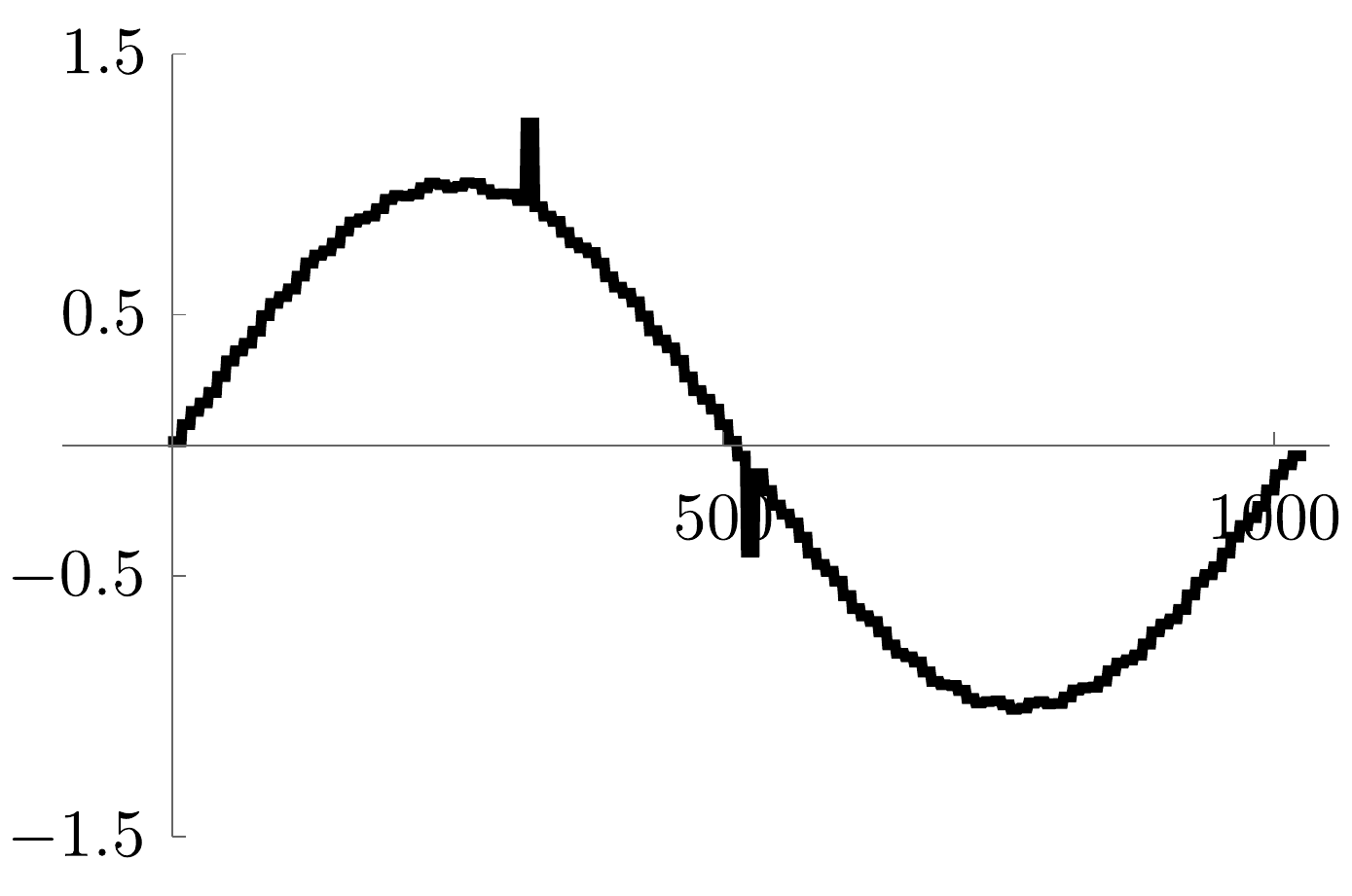}\tabularnewline
(a) $x\left(n\right)$ containing sharp noise\vspace{2em} & (b) $\left(S_{0}^{*}\right)^{3}x\left(n\right)$\tabularnewline
\includegraphics[width=0.42\textwidth]{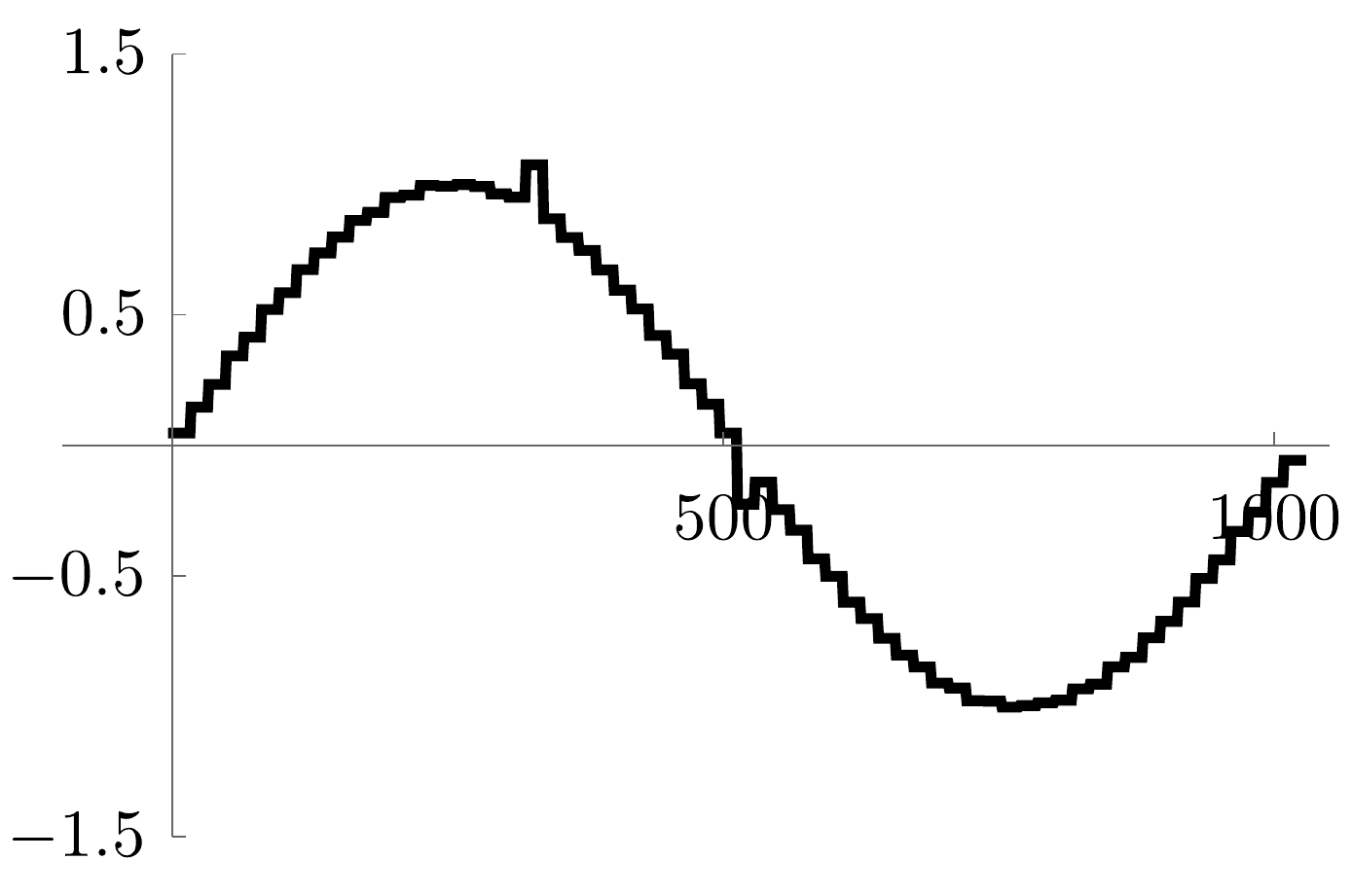} & \includegraphics[width=0.42\textwidth]{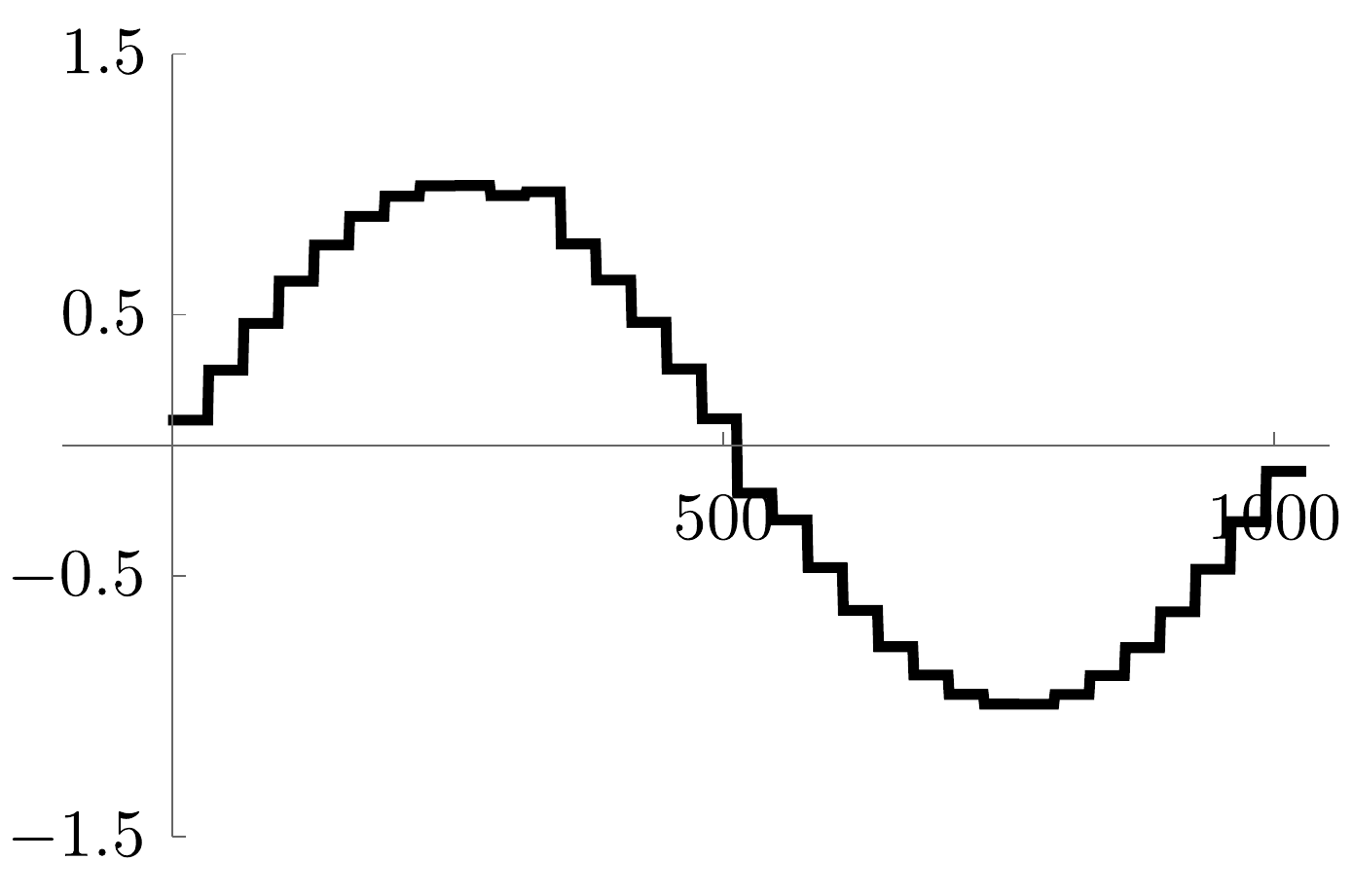}\tabularnewline
(c) $\left(S_{0}^{*}\right)^{4}x\left(n\right)$ & (d) $\left(S_{0}^{*}\right)^{5}x\left(n\right)$\tabularnewline
\end{tabular}

\protect\caption{\label{fig:dn} The outputs of $\left(S_{0}^{*}\right)^{n}$, $n=3,4,5$. }
\end{figure}

\begin{figure}
\begin{tabular}{cc}
\includegraphics[width=0.42\textwidth]{dn1} & \includegraphics[width=0.42\textwidth]{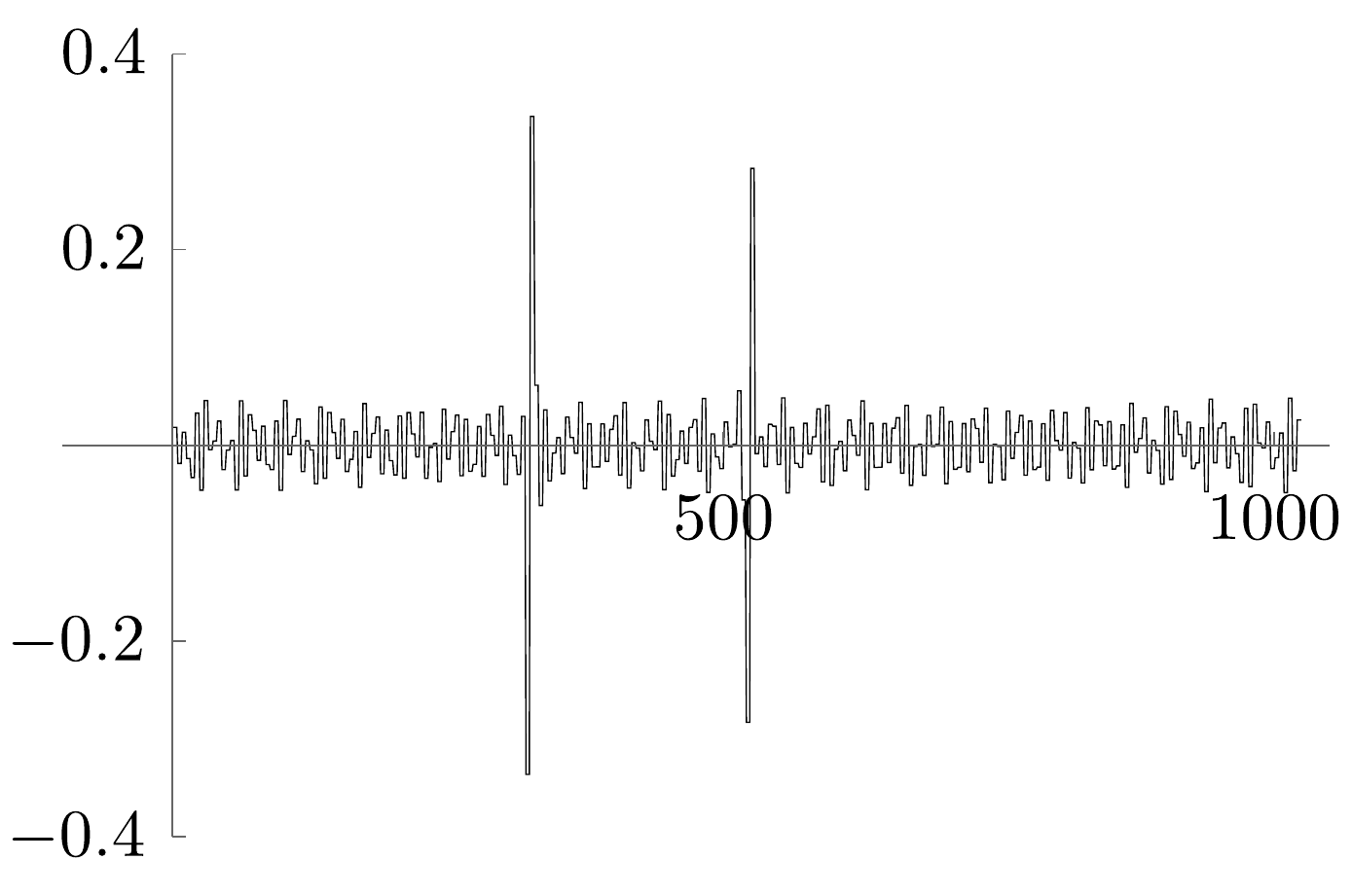}\tabularnewline
(a) $x\left(n\right)$ containing sharp noise\vspace{2em} & (b) $S_{1}^{*}\left(S_{0}^{*}\right)^{2}x\left(n\right)$\tabularnewline
\includegraphics[width=0.42\textwidth]{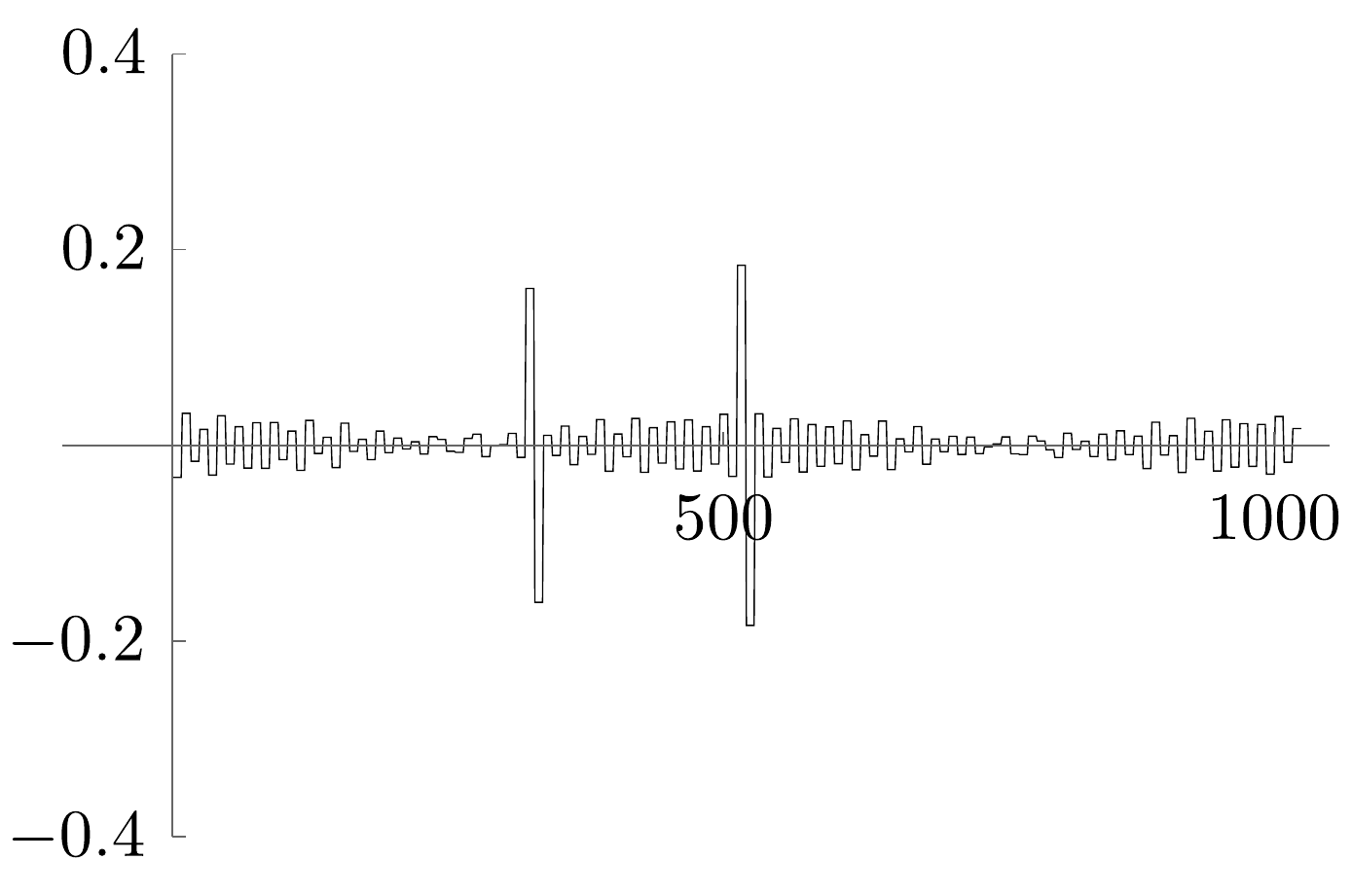} & \includegraphics[width=0.42\textwidth]{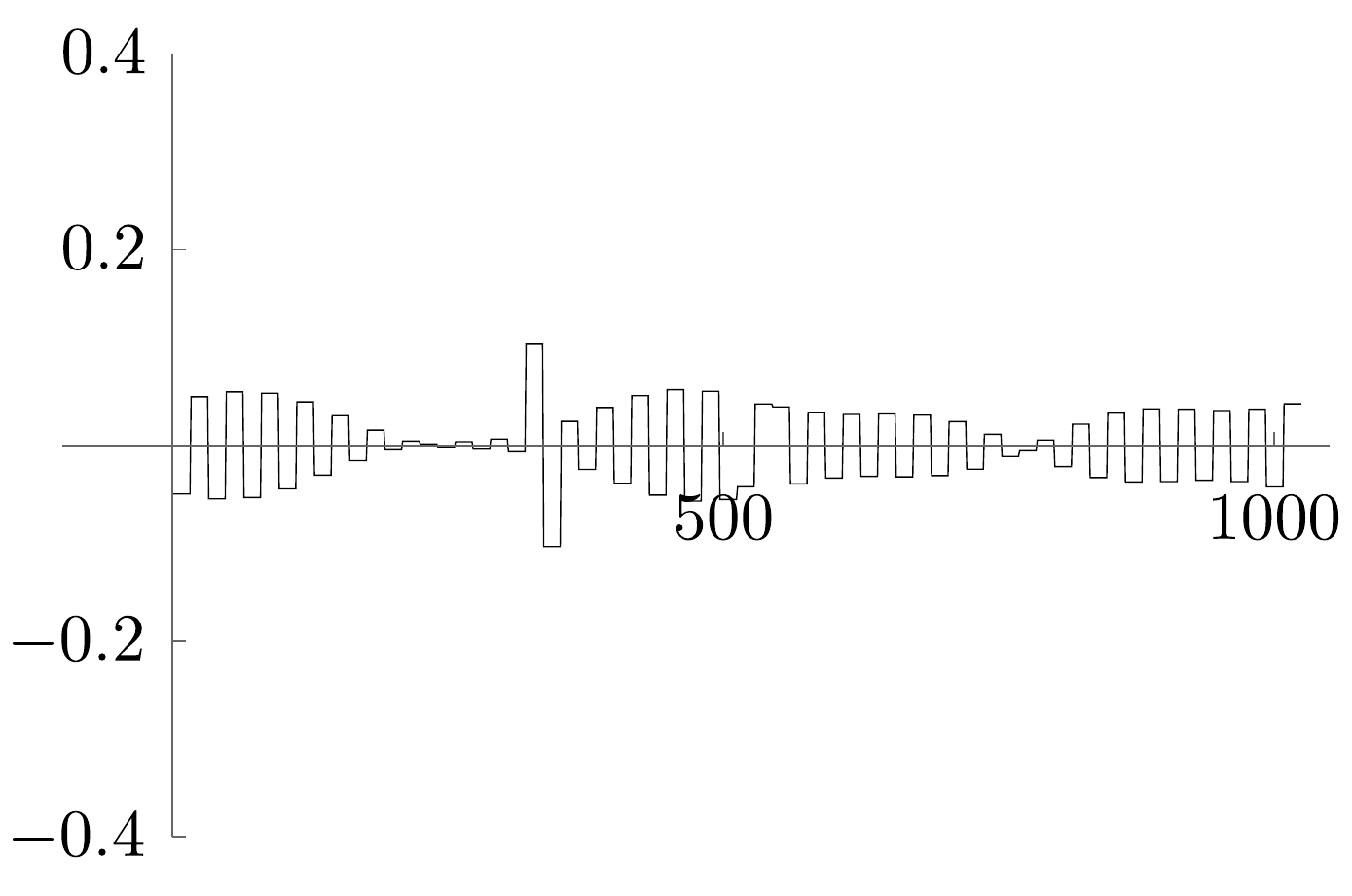}\tabularnewline
(c) $S_{1}^{*}\left(S_{0}^{*}\right)^{3}x\left(n\right)$ & (d) $S_{1}^{*}\left(S_{0}^{*}\right)^{4}x\left(n\right)$\tabularnewline
\end{tabular}

\protect\caption{\label{fig:dn2} The outputs of high-pass filters.}
\end{figure}

\begin{figure}
\noindent \begin{centering}
\includegraphics[width=0.4\textwidth]{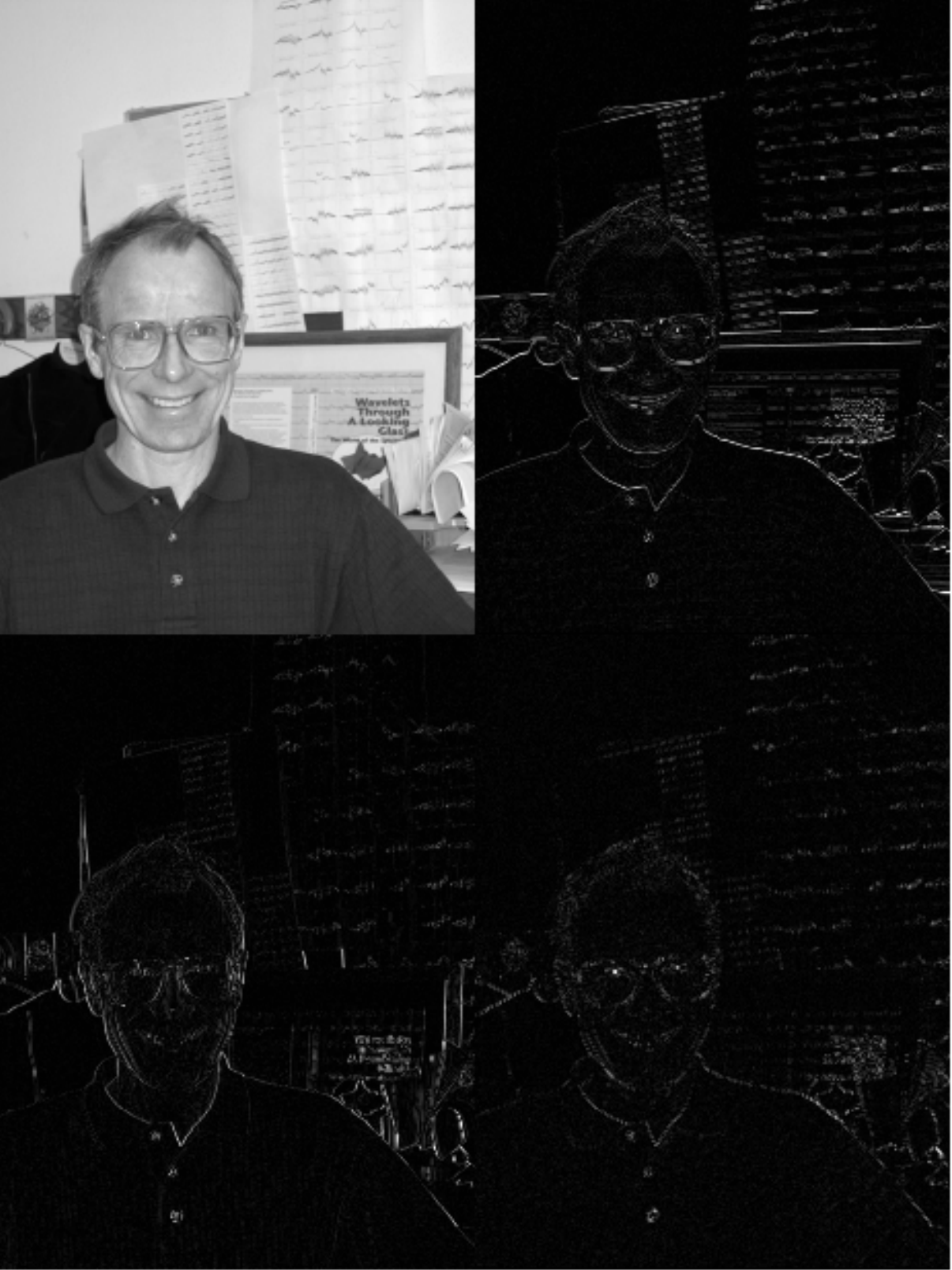}
\par\end{centering}

\protect\caption{\label{fig:2d}A coarser resolution in three directions in the plane,
filtering in directions, $x,y$, and diagonal; -- corresponding dyadic
scaling in each coordinate direction. (Image cited from M.-S. Song,
``\emph{Wavelet Image Compression}'' in \cite{JLH06}.) \index{resolution!multi-}}
\end{figure}

\begin{rem}
While the relations in (\ref{eq:cuntz}), called the Cuntz relations,
of \defref{cuntz} are axioms, they have implications for a host of
applications, and \figref{cuntz} is a graphic representations of
(\ref{eq:cuntz}) stated in a form popular in applications to signal\index{signal}
processing. Effective transmission of signals (speech, or images),
is possible because the transmitted signals can be divided into frequency
sub-bands; \index{frequency band} this is done with filters. A low-pass
filter picks out the band corresponding to frequencies in a \textquotedblleft band\textquotedblright{}
around zero, and similarly with intermediate, and high bands. The
horizontal lines in \figref{cuntz} represent prescribed bands. The
orthogonality part of (\ref{eq:cuntz}) represents non-interference
from one band to the next. Adding the projections on the LHS in (\ref{eq:cuntz})
to recover the identity operator reflects perfect reconstruction,
i.e, signal out equals signal in. The projections on the LHS in (\ref{eq:cuntz})
are projections onto subspaces of a total Hilbert space (of signals
to be transmitted), the subspaces thus representing frequency bands.\index{representation!- of the Cuntz algebra}

Thus \figref{cuntz} represents such a filter design; there are many
such, some good some not. Each one is called a \textquotedblleft filter
bank.\textquotedblright{} And each one corresponds to a representation
of (\ref{eq:cuntz}), or equivalently a representation of the Cuntz
algebra $\mathscr{O}_{N}$ where $N$ is the number of band for the
particular filter design. \index{filter bank}\end{rem}
\begin{xca}[$Rep(\mathscr{O}_{N},\mathscr{H})$]
\myexercise{$Rep(\mathscr{O}_{N},\mathscr{H})$}Fix $N\geq2$, and
let $\mathscr{O}_{N}$ denote the Cuntz-$C^{*}$-algebra (see  \exaref{cuntzN}).
By
\begin{equation}
\pi\in Rep\left(\mathscr{O}_{N},\mathscr{H}\right)\label{eq:cn1}
\end{equation}
we mean a homomorphism
\[
\pi:\mathscr{O}_{N}\longrightarrow\mathscr{B}\left(\mathscr{H}\right),
\]
(in particular, satisfying: $\pi\left(AB\right)=\pi\left(A\right)\pi\left(B\right)$,
$\pi(A^{*})=\pi\left(A\right)^{*}$, $\forall A,B\in\mathscr{O}_{N}$,
and $\pi\left(\mathbf{1}\right)=I_{\mathscr{H}}$.)

Let $\{S_{i}\}_{i=1}^{N}$ be a system of isometries in a Hilbert
space $\mathscr{H}$ satisfying (\ref{eq:cuntz}), called Cuntz-isometries.
For all multi-indices $J=\left(j_{1},j_{2},\cdots,j_{m}\right)$,
$j_{i}\in\left\{ 1,2,\cdots,N\right\} $, set
\begin{eqnarray*}
s_{J} & := & s_{j_{1}}s_{j_{2}}\cdots s_{j_{m}},\;\mbox{and}\\
S_{J} & := & S_{j_{1}}S_{j_{2}}\cdots S_{j_{m}}.
\end{eqnarray*}

Show that, given a (\ref{eq:cuntz})-system $\{S_{i}\}_{i=1}^{N}$
of isometries, there is then a unique $\pi\in Rep(\mathscr{O}_{N},\mathscr{H})$
such that
\begin{equation}
\pi\left(s_{J}s_{K}^{*}\right)=S_{J}S_{K}^{*}\label{eq:cn2}
\end{equation}
holds for all multi-indices $J$, $K$.\end{xca}
\begin{example}[Representation of the Cuntz algebra $\mathscr{O}_{2}$]
. Let $\mathscr{H}=L^{2}\left(\mathbb{T}\right)$. In signal processing
language $\mathscr{H}$ is the $L^{2}$-space of frequency functions.
Set \index{signal processing}\index{representation!- of the Cuntz algebra}
\end{example}
\begin{eqnarray}
\left(S_{0}f\right)\left(x\right) & := & \cos\left(x\right)f\left(2x\right)\label{eq:cu1}\\
\left(S_{1}f\right)\left(x\right) & := & \sin\left(x\right)f\left(2x\right)\label{eq:cu2}
\end{eqnarray}
as the two Cuntz operators, where $f\in\mathscr{H}$, and 
\[
2x:=2x\;\mod\;2\pi\mathbb{Z}\left(=\mbox{multiples of \ensuremath{2\pi}}\right).
\]
(The generators $S_{i}$, $i=1,2$ with up-sampling, and $S_{i}^{*}$
with down-sampling.) \index{sampling!up-sampling}\index{sampling!down-sampling}
\index{up-sampling}\index{down-sampling}
\begin{xca}[Simplest Low/High filter bank]
\myexercise{Simplest Low/High filter bank} Show that (\ref{eq:cu1})-(\ref{eq:cu2})
satisfy the $\mathscr{O}_{2}$-Cuntz relations, i.e., \index{Cuntz relations}
\begin{enumerate}
\item $S_{i}$, $i=0,1$ are isometries in $\mathscr{H}$;
\item $S_{0}^{*}S_{1}=0$ (orthogonality);
\item $S_{0}S_{0}^{*}+S_{1}S_{1}^{*}=I_{\mathscr{H}}$. 
\end{enumerate}

\uline{Hint}: First show that 
\[
\left(S_{0}^{*}f\right)\left(x\right)=\frac{1}{2}\left(\cos\left(\frac{x}{2}\right)f\left(\frac{x}{2}\right)+\cos\left(\frac{x+\pi}{2}\right)f\left(\frac{x+\pi}{2}\right)\right)
\]
are similarly for $S_{1}^{*}f$. Then compute directly that 
\[
\left\Vert S_{0}^{*}f\right\Vert _{\mathscr{H}}^{2}+\left\Vert S_{1}^{*}f\right\Vert _{\mathscr{H}}^{2}=\left\Vert f\right\Vert _{\mathscr{H}}^{2}.
\]

\end{xca}
\begin{example}
Consider the Haar wavelet as in  \exaref{haar}, with $\phi_{0}$
(scaling function), $\varphi_{1}$ and $\psi_{j,k}$, $j,k\in\mathbb{Z}$
be as in (\ref{eq:hw1})-(\ref{eq:hw2}). Set 
\[
h\left(n\right)=\begin{cases}
-\frac{1}{2} & n=-1\\
\frac{1}{2} & n=0\\
0 & \mbox{otherwise}
\end{cases},\quad g\left(n\right)=\begin{cases}
\frac{1}{2} & n=-1\\
\frac{1}{2} & n=0\\
0 & \mbox{otherwise}
\end{cases};
\]
so that $h,g\in l^{2}$; where $g$ is the low-pass filter (averaging
data), and $h$ is the high-pass filter (capturing high-frequency
oscillations). Let $m_{0}$ and $m_{1}$ be Fourier transform of $g$
and $h$ respectively, i.e., 
\begin{eqnarray*}
m_{0}\left(x\right) & = & \sum_{n\in\mathbb{Z}}g\left(n\right)e^{-ixn}\\
m_{1}\left(x\right) & = & \sum_{n\in\mathbb{Z}}h\left(n\right)e^{-ixn}
\end{eqnarray*}
and $m_{0},m_{1}\in L^{2}\left(\mathbb{T}\right)$. Finally set $S_{0}^{*}$
and $S_{1}^{*}$ as in \figref{hbank}. 

An input signal goes through the analysis filter bank (\figref{fba})
and splits into layers (frequency bands) of fine details. The original
signal can be rebuilt through the synthesis filter bank (\figref{fbs}),
i.e., a perfect reconstruction. 

Depending on the applications, the output of the analysis filter bank
will go through other DSP device. For example, in data compression,
insignificant coefficients are dropped; or if the task is to remove
noise in the input signal, the coefficients corresponding to high
frequency components (noise) are removed, and the remaining coefficients
go through the synthesis filter bank. See Figures \ref{fig:dn}-\ref{fig:dn2}
for an illustration, and \figref{2d} for an application in imaging
processing.

We return to a much more detailed discussion of down-sampling and
up-sampling in \chapref{bm}.\index{up-sampling}\index{down-sampling}
\end{example}

\begin{figure}
\begin{minipage}[t]{1\columnwidth}%
\subfloat[Low: $\cos^{2}\left(\frac{x}{2}\right)$]{\noindent \protect\begin{centering}
\protect\includegraphics[width=0.45\textwidth]{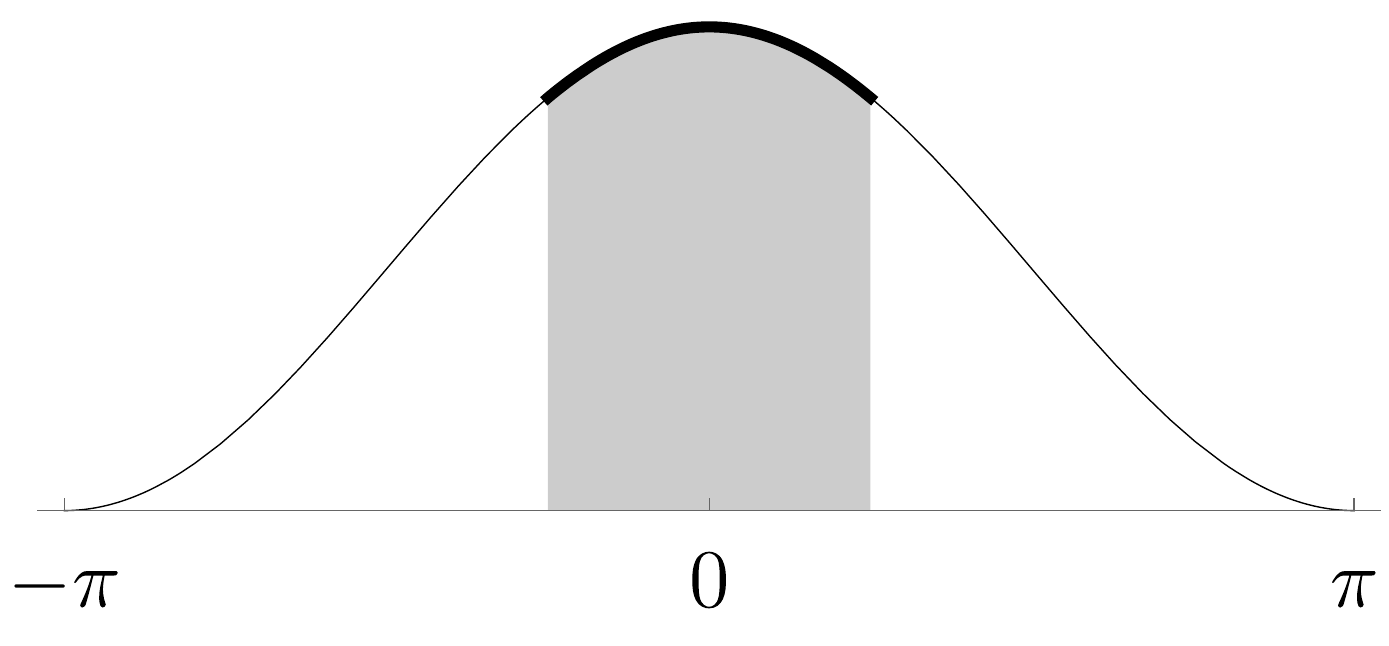}\protect
\par\end{centering}

}\hfill{}\subfloat[High: $\sin^{2}\left(\frac{x}{2}\right)$]{\noindent \protect\begin{centering}
\protect\includegraphics[width=0.45\textwidth]{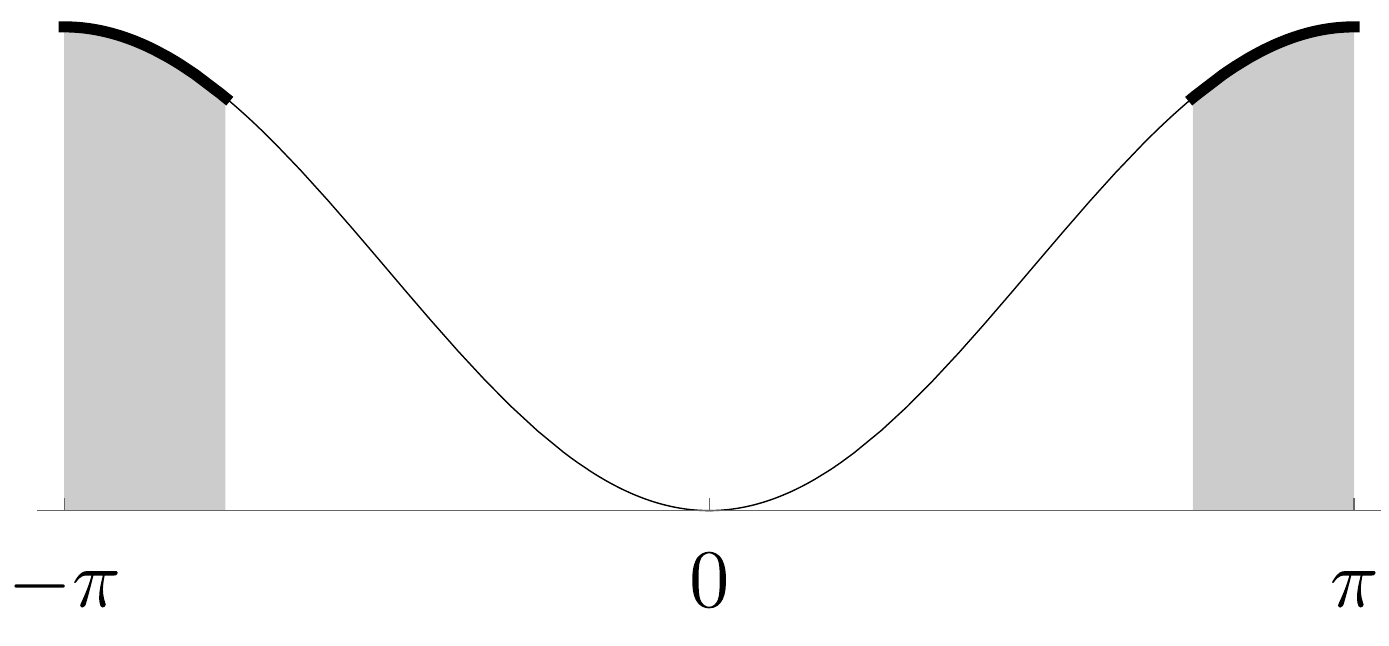}\protect
\par\end{centering}

}%
\end{minipage}

\protect\caption{Low / high pass filters for the Haar wavelet (frequency mod $2\pi$).}
\end{figure}

\begin{figure}
\includegraphics[width=0.5\textwidth]{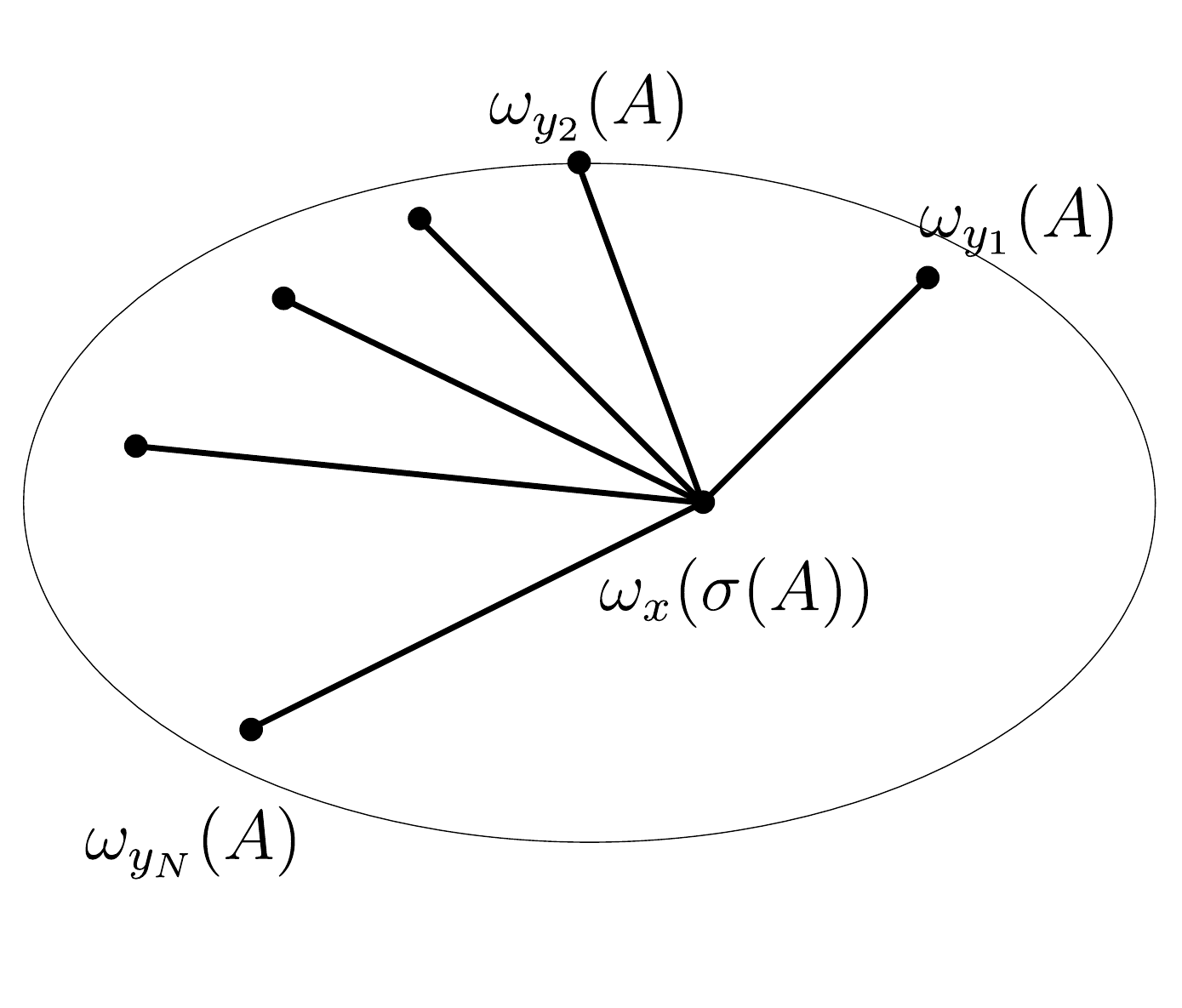}

\protect\caption{\label{fig:NR}Illustration of eq. (\ref{eq:nr}), with $y_{i}:=\frac{S_{i}^{*}x}{\left\Vert S_{i}^{*}x\right\Vert }$,
$i=1,2,\ldots,N$.}
\end{figure}

\begin{xca}[Endomorphism vs representation]
\myexercise{Endomorphism vs representation}\label{exer:cp6}Let
$\mathscr{H}$ be a general separable Hilbert space. The purpose below
is to point out that the study of $Rep(\mathscr{O}_{N},\mathscr{H})$
is essentially equivalent to that of the endomorphisms of $\mathscr{B}(\mathscr{H})$.
\begin{enumerate}
\item \label{enu:on1}Let $\sigma$ be an endomorphism in $\mathscr{B}\left(\mathscr{H}\right)$
of finite index. Show that there is a representation $\left(S_{i}\right)$
of $\mathscr{O}_{N}$ in $\mathscr{H}$ such that
\begin{equation}
\sigma\left(A\right)=\sum_{i=1}^{N}S_{i}AS_{i}^{*}.\label{eq:ed1}
\end{equation}
\emph{Notation}. Given $\sigma\in End\left(\mathscr{B}\left(\mathscr{H}\right)\right)$,
the $N$ in the corresponding representation (\ref{eq:ed1}) is called
\emph{Powers-index} of $\sigma$. It holds that for every $\sigma\in End\left(\mathscr{B}\left(\mathscr{H}\right)\right)$,
the relative commutant
\[
\mathscr{B}\left(\mathscr{H}\right)\cap\sigma\left(\mathscr{B}\left(\mathscr{H}\right)\right)'
\]
is a type $I_{N}$, and this $N$ coincides with the Powers-index.
\index{Powers-index}\index{index!Powers-}
\item Let $\sigma$, $\left\{ S_{i}\right\} $ be as in (\ref{enu:on1}),
and let $A\in\mathscr{B}\left(\mathscr{H}\right)$; then show that
\begin{equation}
NR_{\sigma\left(A\right)}\subseteq NR_{A}.
\end{equation}
In other words, endomorphisms in $\mathscr{B}\left(\mathscr{H}\right)$
contract the numerical range.
\end{enumerate}

\uline{Hint:} Use the following three facts:
\begin{enumerate}[label=(\roman{enumi})]
\item The numerical range $NR_{A}$ is convex; and 
\item if $x\in\mathscr{H}$, $\left\Vert x\right\Vert =1$, then (see \figref{NR})
\begin{equation}
w_{x}\left(\sigma\left(A\right)\right)=\sum_{i=1}^{N}\left\Vert S_{i}^{*}x\right\Vert ^{2}w_{\frac{S_{i}^{*}x}{\left\Vert S_{i}^{*}x\right\Vert }}\left(A\right);\label{eq:nr}
\end{equation}

\item and lastly, 
\begin{equation}
\sum_{i=1}^{N}\left\Vert S_{i}^{*}x\right\Vert ^{2}=1.
\end{equation}

\end{enumerate}
\end{xca}

\begin{xca}[Convex sets that are not numerical ranges]
\label{exer:cp7}\myexercise{Convex sets that are not numerical ranges}
Give an example of a bounded convex subset of the complex plane which
is not $NR_{A}$ for any $A\in\mathscr{B}\left(\mathscr{H}\right)$,
where $\mathscr{H}$ is some Hilbert space. \index{numerical range}
\end{xca}

\section*{A summary of relevant numbers from the Reference List}

For readers wishing to follow up sources, or to go in more depth with
topics above, we suggest: 

The paper \cite{Sti55} is pioneering, starting the study of completely
positive mappings in operator algebra theory. A more comprehensive
list is: \cite{Arv98,BR81,Sti59,MR548728,BJ02,MR2254502,Arv76,Fan10,BJKR84,MR0467330,MR1468230,MR0374334,Sti55,MJD15}.

\chapter{Brownian Motion\label{chap:bm}}
\begin{quotation}
From its shady beginnings devising gambling strategies and counting
corpses in medieval London, probability theory and statistical inference
now emerge as better foundations for scientific models, especially
those of the process of thinking and as essential ingredients of theoretical
mathematics, even the foundations of mathematics itself. 

--- David Mumford\sindex[nam]{Mumford, D., (1937--)}\vspace{1em}

It is intriguing that the mathematics of Brownian motion was discovered
almost simultaneously more than 100 years ago By Bachelier and by
Einstein: In physics (Albert Einstein, 1005,\textquotedbl{}Über die
von der molekularkinetischen Theorie der Wärme geforderte Bewegung
von in ruhenden Flüssigkeiten suspendierten Teilchen;\textquotedbl{}
On the Motion of Small Particles Suspended in a Stationary Liquid,
as Required by the Molecular Kinetic Theory of Heat). And in finance
(Jean-Jacques Bachelier, 1900, \textquotedbl{}The Theory of Speculation\textquotedbl{}).
Einstein\sindex[nam]{Einstein, A., (1879-1955)}

In Einstein's paper, Brownian motion offered one of the first experimental
justification for the atomic theory. The Brownian motion model for
financial markets is a continuous extension of the \textquotedbl{}one-period
market model\textquotedbl{} of H. Markowitz, (in fact much later than
1900). 

Bachelier: Continuous prices of financial asset-markets evolve in
time according to a geometric Brownian motion.\vspace{1em} \\
\textquotedblleft The glory of science is to imagine more than we
can prove.\textquotedblright{} 

--- Freeman Dyson\vspace{1em}\\
\textquotedblleft Not only does God play dice, but... he sometimes
throws them where they cannot be seen.\textquotedblright{}

--- Stephen Hawking\vspace{2em}
\end{quotation}
\sindex[nam]{Bachelier, L., (1870-1946)}

The concept of Brownian motion is not traditional included in Functional
Analysis. Below we offer a presentation which relies on almost all
the big theorems from functional analysis, and especially on $L^{2}$-Hilbert
spaces, built from probability measures on function spaces.

We have included a brief discussion of Brownian motion in order to
illustrate infinite Cartesian products (sect 1.0.1) and unitary one-parameter
group $\left\{ U\left(t\right)\right\} _{t\in\mathbb{R}}$ acting
in Hilbert space. See \cite{Jor14,Nel67,Nel59a,Hi80}.
\begin{defn}
Let $\left(\Omega,\mathcal{F},\mathbb{P}\right)$ be a \emph{probability
space}, i.e.,
\begin{itemize}
\item $\Omega=$ sample space
\item $\mathcal{F}=$ sigma algebra of events
\item $\mathbb{P}=$ probability measure defined on $\mathcal{F}$. 
\end{itemize}

A function $X:\Omega\rightarrow\mathbb{R}$ is called a \emph{random
variable }if it is a measurable function, i.e., if for all intervals
$\left(a,b\right)\subset\mathbb{R}$ the inverse image\emph{ 
\begin{equation}
X^{-1}\left(\left(a,b\right)\right)=\left\{ \omega\in\Omega\:\big|\:X\left(\omega\right)\in\left(a,b\right)\right\} \label{eq:bm1}
\end{equation}
}is in $\mathcal{F}$; and we write $\left\{ a<X\left(\omega\right)<b\right\} \in\mathcal{F}$
in short-hand notation.

We say that $X$ is \emph{Gaussian} if $\exists\:m,\sigma$ (written
$N\left(m,\sigma^{2}\right)$) such that 
\begin{equation}
\mathbb{P}\left(\left\{ a<X\left(\omega\right)<b\right\} \right)=\int_{a}^{b}\frac{1}{\sigma\sqrt{2\pi}}e^{-\left(x-m\right)^{2}/2\sigma^{2}}dx.\label{eq:bm2}
\end{equation}
The function under the integral in (\ref{eq:bm2}) is called the Gaussian
(or normal) distribution. \index{distribution!Gaussian-}\index{random variable}
\index{normal distribution}
\end{defn}

\begin{defn}
Events $A,B\in\mathcal{F}$ are said to be \emph{independent }if\emph{
\[
\mathbb{P}\left(A\cap B\right)=\mathbb{P}\left(A\right)\mathbb{P}\left(B\right).
\]
}Random variables $X$ and $Y$ are said to be independent iff (Def.)
$X^{-1}\left(I\right)$ and $Y^{-1}\left(J\right)$ are independent
for all intervals $I$ and $J$. \index{independent!-events}\index{independent!-increments}
\end{defn}

\begin{defn}
\label{def:bm}A family $\left\{ X_{t}\right\} _{t\in\mathbb{R}}$
of random variables for $\left(\Omega,\mathcal{F},\mathbb{P}\right)$
is said to be a \emph{Brownian motion} iff (Def.) for every $n\in\mathbb{N}$, 
\begin{enumerate}
\item the random variables $X_{t_{1}},X_{t_{2}},\ldots,X_{t_{n}}$ are jointly
Gaussian with 
\[
\mathbb{E}\left(X_{t}\right)=\int_{\Omega}X_{t}\left(\omega\right)d\mathbb{P}\left(\omega\right)=0,\;\forall t\in\mathbb{R};
\]

\item if $t_{1}<t_{2}<\cdots<t_{n}$ then $X_{t_{i+1}}-X_{t_{i}}$ and $X_{t_{i}}-X_{t_{i-1}}$
are independent;
\item \label{enu:bm3}for all $s,t\in\mathbb{R}$,\index{Gaussian process}
\[
\mathbb{E}\left(\left|X_{t}-X_{s}\right|^{2}\right)=\left|t-s\right|.
\]

\end{enumerate}
\end{defn}
\begin{rem}
\begin{flushleft}
There is a list of popular kernels in probability theory (\tabref{pk}).
See any book on probability theory. 
\par\end{flushleft}
\end{rem}
\renewcommand{\arraystretch}{3}

\begin{table}
\begin{tabular}{|>{\centering}p{0.2\textwidth}|>{\centering}p{0.3\textwidth}|>{\centering}m{0.4\textwidth}|}
\hline 
uniform & $a\leq x\leq b$ & $\dfrac{1}{b-1}$\tabularnewline
\hline 
exponential $\left(\lambda\right)$ & $x\geq0$ & $\lambda e^{-\lambda x}$\tabularnewline
\hline 
Gaussian normal $N\left(m,\sigma^{2}\right)$ & $x\in\mathbb{R}$ & $\frac{1}{\sigma\sqrt{2\pi}}e^{-\frac{1}{2}\left(\frac{x-m}{\sigma}\right)^{2}}$\tabularnewline
\hline 
Cauchy & $x\in\mathbb{R}$ & $\dfrac{1}{\pi\left(1+x^{2}\right)}$\tabularnewline
\hline 
$\chi^{2}$ (chi-square) & $x\geq0$ & $\dfrac{e^{-\frac{x}{2}}x^{\frac{\nu}{2}-1}}{2^{\frac{\nu}{2}}\Gamma\left(\frac{\nu}{2}\right)}$\tabularnewline
\hline 
Gamma & $x\geq0$ & $\dfrac{x^{\gamma-1}e^{-x}}{\Gamma\left(\gamma\right)},\;\gamma>0$\tabularnewline
\hline 
\end{tabular}

\protect\caption{\label{tab:pk} Probability kernels (distributions)}
\end{table}

\renewcommand{\arraystretch}{1}
\begin{xca}[Quadratic variation]
\label{exer:qv}\myexercise{Quadratic variation} Let $\left\{ X_{t}\right\} $
be the Brownian motion (\defref{bm}). Fix $T\in\mathbb{R}_{+}$,
and consider partitions $\pi:(t_{i})_{i=0}^{n}$ of $\left[0,T\right]$,
i.e., \index{quadratic variation} 
\begin{equation}
\pi\::\:0=t_{0}<t_{1}<t_{2}<\cdots<t_{n}=T.\label{eq:qv1}
\end{equation}
Set 
\begin{equation}
\mbox{mesh}\left(\pi\right)\big(:=\left|\pi\right|\big)=\max_{i}\left\{ t_{i}-t_{i-1}\right\} .\label{eq:qv2}
\end{equation}
Then show that the limit,
\begin{equation}
\lim_{\text{mesh}\left(\pi\right)\rightarrow0}\:\sum_{i}\left(X_{t_{i}}-X_{t_{i-1}}\right)^{2}=T\label{eq:qv3}
\end{equation}
holds a.e. on $\left(\Omega,\mathcal{F},\mathbb{P}\right)$, where
$\Omega=C\left(\Omega\right)$, $\mathcal{F}=$ Cyl, and $\mathbb{P}=$
the Wiener measure.

\uline{Hint}: Establish that 
\begin{eqnarray}
\mathbb{E}\left(\left|\triangle X_{i}\right|^{2}\right) & = & \triangle t_{i},\label{eq:qv4}\\
\mathbb{E}\left(\left|\triangle X_{i}\right|^{4}\right) & = & 3\left(\triangle t_{i}\right)^{2},\;\mbox{and}\label{eq:qv6}\\
\mathbb{E}\left(\left(\triangle X_{i}\right)^{2n-1}\right) & = & 0,\;n\in\mathbb{N},\label{eq:qv5}
\end{eqnarray}
i.e., all the odd moments vanish; where 
\begin{eqnarray*}
\triangle X_{i} & = & X_{t_{i}}-X_{t_{i-1}},\;\mbox{and}\\
\triangle t_{i} & = & t_{i}-t_{i-1},\;\mbox{for}\:i=1,2,\ldots,n.
\end{eqnarray*}

Note that $\sum_{i}\left(\right)^{2}$ on the LHS in (\ref{eq:qv3})
is a measurable function on $\left(\Omega,\mathcal{F},\mathbb{P}\right)$,
while the RHS is deterministic, i.e., it is the constant function
$T$.
\end{xca}

\begin{rem}
Spectral Theorem and \emph{functional calculus} is about the substitutions
(see (\ref{eq:fc1})). \index{functional calculus}
\begin{equation}
\boxed{\begin{matrix}A\\
\text{selfadjoint operator}
\end{matrix}}\longrightarrow\boxed{\begin{matrix}f:\mathbb{R}\longrightarrow\mathbb{R}\\
\text{scalar function}\\
\hline \hookrightarrow f\left(A\right)
\end{matrix}}\label{eq:fc1}
\end{equation}
By contrast, It\={o}-calculus is about substitutions of Brownian motion
(at least in a special case); as follows:
\begin{equation}
\boxed{\begin{matrix}B_{t}\\
\text{Brownian motion}
\end{matrix}}\longrightarrow\boxed{\begin{matrix}f:\mathbb{R}\longrightarrow\mathbb{R}\\
\text{scalar function}\\
\hline \hookrightarrow f\left(B_{t}\right)
\end{matrix}}\label{eq:fc2}
\end{equation}
See \cite{Sto90, Yos95, Ne69, RS75, DS88b}.
\end{rem}

\begin{figure}
\includegraphics[scale=0.4]{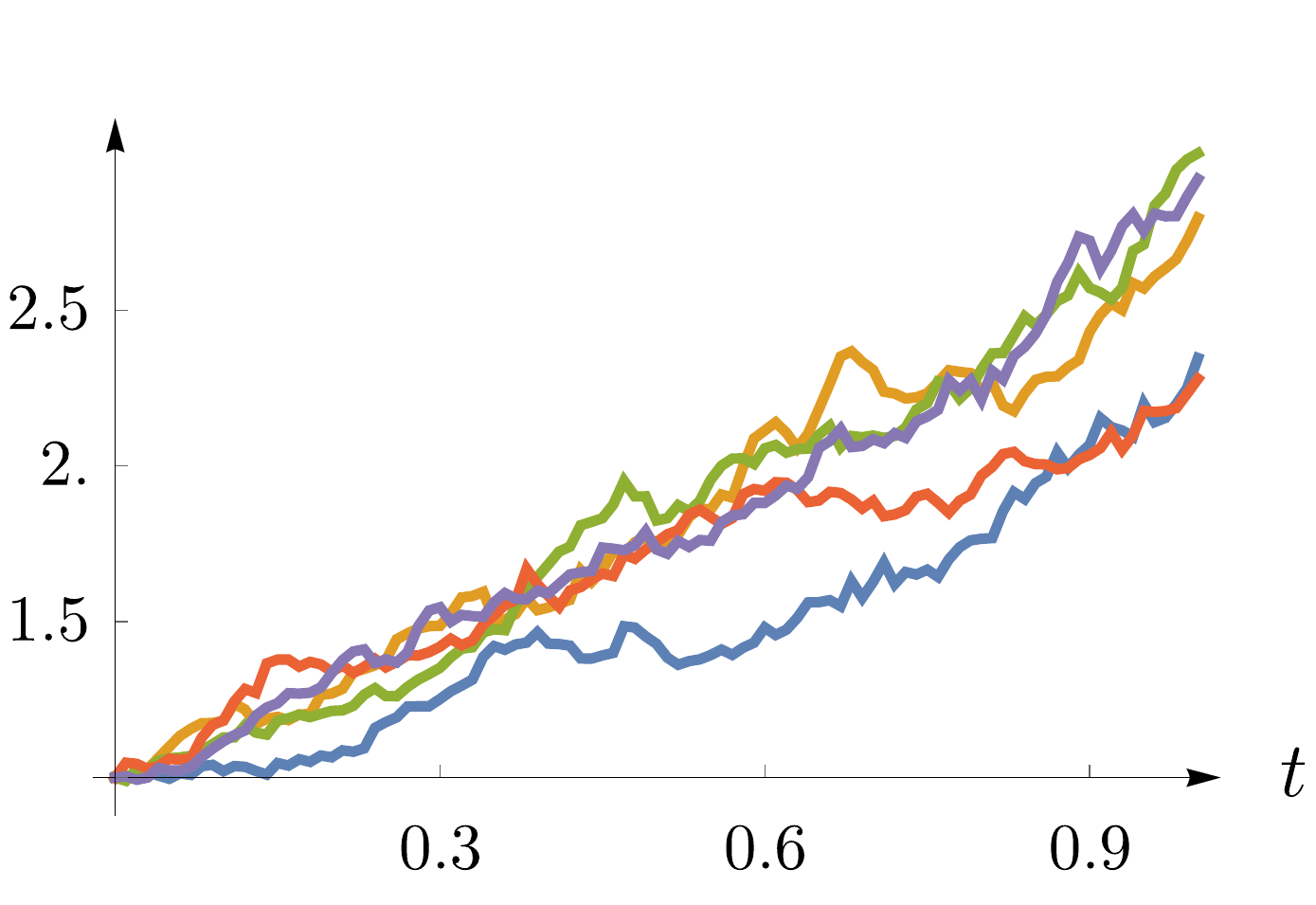}

\protect\caption{\label{fig:geobm}Geometric Brownian motion: 5 sample paths, with
$\mu=1$, $\sigma=0.02$, and $X_{0}=1$. }
\end{figure}

\begin{figure}
\includegraphics[scale=0.4]{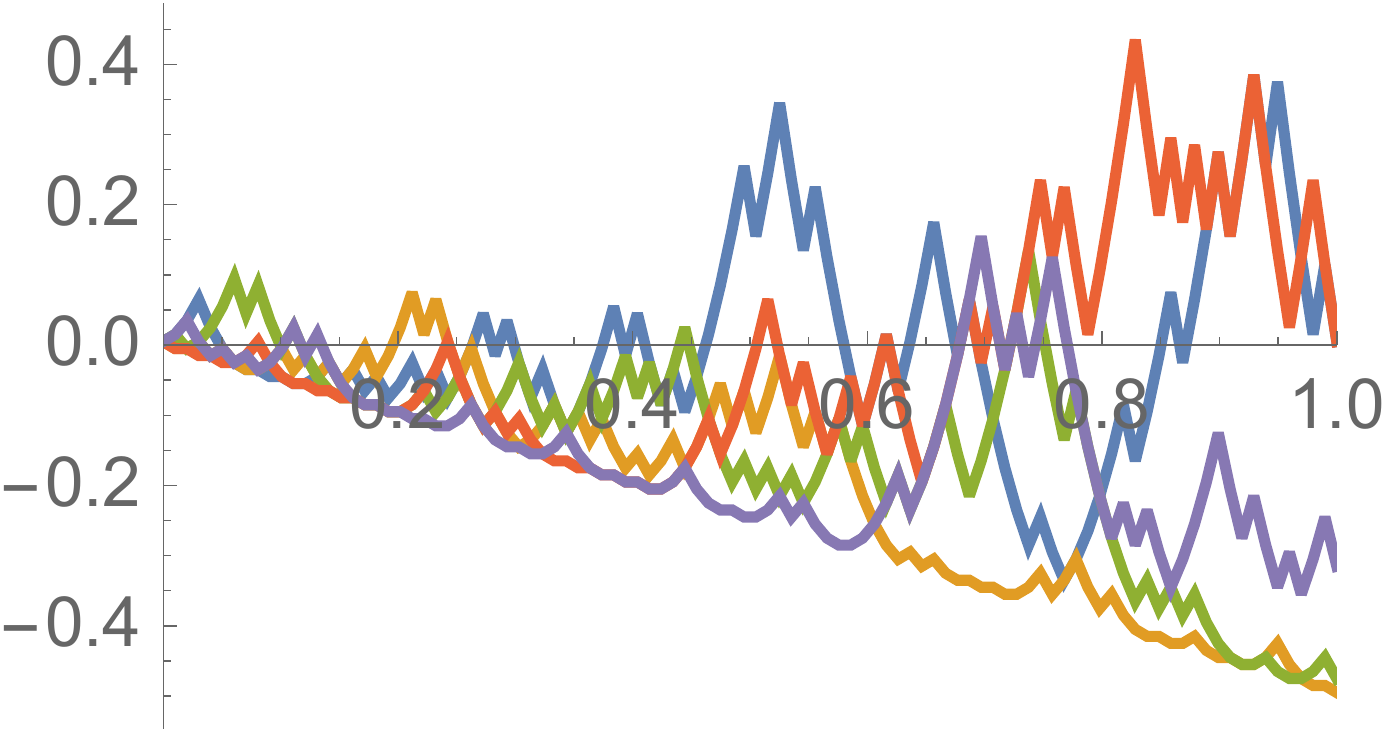}

\protect\caption{\label{fig:geobm2}The process $\frac{1}{2}\left(B_{T}^{2}-T\right)$
in (\ref{eq:fc5}): 5 sample paths, with $T=1$. }
\end{figure}

\index{Geometric Brownian motion} \index{It=o-calculus@It\=o-calculus}
\begin{xca}[Geometric Brownian motion]
\myexercise{Geometric Brownian motion}
\begin{enumerate}
\item Apply (\ref{eq:fc2}) to $f\left(x\right)=\ln x$, $x\in\mathbb{R}_{+}$,
together with (\ref{eq:qv3}) in \exerref{qv}, to show that, for
$T\in\mathbb{R}_{+}$, the process, 
\begin{equation}
X_{T}=X_{0}\exp\left(\left(\mu-\frac{1}{2}\sigma^{2}\right)T+\sigma B_{T}\right)\label{eq:fc3}
\end{equation}
solves the SDE for geometric Brownian motion:
\begin{equation}
dX_{t}=X_{t}\left(\mu dt+\sigma dB_{t}\right).\label{eq:fc4}
\end{equation}
See  \figref{geobm}.
\item Apply (\ref{eq:fc2}) to $f\left(x\right)=x^{2}$, together with (\ref{eq:qv3})
in \exerref{qv} to establish the following:
\begin{equation}
\int_{0}^{T}B_{t}\:dB_{t}=\frac{1}{2}\left(B_{T}^{2}-T\right).\label{eq:fc5}
\end{equation}
See \figref{geobm2}.
\end{enumerate}
\end{xca}

In the previous exercise we explored stochastic processes derived
from standard Brownian motion, but we now return to explore some additional
properties for Brownian motion itself:
\begin{xca}[A unitary one-parameter group]
\myexercise{A unitary one-parameter group}\label{exer:bm1}Using
\defref{bm}, \corref{gns}, \remref{gns-1} and \exaref{bm}, show
that there is a unique strongly continuous unitary one-parameter group
$\left\{ U\left(t\right)\right\} _{t\in\mathbb{R}}$ acting in $L^{2}\left(C\left(\mathbb{R}\right),\mbox{Cyl},\mathbb{P}\right)$,
determined by 
\begin{equation}
U\left(t\right)X_{s}=X_{s+t},\;\forall s,t\in\mathbb{R}.\label{eq:ubm}
\end{equation}
\uline{Hint}: By (\ref{enu:bm3}) in \defref{bm}, we have 
\begin{equation}
\mathbb{E}\left(\left|X_{t_{2}}-X_{t_{1}}\right|^{2}\right)=\mathbb{E}\left(\left|X_{t_{2}+s}-X_{t_{1}+s}\right|^{2}\right),\;\forall s,t_{1},t_{2}\in\mathbb{R}.\label{eq:ubm2}
\end{equation}
Hence, if $U\left(t\right)$ is defined on the generator $\left\{ X_{s}\::\:s\in\mathbb{R}\right\} \subset L^{2}\left(\mathbb{P}\right)$
as in (\ref{eq:ubm}), it follows by (\ref{eq:ubm2}) that it preserves
the $L^{2}\left(\mathbb{P}\right)$-norm. The remaining steps are
left to the reader. 
\end{xca}
\index{strongly continuous}

\index{representation!strongly continuous}

\begin{xca}[Infinitesimal generator]
\myexercise{Infinitesimal generator}\label{exer:bm2}Discuss the
infinitesimal generator of $\left\{ U\left(t\right)\right\} _{t\in\mathbb{R}}$.
\end{xca}

\begin{xca}[An ergodic action]
\myexercise{An ergodic action}\label{exer:bm3}Show that $\left\{ U\left(t\right)\right\} _{t\in\mathbb{R}}$
is induced by an \emph{ergodic action}.\index{ergodic}
\end{xca}

\section{The Path Space\label{sec:pspace}}

\index{Wiener, N.}

\begin{figure}[H]
\includegraphics[scale=0.4]{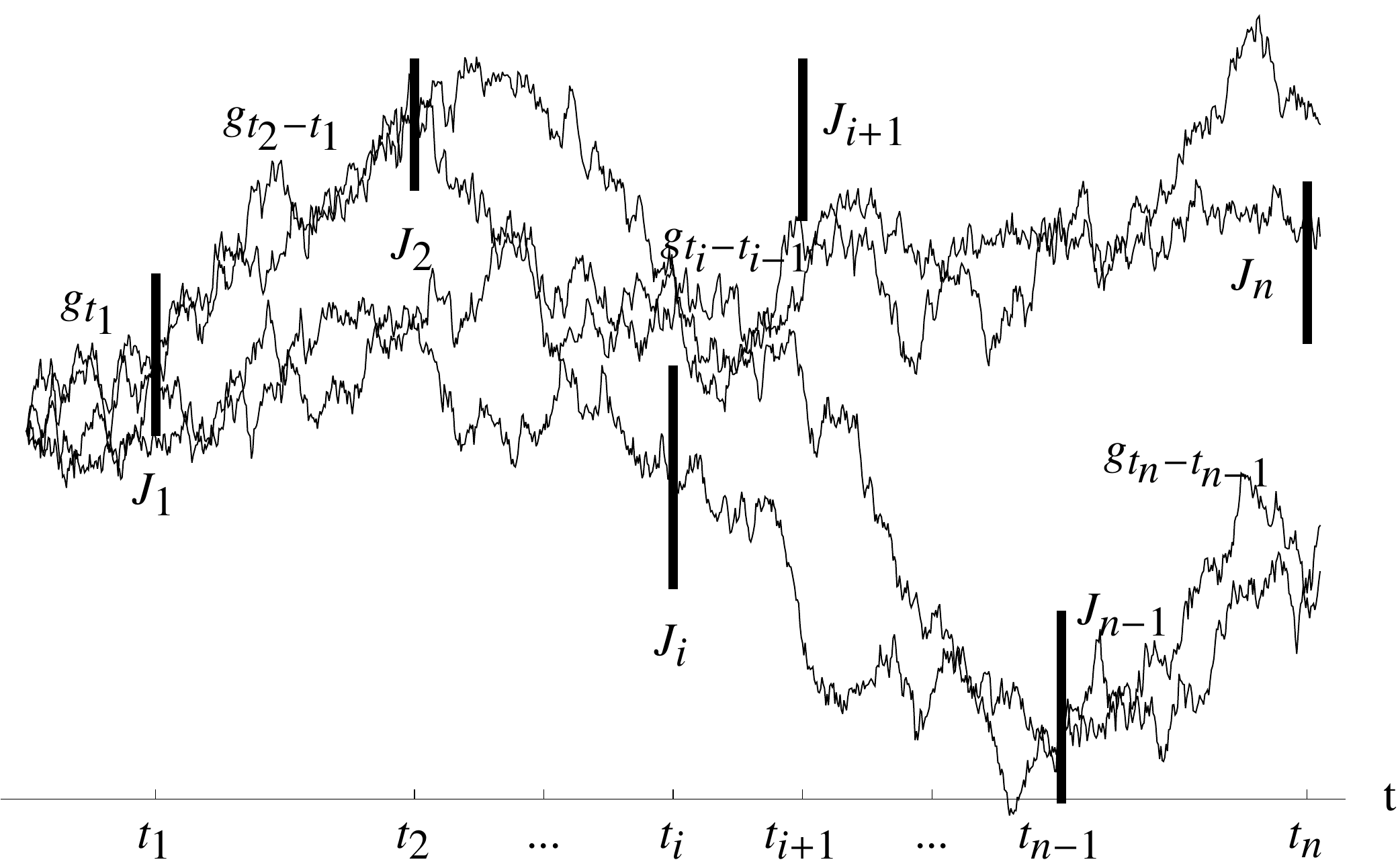}

\protect\caption{\label{fig:cyl}Stochastic processes indexed by time: A cylinder set
$C$ is a special subset of the space of all paths, i.e., functions
of a time variable. A fixed cylinder set $C$ is specified by a finite
set of sample point on the time-axis (horizontal), and a corresponding
set of \textquotedblleft windows\textquotedblright{} (intervals on
the vertical axis). When sample points and intervals are given, we
define the corresponding cylinder set $C$ to be the set of all paths
that pass through the respective windows at the sampled times. In
the figure we illustrate sample points (say future relative to $t=0$).
Imagine the set $C$ of all outcomes with specification at the points
$t_{1},t_{2},\ldots$ etc.}
\end{figure}

\paragraph*{The sample-space as a path-space. }
\begin{thm}[see e.g., \cite{Nel67}]
\label{thm:cnhd}Set $\Omega=C\left(\mathbb{R}\right)=$ (all continuous
real valued function on $\mathbb{R}$), $\mathcal{F}=$ the sigma
algebra generated by cylinder-sets, i.e., determined by finite systems
$t_{1},\ldots,t_{n}$, and intervals $J_{1},\ldots,J_{n}$; 
\begin{equation}
Cyl\left(t_{1},\ldots,t_{n},J_{1},\ldots,J_{n}\right)=\left\{ \omega\in C\left(\mathbb{R}\right)\:\big|\:\omega\left(t_{i}\right)\in J_{i},\;i=1,2\ldots,n\right\} .\label{eq:bm3}
\end{equation}
The measure $\mathbb{P}$ is determined by its value on cylinder sets,
and and an integral over Gaussians; it is called the Wiener-measure.
Set 
\[
X_{t}\left(\omega\right)=\omega\left(t\right),\;\forall t\in\mathbb{R},\omega\in\Omega\left(=C\left(\mathbb{R}\right)\right).
\]

If $0<t_{1}<t_{2}<\cdots<t_{n}$, and the cylinder set is as \eqref{bm3},
then 
\begin{eqnarray*}
 &  & \mathbb{P}\left(Cyl\left(t_{1},\ldots,t_{n},J_{1},\ldots,J_{n}\right)\right)\\
 & = & \int_{J_{1}}\cdots\int_{J_{n}}g_{t_{1}}\left(x_{1}\right)g_{t_{2}-t_{1}}\left(x_{2}-x_{1}\right)\cdots g_{t_{n}-t_{n-1}}\left(x_{n}-x_{n-1}\right)dx_{1}\cdots dx_{n}
\end{eqnarray*}
where 
\[
g_{t}\left(x\right)=\frac{1}{\sqrt{2\pi t}}e^{-x^{2}/2t},\;\forall t>0,
\]
i.e., the $N(0,t)$-Gaussian. See \figref{cyl}. \index{cylinder-set}
\end{thm}

\begin{figure}
\subfloat[$C$ optimistic]{\protect\includegraphics[scale=0.35]{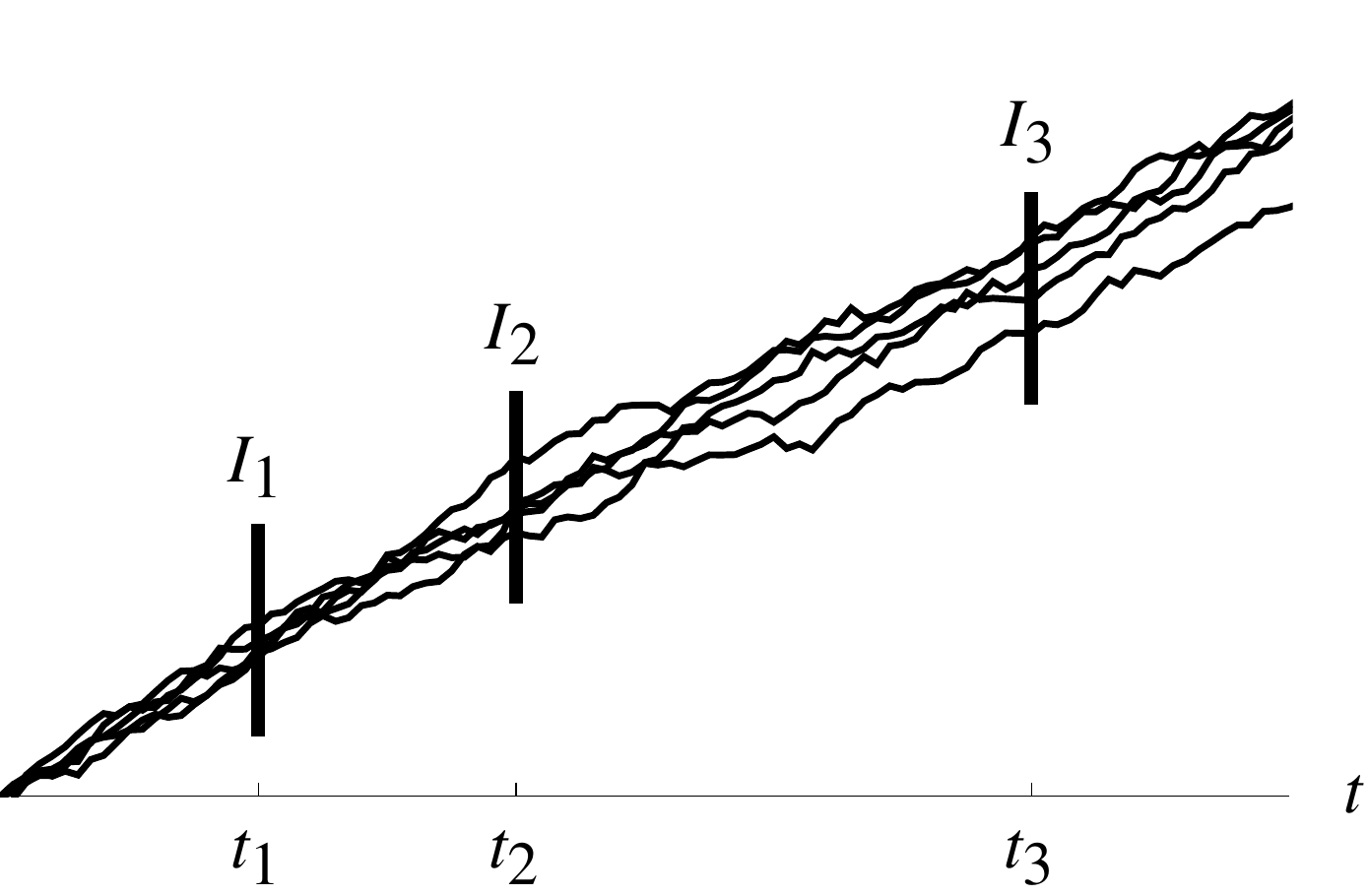}

}\hfill{}\subfloat[$C'$ pessimistic]{\protect\includegraphics[scale=0.35]{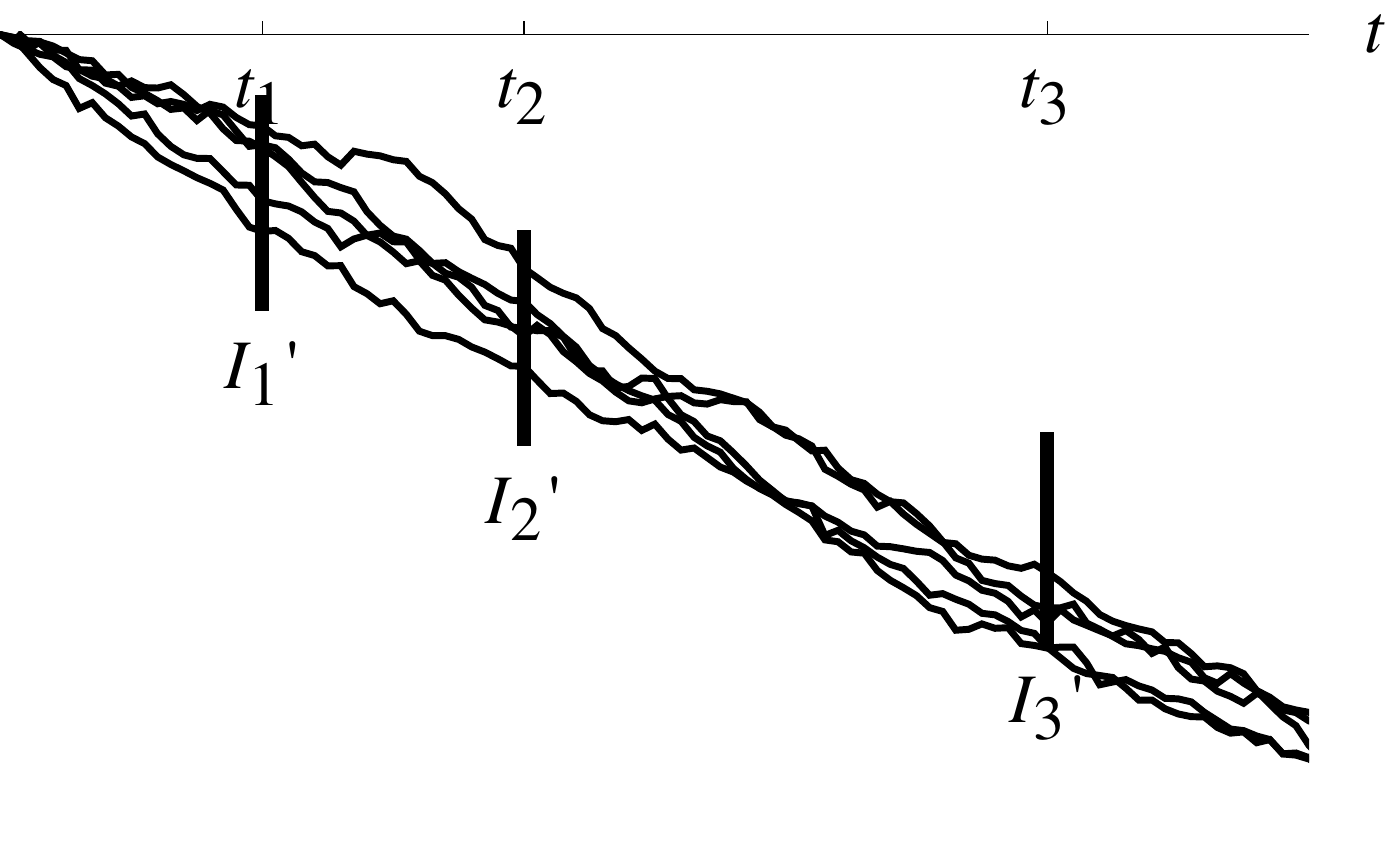}

}\hfill{}\subfloat[$C''$ mixed]{\protect\includegraphics[scale=0.35]{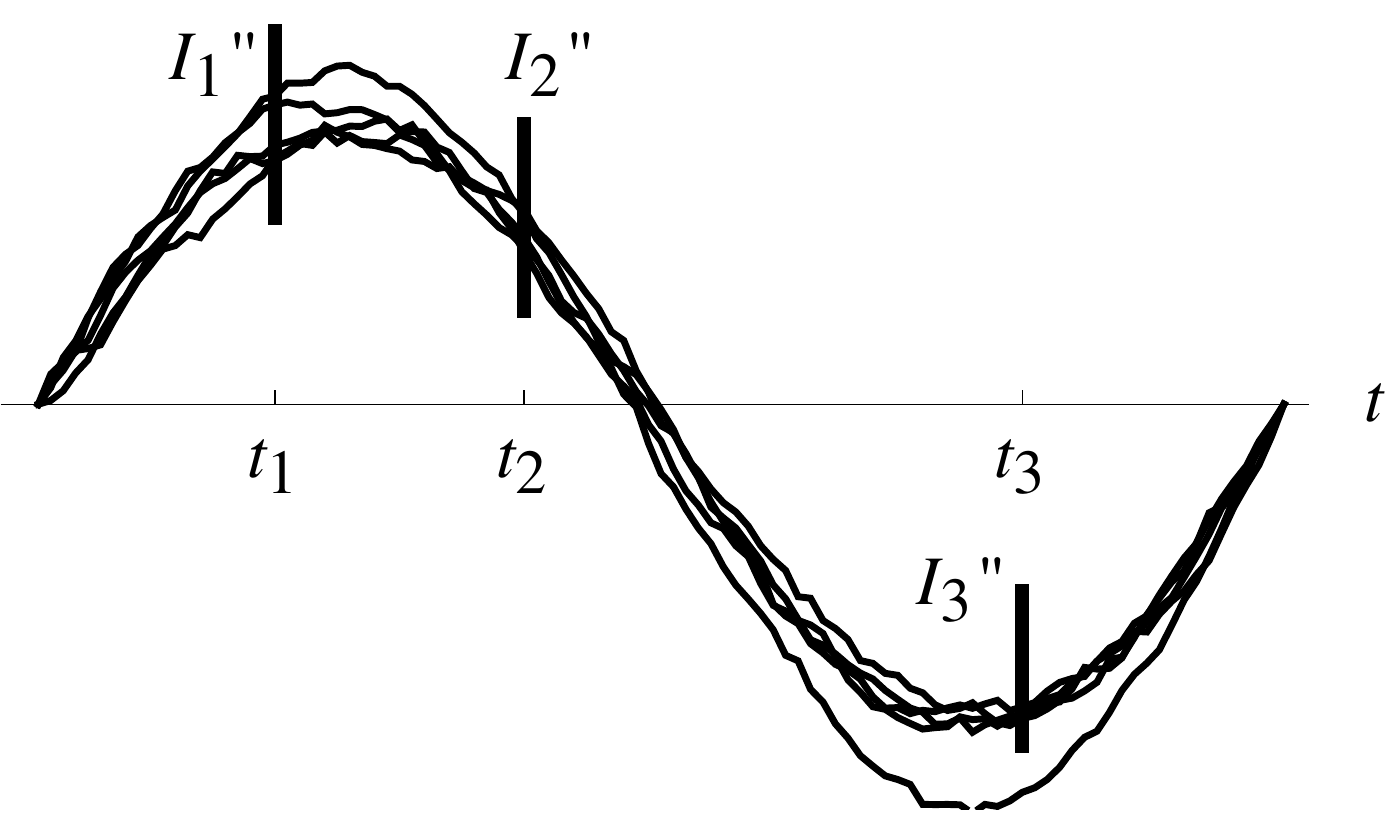}

}

\protect\caption{The cylinder sets $C$, $C'$, and $C''$.}
\end{figure}

\paragraph*{Infinite-product Measure}

Let $\Omega=\prod_{k=1}^{\infty}\{1,-1\}$ be the infinite Cartesian
product of $\{1,-1\}$ with the product topology. $\Omega$ is compact\index{compact}
and Hausdorff by Tychnoff's theorem. \index{product topology}

For each $k\in\mathbb{N}$, let $X_{k}:\Omega\rightarrow\{1,-1\}$
be the $k^{th}$ coordinate projection, and assign probability measures
$\mu_{k}$ on $\Omega$ so that $\mu_{k}\circ X_{k}^{-1}\{1\}=a$
and $\mu_{k}\circ X_{k}^{-1}\{-1\}=1-a$, where $a\in(0,1)$. The
collection of measures $\{\mu_{k}\}$ satisfies the consistency condition,
i.e., $\mu_{k}$ is the restriction of $\mu_{k+1}$ onto the $k^{th}$
coordinate space. By Kolmogorov's extension theorem, there exists
a unique probability\index{measure!probability} measure $P$ on $\Omega$
so that the restriction of $P$ to the $k^{th}$ coordinate is equal
to $\mu_{k}$. 

It follows that $\{X_{k}\}$ is a sequence of independent identically
distributed (i.i.d.) random variables in $L^{2}(\Omega,P)$ with $\mathbb{E}\left(X_{k}\right)=0$
and $Var[X_{k}^{2}]=1$; and $L^{2}(\Omega,P)=\overline{span}\{X_{k}\}$.
\index{independent!-random variables}\index{random variable}
\begin{rem}
Let $\mathscr{H}$ be a separable Hilbert space with an orthonormal
basis $\{u_{k}\}$. The map $\varphi:u_{k}\mapsto X_{k}$ extends
linearly to an isometric embedding of $\mathscr{H}$ into $L^{2}(\Omega,P)$.
Moreover, let $\mathcal{F}_{+}(\mathscr{H})$ be the \emph{symmetric
Fock space}. $\mathcal{F}_{+}(\mathscr{H})$ is the closed span of
the the algebraic tensors $u_{k_{1}}\otimes\cdots\otimes u_{k_{n}}$,
thus $\varphi$ extends to an isomorphism from $\mathcal{F}_{+}(\mathscr{H})$
to $L^{2}(\Omega,P)$. 

\index{orthonormal basis (ONB)}\index{orthogonal!-vectors}\end{rem}
\begin{xca}[The ``fair-coin'' measure]
\myexercise{The fair-coin measure}\label{exer:bm4}Let $\Omega=\prod_{1}^{\infty}\left\{ -1,1\right\} $,
and let $\mu$ be the ``fair-coin'' measure $\left(\frac{1}{2},\frac{1}{2}\right)$
on $\left\{ \pm1\right\} $ (e.g., ``Head v.s. Tail''), let $\mathcal{F}$
be the cylinder sigma-algebra\index{sigma-algebra} of subsets of
$\Omega$. Let $\mathbb{P}=\prod_{1}^{\infty}\mu$ be the infinite-product
measure on $\Omega$. Set $Z_{k}\left(\omega\right)=\omega_{k}$,
$\omega=\left(\omega_{i}\right)\in\Omega$. Finally, let $\left\{ \psi_{j}\right\} _{j\in\mathbb{N}}$
be an ONB in $L^{2}\left(0,1\right)$, and set 
\[
X_{t}\left(\omega\right):=\sum_{j=1}^{\infty}\left(\int_{0}^{t}\psi_{j}\left(s\right)ds\right)Z_{j}\left(\omega\right),\;t\in\left[0,1\right],\omega\in\Omega.
\]
Show that the $\left\{ X_{t}\right\} _{t\in\left[0,1\right]}$ is
Brownian motion, where time ``$t$'' is restricted to $\left[0,1\right]$.
\index{integral!Stochastic-}

\uline{Hint:} Let $t_{1},t_{2}\in\left[0,1\right]$, then 
\begin{eqnarray*}
\mathbb{E}\left(X_{t_{1}}X_{t_{2}}\right) & = & \mathbb{E}\left(\sum_{j}\left(\int_{0}^{t_{1}}\psi_{j}\left(s\right)ds\right)Z_{j}\cdot\sum_{k}\left(\int_{0}^{t_{2}}\psi_{k}\left(s\right)ds\right)Z_{k}\right)\\
 & = & \sum_{j}\sum_{k}\int_{0}^{t_{1}}\psi_{j}\left(s\right)ds\int_{0}^{t_{2}}\psi_{k}\left(s\right)ds\underset{=\delta_{jk}}{\underbrace{\mathbb{E}\left(Z_{j}Z_{k}\right)}}\\
 & = & \sum_{j}\left\langle \chi_{\left[0,t_{1}\right]},\psi_{j}\right\rangle _{L^{2}}\left\langle \chi_{\left[0,t_{2}\right]},\psi_{j}\right\rangle _{L^{2}}\\
 & = & \left\langle \chi_{\left[0,t_{1}\right]},\chi_{\left[0,t_{2}\right]}\right\rangle _{L^{2}\left(0,1\right)}=t_{1}\wedge t_{2}\left(:=\min\left(t_{1},t_{2}\right)\right).
\end{eqnarray*}

\end{xca}

\section{Decomposition of Brownian motion\label{sec:dbm}}

\index{Brownian motion}

The integral kernel $K:[0,1]\times[0,1]\rightarrow\mathbb{R}$, $K(s,t)=s\wedge t$,
is a compact operator on $L^{2}[0,1]$, where
\[
Kf(x)=\int_{0}^{1}(x\wedge y)f(y)dy.
\]
$Kf$ is a solution to the differential equation \index{positive definite!-kernel}
\[
-\frac{d^{2}}{dx^{2}}u=f
\]
with zero boundary conditions.

\index{boundary condition} 

\index{kernel!integral-}

$K$ is also seen as the covariance functions of Brownian motion process.
A stochastic process is a family of measurable functions $\{X_{t}\}$
defined on some sample probability space $\left(\Omega,\mathfrak{B},P\right)$,
where the parameter $t$ usually represents time. $\{X_{t}\}$ is
a Brownian motion process if it is a mean zero Gaussian process such
that \index{Gaussian process} 
\[
E[X_{s}X_{t}]=\int_{\Omega}X_{s}X_{t}dP=s\wedge t.
\]
It follows that the corresponding increment process $\{X_{t}-X_{s}\}\sim N(0,t-s)$.
$P$ is called the Wiener measure. \index{measure!Wiener}\index{stochastic process}

Building $(\Omega,\mathfrak{B},P)$ is a fancy version of Riesz's
representation theorem \cite[Theorem 2.14]{Rud87}. It turns out that
\index{Riesz' theorem}\index{Theorem!Riesz-} 
\[
\Omega=\prod_{t}\bar{\mathbb{R}}
\]
which is a compact Hausdorff space; (where $\overline{\mathbb{R}}=\left(\mathbb{R}\cup\left\{ \infty\right\} \right)^{\sim}$
denotes the one-point compactification of $\mathbb{R}$.) 

Now introduce random variable $X_{t}:\Omega\rightarrow\mathbb{R}$,
defined as
\[
X_{t}(\omega)=\omega(t),\;t\in\mathbb{R};
\]
i.e., $X_{t}$ is the continuous linear functional of evaluation at
$t$ on $\Omega$.

For Brownian motion, the increment of the process $\triangle X_{t}$,
in some statistical sense, is proportional to $\sqrt{\triangle t}$,
i.e., \index{continuous linear functional} 
\[
\triangle X_{t}\sim\sqrt{\triangle t}.
\]
It it this property that makes the set of differentiable functions
have measure zero. In this sense, the trajectory of Brownian motion
is nowhere differentiable. 

A very important application of the spectral theorem of compact operators
is to decompose the Brownian motion process: 
\begin{equation}
B_{t}(\omega)=\sum_{n=1}^{\infty}\frac{\sin(n\pi t)}{n\pi}Z_{n}(\omega)\label{eq:bmd}
\end{equation}
where 
\[
s\wedge t=\sum_{n=1}^{\infty}\frac{\sin\left(n\pi s\right)\sin\left(n\pi t\right)}{\left(n\pi\right)^{2}}
\]
and $Z_{n}\sim N(0,1)$.
\begin{rem}
Consider the Hardy space $\mathbb{H}_{2}$, and the operator $S$
from \exerref{infm}. Writing $f\left(z\right)=\sum_{n=0}^{\infty}x_{n}z^{n}$,
we get 
\begin{equation}
\left(Sf\right)\left(z\right)=f\left(z^{N}\right)=x_{0}+x_{1}z^{N}+x_{2}z^{2N}+\cdots;\label{eq:b1-1}
\end{equation}
and 
\begin{equation}
\left(S^{*}f\right)\left(z\right)=x_{0}+x_{N}z+x_{2N}z^{2}+x_{3N}z^{3}+\cdots;\label{eq:b1-2}
\end{equation}
so in symbol-space, $S^{*}$ acts as follows:
\begin{eqnarray}
 &  & \left(x_{0},x_{1},\cdots,x_{N-1},x_{N},x_{N+1},\cdots,x_{2N},x_{2N+1},\cdots\right)\nonumber \\
S^{*} & \downarrow\label{eq:b1-3}\\
 &  & \left(x_{0},x_{N},x_{2N},x_{3N},\cdots\right);\nonumber 
\end{eqnarray}
so down-sampling $\simeq$ ``decimation'' $\simeq$ killing time-signals
$x_{k}$ when $N+k$, i.e., $k$ is not divisible by $N$. \index{sampling!down-sampling}\index{down-sampling}\index{space!Hardy-}

The projection $SS^{*}$ is:
\begin{align}
\big(x_{0},x_{1},\cdots,x_{N-1},x_{N},x_{N+1},\cdots,x_{2N-1}, & x_{2N},x_{2N+1},\cdots,x_{3N-1},x_{3N},x_{3N+1},\cdots\big)\nonumber \\
SS^{*} & \downarrow\label{eq:b1-4}\\
\big(x_{0},0,\cdots,0,x_{N},0,\cdots, & 0,x_{2N},0\cdots,0,x_{3N},0,\cdots\big)\nonumber 
\end{align}
If the coordinates in $\mathbb{H}_{2}$ label the i.i.d. random variables
$Z_{k}\left(\cdot\right)$ in the expansion (\ref{eq:bmd}) for Brownian
motion, then downsampling corresponds to \emph{conditional expectation};
conditioning on ``less information'', i.e., leaving out the ``decimated
coordinates'' in the expansion (\ref{eq:bmd}) for Brownian motion.

\index{conditional expectation} \index{down-sampling}\end{rem}
\begin{xca}[The Central Limit Theorem]
\myexercise{The Central Limit Theorem}\label{exer:bm5}Look up the
\emph{Central Limit Theorem} (CLT), and prove the following approximation
formula for Brownian motion: \index{Theorem!Central Limit-}\index{Theorem!Spectral-}

Let $\pi$ be the ``fair-coin-measure'' on the two outcomes $\left\{ \pm1\right\} $,
i.e., winning or loosing one unit, and let $\Omega=\vartimes_{\mathbb{N}}\left\{ \pm1\right\} $,
$\mathbb{P}=\vartimes_{\mathbb{N}}\pi$ be the corresponding infinite
product measure. On $\Omega$, set 
\[
W_{k}\left(\omega\right)=\omega_{k},\;\omega=\left(\omega_{1},\omega_{2},\ldots\right)\in\Omega,\;k=1,2,\ldots;\;\mbox{and}
\]
\begin{equation}
S_{n}\left(\cdot\right)=\frac{1}{\sqrt{n}}\sum_{k=1}^{n}W_{k}\left(\cdot\right).\label{eq:bma1}
\end{equation}
Let $X_{t}$ denote Brownian motion. Then show that 
\begin{equation}
X_{t}\left(\cdot\right)=\lim_{n\rightarrow\infty}\frac{1}{\sqrt{n}}\sum_{k=1}^{\left\lfloor n\,t\right\rfloor }W_{k}\left(\cdot\right)\label{eq:bma2}
\end{equation}
where $\left\lfloor n\,t\right\rfloor $ denotes the largest integer
$\leq n\,t$. 

\uline{Hint:} A good reference to the CLT is \cite{CW14}. First
approximate the CLT to the sequence $S_{n}$ in (\ref{eq:bma1}).
 \figref{bm} illustrates the approximation formula in (\ref{eq:bma2}).
($X_{1}=S$.)

\begin{figure}
\subfloat[$n=10$]{\protect\includegraphics[scale=0.4]{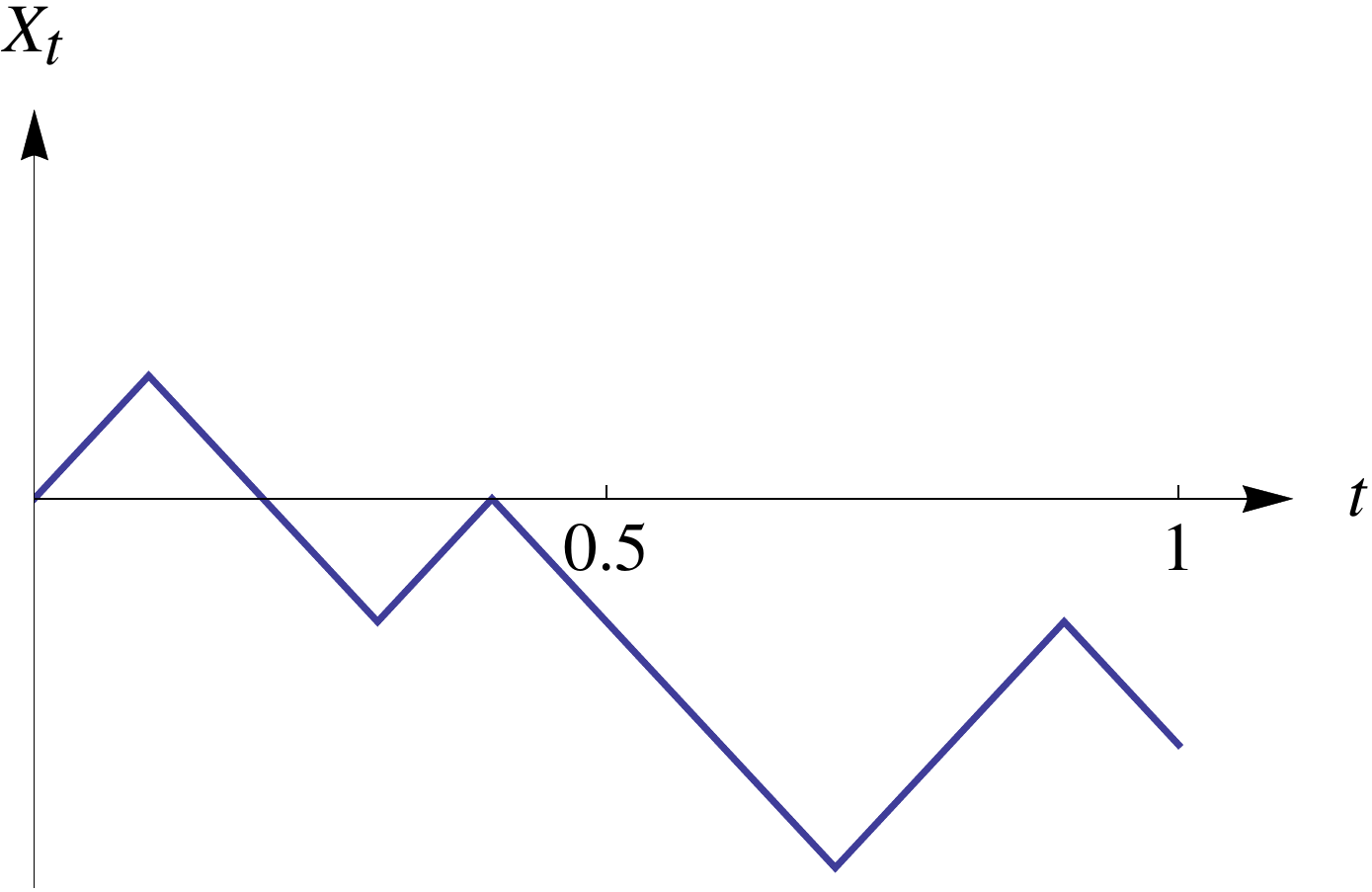}

}\hfill{}\subfloat[$n=50$]{\protect\includegraphics[scale=0.4]{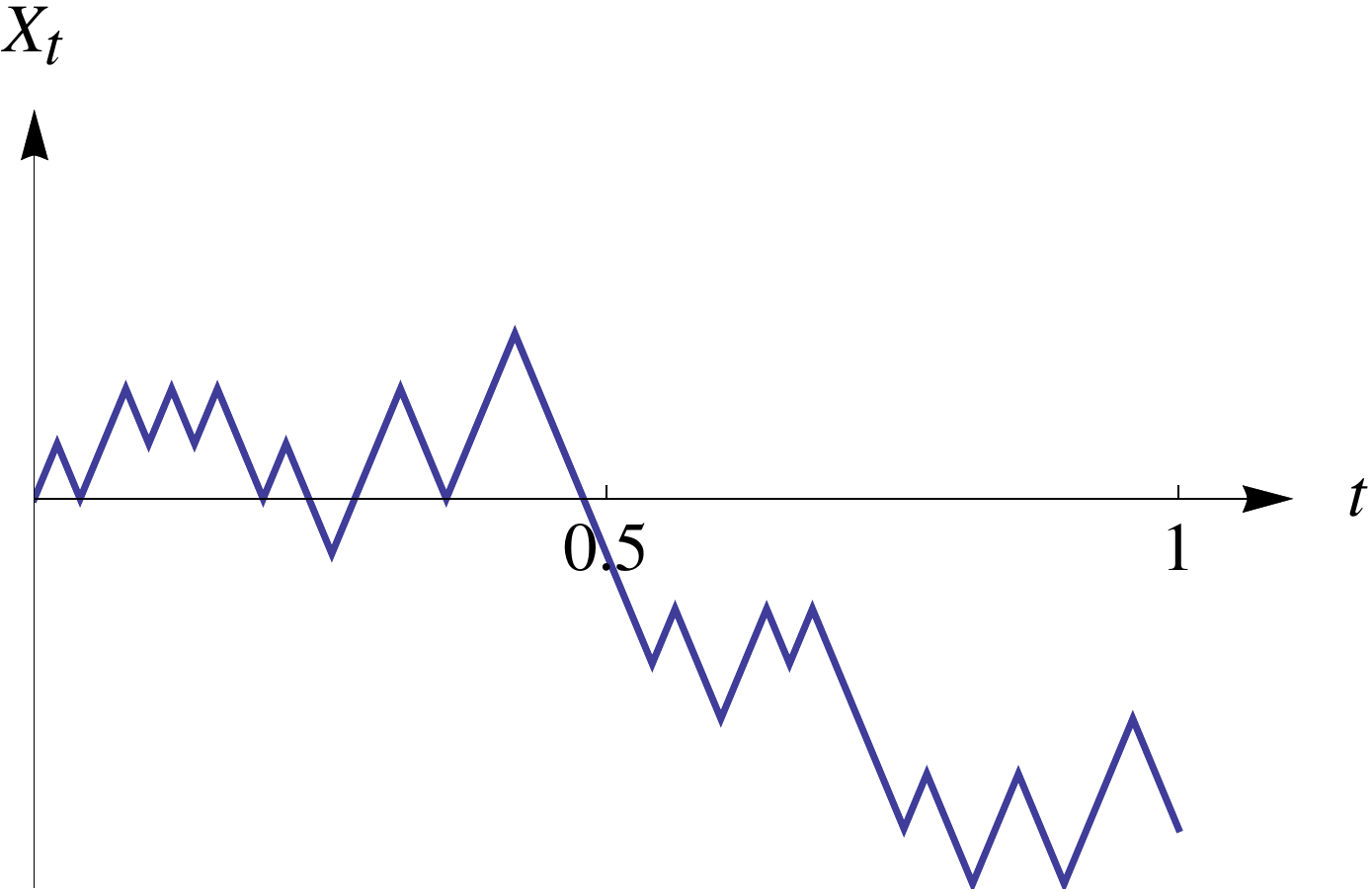}

}\hfill{}\subfloat[$n=100$]{\protect\includegraphics[scale=0.4]{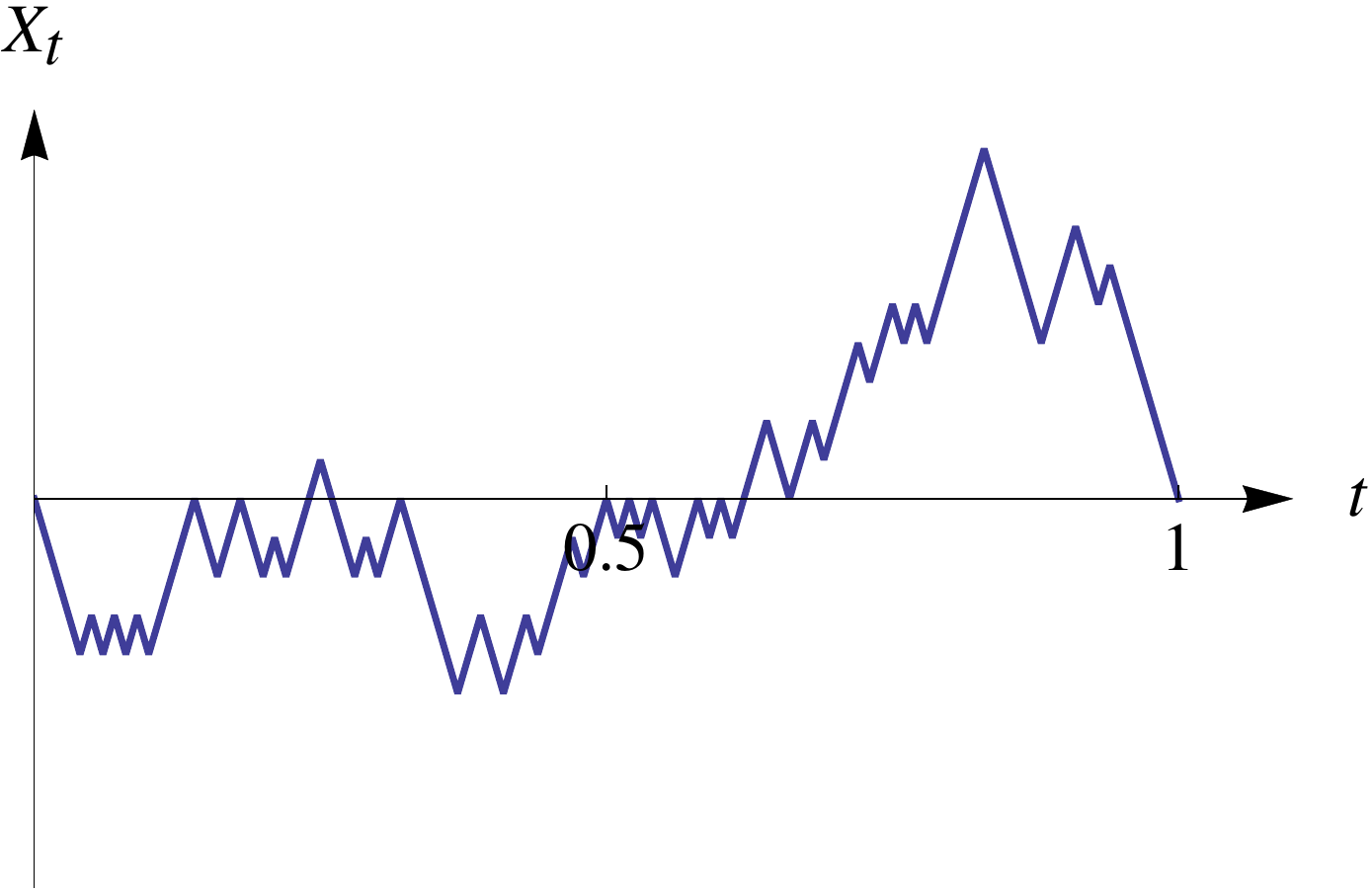}

}\hfill{}\subfloat[$n=500$]{\protect\includegraphics[scale=0.4]{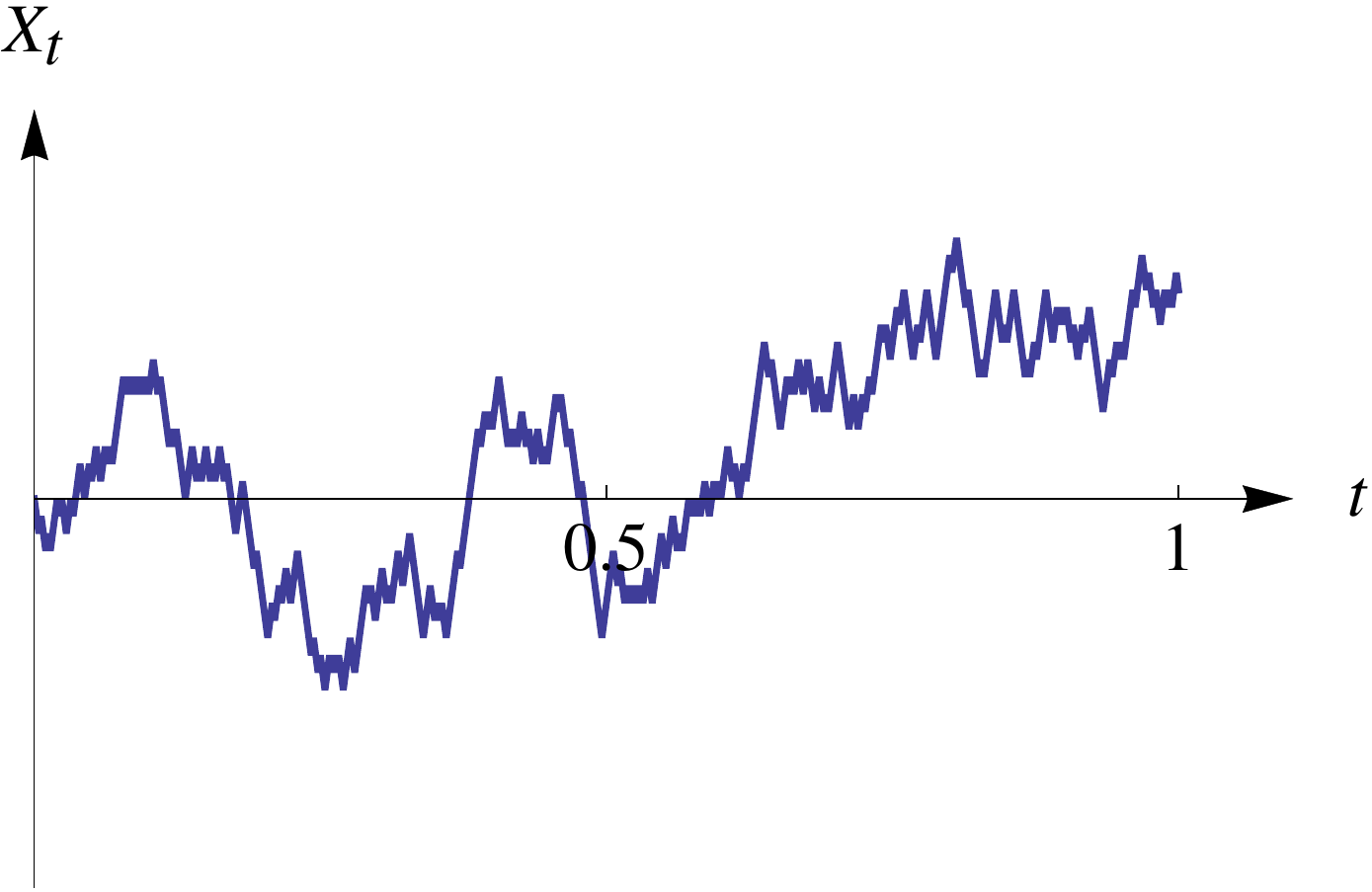}

}

\protect\caption{\label{fig:bm}Monte-Carlo simulation of the standard Brownian motion
process $\left\{ X_{t}\::\:0\leq t\leq1\right\} $, where $\mathbb{E}\left(X_{t}\right)=0$,
and $\mathbb{E}\left(X_{s}X_{t}\right)=s\wedge t=\min\left(s,t\right)$.\protect \\
For $n=10,50,100,500$, set $X_{0}=0$ and $X_{j/n}^{\left(n\right)}=n^{-1/2}\sum_{k=1}^{j}W_{k}$.
Applying linear interpolation between sample points $\left\{ j/n:j=0,\ldots,n\right\} $
yields the $n$-point approximation $X_{t}^{\left(n\right)}$, which
converges in measure to $X_{t}$ (standard BM restricted to the unit
interval), as $n\rightarrow\infty$. }
\end{figure}

The Central Limit Theorem (CLT) states that the limit, $n\rightarrow\infty$,
of the sequence $S_{n}$ in (\ref{eq:bma1}) is a copy of $N\left(0,1\right)$-random
variable; i.e., $\lim_{n\rightarrow\infty}S_{n}\left(\cdot\right)=S\left(\cdot\right)$
exists; and 
\[
\mathbb{P}\left(\left\{ \omega\:\big|\:a\leq S\left(\omega\right)\leq b\right\} \right)=\int_{a}^{b}\frac{1}{\sqrt{2\pi}}e^{-\frac{1}{2}x^{2}}dx,
\]
for all intervals $\left(a,b\right)\subset\mathbb{R}$. Applying this
to (\ref{eq:bma2}), we get existence of $X_{t}$ as a limit, and
$X_{t}\sim N\left(0,t\right)$, i.e., 
\[
\mathbb{P}\left(\left\{ \omega\:\big|\:a\leq X_{t}\left(\omega\right)\leq b\right\} \right)=\int_{a}^{b}\frac{1}{\sqrt{2\pi t}}e^{-\frac{1}{2t}x^{2}}dx.
\]
We claim that
\begin{equation}
\mathbb{E}_{\mathbb{P}}\left(X_{s}X_{t}\right)=s\wedge t.\label{eq:bma3}
\end{equation}
Below we sketch the argument for the assertion in (\ref{eq:bma3}).

Fix $s,t\in\mathbb{R}_{+}$, say $s<t$, and $n\in\mathbb{N}$; then
\begin{eqnarray*}
\mathbb{E}\left(\frac{1}{\sqrt{n}}\left(\sum_{j=1}^{\left\lfloor n\,s\right\rfloor }W_{j}\right)\frac{1}{\sqrt{n}}\left(\sum_{k=1}^{\left\lfloor n\,t\right\rfloor }W_{k}\right)\right) & = & \frac{1}{n}\sum_{j=1}^{\left\lfloor n\,s\right\rfloor }\sum_{k=1}^{\left\lfloor n\,t\right\rfloor }\mathbb{E}\left(W_{j}W_{k}\right)\\
 & = & \frac{1}{n}\sum_{j=1}^{\left\lfloor n\,s\right\rfloor }\sum_{k=1}^{\left\lfloor n\,t\right\rfloor }\delta_{j,k}\\
 & = & \frac{\left\lfloor n\,s\right\rfloor }{n}\rightarrow s,\;\mbox{as}\;n\rightarrow\infty.
\end{eqnarray*}
Hence by the CLT, the desired conclusion in (\ref{eq:bma3}) follows.

This is the key step in proving that the limit $X_{t}$ in (\ref{eq:bma2})
is Brownian motion. The remaining steps are routine left to the readers.
\end{xca}

\section{Large Matrices Revisited}

Since large matrices and limit distributions have played a role in
several topics in the present chapter, we mention here yet a different
one; but now without proofs. Readers will find a complete treatment,
for example in \cite{MR0095527,MR1647832}. 

\index{matrix!random-} \index{limit distribution}

\index{random symmetric matrix}

\index{operators!symmetric-}

\uline{Setting} For all $N\in\mathbb{N}$, consider a symmetric
(random) matrix
\begin{equation}
X=\begin{bmatrix}X_{1,1}^{\left(N\right)} & \cdots & X_{1,N}^{\left(N\right)}\\
\vdots &  & \vdots\\
X_{N,1}^{\left(N\right)} & \cdots & X_{N,N}^{\left(N\right)}
\end{bmatrix}\label{eq:rm1}
\end{equation}
\begin{equation}
X_{i,j}^{\left(N\right)}=X_{j,i}^{\left(N\right)}\;\left(\mbox{real valued}\right)\label{eq:rm2}
\end{equation}
where the entries are i.i.d. random variables (``i.i.d'' is short
for independent identically distributed), mean $0$, and variance
$m^{2}$. And all the moments finite, and with at most exponential
bounds.

Let $a,b\in\mathbb{R}$, $a<b$ be fixed, and set\index{independent!-random variables}\index{random variable}
\begin{eqnarray}
V_{N}^{\left(a,b\right)} & := & \#\;\mbox{of eigenvalues of}\;X^{N}\:\mbox{that}\label{eq:rm3}\\
 &  & \mbox{fall in the interval}\;\left(a\sqrt{N},b\sqrt{N}\right).\nonumber 
\end{eqnarray}
Then the following limit exists, i.e., the semicircle law holds for
the limit distribution of the eigenvalues: \index{semicircle law}\index{eigenvalue}
\begin{equation}
\lim_{N\rightarrow\infty}\frac{\mathbb{E}\left(V_{N}^{\left(a,b\right)}\right)}{N}=\frac{1}{2\pi m^{2}}\int_{a}^{b}\sqrt{4m^{2}-x^{2}}dx.\label{eq:rm4}
\end{equation}

\section*{A summary of relevant numbers from the Reference List}

For readers wishing to follow up sources, or to go in more depth with
topics above, we suggest: \cite{Hi80,MR2053326,Ito06,MR2397781,Nel67,Par09,MR0265548,MR0184066,MR0161189,MR2966130,MR3231624,BoMc13,MR0107312,MR2883397,Jor14}.

\chapter{Lie Groups, and their Unitary Representations\label{chap:groups}}

\chaptermark{Lie groups, and unitary representations}
\begin{quotation}
Every axiomatic (abstract) theory admits, as is well known, an unlimited
number of concrete interpretations besides those from which it was
derived. Thus we find applications in fields of science which have
no relation to the concepts of random event and of probability in
the precise meaning of these words. 

--- A.N. Kolmogorov\sindex[nam]{Kolmogorov, A.N., (1903-1987)} \vspace{1em}\\
The miracle of the appropriateness of the language of mathematics
for the formulation of the laws of physics is a wonderful gift which
we neither understand nor deserve. 

--- Eugene Paul Wigner\sindex[nam]{Wigner, E.P., (1902-1995)}\vspace{1em}\\
\textquotedbl{}Nowadays group theoretical methods--especially those
involving characters and representations, pervade all branches of
quantum mechanics.\textquotedbl{} 

--- George Whitelaw Mackey\sindex[nam]{Mackey, G.W., (1916-2006)}\vspace{1em}\\
\textquotedblleft The universe is an enormous direct product of representations
of symmetry groups.\textquotedblright{} 

--- Hermann Weyl\vspace{2em}
\end{quotation}
As part of our discussion of spectral theory and harmonic analysis,
we had occasion to study unitary one-parameter groups $\mathcal{U}(t)$,
$t\in\mathbb{R}$. Stated differently (see Chapters \ref{chap:lin}-\ref{chap:GNS}),
a unitary one-parameter group acting on a Hilbert space $\mathscr{H}$,
is a strongly continuous unitary representation of the group $\mathbb{R}$
with addition; -- so it is an element in $Rep(\mathbb{R},\mathscr{H})$.
Because of applications to physics, to non-commutative harmonic analysis,
to stochastic processes, and to geometry, it is of interest to generalize
to $Rep(G,\mathscr{H})$ where $G$ is some more general group, other
than $(\mathbb{R},+)$, for example, $G$ may be a matrix group, a
Lie group, both compact and non-compact, or more generally, $G$ may
be a locally compact group.

In Chapters \ref{chap:lin}-\ref{chap:sp} we studied the canonical
commutation-relations for the quantum mechanical momentum and position
operators $P$, respectively $Q$. Below we outline how this problem
can be restated as a result about a unitary representation of the
matrix group $G_{3}$ of all upper triangular $3\times3$ matrices
over $\mathbb{R}$. This is a special unitary irreducible representation
$\mathcal{U}$ in $Rep(G_{3},L^{2}(\mathbb{R}))$. It is called the
Schrödinger representation, and the group $G_{3}$ is called the Heisenberg
group. We shall need the Stone-von Neumann uniqueness theorem, outlined
in the appendix below. (Also see \cite{vN32c,vN31}.) Its proof will
follow from a more general result which is included inside the present
chapter. The Stone-von Neumann uniqueness theorem states that every
unitary irreducible representation of $G_{3}$ is unitarily equivalent
to the Schrödinger representation.\index{representation!Schrödinger-}

We studied operators in Hilbert space of relevance to quantum physics.
A source of examples is relativistic physics. The symmetry group of
Einstein\textquoteright s theory is a particular Lie group, called
the Poincaré group. The study of its unitary representations is central
to relativistic physics. But it turns out that there is a host of
diverse applications (including harmonic analysis) where other groups
arise. Below we offer a glimpse of the theory of unitary representations,
and its many connections to operators in Hilbert space. \index{harmonic}
\index{representation!unitary}\index{Theorem!Stone-von Neumann uniqueness-}

Two pedantic points regarding unbounded operators. The first is the
distinction between \textquotedblleft selfadjoint\textquotedblright{}
vs \textquotedblleft essentially selfadjoint.\textquotedblright{}
An operator is said to be essentially selfadjoint if its closure is
selfadjoint. The second is the distinction between selfadjoint and
skewadjoint. The difference there is just a multiple of $i\left(=\sqrt{-1}\right)$.

\index{operators!essentially selfadjoint-}

\index{essentially selfadjoint operator} \index{selfadjoint operator}

These distinctions plays a role in the study of unitary representations
of Lie groups, but are often swept under the rug, especially in the
physics literature. Every unitary representations of a Lie group has
a derived representation of the corresponding Lie algebra. The individual
operators in a derived representation are skewadjoint; -- but to get
a common dense domain for all these operators, we must resort to essentially
skewadjointness. Nonetheless, indeed there are choices of common dense
domains (e.g., $C^{\infty}$-vectors), but the individual operators
in the derived representation will then only be essentially skewadjoint
there. \index{Lie!algebra}

For more details on this, see e.g., \cite{Pou72}.

\section{\label{sec:groupm}Motivation}

The following non-commutative Lie groups will be of special interest
to us because of their applications to physics, and to a host of areas
within mathematics; they are: the Heisenberg group $G=H_{3}$, the
$ax+b$ group $G=S_{2}$; and $SL_{2}\left(\mathbb{R}\right)$. In
outline:
\begin{itemize}
\item $G=H_{3}$, in real form $\simeq\mathbb{R}^{3}$, with multiplication
\begin{equation}
\left(a,b,c\right)\left(a',b',c'\right)=\left(a+a',b+b',c+c'+ab'\right),\label{eq:gh1}
\end{equation}
$\forall\left(a,b,c\right)$, and $\left(a',b',c'\right)\in\mathbb{R}^{3}$.
This is also matrix-multiplication when $\left(a,b,c\right)$ has
the form 
\[
\begin{pmatrix}1 & a & c\\
0 & 1 & b\\
0 & 0 & 1
\end{pmatrix}.
\]
In complex form, $z\in\mathbb{C}$, $c\in\mathbb{R}$, we have 
\begin{equation}
\left(z,c\right)\left(z',c'\right)=\left(z+z',c+c'+\Im\left(\overline{z}z'\right)\right).\label{eq:gh2}
\end{equation}

\item $G=S_{2}$, the group of transformation $x\mapsto ax+b$ where $a\in\mathbb{R}_{+}$,
and $b\in\mathbb{R}$, with multiplication
\begin{equation}
\left(a,b\right)\left(a',b'\right)=\left(aa',b+ab'\right).\label{eq:gh3}
\end{equation}
 This is also the matrix-multiplication when $\left(a,b\right)$ has
the form 
\[
\begin{pmatrix}a & b\\
0 & 1
\end{pmatrix}.
\]

\item $G=SL_{2}\left(\mathbb{R}\right)=2\times2$ matrices 
\[
\begin{pmatrix}a & b\\
c & d
\end{pmatrix}
\]
over $\mathbb{R}$, with $ad-bc=1$. Note that $SL_{2}\left(\mathbb{R}\right)$
is locally isomorphic to $SU\left(1,1\right)=$ $2\times2$ matrices
over $\mathbb{C}$, 
\[
\begin{pmatrix}\alpha & \beta\\
\overline{\beta} & \overline{\alpha}
\end{pmatrix}
\]
such that $\left|\alpha\right|^{2}-\left|\beta\right|^{2}=1$. In
both cases, the multiplication in $G$ is matrix-multiplication for
$2\times2$ matrices. 
\end{itemize}
The four groups, and their harmonic analysis will be studied in detail
inside this chapter.

An important question in the theory of unitary representations is
the following: For a given Lie group $G$, what are its irreducible
unitary representations, up to unitary equivalence. One aims for lists
of these irreducibles. The question is an important part of non-commutative
harmonic analysis. When answers are available, they have important
implications for physics and for a host of other applications, but
complete lists are hard to come by; and the literature on the subject
is vast. We refer to the book \cite{Tay86}, and its references for
an overview.

To begin with, the tools going into obtaining lists of the equivalence
classes of irreducible representations, differ from one class of Lie
groups to the other. Cases in point are the following classes, nilpotent,
solvable, and semisimple. The Heisenberg group is in the first class,
the $ax+b$ group in the second, and $SL_{2}\left(\mathbb{R}\right)$
in the third. By a theorem of Stone and von Neumann, the classes of
irreducibles for the Heisenberg group are indexed by a real parameter
$h$; they are infinite-dimensional for non-zero values of $h$, and
one dimensional for $h=0$.

For the $ax+b$ group, there are just two classes of unitary irreducibles.
The verification of this can be made with the use of Mackey's theory
of induced representations \cite{Mac88}. But the story is much more
subtle in the semisimple cases, even for $SL_{2}\left(\mathbb{R}\right)$.
The full list is divided up in series of representations (principal,
continuous, discrete, and complementary series representations), and
the paper \cite{JO00} outlines some of their properties. But the
details of this are far beyond the scope of the present book.

We now turn to some:
\begin{xca}[Semidirect product $G\small\textcircled{s}V$]
 \myexercise{Semidirect product $G\small\textcircled{s}V$}Let $V$
be a finite-dimensional vector space, and let $G\subset GL\left(V\right)$
be a sub-group of the corresponding general linear group: Make the
definition
\begin{equation}
\left(g,v\right)\left(g',v'\right)=\left(gg',g\left(v'\right)+v\right)\label{eq:sg1}
\end{equation}
for all $g,g'\in G$, and $v,v'\in V$. 
\begin{enumerate}
\item Show that with (\ref{eq:sg1}) we get a new group; called the semidirect
product .
\item In the group $G\small\textcircled{s}V$, show that 
\[
\left(g,v\right)^{-1}=\left(g^{-1},-g^{-1}\left(v\right)\right),\;g\in G,\:v\in V;
\]
and conclude from this that $V$ identifies as a normal subgroup in
$G\small\textcircled{s}V$. 
\item Show that, if $G$ is a Lie group, then so is $G\small\textcircled{s}V$. 
\item Show that, within the Lie algebra of $G\small\textcircled{s}V$, the
vector space $V$ identifies as an ideal.
\end{enumerate}
\end{xca}

\paragraph{General Considerations}

Every group $G$ is also a $*$-semigroup, with the $*$ operation
\[
g^{*}:=g^{-1}.
\]
$G$ is not a complex $*$ algebra yet, in particular, multiplication
by a complex number is not defined. As a general principal, $G$ can
be embedded into the $*$-algebra 
\[
\mathfrak{A}_{G}=G\otimes\mathbb{C}=\mathbb{C}\mbox{-valued functions on }G.
\]
$\mathfrak{A}_{G}$ has a natural pointwise multiplication and scalar
multiplication, given by 
\begin{alignat*}{1}
(g\otimes c_{g})(h\otimes c_{h}) & =gh\otimes c_{g}c_{h}\\
t(g\otimes c_{g}) & =g\otimes tc_{g}
\end{alignat*}
for all $g,h\in G$ and $c_{g},c_{h},t\in\mathbb{C}$. The $*$ operation
extends from $G$ to $\mathfrak{A}_{G}$ as
\begin{equation}
(g\otimes c_{g})^{*}=g^{-1}\otimes\overline{c_{g^{-1}}}.\label{eq:star}
\end{equation}

\begin{rem}
The $*$ operation so defined (as in (\ref{eq:star})) is the only
way to make it a period-2, conjugate linear, anti-automorphism. That
is, $\left(tA\right)^{*}=\overline{t}A^{*}$, $A^{**}=A$, and $(AB)^{*}=B^{*}A^{*}$,
for all $A,B\in\mathfrak{A}_{G}$, and all $t\in\mathbb{C}$. \index{representation!of group}\index{representation!of Lie group}\index{semigroup}
\end{rem}

There is a bijection between representations of $G$ and representations
of $\mathfrak{A}_{G}$. For if $\pi\in Rep(G,\mathscr{H})$, it then
extends to $\tilde{\pi}=\pi\otimes Id_{\mathbb{C}}\in Rep(\mathfrak{A}_{G},\mathscr{H})$,
by
\[
\tilde{\pi}(g\otimes c_{g})=\pi_{g}\otimes c_{g}.
\]
$id_{\mathbb{C}}$ denotes the identity representation $\mathbb{C}\rightarrow\mathbb{C}$.
Conversely, if $\rho\in Rep(\mathfrak{A}_{G},\mathscr{H})$ then 
\[
\pi=\rho\big|_{G\otimes1}
\]
is a representation of the group $G\simeq G\otimes1$. 
\begin{rem}
The notation $g\otimes c_{g}$ is usually written as $c_{g}$. Thus
pointwise multiplication takes the form 
\[
c_{g}d_{h}=l_{gh}.
\]
Equivalently, we write 
\[
l_{g}=\sum_{h}c_{h}d_{h^{-1}g}
\]
i.e., the usual convolution. In more details: 
\[
\left(\sum_{g}c_{g}\pi\left(g\right)\right)\left(\sum_{h}d_{h}\pi\left(h\right)\right)=\sum_{g}\underset{=l_{g}}{\underbrace{\left(\sum_{h}c_{h}d_{h^{-1}g}\right)}}\pi\left(g\right).
\]
The $*$ operation now becomes 
\[
c_{g}^{*}=\overline{c_{g^{-1}}}.
\]

\end{rem}
More generally, every locally compact\index{groups!locally compact}
group has a left (and right) Haar\index{measure!Haar} measure. Thus
the above construction has a ``continuous'' version. It turns out
that $\mathfrak{A}_{G}$ is the Banach $*$-algebra $L^{1}(G)$. Again,
there is a bijection between $Rep(G,\mathscr{H})$ and $Rep(L^{1}(G),\mathscr{H})$.
In particular, the ``discrete'' version is recovered if the measure
$\mu$ is discrete, in which case $\mathfrak{A}_{G}=l^{1}(G)$. \index{convolution}
\begin{defn}
Let $G$ be a group, and $\psi:G\rightarrow\mathbb{C}$ a function.
We say that $\psi$ is \emph{positive definite} iff (Def.) for all
$n\in\mathbb{N}$, all $g_{1},\ldots,g_{n}\in G$, and all $c_{1},\ldots,c_{n}\in\mathbb{C}$,
we have \index{positive definite!-function} 
\[
\sum_{j=1}^{n}\sum_{k=1}^{n}\overline{c_{j}}c_{k}\psi\left(g_{j}^{-1}g_{k}\right)\geq0.
\]
Often we also assume that $\psi\left(e\right)=1$. \end{defn}
\begin{xca}[Positive definite functions]
\myexercise{Positive definite functions}\label{exer:gp1}Show that
every positive definite function $\psi$ on $G$ extends by linearity
to a positive definite function $\widetilde{\psi}$ on the group algebra
$\mathbb{C}\left[G\right]$, i.e., 
\[
\widetilde{\psi}\left(\sum_{g}c_{g}g\right):=\sum_{g}c_{g}\psi\left(g\right)
\]
on all (finite) linear expressions $\sum_{g}c_{g}g$.
\end{xca}

\begin{xca}[Contractions]
\myexercise{Contractions}\label{exer:gp2}Let $\mathscr{H}$ be
a Hilbert space, and let $T$ be a \emph{contraction}, i.e., $T\in\mathscr{B}\left(\mathscr{H}\right)$,
satisfying one of the two equivalent conditions:
\[
\left\Vert T\right\Vert \leq1\Longleftrightarrow I-T^{*}T\geq0\;\text{in the order on Hermitian operators}.
\]

\begin{enumerate}
\item \label{enu:gp1}For $G=\mathbb{Z}$, define $\psi:\mathbb{Z}\rightarrow\mathbb{C}$
as follows:
\begin{equation}
\psi\left(n\right)=\begin{cases}
T^{n}=\underset{n\:\text{times}}{\underbrace{T\circ\cdots\circ T}} & \;\mbox{if}\;n\geq0\\
\left(T^{*}\right)^{\left|n\right|} & \;\mbox{if}\;n<0.
\end{cases}\label{eq:gpdf1}
\end{equation}
Show that $\psi$ is positive definite.\index{space!Hilbert-}
\item \label{enu:gp2}Conclude that $\widetilde{\psi}$ is completely positive
(see \chapref{cp}). 
\item Apply  (\ref{enu:gp1}) \& (\ref{enu:gp2}) to conclude the existence
of a triple $\left(V,\mathscr{K},\mathcal{U}\right)$, where $\mathscr{K}$
is a Hilbert space, $V:\mathscr{H}\rightarrow\mathscr{K}$ is isometric,
$\mathcal{U}:\mathscr{K}\rightarrow\mathscr{K}$ is a unitary operator;
and we have
\begin{equation}
T^{n}=V^{*}\mathcal{U}^{n}V,\;\forall n\in\mathbb{N}.\label{eq:gpdf2}
\end{equation}

\end{enumerate}

\emph{Historic Note}: The system $\left(V,\mathscr{K},\mathcal{U}\right)$
above, satisfying (\ref{eq:gpdf2}), is called a \emph{unitary dilation}\index{unitary dilation}
\cite{Sch55} (details in  \chapref{cp}.) 

\end{xca}

\begin{xca}[Central Extension]
\myexercise{Central Extension}Let $V$ be a vector space over $\mathbb{R}$,
and let $B:V\times V\rightarrow\mathbb{R}$ be a function. 
\begin{enumerate}
\item Then show that $V\times\mathbb{R}$ turns into a group $G$ when the
operation in $G$ is defined as follows:
\begin{equation}
\left(u,\alpha\right)\left(v,\beta\right)=\left(u+v,\alpha+\beta+B\left(u,v\right)\right),\;\forall\alpha,\beta\in\mathbb{R},\:\forall u,v\in V,\label{eq:ce1}
\end{equation}
if and only if $B$ satisfies
\begin{equation}
\begin{split}B\left(u,0\right)= & B\left(0,u\right)=0,\;\forall u\in V;\;\mbox{and}\\
B\left(u,v\right)+B\left(u+v,w\right)= & B\left(u,v+w\right)+B\left(v,w\right),\;\forall u,v,w\in V.
\end{split}
\label{eq:ce2}
\end{equation}

\item Assuming that (\ref{eq:ce2}) is satisfied; show that the group inverse
under the operation (\ref{eq:ce1}) is
\begin{equation}
\left(u,\alpha\right)^{-1}=\left(-u,-\alpha-B\left(u,-u\right)\right),\;\forall\alpha\in\mathbb{R},u\in V.\label{eq:ce3}
\end{equation}

\end{enumerate}
\end{xca}

\section{\label{sec:ureprep}Unitary One-Parameter Groups}
\begin{quotation}
``The mathematical landscape is full if groups of unitary operators.
The ..., strongly continuous one-parameter groups $U(t)$, $-\infty<t<\infty$,
come mostly from three sources: processes where energy is conserved,
such as those governed by wave equations of all sorts; process where
probability is preserved, for instance, ones governed by Schrödinger
equations; and Hamiltonian and other measure-preserving flows.''

--- Peter Lax, from \cite{MR1892228}\vspace{2em}
\end{quotation}
Let $\mathscr{H}$ be a Hilbert space, and let $P\left(\cdot\right)$
be a projection valued measure (PVM), \index{measure!projection-valued}\sindex[nam]{Lax, P.D., (1926 --)}
\begin{equation}
P:\mathcal{B}\left(\mathbb{R}\right)\longrightarrow Proj\left(\mathscr{H}\right)\label{eq:up1}
\end{equation}
i.e., defined on the sigma-algebra of all Borel subsets in $\mathbb{R}$.
\index{strongly continuous}

Then, as we saw, the integral (operator-valued)
\begin{equation}
U\left(t\right)=\int_{\mathbb{R}}e^{i\lambda t}P\left(d\lambda\right),\;t\in\mathbb{R}\label{eq:up2}
\end{equation}
is well defined, and yield a strongly continuous one-parameter group
acting on $\mathscr{H}$; equivalently, 
\begin{equation}
U\in Rep_{uni}\left(\mathbb{R},\mathscr{H}\right)\label{eq:up3}
\end{equation}
where $U$ is defined by (\ref{eq:up2}).

\index{unitary one-parameter group}

The theorem of M.H. Stone states that the converse holds as well,
i.e., every $U$, as in (\ref{eq:up3}), corresponds to a unique $P$,
a PVM, such that (\ref{eq:up2}) holds.

\index{von Neumann's ergodic theorem} \index{Theorem!Ergodic-}
\begin{xca}[von Neumann's ergodic theorem]
\label{exer:vNergodic}\myexercise{von Neumann's ergodic theorem}
Let $\left\{ U\left(t\right)\right\} _{t\in\mathbb{R}}$ be a strongly
continuous one-parameter group with PVM, $P\left(\cdot\right)$. Denote
by $P\left(\left\{ 0\right\} \right)$ the value of $P$ on the singleton
$\left\{ 0\right\} $. 
\begin{enumerate}
\item Show that 
\begin{equation}
P\left(\left\{ 0\right\} \right)\mathscr{H}=\left\{ h\in\mathscr{H}\::\:U\left(t\right)h=h,\;\forall t\in\mathbb{R}\right\} .\label{eq:up4}
\end{equation}
(We set $\mathscr{H}_{0}:=P\left(\left\{ 0\right\} \right)\mathscr{H}$.)
\item Establish the following limit conclusion:
\begin{equation}
\lim_{T\rightarrow\infty}\frac{1}{2T}\int_{-T}^{T}U\left(t\right)dt=P\left(\left\{ 0\right\} \right).\label{eq:up5}
\end{equation}
\uline{Hint:} 
\[
\frac{1}{2T}\int_{-T}^{T}U\left(t\right)dt=\int_{\mathbb{R}}\frac{\sin\left(\lambda T\right)}{\lambda T}P\left(d\lambda\right)
\]
holds for all $T\in\mathbb{R}_{+}$.
\end{enumerate}
\end{xca}

\section{Group - Algebra - Representations}

In physics, we are interested in representation of symmetry groups,
which preserve inner product or energy, and naturally leads to unitary
representations. Every unitary representation can be decomposed into
irreducible representations; the latter amounts to elementary particles
which can not be broken up further. In practice, quite a lot work
goes into finding irreducible representations of symmetry groups.
Everything we learned about algebras is also true for groups. The
idea is to go from groups to algebras and then to representations.

We summarize the basic definitions:

\index{algebras!group algebra}

\index{representation!unitary}
\begin{itemize}
\item $\pi_{G}\in Rep\left(G,\mathscr{H}\right)$
\[
\begin{cases}
\pi\left(g_{1}g_{2}\right) & =\pi\left(g_{1}\right)\pi\left(g_{2}\right)\\
\pi\left(e_{G}\right) & =I_{\mathscr{H}}\\
\pi\left(g^{-1}\right) & =\pi\left(g\right)^{*}
\end{cases}
\]

\item $\pi_{\mathfrak{A}}\in Rep\left(\mathfrak{A},\mathscr{H}\right)$
\[
\begin{cases}
\pi\left(A_{1}A_{2}\right) & =\pi\left(A_{1}\right)\pi\left(A_{2}\right)\\
\pi\left(1_{\mathfrak{A}}\right) & =I_{\mathscr{H}}\\
\pi\left(A^{*}\right) & =\pi\left(A\right)^{*}
\end{cases}
\]
\end{itemize}
\begin{casenv}
\item \begin{flushleft}
$G$ is discrete $\longrightarrow$ $\mathfrak{A}=G\otimes l^{1}$
\par\end{flushleft}
\end{casenv}
\begin{eqnarray*}
\left(\sum_{g}a(g)g\right)\left(\sum_{g}b(h)h\right) & = & \sum_{g,h}a(g)b(h)gh\\
 & = & \sum_{g'}\sum_{h}a(g'h^{-1})b(h)g'
\end{eqnarray*}
\[
\left(\sum_{g}c(g)g\right)^{*}=\left(\sum_{g}\overline{c(g^{-1})}g\right)
\]
where $c^{*}(g)=\overline{c(g^{-1})}$. The multiplication of functions
in $\mathfrak{A}$ is a generalization of convolutions.

\begin{casenv}
\item $G$ is locally compact $\longrightarrow$ $\mathfrak{A}=G\otimes L^{1}(\mu)\simeq L^{1}(G)$. \end{casenv}
\begin{defn}
Let $\mathcal{B}\left(G\right)$ be the Borel sigma-algebra of $G$.
A regular Borel measure $\lambda$ is said to be \emph{left (resp.
right) invariant}, if $\lambda\left(gE\right)=\lambda\left(E\right)$
(resp. $\lambda\left(Eg\right)=\lambda\left(Eg\right)$), for all
$g\in G$, and $E\in\mathcal{B}\left(G\right)$. $\lambda$ is called
a \emph{left (resp. right) Haar measure} accordingly.
\end{defn}
Note that $\lambda$ is left invariant iff 
\[
\lambda'\left(E\right):=\lambda\left(E^{-1}\right)
\]
is right invariant, where $E^{-1}=\left\{ g\in G:g^{-1}\in E\right\} $,
for all $E\in\mathcal{B}\left(G\right)$. Hence one may choose to
work with either a left or right invariant measure. 
\begin{thm}
\label{thm:haar}Every locally compact group $G$ has a left Haar
measure, unique up to a multiplicative constant. 
\end{thm}
For the existence of Haar measures, one first proves the easy case
when $G$ is compact, and then extends to locally compact cases. For
non compact groups, the left / right Haar measures could be different.
Many non compact groups have no Haar measure. In applications, the
Haar measures are usually constructed explicitly. 

\index{measure!Haar}

\index{compact}

\index{groups!compact}

\index{unimodular}

Given a left Haar measure $\lambda_{L}$, and $g\in G$, then 
\[
E\longmapsto\lambda_{L}\left(Eg\right),\;E\in\mathcal{B}\left(G\right)
\]
is also left invariant. Hence, by  \thmref{haar}, 
\begin{equation}
\lambda_{L}\left(Eg\right)=\triangle_{G}\left(g\right)\lambda_{L}\left(E\right)\label{eq:lm1}
\end{equation}
for some constant $\triangle_{G}\left(g\right)\in\mathbb{R}\backslash\left\{ 0\right\} $.
Note that $\triangle_{G}$ is well-defined, independent of the choice
of $\lambda_{L}$. Moreover, 
\begin{eqnarray*}
\lambda_{L}\left(Egh\right) & = & \triangle_{G}\left(h\right)\lambda_{L}\left(Eg\right)=\triangle_{G}\left(h\right)\triangle_{G}\left(g\right)\lambda_{L}\left(E\right)\\
\lambda_{L}\left(Egh\right) & = & \triangle_{G}\left(gh\right)\lambda_{L}\left(E\right)
\end{eqnarray*}
and it follows that $\triangle_{G}:G\longrightarrow\mathbb{R}_{\times}$
is a homomorphism, i.e., 
\begin{equation}
\triangle_{G}\left(gh\right)=\triangle_{G}\left(h\right)\triangle_{G}\left(g\right),\;\forall g,h\in G.\label{eq:lm2}
\end{equation}

\begin{defn}
$\triangle_{G}$ is called the \emph{modular function} of $G$. $G$
is said to be \emph{unimodular} if $\triangle_{G}\equiv1$. \index{modular function}\end{defn}
\begin{cor}
Every compact group $G$ is unimodular.\end{cor}
\begin{proof}
Since $\triangle_{G}$ is a homomorphism, $\triangle_{G}\left(G\right)$
is compact in $\mathbb{R}_{\times}$, so $\triangle_{G}\equiv1$.
\end{proof}

\begin{cor}
For all $f\in C_{c}\left(G\right)$, and $g\in G$, 
\begin{equation}
\int_{G}f\left(\cdot g\right)d\lambda_{L}=\triangle_{G}\left(g^{-1}\right)\int_{G}fd\lambda_{L}.\label{eq:hm1}
\end{equation}
Equivalently, we get the substitution formula:
\begin{equation}
d\lambda_{L}\left(\cdot g\right)=\triangle_{G}\left(g\right)d\lambda_{L}\left(\cdot\right).\label{eq:hm2}
\end{equation}

Similarly, 
\begin{equation}
\int_{G}f\left(g^{-1}\cdot\right)d\lambda_{R}=\triangle_{G}\left(g^{-1}\right)\int_{G}fd\lambda_{R};\label{eq:hm3}
\end{equation}
i.e., 
\begin{equation}
d\lambda_{R}\left(g\cdot\right)=\triangle_{G}\left(g^{-1}\right)d\lambda_{R}\left(\cdot\right).\label{eq:hm4}
\end{equation}
\end{cor}
\begin{proof}
It suffices to check for characteristic functions. Fix $E\in\mathcal{B}\left(G\right)$,
then 
\begin{eqnarray*}
\int_{G}\chi_{E}\left(\cdot g\right)d\lambda_{L} & = & \int_{G}\chi_{Eg^{-1}}d\lambda_{L}\\
 & = & \lambda_{L}\left(Eg^{-1}\right)\\
 & \underset{\left(\ref{eq:lm1}\right)}{=} & \triangle_{G}\left(g^{-1}\right)\lambda_{L}\left(E\right)\\
 & = & \triangle_{G}\left(g^{-1}\right)\int_{G}\chi_{E}d\lambda_{L}
\end{eqnarray*}
hence (\ref{eq:hm1})-(\ref{eq:hm2}) follow from this and a standard
approximation.

For the right Haar measure, recall that $E\mapsto\lambda_{L}\left(E^{-1}\right)$
is right invariant, and so $\lambda_{R}\left(E\right)=c\lambda_{L}\left(E^{-1}\right)$,
for some constant $c\in\mathbb{R}\backslash\left\{ 0\right\} $. (In
fact, more is true; see  \thmref{hm1} below.) Therefore,
\begin{eqnarray*}
\lambda_{R}\left(gE\right) & = & c\lambda_{L}\left(E^{-1}g^{-1}\right)\\
 & = & c\triangle_{G}\left(g^{-1}\right)\lambda_{L}\left(E^{-1}\right)\\
 & = & \triangle_{G}\left(g^{-1}\right)\lambda_{R}\left(E\right).
\end{eqnarray*}
This yields (\ref{eq:hm3})-(\ref{eq:hm4}).\end{proof}
\begin{thm}
\label{thm:hm1}Let $G$ be a locally compact group, then the two
Haar measures are mutually absolutely continuous, i.e., $\lambda_{L}\ll\lambda_{R}\ll\lambda_{L}$. 

Specifically, fix $\lambda_{L}$, and set 
\[
\lambda_{R}\left(E\right):=\lambda_{L}\left(E^{-1}\right),\;E\in\mathcal{B}\left(G\right);
\]
then 
\[
\frac{d\lambda_{R}}{d\lambda_{L}}\left(g\right)=\triangle_{G}\left(g^{-1}\right)=\mbox{Radon-Nikodym derivative.}
\]
\end{thm}
\begin{proof}
Note that $\triangle_{G}d\lambda_{R}$ is left invariant. Indeed,
\begin{eqnarray*}
\triangle_{G}\left(g\cdot\right)d\lambda_{R}\left(g\cdot\right) & = & \underset{\left(\ref{eq:lm2}\right)}{\underbrace{\left(\triangle_{G}\left(g\right)\triangle_{G}\left(\cdot\right)\right)}}\underset{\left(\ref{eq:hm4}\right)}{\underbrace{\left(\triangle_{G}\left(g^{-1}\right)d\lambda_{R}\left(\cdot\right)\right)}}\\
 & = & \triangle_{G}\left(\cdot\right)d\lambda_{R}\left(\cdot\right).
\end{eqnarray*}
Hence, by the uniqueness of the Haar measure, we have 

\[
\triangle_{G}d\lambda_{R}=c\,d\lambda_{L}
\]
for some constant $c\in\mathbb{R}\backslash\left\{ 0\right\} $. One
then checks that $c\equiv1$.
\end{proof}
\index{homomorphism}
\begin{cor}
If $\lambda_{L}$ is a left Haar measure on $G$, then 
\begin{equation}
d\lambda_{L}\left(g^{-1}\right)=\triangle_{G}\left(g^{-1}\right)d\lambda_{L}\left(g\right).\label{eq:lr1}
\end{equation}
Similarly, if $\lambda_{R}$ is a right Haar measure, then 
\begin{equation}
d\lambda_{R}\left(g^{-1}\right)=\triangle_{G}\left(g\right)d\lambda_{R}\left(g\right).\label{eq:lr2}
\end{equation}
\end{cor}
\begin{rem}
In the case $l^{1}(G)$, $\lambda_{L}=\lambda_{R}=$ the counting
measure, which is unimodular, hence $\triangle_{G}$ does not appear. 
\end{rem}
In $L^{1}(G)$, we define
\begin{eqnarray}
\left(\varphi\star\psi\right)\left(g\right) & := & \int_{G}\varphi\left(gh^{-1}\right)\psi\left(h\right)d\lambda_{R}\left(h\right)\label{eq:L1G}\\
 & = & \int_{G}\varphi\left(h^{-1}\right)\psi\left(hg\right)d\lambda_{R}\left(h\right)\nonumber \\
 & = & \int_{G}\varphi\left(h\right)\psi\left(h^{-1}g\right)\underset{d\lambda_{L}\left(h\right)}{\underbrace{\triangle_{G}\left(h\right)d\lambda_{R}\left(h\right)}}\nonumber 
\end{eqnarray}
and 
\begin{equation}
\varphi^{*}\left(g\right):=\overline{\varphi\left(g^{-1}\right)}\triangle_{G}\left(g\right).\label{eq:L1G1}
\end{equation}
The choice of (\ref{eq:L1G1}) preserves the $L^{1}$-norm. Indeed,
\begin{eqnarray*}
\int_{G}\left|\varphi^{*}\right|d\lambda_{R} & = & \int_{G}\left|\varphi\left(g^{-1}\right)\right|\triangle_{G}\left(g\right)d\lambda_{R}\left(g\right)\\
 & = & \int_{G}\left|\varphi\left(g\right)\right|\triangle_{G}\left(g^{-1}\right)d\lambda_{R}\left(g^{-1}\right)\\
 & = & \int_{G}\left|\varphi\right|d\lambda_{R}
\end{eqnarray*}
where $\triangle_{G}\left(g^{-1}\right)d\lambda_{R}\left(g^{-1}\right)=d\lambda_{R}\left(g\right)$
by (\ref{eq:lr2}).

$L^{1}(G)$ is a Banach {*}-algebra, and $L^{1}\left(G\right)=$ $L^{1}$-completion
of $C_{c}\left(G\right)$. (Fubini's theorem shows that $f\star g\in L^{1}(G)$,
for all $f,g\in L^{1}(G)$.)

We may also use left Haar measure in (\ref{eq:L1G}). Then, we set
\begin{eqnarray}
\left(\varphi\ast\psi\right)\left(g\right) & := & \int_{G}\varphi\left(h\right)\psi\left(h^{-1}g\right)d\lambda_{L}\left(h\right)\label{eq:L1G2}\\
 & = & \int_{G}\varphi\left(gh\right)\psi\left(h^{-1}\right)d\lambda_{L}\left(h\right)\nonumber \\
 & = & \int_{G}\varphi\left(gh^{-1}\right)\psi\left(h\right)\underset{d\lambda_{R}\left(h\right)}{\underbrace{\triangle_{G}\left(h^{-1}\right)d\lambda_{L}\left(h\right)}}\nonumber 
\end{eqnarray}
and set 
\begin{equation}
\varphi^{*}\left(g\right):=\overline{\varphi\left(g^{-1}\right)}\triangle_{G}\left(g^{-1}\right).\label{eq:L1G3}
\end{equation}

There is a bijection between representations of groups and representations
of algebras. 

Given a unitary representation $\pi\in Rep\left(G,\mathscr{H}\right)$,
let $dg$ denote the Haar measure in $L^{1}(G)$, then we get the
group algebra representation $\pi_{L^{1}\left(G\right)}\in Rep\left(L^{1}\left(G\right),\mathscr{H}\right)$,
where 
\begin{eqnarray*}
\pi_{L^{1}\left(G\right)}\left(\varphi\right) & = & \int_{G}\varphi\left(g\right)\pi\left(g\right)dg\\
\pi_{L^{1}\left(G\right)}\left(\varphi^{*}\right) & = & \pi_{L^{1}\left(G\right)}\left(\varphi\right)^{*}
\end{eqnarray*}
Indeed, one checks that 
\[
\pi_{L^{1}\left(G\right)}\left(\varphi_{1}\star\varphi_{2}\right)=\pi_{L^{1}\left(G\right)}\left(\varphi_{1}\right)\pi_{L^{1}\left(G\right)}\left(\varphi_{2}\right).
\]

Conversely, given a representation of $L^{1}(G)$, let $(\varphi_{i})$
be a sequence in $L^{1}$ such that $\varphi_{i}\rightarrow\delta_{g}$
. Then 
\[
\int\varphi_{i}(h)\pi(h)gdh\rightarrow\pi(g),
\]
i.e., the limit is a representation of $G$.

\index{representation!unitary}

\begin{rem}
Let $G$ be a matrix group, then $x^{-1}dx$, $x\in G$, is left translation
invariant. For if $y\in G$ , then 
\[
\left(yx\right)^{-1}d\left(yx\right)=x^{-1}\left(y^{-1}y\right)dx=x^{-1}dx.
\]
Now assume $\dim G=n$, and so $x^{-1}dx$ contains $n$ linearly
independent differential forms, $\sigma_{1},\ldots,\sigma_{n}$; and
each $\sigma_{j}$ is left translation invariant. Thus $\sigma_{1}\wedge\cdots\wedge\sigma_{n}$
is a left invariant volume form, i.e., the left Haar measure. Similarly,
the right Haar measure can be constructed from $dx\cdot x^{-1}$,
$x\in G$.

. 
\end{rem}

\subsection{Example -- $ax+b$ group}

\index{groups!$ax+b$}

\index{representation!of $ax+b$}

Let $G=\left\{ \left[\begin{array}{cc}
a & b\\
0 & 1
\end{array}\right]\::\:a\in\mathbb{R}_{+},\:b\in\mathbb{R}\right\} $. 
\begin{itemize}
\item Multiplication
\[
\left[\begin{array}{cc}
a' & b'\\
0 & 1
\end{array}\right]\left[\begin{array}{cc}
a & b\\
0 & 1
\end{array}\right]=\left[\begin{array}{cc}
a'a & a'b+b'\\
0 & 1
\end{array}\right]
\]

\item Inverse
\[
\left[\begin{array}{cc}
a & b\\
0 & 1
\end{array}\right]^{-1}=\left[\begin{array}{cc}
\frac{1}{a} & -\frac{b}{a}\\
0 & 1
\end{array}\right].
\]

\end{itemize}
$G$ is isomorphic to the transformation group$x\mapsto ax+b$; where
composition gives
\[
x\mapsto ax+b\mapsto a'\left(ax+b\right)+b'=aa'x+\left(a'b+b'\right).
\]

\begin{rem}
Setting $a=e^{t}$, $a'=e^{t'}$, $aa'=e^{t}e^{t'}=e^{t+t'}$, i.e.,
multiplication $aa'$ can be made into addition.
\end{rem}
The left Haar measure is given as follows: 

Let $g=\begin{bmatrix}a & b\\
0 & 1
\end{bmatrix}\in G$, so that 
\begin{eqnarray*}
g^{-1}dg & = & \frac{1}{a}\begin{bmatrix}1\; & -b\\
0\; & a
\end{bmatrix}\left[\begin{array}{cc}
da & db\\
0 & 0
\end{array}\right]\\
 & = & \frac{1}{a}\begin{bmatrix}da & db\\
0 & 0
\end{bmatrix}.
\end{eqnarray*}
Hence we get two left invariant (linear independent) differential
forms:
\[
\frac{da}{a}\quad\mbox{and}\quad\frac{db}{a}.
\]
Set 
\[
d\lambda_{L}\left(g\right)=d\lambda_{L}\left(x,y\right):=\frac{1}{x^{2}}dx\wedge dy;\;\left(g=\begin{bmatrix}x & y\\
0 & 1
\end{bmatrix},\:x\in\mathbb{R}_{+}\right).
\]

Indeed, $\lambda_{L}$ is left invariant. To check this, consider
\[
g=\left[\begin{array}{cc}
a & b\\
0 & 1
\end{array}\right],\:h=\left[\begin{array}{cc}
a' & b'\\
0 & 1
\end{array}\right],\:\mbox{and}
\]
\[
h^{-1}g=\left[\begin{array}{cc}
\frac{1}{a'} & -\frac{b'}{a'}\\
0 & 1
\end{array}\right]\left[\begin{array}{cc}
a & b\\
0 & 1
\end{array}\right]=\left[\begin{array}{cc}
\frac{a}{a'} & \frac{b-b'}{a'}\\
0 & 1
\end{array}\right];
\]
then 
\begin{eqnarray*}
\int_{G}f\left(h^{-1}g\right)d\lambda_{L}\left(g\right) & = & \int_{0}^{\infty}\int_{-\infty}^{\infty}f\left(\frac{a}{a'},\frac{b-b'}{a'}\right)\frac{da\wedge db}{a^{2}}\\
 & = & \int_{0}^{\infty}\int_{-\infty}^{\infty}f\left(s,t\right)\frac{d\left(a's\right)\wedge d\left(a't+b'\right)}{\left(a's\right)\left(a's\right)}\\
 & = & \int_{0}^{\infty}\int_{-\infty}^{\infty}f\left(s,t\right)\frac{ds\wedge dt}{s^{2}}
\end{eqnarray*}
where we set 
\[
s=\frac{a}{a'},\;da=a'ds
\]
\[
t=\frac{b-b'}{a'},\;db=a'dt
\]
so that
\[
\frac{da\wedge db}{a^{2}}=\frac{a'^{2}ds\wedge dt}{a'^{2}s^{2}}=\frac{ds\wedge dt}{s^{2}}.
\]

For the right Haar measure, note that 
\begin{eqnarray*}
\int f\left(gh^{-1}\right)\left(dg\right)g^{-1} & = & \int f\left(g'\right)d\left(g'h\right)\left(g'h\right)^{-1}\\
 & = & \int f\left(g'\right)\left(dg'\right)\left(hh^{-1}\right)g'^{-1}\\
 & = & \int f\left(g'\right)\left(dg'\right)g'^{-1}.
\end{eqnarray*}
Since 
\[
\left(dg\right)g^{-1}=\left[\begin{array}{cc}
da & db\\
0 & 0
\end{array}\right]\left[\begin{array}{cc}
\frac{1}{a} & -\frac{b}{a}\\
0 & 1
\end{array}\right]=\left[\begin{array}{cc}
\frac{da}{a} & -\frac{b\,da}{a}+db\\
0 & 0
\end{array}\right]
\]
we then set 
\[
d\lambda_{R}\left(g\right):=d\lambda_{R}\left(a,b\right)=\frac{da\wedge db}{a}.
\]
Check: 
\[
g=\left[\begin{array}{cc}
a & b\\
0 & 1
\end{array}\right],\:h=\left[\begin{array}{cc}
a' & b'\\
0 & 1
\end{array}\right],
\]
\[
gh^{-1}=\left[\begin{array}{cc}
a & b\\
0 & 1
\end{array}\right]\left[\begin{array}{cc}
\frac{1}{a'} & -\frac{b'}{a'}\\
0 & 1
\end{array}\right]=\left[\begin{array}{cc}
\frac{a}{a'} & -\frac{ab'}{a'}+b\\
0 & 1
\end{array}\right]
\]
and so 
\begin{eqnarray*}
\int_{G}f\left(gh^{-1}\right)d\lambda_{R}(g) & = & \int_{0}^{\infty}\int_{-\infty}^{\infty}f\left(\frac{a}{a'},-\frac{ab'}{a'}+b\right)\frac{da\wedge db}{a}\\
 & = & \int_{0}^{\infty}\int_{-\infty}^{\infty}f\left(s,t\right)\frac{a'ds\wedge dt}{a's}\\
 & = & \int_{0}^{\infty}\int_{-\infty}^{\infty}f\left(s,t\right)\frac{ds\wedge dt}{s}
\end{eqnarray*}
with a change of variable: 
\[
s=\frac{a}{a'},\;da=a'ds
\]
\[
t=-\frac{ab'}{a'}+b,\;db=dt
\]
\[
\frac{da\wedge db}{a}=\frac{\left(a'ds\right)\wedge dt}{a's}=\frac{ds\wedge dt}{s}.
\]

\section{\label{sec:ind}Induced Representations}

Most of the Lie groups considered here fall in the following class:

Let $V$ be a finite-dimensional vector space over $\mathbb{R}$ (or
$\mathbb{C}$). We will consider the real case here, but the modifications
needed for the complex case are straightforward. 

Let $q:V\times V\rightarrow\mathbb{R}$ be a non-degenerate bilinear
form, such that $v\rightarrow q\left(v,\cdot\right)\in V^{*}$ is
1-1. Let $GL\left(V\right)$ be the general linear group for $V$,
i.e., all invertible linear maps $V\rightarrow V$. (If a basis in
$V$ is chosen, this will be a matrix-group.)
\begin{lem}
Set 
\begin{equation}
G\left(q\right)=\left\{ g\in GL\left(V\right)\;\big|\;q\left(gu,gv\right)=q\left(u,v\right),\;\forall u,v\in V\right\} .\label{eq:gr1}
\end{equation}
Then $G\left(q\right)$ is a Lie group, and its Lie algebra consists
of all linear mappings $X:V\rightarrow V$ such that 
\begin{equation}
q\left(Xu,v\right)+q\left(u,Xv\right)=0,\;\forall u,v\in V.\label{eq:gr2}
\end{equation}
\end{lem}
\begin{proof}
Fix a linear mapping $X:V\rightarrow V$, and set 
\[
g_{X}\left(t\right)=\sum_{n=0}^{\infty}\frac{t^{n}}{n!}X^{n}=\exp\left(tX\right)
\]
i.e., the matrix-exponential. Note that $g_{X}\left(t\right)$ satisfies
(\ref{eq:gr1}) for all $t\in\mathbb{R}$ iff $X$ satisfies (\ref{eq:gr2}).
To see this, differentiate: i.e., compute 
\[
\frac{d}{dt}q\left(\exp\left(tX\right),\exp\left(tX\right)v\right)
\]
using that $q$ is assumed bilinear.
\end{proof}
We will address two questions:
\begin{enumerate}
\item How to induce a representation of a group $G$ from a representation
of the a subgroup $\Gamma\subset G$?
\item Given a representation of a group $G$, how to test whether it is
induced from a representation of a subgroup $\Gamma\subset G$? (See
\cite{Mac88}.)
\end{enumerate}
The main examples we will study are the Lie groups of \index{representation!induced}
\index{Lie!group}
\begin{itemize}
\item $ax+b$
\item Heisenberg
\item $SL_{2}(\mathbb{R})$
\item Lorentz
\item Poincaré
\end{itemize}
Among these, the $ax+b$, Heisenberg and Poincaré groups are semi-direct
product groups. Their representations are induced from normal subgroups. 

In more detail, the five groups in the list above are as follows:

The $ax+b$ group is the group of $2\times2$ matrices $\begin{pmatrix}a & b\\
0 & 1
\end{pmatrix}$ where $a\in\mathbb{R}_{+}$, and $b\in\mathbb{R}$. 

The Heisenberg group is the group of upper triangular $3\times3$
real matrices $\begin{pmatrix}1 & x & z\\
0 & 1 & y\\
0 & 0 & 1
\end{pmatrix}$, $x,y,z\in\mathbb{R}$. 

The group $SL_{2}\left(\mathbb{R}\right)$ is the group of all $2\times2$
matrices $\begin{pmatrix}a & b\\
c & d
\end{pmatrix}$ satisfying $a,b,c,d\in\mathbb{R}$, and $ad-bc=1$.

The Lorentz-group is the group $L=G\left(q\right)$ defined by (\ref{eq:gr1})
where $V=\mathbb{R}^{4}$, (space-times in physics) and
\[
q\left(x_{0},x_{1},x_{2},x_{3}\right)=-x_{0}^{2}+x_{1}^{2}+x_{2}^{2}+x_{3}^{2};
\]
so $q$ is the non-degenerate quadratic form with one minus sign,
and three plus signs.

The Poincaré group $P$ is the semi-direct product $P=L\small\textcircled{s}\mathbb{R}^{4}$,
where the group-product in $P$ is as follows:
\[
\left(g,v\right)\left(g',v'\right)=\left(gg',v+gv'\right)
\]
for all $g,g'\in L$, and all $v,v'\in\mathbb{R}^{4}$. 
\begin{xca}[The Heisenberg group as a semidirect product]
\myexercise{The Heisenberg group as a semidirect product}Consider
the following two subgroups $A$ and $B$ in the Heisenberg group
$H$:
\[
A=\begin{bmatrix}1 & x & 0\\
0 & 1 & 0\\
0 & 0 & 1
\end{bmatrix},\quad\mbox{and}\quad B=\begin{bmatrix}1 & 0 & z\\
0 & 1 & y\\
0 & 0 & 1
\end{bmatrix}.
\]

\begin{enumerate}
\item Verify that $A$ and $B$ are both Abelian subgroups under the matrix-multiplication
of $H$.
\item Show that $H$ becomes a semidirect product 
\[
H=A\small\textcircled{s}B
\]
where the action $\alpha_{x}$ of $A$, as a group of automorphisms
in $B$, is as follows:
\[
\alpha_{x}\left(y,z\right)=\left(y,z+xy\right),\quad\forall x,y,z\in\mathbb{R}.
\]

\end{enumerate}
\end{xca}

\begin{xca}[The invariant complex vector fields from the Heisenberg group]
\myexercise{The invariant complex vector fields from the Heisenberg group}In
its complex form, the Heisenberg group takes the form $\mathbb{C}\times\mathbb{R}$,
$z\in\mathbb{C}$, $c\in\mathbb{R}$; and with group multiplication:
\begin{equation}
\left(z,c\right)\left(z',c'\right)=\left(z+z',c+c'+2\Im\left(\overline{z}z'\right)\right)\label{eq:hc1}
\end{equation}
Set $\frac{\partial}{\partial z}=\frac{1}{2}\left(\frac{\partial}{\partial x}-i\frac{\partial}{\partial y}\right)$,
and $\frac{\partial}{\partial\overline{z}}=\frac{1}{2}\left(\frac{\partial}{\partial x}+i\frac{\partial}{\partial y}\right)$,
or in abbreviated form
\begin{equation}
\partial=\partial_{z}=\frac{1}{2}\left(\partial_{x}-i\partial_{y}\right),\quad\mbox{and}\quad\overline{\partial}=\overline{\partial}_{z}=\frac{1}{2}\left(\partial_{x}+i\partial_{y}\right),\label{eq:hc2}
\end{equation}
so 
\begin{equation}
4\left(\partial_{x}^{2}+\partial_{y}^{2}\right)=\partial\overline{\partial}=\overline{\partial}\partial.\label{eq:hc3}
\end{equation}
Show that a basis for the left-invariant vector fields on $H$ is
as follows:
\begin{equation}
\partial_{z}-i\overline{z}\partial_{c},\quad\overline{\partial}_{z}+iz\partial_{c},\quad\mbox{and}\quad i\partial_{c},\label{eq:hc4}
\end{equation}
with commutator
\begin{equation}
\left[\partial_{z}-i\overline{z}\partial_{c},\overline{\partial}_{z}+iz\partial_{c}\right]=2i\partial_{c}.\label{eq:hc5}
\end{equation}

\end{xca}
\textbf{Representation Theory.}

\begin{figure}
\includegraphics[width=0.5\textwidth]{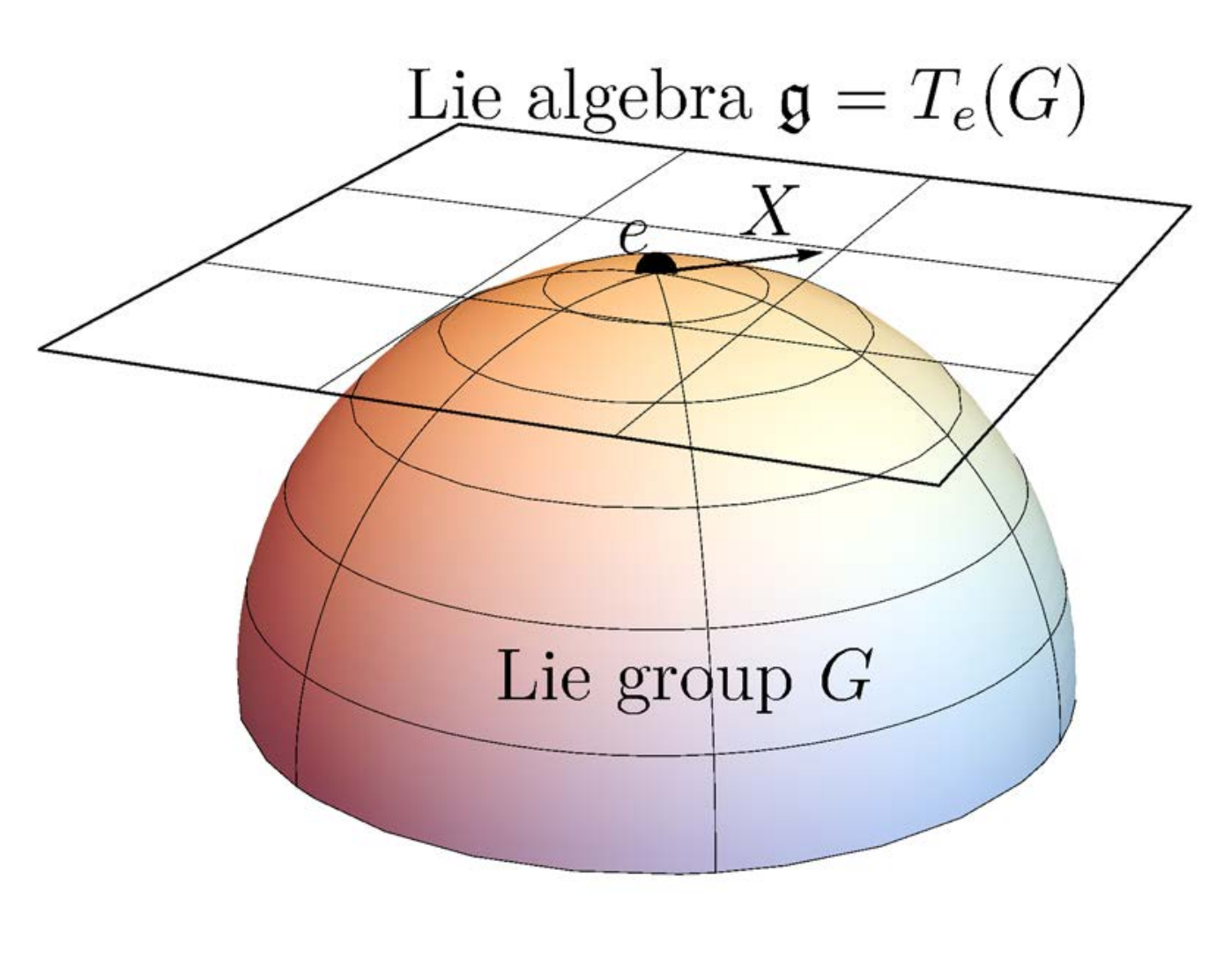}

\protect\caption{\label{fig:BCH}$G$ and $\mathfrak{g}$ (Lie algebra, Lie group,
and exponential mapping).}
\end{figure}

It is extremely easy to find representations of abelian subgroups.
Unitary irreducible representation of abelian subgroups are one-dimensional,
but the induced representation \cite{Mac88} on an enlarged Hilbert
space is infinite dimensional. \index{representation!unitary}
\begin{xca}[The Campbell-Baker-Hausdorff formula]
\myexercise{The Campbell-Baker-Hausdorff formula}\label{exer:gp3}The
exponential function is arguably the most important function in analysis.
The Campbell-Baker-Hausdorff formula (below) illustrates the role
of non-commutativity in this. 

Let $G$ and $\mathfrak{g}$ be as above, and let $\mathfrak{g}\xrightarrow{\;\exp\;}G$
be the exponential mapping. See \figref{BCH}.
\begin{enumerate}
\item Show that there is a convergent series with terms of degree $>1$
being iterated commutators $Z\left(X,Y\right)$ with \index{commutator}
\begin{equation}
\exp X\,\exp Y=\exp Z\left(X,Y\right),\;\mbox{and}\label{eq:com1}
\end{equation}
\begin{equation}
Z\left(X,Y\right)=X+Y+\frac{1}{2}\left[X,Y\right]+\frac{1}{12}\left(\left[X,\left[X,Y\right]\right]+\left[Y,\left[Y,X\right]\right]\right)+\cdots\label{eq:com2}
\end{equation}

\item Use combinatorics and algebra in order to derive an algorithm for
the terms ``$+\cdots$'' in (\ref{eq:com2}). This is the Baker--Campbell--Hausdorff
formula; see, e.g., \cite{HS68}.
\item Show that 
\[
Z\left(X,Y\right)+Z\left(-X,-Y\right)=0.
\]

\end{enumerate}
\end{xca}

\begin{xca}[The Lie algebra of a central extension]
\label{exer:lce}\myexercise{The Lie algebra of a central extension}
Let $V$ be a vector space over $\mathbb{R}$, and $B:V\times V\rightarrow\mathbb{R}$
a \emph{cocycle}, i.e., 
\begin{equation}
B\left(u,v\right)+B\left(u+v,w\right)=B\left(u,v+w\right)+B\left(v,w\right),\;\mbox{for }\forall u,v,w,\in V.\label{eq:lc1}
\end{equation}
Let $G_{B}=V\times\mathbb{R}$ be the corresponding Lie group:
\begin{equation}
\left(u,\alpha\right)\left(v,\beta\right)=\left(u+v,\alpha+\beta+B\left(u,v\right)\right),\;\mbox{for }\forall\alpha,\beta\in\mathbb{R},\forall u,v\in V.\label{eq:lc2}
\end{equation}

\begin{enumerate}
\item Show that the Lie algebra $La\left(G_{B}\right)$ of $G_{B}$ is $V\times\mathbb{R}$
itself with $0\times\mathbb{R}\subseteq$ the center; and with Lie
bracket $\left[\cdot,\cdot\right]$ given by
\begin{equation}
\left[u,v\right]=B\left(u,v\right)-B\left(v,u\right),\;\mbox{for }\forall u,v\in V.\label{eq:lc3}
\end{equation}

\item Show that the exponential mapping $\exp_{G_{B}}$ is the trivial mapping
\[
\exp_{G_{B}}\underset{\in La\left(G_{B}\right)}{\underbrace{\left(u,\alpha\right)}}=\underset{\in G_{B}}{\underbrace{\left(u,\alpha\right)}},\;\mbox{for }\forall u\in V,\alpha\in\mathbb{R}.
\]

\end{enumerate}
\end{xca}
\begin{example}[The $ax+b$ group ($a>0$)]
 $G=\{(a,b)\}$, where $(a,b)=\left[\begin{array}{cc}
a & b\\
0 & 1
\end{array}\right]$. The multiplication rule is given by 
\begin{eqnarray*}
(a,b)(a',b') & = & (aa',ab'+b)\\
(a,b)^{-1} & = & (\frac{1}{a},-\frac{b}{a}).
\end{eqnarray*}
$\Gamma=\{(1,b)\}$ is a one-dimensional abelian, normal subgroup
of $G$. We check that \index{groups!$ax+b$}
\begin{itemize}
\item abelian: $(1,b)(1,c)=(1,c+b)$
\item normal: $(x,y)(1,b)(x,y)^{-1}=(1,xb)$, note that this is also $Ad_{g}$
acting on the normal subgroup $\Gamma$
\item The other subgroup $\{(a,0)\}$ is isomorphic to the multiplicative
group $(\mathbb{R}_{+},\times)$. Because we have 
\[
(a,0)(a',0)=(aa',0)
\]
by the group multiplication rule above. 
\item Notice that $(\mathbb{R}_{+},\times)$ is not a normal subgroup, since
\[
(a,b)(x,0)(\frac{1}{a},-\frac{b}{a})=(ax,b)(\frac{1}{a},-\frac{b}{a})=(x-bx+b).
\]

\end{itemize}
\end{example}
$\Gamma$ is unimodular, hence it is just a copy of $\mathbb{R}$.
Its invariant measure is the Lebesgue measure on $\mathbb{R}$.

The multiplicative group $(\mathbb{R}_{+},\times)$ acts on the additive
group $(\mathbb{R},+)$ by
\begin{eqnarray*}
\varphi:(\mathbb{R}_{+},\times) & \mapsto & Aut((\mathbb{R},+))\\
\varphi_{a}(b) & = & ab
\end{eqnarray*}
check:
\[
(a,b)(a',b')=(aa',b+\varphi_{a}(b'))=(aa',b+ab')
\]
\[
(a,b)^{-1}=(a^{-1},\varphi_{a^{-1}}(b^{-1}))=(a^{-1},a^{-1}(-b))=(\frac{1}{a},-\frac{b}{a})
\]
\begin{eqnarray*}
(a,b)(1,x)(a,b^{-1}) & = & (a,b+\varphi_{a}(x))(a,b^{-1})\\
 & = & (a,b+ax)(\frac{1}{a},-\frac{b}{a})\\
 & = & (1,ax)=\varphi_{a}(x)
\end{eqnarray*}

\begin{example}
The Lie algebra of $G$ is given by $X=\left[\begin{array}{cc}
1 & 0\\
0 & 0
\end{array}\right]$, $=\left[\begin{array}{cc}
0 & 1\\
0 & 0
\end{array}\right]$. We check that
\[
e^{tX}=\left[\begin{array}{cc}
e^{t} & 0\\
0 & 1
\end{array}\right]
\]
which is subgroup $(\mathbb{R}_{+},\times)$; and 
\[
e^{sY}=I+sY+0+\cdots+0=\left[\begin{array}{cc}
1 & s\\
0 & 1
\end{array}\right]
\]
which is subgroup $(\mathbb{R},+)$. We also have $[X,Y]=Y$.\index{algebras!Lie algebra}
\end{example}

\begin{example}
Form $L^{2}(\mu_{L})$ where $\mu_{L}$ is the left Haar measure.
Then $\pi:g\rightarrow\pi(g)f(x)=f(g^{-1}x)$ is a unitary representation
in $L^{2}\left(\mu_{L}\right)$. Specifically, if $g=(a,b)$ then
\index{representation!unitary}
\[
f(g^{-1}x)=f(\frac{x}{a},\frac{y-b}{a}).
\]
Differentiate along the $a$ direction we get
\begin{eqnarray*}
\tilde{X}f & = & \frac{d}{da}\big|_{a=1,b=0}f(\frac{x}{a},\frac{y-b}{a})=(-x\frac{\partial}{\partial x}-y\frac{\partial}{\partial y})f(x,y)\\
\tilde{Y}f & = & \frac{d}{db}\big|_{a=1,b=0}f(\frac{x}{a},\frac{y-b}{a})=-\frac{\partial}{\partial y}f(x,y)
\end{eqnarray*}
therefore we have the vector field 
\begin{eqnarray*}
\tilde{X} & = & -x\frac{\partial}{\partial x}-y\frac{\partial}{\partial y}\\
\tilde{Y} & = & -\frac{\partial}{\partial y}
\end{eqnarray*}
or equivalently we get the Lie algebra representation $d\pi$ on $L^{2}(\mu_{L})$.
Notice that
\begin{eqnarray*}
[\tilde{X},\tilde{Y}] & = & \tilde{X}\tilde{Y}-\tilde{Y}\tilde{X}\\
 & = & (-x\frac{\partial}{\partial x}-y\frac{\partial}{\partial y})(-\frac{\partial}{\partial y})-(-\frac{\partial}{\partial y})(-x\frac{\partial}{\partial x}-y\frac{\partial}{\partial y})\\
 & = & x\frac{\partial^{2}}{\partial x\partial y}+y\frac{\partial^{2}}{\partial y^{2}}-(x\frac{\partial^{2}}{\partial x\partial y}+\frac{\partial}{\partial y}+y\frac{\partial^{2}}{\partial y^{2}})\\
 & = & -\frac{\partial}{\partial y}\\
 & = & \tilde{Y}.
\end{eqnarray*}
Notice that $\tilde{X}$ and $\tilde{Y}$ can be obtained by the exponential
map as well.
\begin{eqnarray*}
\tilde{X}f & = & \frac{d}{dt}\big|_{t=0}f(e^{-tX}x)\\
 & = & \frac{d}{dt}\big|_{t=0}f((e^{-t},1)(x,y))\\
 & = & \frac{d}{dt}\big|_{t=0}f(e^{-t}x,e^{-t}y+1)\\
 & = & (-x\frac{\partial}{\partial x}-y\frac{\partial}{\partial y})f(x,y)\\
\tilde{Yf} & = & \frac{d}{dt}\big|_{t=0}f(e^{-tY}x)\\
 & = & \frac{d}{dt}\big|_{t=0}f((1,-t)(x,y))\\
 & = & \frac{d}{dt}\big|_{t=0}f(x,y-t)\\
 & = & -\frac{\partial}{\partial y}f(x,y)
\end{eqnarray*}

\end{example}

\begin{example}
We may parametrize the Lie algebra of the $ax+b$ group using $(x,y)$
variables. Build the Hilbert space $L^{2}(\mu_{L})$. The unitary
representation $\pi(g)f(\sigma)=f(g^{-1}\sigma)$ induces the follows
representations of the Lie algebra\index{representation!Lie algebra}
\index{Lie!algebra}\index{groups!$ax+b$}
\begin{eqnarray*}
d\pi(\tilde{s})f(\sigma) & = & \frac{d}{dx}\big|_{s=0}f(e^{-sX}\sigma)=\tilde{X}f(\sigma)\\
d\pi(\tilde{t})f(\sigma) & = & \frac{d}{dy}\big|_{t=0}f(e^{-tY}\sigma)=\tilde{Y}f(\sigma).
\end{eqnarray*}
Hence in the parameter space $(s,t)\in\mathbb{R}^{2}$ we have two
usual derivative operators $\partial/\partial s$ and $\partial/\partial t$,
where on the manifold we have
\begin{eqnarray*}
\frac{\partial}{\partial s} & = & -x\frac{\partial}{\partial x}-y\frac{\partial}{\partial y}\\
\frac{\partial}{\partial t} & = & -\frac{\partial}{\partial y}
\end{eqnarray*}
The usual positive Laplacian on $\mathbb{R}^{2}$ translates to 
\begin{eqnarray*}
-\triangle & = & \left(\frac{\partial}{\partial s}\right)^{2}+\left(\frac{\partial}{\partial t}\right)^{2}\\
 & = & (\tilde{X})^{2}+(\tilde{Y})^{2}\\
 & = & \left(-x\frac{\partial}{\partial x}-y\frac{\partial}{\partial y}\right)\left(-x\frac{\partial}{\partial x}-y\frac{\partial}{\partial y}\right)+\left(-\frac{\partial}{\partial y}\right)^{2}\\
 & = & x^{2}\frac{\partial^{2}}{\partial x^{2}}+2xy\frac{\partial^{2}}{\partial x\partial y}+\left(y^{2}+1\right)\frac{\partial^{2}}{\partial y^{2}}+x\frac{\partial}{\partial x}+y\frac{\partial}{\partial y},
\end{eqnarray*}
where we used $\left(x\frac{\partial}{\partial x}\right)^{2}=x^{2}\left(\frac{\partial}{\partial x}\right)^{2}+x\frac{\partial}{\partial x}$.
This is in fact an elliptic operator, since the matrix
\[
\left[\begin{array}{cc}
x^{2} & xy\\
xy & y^{2}+1
\end{array}\right]
\]
has trace $trace=x^{2}+y^{2}+1\geq1$, and $\det=x^{2}\geq0$. If
instead we have ``$y^{2}$'' then the determinant is the constant
zero. \index{trace}

The term ``$y^{2}+1$'' is essential for $\triangle$ being elliptic.
Also note that all the coefficients are analytic functions in the
$(x,y)$ variables.
\end{example}

\begin{example}
\label{exa:hg}Heisenberg group $G=\{a,b,c\}$ where \index{groups!Heisenberg}
\[
(a,b,c)=\left[\begin{array}{ccc}
1 & a & c\\
0 & 1 & b\\
0 & 0 & 1
\end{array}\right]
\]
The multiplication rule is given by
\begin{eqnarray*}
(a,b,c)(a',b',c') & = & (a+a',b+b',c+ab'+c')\\
(a,b,c)^{-1} & = & (-a,-b,-c+ab)
\end{eqnarray*}
The subgroup $\Gamma=\{(0,b,c)\}$ where 
\[
\left(0,b,c\right)=\left[\begin{array}{ccc}
1 & 0 & c\\
0 & 1 & b\\
0 & 0 & 1
\end{array}\right]
\]
is two dimensional, abelian and normal.
\begin{itemize}
\item abelian: $(0,b,c)(0,b',c')=(0,b+b',c+c')$ 
\item normal: 
\begin{eqnarray*}
(a,b,c)(0,x,y)(a,b,c)^{-1} & = & (a,b,c)(0,x,y)(-a,-b,-c+ab)\\
 & = & (a,b+x,c+y+ax)(-a,-b,-c+ab)\\
 & = & (0,x,y+ax+ab-ab)\\
 & = & (0,x,ax+y)
\end{eqnarray*}
Note that this is also $Ad_{g}$ acting on the Lie algebra of $\Gamma$.
\end{itemize}
The additive group $(\mathbb{R},+)$ acts on $\Gamma=\{(0,b,c)\}\simeq(\mathbb{R}^{2},+)$
by
\begin{eqnarray*}
\varphi:(\mathbb{R},+) & \rightarrow & Aut(\Gamma)\\
\varphi(a)\left[\begin{array}{c}
c\\
b
\end{array}\right] & = & \left[\begin{array}{cc}
1 & a\\
0 & 1
\end{array}\right]\left[\begin{array}{c}
c\\
b
\end{array}\right]\\
 & = & \left[\begin{array}{c}
c+ab\\
b
\end{array}\right]
\end{eqnarray*}
check: 
\begin{eqnarray*}
(a,(b,c))(a',(b',c')) & = & (a+a',(b,c)+\varphi(a)(b',c'))\\
 & = & (a+a',(b,c)+(b',c'+ab'))\\
 & = & (a+a',b+b',c+c'+ab)\\
(a,(b,c))^{-1} & = & (-a,\varphi_{a^{-1}}(-b,-c))\\
 & = & (-a,(-b,-c+ab))\\
 & = & (-a,-b,-c+ab)\\
(a,b,c)(0,b',c')(a,b,c)^{-1} & = & (a,b+b',c+c'+ab')(-a,-b,-c+ab)\\
 & = & (0,b',c'+ab')\\
 & = & \varphi_{a}\left[\begin{array}{c}
c'\\
b'
\end{array}\right]
\end{eqnarray*}
\end{example}
\begin{xca}[Weyl's commutation relation]
\myexercise{Weyl's commutation relation}Let $\mathscr{H}$ be a
Hilbert space, and let $P,Q$ be a pair of selfadjoint operators such
that
\begin{equation}
e^{isP}e^{itQ}=e^{ist}e^{itQ}e^{isP}\label{eq:wc1}
\end{equation}
holds for all $s,t\in\mathbb{R}$; then prove that $P$ and $Q$ have
a common dense invariant domain $\mathscr{D}$ such that $P=\overline{P\big|_{\mathscr{D}}}$,
$Q=\overline{Q\big|_{\mathscr{D}}}$; and
\begin{equation}
\left[P,Q\right]=-i\,I\label{eq:wc2}
\end{equation}
holds on $\mathscr{D}$; more precisely, 
\begin{equation}
PQ\varphi-QP\varphi=-i\,\varphi,\label{eq:wc3}
\end{equation}
holds for all $\varphi\in\mathscr{D}$.

\uline{Hint}. One way of proving this is to show that when $P$
and $Q$ satisfy (\ref{eq:wc1}), then we automatically get a unitary
representation of the Heisenberg group (see Sections \ref{sec:groupm},
\ref{sec:ind}, and \exaref{hg}), and for $\mathscr{D}$ we can take
the corresponding Gårding Space (\secref{gaarding}). However there
is also a direct proof from first principles.

\uline{Caution}. The converse implication does \emph{not} hold:
There are selfadjoint solutions to (\ref{eq:wc3}) which do not have
a counterpart relation (\ref{eq:wc1}); see \cite{JM84,MR3210626,JoMu80}.
(Eq (\ref{eq:wc1}) is called the Weyl relation.)
\end{xca}

\subsection{Induced Representations}

This also goes under the name of ``Mackey machine'' \cite{Mac88}.
Its modern formulation is in the context of completely positive map.
\index{completely positive map}\index{modular function}\index{groups!locally compact}
\index{Mackey machine}

Let $G$ be a locally compact group, and $\Gamma\subset G$ be a closed
subgroup. Let $dx$ (resp. $d\xi$) be the right Haar measure on $G$
(resp. $\Gamma$), and $\triangle$ (resp. $\delta$) be the corresponding
modular function. Recall the modular function comes in when the translation
is put on the wrong side, i.e., 
\[
\int_{G}f\left(g^{-1}x\right)dx=\triangle\left(g^{-1}\right)\int_{G}f\left(x\right)dx
\]
or equivalently,
\[
\triangle\left(g\right)\int_{G}f\left(g^{-1}x\right)dx=\int_{G}f\left(x\right)dx.
\]
Form the quotient $M=\Gamma\backslash G$ space, and let $\pi:G\rightarrow\Gamma\backslash G$
be the quotient map (the covering map). $M$ carries a transitive
$G$ action.

\begin{longtable}{c|c|c}
\hline 
group & right Haar measure & modular function\tabularnewline
\hline 
$G$ & $dg$ & $\triangle$\tabularnewline
\hline 
$\Gamma$ & $d\xi$ & $\delta$\tabularnewline
\hline 
\end{longtable}
\begin{note}
$M$ is called \emph{a fundamental domain} or a \emph{homogeneous}
space. $M$ is a group if and only if $\Gamma$ is a normal subgroup
in $G$. In general, $M$ may not be a group, but it is still a very
important manifold.
\end{note}

\begin{note}
$\mu$ is an \emph{invariant measure} on $M$, if $\mu\left(Eg\right)=\mu\left(E\right)$,
$\forall g\in G$. $\mu$ is \emph{quasi-invariant}, if $\mu(E)=0$
$\Leftrightarrow$ $\mu(Eg)=0$, $\forall g$. In general there is
no invariant measures on $M$, but only quasi-invariant measures. 

$G$ has an invariant measure if and only if $G$ is unimodular (e.g.
Heisenberg group.) Not all groups are unimodular. A typical example
is the $ax+b$ group.

\index{measure!quasi-invariant}\index{unimodular}
\end{note}
Define $\tau:C_{c}\left(G\right)\rightarrow C_{c}\left(M\right)$
by
\begin{equation}
\left(\tau\varphi\right)\left(\pi\left(x\right)\right)=\int_{\Gamma}\varphi\left(\xi x\right)d\xi.\label{eq:cexp}
\end{equation}

\begin{note}
Since $\varphi$ has compact support, the integral in (\ref{eq:cexp})
is well-defined. $\tau$ is called \emph{conditional expectation}.
It is the summation of $\varphi$ over the orbit $\Gamma x$. Indeed,
for fixed $x$, if $\xi$ runs over $\Gamma$ then $\xi x$ runs over
$\Gamma x$. We may also say $\tau\varphi$ is a $\Gamma$-periodic
extension, by looking at it as a function defined on $G$. For if
$\xi_{1}\in\Gamma$, we have 
\[
\left(\tau\varphi\right)\left(\xi_{1}x\right)=\int_{\Gamma}\varphi\left(\xi\xi_{1}x\right)d\xi=\tau\varphi\left(x\right)
\]
using the fact that $d\xi$ right-invariant. Thus $\tau\varphi$,
viewed as a function on $G$, is $\Gamma$-periodic, i.e., $\left(\tau\varphi\right)\left(\xi x\right)=\left(\tau\varphi\right)\left(x\right)$,
$\forall\xi\in\Gamma$. \end{note}
\begin{lem}
\label{lem:tau}$\tau$ is surjective.\end{lem}
\begin{proof}
Suppose $f$ is $\Gamma$-periodic, choose $\psi\in C_{c}\left(G\right)$
such that $\tau\psi\equiv1$. Then $\psi f\in C_{c}\left(G\right)$,
and 
\begin{eqnarray*}
\left(\tau\left(\psi f\right)\right)\left(x\right) & = & \int_{\Gamma}\psi\left(\xi x\right)f\left(\xi x\right)d\xi\\
 & = & \int_{\Gamma}\psi\left(\xi x\right)f\left(x\right)d\xi\\
 & = & f\left(x\right)\int_{\Gamma}\psi\left(\xi x\right)d\xi\\
 & = & f\left(x\right)\left(\tau\psi\right)\left(x\right)=f\left(x\right).
\end{eqnarray*}
\end{proof}
\begin{example}
For $G=\mathbb{R}$, $\Gamma=\mathbb{Z}$, $d\xi$ = counting measure
on $\mathbb{Z}$, we have 
\[
\left(\tau\varphi\right)\left(\pi\left(x\right)\right)=\int_{\Gamma}\varphi\left(\xi x\right)d\xi=\sum_{n\in\mathbb{Z}}\varphi\left(n+x\right),\quad\forall\varphi\in C_{c}\left(\mathbb{R}\right).
\]
Since $\varphi$ has compact support, $\varphi\left(n+x\right)$ vanishes
for all but a finite number of $n$. Hence $\tau\varphi$ contains
a finite summation, and so it is well-defined. Moreover, for all $n_{0}\in\mathbb{Z}$,
it follows that 
\[
\left(\tau\varphi\right)\left(n_{0}+x\right)=\sum_{n\in\mathbb{Z}}\varphi\left(n_{0}+n+x\right)=\sum_{n\in\mathbb{Z}}\varphi\left(n+x\right)=\left(\tau\varphi\right)\left(x\right).
\]
Hence $\tau\varphi$ is translation invariant by integers, i.e., $\tau\varphi$
(as a function on $\mathbb{R}$) is $\mathbb{Z}$-periodic. 
\end{example}
Let $L:\Gamma\rightarrow V$ be a unitary representation of $\Gamma$
on a Hilbert space $V$. We now construct a unitary representation
$ind_{\Gamma}^{G}:G\rightarrow\mathscr{H}$ of $G$ on an enlarged
Hilbert space $\mathscr{H}$. 

Let $F_{*}$ be the set of function $f:G\rightarrow V$ so that
\begin{equation}
f\left(\xi g\right)=\rho\left(\xi\right){}^{1/2}L_{\xi}f\left(g\right),\quad\forall\xi\in\Gamma\label{eq:Fstar}
\end{equation}
where $\rho=\frac{\delta}{\Delta}$, i.e., $\rho\left(\xi\right)=\frac{\delta\left(\xi\right)}{\Delta\left(\xi\right)}$,
$\forall\xi\in\Gamma$. For all $f\in F_{*}$, let 
\[
(R_{g}f)(\cdot):=f(\cdot g)
\]
be the right-translation of $f$ by $g\in G$. 
\begin{lem}
$R_{g}f\in F_{*}$. That is, $F_{*}$ is invariant under right-translation
by $g\in G$. \end{lem}
\begin{proof}
To see this, let $f\in F_{*}$, $\xi\in\Gamma$, then 
\[
(R_{g}f)(\xi x)=f(\xi xg)=\rho(\xi)^{1/2}L_{\xi}f(xg)=\rho(\xi)^{1/2}L_{\xi}(R_{g}f)(x)
\]
so that $R_{g}f\in F_{*}$. \end{proof}
\begin{note}
We will defined an inner product on $F_{*}$ so that $\left\Vert f(\xi\cdot)\right\Vert _{new}=\left\Vert f(\cdot)\right\Vert _{new}$,
$\forall\xi\in\Gamma$. Eventually, we will define the induced representation
$U^{ind}:=ind_{\Gamma}^{G}\left(L\right)$ by 
\[
\left(U_{g}^{ind}f\right)\left(\cdot\right):=\left(R_{g}f\right)\left(\cdot\right)
\]
not on $F_{*}$, but pass to a quotient space. The factor $\rho(\xi)^{1/2}$
comes in as we are going to construct a quasi-invariant measure on
$\Gamma\backslash G$. 

To construct $F_{*}$, let $\varphi\in C_{c}\left(G\right)$, and
set 
\[
f\left(g\right):=\int_{\Gamma}\rho^{1/2}\left(\xi^{-1}\right)L\left(\xi^{-1}\right)f\left(\xi g\right)d\xi.
\]
Now if $\xi_{1}\in\Gamma$ then 
\begin{eqnarray*}
f\left(\xi_{1}g\right) & = & \int_{\Gamma}\rho^{1/2}\left(\xi^{-1}\right)L\left(\xi^{-1}\right)f\left(\xi\xi_{1}g\right)d\xi\\
 & = & \int_{\Gamma}\rho^{1/2}\left(\left(\xi\xi_{1}^{-1}\right)^{-1}\right)L\left(\left(\xi\xi_{1}{}^{-1}\right)^{-1}\right)f\left(\xi g\right)d\xi\\
 & = & \rho^{1/2}\left(\xi_{1}\right)L\left(\xi_{1}\right)\int_{\Gamma}\rho^{1/2}\left(\xi^{-1}\right)L\left(\xi^{-1}\right)f\left(\xi g\right)d\xi\\
 & = & \rho^{1/2}\left(\xi_{1}\right)L\left(\xi_{1}\right)f\left(g\right).
\end{eqnarray*}
The proof of \lemref{tau} shows that all functions in $F_{*}$ are
obtained this way. 
\end{note}
\index{representation!unitary}

\begin{note}
Let's ignore the factor $\rho(\xi)^{1/2}$ for a moment. $L_{\xi}$
is unitary implies that for all $f\in F_{*}$, 
\[
\left\Vert f(\xi\cdot)\right\Vert _{V}=\left\Vert L_{\xi}f(\cdot)\right\Vert _{V}=\left\Vert f(\cdot)\right\Vert _{V},\quad\forall\xi\in\Gamma.
\]
Since Hilbert spaces exist up to unitary equivalence, $L_{\xi}f(g)$
and $f(g)$ really are the same function. As $\xi$ running through
$\Gamma$, $\xi g$ running through $\Gamma g$. Thus $\left\Vert f(\xi g)\right\Vert $
is a constant on the orbit $\Gamma g$. It follows that $f(\xi g)$
is in fact a $V$-valued function defined on the quotient $M=\Gamma\backslash G$
(i.e., quasi-$\Gamma$-periodic). We will later use these functions
as multiplication operators.

\index{unitary equivalence}\end{note}
\begin{example}
The Heisenberg group is unimodular, so $\rho\equiv1$. \index{unimodular}
\end{example}

\begin{example}
For the $ax+b$ group, 
\begin{eqnarray*}
d\lambda_{R} & = & \frac{dadb}{a}\\
d\lambda_{L} & = & \frac{dadb}{a^{2}}\\
\triangle & = & \frac{d\lambda_{L}}{d\lambda_{R}}=\frac{1}{a}
\end{eqnarray*}
On the abelian normal subgroup $\Gamma=\{(1,b)\}$, we have $a=1$
and $\triangle(\xi)=1$. $\Gamma$ is unimodular, $\delta(\xi)=1$.
Therefore, $\rho(\xi)=1$, $\forall\xi\in\Gamma$.
\end{example}

For all $f\in F_{*}$, the map $\mu_{f,f}:C_{c}(M)\rightarrow\mathbb{C}$
given by 
\begin{equation}
\mu_{f,f}:\tau\varphi\longmapsto\int_{G}\left\Vert f\left(g\right)\right\Vert _{V}^{2}\varphi\left(g\right)dg,\quad\varphi\in C_{c}\left(G\right)\label{eq:Mlin}
\end{equation}
is a positive linear functional. By Riesz's theorem, there exists
a unique Radon measure $\mu_{f,f}$ on $M$, such that \index{functional}
\index{Riesz' theorem}\index{Theorem!Riesz-} 
\[
\int_{G}\left\Vert f(g)\right\Vert _{V}^{2}\varphi(g)dg=\int_{M}(\tau\varphi)d\mu_{f,f}.
\]

\begin{lem}
(\ref{eq:Mlin}) is a well-defined positive linear functional.
\end{lem}
\allowdisplaybreaks
\begin{proof}
Suppose $\varphi\in C_{c}\left(G\right)$ such that $\tau\varphi\equiv0$
on $M$. It remains to verify that $\mu_{f,f}\left(\tau\varphi\right)=0$,
see (\ref{eq:Mlin}). For this, we choose $\psi\in C_{c}\left(G\right)$
such that $\tau\psi\equiv1$ on $M$, and so 
\begin{eqnarray*}
\int_{G}\left\Vert f\left(g\right)\right\Vert _{V}^{2}\varphi\left(g\right)dg & = & \int_{G}\left\Vert f\left(g\right)\right\Vert _{V}^{2}\left(\tau\psi\right)\left(\pi\left(g\right)\right)\varphi\left(g\right)dg\\
 & = & \int_{G}\left\Vert f\left(g\right)\right\Vert _{V}^{2}\varphi\left(g\right)\left(\int_{\Gamma}\psi\left(\xi g\right)d\xi\right)dg\\
 & \underset{\text{Fubini}}{=} & \int_{\Gamma}\left(\int_{G}\left\Vert f\left(g\right)\right\Vert _{V}^{2}\varphi\left(g\right)\psi\left(\xi g\right)dg\right)d\xi\\
 & = & \int_{\Gamma}\left(\int_{G}\left\Vert f\left(\xi^{-1}g\right)\right\Vert _{V}^{2}\varphi\left(\xi^{-1}g\right)\psi\left(g\right)\triangle\left(\xi\right)dg\right)d\xi\\
 & \underset{\text{Fubini}}{=} & \int_{G}\psi\left(g\right)\left(\int_{\Gamma}\left\Vert f\left(\xi^{-1}g\right)\right\Vert _{V}^{2}\varphi\left(\xi^{-1}g\right)\triangle\left(\xi\right)d\xi\right)dg\\
 & = & \int_{G}\psi\left(g\right)\left(\int_{\Gamma}\left\Vert f\left(\xi g\right)\right\Vert _{V}^{2}\varphi\left(\xi g\right)\triangle\left(\xi^{-1}\right)\delta\left(\xi\right)d\xi\right)dg\\
 & \underset{\left(\ref{eq:Fstar}\right)}{=} & \int_{G}\psi\left(g\right)\left\Vert f\left(g\right)\right\Vert _{V}^{2}\left(\int_{\Gamma}\varphi\left(\xi g\right)d\xi\right)dg\\
 & = & \int_{G}\psi\left(g\right)\left\Vert f\left(g\right)\right\Vert _{V}^{2}\underset{\equiv0}{\underbrace{\left(\tau\varphi\right)\left(\pi\left(g\right)\right)}}dg=0.
\end{eqnarray*}
\end{proof}
\begin{note}
Recall that given a measure space $(X,\mathfrak{M},\mu)$, let $f:X\rightarrow Y$.
Define a linear functional $\Lambda:C_{c}(Y)\rightarrow\mathbb{C}$
by 
\[
\Lambda\varphi:=\int\varphi(f(x))d\mu(x)
\]
$\Lambda$ is positive, hence by Riesz's theorem, there exists a unique
regular Borel\index{measure!Borel} measure $\mu_{f}$ on $Y$ so
that
\[
\Lambda\varphi=\int_{Y}\varphi d\mu_{f}=\int_{X}\varphi(f(x))d\mu(x).
\]
It follows that $\mu_{f}=\mu\circ f^{-1}$. 
\end{note}

\begin{note}
Under current setting, we have a covering map $\pi:G\rightarrow\Gamma\backslash G=:M$,
and the right Haar measure $\mu$ on $G$. Thus we may define a measure
$\mu\circ\pi^{-1}$. However, given $\varphi\in C_{c}(M)$, $\varphi(\pi(x))$
may not have compact support, or equivalently, $\pi^{-1}(E)$ is $\Gamma$
periodic. For example, take $G=\mathbb{R}$, $\Gamma=\mathbb{Z}$,
$M=\mathbb{Z}\backslash\mathbb{R}$. Then $\pi^{-1}([0,1/2))$ is
$\mathbb{Z}$-periodic, which has infinite Lebesgue measure. What
we really need is some map so that the inverse of a subset of $M$
is restricted to a single $\Gamma$ period. This is essentially what
$\tau$ does: from $\tau\varphi\in C_{c}(M)$, get the inverse image
$\varphi\in C_{c}(G)$. Even if $\varphi$ is not restricted to a
single $\Gamma$ period, $\varphi$ always has compact support. 
\end{note}

Hence we get a family of measures indexed by elements in $F_{*}$.
If choosing $f,g\in F_{*}$ then we get complex measures $\mu_{f,g}$
(using polarization identity.)
\begin{itemize}
\item Define $\left\Vert f\right\Vert ^{2}:=\mu_{f,f}(M)$, $\left\langle f,g\right\rangle :=\mu_{f,g}\left(M\right)$. 
\item Complete $F_{*}$ with respect to this norm to get an enlarged Hilbert
space $\mathscr{H}$.
\item Define the induced representation $U^{ind}:=ind_{\Gamma}^{G}\left(L\right)$
on $\mathscr{H}$ as
\[
\left(U_{g}^{ind}f\right)\left(x\right)=f\left(xg\right)
\]
$U^{ind}$ is unitary, in particular, 
\[
\Vert U_{g}^{ind}f\Vert_{\mathscr{H}}=\left\Vert f\right\Vert _{\mathscr{H}},\quad\forall g\in G.
\]
\end{itemize}
\begin{note}
$\mu_{f,g}\left(M\right)=\int_{M}\left(\tau\varphi\right)\left(\xi\right)d\xi$
with $\tau\varphi\equiv1$. What is $\varphi$ then? It turns out
that $\varphi$ could be constant 1 over a single $\Gamma$-period,
or equivalently, $\varphi$ could spread out to a finite number of
$\Gamma$-periods. In the former case, 
\begin{eqnarray*}
\left\Vert f\right\Vert ^{2} & = & \int_{G}\left\Vert f\left(g\right)\right\Vert _{V}^{2}\varphi\left(g\right)dg\\
 & = & \int_{\text{1-period}}\left\Vert f\left(g\right)\right\Vert _{V}^{2}\varphi\left(g\right)dg\\
 & = & \int_{\text{1-period}}\left\Vert f\left(g\right)\right\Vert _{V}^{2}dg\\
 & = & \int_{M}\left\Vert f\left(g\right)\right\Vert _{V}^{2}dg.
\end{eqnarray*}

\end{note}

Define $P(\psi)f(x):=\psi(\pi(x))f(x)$, for $\psi\in C_{c}(M)$,
$f\in\mathscr{H}$, $x\in G$. Note $\left\{ P(\psi)\:|\:\psi\in C_{c}\left(M\right)\right\} $
is the abelian algebra of multiplication operators. 
\begin{lem}
We have
\[
U_{g}^{ind}P(\psi)U_{g^{-1}}^{ind}=P(\psi(\cdot g)).
\]
\end{lem}
\begin{proof}
One checks that 
\begin{eqnarray*}
U_{g}^{ind}P(\psi)f(x) & = & U_{g}^{ind}\psi(\pi(x))f(x)\\
 & = & \psi(\pi(xg))f(xg)\\
P(\psi(\cdot g))U_{g}^{ind}f(x) & = & P(\psi(\cdot g))f(xg)\\
 & = & \psi(\pi(xg))f(xg).
\end{eqnarray*}

\end{proof}
Conversely, how to recognize induced representations? Answer: 
\begin{thm}[Imprimitivity \cite{Ors79}]
 \label{thm:imp}Let $G$ be a locally compact group with a closed
subgroup $\Gamma$. Let $M=\Gamma\backslash G$. Suppose the system
$(U,P)$ satisfies the covariance relation
\[
U_{g}P(\psi)U_{g^{-1}}=P(\psi(\cdot g)),
\]
and $P\left(\cdot\right)$ is non-degenerate. Then, there exists a
unitary representation $L\in Rep(\Gamma,V)$ such that $U\cong ind_{\Gamma}^{G}\left(L\right)$. \end{thm}
\begin{rem}
$P\left(\cdot\right)$ is \emph{non-degenerate} if $P\left(C_{c}\left(M\right)\right)\mathscr{H}=\left\{ P\left(\psi\right)a:\psi\in C_{c}\left(M\right),a\in\mathscr{H}\right\} $
is dense in $\mathscr{H}$. 
\end{rem}
\index{groups!locally compact}\index{representation!unitary} \index{Theorem!imprimitivity-}

\section{\label{sec:Heisenberg}Example - Heisenberg group}

\index{groups!Heisenberg} 

\index{representation!Heisenberg}

\index{Schrödinger representation}

Let $G=\{(a,b,c)\}$ be the \emph{Heisenberg group}, where
\[
(a,b,c)=\left[\begin{array}{ccc}
1 & a & c\\
0 & 1 & b\\
0 & 0 & 1
\end{array}\right]
\]
The multiplication rule is given by
\begin{eqnarray*}
(a,b,c)(a',b',c') & = & (a+a',b+b',c+c'+ab')\\
(a,b,c)^{-1} & = & (-a,-b,-c+ab)
\end{eqnarray*}
The subgroup $\Gamma=\{(0,b,c)\}$ where 
\[
(1,b,c)=\left[\begin{array}{ccc}
1 & 0 & c\\
0 & 1 & b\\
0 & 0 & 1
\end{array}\right]
\]
is two dimensional, abelian and normal.
\begin{itemize}
\item abelian: $(0,b,c)(0,b',c')=(0,b+b',c+c')$ 
\item normal: 
\begin{eqnarray*}
(a,b,c)(0,x,y)(a,b,c)^{-1} & = & (a,b,c)(0,x,y)(-a,-b,-c+ab)\\
 & = & (a,b+x,c+y+ax)(-a,-b,-c+ab)\\
 & = & (0,x,y+ax+ab-ab)\\
 & = & (0,x,ax+y)
\end{eqnarray*}
i.e., $Ad:G\rightarrow GL(\mathfrak{n})$, as 
\begin{eqnarray*}
Ad_{g}\left(n\right) & := & gng^{-1}\\
\left(x,y\right) & \mapsto & \left(x,ax+y\right)
\end{eqnarray*}
the orbit is a 2-d transformation.
\end{itemize}
Fix $h\in\mathbb{R}\backslash\{0\}$. Recall the \emph{Schrödinger
representation} of $G$ on $L^{2}(\mathbb{R})$ \index{Schrödinger representation}\index{representation!Schrödinger-}
\begin{equation}
U_{g}f(x)=e^{ih(c+bx)}f(x+a)\label{eq:sr2}
\end{equation}

\begin{thm}
The Schrödinger representation is induced. 
\end{thm}
\index{representation!unitary}
\begin{proof}
We show that the Schrödinger representation is induced from a unitary
representation $L$ on the subgroup $\Gamma$. \index{representation!Schrödinger-}

Note the Heisenberg group is a non abelian unimodular Lie group ($\triangle=1$,
$\delta=1$, and so $\rho\equiv1$.), The Haar measure on $G$ is
just the product measure $dxdydz$ on $\mathbb{R}^{3}$. Conditional
expectation becomes integrating out the variables correspond to the
subgroup. 

1. Let $L\in Rep(\Gamma,V)$ where $\Gamma=\{(0,b,c)\}$, $V=\mathbb{C}$,
\[
L_{\xi(b,c)}=e^{ihc}.
\]
The complex exponential comes in since we want a unitary representation.
The subgroup $\{(0,0,c)\}$ is the center of $G$. (What is the induced
representation? Is it unitarily equivalent to the Schrödinger representation?)

2. Look for the family $F_{*}$ of functions $f:G\rightarrow\mathbb{C}$
($V$ is the 1-d Hilbert space $\mathbb{C}$), such that
\[
f\left(\xi\left(b,c\right)g\right)=L_{\xi}f\left(g\right).
\]
Since 
\[
f\left(\xi\left(b,c\right)g\right)=f\left(\left(0,b,c\right)\left(x,y,z\right)\right)=f\left(x,b+y,c+z\right),\quad\mbox{and}
\]
\[
L_{\xi\left(b,c\right)}f\left(g\right)=e^{ihc}f\left(x,y,z\right)
\]
so $f$ satisfies
\[
f\left(x,b+y,c+z\right)=e^{ihc}f\left(x,y,z\right).
\]
That is, we may translate the $y,z$ variables by arbitrary amount,
and the only price to pay is the multiplicative factor $e^{ihc}$.
Therefore $f$ is really a function defined on the quotient 
\[
M=\Gamma\backslash G\simeq\mathbb{R}.
\]
The homogeneous space $M=\{(x,0,0)\}$ is identified with $\mathbb{R}$,
and the invariant measure on $M$ is simply the Lebesgue measure.
It is almost clear at this point why the induced representation is
unitarily equivalent to the Schrödinger representation on $L^{2}(\mathbb{R})$. 

3. The positive linear functional $\tau\varphi\mapsto\int_{G}\left\Vert f(g)\right\Vert _{V}^{2}\varphi\left(g\right)dg$
induces a measure $\mu_{f,f}$ on $M$. This can be seen as follows:
\begin{eqnarray*}
\int_{G}\left\Vert f(g)\right\Vert _{V}^{2}\varphi\left(g\right)dg & = & \int_{G\simeq\mathbb{R}^{3}}\left|f\left(x,y,z\right)\right|^{2}\varphi\left(x,y,z\right)dxdydz\\
 & = & \int_{M\simeq\mathbb{R}}\left(\int_{\Gamma\simeq\mathbb{R}^{2}}\left|f\left(x,y,z\right)\right|^{2}\varphi\left(x,y,z\right)dydz\right)dx\\
 & = & \int_{\mathbb{R}}\left|f\left(x,y,z\right)\right|^{2}\left(\int_{\mathbb{R}^{2}}\varphi\left(x,y,z\right)dydz\right)dx\\
 & = & \int_{\mathbb{R}}\left|f\left(x,y,z\right)\right|^{2}\left(\tau\varphi\right)\left(\pi\left(g\right)\right)dx\\
 & = & \int_{\mathbb{R}}\left|f\left(x,0,0\right)\right|^{2}\left(\tau\varphi\right)\left(x\right)dx.
\end{eqnarray*}
Note that
\begin{eqnarray*}
\left(\tau\varphi\right)\left(\pi\left(g\right)\right) & = & \int_{\Gamma}\varphi\left(\xi g\right)d\xi\\
 & = & \int_{\mathbb{R}^{2}}\varphi\left(\left(0,b,c\right)\left(x,y,z\right)\right)dbdc\\
 & = & \int_{\mathbb{R}^{2}}\varphi\left(x,b+y,c+z\right)dbdc\\
 & = & \int_{\mathbb{R}^{2}}\varphi\left(x,b,c\right)dbdc\\
 & = & \left(\tau\varphi\right)\left(x\right),\quad M=\Gamma\backslash G\simeq\mathbb{R}.
\end{eqnarray*}
Hence $\Lambda:C_{c}(M)\rightarrow\mathbb{C}$ given by
\[
\Lambda:\tau\varphi\longmapsto\int_{G}\left\Vert f(g)\right\Vert _{V}^{2}\varphi(g)dg
\]
is a positive linear functional, therefore $\Lambda=\mu_{f,f}$ and
\[
\int_{\mathbb{R}^{3}}\left|f\left(x,y,z\right)\right|^{2}\varphi\left(x,y,z\right)dxdydz=\int_{\mathbb{R}}\left(\tau\varphi\right)\left(x\right)d\mu_{f,f}\left(x\right).
\]

4. Define
\[
\begin{alignedat}{1}\left\Vert f\right\Vert _{ind}^{2} & :=\mu_{f,f}\left(M\right)=\int_{M}\left|f\right|^{2}d\xi=\int_{\mathbb{R}}\left|f\left(x,y,z\right)\right|^{2}dx=\int_{\mathbb{R}}\left|f\left(x,0,0\right)\right|^{2}dx\end{alignedat}
\]
\[
U_{g}^{ind}f\left(g'\right):=f\left(g'g\right)
\]
By definition, if $g=g(a,b,c)$ and $g'=g'(x,y,z)$ then
\begin{eqnarray*}
U_{g}^{ind}f(g') & = & f(g'g)\\
 & = & f((x,y,z)(a,b,c))\\
 & = & f(x+a,y+b,z+c+xb)
\end{eqnarray*}
and $U^{ind}$ is a unitary representation.

5. To see that $U^{ind}$ is unitarily equivalent to the Schrödinger
representation on $L^{2}(\mathbb{R})$, we set\index{representation!Schrödinger-}
\[
W:\mathscr{H}^{ind}\rightarrow L^{2}\left(\mathbb{R}\right),\quad\left(Wf\right)\left(x\right)=f\left(x,0,0\right)
\]
(If put other numbers into $f$, as $f(x,y,z)$, the result is the
same, since $f\in\mathscr{H}^{ind}$ is really defined on the quotient
$M=\Gamma\backslash G\simeq\mathbb{R}$. ) 

$W$ is unitary: 
\[
\left\Vert Wf\right\Vert _{L^{2}}^{2}=\int_{\mathbb{R}}\left|Wf\right|^{2}dx=\int_{\mathbb{R}}\left|f(x,0,0)\right|^{2}dx=\int_{\Gamma\backslash G}\left|f\right|^{2}d\xi=\left\Vert f\right\Vert _{ind}^{2}
\]

The intertwining property: Let $U_{g}$ be the Schrödinger representation,
then 
\begin{eqnarray*}
U_{g}\left(Wf\right) & = & e^{ih\left(c+bx\right)}f\left(x+a,0,0\right)\\
WU_{g}^{ind}f & = & W\left(f\left(\left(x,y,z\right)\left(a,b,c\right)\right)\right)\\
 & = & W\left(f\left(x+a,y+b,z+c+xb\right)\right)\\
 & = & W\left(e^{ih\left(c+bx\right)}f\left(x+a,y,z\right)\right)\\
 & = & e^{ih\left(c+bx\right)}f\left(x+a,0,0\right).
\end{eqnarray*}

6. Since $\left\{ U,L\right\} '\subset\left\{ L\right\} '$, the system
$\left\{ U,L\right\} $ is reducible implies $L$ is reducible. Equivalent,
$\left\{ L\right\} $ is irreducible implies $\left\{ U,L\right\} $
is irreducible. Since $L$ is 1-dimensional, it is irreducible. Consequently,
$U^{ind}$ is irreducible. \end{proof}
\begin{xca}[The Schrödinger representation]
\myexercise{The Schrödinger representation} Prove that for $h\neq0$
fixed, the Schrödinger representation $U^{h}$ (\ref{eq:sr2}) is
irreducible. 

\uline{Hint}: Show that if $A\in\mathscr{B}\left(L^{2}\left(\mathbb{R}\right)\right)$
commutes with $\left\{ U_{g}^{h}\::\:g\in G_{\text{Heis}}\right\} $,
then there exists $\lambda\in\mathbb{C}$ such that $A=\lambda I_{L^{2}\left(\mathbb{R}\right)}$,
i.e., that the commutant of the representation $U^{h}$ is one-dimensional. 
\end{xca}
\index{Schrödinger representation}\index{representation!Schrödinger-}

\subsection{\label{sec:em}$ax+b$ group}

\index{groups!$ax+b$}

$a\in\mathbb{R}_{+}$, $b\in\mathbb{R}$, $g=(a,b)=\left[\begin{array}{cc}
a & b\\
0 & 1
\end{array}\right]$. 
\[
U_{g}f(x)=e^{iax}f(x+b)
\]
could also write $a=e^{t}$, then 
\[
U_{g}f(x)=e^{ie^{t}x}f(x+b)
\]
\[
U_{g(a,b)}f(x)=e^{iae^{x}f(x+b)}
\]
\begin{eqnarray*}
[\frac{d}{dx},ie^{x}] & = & ie^{x}\\
{}[A,B] & = & B
\end{eqnarray*}
or
\[
U_{(e^{t},b)}f=e^{ite^{x}}f(x+b)
\]
\[
\left[\begin{array}{cc}
0 & b\\
0 & 1
\end{array}\right]
\]
1-d representation. $L_{b}=e^{ib}$. Induce $ind_{L}^{G}\simeq$ the
Schrödinger representation.

\section{Co-adjoint Orbits}

It turns out that only a small family of representations are induced.
The question is how to detect whether a representation is induced.
The whole theory is also under the name of ``Mackey machine'' \cite{Mac52,Mac88}.
The notion of ``machine'' refers to something that one can actually
compute in practice. Two main examples are the Heisenberg group and
the $ax+b$ group. \index{co-adjoint orbit} \index{representation!co-adjoint}

What is the mysteries parameter $h$ that comes into the Schrödinger
representation? It is a physical constant, but how to explain it in
mathematical theory?\index{Schrödinger, E. R.}\index{Mackey machine}

\subsection{Review of some Lie theory}
\begin{thm}[Ado]
Every Lie group is diffeomorphic to a matrix group.
\end{thm}
The exponential function $\exp$ maps a neighborhood of $0$ into
a connected component of $G$ containing the identity element. For
example, the Lie algebra of the Heisenberg group is \index{Lie!algebra}
\index{Lie!group}
\[
\left[\begin{array}{ccc}
0 & * & *\\
0 & 0 & *\\
0 & 0 & 0
\end{array}\right]
\]
All the Lie\index{groups!Lie} groups the we will ever encounter come
from a quadratic form. Given a quadratic form \index{quadratic form}
\[
\varphi:V\times V\rightarrow\mathbb{C}
\]
there is an associated group that fixes $\varphi$, i.e. we consider
elements $g$ such that 
\[
\varphi(gx,gy)=\varphi(x,y)
\]
and define $G(\varphi)$ as the collection of these elements. $G(\varphi)$
is clearly a group. Apply the exponential map and the product rule,
\[
\frac{d}{dt}\big|_{t=0}\varphi(e^{tX}x,e^{tX}y)=0\Longleftrightarrow\varphi(Xx,y)+\varphi(x,Xy)=0
\]
hence 
\[
X+X^{tr}=0
\]
The determinant and trace are related so that \index{trace}
\[
\det(e^{tX})=e^{t\cdot trace(X)}
\]
thus $\det=1$ if and only if $trace=0$. It is often stated in differential
geometry that the derivative of the determinant is equal to the trace.
\begin{example}
$\mathbb{R}^{n}$, $\varphi(x,y)=\sum x_{i}y_{i}$. The associated
group is the orthogonal group $O_{n}$.
\end{example}
There is a famous cute little trick to make $O_{n-1}$ into a subgroup
of $O_{n}$. $O_{n-1}$ is not normal in $O_{n}$. We may split the
quadratic form into
\[
\sum_{i=1}^{n-1}x_{i}^{2}+1
\]
where $1$ corresponds to the last coordinate in $O_{n}$. Then we
may identity $O_{n-1}$ as a subgroup of $O_{n}$ 
\[
g\mapsto\left[\begin{array}{cc}
g & 0\\
0 & I
\end{array}\right]
\]
where $I$ is the identity operator. 

Claim: $O_{n}/O_{n-1}\simeq S^{n-1}$. How to see this? Let $u$ be
the unit vector corresponding to the last dimension, look for $g$
that fixes $u$ i.e. $gu=u$. Such $g$ forms a subgroup of $O_{n}$,
and it is called isotropy group. 
\[
I_{n}=\{g:gu=u\}\simeq O_{n-1}
\]
Notice that for all $v\in S^{n-1}$, there exists $g\in O_{n}$ such
that $gu=v$. Hence 
\[
g\mapsto gu
\]
in onto $S^{n-1}$. The kernel of this map is $I_{n}\simeq O_{n-1}$,
thus 
\[
O_{n}/O_{n-1}\simeq S_{n}
\]
Such spaces are called homogeneous spaces.
\begin{example}
visualize this with $O_{3}$ and $O_{2}$.
\end{example}
Other examples of homogeneous spaces show up in number theory all
the time. For example, the Poincaré\index{groups!Poincaré} group
$G/\mbox{discrete subgroup}$.

$G$, $N\subset G$ normal subgroup. The map $g\cdot g^{-1}:G\rightarrow G$
is an automorphism sending identity to identity, hence if we differentiate
it, we get a transformation in $GL(\mathfrak{g})$. i.e. we get a
family of maps $Ad_{g}\in GL(\mathfrak{g})$ indexed by elements in
$G$. $g\mapsto Ad_{g}\in GL(\mathfrak{g})$ is a representation of
$G$, hence if it is differentiated, we get a representation of $\mathfrak{g}$,
$ad_{g}:\mathfrak{g}\mapsto End(\mathfrak{g})$ acting on the vector
space $\mathfrak{g}$.

$gng^{-1}\in N$. $\forall g$, $g\cdot g^{-1}$ is a transformation
from $N$ to $N$, define $Ad_{g}(n)=gng^{-1}$. Differentiate to
get $ad:\mathfrak{n}\rightarrow\mathfrak{n}$. $\mathfrak{n}$ is
a vector space, has a dual. Linear transformation on vector space
passes to the dual space. 
\begin{eqnarray*}
\varphi^{*}(v^{*})(u) & = & v^{*}(\varphi(u))\\
 & \Updownarrow\\
\left\langle \Lambda^{*}v^{*},u\right\rangle  & = & \left\langle v^{*},\Lambda u\right\rangle .
\end{eqnarray*}
In order to get the transformation rules work out, have to pass to
the adjoint or the dual\index{dual} space.
\[
Ad_{g}^{*}:\mathfrak{n}^{*}\rightarrow\mathfrak{n}^{*}
\]
the coadjoint representation of $\mathfrak{n}$. \index{representation!adjoint-}

Orbits of co-adjoint representation amounts precisely to equivalence
classes of irreducible representations.\index{representation!co-adjoint}
\begin{example}
Heisenberg group $G=\{(a,b,c)\}$ with 
\[
(a,b,c)=\left[\begin{array}{ccc}
1 & a & c\\
0 & 1 & b\\
0 & 0 & 1
\end{array}\right]
\]
normal subgroup $N=\{(0,b,c)\}$
\[
(0,b,c)=\left[\begin{array}{ccc}
1 & 0 & c\\
0 & 1 & b\\
0 & 0 & 1
\end{array}\right]
\]
with Lie algebra $\mathfrak{n}=\{(b,c)\}$
\[
(0,\xi,\eta)=\left[\begin{array}{ccc}
1 & 0 & c\\
0 & 1 & b\ \\
0 & 0 & 1
\end{array}\right]
\]
$Ad_{g}:\mathfrak{n}\rightarrow\mathfrak{n}$ given by
\begin{eqnarray*}
gng^{-1} & = & (a,b,c)(0,y,x)(-a,-b,-c+ab)\\
 & = & (a,b+y,c+x+ay)(-a,-b,-c+ab)\\
 & = & (0,y,x+ay)
\end{eqnarray*}
hence $Ad_{g}:\mathbb{R}^{2}\rightarrow\mathbb{R}^{2}$
\[
Ad_{g}:\left[\begin{array}{c}
x\\
y
\end{array}\right]\mapsto\left[\begin{array}{c}
x+ay\\
y
\end{array}\right].
\]
The matrix of $Ad_{g}$ is (before taking adjoint) is
\[
Ad_{g}=\left[\begin{array}{cc}
1 & a\\
0 & 1
\end{array}\right].
\]
The matrix for $Ad_{g}^{*}$ is
\[
Ad_{g}^{*}=\left[\begin{array}{cc}
1 & 0\\
a & 1
\end{array}\right].
\]
We use $[\xi,\eta]^{T}$ for the dual $\mathfrak{n}^{*}$; and use
$[x,y]^{T}$ for $\mathfrak{n}$. Then
\[
Ad_{g}^{*}:\left[\begin{array}{c}
\xi\\
\eta
\end{array}\right]\mapsto\left[\begin{array}{c}
\xi\\
a\xi+\eta
\end{array}\right]
\]
What about the orbit? In the example of $O_{n}/O_{n-1}$, the orbit
is $S^{n-1}$. 

For $\xi\in\mathbb{R}\backslash\{0\}$, the orbit of $Ad_{g}^{*}$
is 
\[
\left[\begin{array}{c}
\xi\\
0
\end{array}\right]\mapsto\left[\begin{array}{c}
\xi\\
\mathbb{R}
\end{array}\right]
\]
i.e. vertical lines with $x$-coordinate $\xi$. $\xi=0$ amounts
to fixed point, i.e. the orbit is a fixed point.

The simplest orbit is when the orbit is a fixed point. i.e.
\[
Ad_{g}^{*}:\left[\begin{array}{c}
\xi\\
\eta
\end{array}\right]\mapsto\left[\begin{array}{c}
\xi\\
\eta
\end{array}\right]\in V^{*}
\]
where if we choose 
\[
\left[\begin{array}{c}
\xi\\
\eta
\end{array}\right]=\left[\begin{array}{c}
0\\
1
\end{array}\right]
\]
it is a fixed point. 

The other extreme is to take any $\xi\neq0$, then 
\[
Ad_{g}^{*}:\left[\begin{array}{c}
\xi\\
0
\end{array}\right]\mapsto\left[\begin{array}{c}
\xi\\
\mathbb{R}
\end{array}\right]
\]
i.e. get vertical lines indexed by the $x$-coordinate $\xi$. In
this example, a cross section is a subset of $\mathbb{R}^{2}$ that
intersects each orbit at precisely one point. Every cross section
in this example is a Borel set in $\mathbb{R}^{2}$.

We don't always get measurable cross sections. An example is the construction
of non-measurable set as was given in Rudin's book. Cross section
is a Borel set that intersects each coset at precisely one point.

Why does it give all the equivalent classes of irreducible representations?
Since we have a unitary representation $L_{n}\in Rep(N,V)$, $L_{n}:V\rightarrow V$
and by construction of the induced representation $U_{g}\in Rep(G,\mathscr{H})$,
$N\subset G$ normal such that
\[
U_{g}L_{n}U_{g^{-1}}=L_{gng^{-1}}
\]
i.e. 
\[
L_{g}\simeq L_{gng^{-1}}
\]
now pass to the Lie algebra and its dual
\[
L_{n}\rightarrow LA\rightarrow LA^{*}.
\]

\end{example}
\index{algebras!Lie algebra} \index{Lie!algebra}

\index{space!dual-}

\index{representation!unitary}

\section{\label{sec:gaarding}Gårding Space}
\begin{defn}
Let $\mathcal{U}$ be a strongly continuous representation of a Lie
group $G$, with Lie algebra $\mathfrak{g}$, and let $\exp:\mathfrak{g}\rightarrow G$
denote the exponential mapping from Lie theory. Fro every $\varphi\in C_{c}^{\infty}\left(G\right)$,
set \index{Gårding!-space} \index{strongly continuous}\index{representation!strongly continuous}
\[
\mathcal{U}\left(\varphi\right)=\int_{G}\varphi\left(g\right)\mathcal{U}_{g}dg
\]
where $dg$ is a left-invariant Haar measure on $G$; and set 
\[
\mathscr{H}_{G\mathring{a}rding}=\left\{ \mathcal{U}\left(\varphi\right)v\;\Big|\;\varphi\in C_{c}^{\infty}\left(G\right),v\in\mathscr{H}\right\} .
\]
\end{defn}
\begin{lem}
Fix $X\in\mathfrak{g}$, set 
\[
d\mathcal{U}\left(X\right)v=\lim_{t\rightarrow0}\frac{\mathcal{U}\left(\exp\left(tX\right)\right)v-v}{t}
\]
then 
\[
\mathscr{H}_{G\mathring{a}rding}\subset\bigcap_{X\in\mathfrak{g}}dom\left(d\mathcal{U}\left(X\right)\right)
\]
and 
\[
d\mathcal{U}\left(X\right)\mathcal{U}\left(\varphi\right)v=\mathcal{U}\left(\widetilde{X}\varphi\right)v,
\]
for all $\varphi\in C_{c}^{\infty}\left(G\right)$, $v\in\mathscr{H}$,
where 
\[
\left(\widetilde{X}\varphi\right)\left(g\right)=\frac{d}{dt}\big|_{t=0}\varphi\left(\exp\left(-tX\right)g\right),\;\forall g\in G.
\]
 \end{lem}
\begin{proof}
(Hint) 
\[
\int_{G}\varphi\left(g\right)\mathcal{U}\left(exp\left(tX\right)\right)\mathcal{U}\left(g\right)dg=\int_{G}\varphi\left(\exp\left(-tX\right)g\right)\mathcal{U}\left(g\right)dg.
\]

\end{proof}
We talked about how to detect whether a representation is induced.
Given a group $G$ with a subgroup $\Gamma$ let $M:=\Gamma\backslash G$.
The map $\pi:G\rightarrow M$ is called a covering map, which sends
$g$ to its equivalent class or the coset $\Gamma g$. $M$ is given
its projective topology, so $\pi$ is continuous. When $G$ is compact,
many things simplify. For example, if $G$ is compact, any irreducible
representation is finite dimensional. But many groups are not compact,
only locally compact. For example, the groups $ax+b$, $H_{3}$, $SL_{n}$.

Specialize to Lie groups. $G$ and subgroup $H$ have Lie algebras
$\mathfrak{g}$ and $\mathfrak{h}$ respectively.\index{algebras!Lie algebra}\index{groups!Lie}\index{groups!$ax+b$}
\index{Lie!group} \index{Lie!algebra} 
\[
\mathfrak{g}=\{X:e^{tX}\in G,\forall t\in\mathbb{R}\}
\]
Almost all Lie algebras we will encounter come from specifying a quadratic
form $\varphi:G\times G\rightarrow\mathbb{C}$. $\varphi$ is then
uniquely determined by a Hermitian matrix $A$ so that \index{quadratic form}
\[
\varphi(x,y)=x^{tr}\cdot Ay
\]
Let $G=G(\varphi)=\{g:\varphi(gx,gy)=\varphi(x,y)\}$, then 
\[
\frac{d}{dt}\big|_{t=0}\varphi(e^{tX}x,e^{tX}y)=0
\]
and with an application of the product rule,
\begin{eqnarray*}
\varphi(Xx,y)+\varphi(x,Xy) & = & 0\\
(Xx)^{tr}\cdot Ay+x^{tr}\cdot AXy & = & 0
\end{eqnarray*}
 
\[
X^{tr}A+AX=0
\]
hence 
\[
\mathfrak{g}=\{X:X^{tr}A+AX=0\}.
\]

Let $U\in Rep(G,\mathscr{H})$, for $X\in\mathfrak{g}$, $U(e^{tX})$
is a one parameter continuous group of unitary operator, hence by
Stone's theorem\index{Stone's theorem} (see \cite{vN32c,Ne69}),
it must have the form\index{Theorem!Stone's-} 
\begin{equation}
U(e^{tX})=e^{itH_{X}}\label{eq:su}
\end{equation}
for some selfadjoint operator $H_{X}$ (possibly unbounded). The RHS
in (\ref{eq:su}) is given by the Spectral Theorem (see \cite{Sto90, Yos95, Ne69, RS75, DS88b}).
We often write 
\[
dU(X):=iH_{X}
\]
to indicate that $dU(X)$ is the directional derivative along the
direction $X$. Notice that $H_{X}^{*}=H_{X}$ but 
\[
(iH_{X})^{*}=-(iH_{X})
\]
i.e. $dU(X)$ is skew adjoint.

\index{Schrödinger representation}\index{representation!Schrödinger-}
\begin{example}
$G=\{(a,b,c)\}$ Heisenberg group. $\mathfrak{g}=\{X_{1}\sim a,X_{2}\sim b,X_{3}\sim c\}$.
Take the Schrödinger representation $U_{g}f(x)=e^{ih(c+bx)}f(x+a)$,
$f\in L^{2}(\mathbb{R})$. 
\begin{itemize}
\item $U(e^{tX_{1}})f(x)=f(x+t)$
\begin{eqnarray*}
\frac{d}{dt}\big|_{t=0}U(e^{tX_{1}})f(x) & = & \frac{d}{dx}f(x)\\
dU(X_{1}) & = & \frac{d}{dx}
\end{eqnarray*}

\item $U(e^{tX_{2}})f(x)=e^{ih(tx)}f(x)$
\begin{eqnarray*}
\frac{d}{dt}\big|_{t=0}U(e^{tX_{2}})f(x) & = & ihxf(x)\\
dU(X_{2}) & = & ihx
\end{eqnarray*}

\item $U(e^{tX_{3}})f(x)=e^{iht}f(x)$
\begin{eqnarray*}
\frac{d}{dt}\big|_{t=0}U(e^{tX_{3}})f(x) & = & ihf(x)\\
dU(X_{2}) & = & ihI
\end{eqnarray*}
Notice that $dU(X_{i})$ are all skew adjoint.
\begin{eqnarray*}
[dU(X_{1}),dU(X_{2})] & = & [\frac{d}{dx},ihx]\\
 & = & ih[\frac{d}{dx},x]\\
 & = & ih
\end{eqnarray*}
In case we want selfadjoint operators, replace $dU(X_{i})$ by$-idU(X_{i})$
and get
\begin{eqnarray*}
-idU(X_{1}) & = & \frac{1}{i}\frac{d}{dx}\\
-idU(X_{2}) & = & hx\\
-idU(X_{3}) & = & hI
\end{eqnarray*}
\[
[\frac{1}{i}\frac{d}{dx},hx]=\frac{h}{i}.
\]

\end{itemize}
\end{example}
Below we answer the following question:

What is the space of functions that $U_{g}$ acts on? L. Gårding /gor-ding/
(Swedish mathematician) looked for one space that always works. It's
now called the Gårding space.

\begin{figure}
\includegraphics[width=0.5\textwidth]{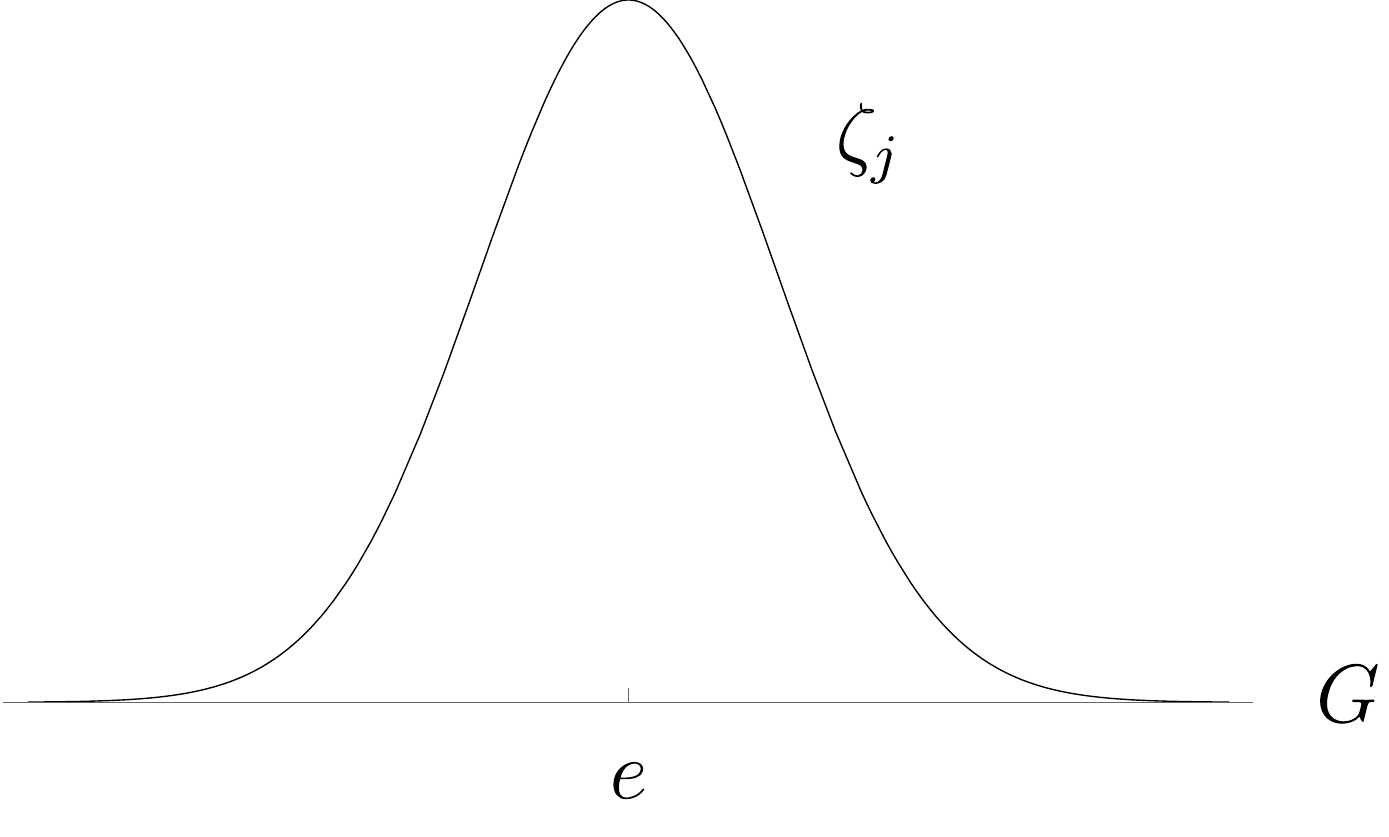}

\protect\caption{\label{fig:approx1}Approximation of identity.}
\end{figure}

Start with $C_{c}(G)$, every $\varphi\in C_{c}(G)$ can be approximated
by the so called Gårding functions, using the convolution argument.
Define convolution as \index{convolution}\index{Gårding!-space}
\begin{eqnarray*}
\varphi\star\psi(g) & = & \int_{G}\varphi(gh)\psi(h)d_{R}h\\
\varphi\star\psi(g) & = & \int_{G}\varphi(h)\psi(g^{-1}h)d_{L}h
\end{eqnarray*}
Take an approximation of identity $\zeta_{j}$ (\figref{approx1}),
so that
\[
\varphi\star\zeta_{j}\rightarrow\varphi,\;j\rightarrow0.
\]

Define Gårding space as the span of the vectors in $\mathscr{H}$,
given by
\[
U(\varphi)v=\int\varphi(h)U(h)vd_{L}h
\]
where $\varphi\in C_{c}(G)$, $v\in\mathscr{H}$, or we say
\[
U(\varphi):=\int_{G}\varphi(h)U(h)\,d_{L}h.
\]
Since $\varphi$ vanishes outside a compact set, and since $U(h)v$
is continuous and bounded in $\left\Vert \cdot\right\Vert $,  it
follows that $U(\varphi)$ is well-defined.

Every representation $\mathcal{U}$ of a Lie group $G$ induces a
representation (also denote $\mathcal{U}$) of the group algebra:
\begin{lem}
$U(\varphi_{1}\star\varphi_{2})=U(\varphi_{1})U(\varphi_{2})$ ($U$
is a representation of the group algebra)\end{lem}
\begin{proof}
Use Fubini,
\begin{eqnarray*}
\int_{G}\varphi_{1}\star\varphi_{2}(g)U(g)dg & = & \iint_{G\times G}\varphi_{1}(h)\varphi(h^{-1}g)U(g)dhdg\\
 & = & \iint_{G\times G}\varphi_{1}(h)\varphi(g)U(hg)dhdg\:(dg\mbox{ is r-Haar}g\mapsto hg)\\
 & = & \iint_{G\times G}\varphi_{1}(h)\varphi(g)U(h)U(g)dhdg\\
 & = & \int_{G}\varphi_{1}(h)U(h)dh\int_{G}\varphi_{2}(g)U(g)dg
\end{eqnarray*}
Choose $\varphi$ to be an approximation of identity, then 
\[
\int_{G}\varphi(g)U(g)vdg\rightarrow U(e)v=v
\]
i.e. any vector $v\in H$ can be approximated by functions in the
Gårding space. It follows that
\[
\{U(\varphi)v\}
\]
is dense in $\mathscr{H}$.\end{proof}
\begin{lem}
\label{lem:dU}$U(\varphi)$ can be differentiated, in the sense that
\[
dU(X)U(\varphi)v=U(\tilde{X}\varphi)v
\]
where we use $\tilde{X}$ to denote the vector field.\end{lem}
\begin{proof}
need to prove 
\[
\lim_{t\rightarrow0}\frac{1}{t}\left[(U(e^{tX})-I)U(\varphi)v\right]=U(\tilde{X}\varphi)v.
\]
Let $v_{\varphi}:=U(\varphi)v$, need to look at in general $U(g)v_{\varphi}$.
\begin{eqnarray*}
U(g)v_{\varphi} & = & U(g)\int_{G}\varphi(h)U(h)vdh\\
 & = & \int_{G}\varphi(h)U(gh)dh\\
 & = & \int_{G}\triangle(g)\varphi(g^{-1}h)U(h)dh
\end{eqnarray*}
set $g=e^{tX}$.\end{proof}
\begin{note}
If assuming unimodular, $\triangle$ does not show up. Otherwise,
$\triangle$ is some correction term which is also differentiable.
$\tilde{X}$ acts on $\varphi$ as $\tilde{X}\varphi$. $\tilde{X}$
is called the derivative of the translation operator $e^{tX}$.\end{note}
\begin{xca}[The Gårding space for the Schrödinger representation]
\myexercise{The Gårding space for the Schrödinger representation}\label{exer:gp4}Show
that Schwartz space $\mathcal{S}$ is the Gårding space for the Schrödinger
representation. \index{Schrödinger representation}\index{Gårding!-space}\index{representation!Schrödinger-}
\end{xca}

\begin{xca}[The Lie bracket]
\myexercise{The Lie bracket}Let $U$ be a representation of a Lie
group $G$, and let $dU\left(\cdot\right)$ be the derived representation,
see  \lemref{dU}. On the dense Gårding space, show that 
\begin{eqnarray*}
dU\left(\left[X,Y\right]\right) & = & \left[dU\left(X\right),dU\left(Y\right)\right]\\
 & = & dU\left(X\right)dU\left(Y\right)-dU\left(Y\right)dU\left(X\right),
\end{eqnarray*}
where $\left[X,Y\right]$ denotes the Lie bracket of the two elements
$X$ and $Y$ in the Lie algebra.
\end{xca}

\section{\label{sec:dreprep}Decomposition of Representations}

We study some examples of duality. 
\begin{itemize}
\item $G=\mathbb{T}$, $\hat{G}=\mathbb{Z}$
\begin{eqnarray*}
\chi_{n}(z) & = & z^{n}\\
\chi_{n}(zw) & = & z^{n}w^{n}=\chi_{n}(z)\chi_{n}(w)
\end{eqnarray*}

\item $G=\mathbb{R}$, $\hat{G}=\mathbb{R}$ 
\begin{eqnarray*}
\chi_{t}(x) & = & e^{itx}
\end{eqnarray*}

\item $G=\mathbb{Z}/n\mathbb{Z}\simeq\{0,1,\cdots,n-1\}$. $\hat{G}=G$.
\\
This is another example where $\hat{G}=G$. \\
Let $\zeta=e^{i2\pi/n}$ be the primitive $n^{th}$-root of unity.
$k\in\mathbb{Z}_{n}$, $l\in\{0,1,\ldots,n-1\}$
\[
\chi_{l}(k)=e^{i\frac{2\pi kl}{n}}
\]

\end{itemize}
If $G$ is a locally compact abelian group, $\hat{G}$ is the set
of 1-dimensional representations. 
\[
\hat{G}=\{\chi:g\mapsto\chi(g)\in\mathbb{T},\chi(gh)=\chi(g)\chi(h),\;\mbox{assumed continuous}\}.
\]
$\hat{G}$ is also a group, with group operation defined by $(\chi_{1}\chi_{2})(g):=\chi_{1}(g)\chi_{2}(g)$.
$\hat{G}$ is called the group characters.
\begin{thm}[Pontryagin]
 If $G$ is a locally compact abelian group, then $G\simeq\hat{\hat{G}}$
(isomorphism between $G$ and the double dual $\hat{\hat{G}}$,) where
``$\simeq$'' means ``natural isomorphism.''\end{thm}
\begin{note}
This result first appears in 1930s in the annals of math, when John
von Neumann was the editor of the journal at the time. The original
paper was hand written. von Neumann rewrote it, since then the theorem
became very popular, see \cite{Ru90}.\index{groups!abelian} \index{Theorem!Pontryagin duality-}
\end{note}
There are many groups that are not abelian. We want to study the duality
question in general. Examples: \index{dual}
\begin{itemize}
\item compact group
\item finite group (abelian, or not)
\item $H_{3}$ locally compact, nonabelian, unimodular\index{unimodular}\index{groups!$ax+b$}
\item $ax+b$ locally compact, nonabelian, non-unimodular
\end{itemize}
If $G$ is not abelian, $\hat{G}$ is not a group. We would like to
decompose $\hat{G}$ into irreducible representations. The big names
in this development are Krein, Peter-Weyl, Weil, Segal. See \cite{AD86,AAR13,BR79,EMC00,JO00,KL14,MR1468230,Rud73,Ru90,Seg50,Sto90}.

Let $G$ be a group (may not be abelian). The right regular representation
is defined as
\[
R_{g}f(\cdot)=f(\cdot g),\;\left(\mbox{translation on the right}\right).
\]
Then $R_{g}$ is a unitary operator acting on $L^{2}(\mu_{R})$, where
$\mu_{R}$ is the right invariant Haar measure.
\begin{thm}[Krein, Weil, Segal]
 Let $G$ be locally compact unimodular (abelian or not). Then the
right regular representation decomposes into a direct integral of
irreducible representations 
\[
R_{g}=\int_{\hat{G}}^{\oplus}"irrep"\:d\mu
\]

\end{thm}
where $\mu$ is called the Plancherel\index{measure!Plancherel} measure.
See \cite{Sti59,Seg50}.
\begin{example}
$G=T$, $\hat{G}=\mathbb{Z}$. Irreducible representations $\{e^{in(\cdot)}\}_{n}\sim\mathbb{Z}$
\begin{eqnarray*}
(U_{y}f)(x) & = & f(x+y)\\
 & = & \sum_{n}\hat{f}(n)\chi_{n}(x+y)\\
 & = & \sum_{n}\hat{f}(n)e^{i2\pi n(x+y)}
\end{eqnarray*}
\[
(U_{y}f)(0)=f(y)=\sum_{n}\hat{f}(n)e^{i2\pi ny}
\]
The Plancherel measure in this case is the counting measure.
\begin{example}
$G=\mathbb{R}$, $\hat{G}=\mathbb{R}$. Irreducible representations
$\{e^{it(\cdot)}\}_{t\in\mathbb{R}}\sim\mathbb{R}$.
\begin{eqnarray*}
(U_{y}f)(x) & = & f(x+y)\\
 & = & \int_{\mathbb{R}}\hat{f}(t)\chi_{t}(x+y)dt\\
 & = & \int_{\mathbb{R}}\hat{f}(t)e^{it(x+y)}dt
\end{eqnarray*}
\[
(U_{y}f)(0)=f(y)=\int_{\mathbb{R}}\hat{f}(t)e^{ity}dt
\]
where the Plancherel measure is the Lebesgue measure on $\mathbb{R}$.
\end{example}
\end{example}
As can be seen that Fourier series and Fourier integrals are special
cases of the decomposition of the right regular representation $R_{g}$
of a unimodular locally compact\index{groups!locally compact} group.
$\int^{\oplus}$ $\Longrightarrow$ $\left\Vert f\right\Vert =\left\Vert \widehat{f}\right\Vert $.
This is a result that was done 30 years earlier before the non abelian
case. Classical function theory studies other types of convergence,
pointwise, uniform, etc. 
\begin{example}
$G=H_{3}$. $G$ is unimodular\index{unimodular}, non abelian. $\hat{G}$
is not a group. 

Irreducible representations: $\mathbb{R}\backslash\{0\}$ Schrödinger
representation, $\{0\}$ 1-d trivial representation

Decomposition: 
\[
R_{g}=\int_{\mathbb{R}\backslash\{0\}}^{\oplus}U_{irrep}^{h}hdh
\]
For all $f\in L^{2}(G)$,
\[
(U_{g}f)(e)=\int^{\oplus}U^{h}f\:hdh,\quad U^{h}\mbox{ irrep}.
\]
Set 
\begin{eqnarray*}
F(g) & = & (R_{g}F)(e)\\
 & = & \int_{\mathbb{R}\backslash\{0\}}^{\oplus}e^{ih(c+bx)}f(x+a)\:hdh;\;\mbox{then}\\
\hat{F}(h) & = & \int_{G}(U_{g}^{h}F)dg
\end{eqnarray*}
\emph{Plancherel measure}: $hdh$ and the point measure $\delta_{0}$
at zero.
\end{example}
\index{Schrödinger representation}\index{representation!Schrödinger-}

\begin{example}
$G$ = $ax+b$ group, non abelian. $\hat{G}$ not a group. 3 irreducible
representations: $+,-,0$ but $G$ is not unimodular. \index{groups!$ax+b$}

The $+$ representation is supported on $\mathbb{R}_{+}$, the $-$
representation on $\mathbb{R}_{-}$, and the $0$ representation is
the trivial one-dimensional representation.
\end{example}
The duality question may also be asked for discrete subgroups. This
leads to remarkable applications in automorphic functions, automorphic
forms, p-adic numbers, compact Riemann surface, hyperbolic geometry,
etc. 
\begin{example}
Cyclic group of order $n$. $G=\mathbb{Z}/n\mathbb{Z}\simeq\{0,1,\cdots,n-1\}$.
$\hat{G}=G$. This is another example where the dual group is identical
to the group itself. Let $\zeta=e^{i2\pi/n}$ be the primitive $n^{th}$-root
of unity. $k\in\mathbb{Z}_{n}$, $l=\{0,1,\ldots,n-1\}$
\[
\chi_{l}(k)=e^{i\frac{2\pi kl}{n}}
\]
In this case, Segal's theorem gives finite Fourier transform. $U:l^{2}(\mathbb{Z})\rightarrow l^{2}(\hat{\mathbb{Z}})$
where
\[
Uf(l)=\frac{1}{\sqrt{N}}\sum_{k}\zeta^{kl}f(k)
\]

\end{example}

\section{Summary of Induced Representations, the Example of $d/dx$}

We study decomposition of group representations. Two cases: abelian
and non abelian. The non abelian case may be induced from the abelian
ones.

non abelian
\begin{itemize}
\item semi product $G=HN$ often $N$ is normal.
\item $G$ simple. $G$ does not have normal subgroups, i.e., the Lie algebra
does not have any ideals. \end{itemize}
\begin{xca}[Normal subgroups]
\myexercise{Normal subgroups}\label{exer:gp5}(1) Find the normal
subgroups in the Heisenberg group. (2) Find the normal subgroups in
the $ax+b$ group.\end{xca}
\begin{example}
$SL_{2}(\mathbb{R})$ (non compact)
\[
\left(\begin{array}{cc}
a & b\\
c & d
\end{array}\right),\;ad-bc=1
\]
with Lie algebra 
\[
sl_{2}(\mathbb{R})=\{X:tr(X)=0\}.
\]
Note that $sl_{2}$ is generated by 
\[
\left(\begin{array}{cc}
0 & 1\\
1 & 0
\end{array}\right),\;\left(\begin{array}{cc}
0 & -1\\
1 & 0
\end{array}\right),\;\left(\begin{array}{cc}
1 & 0\\
0 & -1
\end{array}\right).
\]
In particular, $\left(\begin{array}{cc}
0 & -1\\
1 & 0
\end{array}\right)$ generates the one-parameter group $\left(\begin{array}{cc}
\cos t & -\sin t\\
\sin t & \cos t
\end{array}\right)\simeq\mathbb{T}$ whose dual group is $\mathbb{Z}$, where 
\[
\chi_{n}(g(t))=g(t)^{n}=e^{itn}.
\]
May use this to induce a representation of $G$. This is called principle
series. Need to do something else to get all irreducible representations.

A theorem by Iwasawa states that simple matrix group (Lie group) can
be decomposed into \index{Lie!group} 
\[
G=KAN
\]
where $K$ is compact, $A$ is abelian and $N$ is nilpotent. For
example, in the $SL_{2}$ case, 
\[
SL_{2}(\mathbb{R})=\left(\begin{array}{cc}
\cos t & -\sin t\\
\sin t & \cos t
\end{array}\right)\left(\begin{array}{cc}
e^{s} & 0\\
0 & e^{-s}
\end{array}\right)\left(\begin{array}{cc}
1 & u\\
0 & 1
\end{array}\right).
\]
The simple groups do not have normal subgroups. The representations
are much more difficult.
\end{example}

\subsection{Induced Representations}

Suppose from now on that $G$ has a normal abelian subgroup $N\vartriangleleft G$,
and $G=H\ltimes N$ The $N\simeq\mathbb{R}^{d}$ and $N^{*}\simeq(\mathbb{R}^{d})^{*}=\mathbb{R}^{d}$.
In this case
\[
\chi_{t}(\nu)=e^{it\nu}
\]
for $\nu\in N$ and $t\in\hat{N}=N^{*}$. Notice that $\chi_{t}$
is a 1-d irreducible representation on $\mathbb{C}$. 

Let $\mathscr{H}_{t}$ be the space of functions $f:G\rightarrow\mathbb{C}$
so that
\[
f(\nu g)=\chi_{t}(\nu)f(g).
\]
On $\mathscr{H}_{t}$, define inner product so that
\[
\left\Vert f\right\Vert _{\mathscr{H}_{t}}^{2}:=\int_{G}\left|f(g)\right|^{2}=\int_{G/N}\left\Vert f(g)\right\Vert ^{2}dm
\]
where $dm$ is the invariant measure on $N\backslash G\simeq H$.

Define $U_{t}=ind_{N}^{G}(\chi_{t})\in Rep(G,\mathscr{H}_{t})$. Define
$U_{t}(g)f(x)=f(xg)$, for $f\in\mathscr{H}_{t}$. Notice that the
representation space of $\chi_{t}$ is $\mathbb{C}$, 1-d Hilbert
space; however, the representation space of $U_{t}$ is $\mathscr{H}_{t}$
which is infinite dimensional. $U_{t}$ is a family of irreducible
representations indexed by $t\in N\simeq\hat{N}\simeq\mathbb{R}^{d}$.
\begin{note}
Another way to recognize induced representations is to see these functions
are defined on $H$, not really on $G$.

Define the unitary transformation $W:\mathscr{H}_{t}\rightarrow L^{2}(H)$.
Notice that $H\simeq N\backslash G$ is a group, and it has an invariant
Haar\index{measure!Haar} measure. By uniqueness on the Haar measure,
this has to be $dm$. It would be nice to cook up the same space $L^{2}(H)$
so that all induced representations indexed by $t$ act on it. In
other words, this Hilbert space $L^{2}(H)$ does not depend on $t$.
$W_{t}$ is defined as
\[
WF_{t}(h)=F_{t}(h).
\]
\index{space!Hilbert-}

So what does the induced representation look like in $L^{2}(H)$ then?
Recall by definition that 
\[
U_{t}(g):=W\left(ind_{\chi_{t}}^{G}(g)\right)W^{*}
\]
and the following diagram commutes.

\[\xymatrix{
\mathcal{H}_{t} \ar[r]^{ind_{\chi_{t}}^{G}} \ar[d]_{W} & \mathcal{H}_{t} \ar[d]^{W} \\
L^{2}(H) \ar[r]^{U_{t}} & L^{2}(H) 
}\]

Let $f\in L^{2}(H)$. 
\begin{eqnarray*}
U_{t}\left(g\right)f\left(h\right) & = & W\left(ind_{\chi_{t}}^{G}\left(g\right)\right)W^{*}f\left(h\right)\\
 & = & \left(ind_{\chi_{t}}^{G}(g)W^{*}f\right)\left(h\right)\\
 & = & \left(W^{*}f\right)\left(hg\right).
\end{eqnarray*}
Since $G=H\ltimes N$, $g$ is uniquely decomposed into $g=g_{N}g_{H}$.
Hence $hg=hg_{N}g_{H}=g_{N}g_{N}^{-1}hg_{N}g_{H}=g_{N}\tilde{h}g_{H}$
and 
\begin{eqnarray*}
U_{t}(g)f(h) & = & (W^{*}f)(hg)\\
 & = & (W^{*}f)(g_{N}\tilde{h}g_{H})\\
 & = & \chi_{t}(g_{N})(W^{*}f)(\tilde{h}g_{H})\\
 & = & \chi_{t}(g_{N})(W^{*}f)(g_{N}^{-1}hg_{N}g_{H})
\end{eqnarray*}
This last formula is called the Mackey machine \cite{Mac52,Mac88}.
\index{Mackey machine}

The Mackey machine does not cover many important symmetry groups in
physics. Actually most of these are simple groups. However it can
still be applied. For example, in special relativity theory, we have
the Poincaré group $\mathcal{L}\ltimes\mathbb{R}^{4}$ where $\mathbb{R}^{4}$
is the normal subgroup. The baby version of this is when $\mathcal{L}=SL_{2}(\mathbb{R})$.
V. Bargman formulated this baby version. Wigner pioneered the Mackey
machine, long before Mackey was around. 

Once we get unitary representations, differentiate it and get selfadjoint
algebra of operators (possibly unbounded). These are the observables
in quantum mechanics.

\index{observable}\index{quantum mechanics!observable}\end{note}
\begin{example}
$\mathbb{Z}\subset\mathbb{R}$, $\hat{\mathbb{Z}}=T$. $\chi_{t}\in T$,
$\chi_{t}(n)=e^{itn}$. Let $\mathscr{H}_{t}$ be the space of functions
$f:\mathbb{R}\rightarrow\mathbb{C}$ so that 
\[
f(n+x)=\chi_{t}(n)f(x)=e^{int}f(x).
\]
Define inner product on $\mathscr{H}_{t}$ so that 
\[
\left\Vert f\right\Vert _{\mathscr{H}_{t}}^{2}:=\int_{0}^{1}\left|f(x)\right|^{2}dx.
\]
Define $ind_{\chi_{t}}^{\mathbb{R}}(y)f(x)=f(x+y)$. Claim that $\mathscr{H}_{t}\simeq L^{2}[0,1]$.
The unitary transformation is given by $W:\mathscr{H}_{t}\rightarrow L^{2}[0,1]$
\[
(WF_{t})(x)=F_{t}(x).
\]
Let's see what $ind_{\chi_{t}}^{\mathbb{R}}(y)$ looks like on $L^{2}[0,1]$.
For any $f\in L^{2}[0,1]$, 
\begin{eqnarray*}
\left(W\left(ind_{\chi_{t}}^{G}\left(y\right)\right)W^{*}f\right)\left(x\right) & = & \left(ind_{\chi_{t}}^{G}\left(y\right)W^{*}f\right)\left(x\right)\\
 & = & \left(W^{*}f\right)\left(x+y\right)
\end{eqnarray*}
Since $y\in\mathbb{R}$ is uniquely decomposed as $y=n+x'$ for some
$x'\in[0,1)$, therefore
\begin{eqnarray*}
\left(W\left(ind_{\chi_{t}}^{G}(y)\right)W^{*}f\right)(x) & = & (W^{*}f)(x+y)\\
 & = & (W^{*}f)(x+n+x')\\
 & = & (W^{*}f)(n+(-n+x+n)+x')\\
 & = & \chi_{t}(n)(W^{*}f)((-n+x+n)+x')\\
 & = & \chi_{t}(n)(W^{*}f)(x+x')\\
 & = & e^{itn}(W^{*}f)(x+x')
\end{eqnarray*}
\end{example}
\begin{note}
Are there any functions in $\mathscr{H}_{t}$? Yes, for example, $f(x)=e^{itx}$.
If $f\in\mathscr{H}_{t}$, $\left|f\right|$ is 1-periodic. Therefore
$f$ is really a function defined on $\mathbb{Z}\backslash\mathbb{R}\simeq[0,1]$.
Such a function has the form 
\[
f(x)=(\sum c_{n}e^{i2\pi nx})e^{itx}=\sum c_{n}e^{i(2\pi n+t)x}.
\]
Any 1-periodic function $g$ satisfies the boundary condition $g(0)=g(1)$.
$f\in\mathscr{H}_{t}$ has a modified boundary condition where $f(1)=e^{it}f(0)$. 
\end{note}
\index{boundary condition}

\section{Connections to Nelson's Spectral Theory}

In Nelson's notes \cite{Ne69}, a normal representation has the form
(counting multiplicity\index{multiplicity})\index{Nelson, E}
\[
\rho=\sum^{\oplus}n\pi\big|_{\mathscr{H}_{n}},\;\mathscr{H}_{n}\perp\mathscr{H}_{m}
\]
where 
\[
n\pi=\pi\oplus\cdots\oplus\pi\:(\mbox{n times})
\]
is a representation acting on the Hilbert space 
\[
\sum^{\oplus}K=l_{\mathbb{Z}_{n}}^{2}\otimes K.
\]
In matrix form, this is a diagonal matrix with $\pi$ repeated on
the diagonal $n$ times. $n$ could be $1,2,\ldots,\infty$. We apply
this to group representations.

Locally compact group can be divided into the following types.
\begin{itemize}
\item abelian
\item non-abelian: unimodular, non-unimodular
\item non-abelian: Mackey machine, semidirect product e.g. $H_{3}$, $ax+b$;
simple group $SL_{2}(\mathbb{R})$. Even it's called simple, ironically
its representation is much more difficult than the semidirect\index{product!semi-direct}
product case.
\end{itemize}
We want to apply these to group representations.\index{Spectral Theorem}\index{product!semidirect}

Spectral theorem says that given a normal operator $A$, we may define
$f(A)$ for quite a large class of functions, actually all measurable
functions (see \cite{Sto90, Yos95, Ne69, RS75, DS88b}). One way to
define $f(A)$ is to use the multiplication version of the spectral
theorem, and let \index{normal operator} 
\[
f(A)=\mathcal{F}f(\hat{A})\mathcal{F}^{-1}.
\]
The other way is to use the projection-valued measure version of the
spectral theorem, write 
\begin{eqnarray*}
A & = & \int\lambda P(d\lambda)\\
f(A) & = & \int f(\lambda)P(d\lambda).
\end{eqnarray*}
The effect is $\rho$ is a representation of the abelian algebra of
measurable functions onto operators action on some Hilbert space.
\begin{eqnarray*}
\rho:f\mapsto\rho(f) & = & f(A)\\
\rho(fg) & = & \rho(f)\rho(g)
\end{eqnarray*}
To imitate Fourier\index{transform!Fourier} transform, let's call
$\hat{f}:=\rho(f)$. Notice that $\hat{f}$ is the multiplication
operator.\index{algebras!group algebra}
\begin{example}
$G=(\mathbb{R},+)$, group algebra $L^{1}(\mathbb{R})$. Define Fourier
transform 
\[
\hat{f}(t)=\int f(x)e^{-itx}dx.
\]
$\{e^{itx}\}_{t}$ is a family of 1-dimensional irreducible representation
of $(\mathbb{R},+)$. 
\begin{example}
Fix $t$, $\mathscr{H}=\mathbb{C}$, $\rho(\cdot)=e^{it(\cdot)}\in Rep(G,\mathscr{H})$.
From the group representation $\rho$, we get a group algebra representation
$\tilde{\rho}\in Rep(L^{1}(\mathbb{R}),\mathscr{H})$ defined by
\[
\tilde{\rho}(f)=\int f(x)\rho(x)dx=\int f(x)e^{itx}dx
\]
It follows that 
\begin{eqnarray*}
\hat{f}(\rho) & := & \tilde{\rho}(f)\\
\widehat{f\star g} & = & \widehat{f\star g}=\hat{f}\hat{g}
\end{eqnarray*}
i.e. Fourier transform of $f\in L^{1}(\mathbb{R})$ is a representation
of the group algebra $L^{1}(\mathbb{R})$ on to the 1-dimensional
Hilbert space $\mathbb{C}$. The range of Fourier transform in this
case is 1-d abelian algebra of multiplication operators, multiplication
by complex numbers. 
\begin{example}
$\mathscr{H}=L^{2}(\mathbb{R})$, $\rho\in Rep(G,\mathscr{H})$ so
that 
\[
\rho(y)f(x):=f(x+y)
\]
i.e. $\rho$ is the right regular representation. The representation
space $\mathscr{H}$ in this case is infinite dimensional. From $\rho$,
we get a group algebra representation $\tilde{\rho}\in Rep(L^{1}(\mathbb{R}),\mathscr{H})$
where 
\[
\tilde{\rho}(f)=\int f(y)\rho(y)dy.
\]
Define 
\[
\hat{f}(\rho):=\hat{\rho}(f)
\]
then $\hat{f}(\rho)$ is an operator acting on $\mathscr{H}$.
\begin{alignat*}{1}
\hat{f}(\rho)g=\tilde{\rho}(f)g & =\int f(y)\rho(y)g(\cdot)dy\\
 & =\int f(y)(R_{y}g)(\cdot)dy\\
 & =\int f(y)g(\cdot+y)dy.
\end{alignat*}
If we have used the left regular representation, instead of the right,
then 
\begin{alignat*}{1}
\hat{f}(\rho)g=\tilde{\rho}(f)g & =\int f(y)\rho(y)g(\cdot)dy\\
 & =\int f(y)(L_{y}g)(\cdot)dy\\
 & =\int f(y)g(\cdot-y)dy.
\end{alignat*}
Hence $\hat{f}(\rho)$ is the left or right convolution operator.
\index{convolution}
\end{example}
\end{example}
\end{example}
Back to the general case. Given a locally compact group $G$, form
the group algebra $L^{1}(G)$, and define the left and right convolutions
as\index{operators!convolution-}

\begin{alignat*}{1}
(\varphi\star\psi)(x) & =\int\varphi(g)\psi(g^{-1}x)d_{L}g=\int\varphi(g)(L_{g}\psi)d_{L}g\\
(\varphi\star\psi)(x) & =\int\varphi(xg)\psi(g)d_{R}g=\int(R_{g}\varphi)\psi(g)d_{R}g
\end{alignat*}
Let $\rho(g)\in Rep(G,\mathscr{H})$, define $\tilde{\rho}\in Rep(L^{1}(G),\mathscr{H})$
given by 
\[
\tilde{\rho}(\psi):=\int_{G}\psi(g)\rho(g)dg
\]
and write 
\[
\hat{\psi}(\rho):=\tilde{\rho}(\psi).
\]
$\hat{\psi}$ is an analog of Fourier transform. If $\rho$ is irreducible,
the operators $\hat{\psi}$ forms an abelian algebra. In general,
the range of this generalized Fourier transform gives rise to a non
abelian algebra of operators.

For example, if $\rho(g)=R_{g}$ and $\mathscr{H}=L^{2}(G,d_{R})$,
then 
\[
\tilde{\rho}(\psi)=\int_{G}\psi(g)\rho(g)dg=\int_{G}\psi(g)R_{g}dg
\]
and
\begin{alignat*}{1}
\tilde{\rho}(\psi)\varphi & =\int_{G}\psi(g)\rho(g)\varphi dg=\int_{G}\psi(g)(R_{g}\varphi)dg\\
 & =\int_{G}\psi(g)\varphi(xg)dg\\
 & =(\varphi\star\psi)(x)
\end{alignat*}

\begin{example}
$G=H_{3}\sim\mathbb{R}^{3}$. $\hat{G}=\{\mathbb{R}\backslash\{0\}\}\cup\{0\}$.
$0\in\hat{G}$ corresponds to the trivial representation, i.e. $g\mapsto Id$
for all $g\in G$.
\[
\rho_{h}:G\rightarrow L^{2}(\mathbb{R})
\]
\[
\rho_{h}(g)f(x)=e^{ih(c+bx)}f(x+a)\simeq ind_{H}^{G}(\chi_{h})
\]
where $H$ is the normal subgroup $\{b,c\}$. It is not so nice to
work with $ind_{H}^{G}(\chi_{h})$ directly, so instead, we work with
the equivalent representations, i.e. Schrödinger representation. See
Folland's book on abstract harmonic analysis. \index{harmonic}
\[
\hat{\psi}(h)=\int_{G}\psi(g)\rho_{h}(g)dg
\]
Notice that $\hat{\psi}(h)$ is an operator acting on $L^{2}(\mathbb{R})$.
Specifically,
\begin{eqnarray*}
\hat{\psi}(h) & = & \int_{G}\psi(g)\rho_{h}(g)dg\\
 & = & \iiint\psi(a,b,c)e^{ih(c+bx)}f(x+a)dadbdc\\
 & = & \iint\left(\int\psi(a,b,c)e^{ihc}dc\right)f(x+a)e^{ihbx}dadb\\
 & = & \iint\hat{\psi}(a,b,h)f(x+a)e^{ihbx}dadb\\
 & = & \int\left(\int\hat{\psi}(a,b,h)e^{ihbx}db\right)f(x+a)da\\
 & = & \int\hat{\psi}(a,hx,h)f(x+a)da\\
 & = & \left(\hat{\psi}(\cdot,h\cdot,h)\star f\right)(x)
\end{eqnarray*}
Here the $\hat{\psi}$ on the right hand side in the Fourier transform
of $\psi$ in the usual sense. Therefore the operator $\hat{\psi}(h)$
is the one so that 
\[
L^{2}(\mathbb{R})\ni f\mapsto\left(\hat{\psi}(\cdot,h\cdot,h)\star f\right)(x).
\]
If $\psi\in L^{1}(G)$, $\hat{\psi}$ is not of trace class. But if
$\psi\in L^{1}\cap L^{2}$, then $\hat{\psi}$ is of trace class.
\[
\int_{\mathbb{R}\backslash\{0\}}^{\oplus}tr\left(\hat{\psi}^{*}(h)\hat{\psi}(h)\right)d\mu=\int\left|\psi\right|^{2}dg=\int\bar{\psi}\psi dg
\]
where $\mu$ is the Plancherel measure. 

If the group $G$ is non unimodular, the direct integral is lost (not
orthogonal). These are related to coherent states from physics, which
is about decomposing Hilbert into non orthogonal pieces.
\end{example}
\index{Schrödinger representation}

\index{observable}\index{representation!Schrödinger-}

Important observables in QM come in pairs (dual pairs). For example,
position - momentum; energy - time etc. The Schwartz space $S(\mathbb{R})$
has the property that $\widehat{S(\mathbb{R})}=S(\mathbb{R})$. We
look at the analog of the Schwartz space. $h\mapsto\hat{\psi}(h)$
should decrease faster than any polynomials. 

Take $\psi\in L^{1}(G)$, $X_{i}$ in the Lie algebra, form $\triangle=\sum X_{i}^{2}$.
Require that
\[
\triangle^{n}\psi\in L^{1}(G),\;\psi\in C^{\infty}(G).
\]
For $\triangle^{n}$, see what happens in the transformed domain.
Notice that
\[
\frac{d}{dt}\big|_{t=0}\left(R_{e^{tX}}\psi\right)=\tilde{X}\psi
\]
where $X\mapsto\tilde{X}$ represents the direction vector $X$ as
a vector field.

Let $G$ be any Lie group. $\varphi\in C_{c}^{\infty}(G)$, $\rho\in Rep(G,\mathscr{H})$.
\[
d\rho(X)v=\int(\tilde{X}\varphi)(g)\rho(g)vdg
\]
where 
\[
v=\int\varphi(g)\rho(g)wdg=\rho(\varphi)w.\mbox{ generalized convolution}
\]
If $\rho=R$, the n
\[
v=\int\varphi(g)R(g)
\]
\[
\tilde{X}(\varphi\star w)=(X\varphi)\star w.
\]

\begin{example}
$H_{3}$
\begin{eqnarray*}
a & \rightarrow & \frac{\partial}{\partial a}\\
b & \mapsto & \frac{\partial}{\partial b}\\
c & \mapsto & \frac{\partial}{\partial c}
\end{eqnarray*}
get standard Laplace operator. $\{\rho_{h}(\varphi)w\}\subset L^{2}(\mathbb{R})$
. $"="$ due to Dixmier. $\{\rho_{h}(\varphi)w\}$ is the Schwartz
space. 
\[
\left(\frac{d}{dx}\right)^{2}+(ihx)^{2}+(ih)^{2}=\left(\frac{d}{dx}\right)^{2}-(hx)^{2}-h^{2}
\]
Notice that
\[
-\left(\frac{d}{dx}\right)^{2}+(hx)^{2}+h^{2}
\]
is the Harmonic oscillator. Spectrum = $h\mathbb{Z}_{+}$. 
\end{example}
\index{harmonic}

\index{harmonic oscillator}

\index{operators!Laplace-}

\section{\label{sec:multirep}Multiplicity Revisited}

Let $\mathfrak{A}$ be a $*$-algebra, and let $\pi$ and $\rho$
be representations of $\mathfrak{A}$. To indicate the Hilbert space,
we write $\pi\in Rep\left(\mathfrak{A},\mathscr{H}_{\pi}\right)$,
and $\rho\in Rep\left(\mathfrak{A},\mathscr{H}_{\rho}\right)$.
\begin{defn}
Consider the following space of bounded linear operators $s,t:\mathscr{H}_{\pi}\longrightarrow\mathscr{H}_{\rho}$
which intertwine the respective representations, i.e., we have:
\begin{equation}
s\pi\left(a\right)=\rho\left(a\right)s,\;\forall a\in\mathfrak{A}.\label{eq:mr1}
\end{equation}

\end{defn}
The set of solutions $s$ to (\ref{eq:mr1}) forms a vector space,
and it is denoted $Int\left(\pi,\rho\right)$, the intertwining operators.
We check the following: 
\begin{equation}
s\in Int\left(\pi,\rho\right)\Longleftrightarrow s^{*}\in Int\left(\rho,\pi\right);\label{eq:mr2}
\end{equation}
and therefore, if $s,t\in Int\left(\pi,\rho\right)$, we have:
\begin{equation}
s^{*}t\in Int\left(\pi,\pi\right).\label{eq:mr3}
\end{equation}

Note
\begin{eqnarray}
Int\left(\pi,\pi\right) & = & \pi\left(\mathfrak{A}\right)'\;\left(\mbox{commutant}\right)\nonumber \\
 & = & \left\{ A\in\mathscr{B}\left(\mathscr{H}_{\pi}\right)\::\:\pi\left(a\right)A=A\pi\left(a\right),\;\forall a\in\mathfrak{A}\right\} .\label{eq:mr4}
\end{eqnarray}
If $\pi$ is irreducible, therefore $Int\left(\pi,\pi\right)$ is
one-dimensional; hence, for $\forall s,t\in Int\left(\pi,\rho\right)$,
\begin{equation}
s^{*}t=\left\langle s,t\right\rangle I_{\mathscr{H}_{\pi}},\label{eq:mr5}
\end{equation}
where $\left\langle s,t\right\rangle \in\mathbb{C}$ is uniquely determined.
This form $\left\langle \cdot,\cdot\right\rangle $ is sesquilinear,
and positive definite. We therefore get a Hilbert-completion of $Int\left(\pi,\rho\right)$.
Let $\mathscr{H}\left(\pi,\rho\right)$ be the corresponding Hilbert
space. \index{completion!Hilbert-}
\begin{defn}
\label{def:mr2}Let $\pi$ and $\rho$ be as above, assume that $\pi$
is irreducible, and let $\mathscr{H}\left(\pi,\rho\right)$ be the
corresponding Hilbert space; see (\ref{eq:mr5}). We say that $\pi$
occurs in $\rho$ $m$ times if
\begin{equation}
m=\dim\mathscr{H}\left(\pi,\rho\right).\label{eq:mr6}
\end{equation}
\end{defn}
\begin{xca}[Multiplicity]
\label{exer:mult-1}\myexercise{Multiplicity}Show that the definition
of multiplicity (\defref{mr2}) agrees with the one used inside \chapref{groups}.
\end{xca}

\begin{xca}[A Hilbert space of intertwiners]
\myexercise{A Hilbert space of intertwiners}Let $\pi$ and $\rho$
be as above, $\pi$ irreducible. Show that with the inner product
defined in (\ref{eq:mr5}), $Int\left(\pi,\rho\right)$ is a Hilbert
space.
\end{xca}

\begin{xca}[An ONB in $Int\left(\pi,\rho\right)$]
\myexercise{An ONB in $Int\left(\pi,\rho\right)$} Let $\left(s_{i}\right)$
be an ONB in $Int\left(\pi,\rho\right)$. 
\begin{enumerate}
\item Show that this is a system of isometries, satisfying:
\begin{equation}
s_{i}^{*}s_{j}=\delta_{i,j}I_{\mathscr{H}_{\pi}}.\label{eq:mr7}
\end{equation}

\item For $A\in\mathscr{B}\left(\mathscr{H}_{\pi}\right)$, set 
\begin{equation}
\alpha\left(A\right):=\sum_{i}s_{i}As_{i}^{*}.\label{eq:mr8}
\end{equation}
Show that 
\begin{eqnarray*}
\alpha\left(AB\right) & = & \alpha\left(A\right)\alpha\left(B\right),\;\mbox{and}\\
\alpha\left(A^{*}\right) & = & \alpha\left(A\right)^{*},\;\forall A,B\in\mathscr{B}\left(\mathscr{H}_{\pi}\right).
\end{eqnarray*}

\item What can be said about 
\[
\alpha\left(I_{\mathscr{H}_{\pi}}\right)=\sum_{i}s_{i}s_{i}^{*}\;?
\]

\end{enumerate}
\end{xca}

\begin{xca}[A Hilbert space of intertwining operators]
\myexercise{A Hilbert space of intertwining operators}Verify that
the results above about $Int\left(\pi,\rho\right)$ apply to \emph{unitary
representations} $\pi$, and $\rho$ of some given \emph{group} $G$;
i.e., with \index{representation!unitary}
\[
Int\left(\pi,\rho\right)=\left\{ s:\mathscr{H}_{\pi}\longrightarrow\mathscr{H}_{\rho}\::\:s\pi\left(g\right)=\rho\left(g\right)s,\;\forall g\in G\right\} .
\]
\uline{Hint}: Use the above on the group algebra $\mathfrak{A}_{G}:=\mathbb{C}\left[G\right]$. 
\end{xca}
Now consider the Heisenberg group $G$ of all $3\times3$ matrices\index{space!Hilbert-}
\[
g=\begin{bmatrix}1 & a & c\\
0 & 1 & b\\
0 & 0 & 1
\end{bmatrix},\;\left(a,b,c\right)\in\mathbb{R}^{3}.
\]
Recall its Haar measure is $dg=da\,db\,dc$ = 3-dimensional Lebesgue
measure.

Consider the following two representations $\rho$ and $\pi$ of $G$
(the regular representations RR, and the Schrödinger representation
SR):
\begin{itemize}
\item (RR) $\mathscr{H}_{\rho}=L^{2}\left(G,dg\right)$, Haar measure, and
\[
\left(\rho\left(g\right)f\right)\left(h\right)=f\left(hg\right),\;\forall f\in\mathscr{H}_{\rho},\:\forall g,h\in G.
\]
And the Schrödinger representation ($\hbar=1$):
\item (SR) $\mathscr{H}_{\pi}=L^{2}\left(\mathbb{R}\right)$, Lebesgue measure,
and 
\[
\left(\pi\left(g\right)F\right)\left(x\right)=e^{i\left(c+bx\right)}F\left(x+a\right),\;\forall F\in L^{2}\left(\mathbb{R}\right)=\mathscr{H}_{\pi},\;\forall g\in G,\;x\in\mathbb{R}.
\]
\end{itemize}
\begin{xca}[Specify the operators in $Int\left(\pi,\rho\right)$]
\myexercise{Specify the operators in $Int\left(\pi,\rho\right)$}Let
$G,\pi$, and $\rho$ be as above. What is the Hilbert space $Int\left(\pi,\rho\right)$? 
\end{xca}

\begin{xca}[A formula from Peter-Weyl \cite{Mac92}]
\label{exer:PW}\myexercise{A formula from Peter-Weyl} In case $G$
is a compact group, look up and explain that the Peter-Weyl theorem
states the following: If $\rho$ is the regular representation, and
if $\pi$ is irreducible unitary, then \index{Peter-Weyl theorem}
\[
\dim\left(Int\left(\pi,\rho\right)\right)=\dim\left(\pi\right).
\]
\end{xca}
\begin{rem}
An important class of non-compact, non-commutative, locally compact
groups $G$, and unitary representations $\pi$, for which the intertwining
Hilbert spaces $Int\left(\pi,\rho\right)$ are non-zero is the class
of square-integrable representations: Suppose the representation $\pi$
is irreducible and square-integrable, then $Int\left(\pi,\rho\right)$
is non-zero. Here $\rho$ denotes the regular representation of $G$.
A representation $\pi$ is square-integrable if its matrix coefficients
are in $L^{2}\left(G/Z\right)$, where $Z$ denotes the center of
$G$.
\end{rem}

\section*{A summary of relevant numbers from the Reference List}

For readers wishing to follow up sources, or to go in more depth with
topics above, we suggest: \cite{JO00,Mac52,Mac85,Mac92,JM84,JPS01,JPS05,Jor11,Jor94,Jor88,Szaar,DJ08,Dixmier198101,MR1895530,MR3220588,Nel59,JoMu80,Seg50,Ors79,Pou72,JLH06,DHL09,Tay86,Hal13,Hal15}.

\begin{subappendices}

\section{The Stone-von Neumann Uniqueness Theorem}

The ``uniqueness'' in the title above refers to ``uniqueness up
to \emph{unitary equivalence}.'' \index{Theorem!Stone-von Neumann uniqueness-}
\begin{defn}
Let $\mathscr{H}_{i}$, $i=1,2$ be two Hilbert spaces, and let $S_{1}=\{A_{\alpha}\}\subset\mathscr{B}\left(\mathscr{H}_{1}\right)$,
and $S_{2}=\{B_{\alpha}\}\subset\mathscr{B}\left(\mathscr{H}_{2}\right)$
be systems of bounded operators, where the index set $J=\left\{ \alpha\right\} $
is the same for the two operator systems. 

We say that $S_{1}$ and $S_{2}$ are \emph{unitarily equivalent}
iff (Def) $\exists W:\mathscr{H}_{1}\rightarrow\mathscr{H}_{2}$,
$W$ a unitary isomorphism of $\mathscr{H}_{1}$ onto $\mathscr{H}_{2}$
such that 
\begin{equation}
WA_{\alpha}=B_{\alpha}W,\;\forall\alpha\in J;\;\mbox{see Fig. }\ref{fig:vn}.\label{eq:vu1}
\end{equation}

We say that the system $S_{1}=\{A_{\alpha}\}$ is \emph{irreducible}
iff (Def) the following implication holds
\begin{equation}
\boxed{T\in\mathscr{B}\left(\mathscr{H}_{1}\right),\;TA_{\alpha}=A_{\alpha}T,\;\alpha\in J}\Longrightarrow T=\lambda I_{1},\;\mbox{for some}\;\lambda\in\mathbb{C};\label{eq:vu3}
\end{equation}
i.e., the commutant is one-dimensional. 
\end{defn}
\begin{figure}[H]
\[
\xymatrix{\mathscr{H}_{1}\ar[rr]^{W}\ar[dd]_{A_{\alpha}} &  & \mathscr{H}_{2}\ar[dd]^{B_{\alpha}}\\
 & B_{\alpha}W=WA_{\alpha}\\
\mathscr{H}_{1}\ar[rr]_{W} &  & \mathscr{H}_{2}
}
\]

\protect\caption{\label{fig:vn}$W$ intertwines $S_{1}$ and $S_{2}$.}
\end{figure}

\begin{defn}
The Heisenberg group $G_{3}$ is the matrix group
\[
g=\begin{bmatrix}1 & a & c\\
0 & 1 & b\\
0 & 0 & 1
\end{bmatrix},\quad\left(a,b,c\right)\in\mathbb{R}^{3},
\]
of all upper triangular $3\times3$ matrices. 

Fix $h\in\mathbb{R}\backslash\left\{ 0\right\} $, and set 
\begin{equation}
\left(\mathcal{U}_{h}\left(g\right)f\right)\left(x\right)=e^{ih\left(c+bx\right)}f\left(x+a\right)\label{eq:vu4}
\end{equation}
$\forall g=\left(a,b,c\right)\in G_{3}$, $\forall f\in L^{2}\left(\mathbb{R}\right)$,
$\forall x\in\mathbb{R}$.
\end{defn}
It is easy to see that $\mathcal{U}_{h}$ is a unitary irreducible
representation of $G_{3}$ acting on $L^{2}\left(\mathbb{R}\right)$,
i.e., $\mathcal{U}_{h}\in Rep_{uni}\left(G_{3},L^{2}\left(\mathbb{R}\right)\right)$
for all $h\in\mathbb{R}\backslash\left\{ 0\right\} $. It is called
the \emph{Schrödinger representation}.\index{representation!Schrödinger-}
\begin{thm}[Stone-von Neumann]
 Every unitary irreducible representation of $G_{3}$ in a Hilbert
space (other than the trivial one-dimensional representation) is unitarily
equivalent to the Schrödinger representation $\mathcal{U}_{h}$ for
some $h\in\mathbb{R}\backslash\left\{ 0\right\} $. \end{thm}
\begin{proof}
The proof follows from the more general result, \thmref{imp} above;
the Imprimitivity Theorem. Also see \cite{vN32c,vN31}.\end{proof}
\begin{rem}
The center of $G_{3}$ is the one-dimensional subgroup $g=\left(0,0,c\right)$,
$c\in\mathbb{R}$, and so if $\mathcal{U}_{h}\in Rep_{uni}\left(G_{3},\mathscr{H}\right)$,
$\dim\mathscr{H}>1$, then it follows that $\exists!h\in\mathbb{R}\backslash\left\{ 0\right\} $
such that 
\[
\mathcal{U}\left(0,0,c\right)=e^{i\,c\,h}I_{\mathscr{H}}.
\]
Hence $\mathcal{U}$ is determined by two one-parameter groups
\begin{equation}
\begin{cases}
\mathcal{U}_{1}\left(a\right)=\mathcal{U}\left(a,0,0\right), & a\in\mathbb{R};\:\mbox{and}\\
\mathcal{U}_{2}\left(b\right)=\mathcal{U}\left(0,b,0\right), & b\in\mathbb{R}
\end{cases}\label{eq:vu5}
\end{equation}
such that
\begin{equation}
\mathcal{U}_{1}\left(a\right)\mathcal{U}_{2}\left(b\right)\mathcal{U}_{1}\left(-a\right)=e^{i\,h\,a\,b}\mathcal{U}_{2}\left(b\right),\;\forall a,b\in\mathbb{R}.\label{eq:vu6}
\end{equation}

The system (\ref{eq:vu6}) is called the \emph{Weyl commutation relation}.
It is the integrated form of the corresponding Heisenberg relation
(for unbounded essentially selfadjoint operators). (We omit a systematic
discussion of the interrelationships between the two commutation relations.)
\index{Weyl commutation relation}

Under the unitary equivalence $W:\mathscr{H}\rightarrow L^{2}\left(\mathbb{R}\right)$
from the Stone-von Neumann theorem, we get 
\begin{equation}
\begin{cases}
\left(W\mathcal{U}_{1}\left(a\right)W^{*}f\right)\left(x\right)=f\left(x+a\right), & \mbox{and}\\
\left(W\mathcal{U}_{2}\left(b\right)W^{*}f\right)\left(x\right)=e^{i\,h\,b}f\left(x\right), & \forall a,b,x\in\mathbb{R},\:\forall f\in L^{2}\left(\mathbb{R}\right).
\end{cases}\label{eq:vu7}
\end{equation}

We shall use the following \end{rem}
\begin{lem}
Let $\mathcal{U}_{1}\left(\cdot\right)$ and $\mathcal{U}_{2}\left(\cdot\right)$
be the two one-parameter groups from the Weyl relation (\ref{eq:vu6}),
and let $P_{2}$ be the projection valued measure corresponding to
$\{\mathcal{U}_{2}\left(b\right)\}_{b\in\mathbb{R}}$, i.e., 
\begin{equation}
\mathcal{U}_{2}\left(b\right)=\int_{\mathbb{R}}e^{i\,b\,\lambda}P_{2}\left(d\lambda\right),\;\forall b\in\mathbb{R}.\label{eq:vu9}
\end{equation}
Then the Weyl relation (\ref{eq:vu6}) is equivalent to 
\begin{equation}
\mathcal{U}_{1}\left(a\right)P_{2}\left(\triangle\right)\mathcal{U}_{1}\left(-a\right)=P_{2}\left(\triangle-h\,a\right),\label{eq:vu10}
\end{equation}
$\forall a\in\mathbb{R}$, $\forall\triangle\in\mathcal{B}\left(\mathbb{R}\right)$,
where 
\begin{equation}
\triangle-h\,a=\left\{ s-h\,a\:\big|\:s\in\triangle\right\} .\label{eq:vu11}
\end{equation}
\end{lem}
\begin{proof}
The proof is an easy application of Stone's theorem (\cite{vN32c,Ne69});
see Appendix \ref{sec:stone}.
\end{proof}
\end{subappendices}

\chapter{The Kadison-Singer Problem\label{chap:KS}}
\begin{quotation}
Born wanted a theory which would generalize these matrices or grids
of numbers into something with a continuity comparable to that of
the continuous part of the spectrum. The job was a highly technical
one, and he counted on me for aid.... I had the generalization of
matrices already at hand in the form of what is known as operators.
Born had a good many qualms about the soundness of my method and kept
wondering if Hilbert would approve of my mathematics. Hilbert did,
in fact, approve of it, and operators have since remained an essential
part of quantum theory. 

--- Norbert Wiener\sindex[nam]{Wiener, N., (1894-1964)}\vspace{1em}

In science one tries to tell people, in such a way as to be understood
by everyone, something that no one ever knew before. But in the case
of poetry, it's the exact opposite! 

--- Paul Adrien Maurice Dirac. \vspace{1em}

It seems to be one of the fundamental features of nature that fundamental
physical laws are described in terms of a mathematical equations of
great beauty and power. 

--- Paul Adrien Maurice Dirac\sindex[nam]{Dirac, P.A.M., (1902-1984)}\vspace{2em}
\end{quotation}
\textbf{The Kadison-Singer Problem:} Does every pure state on the
(abelian) von Neumann algebra $\mathbb{D}$ of all bounded diagonal
operators on $l^{2}$ have a unique extension to a pure state on all
$\mathscr{B}(l^{2})$, the von Neumann algebra of all bounded operators
on $l^{2}$? \index{Kadison-Singer problem}

We shall begin by explaining the meaning, and the significance, of
the terms used in the statement of the Kadison-Singer (abbreviated
KS) problem. But on the whole, our discussion of the KS problem (or
conjecture) in the present book will be modest in scope. The first
to say is that it was just solved by Adam Marcus, Dan Spielman, and
N. Srivastava, see \cite{MSS15}. There are several reasons for why
we cannot go into proof-details for the solution: While the formulation
of KS by Kadison and Singer in 1959 was in the language of operator
algebras, as that subject back then was inspired by Dirac's quantum
theory, it turned out that the eventual solution to KS, five decades
later \cite{MSS15}, by Adam Marcus, Dan Spielman, and N. Srivastava,
surprisingly, involves themes that draw on new topics, quite outside
the scope of the present book. And making the connection between tools
from the 2015 solution, back to the original 1959 formulation of KS,
is not at all trivial. The new tools employed by Marcus, Spielman,
and Srivastava are from diverse mathematical areas, and with a heavy
combinatorial component, and also involving mathematical notions which
we have not even defined here; for example: interlacing families,
random vectors, paving, probabilistic frames, discrepancy analysis,
sparsification, \dots{} . We hope readers will find it interesting
to see how apparently disparate areas can meet at the crossroads in
the solution of a famous problem in mathematics\footnote{The many interconnections between the disparate areas of mathematics
coming together in the proof of KS are just emerging in the literature
as of this point; in particular, there is a forthcoming paper by P.
Casazza, M. Bownik, A. Marcus and D. Speegle; which promises to be
an authoritative source. We are grateful to P. Casazza for updates
on KS. On the same theme, see also the paper \textquotedblleft Consequences
of the Marcus/Spielman/Srivastava solution of the Kadison-Singer problem,\textquotedblright{}
By Peter G. Casazza and Janet C. Tremain. arXiv:1407.4768v2.}. In view of this, we stress that even a modest attempt on our part
at going into a detailed discussion of the Marcus-Spielman-Srivastava
solution to KS would take us far afield; and done properly it could
easily become a separate book volume. Since the Marcus-Spielman-Srivastava
paper has just now appeared in the Annals 2015 \cite{MSS15}, it may
also be too early for a proper book presentation. Nonetheless, we
feel that a discussion here, in the present chapter, \emph{of the
original formulation} of KS is in fact appropriate. Indeed, the initial
motivation for KS derives from precisely the topics which central
themes of our present book: Operator theory/algebra, positivity, states,
spectral theory; and with how these mathematical themes intersect
with quantum theory.\index{Kadison-Singer problem}

The authors of \cite{MR3246256,MSS15} have proved the Kadison-Singer
conjecture in an indirect way, by proving instead Weaver\textquoteright s
conjecture \cite{MR2035401,MR3156130}.
\begin{conjecture}[$KS_{2}$]
 There exist universal constants $\eta\geq2$ and $\theta>0$ so
that the following holds. Let $w_{1},\ldots,w_{m}\in\mathbb{C}^{d}$
satisfy $\left\Vert w_{i}\right\Vert \leq1$ for all $i$, and suppose
\begin{equation}
\sum_{i=1}^{m}\left|\left\langle w_{i},u\right\rangle \right|^{2}=\eta\label{eq:ks2-1}
\end{equation}
for every unit vector $u\in\mathbb{C}^{d}$. Then there exists a partition
$S_{1}$, $S_{2}$ of $\left\{ 1,\ldots,m\right\} $ so that 
\[
\sum_{i\in S_{j}}\left|\left\langle w_{i},u\right\rangle \right|^{2}\leq\eta-\theta
\]
for every unit vector $u\in\mathbb{C}^{d}$ and each $j\in\left\{ 1,2\right\} $. 
\end{conjecture}
Akemann and Anderson's projection paving conjecture \cite[Conj. 7.1.3]{MR1086563}
follows directly from $KS_{2}$ (see \cite[p. 229]{MR2035401}) .

Anderson's original paving conjecture says:
\begin{conjecture}[Anderson Paving]
For every $\varepsilon>0$, there is an $r\in\mathbb{N}$ such that
for every $n\times n$ Hermitian matrix $T$ with zero diagonal, there
are diagonal projections $P_{1},\cdots,P_{r}$ with $\sum_{i=1}^{r}P_{i}=I$
such that 
\[
\left\Vert P_{i}TP_{i}\right\Vert \leq\varepsilon\left\Vert T\right\Vert ,\quad\mbox{for \ensuremath{i=1,\ldots,r.}}
\]

\end{conjecture}

The Kadison-Singer problem (KS) lies at the root of how questions
from quantum physics take shape in the language of functional analysis,
and algebras of operators. 

A brief sketch is included below, summarizing some recent advances
(in fact the KS-problem was recently solved.) It is known that the
solution to KS at the same time answers a host of other questions;
all with applications to engineering, especially to signal processing.
The notion from functional analysis here is \textquotedblleft frame.\textquotedblright{}
A frame of vectors in Hilbert space generalizes the notion of orthonormal
basis in Hilbert space. \index{frame} \index{orthonormal basis (ONB)}\index{signal processing}

The Kadison-Singer problem (KS) comes from functional analysis, but
it was resolved (only recently) with tools from areas of mathematics
quite disparate from functional analysis. More importantly, the solution
to KS turned out to have important implications for a host of applied
fields from engineering.\footnote{Atiyah and Singer shared the Abel prize of 2004.}\index{algebras!$C^{*}$-algebra}

This reversal of the usual roles seem intriguing for a number of reasons: 

While the applications considered so far involve problems which in
one way or the other, derive from outside functional analysis itself,
e.g., from physics, from signal processing, or from anyone of a number
of areas of analysis, PDE, probability, statistics, dynamics, ergodic
theory, prediction theory etc.; the Kadison-Singer problem is different.
It comes directly from the foundational framework of functional analysis;
more specifically from the axiomatic formulation of $C^{*}$-algebras.
Then $C^{*}$-algebras are a byproduct of a rigorous formulation of
quantum theory, as proposed by P.A.M. Dirac.\footnote{P.A.M. Dirac gave a lecture at Columbia University in the late 1950's,
in which he claimed without proof that pure states on the algebra
of diagonal operators ($\simeq l^{\infty}$) extends uniquely on $\mathscr{B}(l^{2})$.
Kadison and Singer sitting in the audience were skeptical about whether
Dirac knew what it meant to be an extension. They later formulated
the conjecture in a joint paper, made precise the difference between
MASAs that are continuous vs discrete. They showed that non-uniqueness
holds in the continuous case.}\index{signal}\index{axioms}

\begin{figure}
\begin{tabular}{|>{\centering}p{0.55\textwidth}|>{\centering}p{0.35\textwidth}|}
\hline 
\vspace{0pt}Physics & \vspace{0pt}Mathematics\tabularnewline
\hline 
\vspace{0pt}

\includegraphics[scale=0.6]{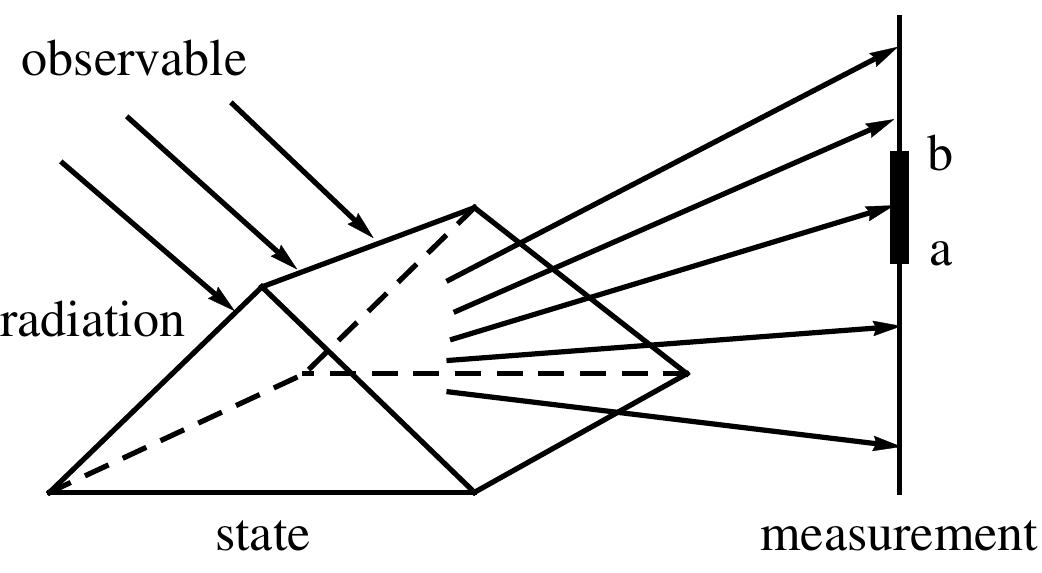} & \vspace{0pt}

$H^{*}=H$, $v\in\mathscr{H}$, $\left\Vert v\right\Vert =1$ (state).
Spectral theorem: $H\sim P_{H}\left(\cdot\right)$ projection-valued\index{measure!projection-valued}
measure.\vspace{1em} \uline{Measurement}:\textbf{ $Prob\left(H\in\left(a,b\right)\right)=\left\Vert P_{H}\left(a,b\right)v\right\Vert ^{2}$} \tabularnewline
\hline 
\end{tabular}

\protect\caption{\label{fig:qmm}Observable, state, measurement. Left column: An idealized
physics experiment. Right: the mathematical counterpart, a selfadjoint
operator $H$, its associated projection-valued measure $P_{H}$,
and a norm-one vector $v$ in Hilbert space.}
\end{figure}

From quantum theory, we have such notions as \emph{state}, \emph{observable}\index{observable},
and \emph{measurement}\index{measurement}. See \figref{qmm}. But
within the framework of $C^{*}$-algebras, each of these same terms,
``state'', ``observable'', and ``measurement'' also has a purely
mathematical definition, see  \secref{qmview} in \chapref{sp}. Indeed
$C^{\ast}$-algebra theory was motivated in part by the desire to
make precise fundamental and conceptual questions in quantum theory,
e.g., the uncertainty\index{uncertainty} principle\index{Heisenberg, W.K.!uncertainty principle}\index{quantum mechanics!uncertainty principle},
measurement, determinacy, hidden variables, to mention a few (see
for example \cite{EMC00}). The interplay between the two sides has
been extraordinarily fruitful since the birth of quantum mechanics
in the 1920ties. 

Cited from \cite{MR0123922}: \index{state!pure-}
\begin{quotation}
\textquotedblleft The main concern of this paper is the problem of
uniqueness of extensions of pure states from maximal abelian self-adjoint
algebras of operators on a Hilbert space to the algebra of all bounded
operators on that space. The answer, as many of us have suspected
for several years, is in negative.'' ... ``We heard of it first
from I.E. Segal and I. Kaplansky, though it is difficult to credit
a problem which stems naturally from the physical interpretation and
the inherent structure of a subject. This problem has arisen, in one
form or another, in our work on several different occasions;...''
\end{quotation}
Now consider the following: (i) the Hilbert space $\mathscr{H}=l^{2}(=l^{2}(\mathbb{N}))$,
all square summable sequences, (ii) the $C^{*}$-algebra $\mathscr{B}(l^{2})$
of all bounded operators on $l^{2}$, and finally (iii) the sub-algebra
$\mathfrak{A}$ of $\mathscr{B}(l^{2})$ consisting of all diagonal
operators, so an isomorphic copy of $l^{\infty}$. 

The Kadison-Singer problem (KS), in the discrete version, is simply
this:

\emph{Does every pure state of $\mathfrak{A}$ have a unique pure-state
extension to $\mathscr{B}(l^{2})$? }

We remark that existence (of a pure-state extension) follows from
the main theorems from functional analysis of Krein and Krein-Milman,
but the uniqueness is difficult. The difficulty lies in the fact that
it's hard to find all states on $l^{\infty}$, i.e., a dual of $l^{\infty}$.
The pure states of $\mathfrak{A}$ are points in the Stone-\v{C}ech
compactification $\beta\left(\mathbb{N}\right)$. The problem was
settled in the affirmatively (uniqueness in the discrete case) only
a year ago, after being open for 50 years. \index{space!state-}

\index{extension!-of state}

\index{quantum mechanics!observable}

\index{quantum mechanics!measurement}

\index{Krein-Milman}

\index{dual}

\index{state!pure-}

\index{space!dual-}
\begin{lem}
\label{lem:ps}Pure normal states on $\mathscr{B}(\mathscr{H})$ are
unit vectors (in fact, the equivalent class of unit vectors.\footnote{Equivalently, pure states sit inside the projective vector space.
If $\mathscr{H}=\mathbb{C}^{n+1}$, pure states is $\mathbb{C}P^{n}$. }) Specifically, let $u\in\mathscr{H}$, $\left\Vert u\right\Vert =1$,
then 
\[
\mathscr{B}(\mathscr{H})\ni A\longmapsto\omega_{u}\left(A\right)=\left\langle u,Au\right\rangle 
\]
is a pure state. All normal pure states on $\mathscr{B}(\mathscr{H})$
are of this form. 
\end{lem}
The pure states on $\mathscr{B}\left(\mathscr{H}\right)$ not of the
form $\omega_{u}$, for $u\in\mathscr{H}$, $\left\Vert u\right\Vert =1$,
are called \emph{singular pure states}.
\begin{rem}
Since $l^{\infty}$ is an abelian algebra Banach $*$-algebra, by
Gelfand's theorem, $l^{\infty}\simeq C(X)$ where $X$ is a compact
Hausdorff space. Indeed, $X=\beta\mathbb{N}$, -- the Stone-\v{C}ech
compactification of $\mathbb{N}$. Points in $\beta\mathbb{N}$ are
called \emph{ultra-filters}. Pure states on $l^{\infty}$ correspond
to pure states on $C(\beta\mathbb{N})$, i.e., Dirac-point measures
on $\beta\mathbb{N}$. 

\index{Stone-v{C}ech compactification@Stone-\v{C}ech compactification}\index{compact}\index{measure!Dirac}

Let $s$ be a pure state on $l^{\infty}$. Using Hahn-Banach theorem
one may extend $s$, as a linear functional, from $l^{\infty}$ to
$\tilde{s}$ on the Banach space $\mathscr{B}(\mathscr{H})$. However,
Hahn-Banach theorem doesn't guarantee the extension remains a \emph{pure
state}. Let $E(s)$ be the set of all states on $\mathscr{B}(\mathscr{H})$
which extend $s$. $E(s)$ is non-empty, compact and convex in the
weak $*$-topology\index{weak{*}-topology}. By Krein-Milman's theorem,
$E(s)=\mbox{closure(Extreme Points)}$. Any extreme point will then
be a pure state extension of $s$; but which one to choose? It's the
uniqueness part that is the famous KS problem.\index{Theorem!Krein-Milman's-}
\end{rem}
\index{Hahn-Banach theorem}\index{measurement}\index{space!Banach-}\index{Theorem!Spectral-}\index{space!state-}\index{selfadjoint operator}
\begin{xca}[Non-normal pure states on $\mathscr{B}\left(l^{2}\right)$]
\myexercise{Non-normal pure states on $\mathscr{B}(l^{2})$}Show
that there are pure states on $\mathscr{B}\left(l^{2}\right)$ which
do not have the form given in \lemref{ps}. 

\begin{flushleft}
\uline{Hint}: \index{pure state!non-normal}
\par\end{flushleft}
\begin{enumerate}
\item[Step 1.]  The states listed in \lemref{ps} have cardinality $c=2^{\aleph_{0}}$.
\item[Step 2.]  The pure states of $C\left(\beta\left(\mathbb{N}\right)\right)$
are given by points in $\beta\left(\mathbb{N}\right)$, and the cardinality
of $\beta\left(\mathbb{N}\right)$ is 
\begin{equation}
2^{2^{\aleph_{0}}}>c.\label{eq:card}
\end{equation}

\item[Step 3.]  By Krien-Milman, every pure state on $\mathscr{D}$ ($\simeq$ $l^{2}\left(\mathbb{N}\right)$)
has a pure state extension to $\mathscr{B}\left(l^{2}\right)$.
\item[Step 4.]  Use (\ref{eq:card}) in step 2 to conclude that some of these pure
state extensions to $\mathscr{B}\left(l^{2}\right)$ are \emph{not}
of the form given in \lemref{ps}. \index{state!pure-}
\end{enumerate}
\end{xca}

\begin{xca}[The Calkin algebra and Non-normal states on $\mathscr{B}\left(l^{2}\right)$]
\myexercise{The Calkin algebra and Non-normal states on $\mathscr{B}(l^{2})$}Let
$\mathscr{K}\subset\mathscr{B}\left(l^{2}\right)$ be the ideal of
all compact operators in $l^{2}$; then the quotient \index{ideal}\index{state!non-normal}
\[
\mathscr{C}:=\mathscr{B}\left(l^{2}\right)/\mathscr{K}
\]
is called the \emph{Calkin algebra}. Show that the quotient is a $C^{*}$-algebra. 

\uline{Hint}: Be careful in defining its $C^{*}$-norm.
\end{xca}
\index{algebra!Calkin-}

\index{algebra!quotient-}

\index{state!normal-}
\begin{xca}[The pure state $\varphi=s\circ\pi$ on $\mathscr{B}\left(l^{2}\right)$]
\myexercise{The pure state $\varphi=s\circ\pi$ on $\mathscr{B}\left(l^{2}\right)$}Let
$\pi:\mathscr{B}\left(l^{2}\right)\longrightarrow\mathscr{C}$ be
the natural quotient mapping, and let $s$ be a pure state on $\mathscr{C}$.
Show that the composition 
\[
\varphi:=s\circ\pi\quad\left(\mbox{see Fig }\ref{fig:psc}.\right)
\]
is a pure state on $\mathscr{B}\left(l^{2}\right)$, and that $\varphi$
does not have the form in  \lemref{ps}. 

\uline{Hint}: Suppose to the contrary, i.e., suppose $\exists x\in l^{2}$,
$\left\Vert x\right\Vert =1$ such that 
\begin{equation}
\varphi\left(A\right)=\left\langle x,Ax\right\rangle =\omega_{x}\left(A\right),\;\forall A\in\mathscr{B}\left(l^{2}\right).\label{eq:psl1}
\end{equation}
We have $\omega_{x}\left(\left|x\left\rangle \right\langle x\right|\right)=1$,
but $\left|x\left\rangle \right\langle x\right|\in\mathscr{K}$, so
$\varphi\left(\left|x\left\rangle \right\langle x\right|\right)=s\left(0\right)=0$;
a contradiction. Hence (\ref{eq:psl1}) cannot hold for any state-vector
$x\in l^{2}$. \index{pure state!normal}
\end{xca}
\begin{figure}
\[
\xymatrix{\mathscr{B}\left(l^{2}\right)\ar@/^{1pc}/[rr]^{\pi}\ar@/_{1pc}/[dr]_{\varphi} &  & \mathscr{C}\ar@/^{1pc}/[dl]^{s\quad\text{(pure)}}\\
 & \mathbb{C}
}
\]

\protect\caption{\label{fig:psc}The pure state $\varphi=s\circ\pi$ on $\mathscr{B}\left(l^{2}\right)$}
\end{figure}

\begin{xca}[The Stone\textendash \v{C}ech compactification]
\myexercise{The Stone-\v{C}ech compactification}\label{exer:ks1}Extend
$+$ on $\mathbb{N}$ to a ``$+$'' on $\beta\mathbb{N}$ (the Stone-\v{C}ech
compactification\index{Stone-v{C}ech compactification@Stone-\v{C}ech compactification}).

\uline{Hint}: 
\begin{enumerate}
\item For subsets $A\subset\mathbb{N}$, and $n\in\mathbb{N}$, set $A-n:=\left\{ k\in\mathbb{N}\:\big|\:k+n\in A\right\} $. 
\item Let $F$ and $G$ be ultra-filters on $\mathbb{N}$, and set
\begin{equation}
F+G:=\left\{ A\subset\mathbb{N}\:\big|\:\left\{ n\in\mathbb{N};A-n\in F\right\} \in G\right\} .\label{eq:ks1}
\end{equation}

\item Show that $F+G$ is an ultra-filter.
\item Show that the ``addition'' operation ``+'' in (\ref{eq:ks1})
is an operation on $\beta\mathbb{N}$, i.e., $\beta\mathbb{N}\times\beta\mathbb{N}\longrightarrow\beta\mathbb{N}$
which is associative, but not commutative, i.e., $F+G\neq G+F$ may
happen.
\item Fix $F\in\beta\mathbb{N}$, and show that 
\[
\beta\mathbb{N}\ni G\longmapsto F+G\in\beta\mathbb{N}
\]
 is continuous, where $F+G$ is defined in (\ref{eq:ks1}).
\end{enumerate}
\end{xca}

\emph{\uline{Ultra-filters}} define pure states of $l^{\infty}$
as follows: If $\left(x_{n}\right)_{n\in\mathbb{N}}\in l^{\infty}$,
and if $F\in\beta\mathbb{N}$, i.e., is an ultra-filter, then there
is a well-defined limit
\[
\lim_{F}x_{n}=\varphi_{F}\left(x\right);
\]
and this defines $\varphi_{F}$ as a state on $l^{\infty}$.

\section{The Dixmier Trace}

A related use of ultra-filters yield the famous \emph{Dixmier-trace}.
For this we need ultra-filters $\omega$ on $\mathbb{N}$ with the
following properties: 

\index{trace}

\index{Dixmier trace}

(i) $x_{n}\geq0\Longrightarrow\lim_{\omega}x_{n}\geq0$.

(ii) If $x_{n}$ is convergent with limit $x$, then $\lim_{\omega}x_{n}=x$.

(iii) For $n\in\mathbb{N}$, set 
\[
\sigma_{N}\left(x\right)=\Big(\underset{N\;\text{times}}{\underbrace{x_{1}\cdots x_{1}}},\underset{N\;\text{times}}{\underbrace{x_{2}\cdots x_{2}}},\underset{N\;\text{times}}{\underbrace{x_{3}\cdots x_{3}}},\cdots\Big)
\]
then $\lim_{\omega}\left(x_{n}\right)=\lim_{\omega}\left(\sigma_{N}\left(x\right)\right)$. 

Let $A$ be a compact operator, and assume the eigenvalues $\lambda_{k}$
of $\left|A\right|=\sqrt{A^{*}A}$ as 
\[
\lambda_{1}\geq\lambda_{2}\geq\cdots,\lambda_{k}=\lambda_{k}\left(A\right)
\]
and set\index{eigenvalue} 
\begin{equation}
tr_{Dix,\omega}\left(A\right)=\lim_{\omega}\frac{1}{\log\left(n+1\right)}\sum_{k=1}^{n}\lambda_{k}\left(A\right).\label{eq:ks2}
\end{equation}

\begin{defn}
We say that $A$ has finite Dixmier trace if the limit in (\ref{eq:ks2})
is finite.\end{defn}
\begin{xca}[The Dixmier trace]
\myexercise{The Dixmier trace}\label{exer:ks2}Show that (\ref{eq:ks2})
is well-defined and that:
\begin{enumerate}
\item $A\longmapsto tr_{Dix,\omega}\left(A\right)$ is linear, and positive.
\item $tr_{Dix,\omega}\left(AB\right)=tr_{Dix,\omega}\left(BA\right)$ holds
if $B$ is bounded, and $A$ has finite Dixmier trace. 
\item If $\sum_{k}\lambda_{k}\left(A\right)<\infty$, then $tr_{Dix,\omega}\left(A\right)=0$.
\end{enumerate}
\end{xca}

\section{Frames in Hilbert Space}

The proof of the KS-problem involves systems of vectors in Hilbert
space called \emph{frames}. For details we refer to \cite{Cas13}.
\index{frame}

Below we include a sketch with some basic fact about frames; also
called \textquotedblleft generalized bases'', see \defref{frame}
below. The general idea is that a \textquotedblleft frame expansion\textquotedblright{}
inherits some (but not all) attractive properties and features of
expansions in ONBs. In frame-analysis, this then offers the desirable
feature of more flexibility in a host of applications; see e.g., \cite{MR2830604}
and \cite{MR1381391,MR2367342}. But we also give up something. For
example, by contrast to what holds for an ONB, non-uniqueness is a
fact of life for frame expansions.

Let $\mathscr{H}$ be a separable Hilbert space, and let $\left\{ u_{k}\right\} _{k\in\mathbb{N}}$
be an ONB, then we have the following unique representation
\begin{equation}
w=\sum_{k\in\mathbb{N}}\left\langle u_{k},w\right\rangle _{\mathscr{H}}u_{k}\label{eq:fr1}
\end{equation}
valid for all $w\in\mathscr{H}$. Moreover, 
\begin{equation}
\left\Vert w\right\Vert _{\mathscr{H}}^{2}=\sum_{k\in\mathbb{N}}\left|\left\langle u_{k},w\right\rangle _{\mathscr{H}}\right|^{2},\label{eq:fr2}
\end{equation}
the Parseval-formula.
\begin{defn}
\label{def:frame}A system $\left\{ v_{k}\right\} _{k\in\mathbb{N}}$
in $\mathscr{H}$ is called a \emph{\uline{frame}} if there are
constants $A,B$ such that $0<A\leq B<\infty$, and 
\begin{equation}
A\left\Vert w\right\Vert _{\mathscr{H}}^{2}\leq\sum_{k\in\mathbb{N}}\left|\left\langle v_{k},w\right\rangle _{\mathscr{H}}\right|^{2}\leq B\left\Vert w\right\Vert _{\mathscr{H}}^{2}\label{eq:fr3}
\end{equation}
holds for all $w\in\mathscr{H}$.
\end{defn}
Note that (\ref{eq:fr3}) generalizes (\ref{eq:fr2}). Below we show
that, for frames, there is also a natural extension of (\ref{eq:fr1}).
\begin{prop}
\label{prop:frame1}Let $\left\{ v_{k}\right\} _{k\in\mathbb{N}}$
be a frame in $\mathscr{H}$; then there is a dual system $\left\{ v_{k}^{*}\right\} _{k\in\mathbb{N}}\subset\mathscr{H}$
such that the following representation holds:
\begin{equation}
w=\sum_{k\in\mathbb{N}}\left\langle v_{k}^{*},w\right\rangle _{\mathscr{H}}v_{k}\label{eq:fr4}
\end{equation}
for all $w\in\mathscr{H}$; absolute convergence.\end{prop}
\begin{proof}
Define the following operator $T:\mathscr{H}\rightarrow l^{2}\left(\mathbb{N}\right)$
by 
\[
Tw=\left(\left\langle v_{k},w\right\rangle _{\mathscr{H}}\right)_{k\in\mathbb{N}},
\]
and show that the adjoint $T^{*}:l^{2}\left(\mathbb{N}\right)\rightarrow\mathscr{H}$
satisfies 
\[
T^{*}\left(\left(x_{k}\right)\right)=\sum_{k\in\mathbb{N}}x_{k}v_{k}.
\]
Hence 
\begin{equation}
T^{*}Tw=\sum_{k\in\mathbb{N}}\left\langle v_{k},w\right\rangle _{\mathscr{H}}v_{k}.\label{eq:fr5}
\end{equation}
It follows that $T^{*}T$ has a bounded inverse, in fact, $A\:I_{\mathscr{H}}\leq T^{*}T\leq B\:I_{\mathscr{H}}$
in the order of selfadjoint operators. As a result, $\left(T^{*}T\right)^{-1}$
and $\left(T^{*}T\right)^{-\frac{1}{2}}$ are well-defined bounded
operators. 

Substitute $\left(T^{*}T\right)^{-1}$ into (\ref{eq:fr5}) yields:
\begin{eqnarray*}
w & = & \sum_{k\in\mathbb{N}}\left\langle v_{k},\left(T^{*}T\right)^{-1}w\right\rangle _{\mathscr{H}}v_{k}\\
 & = & \sum_{k\in\mathbb{N}}\left\langle \left(T^{*}T\right)^{-1}v_{k},w\right\rangle _{\mathscr{H}}v_{k}
\end{eqnarray*}
which is the desired (\ref{eq:fr4}) with $v_{k}^{*}:=\left(T^{*}T\right)^{-1}v_{k}$. \end{proof}
\begin{xca}[Frames from Lax-Milgram]
\myexercise{Frames from Lax-Milgram}\label{exer:ks3}Show that the
conclusion in Proposition \ref{prop:frame1} may also be obtained
from an application of the Lax-Milgram lemma (see \exerref{laxmilgram}.)
\sindex[nam]{Lax, P.D., (1926 --)}\index{Theorem!Lax-Milgram-}

Specifically, starting with a frame and given frame constants (see
(\ref{eq:fr3})), write down the corresponding sesquilinear form $B$
in Lax-Milgram, and verify that it satisfies the premise in Lax-Milgram.
Relate the frame bounds to the constants $b$, and $c$ in Lax-Milgram.
\end{xca}
\index{Lax-Milgram}
\begin{cor}
Let $\left\{ v_{k}\right\} $ be as in Proposition \ref{prop:frame1},
and set $v_{k}^{**}:=\left(T^{*}T\right)^{-\frac{1}{2}}v_{k}$, $k\in\mathbb{N}$;
then 
\[
w=\sum_{k\in\mathbb{N}}\left\langle v_{k}^{**},w\right\rangle _{\mathscr{H}}v_{k}^{**}
\]
holds for all $w\in\mathscr{H}$; absolute convergence.\end{cor}
\begin{rem}
We saw in \chapref{sp} (\secref{PVM}) that, if $\left\{ v_{k}\right\} _{k\in\mathbb{N}}$
is an ONB in some fixed Hilbert space $\mathscr{H}$, then 
\begin{equation}
P\left(\triangle\right)=\sum_{k\in\triangle}\left|v_{k}\left\rangle \right\langle v_{k}\right|,\;\triangle\in\mathcal{B}\left(\mathbb{R}\right),\label{eq:Pdel1}
\end{equation}
is a projection valued measure (PVM) on $\mathbb{R}$.

Suppose now that some $\left\{ v_{k}\right\} _{k\in\mathbb{N}}$ in
the expression (\ref{eq:Pdel1}) is only assumed to be a \emph{frame},
see \defref{frame}. \end{rem}
\begin{xca}[Positive operator valued measures from frames]
\myexercise{Positive operator valued measures from frames} Write
down the modified list of properties for $P\left(\cdot\right)$ in
(\ref{eq:Pdel1}) which generalize the axioms of \defref{PVM} for
PVMs. \end{xca}
\begin{rem}
It is possible to have \emph{uniqueness} for \emph{non-orthogonal}
expansions in Hilbert space. The following theorem of Rota et al.
is a case in point. \end{rem}
\begin{thm}[Rota et al. \cite{MR0125274,MR0056732}]
\label{thm:rota}Let $\mathscr{H}$ be a separable Hilbert space;
let $\left\{ e_{k}\right\} _{k\in\mathbb{N}}$ be an ONB in $\mathscr{H}$;
and let $\left\{ v_{k}\right\} _{k\in\mathbb{N}}$ be a linearly independent
system of vectors in $\mathscr{H}$ such that
\begin{equation}
\sum_{k=1}^{\infty}\left\Vert e_{k}-v_{k}\right\Vert ^{2}<\infty;\label{eq:ro1}
\end{equation}
then every vector $u\in\mathscr{H}$ has a unique representation 
\begin{equation}
u=\sum_{k=1}^{\infty}x_{k}v_{k},\quad x_{k}\in\mathbb{C}.\label{eq:ro2}
\end{equation}
Moreover defining
\begin{equation}
B\left(\sum_{k}x_{k}e_{k}\right):=\sum_{k}x_{k}v_{k},\quad(x_{k})\in l^{2};\label{eq:ro3}
\end{equation}
we get the following conclusions:
\begin{align*}
(i)\quad & B-I\;\mbox{is compact; and}\\
(ii)\quad & \mbox{ran}\left(B\right)=\mathscr{H}.
\end{align*}
\end{thm}
\begin{xca}[The operator $B$]
\myexercise{The operator $B$}Fill in the missing details in the
proof of \thmref{rota}.
\end{xca}

The primary source on KS is the paper by R.V. Kadison and I.M. Singer
\cite{MR0123922}. An important early paper is \cite{MR525951} by
Joel Anderson. 

Since the 1970ties, the KS problem has been studied with the use of
\textquotedblleft pavings;\textquotedblright{} see e.g., \cite{MR3173062,MR2859705,MR2830604,Wea03,MR1124564}.
While this {[}\textquotedblleft pavings\textquotedblright{} and their
equivalents{]} is an extremely interesting area, it is beyond the
scope of the present book.

\section*{A summary of relevant numbers from the Reference List}

For readers wishing to follow up sources, or to go in more depth with
topics above, we suggest: 

The pioneering paper \cite{MR0123922} started the subject, and in
the intervening decades there have been advances, and a discovery
of the relevance of the KS-problem to a host of applied areas, especially
harmonic analysis, frame theory, and signal processing. The problem
was solved two years ago. 

The most current paper concerning the solution to KS appears to be
\cite{MSS15} by Marcus, Spielman, and Strivastava. Paper \cite{1407.4768}
by P. Casazza explains the problem and its implications. A more comprehensive
citation list is: \cite{Arv76,BR81,Cas13,1407.4768,AW14,MSS15,MR3173062,MR525951,MR1124564,MR1381391,MR0123922,Wea03,MR2830604,Dixmier198101,BP44,MJD15}.

\part{Extension of Operators}

\chapter{Selfadjoint Extensions\label{chap:ext}}

\index{extension!-of operator}

\index{extension!selfadjoint-}
\begin{quotation}
Le plus court chemin entre deux vérités dans le domaine réel passe
par le domaine complexe.

--- Jacques Hadamard\sindex[nam]{Hadamard, J., (1865-1963)}\vspace{1em}

It will interest mathematical circles that the mathematical instruments
created by the higher algebra play an essential part in the rational
formulation of the new quantum mechanics. Thus the general proofs
of the conservation theorems in Heisenberg's theory carried out by
Born and Jordan are based on the use of the theory of matrices, which
go back to Cayley and were developed by Hermite. It is to be hoped
that a new era of mutual stimulation of mechanics and mathematics
has commenced. To the physicist it will seem first deplorable that
in atomic problems we have apparently met with such a limitation of
our usual means of visualisation. This regret will, however, have
to give way to thankfulness that mathematics, in this field too, presents
us with the tools to prepare the way for further progress.

--- Niels Bohr\sindex[nam]{Bohr, N., (1885-1962)}\vspace{1em}\\
\textquotedblleft Science is spectral analysis. Art is light synthesis.\textquotedblright{} 

--- Karl Kraus\sindex[nam]{Kraus, K., (1938-1988)}\vspace{2em}
\end{quotation}
Because of dictates from applications (especially quantum physics),
below we stress questions directly related to key-issues for unbounded
linear operators: Some operator from physics may only be \textquotedblleft formally
selfadjoint\textquotedblright{} also called Hermitian; and in such
cases, one ask for selfadjoint extensions (if any). 

The axioms of quantum physics (see e.g., \cite{BoMc13,OdHo13,KS02,CKS79,AAR13,Fan10,Maa10,Par09}
for relevant recent papers), are based on Hilbert space, and selfadjoint
operators. 

\index{operators!unbounded} 

\index{operators!selfadjoint}

\index{space!Hilbert-}

\index{measurement}

\index{selfadjoint operator}\index{axioms}

A quantum mechanical observable is a Hermitian (selfadjoint) linear
operator mapping a Hilbert space, the space of states, into itself.
The values obtained in a physical measurement are in general described
by a probability distribution; and the distribution represents a suitable
``average'' (or ``expectation'') in a measurement of values of
some quantum observable in a state of some prepared system. The states
are (up to phase) unit vectors in the Hilbert space, and a measurement
corresponds to a probability distribution (derived from a projection-valued
spectral measure). The particular probability distribution used depends
on both the state and the selfadjoint operator. The associated spectral
type may be continuous (such as position and momentum; both unbounded)
or discrete (such as spin); this depends on the physical quantity
being measured.

Since the spectral theorem serves as the central tool in quantum measurements
(see \cite{Sto90, Yos95, Ne69, RS75, DS88b}), we must be precise
about the distinction between linear operators with dense domain which
are only Hermitian (formally selfadjoint) as opposed to selfadjoint.
This distinction is accounted for by von  Neumann\textquoteright s
theory of deficiency indices \cite{AG93,DS88b,HdSS12}\footnote{Starting with \cite{vN32a,vN32b,vN32c}, J. von Neumann and M. Stone
did pioneering work in the 1930s on spectral theory for unbounded
operators in Hilbert space; much of it in private correspondence.
The first named author has from conversations with M. Stone, that
the notions \textquotedblleft deficiency-index,\textquotedblright{}
and \textquotedblleft deficiency space\textquotedblright{} are due
to them; suggested by MS to vN as means of translating more classical
notions of \textquotedblleft boundary values\textquotedblright{} into
rigorous tools in abstract Hilbert space: closed subspaces, projections,
and dimension count.}. 

\index{operators!formally selfadjoint}

\index{distribution!probability-}

\index{index!deficiency}

\index{operators!momentum-}

\index{index!von Neumann-}

\section{\label{sec:sa}Extensions of Hermitian Operators}

In order to apply spectral theorem, one must work with self adjoint
operators including the unbounded ones. Some examples first.\index{operators!extension of}\index{operators!Hermitian}

In quantum mechanics \cite{MR1939631,MR965583,MR675039}, to understand
energy levels of atoms and radiation, the energy level comes from
discrete packages. The interactions are given by Coulomb's Law\index{Coulomb's Law}
where
\[
H=-\triangle_{\vec{r}}+\frac{c_{jk}}{\left\Vert r_{j}-r_{k}\right\Vert }
\]
and Laplacian has dimension $3\times\#(\mbox{electrons})$.

\index{operators!Laplace-}

In Schrödinger's wave mechanics, one needs to solve for $\psi(r,t)$
from the equation 
\[
H\psi=\frac{1}{i}\frac{\partial}{\partial t}\psi.
\]
If we apply spectral theorem, then $\psi(t)=e^{itH}\psi(r,t=0)$.
This shows that motion in quantum mechanics is governed by unitary
operators. The two parts in Schrödinger equation are separately selfadjoint,
but justification of the sum being selfadjoint wasn't made rigorous
until 1957 when Kato wrote the book on ``perturbation theory'' \cite{Kat95}.
It is a summary of the sum of selfadjoint operators. \index{Schrödinger equation}

In Heisenberg's matrix mechanics, he suggested that one should look
at two states and the transition probability between them, such that
\[
\left\langle \psi_{1},A\psi_{2}\right\rangle =\left\langle \psi_{1}\left(t\right),A\psi_{2}(t)\right\rangle ,\;\forall t.
\]
If $\psi(t)=e^{itH}\psi$, then it works. In Heisenberg's picture,
one looks at evolution of the observables $e^{-itH}Ae^{itH}$. In
Schrödinger's picture, one looks at evolution of states. The two point
of views are equivalent. 

Everything so far is based on application of the spectral theorem,
which requires the operators being selfadjoint in the first place.

von  Neumann's index theory gives a complete classification of extensions
of single Hermitian unbounded operators with dense domain in a given
Hilbert space. The theory may be adapted to Hermitian representations
of $*$-algebras \cite{Nel59}. 

Let $A$ be a densely defined Hermitian operator on a Hilbert space
$\mathscr{H}$, i.e. $A\subset A^{*}$. If $B$ is any Hermitian extension
of $A$, then 
\begin{equation}
A\subset B\subset B^{*}\subset A^{*}.\label{eq:eso_struc_A_star}
\end{equation}
Since the adjoint operator $A^{*}$ is closed, i.e., $\mathscr{G}\left(A^{*}\right)$
is closed in $\mathscr{H}\oplus\mathscr{H}$, it follows that $\overline{\mathscr{G}\left(A\right)}\subset\mathscr{G}\left(A^{*}\right)$
is a well-defined operator graph, i.e., $A$ is closable and $\overline{\mathscr{G}\left(A\right)}=\mathscr{G}\left(\overline{A}\right)$.
Thus, there is no loss of generality to assume that $A$ is closed
and only consider its closed extensions. \index{closable operator}\index{operators!closable-}\index{operators!closed-}
\index{graph!- of operator}

The containment (\ref{eq:eso_struc_A_star}) suggests a detailed analysis
in $\mathscr{D}\left(A^{*}\right)\setminus\mathscr{D}\left(A\right)$.
Since $\mathscr{D}\left(A\right)$ is dense in $\mathscr{H}$, the
usual structural analysis in $\mathscr{H}$ (orthogonal decomposition,
etc.) is not applicable. However, this structure is brought out naturally
when $\mathscr{D}\left(A^{*}\right)$ is identified with the operator
graph $\mathscr{G}\left(A^{*}\right)$ in $\mathscr{H}\oplus\mathscr{H}$.
That is, $\mathscr{D}\left(A^{*}\right)$ is a Hilbert space under
its graph norm. With this identification, $\mathscr{D}(A)$ becomes
a closed subspace in $\mathscr{D}(A^{*})$, and \index{orthogonal!-decomposition}
\begin{equation}
\mathscr{D}\left(A^{*}\right)=\mathscr{D}\left(A\right)\oplus\left(\mathscr{D}\left(A^{*}\right)\ominus\mathscr{D}\left(A\right)\right).\label{eq:ext-2-1}
\end{equation}

The question of extending $A$ amounts to a further decomposition
\begin{equation}
\mathscr{D}\left(A^{*}\right)\ominus\mathscr{D}\left(A\right)=S\oplus K\label{eq:ext-2-2}
\end{equation}
in such a way that
\begin{eqnarray}
\widetilde{A} & = & A^{*}\big|_{\mathscr{D}\left(\widetilde{A}\right)},\mbox{ where }\label{eq:ext-2-3}\\
\mathscr{D}\big(\widetilde{A}\big) & = & \mathscr{D}\left(A\right)\oplus S\label{eq:ext-2-4}
\end{eqnarray}
defines a (closed) Hermitian operator $\widetilde{A}\supset A$. 

The extension $\widetilde{A}$ in (\ref{eq:ext-2-3})-(\ref{eq:ext-2-4})
is Hermitian iff the closed subspace $S\subset\mathscr{D}\left(A^{*}\right)$
is symmetric, in the sense that\index{operators!symmetric-} 
\begin{equation}
\left\langle A^{*}y,x\right\rangle -\left\langle y,A^{*}x\right\rangle =0,\;\forall x,y\in S.\label{eq:eso_symspace1}
\end{equation}

\begin{lem}
\label{lem:ext-1}Let $S$ be a closed subspace in $\mathscr{D}\left(A^{*}\right)$,
where $\mathscr{D}\left(A^{*}\right)$ is a Hilbert space under the
$A^{*}$-norm. The following are equivalent.
\begin{enumerate}
\item \label{enu:ext-1-1}$\left\langle A^{*}y,x\right\rangle =\left\langle y,A^{*}x\right\rangle $,
for all $x,y\in S$.
\item \label{enu:ext-1-2}$\left\langle x,A^{*}x\right\rangle \in\mathbb{R}$,
for all $x\in S$.
\end{enumerate}
\end{lem}
\begin{proof}
If (\ref{enu:ext-1-1}) holds, setting $x=y$, we get $\left\langle x,A^{*}x\right\rangle =\left\langle A^{*}x,x\right\rangle =\overline{\left\langle x,A^{*}x\right\rangle }$,
which implies that $\left\langle x,A^{*}x\right\rangle $ is real-valued. 

Conversely, assume (\ref{enu:ext-1-2}) is true. Since the mappings
\begin{eqnarray*}
\left(x,y\right) & \mapsto & \left\langle y,A^{*}x\right\rangle \\
\left(x,y\right) & \mapsto & \left\langle A^{*}y,x\right\rangle 
\end{eqnarray*}
are both sesquilinear forms on $S\times S$ (linear in the second
variable, and conjugate linear in the first variable), we apply the
polarization identity: 
\begin{eqnarray*}
\left\langle y,A^{*}x\right\rangle  & = & \frac{1}{4}\sum_{k=0}^{3}i^{k}\left\langle x+i^{k}y,A^{*}\left(x+i^{k}y\right)\right\rangle \\
\left\langle A^{*}y,x\right\rangle  & = & \frac{1}{4}\sum_{k=0}^{3}i^{k}\left\langle A^{*}\left(x+i^{k}y\right),x+i^{k}y\right\rangle 
\end{eqnarray*}
for all $x,y\in\mathscr{D}\left(A^{*}\right)$. Now, since $A$ is
Hermitian, the RHSs of the above equations are equal; therefore, $\left\langle y,A^{*}x\right\rangle =\left\langle A^{*}y,x\right\rangle $,
which is part (\ref{enu:ext-1-2}). \index{sesquilinear form}
\end{proof}
Eqs (\ref{eq:ext-2-3})-(\ref{eq:ext-2-4}) and \lemref{ext-1} set
up a bijection between (closed) Hermitian extensions of $A$ and (closed)
symmetric subspaces in $\mathscr{D}\left(A^{*}\right)\ominus\mathscr{D}\left(A\right)$.
Moreover, by \lemref{ext-1}, condition (\ref{eq:eso_symspace1})
is equivalent to \index{operators!symmetric-} 
\begin{equation}
\left\langle x,A^{*}x\right\rangle \in\mathbb{R},\;\forall x\in\mathscr{D}\left(A\right).\label{eq:eso_symspace2}
\end{equation}

Let $\varphi\in\mathscr{D}\left(A^{*}\right)$, such that $A^{*}\varphi=\lambda\varphi$,
$\Im\left\{ \lambda\right\} \neq0$; then $\left\langle \varphi,A^{*}\varphi\right\rangle =\lambda\left\Vert \varphi\right\Vert ^{2}\notin\mathbb{R}$.
By \lemref{ext-1} and (\ref{eq:eso_symspace2}), $\varphi\notin\mathscr{D}(\widetilde{A})$,
where $\widetilde{A}$ is any possible Hermitian extension of $A$.
This observation is in fact ruling out the ``wrong'' eigenvalues
of $\widetilde{A}$. Indeed, \thmref{esp_sa_criterion} below shows
that $A$ is selfadjoint if and only if ALL the ``wrong'' eigenvalues
of $A^{*}$ are excluded. But first we need the following lemma.\index{eigenvalue}
\begin{lem}
\label{lem:ext-2}Let $A$ be a Hermitian operator in $\mathscr{H}$,
then
\begin{equation}
\left\Vert \left(A-\lambda\right)x\right\Vert ^{2}=\left\Vert \left(A-a\right)x\right\Vert ^{2}+\left|b\right|^{2}\left\Vert x\right\Vert ^{2},\;\forall\lambda=a+ib\in\mathbb{C}.\label{eq:ext-2-5}
\end{equation}
In particular, 
\begin{equation}
\left\Vert \left(A-\lambda\right)x\right\Vert ^{2}\geq\left|\Im\left\{ \lambda\right\} \right|^{2}\left\Vert x\right\Vert ^{2},\;\forall\lambda\in\mathbb{C}.\label{eq:ext-2-6}
\end{equation}
\end{lem}
\begin{proof}
Write $\lambda=a+ib$, $a,b\in\mathbb{R}$; then 
\begin{eqnarray*}
 &  & \left\Vert \left(A-\lambda\right)x\right\Vert ^{2}\\
 & = & \left\langle \left(A-a\right)x-ibx,\left(A-a\right)x-ibx\right\rangle \\
 & = & \left\Vert \left(A-a\right)x\right\Vert ^{2}+\left|b\right|^{2}\left\Vert x\right\Vert ^{2}-i\left(\left\langle \left(A-a\right)x,x\right\rangle -\left\langle x,\left(A-a\right)x\right\rangle \right)\\
 & = & \left\Vert \left(A-a\right)x\right\Vert ^{2}+\left|b\right|^{2}\left\Vert x\right\Vert ^{2}\\
 & \geq & \left|b\right|^{2}\left\Vert x\right\Vert ^{2};
\end{eqnarray*}
where $\left\langle \left(A-a\right)x,x\right\rangle -\left\langle x,\left(A-a\right)x\right\rangle =0$,
since $A-a$ is Hermitian. \end{proof}
\begin{cor}
\label{cor:ext-1}Let $A$ be a closed Hermitian operator acting in
$\mathscr{H}$. Fix $\lambda\in\mathbb{C}$ with $\Im\left\{ \lambda\right\} \neq0$,
then $ran\left(A-\lambda\right)$ is a closed subspace in $\mathscr{H}$.
Consequently, we get the following decomposition\index{operators!closed-}
\begin{equation}
\mathscr{H}=ran\left(A-\lambda\right)\oplus ker\left(A^{*}-\overline{\lambda}\right).\label{eq:ext-1-2}
\end{equation}
\end{cor}
\begin{proof}
Set $B=A-\lambda$; then $B$ is closed, and so is $B^{-1}$, i.e.,
the operator graphs $\mathscr{G}\left(B\right)$ and $\mathscr{G}\left(B^{-1}\right)$
are closed in $\mathscr{H}\oplus\mathscr{H}$. Therefore, $ran\left(B\right)$
$\left(=dom(B^{-1})\right)$ is closed in $\left\Vert \cdot\right\Vert _{B^{-1}}$-norm.
But by (\ref{eq:ext-2-6}), $B^{-1}$ is bounded on $ran\left(B\right)$,
thus the two norms $\left\Vert \cdot\right\Vert $ and $\left\Vert \cdot\right\Vert _{B^{-1}}$
are equivalent on $ran\left(B\right)$. It follows that $ran\left(B\right)$
is also closed in $\left\Vert \cdot\right\Vert $-norm, i.e., it is
a closed subspace in $\mathscr{H}$. The decomposition (\ref{eq:ext-1-2})
follows from this.\end{proof}
\begin{thm}
\label{thm:esp_sa_criterion}Let $A$ be a densely defined, closed,
Hermitian operator in a Hilbert space $\mathscr{H}$; then the following
are equivalent: 
\[
\Big\{\exists\lambda,\:\Im\left\{ \lambda\right\} \neq0,\;ker\left(A^{*}-\lambda\right)=ker\left(A^{*}-\overline{\lambda}\right)=0\Big\}\Longleftrightarrow\Big\{ A=A^{*}\Big\}.
\]
\end{thm}
\begin{proof}
$\Longrightarrow$ By \corref{ext-1}, the hypothesis in the theorem
implies that 
\[
ran\left(A-\lambda\right)=ran\left(A-\overline{\lambda}\right)=\mathscr{H}.
\]
Let $y\in\mathscr{D}\left(A^{*}\right)$, then 
\begin{equation}
\left\langle y,\left(A-\lambda\right)x\right\rangle =\left\langle \left(A^{*}-\overline{\lambda}\right)y,x\right\rangle ,\;\forall x\in\mathscr{D}\left(A\right).\label{eq:ext-1-3}
\end{equation}
Since $ran\left(A-\overline{z}\right)=\mathscr{H}$, $\exists y_{0}\in\mathscr{D}\left(A\right)$
such that 
\[
\left(A^{*}-\overline{\lambda}\right)y=\left(A-\overline{\lambda}\right)y_{0}.
\]
Hence, RHS of (\ref{eq:ext-1-3}) is
\begin{equation}
\left\langle \left(A-\overline{\lambda}\right)y_{0},x\right\rangle =\left\langle y_{0},\left(A-\lambda\right)x\right\rangle .\label{eq:ext-1-4}
\end{equation}
Combining (\ref{eq:ext-1-3})-(\ref{eq:ext-1-4}), we then get 
\[
\left\langle y-y_{0},\left(A-\lambda\right)x\right\rangle =0,\;\forall x\in\mathscr{D}\left(A\right).
\]
Again, since $ran\left(A-\lambda\right)=\mathscr{H}$, the last equation
above shows that $y-y_{0}\perp\mathscr{H}$. In particular, $y-y_{0}\perp y-y_{0}$,
i.e., 
\[
\left\Vert y-y_{0}\right\Vert ^{2}=\left\langle y-y_{0},y-y_{0}\right\rangle =0.
\]
Therefore, $y=y_{0}$, and so $y\in\mathscr{D}\left(A\right)$. This
shows that $A^{*}\subset A$. 

The other containment $A\subset A^{*}$ holds since $A$ is assumed
to be Hermitian. Thus, we conclude that $A=A^{*}$.
\end{proof}
To capture all the ``wrong'' eigenvalues, we consider a family of
closed subspace in $\mathscr{H}$, $ker\left(A^{*}-\lambda\right)$,
where $\Im\left\{ \lambda\right\} \neq0$.
\begin{thm}
\label{thm:ext-1}If $A$ is a closed Hermitian operator in $\mathscr{H}$,
then 
\[
dim\left(ker\left(A^{*}-\lambda\right)\right)
\]
is a constant function on $\Im\left\{ \lambda\right\} >0$, and $\Im\left\{ \lambda\right\} <0$. \end{thm}
\begin{proof}
Fix $\lambda$ with $\Im\left\{ \lambda\right\} >0$. For $\Im\left\{ \lambda\right\} <0$,
the argument is similar. We proceed to verify that if $\eta\in\mathbb{C}$,
close enough to $\lambda$, then $dim\left(ker\left(A^{*}-\eta\right)\right)=dim\left(ker\left(A^{*}-\lambda\right)\right)$.
The desired result then follows immediately.

Since $A$ is closed, we have the following decomposition (by \corref{ext-1}),\index{operators!closed-}
\begin{equation}
\mathscr{H}=ran\left(A-\overline{\lambda}\right)\oplus ker\left(A^{*}-\lambda\right).\label{eq:ext-1-6}
\end{equation}
Now, pick $x\in ker\left(A^{*}-\eta\right)$, and suppose $x\perp ker\left(A^{*}-\lambda\right)$;
assuming $\left\Vert x\right\Vert =1$. By (\ref{eq:ext-1-6}), $\exists x_{0}\in\mathscr{D}\left(A\right)$
such that 
\begin{equation}
x=\left(A-\overline{\lambda}\right)x_{0}.\label{eq:ext-1-7}
\end{equation}
Then, 
\begin{eqnarray}
0=\left\langle \left(A^{*}-\eta\right)x,x_{0}\right\rangle  & = & \left\langle x,\left(A-\overline{\eta}\right)x_{0}\right\rangle \nonumber \\
 & = & \left\langle x,\left(A-\overline{\lambda}\right)x_{0}-\left(\overline{\eta}-\overline{\lambda}\right)x_{0}\right\rangle \nonumber \\
 & = & \left\Vert x\right\Vert ^{2}-\left(\overline{\eta}-\overline{\lambda}\right)\left\langle x,x_{0}\right\rangle \nonumber \\
 & \geq & \left\Vert x\right\Vert ^{2}-\left|\overline{\eta}-\overline{\lambda}\right|\Vert x\Vert^{2}\Vert x_{0}\Vert^{2}\;(\mbox{Cauchy-Schwarz})\nonumber \\
 & = & 1-\left|\eta-\lambda\right|\Vert x_{0}\Vert^{2}\label{eq:ext-1-8}
\end{eqnarray}
Applying \lemref{ext-2} to (\ref{eq:ext-1-7}), we also have 
\[
1=\left\Vert x\right\Vert ^{2}=\left\Vert \left(A-\overline{\lambda}\right)x_{0}\right\Vert ^{2}\geq\left|\Im\left\{ \lambda\right\} \right|^{2}\left\Vert x_{0}\right\Vert ^{2};
\]
substitute this into (\ref{eq:ext-1-8}), we see that 
\[
0\geq1-\left|\eta-\lambda\right|\Vert x_{0}\Vert^{2}\geq1-\left|\eta-\lambda\right|\left|\Im\left\{ \lambda\right\} \right|^{-2}
\]
which would be a contradiction if $\eta$ was close to $\lambda$.

It follows that the projection from $ker\left(A^{*}-\eta\right)$
to $ker\left(A^{*}-\lambda\right)$ is injective. For otherwise, $\exists x\in ker\left(A^{*}-\eta\right)$,
$x\neq0$, and $x\perp ker\left(A^{*}-\lambda\right)$. This is impossible
as shown above. Thus, 
\[
dim\left(ker\left(A^{*}-\eta\right)\right)\leq dim\left(ker\left(A^{*}-\lambda\right)\right).
\]
Similarly, we get the reversed inequality, and so 
\[
dim\left(ker\left(A^{*}-\eta\right)\right)=dim\left(ker\left(A^{*}-\lambda\right)\right).
\]
 
\end{proof}
A complete characterization of Hermitian extensions of a given Hermitian
operator is due to von  Neumann. \thmref{ext-1} suggests the following
definition:
\begin{defn}
Let $A$ be a densely defined, closed, Hermitian operator in $\mathscr{H}$.
The closed subspaces 
\begin{eqnarray}
\mathscr{D}_{\pm}\left(A\right) & = & ker\left(A^{*}\mp i\right)\label{eq:extensionSymOP_def_space}\\
 & = & \left\{ \xi\in\mathscr{D}\left(A^{*}\right):A^{*}\xi=\pm i\,\xi\right\} \nonumber 
\end{eqnarray}
are called the \emph{deficiency spaces} of $A$, and $dim\mathscr{D}_{\pm}\left(A\right)$
are called the \emph{deficiency indices}.
\end{defn}
For illustration, see \figref{def}.

\begin{figure}
$\mathscr{H}\;\begin{Bmatrix}\mathscr{D}_{+} & \xrightarrow{\quad\text{Partial Isometry}\quad} & \mathscr{D}_{-}\\
\oplus &  & \oplus\\
\left(A+i\right)\mathscr{D} & \xrightarrow[C_{A}=\left(A-i\right)\left(A+i\right)^{-1}]{} & \left(A-i\right)\mathscr{D}
\end{Bmatrix}\;\mbox{\ensuremath{\mathscr{H}}}$

\protect\caption{\label{fig:def}$\mathscr{D}_{\pm}=\mbox{Ker}\left(A^{*}\mp i\right)$,
$\mathscr{D}=dom\left(A\right)$}
\end{figure}

\index{deficiency space}

\index{index!von Neumann-}

The role of the \emph{Cayley-transform} $C_{A}:=\left(A-i\right)\left(A+i\right)^{-1}$,
and its extensions by partial isometries $\mathscr{D}_{+}\longrightarrow\mathscr{D}_{-}$,
is illustrated in \figref{def}. The Figure further offers a geometric
account of the conclusion in \thmref{extensionSymOP_VN_decomp_(A_adjoint)}. 

As a result we see that the two subspaces $\mathscr{D}_{\pm}$, also
called \emph{defect-spaces} (or \emph{deficiency-spaces}), are non-zero
precisely when the given symmetric operator $A$ \uline{fails}
to be essentially selfadjoint. The respective dimensions
\begin{equation}
n_{\pm}:=\dim\mathscr{D}_{\pm}\label{eq:vd1}
\end{equation}
are called \emph{deficiency indices}. The pair $\left(n_{+},n_{-}\right)$
in (\ref{eq:vd1}) is called the pair of \emph{von Neumann indices}.
We note that $A$ has selfadjoint extensions if and only if $n_{+}=n_{-}$. 

\index{operators!essentially selfadjoint-}

\index{selfadjoint extensions}
\begin{thm}[von Neumann]
\label{thm:extensionSymOP_VN_decomp_(A_adjoint)} Let $A$ be a densely
defined closed Hermitian operator acting in $\mathscr{H}$. Then 
\begin{equation}
\mathscr{D}\left(A^{*}\right)=\mathscr{D}\left(A\right)\oplus\mathscr{D}_{+}\left(A\right)\oplus\mathscr{D}_{-}\left(A\right);\label{eq:ext-1-9}
\end{equation}
where $\mathscr{D}(A^{*})$ is identified with its graph$\mathscr{G}(A^{*})$,
thus a Hilbert space under the graph inner product; and the decomposition
in (\ref{eq:ext-1-9}) refers to this Hilbert space.\index{space!Hilbert-}\end{thm}
\begin{proof}
By assumption, $A$ is closed, i.e., $\mathscr{D}\left(A\right)$,
identified with $\mathscr{G}\left(A\right)$, is a closed subspace
in $\mathscr{D}\left(A^{*}\right)$. 

Note that $\mathscr{D}_{\pm}\left(A\right)=ker\left(A^{*}\mp i\right)$
are closed subspaces in $\mathscr{H}$. Moreover, 
\[
\left\Vert x\right\Vert _{A^{*}}^{2}=\left\Vert x\right\Vert ^{2}+\left\Vert A^{*}x\right\Vert ^{2}=2\left\Vert x\right\Vert ^{2},\;\forall x\in\mathscr{D}_{\pm}\left(A\right);
\]
and so $\mathscr{D}_{\pm}\left(A\right)$, when identified with the
graph of $A^{*}\Big|_{\mathscr{D_{\pm}}\left(A^{*}\right)}$, are
also closed subspaces in $\mathscr{D}\left(A^{*}\right)$. 

Next, we verify the three subspaces on RHS of (\ref{eq:ext-1-9})
are mutually orthogonal. For all $x\in\mathscr{D}\left(A\right)$,
and all $x_{+}\in\mathscr{D}_{+}\left(A\right)=ker\left(A^{*}-i\right)$,
we have
\begin{eqnarray*}
\left\langle x_{+},x\right\rangle _{A^{*}} & = & \left\langle x_{+},x\right\rangle +\left\langle A^{*}x_{+},A^{*}x\right\rangle \\
 & = & \left\langle x_{+},x\right\rangle -i\left\langle x_{+},Ax\right\rangle \\
 & = & -i\left(\left\langle x_{+},i\,x\right\rangle +\left\langle x_{+},Ax\right\rangle \right)\\
 & = & -i\left\langle x_{+},\left(A+i\right)x\right\rangle =0
\end{eqnarray*}
where the last step follows from $x_{+}\perp ran\left(A+i\right)$
in $\mathscr{H}$, see (\ref{eq:ext-1-2}). Thus, $\mathscr{D}\left(A\right)\perp\mathscr{D}_{+}\left(A\right)$
in $\mathscr{D}\left(A^{*}\right)$. Similarly, $\mathscr{D}\left(A\right)\perp\mathscr{D}_{-}\left(A\right)$
in $\mathscr{D}\left(A^{*}\right)$. 

Moreover, if $x_{+}\in\mathscr{D}_{+}\left(A\right)$ and $x_{-}\in\mathscr{D}_{-}\left(A\right)$,
then 
\begin{eqnarray*}
\left\langle x_{+},x_{-}\right\rangle _{A^{*}} & = & \left\langle x_{+},x_{-}\right\rangle +\left\langle A^{*}x_{+},A^{*}x_{-}\right\rangle \\
 & = & \left\langle x_{+},x_{-}\right\rangle +\left\langle i\,x_{+},-i\,x_{-}\right\rangle \\
 & = & \left\langle x_{+},x_{-}\right\rangle -\left\langle x_{+},x_{-}\right\rangle =0.
\end{eqnarray*}
Hence $\mathscr{D}_{+}\left(A\right)\perp\mathscr{D}_{-}\left(A\right)$
in $\mathscr{D}\left(A^{*}\right)$.

Finally, we show RHS of (\ref{eq:ext-1-9}) yields the entire Hilbert
space $\mathscr{D}\left(A^{*}\right)$. For this, let $x\in\mathscr{D}\left(A^{*}\right)$,
and suppose (\ref{eq:ext-1-9}) holds, say, $x=x_{0}+x_{+}+x_{-}$,
where $x\in\mathscr{D}\left(A\right)$, $x_{\pm}\in\mathscr{D}_{\pm}\left(A\right)$;
then
\begin{eqnarray}
\left(A^{*}+i\right)x & = & \left(A^{*}+i\right)\left(x_{0}+x_{+}+x_{-}\right)\nonumber \\
 & = & \left(A+i\right)x_{0}+2i\,x_{+}.\label{eq:ext-1-10}
\end{eqnarray}
But, by the decomposition $\mathscr{H}=ran\left(A+i\right)\oplus ker\left(A^{*}-i\right)$,
eq. (\ref{eq:ext-1-2}), there exist $x_{0}$ and $x_{+}$ satisfying
(\ref{eq:ext-1-10}). It remains to set $x_{-}:=x-x_{0}-x_{+}$, and
to check $x_{-}\in\mathscr{D}_{-}\left(A\right)$. Indeed, by (\ref{eq:ext-1-10}),
we see that
\[
A^{*}x-Ax_{0}-i\,x_{+}=-i\,x+i\,x_{0}+i\,x_{+};\;\mbox{i.e.,}
\]
\[
A^{*}\left(x-x_{0}-x_{+}\right)=-i\left(x-x_{0}-x_{+}\right)
\]
and so $x_{-}\in\mathscr{D}_{-}\left(A\right)$. Therefore, we get
the desired orthogonal decomposition in (\ref{eq:ext-1-9}).\index{orthogonal!-decomposition}

Another argument: Let $y\in\mathscr{D}\left(A^{*}\right)$ such that
$y\perp\mathscr{D}_{\pm}\left(A\right)$ in $\mathscr{D}\left(A^{*}\right)$.
Then, $y\perp\mathscr{D}_{+}\left(A\right)$ in $\mathscr{D}\left(A^{*}\right)\Longrightarrow$
\begin{eqnarray*}
0 & = & \left\langle y,x_{+}\right\rangle +\left\langle A^{*}y,A^{*}x_{+}\right\rangle \\
 & = & \left\langle y,x_{+}\right\rangle +\left\langle A^{*}y,i\,x_{+}\right\rangle \\
 & = & i\left(\left\langle i\,y,x_{+}\right\rangle +\left\langle A^{*}y,x_{+}\right\rangle \right)\\
 & = & i\left\langle \left(A^{*}+i\right)y,x_{+}\right\rangle ,\;\forall x_{+}\in\mathscr{D}_{+}\left(A\right)=ker\left(A^{*}-i\right)
\end{eqnarray*}
and so $\exists\,x_{1}\in\mathscr{D}\left(A\right)$, and 
\begin{equation}
\left(A^{*}+i\right)y=\left(A+i\right)x_{1}.\label{eq:ext-1-11}
\end{equation}
On the other hand, $y\perp\mathscr{D}_{-}\left(A\right)$ in $\mathscr{D}\left(A^{*}\right)\Longrightarrow$
\begin{eqnarray*}
0 & = & \left\langle y,x_{-}\right\rangle +\left\langle A^{*}y,A^{*}x_{-}\right\rangle \\
 & = & \left\langle y,x_{-}\right\rangle +\left\langle A^{*}y,-i\,x_{-}\right\rangle \\
 & = & -i\left(\left\langle -i\,y,x_{-}\right\rangle +\left\langle A^{*}y,x_{-}\right\rangle \right)\\
 & = & i\left\langle \left(A^{*}-i\right)y,x_{-}\right\rangle ,\;\forall x_{-}\in\mathscr{D}_{-}\left(A\right)=ker\left(A^{*}+i\right);
\end{eqnarray*}
hence $\exists\,x_{2}\in\mathscr{D}\left(A\right)$, and 
\begin{equation}
\left(A^{*}-i\right)y=\left(A-i\right)x_{2}.\label{eq:ext-1-12}
\end{equation}
Subtracting (\ref{eq:ext-1-11})-(\ref{eq:ext-1-12}) then gives 
\[
y=\frac{x_{1}+x_{2}}{2}\in\mathscr{D}\left(A\right).
\]
\end{proof}
\begin{rem}
More generally, there is a family of decompositions 
\begin{equation}
\mathscr{D}\left(A^{*}\right)=\mathscr{D}\left(A\right)+ker\left(A^{*}-z\right)+ker\left(A^{*}-\overline{z}\right),\;\forall z\in\mathbb{C},\Im\left\{ z\right\} \neq0.\label{eq:ext-4-6}
\end{equation}
However, in the general case, we lose orthogonality. \end{rem}
\begin{proof}
Give $z\in\mathbb{C}$, $\Im\left\{ z\right\} \neq0$, suppose $x\in\mathscr{D}\left(A^{*}\right)$
can be written as 
\[
x=x_{0}+x_{+}+x_{-};
\]
where $x_{0}\in\mathscr{D}\left(A\right)$, $x_{+}\in ker\left(A^{*}-z\right)$,
and $x_{-}\in ker\left(A^{*}-\overline{z}\right)$. Then
\begin{eqnarray*}
A^{*}x & = & Ax_{0}+zx_{+}+\overline{z}x_{-}\\
\overline{z}x & = & \overline{z}x_{0}+\overline{z}x_{+}+\overline{z}x_{-}
\end{eqnarray*}
and
\begin{equation}
\left(A^{*}-\overline{z}\right)x=\left(A-\overline{z}\right)x_{0}+\left(z-\overline{z}\right)x_{+}.\label{eq:ext-4-4}
\end{equation}

Now, we start with (\ref{eq:ext-4-4}). By the decomposition 
\[
\mathscr{H}=ran\left(A-\overline{z}\right)\oplus ker\left(A^{*}-z\right),
\]
there exist unique $x_{0}$ and $x_{+}$ such that (\ref{eq:ext-4-4})
holds. This defines $x_{0}$ and $x_{+}$. Then, set 
\[
x_{-}:=x-x_{0}-x_{+};
\]
and it remains to check $x_{-}\in ker\left(A^{*}-\overline{z}\right)$.
Indeed, by (\ref{eq:ext-4-4}), we have
\[
A^{*}x-Ax_{0}-zx_{+}=\overline{z}x-\overline{z}x_{0}-\overline{z}x_{+},\;\mbox{i.e., }
\]
\[
A^{*}\left(x-x_{0}-x_{+}\right)=\overline{z}\left(x-x_{0}-x_{+}\right)
\]
thus, $x_{-}\in ker\left(A^{*}-\overline{z}\right)$.
\end{proof}

\begin{rem}
In the general decomposition (\ref{eq:ext-4-6}), if $f=x+x_{+}+x_{-}$,
$g=y+y_{+}+y_{-}$ where $f,g\in\mathscr{D}\left(A\right)$, $x_{+},y_{+}\in ker\left(A^{*}-z\right)$,
and $x_{-},y_{-}\in ker\left(A^{*}-\overline{z}\right)$; then 

\begin{eqnarray*}
 &  & \left\langle g,A^{*}f\right\rangle -\left\langle A^{*}g,f\right\rangle \\
 & = & \left\langle y+y_{+}+y_{-},A^{*}\left(x+x_{+}+x_{-}\right)\right\rangle -\left\langle A^{*}\left(y+y_{+}+y_{-}\right),x+x_{+}+x_{-}\right\rangle \\
 & = & \left\langle y+y_{+}+y_{-},Ax+zx_{+}+\overline{z}x_{-}\right\rangle -\left\langle Ay+zy_{+}+\overline{z}y_{-},x+x_{+}+x_{-}\right\rangle \\
 & = & \underset{0}{\underbrace{\left\langle y,Ax+zx_{+}+\overline{z}x_{-}\right\rangle -\left\langle Ay,x+x_{+}+x_{-}\right\rangle }}+\\
 &  & \underset{0}{\underbrace{\left\langle y_{+}+y_{-},Ax\right\rangle -\left\langle zy_{+}+\overline{z}y_{-},x\right\rangle }+}\\
 &  & \left\langle y_{+}+y_{-},zx_{+}+\overline{z}x_{-}\right\rangle -\left\langle zy_{+}+\overline{z}y_{-},x_{+}+x_{-}\right\rangle \\
 & = & \left\langle y_{+},zx_{+}\right\rangle -\left\langle zy_{+},x_{+}\right\rangle +\left\langle y_{-},\overline{z}x_{-}\right\rangle -\left\langle \overline{z}y_{-},x_{-}\right\rangle \\
 &  & +\left\langle y_{+},\overline{z}x_{-}\right\rangle +\left\langle y_{-},zx_{+}\right\rangle -\left\langle zy_{+},x_{-}\right\rangle -\left\langle \overline{z}y_{-},x_{+}\right\rangle \\
 & = & \left(z-\overline{z}\right)\left\langle y_{+},x_{+}\right\rangle +\left(\overline{z}-z\right)\left\langle y_{-},x_{-}\right\rangle +\\
 &  & \underset{0}{\underbrace{\overline{z}\left\langle y_{+},x_{-}\right\rangle +z\left\langle y_{-},x_{+}\right\rangle -\overline{z}\left\langle y_{+},x_{-}\right\rangle -z\left\langle y_{-},x_{+}\right\rangle }}\\
 & = & \left(z-\overline{z}\right)\left(\left\langle y_{+},x_{+}\right\rangle -\left\langle y_{-},x_{-}\right\rangle \right).
\end{eqnarray*}
\end{rem}
\begin{thm}[von Neumann]
\label{thm:eso_VN_extension} Let $A$ be a densely defined closed
Hermitian operator in $\mathscr{H}$. \index{operators!closed-}
\begin{enumerate}[leftmargin=*]
\item The (closed) Hermitian extensions of $A$ are indexed by partial
isometries with initial space in $\mathscr{D}_{+}\left(A\right)$
and final space in $\mathscr{D}_{-}\left(A\right)$. \index{isometry!partial-} 
\item Given a partial isometry $U$ as above, the Hermitian extension $\widetilde{A_{U}}\supset A$
is determined as follows: 
\begin{gather*}
\widetilde{A_{U}}\left(x+\left(1+U\right)x_{+}\right)=Ax+i\left(1-U\right)x_{+},\;\mbox{where}\\
\mathscr{D}\left(\widetilde{A_{U}}\right)=\left\{ x+x_{+}+Ux_{+}:x\in\mathscr{D}\left(A\right),x_{+}\in\mathscr{D}_{+}\left(A\right)\right\} 
\end{gather*}
\begin{equation}
\end{equation}

\end{enumerate}
\end{thm}
\begin{proof}
By the discussion in (\ref{eq:eso_symspace1}) and (\ref{eq:eso_symspace2}),
and \lemref{ext-1}, it remains to characterize the closed symmetric
subspaces $S$ in $\mathscr{D}_{+}\left(A\right)\oplus\mathscr{D}_{-}\left(A\right)$
$\left(\subset\mathscr{D}\left(A^{*}\right)\right)$. For this, let
$x=x_{+}+x_{-}$, $x_{\pm}\in\mathscr{D}_{\pm}\left(A\right)$, then
\index{operators!symmetric-}
\begin{eqnarray}
\left\langle x,A^{*}x\right\rangle  & = & \left\langle x_{+}+x_{-},A\left(x_{+}+x_{-}\right)\right\rangle \nonumber \\
 & = & \left\langle x_{+}+x_{-},i\left(x_{+}-x_{-}\right)\right\rangle \nonumber \\
 & = & i\left(\left\Vert x_{+}\right\Vert ^{2}-\left\Vert x_{-}\right\Vert ^{2}-2i\Im\left\{ \left\langle x_{+},x_{-}\right\rangle \right\} \right)\nonumber \\
 & = & i\left(\left\Vert x_{+}\right\Vert ^{2}-\left\Vert x_{-}\right\Vert ^{2}\right)+2\Im\left\{ \left\langle x_{+},x_{-}\right\rangle \right\} .\label{eq:eso_boundary_form}
\end{eqnarray}
Thus, 
\begin{align*}
\left\langle x,A^{*}x\right\rangle  & \in\mathbb{R},\;\forall x\in S\\
 & \Updownarrow\\
S= & \left\{ \left(x_{+},x_{-}\right):\left\Vert x_{+}\right\Vert =\left\Vert x_{-}\right\Vert ,\;x_{\pm}\in\mathscr{D}_{\pm}\left(A\right)\right\} 
\end{align*}
i.e., $S$ is identified with the graph of a partial isometry, say
$U$, with initial space in $\mathscr{D}_{+}\left(A\right)$ and final
space in $\mathscr{D}_{-}\left(A\right)$.\end{proof}
\begin{cor}
\label{cor:ext-2}Let $A$ be a densely defined, closed, Hermitian
operator on $\mathscr{H}$, and set $d_{\pm}=dim\left(\mathscr{D}_{\pm}\left(A\right)\right)$;
then 
\begin{enumerate}
\item $A$ is maximally Hermitian if and only if one of the deficiency indices
is $0$;
\item $A$ has a selfadjoint extension if and only if $d_{+}=d_{-}\neq0$;
\item $\overline{A}$ is selfadjoint if and only if $d_{+}=d_{-}=0$.
\end{enumerate}
\end{cor}
\begin{proof}
Immediate from \thmref{extensionSymOP_VN_decomp_(A_adjoint)} and
\thmref{eso_VN_extension}.\end{proof}
\begin{example}
$d_{+}=d_{-}=1$. Let $e_{\pm}$ be corresponding eigenvalues. $e_{+}\mapsto ze_{-}$
is the unitary operator sending one to the other eigenvalue. It is
clear that $\left|z\right|=1$. Hence the self adjoint extension is
indexed by $U_{1}(\mathbb{C})$.\index{eigenvalue}
\end{example}

\begin{example}
$d_{+}=d_{-}=2$, get a family of extensions indexed by $U_{2}(\mathbb{C})$.\end{example}
\begin{rem}
M. Stone and von Neumann are the two pioneers who worked at the same
period. They were born at about the same time. Stone died at 1970's
and von Neumann died in the 1950's.\index{Stone, M. H.}
\end{rem}
There is a simple criterion to test whether a Hermitian operator has
equal deficiency indices.
\begin{defn}
An operator $J:\mathscr{H}\rightarrow\mathscr{H}$ is called a \emph{conjugation}
if 
\begin{itemize}
\item $J$ is conjugate linear, i.e., $J\left(cx\right)=\overline{c}x$,
for all $x\in\mathscr{H}$, and all $c\in\mathbb{C}$, 
\item $J^{2}=1$, and
\item $\left\langle Jx,Jy\right\rangle =\left\langle y,x\right\rangle $,
for all $x,y\in\mathscr{H}$.
\end{itemize}
\end{defn}
\begin{thm}[von Neumann]
\label{thm:eso_VN_conjugation} Let $A$ be a densely defined closed
Hermitian operator in $\mathscr{H}$. Set $d_{\pm}=dim\left(\mathscr{D}_{\pm}\left(A\right)\right)$.
Suppose $AJ=JA$, where $J$ is a conjugation, then $d_{+}=d_{-}$.
In particular,  $A$ has selfadjoint extensions. \index{selfadjoint extensions}\end{thm}
\begin{proof}
Note that, by definition, we have $\left\langle Jx,y\right\rangle =\left\langle Jx,J^{2}y\right\rangle =\left\langle Jy,x\right\rangle $,
for all $x,y\in\mathscr{H}$. 

We proceed to show that $J$ commutes with $A^{*}$. For this, let
$x\in\mathscr{D}\left(A\right)$, $y\in\mathscr{D}\left(A^{*}\right)$,
then
\begin{equation}
\left\langle JA^{*}y,x\right\rangle =\left\langle Jx,A^{*}y\right\rangle =\left\langle AJx,y\right\rangle =\left\langle JAx,y\right\rangle =\left\langle Jy,Ax\right\rangle .\label{eq:ext-4-1}
\end{equation}
It follows that $x\mapsto\left\langle Jy,Ax\right\rangle $ is bounded,
and $Jy\in\mathscr{D}\left(A^{*}\right)$. Thus, $J\mathscr{D}\left(A^{*}\right)\subset\mathscr{D}\left(A^{*}\right)$.
Since $J^{2}=1$, $\mathscr{D}\left(A^{*}\right)=J^{2}\mathscr{D}\left(A^{*}\right)\subset J\mathscr{D}\left(A^{*}\right)$;
therefore, $J\mathscr{D}\left(A^{*}\right)=\mathscr{D}\left(A^{*}\right)$.
Moreover, (\ref{eq:ext-4-1}) shows that $JA^{*}=A^{*}J$. 

Now if $x\in\mathscr{D}_{+}\left(A\right)$, then 
\[
A^{*}Jx=JA^{*}x=J\left(ix\right)=-iJx
\]
i.e., $J\mathscr{D}_{+}\left(A\right)\subset\mathscr{D}_{-}\left(A\right)$.
Similarly, $J\mathscr{D}_{-}\left(A\right)\subset\mathscr{D}_{+}\left(A\right)$. 

Using $J^{2}=1$ again, $\mathscr{D}_{-}\left(A\right)=J^{2}\mathscr{D}_{-}\left(A\right)\subset J\mathscr{D}_{+}\left(A\right)$;
and we conclude that $J\mathscr{D}_{+}\left(A\right)=\mathscr{D}_{-}\left(A\right)$. 

Since the restriction of $J$ to $\mathscr{D}_{+}\left(A\right)$
preserves orthonormal basis, we then get $dim\left(\mathscr{D}_{+}\left(A\right)\right)=dim\left(\mathscr{D}_{-}\left(A\right)\right)$. 
\end{proof}

\section{Cayley Transform\label{sec:Cayley}}

There is an equivalent characterization of Hermitian extensions, taking
place entirely in $\mathscr{H}$ and without the identification of
$\mathscr{D}\left(A^{*}\right)\simeq\mathscr{G}\left(A^{*}\right)$,
where $\mathscr{G}\left(A^{*}\right)$ is seen as a Hilbert space
under its graph inner product. This is the result of the following
observation.\index{transform!Cayley}
\begin{lem}
\label{lem:ext-3}Let $A$ be a Hermitian operator acting in $\mathscr{H}$;
then 
\begin{equation}
\left\Vert \left(A\pm i\right)x\right\Vert ^{2}=\left\Vert x\right\Vert ^{2}+\left\Vert Ax\right\Vert ^{2},\;\forall x\in\mathscr{D}\left(A\right).\label{eq:ext-1-13}
\end{equation}
\end{lem}
\begin{proof}
See \lemref{ext-2}. Or, a direct computation shows that 
\begin{align*}
\left\Vert \left(A+i\right)x\right\Vert ^{2} & =\left\langle \left(A+i\right)x,\left(A+i\right)x\right\rangle \\
 & =\left\Vert x\right\Vert ^{2}+\left\Vert Ax\right\Vert ^{2}+i\left(\left\langle Ax,x\right\rangle -\left\langle x,Ax\right\rangle \right)\\
 & =\left\Vert x\right\Vert ^{2}+\left\Vert Ax\right\Vert ^{2};
\end{align*}
where $\left\langle Ax,x\right\rangle -\left\langle x,Ax\right\rangle =0$
since $A$ is Hermitian. \end{proof}
\begin{thm}[Cayley transform]
\label{thm:cayley}Let $A$ be a densely defined, closed, Hermitian
operator in $\mathscr{H}$. 
\begin{enumerate}
\item The following subspaces in $\mathscr{H}$ are isometrically isomorphic:
\[
ran\left(A\pm i\right)\simeq\mathscr{G}\left(A\right)\simeq\mathscr{D}\left(A\right).
\]
In particular, $ran\left(A\pm i\right)$ are closed subspace in $\mathscr{H}$.
\item The map $C_{A}:ran\left(A+i\right)\rightarrow ran\left(A-i\right)$
by 
\begin{equation}
\left(A+i\right)x\mapsto\left(A-i\right)x,\;\forall x\in\mathscr{D}\left(A\right)\label{eq:ext-4-2}
\end{equation}
is isometric. Equivalently, 
\begin{equation}
C_{A}x=\left(A-i\right)\left(A+i\right)^{-1}x\label{eq:ext-14}
\end{equation}
for all $x\in ran\left(A+i\right)$. 
\item Moreover, 
\begin{equation}
A=i\left(1+C_{A}\right)\left(1-C_{A}\right)^{-1}.\label{eq:ext-15}
\end{equation}

\end{enumerate}
\end{thm}
\begin{proof}
By (\ref{lem:ext-3}), $ran\left(A\pm i\right)$ are isometric to
the graph of $A$, and the latter is closed (as a subset in $\mathscr{H}\oplus\mathscr{H}$)
since $A$ is closed (i.e., $\mathscr{G}\left(A\right)$ is closed).
Thus, $ran\left(A\pm i\right)$ are closed in $\mathscr{H}$. Note
this is also a result of \corref{ext-1}. 

The mapping (\ref{eq:ext-14}) being isometric follows from (\ref{eq:ext-1-13}).

By (\ref{eq:ext-4-2}), we have 
\begin{align*}
\left(1-C_{A}\right)\left(\left(A+i\right)x\right) & =\left(A+i\right)x-\left(A-i\right)x=2ix\\
\left(1+C_{A}\right)\left(\left(A+i\right)x\right) & =\left(A+i\right)x+\left(A-i\right)x=2Ax
\end{align*}
for all $x\in\mathscr{D}\left(A\right)$. It follows that 
\[
\left(1+C_{A}\right)\left(1-C_{A}\right)^{-1}\left(2ix\right)=\left(1+C_{A}\right)\left(\left(A+i\right)x\right)=2Ax;\mbox{ i.e.,}
\]
\[
Ax=i\left(1+C_{A}\right)\left(1-C_{A}\right)^{-1}x,\;\forall x\in\mathscr{D}\left(A\right)
\]
which is (\ref{eq:ext-15}).\end{proof}
\begin{thm}
Suppose $A$ is densely defined, closed, and Hermitian in $\mathscr{H}$.
Then the family of (closed) Hermitian extensions of $A$ is indexed
by partial isometries $U$ with initial space in $\mathscr{D}_{+}\left(A\right)$
and final space in $\mathscr{D}_{-}\left(A\right)$. Given $U$, the
corresponding extension $\widetilde{A}_{U}\supset A$ is determined
by
\begin{gather*}
\widetilde{A}_{U}(x+\left(1-U\right)x_{+})=x+i\left(1+U\right)x_{+},\;\mbox{where}\\
dom(\widetilde{A}_{U})=\left\{ x+\left(1-U\right)x_{+}:x\in\mathscr{D}\left(A\right),x_{+}\in\mathscr{D}_{+}\left(A\right)\right\} 
\end{gather*}
Moreover, $\widetilde{A}_{U}$ is selfadjoint if and only if $U$
is unitary from $\mathscr{D}_{+}\left(A\right)$ onto $\mathscr{D}_{-}\left(A\right)$.\index{extension!selfadjoint-}\end{thm}
\begin{proof}
Since $A$ is closed, we get the following decompositions (\corref{ext-1})
\begin{align*}
\mathscr{H} & =ran\left(A+i\right)\oplus ker\left(A^{*}-i\right)\\
 & =ran\left(A-i\right)\oplus ker\left(A^{*}+i\right).
\end{align*}
By \thmref{cayley}, $C_{A}:ran\left(A+i\right)\rightarrow ran\left(A-i\right)$
is isometric. Consequently, getting a Hermitian extension of $A$
amounts to choosing a partial isometry $U$ with initial space in
$ker\left(A^{*}-i\right)\left(=\mathscr{D}_{+}\left(A\right)\right)$
and final space in $ker\left(A^{*}+i\right)\left(=\mathscr{D}_{-}\left(A\right)\right)$,
such that 
\[
C_{\widetilde{A}_{U}}:=C_{A}\oplus U
\]
is the Cayley transform of $\widetilde{A}_{U}\supset A$. 

Given $U$ as above, for all $x\in\mathscr{D}\left(A\right)$, $x_{+}\in\mathscr{D}_{+}\left(A\right)$,
we have 
\[
C_{\widetilde{A}_{U}}\left(\left(A+i\right)x\oplus x_{+}\right)=\left(A-i\right)x\oplus Ux_{+}.
\]
 Then, 
\begin{align*}
(1-C_{\widetilde{A}_{U}})\left(\left(A+i\right)x\oplus x_{+}\right) & =\left(\left(A+i\right)x+x_{+}\right)-\left(\left(A-i\right)x+Ux_{+}\right)\\
 & =2ix+\left(1-U\right)x_{+}\\
(1+C_{\widetilde{A}_{U}})\left(\left(A+i\right)x\oplus x_{+}\right) & =\left(\left(A+i\right)x+x_{+}\right)+\left(\left(A-i\right)x+Ux_{+}\right)\\
 & =2Ax+\left(1+U\right)x_{+};
\end{align*}
and so 
\[
i(1+C_{\widetilde{A}_{U}})(1-C_{\widetilde{A}_{U}})^{-1}\left(x+\frac{1}{2i}\left(1-U\right)x_{+}\right)=Ax+\frac{1}{2}\left(1+U\right)x_{+}.
\]
The theorem follows by setting $x_{+}:=2iy_{+}$. \index{isometry!Cayley-}
\end{proof}

\section{Boundary Triple}

\index{extension!selfadjoint-}

In applications, especially differential equations, it is convenient
to characterize selfadjoint extensions using boundary conditions.
For recent applications, see \cite{JPT11-1,JPT12,JPT12-1}. A slightly
modified version can be found in \cite{dO09}. \index{boundary condition}
\index{boundary triple}

Let $A$ be a densely defined, closed, Hermitian operator acting in
a Hilbert space $\mathscr{H}$. Assume $A$ has deficiency indices
$\left(d,d\right)$, $d>0$, and so $A$ has non-trivial selfadjoint
extensions. By von Neumann's theorem (\thmref{extensionSymOP_VN_decomp_(A_adjoint)}),
for all $x,y\in\mathscr{D}\left(A^{*}\right)$, we have the following
decomposition, 
\begin{eqnarray*}
x & = & x_{0}+x_{+}+x_{-}\\
y & = & y_{0}+y_{+}+y_{-}
\end{eqnarray*}
where $x_{0},y_{0}\in\mathscr{D}\left(A\right)$, $x_{+},y_{+}\in\mathscr{D}_{+}\left(A\right)$,
and $x_{-},y_{-}\in\mathscr{D}_{-}\left(A\right)$. Then, 
\begin{eqnarray}
 &  & \left\langle y,A^{*}x\right\rangle -\left\langle A^{*}y,x\right\rangle \nonumber \\
 & = & \left\langle y_{0}+y_{+}+y_{-},Ax_{0}+i\left(x_{+}-x_{-}\right)\right\rangle -\nonumber \\
 &  & \left\langle Ay_{0}+i\left(y_{+}-y_{-}\right),x_{0}+x_{+}+x_{-}\right\rangle \nonumber \\
 & = & \underset{0}{\underbrace{\left\langle y_{0},Ax_{0}\right\rangle -\left\langle Ay_{0},x_{0}\right\rangle }}+\underset{0}{\underbrace{\left\langle y_{0},i\left(x_{+}-x_{-}\right)\right\rangle -\left\langle Ay_{0},x_{+}+x_{-}\right\rangle }}+\nonumber \\
 &  & \underset{0}{\underbrace{\left\langle y_{+}+y_{-},Ax_{0}\right\rangle -\left\langle i\left(y_{+}-y_{-}\right),x_{0}\right\rangle }}+\nonumber \\
 &  & \left\langle y_{+}+y_{-},i\left(x_{+}-x_{-}\right)\right\rangle -\left\langle i\left(y_{+}-y_{-}\right),x_{+}+x_{-}\right\rangle \nonumber \\
 & = & 2i\left\{ \left\langle y_{+},x_{+}\right\rangle -\left\langle y_{-},x_{-}\right\rangle \right\} .\label{eq:ext-3-1}
\end{eqnarray}
Therefore, we see that 
\[
\left[x,y\in\mathscr{D}\big(\widetilde{A}\big),\;\widetilde{A}\supset A,\;\mbox{Hermitian extension}\right]\Longleftrightarrow\left[\mbox{RHS of \ensuremath{\left(\ref{eq:ext-3-1}\right)} vanishes}\right]
\]
For selfadjoint extensions, this is equivalent to choosing a partial
isometry $U$ from $\mathscr{D}_{+}\left(A\right)$ onto $\mathscr{D}_{-}\left(A\right)$,
and setting 
\[
x_{-}=Ux_{+},\;y_{-}=Uy_{+};\mbox{ so that}
\]
\begin{eqnarray*}
\left\langle y,A^{*}x\right\rangle -\left\langle A^{*}y,x\right\rangle  & = & 2i\left\{ \left\langle y_{+},x_{+}\right\rangle -\left\langle Uy_{+},Ux_{+}\right\rangle \right\} \\
 & = & 2i\left\{ \left\langle y_{+},x_{+}\right\rangle -\left\langle y_{+},x_{+}\right\rangle \right\} =0.
\end{eqnarray*}
The discussion above leads to the following definition:
\begin{defn}
Let $A$ be a densely defined, closed, Hermitian operator in $\mathscr{H}$.
Suppose $A$ has deficiency indices $\left(d,d\right)$, $d>0$. A
boundary space for $A$ is a triple $\left(\mathscr{H}_{b},\rho_{1},\rho_{2}\right)$
consisting of a Hilbert space $\mathscr{H}_{b}$ and two linear maps
$\rho_{1},\rho_{2}:\mathscr{D}\left(A^{*}\right)\rightarrow\mathscr{H}_{b}$,
such that
\begin{enumerate}
\item $\rho_{i}\left(\mathscr{D}\left(A^{*}\right)\right)$ is dense in
$\mathscr{H}_{b}$, $i=1,2$; and
\item for all $x,y\in\mathscr{D}\left(A^{*}\right)$, $\exists\,c\neq0$,
such that 
\begin{equation}
\left\langle y,A^{*}x\right\rangle -\left\langle A^{*}y,x\right\rangle =c\left[\left\langle \rho_{1}\left(y\right),\rho_{1}\left(x\right)\right\rangle _{b}-\left\langle \rho_{2}\left(y\right),\rho_{2}\left(x\right)\right\rangle _{b}\right].\label{eq:eso_boundary_form1}
\end{equation}

\end{enumerate}
\end{defn}
\begin{rem}
In (\ref{eq:ext-3-1}), we set 
\begin{flalign*}
\mathscr{H}_{b} & =\mathscr{D}_{+}\left(A\right)\\
\rho_{1}\left(x_{0}+x_{+}+x_{-}\right) & =x_{+}\\
\rho_{2}\left(x_{0}+x_{+}+x_{-}\right) & =Ux_{+}
\end{flalign*}
for any $x=x_{0}+x_{+}+x_{-}$ in $\mathscr{D}\left(A^{*}\right)$.
Then $\left(\mathscr{H}_{b},\rho_{1},\rho_{2}\right)$ is a boundary
space for $A$. In this special case, $\rho_{1},\rho_{2}$ are surjective.
It is clear that the choice of a boundary triple is not unique. In
applications, $\mathscr{H}_{b}$ is usually chosen to have the same
dimension as $\mathscr{D}_{\pm}\left(A\right)$.\index{boundary triple}
\end{rem}
Consequently, \thmref{eso_VN_extension} can be restated as follows.
\begin{thm}
\label{thm:ext-3-1}Let $A$ be a densely defined, closed, Hermitian
operator in $\mathscr{H}$. Suppose $A$ has deficiency indices $\left(d,d\right)$,
$d>0$. Let $\left(\mathscr{H}_{b},\rho_{1},\rho_{2}\right)$ be a
boundary triple. Then the selfadjoint extensions of $A$ are indexed
by unitary operators $U:\mathscr{H}_{b}\rightarrow\mathscr{H}_{b}$,
such that given $U$, the corresponding selfadjoint extension $\widetilde{A_{U}}\supset A$
is determined by 
\begin{gather*}
\widetilde{A_{U}}=A^{*}\Big|_{\mathscr{D}\left(\widetilde{A_{U}}\right)},\;\mbox{where}\\
\mathscr{D}\left(\widetilde{A_{U}}\right)=\left\{ x\in\mathscr{D}\left(A^{*}\right):U\rho_{1}\left(x\right)=\rho_{2}\left(x\right)\right\} .
\end{gather*}

\end{thm}
Certain variations of \thmref{ext-3-1} are convenient in the boundary
value problems (BVP) of differential equations. In \cite{DM91,GG91},
a boundary triple $\left(\mathscr{H}_{b},\beta_{1},\beta_{2}\right)$
is defined to satisfy
\begin{equation}
\left\langle x,A^{*}y\right\rangle =\left\langle A^{*}x,y\right\rangle =c'\left[\left\langle \beta_{1}\left(x\right),\beta_{2}\left(y\right)\right\rangle _{b}-\left\langle \beta_{2}\left(x\right),\beta_{1}\left(y\right)\right\rangle _{b}\right]\label{eq:eso_boundary_form2}
\end{equation}
for all $x,y\in\mathscr{D}\left(A^{*}\right)$; and $c'$ is some
nonzero constant. Also, see \cite{JPT11-1,JPT12,JPT12-1}.\index{selfadjoint extensions}

The connection between (\ref{eq:eso_boundary_form1}) and (\ref{eq:eso_boundary_form2})
is via the bijection
\begin{equation}
\begin{Bmatrix}\rho_{1} & = & \beta_{1}+i\beta_{2}\\
\rho_{2} & = & \beta_{1}-i\beta_{2}
\end{Bmatrix}\Longleftrightarrow\begin{Bmatrix}\beta_{1} & = & \dfrac{\rho_{1}+\rho_{2}}{2}\\
\beta_{2} & = & \dfrac{\rho_{1}-\rho_{2}}{2i}
\end{Bmatrix}.\label{eq:ext-3-2}
\end{equation}

\begin{lem}
Under the bijection (\ref{eq:ext-3-2}), we have 
\[
\left\langle \rho_{1}\left(x\right),\rho_{1}\left(y\right)\right\rangle _{b}-\left\langle \rho_{2}\left(x\right),\rho_{2}\left(y\right)\right\rangle _{b}=2i\left(\left\langle \beta_{1}\left(x\right),\beta_{2}\left(y\right)\right\rangle _{b}-\left\langle \beta_{2}\left(x\right),\beta_{1}\left(y\right)\right\rangle _{b}\right)
\]
\end{lem}
\begin{proof}
For convenience, we suppress the variables $x,y$. Then a direct computation
shows that,
\begin{eqnarray*}
 &  & \left\langle \rho_{1},\rho_{1}\right\rangle _{b}-\left\langle \rho_{2},\rho_{2}\right\rangle _{b}\\
 & = & \left\langle \beta_{1}+i\beta_{2},\beta_{1}+i\beta_{2}\right\rangle _{b}-\left\langle \beta_{1}-i\beta_{2},\beta_{1}-i\beta_{2}\right\rangle _{b}\\
 & = & i\left\langle \beta_{1},\beta_{2}\right\rangle _{b}-i\left\langle \beta_{2},\beta_{1}\right\rangle _{b}+i\left\langle \beta_{1},\beta_{2}\right\rangle _{b}-i\left\langle \beta_{2},\beta_{1}\right\rangle _{b}\\
 & = & 2i\left(\left\langle \beta_{1},\beta_{2}\right\rangle _{b}-\left\langle \beta_{2},\beta_{1}\right\rangle _{b}\right)
\end{eqnarray*}
which is the desired conclusion. \end{proof}
\begin{thm}
\label{thm:bt}Given a boundary triple $\left(\mathscr{H}_{b},\beta_{1},\beta_{2}\right)$
satisfying (\ref{eq:eso_boundary_form2}), the family of selfadjoint
extensions $\widetilde{A_{U}}\supset A$ is indexed by unitary operators
$U:\mathscr{H}_{b}\rightarrow\mathscr{H}_{b}$, such that \index{boundary triple}
\begin{gather}
\widetilde{A_{U}}=A^{*}\Big|_{\mathscr{D}\left(\widetilde{A_{U}}\right)},\;\mbox{where}\\
\mathscr{D}\left(\widetilde{A_{U}}\right)=\left\{ x\in\mathscr{D}\left(A^{*}\right):\left(1-U\right)\beta_{1}\left(x\right)=i\left(1+U\right)\beta_{2}\left(x\right)\right\} .
\end{gather}
\end{thm}
\begin{proof}
By \thmref{ext-3-1}, we need only pick a unitary operator $U:\mathscr{H}_{b}\rightarrow\mathscr{H}_{b}$,
such that $\rho_{2}=U\rho_{1}$. In view of the bijection (\ref{eq:ext-3-2}),
this yields
\[
\beta_{1}-i\beta_{2}=U\left(\beta_{1}+i\beta_{2}\right)\Longleftrightarrow\left(1-U\right)\beta_{1}=i\left(1+U\right)\beta_{2}
\]
and the theorem follows.\end{proof}
\begin{example}
Let $A=-i\frac{d}{dx}\Big|_{\mathscr{D}\left(A\right)}$, and 
\[
\mathscr{D}\left(A\right)=\left\{ f\in L^{2}\left(0,1\right):f'\in L^{2}\left(0,1\right),f\left(0\right)=f\left(1\right)=0\right\} .
\]
Then $A^{*}=-i\frac{d}{dx}\Big|_{\mathscr{D}\left(A^{*}\right)}$,
where 
\[
\mathscr{D}\left(A^{*}\right)=\left\{ f:f,f'\in L^{2}\left(0,1\right)\right\} .
\]
For all $f,g\in\mathscr{D}\left(A^{*}\right)$, using integration
by parts, we get 
\[
\left\langle g,A^{*}f\right\rangle -\left\langle A^{*}g,f\right\rangle =-i\overline{g\left(x\right)}f\left(x\right)\Big|_{0}^{1}=-i\left(\overline{g\left(1\right)}f\left(1\right)-\overline{g\left(0\right)}f\left(0\right)\right).
\]
Let $\mathscr{H}_{b}=\mathbb{C}$, i.e., one-dimensional, and set
\[
\rho_{1}\left(f\right)=f\left(1\right),\;\rho_{2}\left(f\right)=f\left(0\right);\;\mbox{then}
\]
\[
\left\langle g,A^{*}f\right\rangle -\left\langle A^{*}g,f\right\rangle =-i\left(\left\langle \rho_{1}\left(g\right),\rho_{1}\left(f\right)\right\rangle _{b}-\left\langle \rho_{2}\left(g\right),\rho_{2}\left(f\right)\right\rangle _{b}\right).
\]
Therefore, $\left(\mathscr{H}_{b},\rho_{1},\rho_{2}\right)$ is a
boundary triple. 

The family of selfadjoint extensions of $A$ is given by the unitary
operator 
\[
e^{i\theta}:\mathscr{H}_{b}\rightarrow\mathscr{H}_{b},\;s.t.\;\rho_{2}=e^{i\theta}\rho_{1};
\]
i.e.,
\[
\widetilde{A_{\theta}}=-i\frac{d}{dx}\Big|_{\left\{ f\in\mathscr{D}\left(A^{*}\right):f\left(0\right)=e^{i\theta}f\left(1\right)\right\} }.
\]

\end{example}

\begin{example}
\label{exa:bt1}Let $Af=-f''$, with $\mathscr{D}\left(A\right)=C_{c}^{\infty}\left(0,\infty\right)$.
Since $A$ is Hermitian and $A\geq0$, it follows that it has equal
deficiency indices. Also, $\mathscr{D}\left(A^{*}\right)=\left\{ f,f''\in L^{2}\left(0,\infty\right)\right\} $,
and $A^{*}f=-f''$, $\forall f\in\mathscr{D}\left(A^{*}\right)$. 

For $f,g\in\mathscr{D}^{*}\left(A\right)$, we have 
\begin{align*}
\left\langle g,A^{*}f\right\rangle =-\int_{0}^{\infty}\overline{g}f'' & =-\left(\left[\overline{g}f'-\overline{g}'f\right]_{0}^{\infty}+\int_{0}^{\infty}\overline{g''}f\right)\\
 & =\left(\overline{g}f'\right)\left(0\right)-\left(\overline{g}'f\right)\left(0\right)-\int_{0}^{\infty}\overline{g''}f\\
 & =\left(\overline{g}f'\right)\left(0\right)-\left(\overline{g}'f\right)\left(0\right)+\left\langle A^{*}g,f\right\rangle 
\end{align*}
and so 
\[
\left\langle g,A^{*}f\right\rangle -\left\langle A^{*}g,f\right\rangle =\left(\overline{g}f'\right)\left(0\right)-\left(\overline{g}'f\right)\left(0\right).
\]

Now, set $\mathscr{H}_{b}=\mathbb{C}$, i.e., one-dimensional, and
\[
\beta_{1}\left(\varphi\right)=\varphi\left(0\right),\quad\beta_{2}\left(\varphi\right)=\varphi'\left(0\right);\;\mbox{then}
\]
\[
\left\langle g,A^{*}f\right\rangle -\left\langle A^{*}g,f\right\rangle =\left\langle \beta_{1}\left(g\right),\beta_{2}\left(f\right)\right\rangle _{b}-\left\langle \beta_{2}\left(g\right),\beta_{1}\left(f\right)\right\rangle _{b}.
\]
This defines the boundary triple.

The selfadjoint extensions are parameterized by $e^{i\theta}$, where
\[
\left(1-e^{i\theta}\right)\beta_{1}\left(f\right)=i\left(1+e^{i\theta}\right)\beta_{2}\left(f\right);\;\mbox{i.e.,}
\]
\[
f\left(0\right)=zf'\left(0\right),\;f\in\mathscr{D}\left(A^{*}\right)
\]
where 
\[
z=i\frac{1+e^{i\theta}}{1-e^{i\theta}}.
\]
We take the convention that $z=\infty\Longleftrightarrow f'\left(0\right)=0$,
i.e., the Neumann boundary condition. \index{boundary condition}\index{selfadjoint extensions}
\end{example}

\begin{example}
\label{exa:bt2}$Af=-f''$, $\mathscr{D}\left(A\right)=C_{c}^{\infty}\left(0,1\right)$;
then 
\[
\mathscr{D}\left(A^{*}\right)=\left\{ f,f''\in L^{2}\left(0,1\right)\right\} .
\]
Integration by parts gives 
\begin{align*}
\left\langle g,A^{*}f\right\rangle  & =-\int_{0}^{1}\overline{g}f''\\
 & =-\left[\overline{g}f'-\overline{g}'f\right]_{0}^{1}-\int_{0}^{1}\overline{g}''f\\
 & =-\left[\overline{g}f'-\overline{g}'f\right]_{0}^{1}+\left\langle A^{*}g,f\right\rangle .
\end{align*}
Thus, 
\begin{align*}
\left\langle g,A^{*}f\right\rangle -\left\langle A^{*}g,f\right\rangle  & =\left[\left(\overline{g}f'\right)\left(0\right)+\left(\overline{g}'f\right)\left(1\right)\right]-\left[\left(\overline{g}'f\right)\left(0\right)+\left(\overline{g}f'\right)\left(1\right)\right]\\
 & =\left\langle \beta_{1}\left(g\right),\beta_{2}\left(f\right)\right\rangle _{b}-\left\langle \beta_{2}\left(g\right),\beta_{1}\left(f\right)\right\rangle _{b}
\end{align*}
where 
\[
\beta_{1}\left(\varphi\right)=\begin{bmatrix}\varphi\left(0\right)\\
\varphi'\left(1\right)
\end{bmatrix},\;\beta_{2}\left(\varphi\right)=\begin{bmatrix}\varphi'\left(0\right)\\
\varphi\left(1\right)
\end{bmatrix}.
\]
The boundary space is $\mathscr{H}_{b}=\mathbb{C}^{2}$, i.e., 2-dimensional
\end{example}
The family of selfadjoint extensions is parameterized by $U\in M\left(2,\mathbb{C}\right)$.
Given $U$, the corresponding extension $\widetilde{A_{U}}$ is determined
by 
\begin{gather*}
\widetilde{A_{U}}=A^{*}\Big|_{\mathscr{D}\left(\widetilde{A_{U}}\right)},\;\mbox{where}\\
\mathscr{D}\left(\widetilde{A_{U}}\right)=\left\{ f\in\mathscr{D}\left(A^{*}\right):\left(1-U\right)\beta_{1}\left(f\right)=i\left(1+U\right)\beta_{2}\left(f\right)\right\} .
\end{gather*}

\begin{rem}
Another choice of the boundary map: 
\begin{align*}
\left\langle g,A^{*}f\right\rangle -\left\langle A^{*}g,f\right\rangle  & =\left[\left(\overline{g}f'\right)\left(0\right)-\left(\overline{g}f'\right)\left(1\right)\right]-\left[\left(\overline{g}'f\right)\left(0\right)-\left(\overline{g}'f\right)\left(1\right)\right]\\
 & =\left\langle \beta_{1}\left(g\right),\beta_{2}\left(f\right)\right\rangle _{b}-\left\langle \beta_{2}\left(g\right),\beta_{1}\left(f\right)\right\rangle _{b}
\end{align*}
where 
\[
\beta_{1}\left(\varphi\right)=\begin{bmatrix}\varphi\left(0\right)\\
\varphi\left(1\right)
\end{bmatrix},\quad\beta_{2}\left(\varphi\right)=\begin{bmatrix}\varphi'\left(0\right)\\
-\varphi'\left(1\right)
\end{bmatrix}.
\]
The selfadjoint boundary condition leads to 
\[
\left(1-U\right)\begin{bmatrix}f\left(0\right)\\
f\left(1\right)
\end{bmatrix}=i\left(1+U\right)\begin{bmatrix}f'\left(0\right)\\
-f'\left(1\right)
\end{bmatrix}.
\]
For $U=1$, we get the Neumann boundary condition: 
\[
f'\left(0\right)=f'\left(1\right)=0.
\]
For $U=-1$, we get the Dirichlet boundary condition: 
\[
f\left(0\right)=f\left(1\right)=0.
\]

\end{rem}
\index{boundary condition}
\begin{xca}[From selfadjoint extension to unitary one-parameter group]
\myexercise{From selfadjoint extension to unitary one-parameter group}~
\begin{enumerate}
\item For each of the selfadjoint extensions from \thmref{bt}, write down
the corresponding unitary one-parameter group; and identify it as
an induced representation; induction $\mathbb{Z}\longrightarrow\mathbb{R}$;
see \secref{ind}.
\item Same question for the selfadjoint extension operators computed in
Examples \ref{exa:bt1} and \ref{exa:bt2}.
\end{enumerate}
\end{xca}

\section{\label{sec:fried}The Friedrichs Extension }

\begin{figure}[h]
\[
\xymatrix{\mathscr{H}_{A}\ar@/^{1pc}/[rr]^{J} &  & \mathscr{H}\ar@/^{1pc}/[ll]^{J^{*}}}
\]

\protect\caption{\label{fig:fj}The operator $J$ and its adjoint.}
\end{figure}
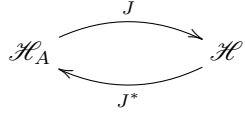

Let $A:\mathscr{D}\rightarrow\mathscr{H}$ be an operator with dense
domain $dom\left(A\right):=\mathscr{D}$ in $\mathscr{H}$, such that
\begin{equation}
\left\langle \varphi,A\varphi\right\rangle \geq\left\Vert \varphi\right\Vert ^{2},\quad\forall\varphi\in\mathscr{D}.\label{eq:fj1}
\end{equation}
Set $\mathscr{H}_{A}:=$ Hilbert completion of $\mathscr{D}$ with
respect to the 
\begin{equation}
\left\Vert \varphi\right\Vert _{A}:=\left\langle \varphi,A\varphi\right\rangle ^{\frac{1}{2}}.\label{eq:fj2}
\end{equation}
Then $\varphi\rightarrow\varphi$ defines a contraction $J:\mathscr{H}_{A}\rightarrow\mathscr{H}$,
extending $J\varphi=\varphi$, for $\varphi\in\mathscr{D}$. Note
that (\ref{eq:fj1}) $\Leftrightarrow$ 
\[
\left\Vert J\varphi\right\Vert \leq\left\Vert \varphi\right\Vert _{A},\quad\forall\varphi\in\mathscr{H}_{A},
\]
(see (\ref{eq:fj2}).) 
\begin{rem}
We will make use of two inner products: $\left\langle \cdot,\cdot\right\rangle $
in $\mathscr{H}$, and $\left\langle \cdot,\cdot\right\rangle _{A}$
(with subscript $A$) in $\mathscr{H}_{A}$. 
\end{rem}
We have
\begin{equation}
\left\langle J\varphi,f\right\rangle =\left\langle \varphi,J^{*}f\right\rangle _{A}\label{eq:jf3}
\end{equation}
see \figref{fj}: $\varphi\in\mathscr{H}_{A}$, $f\in\mathscr{H}$,
$J^{*}f\in\mathscr{H}_{A}$.

Note both $J$ and $J^{*}$ are contractions with respect to the respective
norms, so
\begin{equation}
\left\Vert J^{*}f\right\Vert _{A}\leq\left\Vert f\right\Vert ,\quad\forall f\in\mathscr{H}.\label{eq:jf4}
\end{equation}
So $JJ^{*}:\mathscr{H}\rightarrow\mathscr{H}$ is a contractive selfadjoint
operator in $\mathscr{H}$, and $\left(JJ^{*}\right)^{\frac{1}{2}}$
is well defined by the Spectral Theorem.
\begin{thm}
Let $A,J,\mathscr{D},\mathscr{H},\mathscr{H}_{A}$ be as above. Then
there is a selfadjoint extension $\widetilde{A}\supset A$ in $\mathscr{H}$
such that 
\[
\langle\widetilde{A}x,y\rangle=\left\langle x,y\right\rangle _{1},\;\forall x\in dom(\widetilde{A}),\forall y\in dom\left(A\right)\left(=\mathscr{D}\right).
\]
\end{thm}
\begin{proof}
The theorem is established in three steps:

Step 1. $JJ^{*}A\varphi=\varphi$, $\forall\varphi\in\mathscr{D}$.

Step 2. $JJ^{*}$ is invertible (easy from (\ref{eq:fj1}).)

Step 3. $A\subset\left(JJ^{*}\right)^{-1}$, where $\left(JJ^{*}\right)^{-1}$
is selfadjoint; it is the Friedrichs extension of $A$. Note Step
3 is immediate from step 1 by definition.

\emph{Proof of Step 1}. Since $\mathscr{D}$ is dense in $\mathscr{H}$,
it is enough to prove that
\begin{equation}
\left\langle \psi,JJ^{*}A\varphi\right\rangle =\left\langle \psi,\varphi\right\rangle ,\quad\forall\varphi,\psi\in\mathscr{D}.\label{eq:jf8}
\end{equation}
Let $\varphi,\psi\in\mathscr{D}$, then 
\begin{eqnarray*}
\mbox{LHS}_{\left(\ref{eq:jf8}\right)} & = & \left\langle JJ^{*}\psi,A\varphi\right\rangle \\
 & = & \left\langle J^{*}\psi,A\varphi\right\rangle \\
 & = & \left\langle J^{*}\psi,\varphi\right\rangle _{A}\;\left(\mbox{by \ensuremath{\left(\ref{eq:fj1}\right)} \& \ensuremath{\left(\ref{eq:fj2}\right)}}\right)\\
 & = & \left\langle \psi,J\varphi\right\rangle \;\left(\mbox{by \ensuremath{\left(\ref{eq:jf3}\right)} and use \ensuremath{J^{**}=J}}\right)\\
 & = & \left\langle \psi,\varphi\right\rangle =\mbox{RHS}_{\left(\ref{eq:jf8}\right)}.
\end{eqnarray*}
Hence Step 1 follows.
\end{proof}
\index{extension!Friedrich's-}

Let $A$ be a densely defined Hermitian operator in a Hilbert space
$\mathscr{H}$. $A$ is \emph{semi-bounded} if $A\geq c>-\infty$,
in the sense that, $\left\langle x,Ax\right\rangle \geq c\left\langle x,x\right\rangle $,
$\forall x\in dom\left(A\right)$. Set\index{Friedrichs extension}
\begin{equation}
L_{A}:=\inf\left\{ \left\langle x,Ax\right\rangle :x\in dom\left(A\right),\:\left\Vert x\right\Vert =1\right\} \label{eq:tr-low}
\end{equation}
and $L_{A}$ is called the \emph{lower bound} of $A$. 

In the following discussion, we first assume $A\geq1$ and eventually
drop the constraint. 

Let $\mathscr{H}_{1}=$ completion of $dom\left(A\right)$ with respect
to the inner product
\begin{align}
\left\langle x,y\right\rangle _{1} & :=\left\langle x,Ay\right\rangle \label{eq:tr-1}
\end{align}

\begin{thm}
\label{thm:AH}Let $A\geq I$, i.e., $\left\langle x,Ax\right\rangle \geq\left\Vert x\right\Vert ^{2}$,
$\forall x\in dom\left(A\right)$; and let $\mathscr{H}$, and $\mathscr{H}_{1}$
be as above. Then
\begin{enumerate}
\item \label{enu:fri1}$\left\Vert x\right\Vert \leq\left\Vert x\right\Vert _{1},\;\forall x\in dom\left(A\right)$.
\item \label{enu:fri2}$\left\Vert \cdot\right\Vert $ and $\left\Vert \cdot\right\Vert _{1}$
are topologically consistent, i.e., the identity map 
\[
\varphi:dom\left(A\right)\rightarrow dom\left(A\right)
\]
extends by continuity to 
\[
\widetilde{\varphi}:\mathscr{H}_{1}\hookrightarrow\mathscr{H}
\]
such that 
\begin{equation}
\left\Vert x\right\Vert \leq\left\Vert x\right\Vert _{1},\;\forall x\in\mathscr{H}_{1}.\label{eq:tr-3-1}
\end{equation}
Therefore, $\mathscr{H}_{1}$ is identified as a dense subspace in
$\mathscr{H}$. 
\item \label{enu:fri3}Moreover, 
\begin{equation}
\left\langle y,x\right\rangle _{1}=\left\langle y,Ax\right\rangle ,\;\forall x\in dom\left(A\right),\forall y\in\mathscr{H}_{1}.\label{eq:tr-2-1}
\end{equation}

\item \label{enu:fri4}Define 
\begin{gather*}
\widetilde{A}:=A^{*}\Big|_{dom(\widetilde{A})},\mbox{ where }\\
dom(\widetilde{A}):=dom(\widetilde{A}^{*})\cap\mathscr{H}_{1}.
\end{gather*}
Then $\widetilde{A}=\widetilde{A}^{*}$, and $L_{\widetilde{A}}=L_{A}$.
\end{enumerate}
\end{thm}
\begin{proof}
(\ref{enu:fri1})-(\ref{enu:fri2}) The assumption $A\geq1$ implies
that 
\[
\left\Vert x\right\Vert ^{2}=\left\langle x,x\right\rangle \leq\left\langle x,Ax\right\rangle =\left\Vert x\right\Vert _{1}^{2},\;\forall x\in dom\left(A\right).
\]
Hence $\varphi$ is continuous and the norm ordering passes to the
completions of $dom\left(A\right)$ with respect to $\left\Vert \cdot\right\Vert _{1}$
and $\left\Vert \cdot\right\Vert $. Therefore (\ref{eq:tr-3-1})
holds. 

Next, we verify that $\widetilde{\varphi}$ is injective (i.e., $ker\widetilde{\varphi}=0$.)
Suppose $\left(x_{n}\right)\subset dom\left(A\right)$ such that $x_{n}\xrightarrow{\left\Vert \cdot\right\Vert _{1}}x\in\mathscr{H}_{1}$,
and $x_{n}\xrightarrow{\left\Vert \cdot\right\Vert }0$. We must show
that $\left\Vert x\right\Vert _{1}=0$. But 
\begin{align*}
\left\Vert x\right\Vert _{1}^{2} & =\lim_{m,n\rightarrow\infty}\left\langle x_{m},x_{n}\right\rangle _{1}\qquad\left(\mbox{limit exists by assumption}\right)\\
 & =\lim_{m,n\rightarrow\infty}\left\langle x_{m},Ax_{n}\right\rangle \\
 & =\lim_{n\rightarrow\infty}\left\langle 0,Ax_{n}\right\rangle =0.
\end{align*}
In the computation, we used the fact that 
\begin{align*}
\left|\left\langle x_{m}-x,Ax_{n}\right\rangle \right| & \leq\left\Vert x_{m}-x\right\Vert \left\Vert Ax_{n}\right\Vert \\
 & \leq\left\Vert x_{m}-x\right\Vert _{1}\left\Vert Ax_{n}\right\Vert \rightarrow0,\;\mbox{as \ensuremath{m\rightarrow\infty}.}
\end{align*}

(\ref{enu:fri3}) Let $\left(y_{n}\right)\subset dom\left(A\right)$,
and $\left\Vert y_{n}-y\right\Vert _{1}\rightarrow0$. For all $x\in dom\left(A\right)$,
we have 
\[
\left\langle y,x\right\rangle _{1}=\lim_{n\rightarrow\infty}\left\langle y_{n},x\right\rangle _{1}=\lim_{n\rightarrow\infty}\left\langle y_{n},Ax\right\rangle =\left\langle y,Ax\right\rangle .
\]
Equivalently, 
\begin{align*}
\left|\left\langle y_{n},Ax\right\rangle -\left\langle y,Ax\right\rangle \right| & =\left|\left\langle y_{n}-y,Ax\right\rangle \right|\\
 & \leq\left\Vert y_{n}-y\right\Vert \left\Vert Ax\right\Vert \\
 & \leq\left\Vert y_{n}-y\right\Vert _{1}\left\Vert Ax\right\Vert \rightarrow0,\;\mbox{as \ensuremath{n\rightarrow\infty}.}
\end{align*}

(\ref{enu:fri4}) For all $x,y\in dom(\widetilde{A})$, $\exists\:\left(x_{n}\right),\left(y_{n}\right)\subset dom\left(A\right)$
such that $\left\Vert x_{n}-x\right\Vert _{1}\rightarrow0$ and $\left\Vert y_{n}-y\right\Vert _{1}\rightarrow0$.
Hence the following limit exists: 
\[
\lim_{m,n\rightarrow\infty}\left\langle x_{m},Ay_{n}\right\rangle \left(=\lim_{m,n\rightarrow\infty}\left\langle x_{m},y_{n}\right\rangle _{1}\right).
\]
Consequently, 
\begin{align*}
\lim_{m\rightarrow\infty}\lim_{n\rightarrow\infty}\left\langle x_{m},Ay_{n}\right\rangle  & =\lim_{m\rightarrow\infty}\lim_{n\rightarrow\infty}\left\langle Ax_{m},y_{n}\right\rangle \\
 & =\lim_{m\rightarrow\infty}\left\langle Ax_{m},y\right\rangle \\
 & =\lim_{m\rightarrow\infty}\left\langle x_{m},A^{*}y\right\rangle \\
 & =\lim_{m\rightarrow\infty}\langle x_{m},\widetilde{A}y\rangle=\langle x,\widetilde{A}y\rangle
\end{align*}
and 
\begin{align*}
\lim_{n\rightarrow\infty}\lim_{m\rightarrow\infty}\left\langle x_{m},Ay_{n}\right\rangle  & =\lim_{n\rightarrow\infty}\left\langle x,Ay_{n}\right\rangle \\
 & =\lim_{n\rightarrow\infty}\left\langle A^{*}x,y_{n}\right\rangle \\
 & =\lim_{n\rightarrow\infty}\langle\widetilde{A}x,y_{n}\rangle=\langle\widetilde{A}x,y\rangle.
\end{align*}
Thus, $\widetilde{A}$ is Hermitian. 

Fix $y\in\mathscr{H}$. The map $x\longmapsto\left\langle y,x\right\rangle $,
$\forall x\in dom\left(A\right)\subset\mathscr{H}_{1}$, is linear
and satisfies 
\[
\left|\left\langle y,x\right\rangle \right|\leq\left\Vert y\right\Vert \left\Vert x\right\Vert \leq\left\Vert y\right\Vert \left\Vert x\right\Vert _{1}.
\]
Hence it extends to a unique bounded linear functional on $\mathscr{H}_{1}$,
as $dom\left(A\right)$ is dense in $\mathscr{H}_{1}$. \index{functional}

By Riesz's theorem, there exists unique $h_{y}\in\mathscr{H}_{1}$
such that 
\begin{equation}
\left\langle y,x\right\rangle =\left\langle h_{y},x\right\rangle _{1},\;\forall x\in\mathscr{H}_{1}.\label{eq:psi}
\end{equation}
In particular,\index{Theorem!Riesz-} 
\[
\left\langle y,x\right\rangle =\left\langle h_{y},x\right\rangle _{1}=\left\langle h_{y},Ax\right\rangle ,\;\forall x\in dom\left(A\right).
\]
Then, $h_{y}\in\mathscr{H}_{1}\cap dom(A^{*})=dom(\widetilde{A})$,
and $\widetilde{A}h_{y}=y$. Therefore, $ran(\widetilde{A})=\mathscr{H}$.
Note we have established the identity
\begin{equation}
\langle\widetilde{A}y,x\rangle=\left\langle y,x\right\rangle _{1},\:\forall y\in dom(\widetilde{A}),\forall x\in dom\left(A\right).\label{eq:defeq}
\end{equation}
\end{proof}
\begin{claim}
$ran(\widetilde{A})=\mathscr{H}$ implies that $\widetilde{A}$ is
selfadjoint. In fact, for all $x\in dom(\widetilde{A})$ and $y\in dom(\widetilde{A}^{*})$,
we have 
\[
\langle y,\widetilde{A}x\rangle=\langle\widetilde{A}^{*}y,x\rangle=\langle\widetilde{A}h,x\rangle=\langle h,\widetilde{A}x\rangle
\]
where $\widetilde{A}^{*}y=\widetilde{A}h$, for some $h\in dom(\widetilde{A})$,
using the assumption $ran(\widetilde{A})=\mathscr{H}$. Thus, 
\[
\langle y-h,\widetilde{A}x\rangle=0,\;\forall x\in dom(\widetilde{A});
\]
i.e., $y-h\perp ran(\widetilde{A})=\mathscr{H}$. Therefore, $y=h$,
and so $y\in dom(\widetilde{A})$.\end{claim}
\begin{proof}
Finally, we show $L_{\widetilde{A}}=L_{A}$. By the definition of
lower bound, $dom(A)\subset dom(\widetilde{A})$ implies $L_{\widetilde{A}}\leq L_{A}$.
On the other hand, let $\left(x_{n}\right)\subset dom\left(A\right)$
such that $x_{n}\xrightarrow{\left\Vert \cdot\right\Vert _{1}}x\in dom(\widetilde{A})$,
then 
\begin{align*}
\langle x,\widetilde{A}x\rangle & =\lim_{n\rightarrow\infty}\langle x,\widetilde{A}x_{n}\rangle=\lim_{n\rightarrow\infty}\langle x,Ax_{n}\rangle\\
 & =\lim_{n\rightarrow\infty}\left\langle x,x_{n}\right\rangle _{1}=\left\langle x,x\right\rangle _{1}\geq L_{A}\left\langle x,x\right\rangle 
\end{align*}
which shows that $L_{\widetilde{A}}\geq L_{A}$.\end{proof}
\begin{rem}
In the proof of \thmref{AH}, we established an embedding $\psi:\mathscr{H}\hookrightarrow\mathscr{H}_{1}^{*}$
by $\psi:y\mapsto h_{y}$ with the defining equation (\ref{eq:psi}).
And we define $\widetilde{A}h_{y}=y$, i.e., $h_{y}=\widetilde{A}^{-1}y$.
It follows that $dom(\widetilde{A})=ran\left(\psi\right)$, and $ran(\widetilde{A})=\mathscr{H}$. \end{rem}
\begin{thm}
The Friedrichs extension of $A$ is the unique selfadjoint operator
satisfying 
\[
\langle\widetilde{A}x,y\rangle=\left\langle x,y\right\rangle _{1},\;\forall x\in dom(\widetilde{A}),\forall y\in dom\left(A\right).
\]
See (\ref{eq:defeq}). \index{selfadjoint operator}\end{thm}
\begin{proof}
Suppose $A\subset B,C\subset A^{*}$, and $B,C$ selfadjoint, satisfying
\begin{align*}
\left\langle Bx,y\right\rangle  & =\left\langle x,y\right\rangle _{1},\;\forall x\in dom\left(B\right),\forall y\in dom\left(A\right)\\
\left\langle Cx,y\right\rangle  & =\left\langle x,y\right\rangle _{1},\;\forall x\in dom\left(C\right),\forall y\in dom\left(A\right).
\end{align*}
Then, for all $x,y\in dom\left(A\right)$, we have 
\[
\left\langle Bx,y\right\rangle =\left\langle Cx,y\right\rangle =\left\langle x,Cy\right\rangle \left(=\left\langle x,Ay\right\rangle \right).
\]
Fix $x\in dom\left(A\right)$, and the above identify passes to $y\in dom\left(C\right)$.
Therefore, $y\in B^{*}=B$, and $By=Cy$. This shows $C\subset B$.
Since 
\[
C=C^{*}\supset B^{*}=B
\]
i.e., $C\supset B$, it then follows that $B=C$.\end{proof}
\begin{rem}
If $A$ is only assumed to be semi-bounded, i.e., $A\geq c>-\infty$,
then $B:=A-c+1\geq1$, and we get the Friedrichs extension $\widetilde{B}$
of $B$; and $\widetilde{B}-1+c$ is the Friedrichs extension of $A$.
\end{rem}

\section{\label{sec:rhs}Rigged Hilbert Space}

In the construction of Friedrichs extensions of semi-bounded operators,
we have implicitly used the idea of rigged Hilbert spaces. We study
this method systematically and recover the Friedrichs extension as
a special case. \index{Hilbert space!rigged}

Let $\mathscr{H}_{0}$ be a Hilbert space with inner product $\left\langle \cdot,\cdot\right\rangle _{0}$,
and $\mathscr{H}_{1}$ be a dense subspace in $\mathscr{H}_{0}$,
which by itself, is a Hilbert space with respect to $\left\langle \cdot,\cdot\right\rangle _{1}$.
Further, assume the ordering\index{space!Hilbert-} 
\begin{equation}
\left\Vert x\right\Vert _{0}\leq\left\Vert x\right\Vert _{1},\;\forall x\in\mathscr{H}_{1}.\label{eq:fr-1}
\end{equation}
Hence, the identity map
\begin{equation}
id:\mathscr{H}_{1}\hookrightarrow\mathscr{H}_{0}\label{eq:fr-2}
\end{equation}
is continuous with a dense image.

Let $\mathscr{H}_{-1}$ be the space of bounded \emph{conjugate} linear
functionals on $\mathscr{H}_{1}$. By Riesz's theorem, $\mathscr{H}_{-1}$
is identified with $\mathscr{H}_{1}$ via the map \index{Riesz' theorem}\index{Theorem!Riesz-}
\begin{gather}
\mathscr{H}_{-1}\rightarrow\mathscr{H}_{1},\;f\mapsto\xi_{f},\;\mbox{s,t.}\label{eq:fr-3}\\
f\left(x\right)=\left\langle x,\xi_{f}\right\rangle _{1},\;\forall x\in\mathscr{H}_{1}.\label{eq:fr-4}
\end{gather}
Then $\mathscr{H}_{-1}$ is a Hilbert space with respect to the inner
product 
\begin{equation}
\left\langle f,g\right\rangle _{-1}=\left\langle \xi_{f},\xi_{g}\right\rangle _{1},\;\forall f,g\in\mathscr{H}_{-1}.\label{eq:fr-5}
\end{equation}

\begin{rem}
\label{rem:H-1}The map $f\mapsto\xi_{f}$ in (\ref{eq:fr-3}) is
linear. For if $c\in\mathbb{C}$, then $\left\langle x,\xi_{cf}\right\rangle _{1}=cf\left(x\right)=c\left\langle x,\xi_{f}\right\rangle _{1}=\left\langle x,c\xi_{f}\right\rangle _{1}$,
for all $x\in\mathscr{H}_{1}$; i.e., $\left(\xi_{cf}-c\xi_{f}\right)\perp\mathscr{H}_{1}$,
and so $\xi_{cf}=c\xi_{f}$. \end{rem}
\begin{thm}
\label{thm:fe-emb}The mapping 
\begin{equation}
\mathscr{H}_{0}\hookrightarrow\mathscr{H}_{-1},\;x\mapsto\left\langle \cdot,x\right\rangle _{0},\;\forall x\in\mathscr{H}_{0}\label{eq:fr-6}
\end{equation}
is linear, injective, continuous, and having a dense image in $\mathscr{H}_{-1}$.\end{thm}
\begin{proof}
Since $cx\mapsto\left\langle \cdot,cx\right\rangle =c\left\langle \cdot,x\right\rangle $,
the mapping in (\ref{eq:fr-6}) is linear.

For all $x\in\mathscr{H}_{0}$, we have \ref{eq:fr1}
\begin{equation}
\left|\left\langle y,x\right\rangle _{0}\right|\leq\left\Vert x\right\Vert _{0}\left\Vert y\right\Vert _{0}\overset{\left(\ref{eq:fr-1}\right)}{\leq}\left\Vert x\right\Vert _{0}\left\Vert y\right\Vert _{1},\;\forall y\in\mathscr{H}_{1};\label{eq:fr-7}
\end{equation}
hence $\left\langle \cdot,x\right\rangle _{0}$ is a bounded conjugate
linear functional on $\mathscr{H}_{1}$, i.e., $\left\langle \cdot,x\right\rangle _{0}\in\mathscr{H}_{-1}$.
Moreover, by (\ref{eq:fr-7}), 
\begin{equation}
\left\Vert \left\langle \cdot,x\right\rangle _{0}\right\Vert _{-1}\leq\left\Vert x\right\Vert _{0}.\label{eq:fr-8}
\end{equation}

If $\left\langle \cdot,x\right\rangle _{0}\equiv0$ in $\mathscr{H}_{-1}$,
then $\left\langle y,x\right\rangle _{0}=0$, for all $y\in\mathscr{H}_{1}$.
Since $\mathscr{H}_{1}$ is dense in $\mathscr{H}_{0}$, it follows
that $x=0$ in $\mathscr{H}_{0}$. Thus, (\ref{eq:fr-6}) is injective. 

Now, if $f\perp\left\{ \left\langle \cdot,x\right\rangle _{0}:x\in\mathscr{H}_{1}\right\} $
in $\mathscr{H}_{-1}$, then 
\begin{equation}
\left\langle f,\left\langle \cdot,x\right\rangle _{0}\right\rangle _{-1}\overset{\left(\ref{eq:fr-5}\right)}{=}\left\langle \xi_{f},\xi_{\left\langle \cdot,x\right\rangle _{0}}\right\rangle _{1}\overset{\left(\ref{eq:fr-3}\right)}{=}\left\langle \xi_{f},x\right\rangle _{0}=0,\;\forall x\in\mathscr{H}_{1}.\label{eq:fr-9}
\end{equation}
Thus, $\left\Vert \xi_{f}\right\Vert _{0}=0$, since $\mathscr{H}_{1}$
is dense in $\mathscr{H}_{0}$. Since $id:\mathscr{H}_{1}\hookrightarrow\mathscr{H}_{0}$
is injective, and $\xi_{f}\in\mathscr{H}_{1}$, it follows that $\left\Vert \xi_{f}\right\Vert _{1}=0$.
This, in turn, implies $\left\Vert f\right\Vert _{-1}=0$, and so
$f=0$ in $\mathscr{H}_{-1}$. Consequently, the image of $\mathscr{H}_{1}$
(resp. $\mathscr{H}_{0}$ as it contains $\mathscr{H}_{1}$) under
(\ref{eq:fr-6}) is dense in $\mathscr{H}_{-1}$. 
\end{proof}
Combining (\ref{eq:fr-1}) and \thmref{fe-emb}, we get the triple
of Hilbert spaces
\begin{equation}
\mathscr{H}_{1}\xhookrightarrow[]{\left(\ref{eq:fr-2}\right)}\mathscr{H}_{0}\xhookrightarrow[]{\left(\ref{eq:fr-6}\right)}\mathscr{H}_{-1}.\label{eq:fr-10}
\end{equation}

The following are immediate:
\begin{enumerate}
\item All mappings in (\ref{eq:fr-10}) are injective, continuous (in fact,
contractive), having dense images. 
\item The map $x\mapsto\xi_{\left\langle \cdot,x\right\rangle _{0}}$ is
a contraction from $\mathscr{H}_{1}\subset\mathscr{H}_{0}$ into $\mathscr{H}_{0}$.
This follows from the estimate: 
\begin{align}
\left\Vert \xi_{\left\langle \cdot,x\right\rangle _{0}}\right\Vert _{0} & \overset{\left(\ref{eq:fr-1}\right)}{\leq}\left\Vert \xi_{\left\langle \cdot,x\right\rangle _{0}}\right\Vert _{1}\nonumber \\
 & \overset{\left(\ref{eq:fr-5}\right)}{=}\left\Vert \left\langle \cdot,x\right\rangle _{0}\right\Vert _{-1}\nonumber \\
 & \overset{\left(\ref{eq:fr-8}\right)}{\leq}\left\Vert x\right\Vert _{0}\overset{\left(\ref{eq:fr-1}\right)}{\leq}\left\Vert x\right\Vert _{1},\;\forall x\in\mathscr{H}_{1}.\label{eq:fr-11}
\end{align}
In particular, for $x\in\mathscr{H}_{1}$, $x\neq\xi_{\left\langle \cdot,x\right\rangle _{0}}$
in general. 
\item The canonical bilinear form $\left\langle \cdot,\cdot\right\rangle :\mathscr{H}_{1}\times\mathscr{H}_{-1}\rightarrow\mathbb{C}$
is given by 
\begin{equation}
f\left(x\right)=\left\langle x,\xi_{f}\right\rangle _{1},\;\forall x\in\mathscr{H}_{1},\forall f\in\mathscr{H}_{-1}.\label{eq:fr-12}
\end{equation}
In particular, if $f=\left\langle \cdot,y\right\rangle _{0}$, $y\in\mathscr{H}_{0}$,
then 
\[
\left\langle x,\xi_{f}\right\rangle _{1}=\left\langle x,\xi_{\left\langle \cdot,y\right\rangle _{0}}\right\rangle _{1}=\left\langle x,y\right\rangle _{0},\;\forall x\in\mathscr{H}_{1}.
\]

\item By \thmref{fe-emb}, $\mathscr{H}_{0}$ is dense in $\mathscr{H}_{-1}$,
and 
\begin{align}
\left\langle x,y\right\rangle _{-1} & :=\left\langle \left\langle \cdot,x\right\rangle _{0},\left\langle \cdot,y\right\rangle _{0}\right\rangle _{-1}\nonumber \\
 & =\left\langle \xi_{\left\langle \cdot,x\right\rangle _{0}},\xi_{\left\langle \cdot,y\right\rangle _{0}}\right\rangle _{1}\nonumber \\
 & =\left\langle x,\xi_{\left\langle \cdot,y\right\rangle _{0}}\right\rangle _{0}=\left\langle \xi_{\left\langle \cdot,x\right\rangle _{0}},y\right\rangle _{0},\;\forall x,y\in\mathscr{H}_{0}.\label{eq:fr-12-1}
\end{align}
Combined with the order relation $\left\Vert x\right\Vert _{-1}\leq\left\Vert x\right\Vert _{0}$
for all $x\in\mathscr{H}_{0}$, we see that 
\begin{align*}
\left|\left\langle x,y\right\rangle _{-1}\right| & =\left|\left\langle x,\xi_{\left\langle \cdot,y\right\rangle _{0}}\right\rangle _{0}\right|\\
 & \leq\left\Vert x\right\Vert _{0}\left\Vert \xi_{\left\langle \cdot,y\right\rangle _{0}}\right\Vert _{-1}\\
 & =\left\Vert x\right\Vert _{0}\left\Vert y\right\Vert _{-1}.\\
 & \leq\left\Vert x\right\Vert _{0}\left\Vert y\right\Vert _{0},\,,\forall x,y\in\mathscr{H}_{0}.
\end{align*}
Thus $\left\langle \cdot,\cdot\right\rangle _{-1}$ is a continuous
extension of $\left\langle \cdot,\cdot\right\rangle _{0}$. \end{enumerate}
\begin{thm}
\label{thm:fr-B}Let $\mathscr{H}_{1}\hookrightarrow\mathscr{H}_{0}\hookrightarrow\mathscr{H}_{-1}$
be the triple in (\ref{eq:fr-10}). Define $B:\mathscr{H}_{0}\rightarrow\mathscr{H}_{0}$
by
\begin{equation}
B:\underset{\mathscr{H}_{0}}{x}\xmapsto[]{\left(\ref{eq:fr-6}\right)}\underset{\mathscr{H}_{-1}}{\left\langle \cdot,x\right\rangle _{0}}\xmapsto[]{\left(\ref{eq:fr-3}\right)}\underset{\mathscr{H}_{1}}{\xi_{\left\langle \cdot,x\right\rangle _{0}}}\xmapsto[]{\left(\ref{eq:fr-2}\right)}\underset{\mathscr{H}_{0}}{\xi_{\left\langle \cdot,x\right\rangle _{0}}}.\label{eq:fr-13}
\end{equation}
Then, 
\begin{enumerate}
\item For all $x\in\mathscr{H}_{1}$, and all $y\in\mathscr{H}_{0}$, 
\begin{equation}
\left\langle x,By\right\rangle _{1}=\left\langle x,y\right\rangle _{0}.\label{eq:fr-15}
\end{equation}
In particular, 
\begin{align}
\left\langle x,y\right\rangle _{-1} & =\left\langle Bx,By\right\rangle _{1}\nonumber \\
 & =\left\langle x,By\right\rangle _{0}=\left\langle Bx,y\right\rangle _{0},\;\forall x,y\in\mathscr{H}_{0};\label{eq:fr-15-1}
\end{align}
where $\left\langle x,y\right\rangle _{-1}:=\left\langle \xi_{\left\langle \cdot,x\right\rangle _{0}},\xi_{\left\langle \cdot,y\right\rangle _{0}}\right\rangle _{1}$,
as defined in (\ref{eq:fr-12-1}).
\item $B$ is invertible. 
\item $ran\left(B\right)$ is dense in both $\mathscr{H}_{1}$ and $\mathscr{H}_{0}$. 
\item $0\leq B\leq1$. In particular, $B$ is a bounded selfadjoint operator
on $\mathscr{H}_{0}$.
\end{enumerate}
\end{thm}
\begin{proof}
~\end{proof}
\begin{enumerate}
\item For $y\in\mathscr{H}_{0}$, $By=\xi_{\left\langle \cdot,y\right\rangle _{0}}\in\mathscr{H}_{1}$,
where $\xi_{\left\langle \cdot,y\right\rangle _{0}}$ is given in
(\ref{eq:fr-3})-(\ref{eq:fr-4}). Thus, 
\[
\left\langle x,By\right\rangle _{1}=\left\langle x,\xi_{\left\langle \cdot,y\right\rangle _{0}}\right\rangle _{1}=\left\langle x,y\right\rangle _{0},\;\forall x\in\mathscr{H}_{1}.
\]
(\ref{eq:fr-15-1}) follows from this.

\begin{enumerate}
\item If $\left\Vert Bx\right\Vert _{0}=0$, then $\left\Vert Bx\right\Vert _{1}=0$,
since $Bx\in\mathscr{H}_{1}$ and $\mathscr{H}_{1}\hookrightarrow\mathscr{H}_{0}$
is injective. But then 
\[
\left\Vert Bx\right\Vert _{1}\overset{\left(\ref{eq:fr-15-1}\right)}{=}\left\Vert x\right\Vert _{-1}=0\Longrightarrow\left\Vert x\right\Vert _{0}=0
\]
since $\mathscr{H}_{0}\hookrightarrow\mathscr{H}_{-1}$ is injective.
This shows that $B$ is injective.
\item Since $\mathscr{H}_{0}\hookrightarrow\mathscr{H}_{-1}$ is dense,
and $\mathscr{H}_{-1}\simeq\mathscr{H}_{1}$, it follows that $ran\left(B\right)$
is dense in $\mathscr{H}_{1}$. Now if $y\in\mathscr{H}_{0}$, and
$\left\langle y,Bx\right\rangle _{0}=0$, for all $x\in\mathscr{H}_{0}$,
then 
\[
\left\langle By,Bx\right\rangle _{1}=0,\;\forall x\in\mathscr{H}_{0};
\]
equivalently, 
\[
\left\langle y,x\right\rangle _{-1}=0,\;\forall x\in\mathscr{H}_{0}.
\]
Since $\mathscr{H}_{0}\hookrightarrow\mathscr{H}_{-1}$ is dense,
we have $\left\Vert y\right\Vert _{-1}=0$. But $y\in\mathscr{H}_{0}$
and $\mathscr{H}_{0}\hookrightarrow\mathscr{H}_{-1}$ is injective,
it follows that $\left\Vert y\right\Vert _{0}=0$, i.e., $y=0$ in
$\mathscr{H}_{0}$. Therefore, $ran\left(B\right)$ is also dense
in $\mathscr{H}_{0}$. 
\item For all $x\in\mathscr{H}_{0}$, we have 
\[
\left\langle x,Bx\right\rangle _{0}\overset{\left(\ref{eq:fr-15}\right)}{=}\left\langle Bx,Bx\right\rangle _{1}\geq0\Longrightarrow B\geq0.
\]
On the other hand, 
\[
\left\langle x,Bx\right\rangle _{0}=\left\langle Bx,Bx\right\rangle _{1}\overset{\left(\ref{eq:fr-15-1}\right)}{=}\left\langle x,x\right\rangle _{-1}\overset{\left(\ref{eq:fr-8}\right)}{\leq}\left\langle x,x\right\rangle _{0};
\]
and so $B\leq1$. Since $B$ is positive and bounded, it is selfadjoint.\\
\\
Another argument: 
\[
\left\Vert Bx\right\Vert _{0}\leq\left\Vert Bx\right\Vert _{1}=\left\Vert x\right\Vert _{-1}\leq\left\Vert x\right\Vert _{0},\;\forall x\in\mathscr{H}_{0}.
\]

\end{enumerate}
\end{enumerate}
In view of applications, it is convenient to reformulate the previous
theorem in terms of $B^{-1}$. 
\begin{thm}
\label{thm:fr-A}Let $B$, $\mathscr{H}_{1}$, $\mathscr{H}_{0}$,
$\mathscr{H}_{-1}$, be as in \thmref{fr-B}, set $A:=B^{-1}$, then
\begin{enumerate}
\item \label{enu:fe-2-1}$A=A^{*}$, $A\geq1$.
\item \label{enu:fe-2-2}$dom\left(A\right)$ is dense in $\mathscr{H}_{1}$
and $\mathscr{H}_{0}$ , and $ran\left(A\right)=\mathscr{H}_{0}$. 
\item \label{enu:fe-2-3}For all $y\in dom\left(A\right)$, $x\in\mathscr{H}_{1}$,
\begin{equation}
\left\langle x,y\right\rangle _{1}=\left\langle x,Ay\right\rangle _{0}.\label{eq:fr-16}
\end{equation}
In particular, 
\begin{align}
\left\langle x,y\right\rangle _{1} & =\left\langle Ax,Ay\right\rangle _{-1}\nonumber \\
 & =\left\langle Ax,y\right\rangle _{0}=\left\langle x,Ay\right\rangle _{0},\;\forall x,y\in dom\left(A\right).\label{eq:fr-16-1}
\end{align}

\item \label{enu:fe-unique}There is a unique selfadjoint operator in $\mathscr{H}_{0}$
satisfying (\ref{eq:fr-16}).
\end{enumerate}
\end{thm}
\begin{proof}
Part (\ref{enu:fe-2-1})-(\ref{enu:fe-2-3}) are immediate by \thmref{fr-B}.
For (\ref{enu:fe-unique}), suppose $A,B$ are selfadjoint in $\mathscr{H}_{0}$
such that 

(i) $dom\left(A\right)$, $dom\left(B\right)$ are contained in $\mathscr{H}_{1}$,
dense in $\mathscr{H}_{0}$;

(ii) 
\begin{align*}
\left\langle x,y\right\rangle _{1} & =\left\langle x,Ay\right\rangle _{0},\;\forall x\in\mathscr{H}_{1},y\in dom\left(A\right)\\
\left\langle x,y\right\rangle _{1} & =\left\langle x,By\right\rangle _{0},\;\forall x\in\mathscr{H}_{1},y\in dom\left(B\right).
\end{align*}
Then, for all $x\in dom\left(A\right)$ and $y\in dom\left(B\right)$,
\[
\left\langle x,By\right\rangle _{0}=\left\langle x,y\right\rangle _{1}=\left\langle Ax,y\right\rangle _{0}.
\]
Thus, $x\mapsto\left\langle Ax,y\right\rangle _{0}$ is a bounded
linear functional on $dom\left(A\right)$, and so $y\in dom\left(A^{*}\right)=dom\left(A\right)$
and $A^{*}y=Ay=By$; i.e., $A\supset B$. Since $A,B$ are selfadjoint,
then 
\[
B=B^{*}\subset A^{*}=A.
\]
Therefore $A=B$.\end{proof}
\begin{thm}
Let $\mathscr{H}_{1},\mathscr{H}_{0},A$ as in \thmref{fr-A}. Then 
\begin{enumerate}
\item \label{enu:f2-2-4}$\mathscr{H}_{1}=dom\left(A^{1/2}\right)$, and
\begin{equation}
\left\langle x,y\right\rangle _{1}=\left\langle A^{1/2}x,A^{1/2}y\right\rangle _{0},\;\forall x,y\in\mathscr{H}_{1}.\label{eq:fr-17}
\end{equation}

\item \label{enu:fe-2-5}For all $x,y\in\mathscr{H}_{0}$, 
\begin{equation}
\left\langle x,y\right\rangle _{-1}=\left\langle A^{-1/2}x,A^{-1/2}y\right\rangle _{0}.\label{eq:fr-18}
\end{equation}
Since $\mathscr{H}_{0}$ is dense in $\mathscr{H}_{-1}$, then $\mathscr{H}_{-1}=$
completion of $\mathscr{H}_{0}$ under the $\left\Vert A^{-1/2}\cdot\right\Vert _{0}$-norm.
\item \label{enu:fe-2-6}For all $x\in dom\left(A\right)$, 
\begin{equation}
\left\Vert Ax\right\Vert _{-1}=\left\Vert x\right\Vert _{1}\left(=\left\Vert A^{1/2}x\right\Vert _{0}\right).\label{eq:fr-19}
\end{equation}
Consequently, the map $dom\left(A\right)\ni x\mapsto Ax\in\mathscr{H}_{0}$
extends by continuity to a unitary operator from $\mathscr{H}_{1}\left(=dom\left(A^{1/2}\right)\right)$
onto $\mathscr{H}_{-1}$, which is precisely the inverse of (\ref{eq:fr-3})-(\ref{eq:fr-4}).
\end{enumerate}
\end{thm}
\begin{proof}
\begin{flushleft}
(\ref{enu:f2-2-4}) This is the result of the following observations:
\par\end{flushleft}\end{proof}
\begin{enumerate}[label=\alph{enumi}.]
\item  $dom\left(A\right)\subset dom\left(A^{1/2}\right)$. With the assumption
$A\geq1$, the containment is clear. The assertion also holds in general:
By spectral theorem, 
\[
x\in dom\left(A\right)\Longleftrightarrow\int\left(1+\left|\lambda\right|^{2}\right)\left\Vert P\left(d\lambda\right)x\right\Vert _{0}^{2}<\infty;
\]
where $P\left(\cdot\right)$ is the projection-valued measure (PVM)
of $A$, defined on the set of all Borel sets in $\mathbb{R}$, and
$d\mu_{x}:=\left\Vert P\left(d\lambda\right)x\right\Vert _{0}^{2}$
is a finite positive Borel measure on $\mathbb{R}$. Thus,\index{Borel measure}
\[
\int\left(1+\left|\lambda\right|\right)\left\Vert P\left(d\lambda\right)x\right\Vert _{0}^{2}<\infty
\]
since $L^{2}\subset L^{1}$ when the measure is finite. But this is
equivalent to $x\in dom\left(A^{1/2}\right)$. 

\begin{enumerate}
\item For any Hermitian operator $T$ satisfying $T\geq c>0$, we have the
estimate:
\[
\left\Vert Tx\right\Vert \leq\left\Vert x\right\Vert +\left\Vert Tx\right\Vert \leq\left(1+c\right)\left\Vert Tx\right\Vert 
\]
Thus, the graph norm of $T$ is equivalent to $\left\Vert T\cdot\right\Vert $.
\index{graph!- of operator}
\item $dom\left(A\right)$ is dense in $dom\left(A^{1/2}\right)$. Note
$dom\left(A^{1/2}\right)$ is a Hilbert space with respect to the
$A^{1/2}$-graph norm. By the discussion above, $\left\Vert \cdot\right\Vert _{A^{1/2}}\simeq\left\Vert A^{1/2}\cdot\right\Vert _{0}$.
\\
Let $y\in dom\left(A^{1/2}\right)$, then 
\begin{align*}
y\perp dom\left(A\right) & \mbox{ in }\ensuremath{dom\left(A^{1/2}\right)}\\
 & \Updownarrow\\
\left\langle A^{1/2}y,A^{1/2}x\right\rangle _{0} & =0,\;\forall x\in dom\left(A^{1/2}\right)\\
 & \Updownarrow\\
\left\langle y,Ax\right\rangle _{0} & =0,\;\forall x\in dom\left(A^{1/2}\right)\\
 & \Updownarrow\\
y=0 & \mbox{ in }\ensuremath{\mathscr{H}_{0}}
\end{align*}

\item $dom\left(A\right)$ is dense in $\mathscr{H}_{1}$. See theorems
\ref{thm:fr-B}-\ref{thm:fr-A}.
\item $\left\Vert A^{1/2}x\right\Vert _{0}=\left\Vert x\right\Vert _{1}$,
$\forall x\in dom\left(A\right)$. Indeed, 
\[
\left\langle x,x\right\rangle _{1}\overset{\left(\ref{eq:fr-16}\right)}{=}\left\langle x,Ax\right\rangle _{0}=\left\langle A^{1/2}x,A^{1/2}x\right\rangle _{0},\;\forall x\in dom\left(A\right).
\]

\end{enumerate}

\textbf{Conclusion:} (i) $dom\left(A\right)$ is dense in $\mathscr{H}_{1}$
and $dom\left(A^{1/2}\right)$; (ii) $\left\Vert \cdot\right\Vert _{1}$
and $\left\Vert A^{1/2}\cdot\right\Vert _{0}$ agree on $dom\left(A\right)$.
Therefore the closures of $dom\left(A\right)$ in $\mathscr{H}_{1}$
and $dom\left(A^{1/2}\right)$ are identical. This shows $\mathscr{H}_{1}=dom\left(A^{1/2}\right)$.
(\ref{eq:fr-17}) is immediate. 

\end{enumerate}
\begin{proof}
\begin{flushleft}
(\ref{enu:fe-2-5}) Given $x,y\in\mathscr{H}_{0}$, 
\[
\left\langle x,y\right\rangle _{-1}\overset{\left(\ref{eq:fr-15-1}\right)}{=}\left\langle A^{-1}x,A^{-1}y\right\rangle _{1}\overset{\left(\ref{eq:fr-17}\right)}{=}\left\langle A^{-1/2}x,A^{-1/2}y\right\rangle _{0}.
\]

\par\end{flushleft}

\begin{flushleft}
(\ref{enu:fe-2-6}) Given $x\in dom\left(A\right)$, we have 
\[
\left\Vert Ax\right\Vert _{-1}\overset{\left(\ref{eq:fr-15-1}\right)}{=}\left\Vert A^{-1}\left(Ax\right)\right\Vert _{1}=\left\Vert x\right\Vert _{1}\overset{\left(\ref{eq:fr-17}\right)}{=}\left\Vert A^{1/2}x\right\Vert _{0}.
\]

\par\end{flushleft}
\end{proof}
\textbf{Application: The Friedrichs extension revisited. }\index{Friedrichs extension}
\begin{thm}[Friedrichs]
\label{thm:Friedrichs} Let $A$ be a densely defined Hermitian operator
acting in $\mathscr{H}_{0}$, and assume $A\geq1$. There exists a
selfadjoint extension $S\supset A$, such that $L_{S}=L_{A}$, i.e.,
$A$ and $S$ have the same lower bound. \end{thm}
\begin{proof}
Given $A$, construct the triple $\mathscr{H}_{1}\hookrightarrow\mathscr{H}_{0}\hookrightarrow\mathscr{H}_{-1}$
as in \thmref{AH}, so that $\mathscr{H}_{1}=cl_{\left\langle \cdot,A\cdot\right\rangle _{0}}\left(dom\left(A\right)\right)$,
and 
\begin{equation}
\left\langle y,x\right\rangle _{1}=\left\langle y,Ax\right\rangle _{0},\;\forall x\in dom\left(A\right),\forall y\in\mathscr{H}_{1}.\label{eq:tr-3}
\end{equation}
By \thmref{fr-A}, there is a densely defined selfadjoint operator
$S$ in $\mathscr{H}_{0}$, such that (i) $\mathscr{H}_{1}=dom\left(S^{1/2}\right)\supset dom\left(S\right)$;
(ii) \index{selfadjoint operator} 
\begin{equation}
\left\langle y,x\right\rangle _{1}=\left\langle S^{1/2}y,S^{1/2}x\right\rangle _{0},\;\forall x,y\in\mathscr{H}_{1}.\label{eq:tr-4}
\end{equation}
Combing (\ref{eq:tr-3})-(\ref{eq:tr-4}), we get 
\[
\left\langle y,Ax\right\rangle _{0}=\left\langle S^{1/2}y,S^{1/2}x\right\rangle _{0},\;\forall x\in dom\left(A\right),\forall y\in dom\left(S^{1/2}\right).
\]
Therefore, $S^{1/2}x\in dom\left(S^{1/2}\right)$, i.e., $x\in dom\left(S\right)$,
and 
\[
\left\langle y,Ax\right\rangle _{0}=\left\langle y,Sx\right\rangle _{0},\;\forall x\in dom\left(A\right),\forall y\in dom\left(S^{1/2}\right).
\]
Since $\mathscr{H}_{1}=dom\left(S^{1/2}\right)$ is dense in $\mathscr{H}_{0}$,
we conclude that $Sx=Ax$, for all $x\in dom\left(A\right)$. Thus,
$S\supset A$. 

Clearly, $S\supset A$ implies $L_{A}\geq L_{S}$. On the other hand,
\[
\left\Vert x\right\Vert _{1}^{2}=\left\langle x,Ax\right\rangle _{0}\geq L_{A}\left\langle x,x\right\rangle _{0}=L\left\Vert x\right\Vert _{0}^{2},\;\forall x\in dom\left(A\right)
\]
and the inequality passes by continuity to all $x\in\mathscr{H}_{1}$;
i.e., 
\[
\left\Vert x\right\Vert _{1}^{2}\geq L_{A}\left\Vert x\right\Vert _{0}^{2},\;\forall x\in\mathscr{H}_{1}.
\]
This is equivalent (by \thmref{fr-A}) to 
\[
\left\langle S^{1/2}x,S^{1/2}x\right\rangle _{0}\geq L_{A}\left\langle x,x\right\rangle _{0},\;\forall x\in dom\left(S^{1/2}\right)\left(=\mathscr{H}_{1}\right).
\]
In particular, 
\[
\left\langle x,Sx\right\rangle _{0}\geq L_{A}\left\langle x,x\right\rangle _{0},\;\forall x\in dom\left(S\right)
\]
and so $L_{S}\geq L_{A}$. Therefore, $L_{A}=L_{S}$.
\end{proof}

\section*{A summary of relevant numbers from the Reference List}

For readers wishing to follow up sources, or to go in more depth with
topics above, we suggest: \cite{AG93,Ne69,Dev72,DS88b,Kre46,Jor08,Rud73,Sto51,Sto90,FL28,Fug82,MR887102,JoMu80,VN35,JPT11-1,JP14,RS75,Hel13}.

\chapter{Unbounded Graph-Laplacians\label{chap:gLap}}
\begin{quotation}
Mathematics is an experimental science, and definitions do not come
first, but later on. 

--- Oliver Heaviside\sindex[nam]{Heaviside, O., (1850-1925)}\vspace{1em}\\
It is nice to know that the computer understands the problem. But
I would like to understand it too. 

--- Eugene Wigner\sindex[nam]{Wigner, E.P., (1902-1995)}\vspace{2em}\\
Knowing is not enough; we must apply. 

---Göthe \vspace{2em}\\
Chance is a more fundamental conception than causality.

--- Max Born\vspace{2em}
\end{quotation}
Below we study selfadjoint operators, and extensions in a particular
case arising in the study of infinite graphs; the operators here are
infinite discrete Laplacians. 

As an application of the previous chapter, we consider the Friedrichs
extension of discrete Laplacian in infinite networks \cite{JoPe10,JoPe13,JoPe13b,JT14,JT14c}. 

\index{discrete Laplacian}

\index{infinite network}

\index{extension!Friedrich's-}

\index{operators!Laplace-}

\index{Laplacian!graph-}\index{graph!network-}

By an electrical network we mean a graph $G$ of vertices and edges
satisfying suitable conditions which allow for computation of voltage
distribution from a network of prescribed resistors assigned to the
edges in $G$. The mathematical axioms are prescribed in a way that
facilitates the use of the laws of Kirchhoff and Ohm in computing
voltage distributions and resistance distances in $G$. It will be
more convenient to work with prescribed conductance functions $c$
on $G$. Indeed with a choice of conductance function $c$ specified
we define two crucial tools for our analysis, a graph Laplacian $\Delta\left(=\Delta_{c},\right)$
a discrete version of more classical notions of Laplacians, and an
energy Hilbert space $\mathscr{H}_{E}$. \index{graph Laplacian}\index{conductance}
\index{resistance distance}

Because of statistical consideration, and our use of random walk\index{random walk}
models, we focus our study on infinite electrical networks, i.e.,
the case when the graph $G$ is countable infinite. In this case,
for realistic models the graph Laplacian $\Delta_{c}$ will then be
an unbounded operator with dense domain in $\mathscr{H}_{E}$, Hermitian
and semibounded. Hence it has a unique Friedrichs extension.

\index{Friedrichs extension}

\index{operators!Laplace-}

\index{infinite electrical network}

\index{operators!semibounded-}

Large networks arise in both pure and applied mathematics, e.g., in
graph theory (the mathematical theory of networks), and more recently,
they have become a current and fast developing research area; with
applications including a host of problems coming from for example
internet search, and social networks. Hence, of the recent applications,
there is a change in outlook from finite to infinite.

More precisely, in traditional graph theoretical problems, the whole
graph is given exactly, and we are then looking for relationships
between its parameters, variables and functions; or for efficient
algorithms for computing them. By contrast, for very large networks
(like the Internet), variables are typically not known completely;
-- in most cases they may not even be well defined. In such applications,
data about them can only be collected by indirect means; hence random
variables and local sampling must be used as opposed to global processes.\index{random variable}\index{graph!network-}

Although such modern applications go far beyond the setting of large
electrical networks (even the case of infinite sets of vertices and
edges), it is nonetheless true that the framework of large electrical
networks is helpful as a basis for the analysis we develop below;
and so our results will be phrased in the setting of large electrical
networks, even though the framework is much more general.

The applications of \textquotedblleft large\textquotedblright{} or
infinite graphs are extensive, counting just physics; see for example
\cite{BCD06,RAKK05,KMRS05,BC05,TD03,VZ92}.

In discrete harmonic analysis, two operations play a key role, the
Laplacian $\Delta$, and the Markov operator $P$. An infinite network
is a pair of sets, $V$ vertices, and $E$, edges. In addition to
this, one specifies a conductance function c . This is a function
c defined on the edge set $E$. There are then two associated operators
$\Delta$ and $P$ are defined from, and they depend on the entire
triple $(V,E,c)$. For many problems one of the two operators is even
used in the derivation of properties of the other. Both represent
actions (operations) on appropriate spaces of functions, i.e., functions
defined on the infinite set of vertices $V$. For the networks of
interest to us, the vertex set $V$ will be infinite, and we are therefore
faced with a variety of choices of infinite-dimensional function spaces.
Because of spectral theory, the useful choices will be Hilbert spaces.\index{Laplacian!graph-}

But even restricting to Hilbert spaces, there are at least three natural
candidates: (i) the plain $l^{2}$ sequence space, so an $l^{2}$-space
of functions on $V$, (ii) a suitably weighted $l^{2}$-space, and
finally (iii), an energy Hilbert space $\mathscr{H}_{E}$. (The latter
is an abstraction of more classical notions of Dirichlet spaces.)
Which one of the three to use depends on the particular operator considered,
and also on the questions asked.

We note that in infinite network models, both the Laplacian $\Delta$,
and the Markov operator $P$ will have infinite by infinite matrix
representations. Each of these infinite by infinite matrices is special,
in that, as an infinite by infinite matrix, it will have non-zero
entries localized only in finite bands containing the infinite matrix-diagonal
(i.e., they are infinite banded matrices). This makes the algebraic
matrix operations well defined.

Functional analytic and spectral theoretic tools enter the analysis
as follows: In passing to appropriate Hilbert spaces, we arrive at
classes of Hilbert space-operators, and the operators in question
will be Hermitian. But the Laplacian $\Delta$ will be typically be
an unbounded operator, albeit semibounded. By contrast we show that
there is a weighted $l^{2}$-space such that the corresponding Markov
operator $P$ is a bounded, selfadjoint operator, and that it has
its spectrum contained in the finite interval $[-1,1]$. We caution,
that in general this spectrum may be continuous, or have a mix of
spectral types, continuous (singular or Lebesgue), and discrete.

\section{Basic Setting}

Let $V$ be a countable discrete set, and let $E\subset V\times V$
be a subset such that:
\begin{enumerate}
\item $\left(x,y\right)\in E\Longleftrightarrow\left(y,x\right)\in E$;
$x,y\in V$;
\item $\#\left\{ y\in V\:|\:\left(x,y\right)\in E\right\} $ is finite,
and $>0$ for all $x\in V$; 
\item $\left(x,x\right)\notin E$; and
\item $\exists\,o\in V$ such that for all $y\in V$ $\exists\,x_{0},x_{1},\ldots,x_{n}\in V$
with $x_{0}=o$, $x_{n}=y$, $\left(x_{i-1},x_{i}\right)\in E$, $\forall i=1,\ldots,n$.
(This property is called connectedness.)
\item If a conductance function $c$ is given we require $c_{x_{i-1}x_{i}}>0$.
See \defref{cond} below.\end{enumerate}
\begin{defn}
\label{def:cond}A function $c:E\rightarrow\mathbb{R}_{+}\cup\left\{ 0\right\} $
is called \emph{conductance function} if \index{conductance}
\begin{enumerate}
\item $c\left(e\right)\geq0$, $\forall e\in E$; and 
\item Given $x\in V$, $c_{xy}>0$, $c_{xy}=c_{yx}$, for all $\left(xy\right)\in E$.
\end{enumerate}

If $x\in V$, we set 
\begin{equation}
c\left(x\right):=\sum_{\left(xy\right)\in E}c_{xy}.\label{eq:cond}
\end{equation}
The summation in (\ref{eq:cond}) is denoted $x\sim y$; i.e., $x\sim y$
if $\left(xy\right)\in E$. 

\end{defn}

\begin{defn}
\label{def:lap}When $c$ is a conductance function (see also \defref{cond})
we set $\Delta=\Delta_{c}$ (the corresponding graph Laplacian \index{graph Laplacian}
\index{operators!Laplace-}\index{Laplacian!graph-}
\begin{equation}
\left(\Delta u\right)\left(x\right)=\sum_{y\sim x}c_{xy}\left(u\left(x\right)-u\left(y\right)\right)=c\left(x\right)u\left(x\right)-\sum_{y\sim x}c_{xy}u\left(y\right).\label{eq:lap}
\end{equation}

\end{defn}
\label{rem:lapl2}Given $G=\left(V,E,c\right)$ as above, and let
$\Delta=\Delta_{c}$ be the corresponding graph Laplacian. With a
suitable ordering on $V$, we obtain the following banded $\infty\times\infty$
matrix-representation for $\Delta$ (eq. (\ref{eq:lapm})). We refer
to \cite{GLS12} for a number of applications of infinite banded matrices.
\index{matrix!inftytimesinfty
@$\infty\times\infty$}\index{matrix!banded-}\index{graph!network-} 
\begin{equation}
\begin{bmatrix}c\left(x_{1}\right) & -c_{x_{1}x_{2}} & 0 & \cdots & \cdots & \cdots & \cdots & 0 & \cdots\\
-c_{x_{2}x_{1}} & c\left(x_{2}\right) & -c_{x_{2}x_{3}} & 0 & \cdots & \cdots & \cdots & \vdots & \cdots\\
0 & -c_{x_{3}x_{2}} & c\left(x_{3}\right) & -c_{x_{3}x_{4}} & 0 & \cdots & \cdots & \huge\mbox{0} & \cdots\\
\vdots & 0 & \ddots & \ddots & \ddots & \ddots & \vdots & \vdots & \cdots\\
\vdots & \vdots & \ddots & \ddots & \ddots & \ddots & 0 & \vdots & \cdots\\
\vdots & \huge\mbox{0} & \cdots & 0 & -c_{x_{n}x_{n-1}} & c\left(x_{n}\right) & -c_{x_{n}x_{n+1}} & 0 & \cdots\\
\vdots & \vdots & \cdots & \cdots & 0 & \ddots & \ddots & \ddots & \ddots
\end{bmatrix}\label{eq:lapm}
\end{equation}

\index{matrix!banded-}
\begin{rem}[\textbf{Random Walk}]
\label{rem:rw} If $\left(V,E,c\right)$ is given as in \defref{lap},
then for $\left(x,y\right)\in E$, set \index{random walk} 
\begin{equation}
p_{xy}:=\frac{c_{xy}}{c\left(x\right)}\label{eq:pxy}
\end{equation}
and note then $\left\{ p_{xy}\right\} $ in (\ref{eq:pxy}) is a system
of transition probabilities\index{transition probabilities}, i.e.,
$\sum_{y}p_{xy}=1$, $\forall x\in V$, see  \figref{tp}.

\begin{figure}
\includegraphics[scale=0.5]{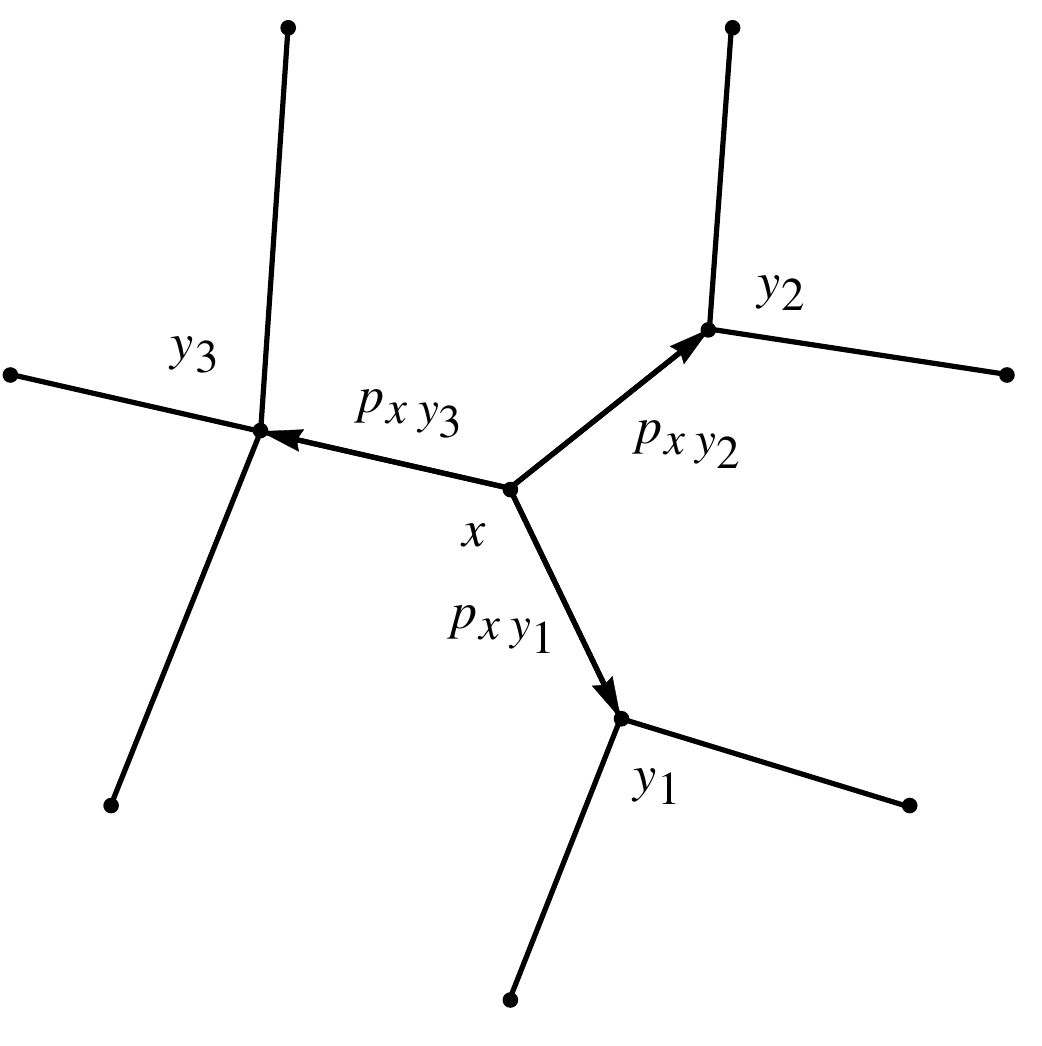}

\protect\caption{\label{fig:tp}Transition probabilities $p_{xy}$ at a vertex $x$
$\left(\mbox{in }V\right)$. }
\end{figure}

\end{rem}
A Markov-random walk on $V$ with transition probabilities $\left(p_{xy}\right)$
is said to be \emph{\uline{reversible}} iff $\exists$ a positive
function $\widetilde{c}$ on $V$ such that \index{reversible}
\begin{equation}
\widetilde{c}\left(x\right)p_{xy}=\widetilde{c}\left(y\right)p_{yx},\;\forall\left(xy\right)\in E.\label{eq:mar1}
\end{equation}
 
\begin{thm}[\cite{Jor08,MR2542093,MR2711706}]
Let $G=\left(V,E,c\right)$ be as above, $V$: vertices, $E$: edges,
and $c:E\longrightarrow\mathbb{R}_{+}$ a given conductance function;
we assume finite range so that the $\infty\times\infty$ matrix in
(\ref{eq:lapm}) is banded.

Then the banded $\infty\times\infty$ matrix in (\ref{eq:lapm}) defines
an essentially selfadjoint operator in $l^{2}\left(V\right)$, with
dense domain equal to all finitely supported functions on $V$.
\end{thm}
\index{essentially selfadjoint operator}

\section{\label{sec:enSp}The Energy Hilbert Spaces $\mathscr{H}_{E}$ }

Let $G=\left(V,E,c\right)$ be an infinite connected network introduced
above. Set $\mathscr{H}_{E}:=$ completion of the space of all compactly
supported functions $u:V\rightarrow\mathbb{C}$ with respect to\index{Hilbert space!energy}\index{completion!Hilbert-}
\begin{align}
\left\langle u,v\right\rangle _{\mathscr{H}_{E}} & :=\frac{1}{2}\underset{\left(x,y\right)\in E}{\sum\sum}c_{xy}(\overline{u\left(x\right)}-\overline{u\left(y\right)})\left(v\left(x\right)-v\left(y\right)\right)\label{eq:Einn}\\
\left\Vert u\right\Vert _{\mathscr{H}_{E}}^{2}: & =\frac{1}{2}\underset{\left(x,y\right)\in E}{\sum\sum}c_{xy}\left|u\left(x\right)-u\left(y\right)\right|^{2}\label{eq:Enorm}
\end{align}
then $\mathscr{H}_{E}$ is a Hilbert space \cite{JoPe10,JT14c}.\index{space!Hilbert-}
\begin{lem}
\label{lem:dipole}For all $x,y\in V$, there is a unique real-valued
\uline{dipole} vector $v_{xy}\in\mathscr{H}_{E}$ such that \index{Riesz' theorem}
\[
\left\langle v_{xy},u\right\rangle _{\mathscr{H}_{E}}=u\left(x\right)-u\left(y\right),\;\forall u\in\mathscr{H}_{E}.
\]
\end{lem}
\begin{proof}
Apply Riesz' theorem. \index{Theorem!Riesz-}
\end{proof}
\index{Gaussian process}

\index{Gaussian free field (GFF)}

\begin{figure}[H]
\begin{tabular}[t]{>{\centering}p{0.45\textwidth}>{\centering}p{0.45\textwidth}}
\includegraphics[scale=0.5]{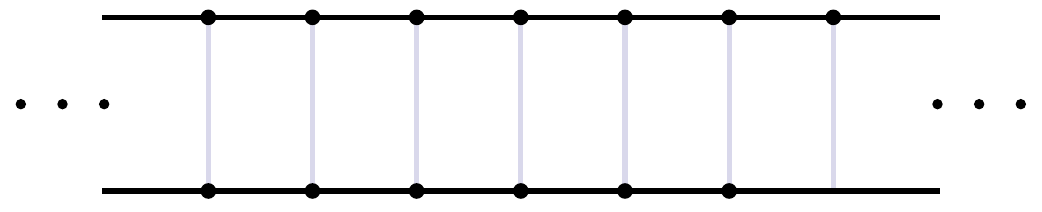}\\
 & \includegraphics[scale=0.45]{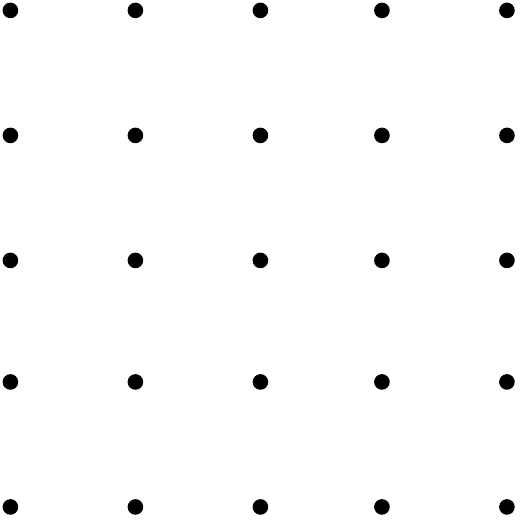}\\
\tabularnewline
$V=\mbox{Band}$ & $V=\mathbb{Z}^{2}$\tabularnewline
\end{tabular}

\protect\caption{\label{fig:frame}non-linear system of vertices}
\end{figure}

\index{dipole}
\begin{xca}[Gaussian free field (GFF) \cite{MR2883397}]
\myexercise{Gaussian free field (GFF)} Let $G=\left(V,E,c\right)$
be as in the setting of \secref{enSp}, and let $\mathscr{H}_{E}$
be the corresponding energy Hilbert space with inner product, and
$\mathscr{H}_{E}$-norm as in (\ref{eq:Einn})-(\ref{eq:Enorm}). 
\begin{enumerate}
\item \label{enu:gf1}Show that there is a probability space $\left(\Omega,\mathcal{F},\mathbb{P}^{\left(G\right)}\right)$
and a Gaussian field $X_{\varphi}$, indexed by $\varphi\in\mathscr{H}_{E}$
(real valued) such that $E_{\mathbb{P}}\left(X_{\varphi}\right)=0$,
and 
\begin{equation}
\mathbb{E}_{\mathbb{P}^{\left(G\right)}}\left(e^{iX_{\varphi}}\right)=e^{-\frac{1}{2}\left\Vert \varphi\right\Vert _{\mathscr{H}_{E}}^{2}};\label{eq:gf1}
\end{equation}
in particular, 
\begin{equation}
\mathbb{E}_{\mathbb{P}^{\left(G\right)}}\left(X_{\varphi}X_{\psi}\right)=\left\langle \varphi,\psi\right\rangle _{\mathscr{H}_{E}}\label{eq:gf2}
\end{equation}
for all $\varphi,\psi\in\mathscr{H}_{E}$.
\item \label{enu:gf2}Show that $X$ arises from a Gaussian point process
$\left\{ X_{x}\right\} _{x\in V}$ such that 
\begin{equation}
\mathbb{E}_{\mathbb{P}^{\left(G\right)}}\left(X_{x}X_{\varphi}\right)=\left\langle v_{xo},\varphi\right\rangle _{\mathscr{H}_{E}}=\varphi\left(x\right)\label{eq:gf3}
\end{equation}
for all $\varphi\in\mathscr{H}_{E}$, and all $x\in V$, where $o$
is a fixed base-point in the vertex set $V$, and we normalize in
(\ref{eq:gf3}) such that $\varphi\left(o\right)=0$.
\end{enumerate}

\uline{Hint}: For (\ref{enu:gf1}), use  \corref{gns}; and for
(\ref{enu:gf2}), use  \lemref{dipole}.

\end{xca}
\begin{defn}
Let $\mathscr{H}$ be a Hilbert space with inner product denoted $\left\langle \cdot,\cdot\right\rangle $,
or $\left\langle \cdot,\cdot\right\rangle _{\mathscr{H}}$ when there
is more than one possibility to consider. Let $J$ be a countable
index set, and let $\left\{ w_{j}\right\} _{j\in J}$ be an indexed
family of non-zero vectors in $\mathscr{H}$. We say that $\left\{ w_{j}\right\} _{j\in J}$
is a \emph{\uline{frame}}\emph{ }for $\mathscr{H}$ iff (Def.)
there are two finite positive constants $b_{1}$ and $b_{2}$ such
that \index{frame}
\begin{equation}
b_{1}\left\Vert u\right\Vert _{\mathscr{H}}^{2}\leq\sum_{j\in J}\left|\left\langle w_{j},u\right\rangle _{\mathscr{H}}\right|^{2}\leq b_{2}\left\Vert u\right\Vert _{\mathscr{H}}^{2}\label{eq:en1}
\end{equation}
holds for all $u\in\mathscr{H}$. We say that it is a \emph{\uline{Parseval}}
frame if $b_{1}=b_{2}=1$. 

For references to the theory and application of \emph{\uline{frames}},
see e.g., \cite{HJL13,KLZ09,CM13,SD13,KOPT13,EO13}.\end{defn}
\begin{lem}
\label{lem:eframe}If $\left\{ w_{j}\right\} _{j\in J}$ is a Parseval
frame in $\mathscr{H}$, then the (analysis) operator $A=A_{\mathscr{H}}:\mathscr{H}\longrightarrow l^{2}\left(J\right)$,
\begin{equation}
Au=\left(\left\langle w_{j},u\right\rangle _{\mathscr{H}}\right)_{j\in J}\label{eq:en2}
\end{equation}
is well-defined and isometric. Its adjoint $A^{*}:l^{2}\left(J\right)\longrightarrow\mathscr{H}$
is given by 
\begin{equation}
A^{*}\left(\left(\gamma_{j}\right)_{j\in J}\right):=\sum_{j\in J}\gamma_{j}w_{j}\label{eq:en3}
\end{equation}
and the following hold:
\begin{enumerate}
\item The sum on the RHS in (\ref{eq:en3}) is norm-convergent;
\item $A^{*}:l^{2}\left(J\right)\longrightarrow\mathscr{H}$ is co-isometric;
and for all $u\in\mathscr{H}$, we have 
\begin{equation}
u=A^{*}Au=\sum_{j\in J}\left\langle w_{j},u\right\rangle w_{j}\label{eq:en4}
\end{equation}
where the RHS in (\ref{eq:en4}) is norm-convergent. 
\end{enumerate}
\end{lem}
\begin{proof}
The details are standard in the theory of frames; see the cited papers
above. Note that (\ref{eq:en1}) for $b_{1}=b_{2}=1$ simply states
that $A$ in (\ref{eq:en2}) is isometric, and so $A^{*}A=I_{\mathscr{H}}=$
the identity operator in $\mathscr{H}$, and $AA^{*}=$ the projection
onto the range of $A$.\end{proof}
\begin{thm}
\label{thm:eframe}Let $G=\left(V,E,c\right)$ be an infinite network.
Choose an orientation on the edges, denoted by $E^{\left(ori\right)}$.
Then the system of vectors
\begin{equation}
\left\{ w_{xy}:=\sqrt{c_{xy}}v_{xy},\;\ensuremath{\left(xy\right)\in E^{\left(ori\right)}}\right\} \label{eq:en8}
\end{equation}
is a Parseval frame for the energy Hilbert space $\mathscr{H}_{E}$.
For all $u\in\mathscr{H}_{E}$, we have the following representation
\begin{align}
u & =\sum_{\left(xy\right)\in E^{\left(ori\right)}}c_{xy}\left\langle v_{xy},u\right\rangle v_{xy},\;\mbox{and}\label{eq:frep}\\
\left\Vert u\right\Vert _{\mathscr{H}_{E}}^{2} & =\sum_{\left(xy\right)\in E^{\left(ori\right)}}c_{xy}\left|\left\langle v_{xy},u\right\rangle \right|^{2}\label{eq:fnorm}
\end{align}
\end{thm}
\begin{proof}
See \cite{JT14,CH08}. \index{frame}
\end{proof}
Frames in $\mathscr{H}_{E}$ consisting of our system (\ref{eq:en8})
are not ONBs when resisters are configured in non-linear systems of
vertices, for example, resisters in parallel. See \figref{frame},
and \exref{tri}. \index{infinite network}

\begin{figure}[H]
\includegraphics[scale=0.4]{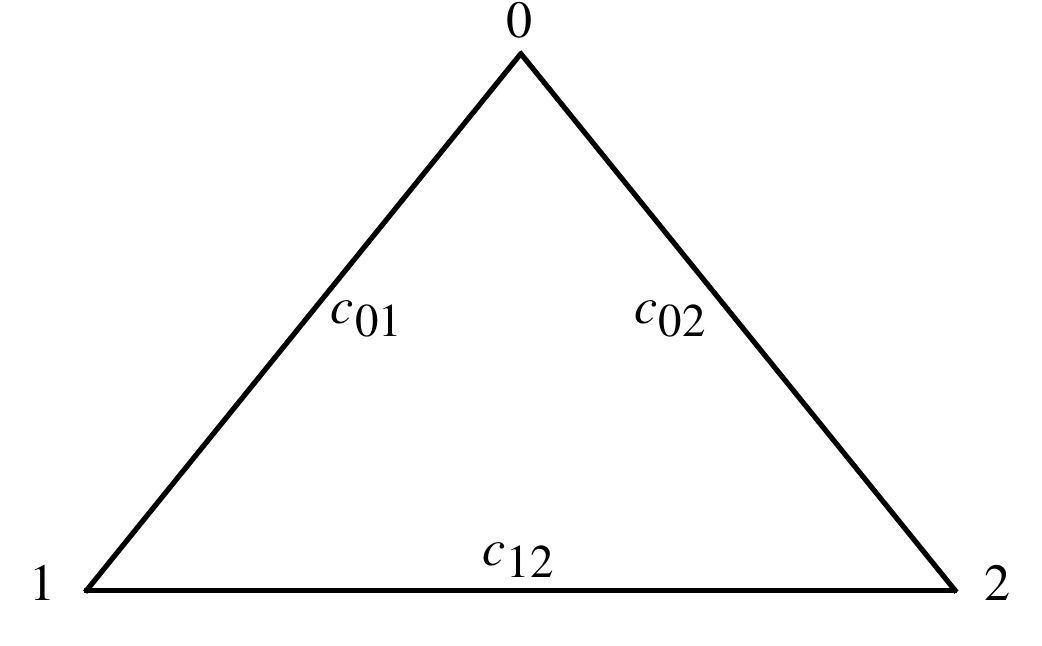}

\protect\caption{\label{fig:tri}The set $\left\{ v_{xy}:\left(xy\right)\in E\right\} $
is not orthogonal.}
\end{figure}

\begin{example}
\label{ex:tri}Let $c_{01},c_{02},c_{12}$ be positive constants,
and assign conductances on the three edges (see \figref{tri}) in
the triangle network. \index{conductance}
\end{example}

In this case, $w_{ij}=\sqrt{e_{ij}}v_{ij}$, $i<j$, in the cyclic
order is a Parseval frame but not an ONB in $\mathscr{H}_{E}$ \cite{JT14}.
\index{frame}

Note the corresponding Laplacian $\Delta\left(=\Delta_{c}\right)$
has the following matrix representation 
\begin{equation}
M:=\begin{bmatrix}c\left(0\right) & -c_{01} & -c_{02}\\
-c_{01} & c\left(1\right) & -c_{12}\\
-c_{02} & -c_{12} & c\left(2\right)
\end{bmatrix}\label{eq:LM}
\end{equation}
The dipoles $\left\{ v_{xy}:\left(xy\right)\in E^{\left(ori\right)}\right\} $
as 3-D vectors are the solutions to the equation \index{dipole} 
\[
\Delta v_{xy}=\delta_{x}-\delta_{y}.
\]
Hence, 
\begin{align*}
Mv_{01} & =\begin{bmatrix}1 & -1 & 0\end{bmatrix}^{tr}\\
Mv_{02} & =\begin{bmatrix}1 & 0 & -1\end{bmatrix}^{tr}\\
Mv_{12} & =\begin{bmatrix}0 & 1 & -1\end{bmatrix}^{tr}
\end{align*}
The Parseval frame from \lemref{eframe} is 
\begin{align*}
w_{01} & =\sqrt{c_{01}}v_{01}=\left[\frac{\sqrt{c_{01}}\,c_{12}}{c_{01}c_{02}+c_{01}c_{12}+c_{02}c_{12}},-\frac{\sqrt{c_{01}}\,c_{02}}{c_{01}c_{02}+c_{01}c_{12}+c_{02}c_{12}},0\right]^{tr}\\
w_{12} & =\sqrt{c_{12}}v_{12}=\left[0,\frac{\sqrt{c_{12}}\,c_{02}}{c_{01}c_{02}+c_{01}c_{12}+c_{02}c_{12}},-\frac{\sqrt{c_{12}}\,c_{01}}{c_{01}c_{02}+c_{01}c_{12}+c_{02}c_{12}}\right]^{tr}\\
w_{20} & =\sqrt{c_{20}}v_{20}=\left[\frac{-\sqrt{c_{20}}\,c_{12}}{c_{01}c_{02}+c_{01}c_{12}+c_{02}c_{12}},0,\frac{\sqrt{c_{20}}\,c_{01}}{c_{01}c_{02}+c_{01}c_{12}+c_{02}c_{12}}\right]^{tr}.
\end{align*}

\begin{rem}
The dipole $v_{xy}$ is unique in $\mathscr{H}_{E}$ as an equivalence
class, not a function on $V$. Note $\ker M$ = harmonic functions
= constant (see (\ref{eq:LM})), and so $v_{xy}+\mbox{const}=v_{xy}$
in $\mathscr{H}_{E}$. Thus, the above frame vectors have non-unique
representations as functions on $V$. \index{harmonic}
\end{rem}

\section{The Graph-Laplacian\index{Laplacian!graph-}}

Here we include some technical lemmas for graph Laplacian in the energy
Hilbert space $\mathscr{H}_{E}$ . \index{graph Laplacian} \index{operators!Laplace-}

Let $G=\left(V,E,c\right)$ be as above; assume $G$ is connected;
i.e., there is a base point $o$ in $V$ such that every $x\in V$
is connected to $o$ via a finite path of edges. 

If $x\in V$, we set 
\begin{equation}
\delta_{x}\left(y\right)=\begin{cases}
1 & \mbox{if }y=x\\
0 & \mbox{if }y\neq x
\end{cases}\label{eq:delx}
\end{equation}

\begin{defn}
\label{def:D}Let $\left(V,E,c,o,\Delta\right)$ be as above. Let
$V':=V\backslash\left\{ o\right\} $, and set
\[
v_{x}:=v_{x,o},\;\forall x\in V'.
\]
Further, let 
\begin{align}
\mathscr{D}_{2} & :=span\left\{ \delta_{x}\:\big|\:x\in V\right\} ,\;\mbox{and}\label{eq:D1}\\
\mathscr{D}_{E} & :=\left\{ \sum\nolimits _{x\in V'}\xi_{x}v_{x}\:\big|\:\mbox{finite support}\right\} ;\label{eq:D2}
\end{align}
where by ``span'' we mean of all \emph{finite} linear combinations.
\end{defn}
\lemref{Delta} below summarizes the key properties of $\Delta$ as
an operator, both in $l^{2}(V)$ and in $\mathscr{H}_{E}$. 
\begin{lem}
\label{lem:Delta}The following hold:
\begin{enumerate}
\item $\left\langle \Delta u,v\right\rangle _{l^{2}}=\left\langle u,\Delta v\right\rangle _{l^{2}}$,
$\forall u,v\in\mathscr{D}_{2}$;
\item \label{enu:D2}$\left\langle \Delta u,v\right\rangle _{\mathscr{H}_{E}}=\left\langle u,\Delta v\right\rangle _{\mathscr{H}_{E}},$
$\forall u,v\in\mathscr{D}_{E}$;
\item \label{enu:D3}$\left\langle u,\Delta u\right\rangle _{l^{2}}\geq0$,
$\forall u\in\mathscr{D}_{2}$, and
\item $\left\langle u,\Delta u\right\rangle _{\mathscr{H}_{E}}\geq0$, $\forall u\in\mathscr{D}_{E}$.
\end{enumerate}

Moreover, we have
\begin{enumerate}[resume]
\item \label{enu:D5}$\left\langle \delta_{x},u\right\rangle _{\mathscr{H}_{E}}=\left(\Delta u\right)\left(x\right)$,
$\forall x\in V$, $\forall u\in\mathscr{H}_{E}$.
\item $\Delta v_{xy}=\delta_{x}-\delta_{y}$, $\forall v_{xy}\in\mathscr{H}_{E}$.
In particular, $\Delta v_{x}=\delta_{x}-\delta_{o}$, $x\in V'=V\backslash\left\{ o\right\} $. 
\item \label{enu:D7}
\[
\delta_{x}\left(\cdot\right)=c\left(x\right)v_{x}\left(\cdot\right)-\sum_{y\sim x}c_{xy}v_{y}\left(\cdot\right),\;\forall x\in V'.
\]

\item 
\[
\left\langle \delta_{x},\delta_{y}\right\rangle _{\mathscr{H}_{E}}=\begin{cases}
c\left(x\right)=\sum_{t\sim x}c_{xt} & \mbox{if \ensuremath{y=x}}\\
-c_{xy} & \mbox{if \ensuremath{\left(xy\right)\in E}}\\
0 & \mbox{if }\left(xy\right)\notin E,\;x\neq y
\end{cases}
\]

\end{enumerate}
\end{lem}
\begin{proof}
See \cite{JoPe10,JoPe11a,JT14}.
\end{proof}

\section{The Friedrichs Extension}

Fix a conductance function $c$. In this section we turn to some technical
lemmas we will need for the Friedrichs extension of $\Delta\left(=\Delta_{c}\right)$.\index{Friedrichs extension}\index{extension!Friedrich's-}

It is known the graph-Laplacian $\Delta$ is automatically essentially
selfadjoint as a densely defined operator in $l^{2}(V)$, but not
as a $\mathscr{H}_{E}$ operator \cite{Jor08,JoPe11b}. Since $\Delta$
defined on $\mathscr{D}_{E}$ is semibounded, it has the Friedrichs
extension $\Delta_{Fri}$ (in $\mathscr{H}_{E}$).\index{operators!semibounded-}

\index{essentially selfadjoint operator}
\begin{lem}
\label{lem:FriedDomain}Consider $\Delta$ with $dom\left(\Delta\right):=span\left\{ v_{xy}:x,y\in V\right\} $,
then 
\[
\left\langle \varphi,\Delta\varphi\right\rangle _{\mathscr{H}_{E}}=\sum_{\left(xy\right)\in E}c_{xy}^{2}\left|\left\langle v_{xy},\varphi\right\rangle _{\mathscr{H}_{E}}\right|^{2}.
\]
\end{lem}
\begin{proof}
Suppose $\varphi=\sum\varphi_{xy}v_{xy}\in dom(\Delta)$. Note the
edges are not oriented, and a direct computation shows that
\[
\left\langle \varphi,\Delta\varphi\right\rangle _{\mathscr{H}_{E}}=4\sum_{x,y}\left|\varphi_{xy}\right|^{2}.
\]

Using the Parseval frames in \thmref{eframe}, we have the following
representation 
\[
\varphi=\sum_{\left(xy\right)\in E}\underset{=:\varphi_{xy}}{\underbrace{\frac{1}{2}c_{xy}\left\langle v_{xy},\varphi\right\rangle _{\mathscr{H}_{E}}}}v_{xy}
\]
Note $\varphi\in span\left\{ v_{xy}:x,y\in V\right\} $, so the above
equation contains a finite sum. 

It follows that 
\[
\left\langle \varphi,\Delta\varphi\right\rangle _{\mathscr{H}_{E}}=4\sum_{\left(xy\right)\in E}\left|\varphi_{xy}\right|^{2}=\sum_{\left(xy\right)\in E}c_{xy}^{2}\left|\left\langle v_{xy},\varphi\right\rangle _{\mathscr{H}_{E}}\right|^{2}
\]
which is the assertion.\end{proof}
\begin{thm}
\label{thm:domF}Let $G=\left(V,E,c\right)$ be an infinite network.
If the deficiency indices of $\Delta\left(=\Delta_{c}\right)$ are
$\left(k,k\right)$, $k>0$, where $dom(\Delta)=span\left\{ v_{xy}\right\} $,
then the Friedrichs extension $\Delta_{Fri}\supset\Delta$ is the
restriction of $\Delta^{*}$ to 
\begin{equation}
dom(\Delta_{Fri}):=\left\{ u\in\mathscr{H}_{E}\:\big|\:\sum\nolimits _{\left(xy\right)\in E}c_{xy}^{2}\left|\left\langle v_{xy},\varphi\right\rangle _{E}\right|^{2}<\infty\right\} .\label{eq:FrieDomain}
\end{equation}
 \end{thm}
\begin{proof}
Follows from \lemref{FriedDomain}, and the characterization of Friedrichs
extensions of semibounded Hermitian operators (\chapref{ext}); see,
e.g., \cite{DS88b,AG93,RS75}. \index{index!deficiency}
\end{proof}

\section{\label{ex:1d}A 1D Example}

Consider $G=\left(V,E,c\right)$, where $V=\left\{ 0\right\} \cup\mathbb{Z}_{+}$.
Observation: Every sequence $a_{1},a_{2},\ldots$ in $\mathbb{R}_{+}$
defines a conductance $c_{n-1,n}:=a_{n}$, $n\in\mathbb{Z}_{+}$,
i.e., 
\[
\xymatrix{0\ar@{<->}[r]_{a_{1}} & 1\ar@{<->}[r]_{a_{2}} & 2\ar@{<->}[r]_{a_{3}} & 3 & \cdots & n\ar@{<->}[r]_{a_{n+1}} & n+1}
\cdots
\]

\begin{figure}[H]
\includegraphics[scale=0.5]{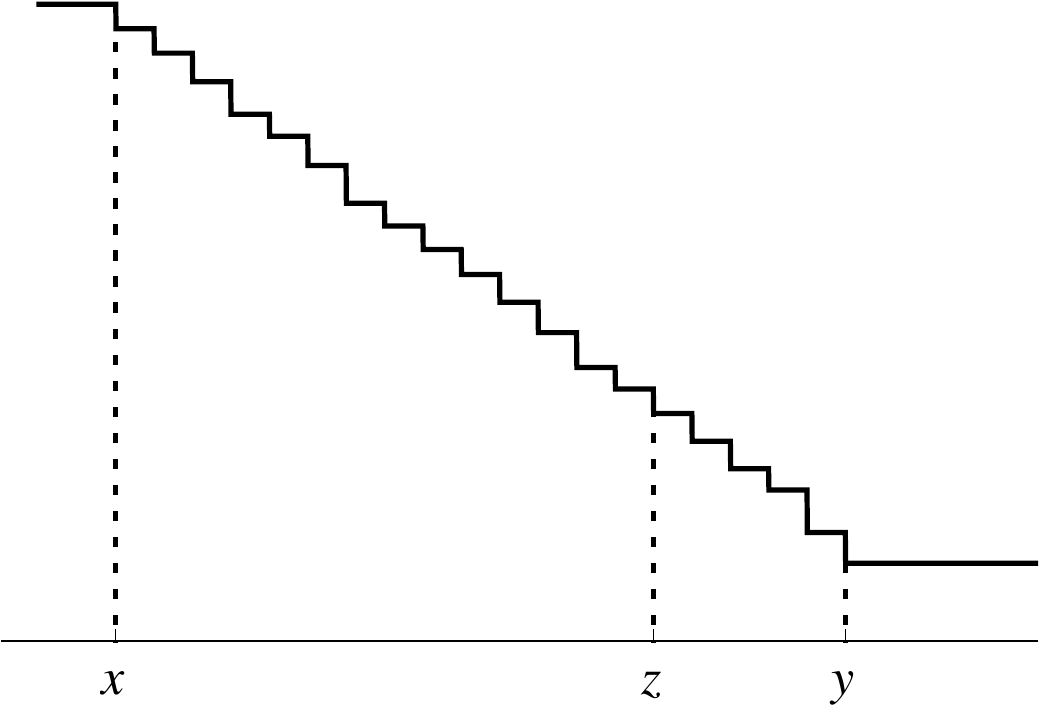}

\protect\caption{\label{fig:1dipole}The dipole $v_{xy}$.}
\end{figure}

The dipole vectors $v_{xy}$ (for $x,y\in\mathbb{N}$) are given by
\[
v_{xy}\left(z\right)=\begin{cases}
0 & \mbox{if \ensuremath{z\leq x}}\\
-\sum_{k=x+1}^{z}\frac{1}{a_{k}} & \mbox{if \ensuremath{x<z<y}}\\
-\sum_{k=x+1}^{y}\frac{1}{a_{k}} & \mbox{if \ensuremath{z\geq y}}
\end{cases}
\]
See \figref{1dipole}. \index{conductance}\index{dipole}

The corresponding graph Laplacian has the following matrix representation:
\index{Laplacian!graph-}\index{graph!network-}

\begin{equation}
\begin{bmatrix}a_{1} & -a_{1}\\
-a_{1} & a_{1}+a_{2} & -a_{2}\\
 & -a_{2} & a_{2}+a_{3} & -a_{3} &  &  & \Huge\mbox{0}\\
 &  & -a_{3} & a_{3}+a_{4} & \ddots\\
 &  &  & \ddots & \ddots & -a_{n}\\
 &  &  &  & -a_{n} & a_{n}+a_{n+1} & -a_{n+1}\\
 & \Huge\mbox{0} &  &  &  & -a_{n+1} & \ddots & \ddots\\
 &  &  &  &  &  & \ddots & \ddots\\
 &  &  &  &  &  & \ddots & \ddots
\end{bmatrix}\label{eq:Lma1}
\end{equation}
That is, \index{matrix!tri-diagonal-} 
\begin{equation}
\begin{cases}
\left(\Delta u\right)_{0} & =a_{1}\left(u_{0}-u_{1}\right)\\
\left(\Delta u\right)_{n} & =a_{n}\left(u_{n}-u_{n-1}\right)+a_{n+1}\left(u_{n}-u_{n+1}\right)\\
 & =\left(a_{n}+a_{n+1}\right)u_{n}-a_{n}u_{n-1}-a_{n+1}u_{n+1},\;\forall n\in\mathbb{Z}_{+}.
\end{cases}\label{eq:Le1}
\end{equation}

\begin{lem}
Let $G=\left(V,c,E\right)$ be as above, where $a_{n}:=c_{n-1,n}$,
$n\in\mathbb{Z}_{+}$. Then $u\in\mathscr{H}_{E}$ is the solution
to $\Delta u=-u$ (i.e., $u$ is a defect vector of $\Delta$) if
and only if $u$ satisfies the following equation:
\begin{equation}
\sum_{n=1}^{\infty}a_{n}\left\langle v_{n-1,n},u\right\rangle _{\mathscr{H}_{E}}\left(\delta_{n-1}\left(s\right)-\delta_{n}\left(s\right)+v_{n-1,n}\left(s\right)\right)=0,\;\forall s\in\mathbb{Z}_{+};\label{eq:1ddef}
\end{equation}
where 
\begin{equation}
\left\Vert u\right\Vert _{\mathscr{H}_{E}}^{2}=\sum_{n=1}^{\infty}a_{n}\left|\left\langle v_{n-1,n},u\right\rangle _{\mathscr{H}_{E}}\right|^{2}<\infty.\label{eq:1ddef-1}
\end{equation}
\end{lem}
\begin{proof}
By  \thmref{eframe}, the set $\left\{ \sqrt{a_{n}}v_{n-1,n}\right\} _{n=1}^{\infty}$
forms a Parseval frame in $\mathscr{H}_{E}$. In fact, the dipole
vectors are \index{frame} 
\begin{equation}
v_{n-1,n}\left(s\right)=\begin{cases}
0 & s\leq n-1\\
-\frac{1}{a_{n}} & s\geq n
\end{cases};n=1,2,\ldots\label{eq:2}
\end{equation}
and so $\left\{ \sqrt{a_{n}}v_{n-1,n}\right\} _{n=1}^{\infty}$ forms
an ONB in $\mathscr{H}_{E}$; and $u\in\mathscr{H}_{E}$ has the representation\index{dipole}
\[
u=\sum_{n=1}^{\infty}a_{n}\left\langle v_{n-1,n},u\right\rangle _{\mathscr{H}_{E}}v_{n-1,n}
\]
see (\ref{eq:en4}). Therefore, $\Delta u=-u$ if and only if 
\begin{align*}
\sum_{n=1}^{\infty}a_{n}\left\langle v_{n-1,n},u\right\rangle _{\mathscr{H}_{E}}\left(\delta_{n-1}\left(s\right)-\delta_{n}\left(s\right)\right) & =-\sum_{n=1}^{\infty}a_{n}\left\langle v_{n-1,n},u\right\rangle _{\mathscr{H}_{E}}v_{n-1,n}\left(s\right)
\end{align*}
for all $s\in\mathbb{Z}_{+}$, which is the assertion.
\end{proof}
Below we compute the deficiency space in an example with index values
$\left(1,1\right)$.\index{deficiency space}
\begin{lem}
\label{lem:Q11}Let $\left(V,E,c=\left\{ a_{n}\right\} \right)$ be
as above. Let $Q>1$ and set $a_{n}:=Q^{n}$, $n\in\mathbb{Z}_{+}$;
then $\Delta$ has deficiency indices $\left(1,1\right)$.\end{lem}
\begin{proof}
Suppose $\Delta u=-u$, $u\in\mathscr{H}_{E}$. Then, 
\begin{align*}
-u_{1} & =Q\left(u_{1}-u_{0}\right)+Q^{2}\left(u_{1}-u_{2}\right)\Longleftrightarrow u_{2}=\left(\frac{1}{Q^{2}}+\frac{1+Q}{Q}\right)u_{1}-\frac{1}{Q}u_{0}\\
-u_{2} & =Q^{2}\left(u_{2}-u_{1}\right)+Q^{3}\left(u_{2}-u_{3}\right)\Longleftrightarrow u_{3}=\left(\frac{1}{Q^{3}}+\frac{1+Q}{Q}\right)u_{2}-\frac{1}{Q}u_{1}
\end{align*}
and by induction, 
\[
u_{n+1}=\left(\frac{1}{Q^{n+1}}+\frac{1+Q}{Q}\right)u_{n}-\frac{1}{Q}u_{n-1},\;n\in\mathbb{Z}_{+}
\]
i.e., $u$ is determined by the following matrix equation:

\[
\begin{bmatrix}u_{n+1}\\
u_{n}
\end{bmatrix}=\begin{bmatrix}\frac{1}{Q^{n+1}}+\frac{1+Q}{Q} & -\frac{1}{Q}\\
1 & 0
\end{bmatrix}\begin{bmatrix}u_{n}\\
u_{n-1}
\end{bmatrix}
\]

The eigenvalues of the coefficient matrix are 
\begin{align*}
\lambda_{\pm} & =\frac{1}{2}\left(\frac{1}{Q^{n+1}}+\frac{1+Q}{Q}\pm\sqrt{\left(\frac{1}{Q^{n+1}}+\frac{1+Q}{Q}\right)^{2}-\frac{4}{Q}}\right)\\
 & \sim\frac{1}{2}\left(\frac{1+Q}{Q}\pm\left(\frac{Q-1}{Q}\right)\right)=\begin{cases}
1\\
\dfrac{1}{Q}
\end{cases}\mbox{as \ensuremath{n\rightarrow\infty}.}
\end{align*}
Equivalently, as $n\rightarrow\infty$, we have
\[
u_{n+1}\sim\left(\frac{1+Q}{Q}\right)u_{n}-\frac{1}{Q}u_{n-1}=\left(1+\frac{1}{Q}\right)u_{n}-\frac{1}{Q}u_{n-1}
\]
and so 
\[
u_{n+1}-u_{n}\sim\frac{1}{Q}\left(u_{n}-u_{n-1}\right).
\]
Therefore, for the tail-summation, we have:
\[
\sum_{n}Q^{n}\left(u_{n+1}-u_{n}\right)^{2}=\mbox{const}\sum_{n}\frac{\left(Q-1\right)^{2}}{Q^{n+2}}<\infty
\]
which implies $\left\Vert u\right\Vert _{\mathscr{H}_{E}}<\infty$. 
\end{proof}
Next, we give a random walk interpretation of \lemref{Q11}. See \remref{rw},
and \figref{tp}. \index{random walk}
\begin{rem}[Harmonic functions in $\mathscr{H}_{E}$]
 \label{rem:Qharm}Note that in \exref{1d} (\lemref{Q11}), the
space of harmonic functions in $\mathscr{H}_{E}$ is one-dimensional;
in fact if $Q>1$ is fixed, then 
\[
\left\{ u\in\mathscr{H}_{E}\:\big|\:\Delta u=0\right\} 
\]
is spanned by $u=\left(u_{n}\right)_{n=0}^{\infty}$, $u_{n}=\frac{1}{Q^{n}}$,
$n\in\mathbb{N}$; and of course $\Vert1/Q^{n}\Vert_{\mathscr{H}_{E}}^{2}<\infty$.
\end{rem}
\index{harmonic}

\begin{figure}[H]
\includegraphics[scale=0.6]{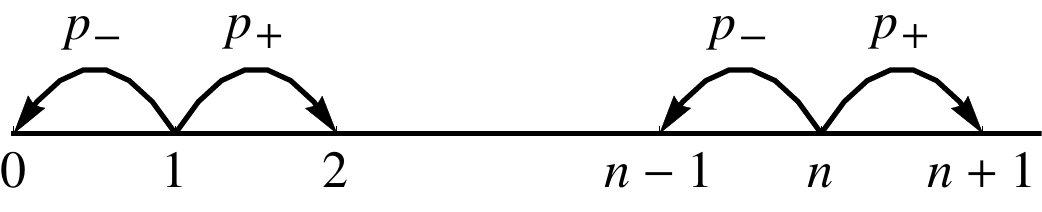}

\protect\caption{\label{fig:tp1}The transition probabilities $p_{+},p_{-}$, in the
case of constant transition probabilities, i.e., $p_{+}\left(n\right)=p_{+}$,
and $p_{-}\left(n\right)=p_{-}$ for all $n\in\mathbb{Z}_{+}$.}
\end{figure}

\begin{rem}
\label{rem:domFri}For the domain of the Friedrichs extension\index{Friedrichs extension}
$\Delta_{Fri}$, we have:
\begin{equation}
dom(\Delta_{Fri})=\left\{ f\in\mathscr{H}_{E}\:|\:\left(f\left(x\right)-f\left(x+1\right)\right)Q^{x}\in l^{2}\left(\mathbb{Z}_{+}\right)\right\} \label{eq:1d-2-1}
\end{equation}
i.e., 
\[
dom(\Delta_{Fri})=\left\{ f\in\mathscr{H}_{E}\:|\:\sum_{x=0}^{\infty}\left|f\left(x\right)-f\left(x+1\right)\right|^{2}Q^{2x}<\infty\right\} .
\]
\end{rem}
\begin{proof}
By \thmref{eframe}, we have the following representation, valid for
all $f\in\mathscr{H}_{E}$: 
\begin{align*}
f & =\sum_{x}\left\langle f,Q^{\frac{x}{2}}v_{\left(x,x+1\right)}\right\rangle _{\mathscr{H}_{E}}Q^{\frac{x}{2}}v_{\left(x,x+1\right)}\\
 & =\sum_{x}\left(f\left(x\right)-f\left(x+1\right)\right)Q^{x}v_{\left(x,x+1\right)};
\end{align*}
and
\[
\left\langle f,\Delta f\right\rangle _{\mathscr{H}_{E}}=\sum_{x}\left|f\left(x\right)-f\left(x+1\right)\right|^{2}Q^{2x}.
\]
The desired conclusion (\ref{eq:1d-2-1}) now follows from \thmref{domF}.
Also see e.g. \cite{DS88b,AG93}.\end{proof}
\begin{defn}
Let $G=\left(V,E,c\right)$ be a connected graph. The set of transition
probabilities\index{transition probabilities} $\left(p_{xy}\right)$
is said to be reversible if there exists $c:V\rightarrow\mathbb{R}_{+}$
such that 
\begin{equation}
c\left(x\right)p_{xy}=c\left(y\right)p_{yx};\label{eq:tp1}
\end{equation}
and then \index{reversible} 
\begin{equation}
c_{xy}:=c\left(x\right)p_{xy}\label{eq:tpc}
\end{equation}
is a system of conductance. Conversely, for a system of conductance
$\left(c_{xy}\right)$ we set 
\begin{align}
c\left(x\right) & :=\sum_{y\sim x}c_{xy},\;\mbox{and}\label{eq:cond1}\\
p_{xy} & :=\frac{c_{xy}}{c\left(x\right)}\label{eq:tp-1}
\end{align}
and so $\left(p_{xy}\right)$ is a set of transition probabilities.
See \figref{tp2}. \index{conductance}

\begin{figure}
\begin{tabular}{cc}
\includegraphics[scale=0.45]{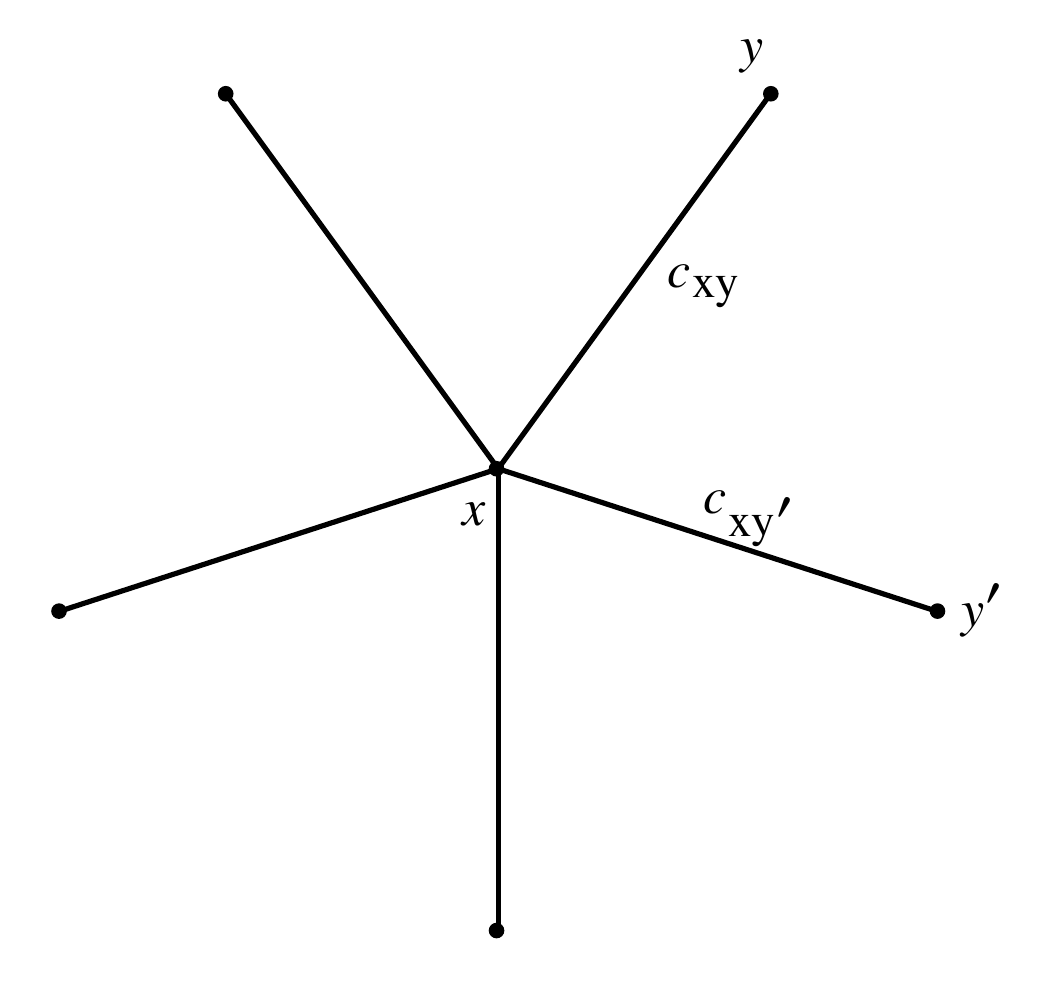} & \includegraphics[scale=0.45]{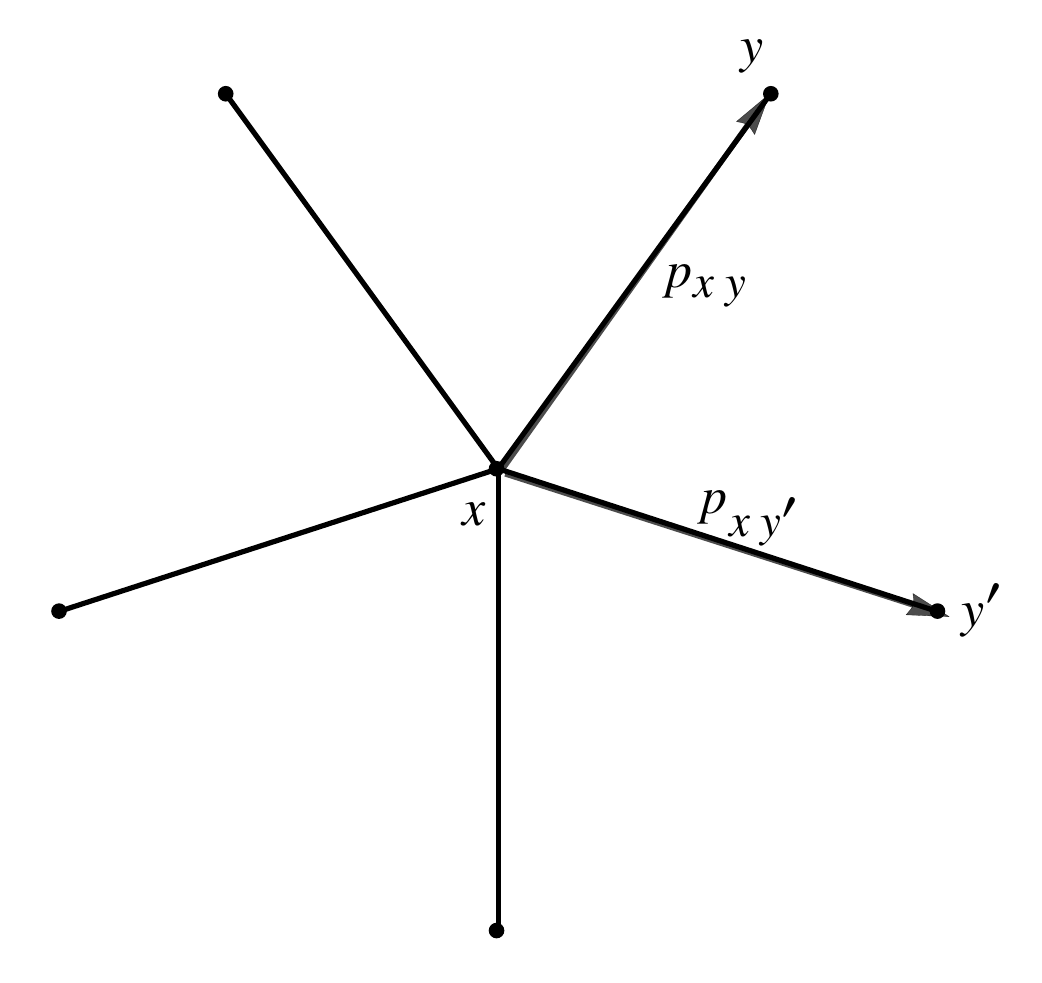}\tabularnewline
$c_{xy}$, $y\sim x$ & transition probabilities\tabularnewline
\end{tabular}

\protect\caption{\label{fig:tp2}neighbors of $x$}
\end{figure}

\end{defn}
Recall the graph Laplacian in (\ref{eq:Le1}) can be written as 
\begin{equation}
\left(\Delta u\right)_{n}=c\left(n\right)\left(u_{n}-p_{-}\left(n\right)u_{n-1}-p_{+}\left(n\right)u_{n+1}\right),\;\forall n\in\mathbb{Z}_{+};\label{eq:Le2}
\end{equation}
where 
\begin{equation}
c\left(n\right):=a_{n}+a_{n+1}\label{eq:Lec}
\end{equation}
and 
\begin{equation}
p_{-}\left(n\right):=\frac{a_{n}}{c\left(n\right)},\;p_{+}\left(n\right):=\frac{a_{n+1}}{c\left(n\right)}\label{eq:Lep}
\end{equation}
are the left/right transition probabilities\index{transition probabilities},
as shown in \figref{tp1}.

\index{operators!Laplace-}

In the case $a_{n}=Q^{n}$, $Q>1$, as in \lemref{Q11}, we have 
\begin{equation}
c\left(n\right):=Q^{n}+Q^{n+1},\;\mbox{and }\label{eq:cn}
\end{equation}
\begin{align}
p_{+} & :=p_{+}\left(n\right)=\frac{Q^{n+1}}{Q^{n}+Q^{n+1}}=\frac{Q}{1+Q}\label{eq:pplus}\\
p_{-} & :=p_{-}\left(n\right)=\frac{Q^{n}}{Q^{n}+Q^{n+1}}=\frac{1}{1+Q}\label{eq:pmimus}
\end{align}
For all $n\in\mathbb{Z}_{+}\cup\left\{ 0\right\} $, set 
\begin{equation}
\left(Pu\right)_{n}:=p_{-}u_{n-1}+p_{+}u_{n+1}.\label{eq:cp}
\end{equation}
Note $\left(Pu\right)_{0}=u_{1}$. By (\ref{eq:Le2}), we have 
\begin{equation}
\Delta=c\left(1-P\right).\label{eq:dp}
\end{equation}

In particular, $p_{+}>\frac{1}{2}$, i.e., a random walker has probability
$>\frac{1}{2}$ of moving to the right. It follows that 
\[
\underset{=\mbox{ dist to \ensuremath{\infty}}}{\underbrace{\mbox{travel time}\left(n,\infty\right)}}<\infty;
\]
and so $\Delta$ is not essentially selfadjoint, i.e., indices $\left(1,1\right)$.

\index{operators!essentially selfadjoint-}
\begin{lem}
Let $\left(V,E,\Delta(=\Delta_{c})\right)$ be as above, where the
conductance $c$ is given by $c_{n-1,n}=Q^{n}$, $n\in\mathbb{Z}_{+}$,
$Q>1$ (see \lemref{Q11}). For all $\lambda>0$, there exists $f_{\lambda}\in\mathscr{H}_{E}$
satisfying $\Delta f_{\lambda}=\lambda f_{\lambda}$.\end{lem}
\begin{proof}
By (\ref{eq:dp}), we have $\Delta f_{\lambda}=\lambda f_{\lambda}\Longleftrightarrow Pf_{\lambda}=\left(1-\frac{\lambda}{c}\right)f_{\lambda}$,
i.e., 
\[
\frac{1}{1+Q}f_{\lambda}\left(n-1\right)+\frac{Q}{1+Q}f_{\lambda}\left(n+1\right)=\left(1-\frac{\lambda}{Q^{n-1}\left(1+Q\right)}\right)f_{\lambda}\left(n\right)
\]
and so 
\begin{equation}
f_{\lambda}\left(n+1\right)=\left(\frac{1+Q}{Q}-\frac{\lambda}{Q^{n}}\right)f_{\lambda}\left(n\right)-\frac{1}{Q}f_{\lambda}\left(n-1\right).\label{eq:Fe0}
\end{equation}
This corresponds to the following matrix equation:
\begin{align*}
\begin{bmatrix}f\left(n+1\right)\\
f\left(n\right)
\end{bmatrix} & =\begin{bmatrix}\frac{1+Q}{Q}-\frac{\lambda}{Q^{n}} & -\frac{1}{Q}\\
1 & 0
\end{bmatrix}\begin{bmatrix}f\left(n\right)\\
f\left(n-1\right)
\end{bmatrix}\\
 & \sim\begin{bmatrix}\frac{1+Q}{Q} & -\frac{1}{Q}\\
1 & 0
\end{bmatrix}\begin{bmatrix}f\left(n\right)\\
f\left(n-1\right)
\end{bmatrix},\;\mbox{as \ensuremath{n\rightarrow\infty}.}
\end{align*}
The eigenvalues of the coefficient matrix are given by\index{eigenvalue}
\[
\lambda_{\pm}\sim\frac{1}{2}\left(\frac{1+Q}{Q}\pm\left(\frac{Q-1}{Q}\right)\right)=\begin{cases}
1\\
\dfrac{1}{Q}
\end{cases}\mbox{as \ensuremath{n\rightarrow\infty}.}
\]
That is, as $n\rightarrow\infty$, 
\[
f_{\lambda}\left(n+1\right)\sim\left(\frac{1+Q}{Q}\right)f_{\lambda}\left(n\right)-\frac{1}{Q}f_{\lambda}\left(n-1\right);
\]
i.e.,
\begin{equation}
f_{\lambda}\left(n+1\right)\sim\frac{1}{Q}f_{\lambda}\left(n\right);\label{eq:Fe1}
\end{equation}
and so the tail summation of $\left\Vert f_{\lambda}\right\Vert _{\mathscr{H}_{E}}^{2}$
is finite. (See the proof of \lemref{Q11}.) We conclude that $f_{\lambda}\in\mathscr{H}_{E}$.\end{proof}
\begin{cor}
Let $\left(V,E,\Delta\right)$ be as in the lemma. The Friedrichs
extension $\Delta_{Fri}$ has continuous spectrum $[0,\infty)$. \index{spectrum!continuous}\end{cor}
\begin{proof}
Fix $\lambda\geq0$. We prove that if $\Delta f_{\lambda}=\lambda f_{\lambda}$,
$f\in\mathscr{H}_{E}$, then $f_{\lambda}\notin dom(\Delta_{Fri})$. 

Note for $\lambda=0$, $f_{0}$ is harmonic\index{harmonic}, and
so $f_{0}=k\left(\frac{1}{Q^{n}}\right)_{n=0}^{\infty}$ for some
constant $k\neq0$. See \remref{Qharm}. It follows from (\ref{eq:1d-2-1})
that $f_{0}\notin dom(\Delta_{Fri})$. 

The argument for $\lambda>0$ is similar. Since as $n\rightarrow\infty$,
$f_{\lambda}\left(n\right)\sim\frac{1}{Q^{n}}$ (eq. (\ref{eq:Fe1})),
so by (\ref{eq:1d-2-1}) again, $f_{\lambda}\notin dom(\Delta_{Fri})$.

However, if $\lambda_{0}<\lambda_{1}$ in $[0,\infty)$ then 
\begin{equation}
\int_{\lambda_{0}}^{\lambda_{1}}f_{\lambda}\left(\cdot\right)d\lambda\in dom(\Delta_{Fri})\label{eq:Fe2}
\end{equation}
and so every $f_{\lambda}$, $\lambda\in[0,\infty)$, is a generalized
eigenfunction\index{eigenfunction}\index{generalized eigenfunction},
i.e., the spectrum of $\Delta_{Fri}$ is purely continuous with Lebesgue\index{measure!Lebesgue}
measure, and multiplicity one. 

The verification of (\ref{eq:Fe2}) follows from (\ref{eq:Fe0}),
i.e., 
\begin{equation}
f_{\lambda}\left(n+1\right)=\left(\frac{1+Q}{Q}-\frac{\lambda}{Q^{n}}\right)f_{\lambda}\left(n\right)-\frac{1}{Q}f_{\lambda}\left(n-1\right).\label{eq:Fe3}
\end{equation}
Set 
\begin{equation}
F_{\left[\lambda_{0},\lambda_{1}\right]}:=\int_{\lambda_{0}}^{\lambda_{1}}f_{\lambda}\left(\cdot\right)d\lambda.\label{eq:Fe4}
\end{equation}
Then by (\ref{eq:Fe3}) and (\ref{eq:Fe4}),
\[
F_{\left[\lambda_{0},\lambda_{1}\right]}\left(n+1\right)=\frac{1+Q}{Q}F_{\left[\lambda_{0},\lambda_{1}\right]}\left(n\right)-\frac{1}{Q^{n}}\int_{\lambda_{0}}^{\lambda_{1}}\lambda f_{\lambda}\left(n\right)d\lambda-\frac{1}{Q}F_{\left[\lambda_{0},\lambda_{1}\right]}\left(n-1\right)
\]
and $\int_{\lambda_{0}}^{\lambda_{1}}\lambda f_{\lambda}d\lambda$
is computed using integration by parts.
\end{proof}

\section*{A summary of relevant numbers from the Reference List}

For readers wishing to follow up sources, or to go in more depth with
topics above, we suggest: \cite{Jor08,JoPe10,JoPe11a,JoPe11b,JoPe13,JoPe13b,JT14,RAKK05,Yos95,MR3042410,MR2921088,JP14,LPW13,JP14,MR2354721,AJSV13,MR2883397,JP12,JT14c}.

\chapter{Reproducing Kernel Hilbert Space\label{chap:RKHS}}
\begin{quotation}
The simplicities of natural laws arise through the complexities of
the language we use for their expression. 

--- Eugene Wigner\sindex[nam]{Wigner, E.P., (1902-1995)}\vspace{1em}\\
... an apt comment on how science, and indeed the whole of civilization,
is a series of incremental advances, each building on what went before. 

--- Stephen Hawking\sindex[nam]{Hawking, S., (1942--)}\vspace{2em}\\
At a given moment there is only a fine layer between the \textquoteleft trivial\textquoteright{}
and the impossible. Mathematical discoveries are made in this layer.

--- Andrey Kolmogorov \vspace{2em}\\
\dots{} the distinction between two sorts of truths, profound truths
recognized by the fact that the opposite is also a profound truth;
by contrast to trivialities, where opposites are obviously absurd.

--- Niels Bohr
\end{quotation}
A special family of Hilbert spaces $\mathscr{H}$ are reproducing
kernel Hilbert spaces (RKHSs). We say that $\mathscr{H}$ is a RKHS
if $\mathscr{H}$ is a Hilbert space of functions on some set $X$
such that for every $x$ in $X$, the linear mapping $f\longmapsto f\left(x\right)$
is continuous in the norm of $\mathscr{H}$. \index{kernel!reproducing-}

We begin with a general 
\begin{defn}
\label{def:RKHS}Let $S$ be a set. We say that $\mathscr{H}$ is
a $S$-reproducing kernel Hilbert space if:
\begin{enumerate}
\item \label{enu:RKHS1}$\mathscr{H}$ is a Hilbert space of functions on
$S$; and
\item \label{enu:RKHS2}For all $s\in S$, the mapping $\mathscr{H}\longrightarrow\mathbb{C}$,
by $E_{s}:h\longmapsto h\left(s\right)$ is continuous on $\mathscr{H}$,
i.e., by Riesz, there is a $K_{s}\in\mathscr{H}$ such that $h\left(s\right)=\left\langle K_{s},h\right\rangle $,
$\forall h\in\mathscr{H}$.
\end{enumerate}

\uline{Notation}. The system of functions $\left\{ K_{s}:s\in S\right\} \subset\mathscr{H}$
is called the associated reproducing kernel. 

\end{defn}
Let $\mathscr{H}$ be a RKHS, see \defref{RKHS}, hence $S$ is a
set, and $\mathscr{H}$ is a Hilbert space of functions on $S$ such
that (\ref{enu:RKHS2}) holds. Fix $s\in S$, and note that $E_{s}:\mathscr{H}\longrightarrow\mathbb{C}$
is then a bounded linear operator, where $\mathbb{C}$ is a 1-dimensional
Hilbert space. Hence its adjoint $E_{s}^{*}:\mathbb{C}\longrightarrow\mathscr{H}^{*}\simeq\mathscr{H}$
is well-defined.

\index{reproducing kernel Hilbert space (RKHS)}
\begin{claim}
The kernel $K_{s}$ in  \defref{RKHS} is $E_{s}^{*}\left(1\right)=K_{s}$. \end{claim}
\begin{proof}
For $\lambda\in\mathbb{C}$ and $h\in\mathscr{H}$, we have 
\[
\left\langle E_{s}^{*}\left(\lambda\right),h\right\rangle _{\mathscr{H}}=\overline{\lambda}E_{s}\left(h\right)=\overline{\lambda}h\left(s\right)\underset{\left(\text{by}\;\left(\ref{enu:RKHS2}\right)\right)}{=}\overline{\lambda}\left\langle K_{s},h\right\rangle _{\mathscr{H}}=\left\langle \lambda K_{s},h\right\rangle _{\mathscr{H}};
\]
and therefore $E_{s}^{*}\left(\lambda\right)=\lambda K_{s}$ as desired.\end{proof}
\begin{example}
For $s,t\in\left[0,1\right]$, set 
\begin{equation}
K\left(s,t\right)=s\wedge t=\min\left(s,t\right),\label{eq:rk1-1}
\end{equation}
i.e., the covariance-kernel for Brownian motion on the interval $\left[0,1\right]$,
see Chapter \ref{chap:bm}. \index{positive definite!-kernel}\end{example}
\begin{xca}[The Fundamental Theorem of Calculus and a RKHS]
\myexercise{The Fundamental Theorem of Calculus and a RKHS}\label{exer:rkhs1}Show
that the RKHS for the kernel $K$ in (\ref{eq:rk1-1}) is
\[
\mathscr{H}=\begin{Bmatrix}f\::\: & \mbox{locally integrable with distribution-}\\
 & \mbox{derivative}\;f'\in L^{2}\left(0,1\right)\;\mbox{and}\;f\left(0\right)=0
\end{Bmatrix}
\]
with the inner product 
\[
\left\langle f,g\right\rangle _{\mathscr{H}}:=\int_{0}^{1}\overline{f'\left(x\right)}g'\left(x\right)dx.
\]
\uline{Hint}: 

\uline{Step 1} Show that $K_{t}=K\left(t,\cdot\right)$ is in $\mathscr{H}$

\uline{Step 2} Show that for all $f\in\mathscr{H}$, we have 
\[
\left\langle K_{t},f\right\rangle _{\mathscr{H}}=\int_{0}^{t}f'\left(x\right)dx=f\left(t\right).
\]

\end{xca}

Suppose $K:\left[0,1\right]\times\left[0,1\right]\longrightarrow\mathbb{C}$
is positive definite \uline{and continuous}, i.e., for all finite
sums:
\begin{equation}
\sum_{j}\sum_{k}\overline{c_{j}}c_{k}K\left(t_{j},t_{k}\right)\geq0,\label{eq:rkh1}
\end{equation}
$\left\{ c_{j}\right\} \subset\mathbb{C}$, $\left\{ t_{j}\right\} \subset\left[0,1\right]$.
\begin{xca}[Mercer \cite{MR1356475,MR0380303}]
\label{exer:mer}\myexercise{Mercer}Show that the operator $T_{K}$,
\begin{equation}
\left(T_{K}\varphi\right)\left(t\right)=\int_{0}^{1}K\left(t,s\right)\varphi\left(s\right)ds\label{eq:rkh2}
\end{equation}
in $L^{2}\left(0,1\right)$ is trace-class, and \index{Theorem!Mercer's-}
\begin{equation}
trace\left(T_{K}\right)=\int_{0}^{1}K\left(t,t\right)dt.\label{eq:rkh3}
\end{equation}

\uline{Hint}: Apply weak-compactness, and the Spectral Theorem
for compact selfadjoint operators,  \thmref{spcpt}. Combine this
with a choice of an ONB in the RKHS defined from (\ref{eq:rkh1}).
\end{xca}
\index{Mercer's Theorem}

\index{selfadjoint operator}

\begin{xca}[The Szegö-kernel]
\label{exer:zego}\myexercise{Szeg\"{o}-kernel}Let $\mathbb{H}_{2}$
be the Hardy space of the disk $\mathbb{D}=\left\{ z\in\mathbb{C}\::\:\left|z\right|<1\right\} $,
see \secref{egCalgebras}. Show that $\mathbb{H}_{2}$ is a RKHS with
reproducing kernel (the Szegö-kernel): \index{kernel!reproducing-}\index{space!Hardy-}
\begin{equation}
K\left(z,w\right)=\frac{1}{1-\overline{z}w}\label{eq:zk1}
\end{equation}
i.e., that we have: \index{positive definite!-kernel}\index{Theorem!Spectral-}
\begin{equation}
\left\langle K\left(z,\cdot\right),f\right\rangle _{\mathbb{H}_{2}}=f\left(z\right),\;\forall f\in\mathbb{H}_{2},\:\forall z\in\mathbb{D}.\label{eq:zk2}
\end{equation}

\uline{Hint}: Substitute (\ref{eq:zk1}) into the formula from
$\left\langle \cdot,\cdot\right\rangle _{\mathbb{H}_{2}}$ inner product
on the LHS in (\ref{eq:zk2}), and recall our convention: Inner products
are linear in the second variable.
\end{xca}
\index{multiplier}
\begin{defn}
Let $X$ be a set and $K:X\times X\rightarrow\mathbb{C}$ a fixed
positive definite function; and let $\mathscr{H}_{K}$ be the corresponding
RKHS; (see \defref{RKHS}). Let $\varphi:X\rightarrow\mathbb{C}$
be a function on $X$; then we say $\varphi$ is a \emph{multiplier},
written ``$\varphi\in\mbox{Multp}(\mathscr{H}_{K})$'' if multiplication
by $\varphi$ defines a bounded linear operator in $\mathscr{H}_{K}$,
so
\begin{equation}
\begin{matrix}M_{\varphi}:\mathscr{H}_{K}\longrightarrow\mathscr{H}_{K}\\
\left(M_{\varphi}f\right)\left(x\right)=\varphi\left(x\right)f\left(x\right),\;\forall f\in\mathscr{H}_{K},\:\forall x\in X.
\end{matrix}\label{eq:ml1}
\end{equation}
\end{defn}
\begin{xca}[Multp ($\mathscr{H}_{K}$)]
\myexercise{Multp. ($\mathscr{H}_{K}$)} ~
\begin{enumerate}
\item \label{enu:ml1}Show the following equivalence:
\begin{eqnarray}
\varphi & \in & \mbox{Multp(\ensuremath{\mathscr{H}_{K}})}\label{eq:ml1-1}\\
 & \Updownarrow\nonumber 
\end{eqnarray}
$\exists$ constant $B<\infty$ such that we have the following estimate
for all finite sums:
\begin{equation}
\sum_{i}\sum_{j}\overline{c_{i}}c_{j}\left(B-\overline{\varphi\left(x_{i}\right)}\varphi\left(x_{j}\right)\right)K\left(x_{i},x_{j}\right);\label{eq:ml1-2}
\end{equation}
i.e., computed for all systems $\{c_{i}\}\subset\mathbb{C}$, and
$\{x_{i}\}\subset X$. 
\item \label{enu:ml2}Let $\varphi\in\mbox{Multp}(\mathscr{H}_{K})$, and
let $\{K_{x}\}_{x\in X}$ be the kernel functions. Prove that
\begin{equation}
M_{\varphi}^{*}(K_{x})=\overline{\varphi\left(x\right)}K_{x},\;\forall x\in X,\label{eq:ml2-1}
\end{equation}
where $M_{\varphi}^{*}$ denotes the adjoint operator.

\uline{Hint}: Verify the following identity:
\begin{equation}
\left\langle \overline{\varphi\left(x\right)}K_{x},f\right\rangle _{\mathscr{H}_{K}}=\left\langle K_{x},M_{\varphi}f\right\rangle _{\mathscr{H}_{K}},\;\forall x\in X,\:f\in\mathscr{H}_{K}.\label{eq:ml2-2}
\end{equation}

\item Show directly from (\ref{enu:ml1}) and (\ref{enu:ml2}) that $\mbox{Multp}(\mathscr{H}_{K})$
is an algebra.
\item Apply (\ref{enu:ml1}) to the Hardy space $\mathbb{H}_{2}$ of the
disk to conclude that 
\begin{equation}
\mbox{Multp}(\mathbb{H}_{2})=\mathbb{H}_{\infty}.\label{eq:ml1-3}
\end{equation}

\end{enumerate}
\end{xca}
\uline{Contents of the Chapter.}

In this chapter, we study two extension problems, and their interconnections.
The first class of extension problems concerns (i) positive definite
(p.d.) continuous functions on Lie groups $G$, and the second deals
with (ii) Lie algebras of unbounded skew-Hermitian operators in a
certain family of reproducing kernel Hilbert spaces (RKHS). The analysis
is non-trivial even if $G=\mathbb{R}^{n}$, and even if $n=1$. If
$G=\mathbb{R}^{n}$, we are concerned in (ii) with the study of systems
of $n$ skew-Hermitian operators $\left\{ S_{i}\right\} $ on a common
dense domain in Hilbert space, and in deciding whether it is possible
to find a corresponding system of strongly commuting selfadjoint operators
$\left\{ T_{i}\right\} $ such that, for each value of $i$, the operator
$T_{i}$ extends $S_{i}$.\index{Hilbert space!reproducing kernel}
\index{positive definite!-function}\index{space!Hardy-}

From the postulates of quantum physics, we know that measurements
of observable\index{observable}s are computed from associated selfadjoint
operators---observables. From the corresponding spectral resolutions,
we get probability measures, and of course uncertainty. There are
many philosophical issues (which we bypass here), and we do not yet
fully understand quantum reality. See for example, \cite{Sla03,CJK+12}.
\index{selfadjoint operator}

The axioms are as follows: An \emph{observable} is a Hermitian (selfadjoint)
linear operator mapping a Hilbert space, the space of \emph{states},
into itself. The values obtained in a \emph{physical measurement}\index{measurement}
are, in general, described by a probability distribution; and the
distribution represents a suitable \textquotedblleft average\textquotedblright{}
(or \textquotedblleft \emph{expectation}\textquotedblright ) in a
measurement of values of some quantum observable in a state of some
prepared system. The states are (up to phase) unit vectors in the
Hilbert space, and a measurement corresponds to a \emph{probability
distribution} (derived from a projection-valued spectral measure).
The spectral type may be \emph{continuous} (such as position and momentum)
or \emph{discrete} (such as spin).

Information about the measures $\mu$ are computed with the use of
\emph{generating functions} (on $\mathbb{R}$), i.e., spectral (Bochner/Fourier)
transforms of the corresponding measure. Generating functions are
positive definite continuous functions $F\left(=F_{\mu}\right)$ on
$\mathbb{R}$. One then tries to recover $\mu$ from information about
$F$. In this chapter we explore the cases when information about
$F\left(x\right)$ is only available for $x$ in a bounded interval.

\index{generating function}

\index{spectrum!continuous}

\index{spectrum!discrete}

\index{operators!momentum-}

In probability theory, normalized continuous positive definite functions
$F$, i.e., $F(0)=1$, arise as \emph{generating functions} for probability
measures, and one passes from information about one to the other;
-- from generating function to probability measure is called \textquotedblleft the
inverse problem\textquotedblright , see e.g., \cite{DM85}. Hence
the study of partially defined p.d. functions addresses the inverse
question: ambiguity of measures when only partial information for
a possible generating function is available. \index{positive definite!-function}

\section{\label{sec:stoch}A Digression: Stochastic Processes}

Below we continue the discussion of stochastic processes started in
 \secref{Hilbert}.

The interest in positive definite functions has at least three roots:
(i) Fourier analysis, and harmonic\index{harmonic} analysis more
generally, including the non-commutative variant where we study unitary
representations of groups; (ii) optimization and approximation problems,
involving for example spline approximations as envisioned by I. Schöenberg;
and (iii) the study of stochastic (random) processes. \index{stochastic process}

A stochastic process is an indexed family of random variables based
on a fixed probability space; in our present analysis, the processes
will be indexed by some group $G$; for example $G=\mathbb{R}$, or
$G=\mathbb{Z}$ correspond to processes indexed by real time, respectively
discrete time. A main tool in the analysis of stochastic processes
is an associated covariance function, see (\ref{eq:stat1}).

A process $\left\{ X_{g}\:\big|\:g\in G\right\} $ is called Gaussian
if each random variable $X_{g}$ is Gaussian, i.e., its distribution
is Gaussian. For Gaussian processes we only need two moments. So if
we normalize, setting the mean equal to $0$, then the process is
determined by the covariance function. In general the covariance function
is a function on $G\times G$, or on a subset, but if the process
is stationary, the covariance function will in fact be a positive
definite function defined on $G$, or a subset of $G$. We will be
using three stochastic processes in this book, Brownian motion, Brownian
Bridge, and the Ornstein-Uhlenbeck process, all Gaussian, or It\={o}
integrals. \index{It=o, K.@It\=o, K.!integral} \index{distribution!Gaussian-}
\index{integral!It-@It\=o-} \index{Ornstein-Uhlenbeck process}

We outline a brief sketch of these facts below.\index{random variable}

Let $G$ be a locally compact group, and let $\left(\Omega,\mathscr{F},\mathbb{P}\right)$
be a probability space, $\mathscr{F}$ a sigma-algebra\index{sigma-algebra},
and $\mathbb{P}$ a probability\index{measure!probability} measure
defined on $\mathscr{F}$. A stochastic $L^{2}$-process is a system
of random variables $\left\{ X_{g}\right\} _{g\in G}$, $X_{g}\in L^{2}\left(\Omega,\mathscr{F},\mathbb{P}\right)$.
The covariance function $c_{X}$ of the process is the function $G\times G\rightarrow\mathbb{C}$
given by 
\begin{equation}
c_{X}\left(g_{1},g_{2}\right)=\mathbb{E}\left(\overline{X}_{g_{1}}X_{g_{2}}\right),\;\forall\left(g_{1},g_{2}\right)\in G\times G.\label{eq:stat1}
\end{equation}
To simplify will assume that the mean $\mathbb{E}\left(X_{g}\right)=\int_{\Omega}X_{g}d\mathbb{P}\left(\omega\right)=0$
for all $g\in G$. 

We say that $\left(X_{g}\right)$ is stationary iff 
\begin{equation}
c_{X}\left(hg_{1},hg_{2}\right)=c_{X}\left(g_{1},g_{2}\right),\;\forall h\in G.\label{eq:stat2}
\end{equation}
In this case $c_{X}$ is a function of $g_{1}^{-1}g_{2}$, i.e., 
\begin{equation}
\mathbb{E}\left(X_{g_{1},}X_{g_{2}}\right)=c_{X}\left(g_{1}^{-1}g_{2}\right),\;\forall g_{1},g_{2}\in G.\label{eq:stat3}
\end{equation}
(Just take $h=g_{1}^{-1}$ in (\ref{eq:stat2}).)

We now recall the following theorem of Kolmogorov (see \cite{PaSc75}).
One direction is easy, and the other is the deep part:
\begin{defn}
A function $c$ defined on a subset of $G$ is said to be \emph{positive
definite} iff
\[
\sum_{i}\sum_{j}\overline{\lambda_{i}}\lambda_{j}c\left(g_{i}^{-1}g_{j}\right)\geq0
\]
for all finite summation, where $\lambda_{i}\in\mathbb{C}$ and $g_{i}^{-1}g_{j}$
in the domain of $c$. \end{defn}
\begin{thm}[Kolmogorov]
 A function $c:G\rightarrow\mathbb{C}$ is positive definite if and
only if there is a stationary Gaussian process $\left(\Omega,\mathscr{F},\mathbb{P},X\right)$
with mean zero, such that $c=c_{X}$. \end{thm}
\begin{proof}
To stress the idea, we include the easy part of the theorem, and we
refer to \cite{PaSc75} for the non-trivial direction: 

Let $\lambda_{1},\lambda_{2},\ldots,\lambda_{n}\in\mathbb{C}$, and
$\left\{ g_{i}\right\} _{i=1}^{N}\subset G$, then for all finite
summations, we have:
\[
\sum_{i}\sum_{j}\overline{\lambda_{i}}\lambda_{j}c_{X}\left(g_{i}^{-1}g_{j}\right)=\mathbb{E}\left(\left|\sum_{i=1}^{N}\lambda_{i}X_{g_{i}}\right|^{2}\right)\geq0.
\]

\end{proof}

\section{\label{sec:2ext}Two Extension Problems}

While each of the two extension problems has received a considerable
amount of attention in the literature, our emphasis here will be the
interplay between the two problems: Our aim is a duality theory; and,
in the case $G=\mathbb{R}^{n}$, and $G=\mathbb{T}^{n}=\mathbb{R}^{n}/\mathbb{Z}^{n}$,
we will state our theorems in the language of Fourier duality of abelian
groups: With the time frequency duality formulation of Fourier duality
for $G=\mathbb{R}^{n}$ we have that both the time domain and the
frequency domain constitute a copy of $\mathbb{R}^{n}$. We then arrive
at a setup such that our extension questions (i) are in time domain,
and extensions from (ii) are in frequency domain. Moreover we show
that each of the extensions from (i) has a variant in (ii). Specializing
to $n=1$, we arrive of a spectral theoretic characterization of all
skew-Hermitian operators with dense domain in a separable Hilbert
space, having deficiency-indices $\left(1,1\right)$. \index{index!deficiency}\index{deficiency indices}\index{index!von Neumann-}

A systematic study of densely defined Hermitian operators with deficiency
indices $\left(1,1\right)$, and later $\left(d,d\right)$, was initiated
by M. Krein \cite{Kre46}, and is also part of de Branges\textquoteright{}
model theory; see \cite{deB68,BrRo66}. The direct connection between
this theme and the problem of extending continuous positive definite
(p.d.) functions $F$ when they are only defined on a fixed open subset
to $\mathbb{R}^{n}$ was one of our motivations. One desires continuous
p.d. extensions to $\mathbb{R}^{n}$. 

If $F$ is given, we denote the set of such extensions $Ext\left(F\right)$.
If $n=1$, $Ext\left(F\right)$ is always non-empty, but for $n=2$,
Rudin gave examples in \cite{Ru70,Ru63} when $Ext\left(F\right)$
may be empty. Here we extend these results, and we also cover a number
of classes of positive definite functions on locally compact groups
in general; so cases when $\mathbb{R}^{n}$ is replaced with other
groups, both Abelian and non-abelian.

The results in the framework of locally compact\index{groups!locally compact}
Abelian\index{groups!abelian} groups are more complete than their
counterparts for non-Abelian Lie groups, one reason is the availability
of Bochner\textquoteright s duality theorem for locally compact Abelian
groups; -- not available for non-Abelian Lie groups. \index{Lie!group}
\begin{rem}
Even in one dimension the extension problem for locally defined positive
definite functions is interesting. One reason is that among the Fourier
transforms (generating functions) for finite positive Borel measures
$P$ on $\mathbb{R}$, 
\begin{equation}
g_{P}\left(u\right)=\int_{\mathbb{R}}e^{iux}dP\left(x\right),\;u\in\mathbb{R};\label{eq:lk1}
\end{equation}
one wishes to identify the \emph{infinitely divisible distributions}.
We have:\end{rem}
\begin{thm}[Lévy-Khinchin \cite{MR940528}]
Infinite divisibility holds if and only if $g_{P}$ has the following
representation: $g_{P}\left(u\right)=e^{\eta\left(u\right)}$, such
that for some $a\in\mathbb{R}$, $\sigma\in\mathbb{R}_{+}$, and Borel
measure $L$ on $\mathbb{R}\backslash\left\{ 0\right\} $, we have:
\begin{equation}
\eta\left(u\right)=i\,a\,u-\frac{\sigma^{2}}{2}u^{2}+\int_{\mathbb{R}\backslash\left\{ 0\right\} }\left(e^{i\,u\,x}-1-\frac{i\,u\,x}{1+x^{2}}\right)L\left(dx\right)\label{eq:lk2}
\end{equation}
and the measure $L$ satisfying
\begin{equation}
\int_{\mathbb{R}\backslash\left\{ 0\right\} }\left(1\wedge x^{2}\right)L\left(dx\right)<\infty.\label{eq:lk3}
\end{equation}

\end{thm}
\index{Lévy-Khinchin}

\section{The Reproducing Kernel Hilbert Space $\mathscr{H}_{F}$}

Reproducing kernel Hilbert spaces were pioneered by Aronszajn\index{Aronszajn}
\cite{Aro50}, and subsequently they have been used in a host of applications;
e.g., \cite{Sza04,SZ09,SZ07}. The reproducing kernel property appeared
for the first time in Zaremba's paper \cite{Zar07}.

As for positive definite functions, their use and applications are
extensive and includes such areas as stochastic processes, see e.g.,
\cite{JoPe13,AJSV13,JP12,MR2966130}; harmonic analysis (see \cite{JO00}),
and the references there); potential theory \cite{Fu74b,KL14}; operators
in Hilbert space \cite{Al92,AD86}; and spectral theory \cite{AH13,Nus75,Dev72,Dev59}.
We stress that the literature is vast, and the above list is only
a small sample.

Associated to a pair $\left(\Omega,F\right)$, where $F$ is a prescribed
continuous positive definite function defined on $\Omega$, we outline
a reproducing kernel Hilbert space $\mathscr{H}_{F}$ which will serve
as a key tool in our analysis. The particular RKHSs we need here will
have additional properties (as compared to a general framework); which
allow us to give explicit formulas for our solutions. 

\index{reproducing kernel Hilbert space (RKHS)}

\index{kernel!reproducing-}
\begin{defn}
Let $G$ be a Lie group\index{groups!Lie}. Fix $\Omega\subset G$,
non-empty, open and connected. A continuous function 
\begin{equation}
F:\Omega^{-1}\cdot\Omega\rightarrow\mathbb{C}\label{eq:li-1}
\end{equation}
is \emph{positive definite }(p.d.) if \index{positive definite!-function}
\begin{equation}
\sum_{i}\sum_{j}\overline{c_{i}}c_{j}F\left(x_{i}^{-1}x_{j}\right)\geq0,\label{eq:li-2}
\end{equation}
for all finite systems $\left\{ c_{i}\right\} \subset\mathbb{C}$,
and points $\left\{ x_{i}\right\} \subset\Omega$. 

Equivalently,
\begin{equation}
\int_{\Omega}\int_{\Omega}\overline{\varphi\left(x\right)}\varphi\left(y\right)F\left(x^{-1}y\right)dxdy\geq0,\label{eq:li-3}
\end{equation}
for all $\varphi\in C_{c}\left(\Omega\right)$; where $dx$ denotes
a choice of left-invariant Haar\index{measure!Haar} measure on $G$.
\end{defn}
For simplicity we focus on the case $G=\mathbb{R},$ indicating the
changes needed for general Lie groups.
\begin{defn}
Fix $0<a<\infty$, set $\Omega:=\left(0,a\right)$. Let $F:\Omega-\Omega\rightarrow\mathbb{C}$
be a continuous p.d. function. The \emph{reproducing kernel Hilbert
space (RKHS),} $\mathscr{H}_{F}$, is the completion of the space
of functions \index{completion!Hilbert-}\index{kernel!reproducing-}
\begin{equation}
\sum_{\text{finite}}c_{j}F\left(\cdot-x_{j}\right):c_{j}\in\mathbb{C}\label{eq:H1}
\end{equation}
with respect to the inner product
\[
\left\langle F\left(\cdot-x\right),F\left(\cdot-y\right)\right\rangle _{\mathscr{H}_{F}}=F\left(x-y\right),\;\forall x,y\in\Omega,\;\mbox{and}
\]
\begin{equation}
\big\langle\sum_{i}c_{i}F\left(\cdot-x_{i}\right),\sum_{j}c_{j}F\left(\cdot-x_{j}\right)\big\rangle_{\mathscr{H}_{F}}=\sum_{i}\sum_{j}\overline{c_{i}}c_{j}F\left(x_{i}-x_{j}\right),\label{eq:ip-discrete}
\end{equation}
\end{defn}
\begin{rem}
Throughout, we use the convention that the inner product is conjugate
linear in the first variable, and linear in the second variable. When
more than one inner product is used, subscripts will make reference
to the Hilbert space. 

\textbf{Notation.} Inner product and norms will be denoted $\left\langle \cdot,\cdot\right\rangle $,
and $\left\Vert \cdot\right\Vert $ respectively. Often more than
one inner product is involved, and subscripts are used for identification.\end{rem}
\begin{lem}
\label{lem:RKHS-def-by-integral}The reproducing kernel Hilbert space
(RKHS), $\mathscr{H}_{F}$, is the Hilbert completion of the space
of functions 
\begin{equation}
F_{\varphi}\left(x\right)=\int_{\Omega}\varphi\left(y\right)F\left(x-y\right)dy,\;\forall\varphi\in C_{c}^{\infty}\left(\Omega\right),x\in\Omega\label{eq:H2}
\end{equation}
with respect to the inner product
\begin{equation}
\left\langle F_{\varphi},F_{\psi}\right\rangle _{\mathscr{H}_{F}}=\int_{\Omega}\int_{\Omega}\overline{\varphi\left(x\right)}\psi\left(y\right)F\left(x-y\right)dxdy,\;\forall\varphi,\psi\in C_{c}^{\infty}\left(\Omega\right).\label{eq:hi2}
\end{equation}
In particular, 
\begin{equation}
\left\Vert F_{\varphi}\right\Vert _{\mathscr{H}_{F}}^{2}=\int_{\Omega}\int_{\Omega}\overline{\varphi\left(x\right)}\varphi\left(y\right)F\left(x-y\right)dxdy,\;\forall\varphi\in C_{c}^{\infty}\left(\Omega\right)\label{eq:hn2}
\end{equation}
and 
\begin{equation}
\left\langle F_{\varphi},F_{\psi}\right\rangle _{\mathscr{H}_{F}}=\int_{\Omega}\overline{\varphi\left(x\right)}F_{\psi}\left(x\right)dx,\;\forall\phi,\psi\in C_{c}^{\infty}(\Omega).
\end{equation}
\end{lem}
\begin{proof}
Apply standard approximation, see \lemref{dense-1} below. \index{positive definite!-kernel}
\end{proof}

\begin{figure}[H]
\includegraphics[width=0.5\textwidth]{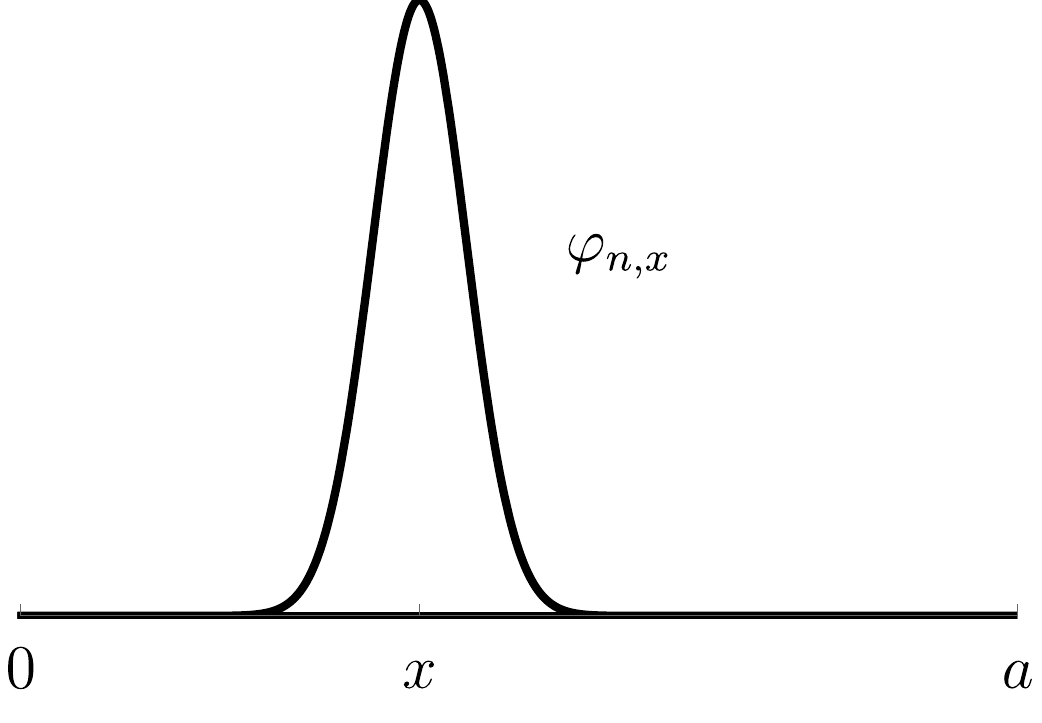}

\protect\caption{\label{fig:approx}The approximate identity $\varphi_{n,x}$}
\end{figure}

The remaining of this section is devoted to a number of technical
lemmas which will be used throughout the chapter. Given a locally
defined continuous positive definite function $F$, the issues addressed
below are: approximation (\lemref{dense-1}), a reproducing kernel
Hilbert space (RKHS) $\mathscr{H}_{F}$ built from $F$, an integral
transform, and a certain derivative operator $D^{\left(F\right)}$,
generally unbounded in the RKHS $\mathscr{H}_{F}$. We will be concerned
with boundary value problems for $D^{\left(F\right)}$, and in order
to produce suitable orthonormal bases in $\mathscr{H}_{F}$, we be
concerned with an explicit family of skew-adjoint extensions of $D^{\left(F\right)}$,
as well as the associated spectra, see Corollaries \ref{cor:defg}
and \ref{cor:ptspec}. \index{extension!selfadjoint-}\index{orthogonal!-vectors}
\begin{lem}
\label{lem:dense-1}Let $\varphi$ be a function such that
\begin{enumerate}
\item $\mathrm{supp}\left(\varphi\right)\subset\left(0,a\right)$;
\item $\varphi\in C_{c}^{\infty}\left(0,a\right)$, $\varphi\geq0$;
\item $\int_{0}^{a}\varphi\left(t\right)dt=1$. 
\end{enumerate}

Fix $x\in\left(0,a\right)$, and set $\varphi_{n,x}\left(t\right):=n\varphi\left(n\left(t-x\right)\right)$.
Then $\lim_{n\rightarrow\infty}\varphi_{n,x}=\delta_{x}$, i.e., the
Dirac measure at $x$; and 
\begin{equation}
\left\Vert F_{\varphi_{n,x}}-F\left(\cdot-x\right)\right\Vert _{\mathscr{H}_{F}}\rightarrow0,\;\mbox{as }n\rightarrow\infty.\label{eq:approx}
\end{equation}
Hence $\left\{ F_{\varphi}:\varphi\in C_{c}^{\infty}\left(0,a\right)\right\} $
spans a dense subspace in $\mathscr{H}_{F}$. See \figref{approx}.

\end{lem}

Recall, the following facts about $\mathscr{H}_{F},$ which follow
from the general theory \cite{Aro50} of RKHS:
\begin{itemize}
\item $F(0)>0,$ so we can always arrange $F(0)=1.$
\item $F(-x)=\overline{F(x)}$
\item $\mathscr{H}_{F}$ consists of continuous functions $\xi:\Omega-\Omega\rightarrow\mathbb{C}.$
\item The reproducing property:
\[
\left\langle F\left(\cdot-x\right),\xi\right\rangle _{\mathscr{H}_{F}}=\xi\left(x\right),\;\forall\xi\in\mathscr{H}_{F},\forall x\in\Omega,
\]
is a direct consequence of (\ref{eq:ip-discrete}).\end{itemize}
\begin{rem}
It follows from the reproducing property that if $F_{\phi_{n}}\to\xi$
in $\mathscr{H}_{F},$ then $F_{\phi_{n}}$ converges uniformly to
$\xi$ in $\Omega.$ In fact 
\begin{align*}
\left|F_{\phi_{n}}\left(x\right)-\xi\left(x\right)\right| & =\left|\left\langle F\left(\cdot-x\right),F_{\phi_{n}}-\xi\right\rangle _{\mathscr{H}_{F}}\right|\\
 & \leq\left\Vert F\left(\cdot-x\right)\right\Vert _{\mathscr{H}_{F}}\left\Vert F_{\phi_{n}}-\xi\right\Vert _{\mathscr{H}_{F}}\\
 & =F\left(0\right)^{1/2}\left\Vert F_{\phi_{n}}-\xi\right\Vert _{\mathscr{H}_{F}}.
\end{align*}
\end{rem}
\begin{lem}
\label{lem:TFvarphi}Let $F:\left(-a,a\right)\rightarrow\mathbb{C}$
be a continuous and p.d. function, and let $\mathscr{H}_{F}$ be the
corresponding RKHS. Then:
\begin{enumerate}
\item \label{enu:F1-1}the integral $F_{\varphi}:=\int_{0}^{a}\varphi\left(y\right)F\left(\cdot-y\right)dy$
is convergent in $\mathscr{H}_{F}$ for all $\varphi\in C_{c}\left(0,a\right)$;
and 
\item \label{enu:F1-2}for all $\xi\in\mathscr{H}_{F}$, we have:
\begin{equation}
\left\langle F_{\varphi},\xi\right\rangle _{\mathscr{H}_{F}}=\int_{0}^{a}\overline{\varphi\left(x\right)}\xi\left(x\right)dx.\label{eq:F1-1}
\end{equation}

\end{enumerate}
\end{lem}
\begin{proof}
For simplicity, we assume the following normalization $F\left(0\right)=1$;
then for all $y_{1},y_{2}\in\left(0,1\right)$, we have 
\begin{equation}
\left\Vert F\left(\cdot-y_{1}\right)-F\left(\cdot-y_{2}\right)\right\Vert _{\mathscr{H}_{F}}^{2}=2\left(1-\Re\left\{ F\left(y_{1}-y_{2}\right)\right\} \right).\label{eq:F1-2}
\end{equation}
Now, view the integral in (\ref{enu:F1-1}) as a $\mathscr{H}_{F}$-vector
valued integral. If $\varphi\in C_{c}\left(0,a\right)$, this integral
$\int_{0}^{a}\varphi\left(y\right)F\left(\cdot-y\right)dy$ is the
$\mathscr{H}_{F}$-norm convergent. Since $\mathscr{H}_{F}$ is a
RKHS, $\left\langle \cdot,\xi\right\rangle _{\mathscr{H}_{F}}$ is
continuous on $\mathscr{H}_{F}$, and it passes under the integral
in (\ref{enu:F1-1}). Using 
\begin{equation}
\left\langle F\left(y-\cdot\right),\xi\right\rangle _{\mathscr{H}_{F}}=\xi\left(y\right)\label{eq:F1-3}
\end{equation}
the desired conclusion (\ref{eq:F1-1}) follows.
\end{proof}

\begin{cor}
Let $F:\left(-a,a\right)\rightarrow\mathbb{C}$ be as above, and let
$\mathscr{H}_{F}$ be the corresponding RKHS. For $\varphi\in C_{c}^{1}\left(0,a\right)$,
set 
\begin{equation}
F_{\varphi}\left(x\right)=\left(T_{F}\varphi\right)\left(x\right)=\int_{0}^{a}\varphi\left(y\right)F\left(x-y\right)dy;\label{eq:F1-4}
\end{equation}
then $F_{\varphi}\in C^{1}\left(0,a\right)$, and 
\begin{equation}
\frac{d}{dx}F_{\varphi}\left(x\right)=\left(T_{F}\left(\varphi'\right)\right)\left(x\right),\;\forall x\in\left(0,a\right).\label{eq:F1-5}
\end{equation}
\end{cor}
\begin{proof}
Since $F_{\varphi}\left(x\right)=\int_{0}^{a}\varphi\left(y\right)F\left(x-y\right)dy$,
$x\in\left(0,a\right)$; the desired assertion (\ref{eq:F1-5}) follows
directly from the arguments in the proof of \lemref{TFvarphi}.\end{proof}
\begin{thm}
\label{thm:HF}Fix $0<a<\infty$. A continuous function $\xi:\left(0,a\right)\rightarrow\mathbb{C}$
is in $\mathscr{H}_{F}$ if and only if there exists a finite constant
$A>0$, such that
\begin{equation}
\sum_{i}\sum_{j}\overline{c_{i}}c_{j}\overline{\xi\left(x_{i}\right)}\xi\left(x_{j}\right)\leq A\sum_{i}\sum_{j}\overline{c_{i}}c_{j}F\left(x_{i}-x_{j}\right)\label{eq:bdd}
\end{equation}
for all finite system $\left\{ c_{i}\right\} \subset\mathbb{C}$ and
$\left\{ x_{i}\right\} \subset\left(0,a\right)$. Equivalently, for
all $\varphi\in C_{c}^{\infty}\left(\Omega\right)$, 
\begin{align}
\left|\int_{0}^{a}\varphi\left(y\right)\xi\left(y\right)dy\right|^{2} & \leq A\int_{0}^{a}\int_{0}^{a}\overline{\varphi\left(x\right)}\varphi\left(y\right)F\left(x-y\right)dxdy\label{eq:bdd2}
\end{align}

\end{thm}
We will use these two conditions (\ref{eq:bdd})($\Leftrightarrow$(\ref{eq:bdd2}))
when considering for example the von Neumann deficiency-subspaces
for skew Hermitian operators with dense domain in $\mathscr{H}_{F}$.
\begin{proof}[Proof of \thmref{HF}]
Note, if $\xi\in\mathscr{H}_{F}$, then 
\[
\mbox{LHS}_{\left(\ref{eq:bdd2}\right)}=\big|\left\langle F_{\varphi},\xi\right\rangle _{\mathscr{H}_{F}}\big|^{2},
\]
and so (\ref{eq:bdd2}) holds, since $\left\langle \cdot,\xi\right\rangle _{\mathscr{H}_{F}}$
is continuous on $\mathscr{H}_{F}$.

If $\xi$ is continuous on $\left[0,a\right]$, and if (\ref{eq:bdd2})
holds, then 
\[
\mathscr{H}_{F}\ni F_{\varphi}\longmapsto\int_{0}^{a}\varphi\left(y\right)\xi\left(y\right)dy
\]
is well-defined, continuous, linear; and extends to $\mathscr{H}_{F}$
by density (see \lemref{dense-1}). Hence, by Riesz' theorem, $\exists!$
$k_{\xi}\in\mathscr{H}_{F}$ such that \index{Riesz' theorem}\index{Theorem!Riesz-}
\[
\int_{0}^{a}\varphi\left(y\right)\xi\left(y\right)dy=\left\langle F_{\varphi},k_{\xi}\right\rangle _{\mathscr{H}_{F}}.
\]
But using the reproducing property in $\mathscr{H}_{F}$, and $F_{\varphi}\left(x\right)=\int_{0}^{a}\varphi\left(x\right)F\left(x-y\right)dy$,
we get 
\[
\int_{0}^{a}\overline{\varphi\left(x\right)}\xi\left(x\right)dx=\int_{0}^{a}\overline{\varphi\left(x\right)}k_{\xi}\left(x\right)dx,\;\forall\varphi\in C_{c}\left(0,a\right)
\]
so 
\[
\int_{0}^{a}\varphi\left(x\right)\left(\xi\left(x\right)-k_{\xi}\left(x\right)\right)dx=0,\;\forall\varphi\in C_{c}\left(0,a\right);
\]
it follows that $\xi-k_{\xi}=0$ on $\left(0,a\right)$ $\Longrightarrow$
$\xi-k_{\xi}=0$ on $\left[0,a\right]$. \end{proof}
\begin{defn}[The operator $D_{F}$]
\label{def:DF} Let $D_{F}\left(F_{\varphi}\right)=F_{\varphi'}$,
for all $\varphi\in C_{c}^{\infty}\left(0,a\right)$, where $\varphi'=\frac{d\varphi}{dt}$
and $F_{\varphi}$ is as in (\ref{eq:H2}). \end{defn}
\begin{lem}
\label{lem:DF}The operator $D_{F}$ defines a skew-Hermitian operator
with dense domain in $\mathscr{H}_{F}$.\end{lem}
\begin{proof}
By \lemref{dense-1}, $dom\left(D_{F}\right)$ is dense in $\mathscr{H}_{F}.$
If $\psi\in C_{c}^{\infty}\left(0,a\right)$ and 
\[
\left|t\right|<\mathrm{dist}\left(\mathrm{supp}\left(\psi\right),\mbox{endpoints}\right),
\]
then
\begin{equation}
\left\Vert F_{\psi\left(\cdot+t\right)}\right\Vert _{\mathscr{H}_{F}}^{2}=\left\Vert F_{\psi}\right\Vert _{\mathscr{H}_{F}}^{2}=\int_{0}^{a}\int_{0}^{a}\overline{\psi\left(x\right)}\psi\left(y\right)F\left(x-y\right)dxdy
\end{equation}
see (\ref{eq:hn2}), so 
\[
\frac{d}{dt}\left\Vert F_{\psi\left(\cdot+t\right)}\right\Vert _{\mathscr{H}_{F}}^{2}=0
\]
which is equivalent to 
\begin{equation}
\left\langle D_{F}F_{\psi},F_{\psi}\right\rangle _{\mathscr{H}_{F}}+\left\langle F_{\psi},D_{F}F_{\psi}\right\rangle _{\mathscr{H}_{F}}=0.
\end{equation}
It follows that $D_{F}$ is well-defined and skew-Hermitian in $\mathscr{H}_{F}$. \end{proof}
\begin{lem}
\label{lem:Dadj}Let $F$ be a positive definite function on $\left(-a,a\right)$,
$0<a<\infty$ fixed. Let $D_{F}$ be as in \defref{DF}, so that $D_{F}\subset D_{F}^{*}$
(\lemref{DF}), where $D_{F}^{*}$ is the adjoint relative to the
$\mathscr{H}_{F}$ inner product. 

Then $\xi\in\mathscr{H}_{F}$ (as a continuous function on $\left[0,a\right]$)
is in $dom\left(D_{F}^{*}\right)$ iff 
\begin{gather}
\xi'\in\mathscr{H}_{F}\;\mbox{where }\xi'=\mbox{distribution derivative, and}\label{eq:Dadj1}\\
D_{F}^{*}\xi=-\xi'\label{eq:Dadj2}
\end{gather}
\end{lem}
\begin{proof}
By \thmref{HF}, a fixed $\xi\in\mathscr{H}_{F}$, i.e., $x\mapsto\xi\left(x\right)$
is a continuous function on $\left[0,a\right]$ such that $\exists C$,
$\left|\int_{0}^{a}\varphi\left(x\right)\xi\left(x\right)dx\right|^{2}\leq C\left\Vert F_{\varphi}\right\Vert _{\mathscr{H}_{F}}^{2}$. 

$\xi$ is in $dom\left(D_{F}^{*}\right)$ $\Longleftrightarrow$ $\exists C=C_{\xi}<\infty$
such that 
\begin{equation}
\left|\left\langle D_{F}\left(F_{\varphi}\right),\xi\right\rangle _{\mathscr{H}_{F}}\right|^{2}\leq C\left\Vert F_{\varphi}\right\Vert _{\mathscr{H}_{F}}^{2}=C\int_{0}^{a}\int_{0}^{a}\overline{\varphi\left(x\right)}\varphi\left(y\right)F\left(x-y\right)dxdy\label{eq:Dadj3}
\end{equation}
But LHS of (\ref{eq:Dadj3}) under $\left|\left\langle \cdot,\cdot\right\rangle \right|^{2}$
is:
\begin{equation}
\left|\left\langle D_{F}\left(F_{\varphi}\right),\xi\right\rangle _{\mathscr{H}_{F}}\right|^{2}=\left\langle F_{\varphi'},\xi\right\rangle _{\mathscr{H}_{F}}\overset{\left(\ref{eq:F1-1}\right)}{=}\int_{0}^{a}\overline{\varphi'\left(x\right)}\xi\left(x\right)dx,\;\forall\varphi\in C_{c}^{\infty}\left(0,a\right)\label{eq:Dadj4}
\end{equation}
So (\ref{eq:Dadj3}) holds $\Longleftrightarrow$ 
\[
\left|\int_{0}^{a}\overline{\varphi'\left(x\right)}\xi\left(x\right)dx\right|^{2}\leq C\left\Vert F_{\varphi}\right\Vert _{\mathscr{H}_{F}}^{2},\;\forall\varphi\in C_{c}^{\infty}\left(0,a\right)
\]
i.e., 
\[
\left|\int_{0}^{a}\overline{\varphi\left(x\right)}\xi'\left(x\right)dx\right|^{2}\leq C\left\Vert F_{\varphi}\right\Vert _{\mathscr{H}_{F}}^{2},\;\forall\varphi\in C_{c}^{\infty}\left(0,a\right),\;\mbox{and}
\]
$\xi'$ as a distribution is in $\mathscr{H}_{F}$, and 
\[
\int_{0}^{a}\overline{\varphi\left(x\right)}\xi'\left(x\right)dx=\left\langle F_{\varphi},\xi'\right\rangle _{\mathscr{H}_{F}}
\]
where we use the characterization of $\mathscr{H}_{F}$ in (\ref{eq:bdd2}),
i.e., a function $\eta:\left[0,a\right]\rightarrow\mathbb{C}$ is
in $\mathscr{H}_{F}$ $\Longleftrightarrow$ $\exists C<\infty$,
$\left|\int_{0}^{a}\overline{\varphi\left(x\right)}\eta\left(x\right)dx\right|\leq C\left\Vert F_{\varphi}\right\Vert _{\mathscr{H}_{F}}$,
$\forall\varphi\in C_{c}^{\infty}\left(0,a\right)$, and then $\int_{0}^{a}\overline{\varphi\left(x\right)}\eta\left(x\right)dx=\left\langle F_{\varphi},\eta\right\rangle _{\mathscr{H}_{F}}$,
$\forall\varphi\in C_{c}^{\infty}\left(0,a\right)$. See \thmref{HF}.\end{proof}
\begin{cor}
$h\in\mathscr{H}_{F}$ is in $dom\left(\left(D_{F}^{2}\right)^{*}\right)$
iff $h''\in\mathscr{H}_{F}$ ($h''$ distribution derivative) and
$\left(D_{F}^{*}\right)^{2}h=\left(D_{F}^{2}\right)^{*}h=h''$. \end{cor}
\begin{proof}
Application of (\ref{eq:Dadj4}) to $D_{F}\left(F_{\varphi}\right)=F_{\varphi'}$,
we have $D_{F}^{2}\left(F_{\varphi}\right)=F_{\varphi''}=\left(\frac{d}{dx}\right)^{2}F_{\varphi}$,
$\forall\varphi\in C_{c}^{\infty}\left(0,a\right)$, and 
\begin{align*}
\left\langle D_{F}^{2}\left(F_{\varphi}\right),h\right\rangle _{\mathscr{H}_{F}} & =\left\langle F_{\varphi''},h\right\rangle _{\mathscr{H}_{F}}=\int_{0}^{a}\overline{\varphi''\left(x\right)}h\left(x\right)dx\\
 & =\int_{0}^{a}\overline{\varphi\left(x\right)}h''\left(x\right)dx=\left\langle F_{\varphi},\left(D_{F}^{2}\right)^{*}h\right\rangle _{\mathscr{H}_{F}}.
\end{align*}

\end{proof}
\index{distribution!Schwartz-}
\begin{defn}
\cite{DS88b}Let $D_{F}^{*}$ be the adjoint of $D_{F}$ relative
to $\mathscr{H}_{F}$ inner product. The deficiency spaces $DEF^{\pm}$
consists of $\xi_{\pm}\in dom\left(D_{F}^{*}\right)$, such that $D_{F}^{*}\xi_{\pm}=\pm\xi_{\pm}$,
i.e.,\index{deficiency space} 
\[
DEF^{\pm}=\left\{ \xi_{\pm}\in\mathscr{H}_{F}:\left\langle F_{\psi'},\xi_{\pm}\right\rangle _{\mathscr{H}_{F}}=\left\langle F_{\psi},\pm\xi_{\pm}\right\rangle _{\mathscr{H}_{F}},\forall\psi\in C_{c}^{\infty}\left(\Omega\right)\right\} .
\]
\end{defn}
\begin{cor}
\label{cor:defg}If $\xi\in DEF^{\pm}$ then $\xi(x)=\mathrm{constant}\,e^{\mp x}.$\end{cor}
\begin{proof}
Immediate from \lemref{Dadj}.
\end{proof}
The role of deficiency indices for the canonical skew-Hermitian operator
$D_{F}$ (\defref{DF}) in the RKHS $\mathscr{H}_{F}$ is as follows:
using von Neumann's conjugation trick \cite{DS88b}, we see that the
deficiency indices can be only $\left(0,0\right)$ or $\left(1,1\right)$.\index{deficiency indices}

We conclude that there exists proper skew-adjoint extensions $A\supset D_{F}$
in $\mathscr{H}_{F}$ (in case $D_{F}$ has indices $\left(1,1\right)$).
Then 
\begin{equation}
D_{F}\subseteq A=-A^{*}\subseteq-D_{F}^{*}\label{eq:Dadj7}
\end{equation}
(If the indices are $\left(0,0\right)$ then $\overline{D_{F}}=-D_{F}^{*}$;
see \cite{DS88b}.)

Hence, set $U\left(t\right)=e^{tA}:\mathscr{H}_{F}\rightarrow\mathscr{H}_{F}$,
and get the strongly continuous unitary one-parameter group 
\[
\left\{ U\left(t\right):t\in\mathbb{R}\right\} ,\;U\left(s+t\right)=U\left(s\right)U\left(t\right),\:\forall s,t\in\mathbb{R};
\]
and if 
\[
\xi\in dom\left(A\right)=\left\{ \xi\in\mathscr{H}_{F}:\:\mbox{s.t.}\lim_{t\rightarrow0}\frac{U\left(t\right)\xi-\xi}{t}\:\mbox{exists}\right\} 
\]
then 
\begin{equation}
A\xi=\mbox{s.t.}\lim_{t\rightarrow0}\frac{U\left(t\right)\xi-\xi}{t}.
\end{equation}

\index{strongly continuous}

Now use $F_{x}(\cdot)=F\left(x-\cdot\right)$ defined in $\left(0,a\right)$;
and set 
\begin{equation}
F_{A}\left(t\right):=\left\langle F_{0},U\left(t\right)F_{0}\right\rangle _{\mathscr{H}_{F}},\;\forall t\in\mathbb{R}\label{eq:Fext}
\end{equation}
then using (\ref{eq:approx}), we see that $F_{A}$ is a continuous
positive definite extension of $F$ on $\left(-a,a\right)$. This
extension is in $Ext_{1}\left(F\right)$. 
\begin{cor}
\label{cor:ptspec}Assume $\lambda\in\mathbb{R}$ is in the point
spectrum of $A$, i.e., $\exists\xi_{\lambda}\in dom\left(A\right)$,
$\xi_{\lambda}\neq0$, such that $A\xi_{\lambda}=i\lambda\xi_{\lambda}$
holds in $\mathscr{H}_{F}$, then $\xi_{\lambda}=\mbox{const}\cdot e_{\lambda}$,
i.e., 
\begin{equation}
\xi_{\lambda}\left(x\right)=\mbox{const}\cdot e^{i\lambda x},\;\forall x\in\left[0,a\right].\label{eq:Dadj8}
\end{equation}
\end{cor}
\begin{proof}
Assume $\lambda$ is in $spec_{pt}\left(A\right)$, and $\xi_{\lambda}\in dom\left(A\right)$
satisfying 
\begin{equation}
\left(A\xi_{\lambda}\right)\left(x\right)=i\lambda\xi_{\lambda}\left(x\right)\;\mbox{in }\mathscr{H}_{F},\label{eq:Dadj9}
\end{equation}
then since $A\subset-D_{F}^{*}$, we get $\xi\in dom\left(D_{F}^{*}\right)$
by \lemref{Dadj} and (\ref{eq:Dadj7}), and $D_{F}^{*}\xi_{\lambda}=-\xi_{\lambda}'$
where $\xi'$ is the distribution derivative (see (\ref{eq:Dadj2}));
and by (\ref{eq:Dadj7})
\begin{equation}
\left(A\xi_{\lambda}\right)\left(x\right)=-\left(D_{F}^{*}\xi_{\lambda}\right)\left(x\right)=\xi'_{\lambda}\left(x\right)\overset{\left(\ref{eq:Dadj9}\right)}{=}i\lambda\xi_{\lambda}\left(x\right),\;\forall x\in\left(0,a\right)\label{eq:Dadj10}
\end{equation}
so $\xi_{\lambda}$ is the distribution derivative solution to 
\begin{eqnarray}
\xi'_{\lambda}\left(x\right) & = & i\lambda\xi_{\lambda}\left(x\right)\label{eq:Dadj11}\\
 & \Updownarrow\nonumber \\
-\int_{0}^{a}\overline{\varphi'\left(x\right)}\xi_{\lambda}\left(x\right)dx & = & i\lambda\int_{0}^{a}\overline{\varphi\left(x\right)}\xi_{\lambda}\left(x\right)dx,\;\forall\varphi\in C_{c}^{\infty}\left(0,a\right)\nonumber \\
 & \Updownarrow\nonumber \\
-\left\langle D_{F}\left(F_{\varphi}\right),\xi_{\lambda}\right\rangle _{\mathscr{H}_{F}} & = & i\lambda\left\langle F_{\varphi},\xi_{\lambda}\right\rangle _{\mathscr{H}_{F}},\;\forall\varphi\in C_{c}^{\infty}\left(0,a\right).\nonumber 
\end{eqnarray}
But by Schwartz, the distribution solutions to (\ref{eq:Dadj11})
are $\xi_{\lambda}\left(x\right)=\mbox{const}\cdot e_{\lambda}\left(x\right)=\mbox{const}\cdot e^{i\lambda x}$. 
\end{proof}

In the considerations below, we shall be primarily concerned with
the case when a fixed continuous p.d. function $F$ is defined on
a finite interval $\left(-a,a\right)\subset\mathbb{R}$. In this case,
by a Mercer operator, we mean an operator $T_{F}$ in $L^{2}\left(0,a\right)$
where $L^{2}\left(0,a\right)$ is defined from Lebesgue measure on
$\left(0,a\right)$, given by
\begin{equation}
\left(T_{F}\varphi\right)\left(x\right):=\int_{0}^{a}\varphi\left(y\right)F\left(x-y\right)dy,\;\forall\varphi\in L^{2}\left(0,a\right),\forall x\in\left(0,a\right).\label{eq:mer-1}
\end{equation}

\begin{lem}
\label{lem:mer-1}Under the assumptions stated above, the Mercer operator
$T_{F}$ is trace class in $L^{2}\left(0,a\right)$; and if $F\left(0\right)=1$,
then \index{trace} 
\begin{equation}
trace\left(T_{F}\right)=a.\label{eq:mer-2}
\end{equation}
\end{lem}
\begin{proof}
This is an application of Mercer's theorem \cite{LP89,FR42,FM13}
to the integral operator $T_{F}$ in (\ref{eq:mer-1}). But we must
check that $F$, on $\left(-a,a\right)$, extends uniquely by limit
to a continuous p.d. function $F_{ex}$ on $\left[-a,a\right]$, the
closed interval. This is true, and easy to verify, see e.g. \cite{JPT14}.
\index{operators!integral-}\index{Theorem!Mercer's-}\end{proof}
\begin{cor}
\label{cor:mer1}Let $F$ and $\left(-a,a\right)$ be as in \lemref{mer-1}.
Then there is a sequence $\left(\lambda_{n}\right)_{n\in\mathbb{N}}$,
$\lambda_{n}>0$, such that $\sum_{n\in\mathbb{N}}\lambda_{n}=a$,
and a system of orthogonal functions $\left\{ \xi_{n}\right\} \subset L^{2}\left(0,a\right)\cap\mathscr{H}_{F}$
such that
\begin{equation}
F\left(x-y\right)=\sum_{n\in\mathbb{N}}\lambda_{n}\xi_{n}\left(x\right)\overline{\xi_{n}\left(y\right)},\mbox{ and}\label{eq:mer-3}
\end{equation}
\begin{equation}
\int_{0}^{a}\overline{\xi_{n}\left(x\right)}\xi_{m}\left(x\right)dx=\delta_{n,m},\;n,m\in\mathbb{N}.\label{eq:mer-4}
\end{equation}
\end{cor}
\begin{proof}
An application of Mercer's theorem \cite{LP89,FR42,FM13}. See also
 \exerref{mer}.\end{proof}
\begin{cor}
\label{cor:merinn}For all $\psi,\varphi\in C_{c}^{\infty}\left(0,a\right)$,
we have 
\begin{equation}
\left\langle F_{\psi},F_{\varphi}\right\rangle _{\mathscr{H}_{F}}=\left\langle F_{\psi},T_{F}^{-1}F_{\varphi}\right\rangle _{2}.\label{eq:mnorm1}
\end{equation}
Consequently, 
\begin{equation}
\left\Vert h\right\Vert _{\mathscr{H}_{F}}=\Vert T_{F}^{-1/2}h\Vert_{2},\;\forall h\in\mathscr{H}_{F}.\label{eq:mnorm2}
\end{equation}
\end{cor}
\begin{proof}
Note 
\begin{align*}
\left\langle F_{\psi},T_{F}^{-1}F_{\varphi}\right\rangle _{2} & =\left\langle F_{\psi},T_{F}^{-1}T_{F}\varphi\right\rangle _{2}=\left\langle F_{\psi},\varphi\right\rangle _{2}\\
 & =\int_{0}^{a}\overline{\left(\int_{0}^{a}\psi\left(x\right)F\left(y-x\right)dx\right)}\,\varphi\left(y\right)dy\\
 & =\int_{0}^{a}\int_{0}^{a}\overline{\psi\left(x\right)}\varphi\left(y\right)F\left(x-y\right)dxdy=\left\langle F_{\psi},F_{\varphi}\right\rangle _{\mathscr{H}_{F}}.
\end{align*}
\end{proof}
\begin{cor}
Let $\left\{ \xi_{n}\right\} $ be the ONB in $L^{2}\left(0,a\right)$
as in \corref{mer1}; then $\left\{ \sqrt{\lambda_{n}}\xi_{n}\right\} $
is an ONB in $\mathscr{H}_{F}$. \end{cor}
\begin{proof}
The functions $\xi_{n}$ are in $\mathscr{H}_{F}$ by \thmref{HF}.
We check directly (\corref{merinn}) that 
\begin{align*}
\left\langle \sqrt{\lambda_{n}}\xi_{n},\sqrt{\lambda_{m}}\xi_{m}\right\rangle _{\mathscr{H}_{F}} & =\sqrt{\lambda_{n}\lambda_{m}}\left\langle \xi_{n},T^{-1}\xi_{m}\right\rangle _{2}\\
 & =\sqrt{\lambda_{n}\lambda_{m}}\lambda_{m}^{-1}\left\langle \xi_{n},\xi_{m}\right\rangle _{2}=\delta_{n,m}.
\end{align*}

\end{proof}

\section{\label{sec:types}Type I v.s. Type II Extensions}

When a pair $\left(\Omega,F\right)$ is given, where $F$ is a prescribed
continuous positive definite function defined on $\Omega$, we consider
the possible continuous positive definite extensions to all of $\mathbb{R}^{n}$.
The reproducing kernel Hilbert space $\mathscr{H}_{F}$ will play
a key role in our analysis. In constructing various classes of continuous
positive definite extensions to $\mathbb{R}^{n}$, we introduce operators
in $\mathscr{H}_{F}$, and their \emph{dilation} to operators, possibly
acting in an enlargement Hilbert space \cite{JPT14,KL14}. Following
techniques from dilation theory we note that every dilation contains
a minimal one. If a continuous positive definite extensions to $\mathbb{R}^{n}$
has its minimal dilation Hilbert space equal to $\mathscr{H}_{F}$,
we say it is type 1, otherwise we say it is type 2. \index{positive definite!-function}
\begin{defn}
Let $G$ be a locally compact group\index{groups!locally compact},
and let $\Omega$ be an open connected subset of $G$. Let $F:\Omega^{-1}\cdot\Omega\rightarrow\mathbb{C}$
be a continuous positive definite function.
\end{defn}
\index{representation!unitary}

\index{representation!strongly continuous}

\index{kernel!reproducing-}
\begin{defn}
Consider a strongly continuous unitary representation $U$ of $G$
acting in some Hilbert space $\mathscr{K}$, containing the RKHS $\mathscr{H}_{F}$.
We say that $\left(U,\mathscr{K}\right)\in Ext\left(F\right)$ iff
there is a vector $k_{0}\in\mathscr{K}$ such that
\begin{equation}
F\left(g\right)=\left\langle k_{0},U\left(g\right)k_{0}\right\rangle _{\mathscr{K}},\;\forall g\in\Omega^{-1}\cdot\Omega.\label{eq:ext-1-1}
\end{equation}

\begin{enumerate}[leftmargin=*,label=\arabic{enumi}.]
\item The subset of $Ext\left(F\right)$ consisting of $\left(U,\mathscr{H}_{F},k_{0}=F_{e}\right)$
with 
\begin{equation}
F\left(g\right)=\left\langle F_{e},U\left(g\right)F_{e}\right\rangle _{\mathscr{H}_{F}},\;\forall g\in\Omega^{-1}\cdot\Omega\label{eq:ext-1-2-1}
\end{equation}
is denoted $Ext_{1}\left(F\right)$; and we set 
\[
Ext_{2}\left(F\right):=Ext\left(F\right)\backslash Ext_{1}\left(F\right);
\]
i.e., $Ext_{2}\left(F\right)$, consists of the solutions to problem
(\ref{eq:ext-1-1}) for which $\mathscr{K}\supsetneqq\mathscr{H}_{F}$,
i.e., unitary representations realized in an enlargement Hilbert space.
\\
(We write $F_{e}\in\mathscr{H}_{F}$ for the vector satisfying $\left\langle F_{e},\xi\right\rangle _{\mathscr{H}_{F}}=\xi\left(e\right)$,
$\forall\xi\in\mathscr{H}_{F}$, where $e$ is the neutral (unit)
element in $G$, i.e., $e\,g=g$, $\forall g\in G$.)
\item In the special case, where $G=\mathbb{R}^{n}$, and $\Omega\subset\mathbb{R}^{n}$
is open and connected, we consider 
\[
F:\Omega-\Omega\rightarrow\mathbb{C}
\]
continuous and positive definite. In this case,
\begin{align}
Ext\left(F\right)= & \Bigl\{\mu\in\mathscr{M}_{+}\left(\mathbb{R}^{n}\right)\:\big|\:\widehat{\mu}\left(x\right)=\int_{\mathbb{R}^{n}}e^{i\lambda\cdot x}d\mu\left(\lambda\right)\label{eq:ext-1-4-1}\\
 & \mbox{ is a p.d. extensiont of \ensuremath{F}}\Bigr\}.\nonumber 
\end{align}

\end{enumerate}
\end{defn}
\begin{rem}
Note that (\ref{eq:ext-1-4-1}) is consistent with (\ref{eq:ext-1-1}):
For if $\left(U,\mathscr{K},k_{0}\right)$ is a unitary representation
of $G=\mathbb{R}^{n}$, such that (\ref{eq:ext-1-1}) holds; then,
by a theorem of Stone, there is a projection-valued\index{measure!projection-valued}
measure (PVM) $P_{U}\left(\cdot\right)$, defined on the Borel subsets
of $\mathbb{R}^{n}$ such that 
\begin{equation}
U\left(x\right)=\int_{\mathbb{R}^{n}}e^{i\lambda\cdot x}P_{U}\left(d\lambda\right),\;x\in\mathbb{R}^{n}.\label{eq:ex-1-5}
\end{equation}
Setting 
\begin{equation}
d\mu\left(\lambda\right):=\left\Vert P_{U}\left(d\lambda\right)k_{0}\right\Vert _{\mathscr{K}}^{2},\label{eq:ext-1-6-7}
\end{equation}
it is then immediate that we have: $\mu\in\mathscr{M}_{+}\left(\mathbb{R}^{n}\right)$,
and that the finite measure $\mu$ satisfies 
\begin{equation}
\widehat{\mu}\left(x\right)=F\left(x\right),\;\forall x\in\Omega-\Omega.\label{eq:ext-1-6-1}
\end{equation}

\end{rem}
Set $n=1$: Start with a local p.d. continuous function $F$, and
let $\mathscr{H}_{F}$ be the corresponding RKHS. Let $Ext(F)$ be
the compact convex set of probability measures on $\mathbb{R}$ defining
extensions of $F$.

We now divide $Ext(F)$ into two parts, say $Ext_{1}\left(F\right)$
and $Ext_{2}\left(F\right)$. 

All continuous p.d. extensions of $F$ come from strongly continuous
unitary representations. So in the case of 1D, from unitary one-parameter
groups of course, say $U(t)$. 

\index{strongly continuous}

\index{representation!unitary}

\index{representation!strongly continuous}

Let $Ext_{1}\left(F\right)$ be the subset of $Ext(F)$ corresponding
to extensions when the unitary representation $U(t)$ acts in $\mathscr{H}_{F}$
(internal extensions), and $Ext_{2}\left(F\right)$ denote the part
of $Ext(F)$ associated to unitary representations $U(t)$ acting
in a proper enlargement Hilbert space $\mathscr{K}$ (if any), i.e.,
acting in a Hilbert space $\mathscr{K}$ corresponding to a proper
dilation of $\mathscr{H}_{F}$.

\section{\label{sec:exp}The Case of $e^{-\left|x\right|}$, $\left|x\right|<1$}

Our emphasis is von Neumann indices, and explicit formulas for partially
defined positive definite functions $F$, defined initially only on
a symmetric interval $\left(-a,a\right)$. Among the cases of partially
defined positive definite functions, the following example $F\left(x\right)=e^{-\left|x\right|}$,
in the symmetric interval $\left(-1,1\right)$, will play a special
role. The present section is devoted to this example.

There are many reasons for this:

\index{reproducing kernel Hilbert space (RKHS)}

\index{kernel!reproducing-}
\begin{enumerate}[leftmargin=*,label=(\roman{enumi})]
\item It is of independent interest, and its type 1 extensions (see \secref{types})
can be written down explicitly.
\item Its applications include stochastic analysis \cite{Ito06} as follows.
Given a random variable $X$ in a process; if $\mu$ is its distribution,
then there are two measures of concentration for $\mu$, one called
\textquotedblleft degree of concentration,\textquotedblright{} and
the other \textquotedblleft dispersion,\textquotedblright{} both computed
directly from $F\left(x\right)=e^{-\left|x\right|}$ applied to $\mu$.
\index{distribution!probability-}
\item In addition, there are analogous relative notions for comparing different
samples in a fixed stochastic process. These notions are defined with
the use of example $F\left(x\right)=e^{-\left|x\right|}$, and it
will frequently be useful to localize the $x$-variable in a compact
interval.
\item Additional reasons for special attention to example $F\left(x\right)=e^{-\left|x\right|}$,
for $x\in\left(-1,1\right)$ is its use in sampling theory, and analysis
of de Branges spaces \cite{DM85}, as well as its role as a Greens
function for an important boundary value problem. 
\item Related to this, the reproducing kernel Hilbert space $\mathscr{H}_{F}$
associated to this p.d. function $F$ has a number of properties that
also hold for wider families of locally defined positive definite
function of a single variable. In particular, $\mathscr{H}_{F}$ has
Fourier bases: The RKHS $\mathscr{H}_{F}$ has orthogonal bases of
complex exponentials $e_{\lambda}$ with aperiodic frequency distributions,
i.e., frequency points $\left\{ e_{\lambda}\right\} $ on the real
line which do not lie on any arithmetic progression, see \figref{expExt}.
For details on this last point, see Corollaries \ref{cor:elambda},
\ref{cor:HFinner2}, \ref{cor:elambda1}, and \ref{cor:Fext}. 
\end{enumerate}

\subsection{\label{sub:saext}The selfadjoint Extensions $A_{\theta}\supset-iD_{F}$}

The notation ``$\supseteq$'' above refers to containment of operators,
or rather of the respective graphs of the two operators; see \cite{DS88b}.
\index{selfadjoint extensions}
\begin{lem}
\label{lem:dev1}Let $F\left(x\right)=e^{-\left|x\right|}$, $\left|x\right|<1$.
Set $F_{x}\left(y\right):=F\left(x-y\right)$, $\forall x,y\in\left(0,1\right)$;
and $F_{\varphi}\left(x\right)=\int_{0}^{1}\varphi\left(y\right)F\left(x-y\right)dy$,
$\forall\varphi\in C_{c}^{\infty}\left(0,1\right)$. Define $D_{F}\left(F_{\varphi}\right)=F_{\varphi'}$
on the dense subset 
\begin{equation}
dom\left(D_{F}\right)=\left\{ F_{\varphi}:\varphi\in C_{c}^{\infty}\left(0,1\right)\right\} \subset\mathscr{H}_{F}.\label{eq:D}
\end{equation}
Then the skew-Hermitian operator $D_{F}$ has deficiency indices $\left(1,1\right)$
in $\mathscr{H}_{F}$, where the defect vectors are 
\begin{align}
\xi_{+}\left(x\right) & =F_{0}\left(x\right)=e^{-x}\label{eq:edev1}\\
\xi_{-}\left(x\right) & =F_{1}\left(x\right)=e^{x-1};\label{eq:edev2}
\end{align}
moreover, 
\begin{equation}
\left\Vert \xi_{+}\right\Vert _{\mathscr{H}_{F}}=\left\Vert \xi_{+}\right\Vert _{\mathscr{H}_{F}}=1.\label{eq:edev3}
\end{equation}
\end{lem}
\begin{proof}
(Note if $\Omega$ is any bounded, open and connected domain in $\mathbb{R}^{n}$,
then a locally defined continuous p.d. function, $F:\Omega-\Omega:\rightarrow\mathbb{C}$,
extends uniquely to the boundary $\partial\Omega:=\overline{\Omega}\backslash\Omega$
by continuity \cite{JPT14}.)

In our current settings, $\Omega=\left(0,1\right)$, and $F_{x}\left(y\right):=F\left(x-y\right)$,
$\forall x,y\in\left(0,1\right)$. Thus, $F_{x}\left(y\right)$ extends
to all $x,y\in\left[0,1\right]$. In particular, 
\[
F_{0}\left(x\right)=e^{-x},\;F_{1}\left(x\right)=e^{x-1}
\]
are the two defect vectors, as shown in \corref{defg}. Moreover,
using the reproducing property, we have 
\begin{align*}
\left\Vert F_{0}\right\Vert _{\mathscr{H}_{F}}^{2} & =\left\langle F_{0},F_{0}\right\rangle _{\mathscr{H}_{F}}=F_{0}\left(0\right)=F\left(0\right)=1\\
\left\Vert F_{1}\right\Vert _{\mathscr{H}_{F}}^{2} & =\left\langle F_{1},F_{1}\right\rangle _{\mathscr{H}_{F}}=F_{1}\left(1\right)=F\left(0\right)=1
\end{align*}
and (\ref{eq:edev3}) follows. For more details, see \cite[lemma 2.10.14]{JPT14}.\end{proof}
\begin{lem}
\label{lem:dev}Let $F$ be any continuous p.d. function on $\left(-1,1\right)$.
Set 
\[
h\left(x\right)=\int_{0}^{1}\varphi\left(y\right)F\left(x-y\right)dy,\;\forall\varphi\in C_{c}^{\infty}\left(0,1\right);
\]
then 
\begin{align}
h\left(0\right) & =\int_{0}^{1}\varphi\left(y\right)F\left(-y\right)dy,\qquad h\left(1\right)=\int_{0}^{1}\varphi\left(y\right)F\left(1-y\right)dy\label{eq:dev-1-1}\\
h'\left(0\right) & =\int_{0}^{1}\varphi\left(y\right)F'\left(-y\right)dy,\quad\,\,\,h'\left(1\right)=\int_{0}^{1}\varphi\left(y\right)F'\left(1-y\right)dy;\label{eq:dev-1-2}
\end{align}
where the derivatives $F'$ in (\ref{eq:dev-1-1})-(\ref{eq:dev-1-2})
are in the sense of distribution.\end{lem}
\begin{proof}
Note that
\begin{align*}
h\left(x\right) & =\int_{0}^{x}\varphi\left(y\right)F\left(x-y\right)dy+\int_{x}^{1}\varphi\left(y\right)F\left(x-y\right)dy;\\
h'\left(x\right) & =\int_{0}^{x}\varphi\left(y\right)F'\left(x-y\right)dy+\int_{x}^{1}\varphi\left(y\right)F'\left(x-y\right)dy.
\end{align*}
and so (\ref{eq:dev-1-1})-(\ref{eq:dev-1-2}) follow.
\end{proof}
We now specialize to the function $F\left(x\right)=e^{-\left|x\right|}$
defined in $\left(-1,1\right)$.
\begin{cor}
\label{cor:dev1}For $F\left(x\right)=e^{-\left|x\right|}$, $\left|x\right|<1$,
set $h=T_{F}\varphi$, i.e., 
\[
h:=F_{\varphi}=\int_{0}^{1}\varphi\left(y\right)F\left(\cdot-y\right)dy,\;\forall\varphi\in C_{c}^{\infty}\left(0,1\right);
\]
then 
\begin{align}
h\left(0\right) & =\int_{0}^{1}\varphi\left(y\right)e^{-y}dy,\qquad h\left(1\right)=\int_{0}^{1}\varphi\left(y\right)e^{y-1}dy\label{eq:edev4}\\
h'\left(0\right) & =\int_{0}^{1}\varphi\left(y\right)e^{-y}dy,\quad\,\,\,\,h'\left(1\right)=-\int_{0}^{1}\varphi\left(y\right)e^{y-1}dy\label{eq:edev5}
\end{align}
In particular, 
\begin{align}
h\left(0\right)-h'\left(0\right) & =0\label{eq:dev6}\\
h\left(1\right)+h'\left(1\right) & =0.\label{eq:dev7}
\end{align}
\end{cor}
\begin{proof}
Immediately from \lemref{dev}. Specifically, 
\begin{align*}
h\left(x\right) & =e^{-x}\int_{0}^{x}\varphi\left(y\right)e^{y}dy+e^{x}\int_{x}^{1}\varphi\left(y\right)e^{-y}dy\\
h'\left(x\right) & =-e^{-x}\int_{0}^{x}\varphi\left(y\right)e^{y}dy+e^{x}\int_{x}^{1}\varphi\left(y\right)e^{-y}dy.
\end{align*}
Setting $x=0$ and $x=1$ gives the desired conclusions.\end{proof}
\begin{rem}
The space 
\[
\Big\{ h\in\mathscr{H}_{F}\:\big|\:h\left(0\right)-h'\left(0\right)=0,\;h\left(1\right)+h'\left(1\right)=0\Big\}
\]
is dense in $\mathscr{H}_{F}$. This is because it contains $\left\{ F_{\varphi}\:\big|\:\varphi\in C_{c}^{\infty}\left(0,1\right)\right\} $.
Note
\begin{align*}
F_{0}+F_{0}' & =-\delta_{0},\;\mbox{and}\\
F_{1}-F_{1}' & =-\delta_{1};
\end{align*}
however,  $\delta_{0},\delta_{1}\notin\mathscr{H}_{F}$.
\end{rem}
By von  Neumann's theory \cite{DS88b} and \lemref{Dadj}, the family
of selfadjoint extensions of the Hermitian operator $-iD_{F}$ is
characterized by
\begin{gather}
\begin{split}A_{\theta}\left(h+c\left(e^{-x}+e^{i\theta}e^{x-1}\right)\right)=-i\,h'+i\,c\left(e^{-x}-e^{i\theta}e^{x-1}\right),\;\mbox{where}\\
dom\left(A_{\theta}\right):=\left\{ h+c\left(e^{-x}+e^{i\theta}e^{x-1}\right)\:\big|\:h\in dom\left(D_{F}\right),c\in\mathbb{C}\right\} .
\end{split}
\label{eq:saextA}
\end{gather}

\begin{figure}
\includegraphics[scale=0.6]{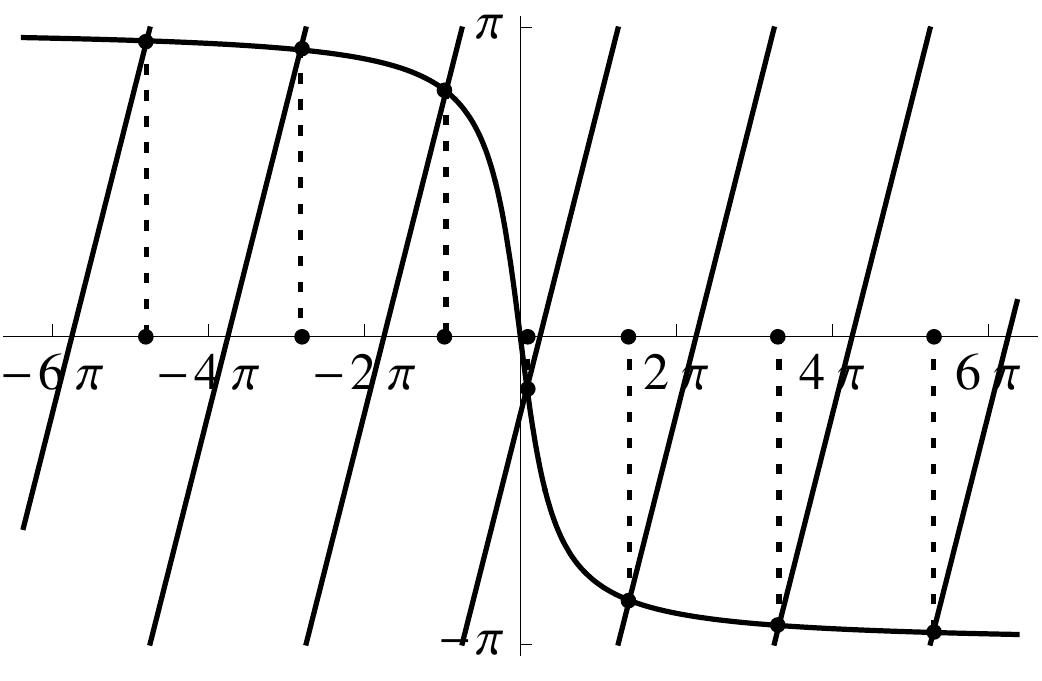}

\protect\caption{\label{fig:expExpSp}Fix $\theta=0.8$, $\Lambda_{\theta}=\left\{ \lambda_{n}\left(\theta\right)\right\} $
= intersections of two curves. (spectrum from curve intersections)}
\end{figure}

\begin{rem}
In (\ref{eq:saextA}), $h\in dom\left(D_{F}\right)$ (see (\ref{eq:D})),
and by \corref{dev1}, $h$ satisfies the boundary conditions (\ref{eq:dev6})-(\ref{eq:dev7}).
Also, by \lemref{dev1}, $\xi_{+}=F_{0}=e^{-x}$, $\xi_{-}=F_{1}=e^{x-1}$,
and $\left\Vert \xi_{+}\right\Vert _{\mathscr{H}_{F}}=\left\Vert \xi_{-}\right\Vert _{\mathscr{H}_{F}}=1$.
\end{rem}
\index{boundary condition}
\begin{prop}
\label{prop:dev1}Let $A_{\theta}$ be a selfadjoint extension of
$-iD$ as in (\ref{eq:saextA}). Then, 
\begin{equation}
\psi\left(1\right)+\psi'\left(1\right)=e^{i\theta}\left(\psi\left(0\right)-\psi'\left(0\right)\right),\;\forall\psi\in dom\left(A_{\theta}\right).\label{eq:dev-bd}
\end{equation}
\end{prop}
\begin{proof}
Any $\psi\in dom\left(A_{\theta}\right)$ has the decomposition 
\[
\psi\left(x\right)=h\left(x\right)+c\left(e^{-x}+e^{i\theta}e^{x-1}\right)
\]
where $h\in dom\left(D_{F}\right)$, and $c\in\mathbb{C}$. An application
of \corref{ptspec} gives 
\begin{align*}
\psi\left(1\right)+\psi'\left(1\right) & =\underset{=0\:\left(\text{by }\left(\ref{eq:dev7}\right)\right)}{\underbrace{h\left(1\right)+h'\left(1\right)}}+c\left(e^{-1}+e^{i\theta}\right)+c\left(-e^{-1}+e^{i\theta}\right)=2c\,e^{i\theta}\\
\psi\left(0\right)-\psi'\left(0\right) & =\underset{=0\;\left(\text{by }\left(\ref{eq:dev6}\right)\right)}{\underbrace{h\left(0\right)-h'\left(0\right)}}+c\left(1+e^{-1}e^{i\theta}\right)-c\left(-1+e^{-1}e^{i\theta}\right)=2c
\end{align*}
which is the assertion in (\ref{eq:dev-bd}).\end{proof}
\begin{cor}
\label{cor:spext}Let $A_{\theta}$ be a selfadjoint extension of
$-iD_{F}$ as in (\ref{eq:saextA}). Fix $\lambda\in\mathbb{R}$,
then $\lambda\in spec_{pt}\left(A_{\theta}\right)$ $\Longleftrightarrow$
$e_{\lambda}\left(x\right):=e^{i\lambda x}\in dom\left(A_{\theta}\right)$,
and $\lambda$ is a solution to the following equation: 
\begin{equation}
\lambda=\theta+\tan^{-1}\left(\frac{2\lambda}{\lambda^{2}-1}\right)+2n\pi,\;n\in\mathbb{Z}.\label{eq:dev9}
\end{equation}
\end{cor}
\begin{proof}
By assumption, $e^{i\lambda x}\in dom\left(A_{\theta}\right)$, so
$\exists h_{\lambda}\in dom\left(D_{F}\right)$, and $\exists c_{\lambda}\in\mathbb{C}$
such that 
\begin{equation}
e^{i\lambda x}=h_{\lambda}\left(x\right)+c_{\lambda}\left(e^{x}+e^{i\theta}e^{x-1}\right).\label{eq:dev8}
\end{equation}
Applying the boundary condition in \propref{dev1}, we have 
\[
e^{i\lambda}+i\lambda e^{i\lambda}=e^{i\theta}\left(1-i\lambda\right);\;\mbox{i.e.,}
\]
\begin{equation}
e^{i\lambda}=e^{i\theta}\frac{1-i\lambda}{1+i\lambda}=e^{i\theta}e^{i\arg\left(\frac{1-i\lambda}{1+i\lambda}\right)}\label{eq:dev8-1}
\end{equation}
where 
\[
\arg\left(\frac{1-i\lambda}{1+i\lambda}\right)=\tan^{-1}\left(\frac{2\lambda}{\lambda^{2}-1}\right)
\]
and (\ref{eq:dev9}) follows. For a discrete set of solutions, see
\figref{expExpSp}.
\end{proof}

\index{spectrum!discrete}
\begin{cor}
If $A_{\theta}\supset-iD_{F}$ is a selfadjoint extension in $\mathscr{H}_{F}$,
then
\begin{align*}
spect\left(A_{\theta}\right)= & \left\{ \lambda\in\mathbb{R}\:\big|\:e_{\lambda}\in\mathscr{H}_{F}\:\mbox{satisfying }\left(\ref{eq:dev-bd}\right)\right\} \\
= & \big\{\lambda\in\mathbb{R}\:\big|\:e_{\lambda}\in\mathscr{H}_{F},\:e_{\lambda}=h_{\lambda}+c_{\lambda}\left(e^{x}+e^{i\theta}e^{x-1}\right),\\
 & \quad h_{\lambda}\in dom\left(D_{F}\right),\;c_{\lambda}\in\mathbb{C}\big\}.
\end{align*}
\end{cor}
\begin{rem}
The corollary holds for \emph{all} continuous p.d. functions $F:\left(-a,a\right)\rightarrow\mathbb{C}$.\end{rem}
\begin{cor}
\label{cor:spdiscrete}All selfadjoint extensions $A_{\theta}\supset-iD_{F}$
have purely atomic spectrum; i.e., 
\begin{equation}
\Lambda_{\theta}:=spect\left(A_{\theta}\right)=\mbox{discrete subset in }\mathbb{R}.\label{eq:Aspect1}
\end{equation}
And for all $\lambda\in\Lambda_{\theta}$, 
\begin{equation}
ker\left(A_{\theta}-\lambda I_{\mathscr{H}_{F}}\right)=\mathbb{C}e_{\lambda},\;\mbox{where }e_{\lambda}\left(x\right)=e^{i\lambda x}\label{eq:Aspect2}
\end{equation}
i.e., all eigenvalues have multiplicity $1$. (The set $\Lambda_{\theta}$
will be denoted $\left\{ \lambda_{n}\left(\theta\right)\right\} _{n\in\mathbb{Z}}$
following \figref{expExpSp}. )\index{eigenvalue}\end{cor}
\begin{proof}
This follows by solving eq. (\ref{eq:dev9}).\end{proof}
\begin{cor}
Let $A$ be a selfadjoint extension of $-iD_{F}$ as before. Suppose
$\lambda_{1},\lambda_{2}\in spec\left(A\right)$, $\lambda_{1}\neq\lambda_{2}$,
then $e_{\lambda_{i}}\in\mathscr{H}_{F}$, $i=1,2$; and $\left\langle e_{\lambda_{1}},e_{\lambda_{2}}\right\rangle _{\mathscr{H}_{F}}=0$. \end{cor}
\begin{proof}
Let $\lambda_{1},\lambda_{2}$ be as in the statement, then 
\[
\left(\lambda_{1}-\lambda_{2}\right)\left\langle e_{\lambda_{1}},e_{\lambda_{2}}\right\rangle _{\mathscr{H}_{F}}=\left\langle Ae_{\lambda_{1}},e_{\lambda_{2}}\right\rangle _{\mathscr{H}_{F}}-\left\langle e_{\lambda_{1}},Ae_{\lambda_{2}}\right\rangle _{\mathscr{H}_{F}}=0;
\]
so since $\lambda_{1}-\lambda_{2}\neq0$, we get $\left\langle e_{\lambda_{1}},e_{\lambda_{2}}\right\rangle _{\mathscr{H}_{F}}=0$. 
\end{proof}
For explicit computations regarding these points, see also Corollaries
\ref{cor:expinner}, \ref{cor:exporg}, and \ref{cor:Fext} below.

\subsection{The Spectra of the s.a. Extensions $A_{\theta}\supset-iD_{F}$}

Let $F\left(x\right)=e^{-\left|x\right|}$, $\left|x\right|<1$. Define
$D_{F}$$\left(F_{\varphi}\right)=F_{\varphi'}$ as before, where
\begin{align*}
F_{\varphi}\left(x\right) & =\int_{0}^{1}\varphi\left(y\right)F\left(x-y\right)dy\\
 & =\int_{0}^{1}\varphi\left(y\right)e^{-\left|x-y\right|}dy,\;\forall\varphi\in C_{c}^{\infty}\left(0,1\right).
\end{align*}
And let $\mathscr{H}_{F}$ be the RKHS of $F$. 
\begin{lem}
\label{lem:spInner}For all $\varphi\in C_{c}^{\infty}\left(0,1\right)$,
and all $h,h''\in\mathscr{H}_{F}$, we have 
\begin{equation}
\left\langle F_{\varphi},h\right\rangle _{\mathscr{H}_{F}}=\left\langle F_{\varphi},\tfrac{1}{2}\left(h-h''\right)\right\rangle _{2}-\tfrac{1}{2}\left[W\right]_{0}^{1}\label{eq:sp1}
\end{equation}
where 
\begin{equation}
W=\det\begin{bmatrix}h & F_{\varphi}\\
h' & F_{\varphi'}
\end{bmatrix}.\label{eq:sp1-2}
\end{equation}
Setting $l:=F_{\varphi}$, we have 
\begin{equation}
\left[W\right]_{0}^{1}=-\overline{l}\left(1\right)\left(h\left(1\right)+h'\left(1\right)\right)-\overline{l}\left(0\right)\left(h\left(0\right)-h'\left(0\right)\right).\label{eq:sp1-3}
\end{equation}
\end{lem}
\begin{proof}
Note 
\begin{align*}
\left\langle F_{\varphi},h\right\rangle _{\mathscr{H}_{F}} & =\int_{0}^{1}\varphi\left(x\right)h\left(x\right)dx\quad(\mbox{reproducing property})\\
 & =\left\langle \tfrac{1}{2}\left(I-\left(\tfrac{d}{dx}\right)^{2}\right)F_{\varphi},h\right\rangle _{2}\\
 & =\left\langle F_{\varphi},\tfrac{1}{2}\left(h-h''\right)\right\rangle _{2}-\tfrac{1}{2}\left[W\right]_{0}^{1}.
\end{align*}
Set $l:=F_{\varphi}\in\mathscr{H}_{F}$, $\varphi\in C_{c}^{\infty}\left(0,1\right)$.
Recall the boundary condition in \corref{dev1}:
\[
l\left(0\right)-l'\left(0\right)=l\left(1\right)+l'\left(1\right)=0.
\]
Then 
\begin{align*}
\left[W\right]_{0}^{1} & =\left(\overline{l'}h-\overline{l}h'\right)\left(1\right)-\left(\overline{l'}h-\overline{l}h'\right)\left(0\right)\\
 & =-\overline{l}\left(1\right)h\left(1\right)-\overline{l}\left(1\right)h'\left(1\right)-\overline{l}\left(0\right)h\left(0\right)+\overline{l}\left(0\right)h'\left(0\right)\\
 & =-\overline{l}\left(1\right)\left(h\left(1\right)+h'\left(1\right)\right)-\overline{l}\left(0\right)\left(h\left(0\right)-h'\left(0\right)\right)
\end{align*}
which is (\ref{eq:sp1-3}).\end{proof}
\begin{cor}
\label{cor:elambda}$e_{\lambda}\in\mathscr{H}_{F}$, $\forall\lambda\in\mathbb{R}$.\end{cor}
\begin{proof}
By \thmref{HF}, we need the following estimate: $\exists C<\infty$
such that 
\begin{equation}
\left|\int_{0}^{1}\varphi\left(x\right)e_{\lambda}\left(x\right)dx\right|^{2}\leq C\left\Vert F_{\varphi}\right\Vert _{\mathscr{H}_{F}}^{2}.\label{eq:sp1-5}
\end{equation}
But
\begin{eqnarray*}
 &  & \int_{0}^{1}\varphi\left(x\right)e_{\lambda}\left(x\right)dx\\
 & = & \left\langle \tfrac{1}{2}\left(I-\left(\tfrac{d}{dx}\right)^{2}\right)F_{\varphi},e_{\lambda}\right\rangle _{2}\\
 & = & \left\langle F_{\varphi},\tfrac{1}{2}\left(e_{\lambda}-e_{\lambda}''\right)\right\rangle _{2}-\frac{1}{2}\left[W\right]_{0}^{1}\\
 & = & \tfrac{1}{2}\left(1+\lambda^{2}\right)\left\langle F_{\varphi},e_{\lambda}\right\rangle _{2}-\tfrac{1}{2}\left(-l\left(1\right)\left(1+i\lambda\right)e^{i\lambda}-l\left(0\right)\left(1-i\lambda\right)\right);
\end{eqnarray*}
see (\ref{eq:sp1})-(\ref{eq:sp1-3}). Here, $l:=F_{\varphi}$. 

It suffices to show

(i) $\exists C_{1}<\infty$ such that 
\[
\left|l\left(0\right)\right|^{2}\mbox{ and }\left|l\left(1\right)\right|^{2}\leq C_{1}\left\Vert F_{\varphi}\right\Vert _{\mathscr{H}_{F}}^{2}.
\]

(ii) $\exists C_{2}<\infty$ such that 
\[
\left|\left\langle F_{\varphi},e_{\lambda}\right\rangle _{2}\right|^{2}\leq C_{2}\left\Vert F_{\varphi}\right\Vert _{\mathscr{H}_{F}}^{2}.
\]

For (i), note that 
\begin{align*}
\left|l\left(0\right)\right| & =\left|\left\langle F_{0},l\right\rangle _{\mathscr{H}_{F}}\right|\leq\left\Vert F_{0}\right\Vert _{\mathscr{H}_{F}}\left\Vert l\right\Vert _{\mathscr{H}_{F}}=\left\Vert F_{0}\right\Vert _{\mathscr{H}_{F}}\left\Vert F_{\varphi}\right\Vert _{\mathscr{H}_{F}}\\
\left|l\left(1\right)\right| & =\left|\left\langle F_{1},l\right\rangle _{\mathscr{H}_{F}}\right|\leq\left\Vert F_{1}\right\Vert _{\mathscr{H}_{F}}\left\Vert l\right\Vert _{\mathscr{H}_{F}}=\left\Vert F_{1}\right\Vert _{\mathscr{H}_{F}}\left\Vert F_{\varphi}\right\Vert _{\mathscr{H}_{F}}
\end{align*}
and we have 
\begin{align*}
\left\Vert F_{0}\right\Vert _{\mathscr{H}_{F}} & =\left\Vert F_{1}\right\Vert _{\mathscr{H}_{F}}=1\\
\left\Vert l\right\Vert _{\mathscr{H}_{F}}^{2} & =\left\Vert F_{\varphi}\right\Vert _{\mathscr{H}_{F}}^{2}=\left\Vert T_{F}\varphi\right\Vert _{2}^{2}\leq\lambda_{1}^{2}\left\Vert \varphi\right\Vert _{2}^{2}<\infty
\end{align*}
where $\lambda_{1}$ is the top eigenvalue of the Mercer operator
$T_{F}$ (\lemref{mer-1}). 

For (ii),\index{eigenvalue} 
\begin{eqnarray*}
\left|\left\langle F_{\varphi},e_{\lambda}\right\rangle _{2}\right|^{2} & = & \left|\left\langle T_{F}\varphi,e_{\lambda}\right\rangle _{2}\right|^{2}\\
 & = & \left|\left\langle T_{F}^{1/2}\varphi,T_{F}^{1/2}e_{\lambda}\right\rangle _{2}\right|^{2}\\
 & \leq & \left\Vert T_{F}^{1/2}\varphi\right\Vert _{2}^{2}\left\Vert T_{F}^{1/2}e_{\lambda}\right\Vert _{2}^{2}\;\left(\mbox{by Cauchy-Schwarz}\right)\\
 & = & \left\langle \varphi,T_{F}\varphi\right\rangle _{2}\left\Vert T_{F}^{1/2}e_{\lambda}\right\Vert _{2}^{2}\\
 & \leq & \left\Vert F_{\varphi}\right\Vert _{\mathscr{H}_{F}}^{2}\left\Vert e_{\lambda}\right\Vert _{2}^{2}=\left\Vert F_{\varphi}\right\Vert _{\mathscr{H}_{F}}^{2};
\end{eqnarray*}
where we used the fact that $\left\Vert T_{F}^{1/2}e_{\lambda}\right\Vert _{2}^{2}\leq\lambda_{1}\left\Vert e_{\lambda}\right\Vert _{2}^{2}\leq1$,
since $\lambda_{1}<1$ = the right endpoint of the interval $\left[0,1\right]$
(see \lemref{mer-1}), and $\left\Vert e_{\lambda}\right\Vert _{2}=1$.

Therefore, the corollary follows.\end{proof}
\begin{cor}
\label{cor:HFinner2}For all $\lambda\in\mathbb{R}$, and all $F_{\varphi}$,
$\varphi\in C_{c}^{\infty}\left(0,1\right)$, we have 
\begin{align}
\left\langle F_{\varphi},e_{\lambda}\right\rangle _{\mathscr{H}_{F}} & =\tfrac{1}{2}\left(1+\lambda^{2}\right)\left\langle F_{\varphi},e_{\lambda}\right\rangle _{2}\label{eq:sp1-4}\\
 & \quad+\tfrac{1}{2}\left(\overline{l}\left(1\right)\left(1+i\lambda\right)e^{i\lambda}+\overline{l}\left(0\right)\left(1-i\lambda\right)\right).\nonumber 
\end{align}
\end{cor}
\begin{proof}
By \lemref{spInner}, 
\[
\left\langle F_{\varphi},e_{\lambda}\right\rangle _{\mathscr{H}_{F}}=\left\langle F_{\varphi},\tfrac{1}{2}\left(e_{\lambda}-e_{\lambda}''\right)\right\rangle _{2}-\tfrac{1}{2}\left[W\right]_{0}^{1}.
\]
where 
\[
\tfrac{1}{2}\left(e_{\lambda}-e_{\lambda}''\right)=\tfrac{1}{2}\left(1+\lambda^{2}\right)e_{\lambda};\;\mbox{and}
\]
\[
\left[W\right]_{0}^{1}\overset{\left(\ref{eq:sp1-3}\right)}{=}-\overline{l}\left(1\right)\left(1+i\lambda\right)e^{i\lambda}-\overline{l}\left(0\right)\left(1-i\lambda\right),\;l:=F_{\varphi}.
\]
\end{proof}
\begin{lem}
For all $F_{\varphi}$, $\varphi\in C_{c}^{\infty}\left(0,1\right)$,
and all $\lambda\in\mathbb{R}$, 
\begin{equation}
\left\langle F_{\varphi},e_{\lambda}\right\rangle _{\mathscr{H}_{F}}=\left\langle \varphi,e_{\lambda}\right\rangle _{2}.\label{eq:sp2}
\end{equation}
In particular, set $\lambda=0$, we get 
\begin{align*}
\left\langle F_{\varphi},\mathbf{1}\right\rangle _{\mathscr{H}_{F}} & =\int_{0}^{1}\varphi\left(x\right)dx=\frac{1}{2}\int_{0}^{1}\left(F_{\varphi}-F_{\varphi}''\right)\left(x\right)dx\\
 & =\frac{1}{2}\left(\left\langle F_{\varphi},\mathbf{1}\right\rangle _{2}-\left\langle F_{\varphi}'',\mathbf{1}\right\rangle _{2}\right)\\
 & \leq C\left\Vert F_{\varphi}\right\Vert _{\mathscr{H}}
\end{align*}
\end{lem}
\begin{proof}
Eq. (\ref{eq:sp2}) follows from basic fact of the Mercer operator.
See  \lemref{mer-1} and its corollaries. It suffices to note the
following estimate:
\begin{align*}
\int_{0}^{1}F''_{\varphi}\left(x\right)dx & =F'{}_{\varphi}\left(1\right)-F'_{\varphi}\left(0\right)\\
 & =-e^{-1}\int_{0}^{1}e^{y}\varphi\left(y\right)dy-\int_{0}^{1}e^{-y}\varphi\left(y\right)dy\\
 & =-F_{\varphi}\left(1\right)-F_{\varphi}\left(0\right)\leq2\left\Vert F_{\varphi}\right\Vert _{\mathscr{H}}.
\end{align*}
\end{proof}
\begin{cor}
\label{cor:elambda1}For all $\lambda\in\mathbb{R}$, 
\begin{equation}
\left\langle e_{\lambda},e_{\lambda}\right\rangle _{\mathscr{H}_{F}}=\frac{\lambda^{2}+3}{2}.\label{eq:sp3}
\end{equation}
\end{cor}
\begin{proof}
By \corref{HFinner2}, we see that 
\begin{align}
\left\langle F_{\varphi},e_{\lambda}\right\rangle _{\mathscr{H}_{F}} & =\frac{1}{2}\left(1+\lambda^{2}\right)\left\langle F_{\varphi},e_{\lambda}\right\rangle _{2}\nonumber \\
 & +\frac{1}{2}\left(\overline{l}\left(1\right)\left(1+i\lambda\right)e^{i\lambda}+\overline{l}\left(0\right)\left(1-i\lambda\right)\right);\;l:=F_{\varphi}.\label{eq:sp4-1}
\end{align}
Since $\left\{ F_{\varphi}:\varphi\in C_{c}^{\infty}\left(0,1\right)\right\} $
is dense in $\mathscr{H}_{F}$, $\exists F_{\varphi_{n}}\rightarrow e_{\lambda}$
in $\mathscr{H}_{F}$, so that
\begin{align*}
\left\langle F_{\varphi_{n}},e_{\lambda}\right\rangle _{\mathscr{H}_{F}}\rightarrow & \left\langle e_{\lambda},e_{\lambda}\right\rangle _{\mathscr{H}_{F}}\\
= & \frac{1}{2}\left(1+\lambda^{2}\right)+\frac{1}{2}\left(e^{-i\lambda}\left(1+i\lambda\right)e^{i\lambda}+\left(1-i\lambda\right)\right)\\
= & \frac{1}{2}\left(1+\lambda^{2}\right)+1=\frac{\lambda^{2}+3}{2}.
\end{align*}

The approximation is justified since all the terms in the RHS of (\ref{eq:sp4-1})
satisfy the estimate $\left|\cdots\right|^{2}\leq C\left\Vert F_{\varphi}\right\Vert _{\mathscr{H}_{F}}^{2}$.
See the proof of \corref{elambda} for details.
\end{proof}
Note \lemref{spInner} is equivalent to the following:
\begin{cor}
\label{cor:expinner}For all $h\in\mathscr{H}_{F}$, and all $k\in dom\left(T_{F}^{-1}\right)$,
i.e., $k\in\left\{ F_{\varphi}:\varphi\in C_{c}^{\infty}\left(0,1\right)\right\} $,
we have 
\begin{equation}
\left\langle h,k\right\rangle _{\mathscr{H}}=\frac{1}{2}\left(\left\langle h,k\right\rangle _{0}+\left\langle h',k'\right\rangle _{0}\right)+\frac{1}{2}\left(\overline{h\left(0\right)}k\left(0\right)+\overline{h\left(1\right)}k\left(1\right)\right)\label{eq:HFinner1-1}
\end{equation}
and eq. (\ref{eq:HFinner1-1}) extends to all $k\in\mathscr{H}_{F}$,
since $dom\left(T_{F}^{-1}\right)$ is dense in $\mathscr{H}_{F}$.\end{cor}
\begin{example}
Take $h=k=e_{\lambda}$, $\lambda\in\mathbb{R}$, then (\ref{eq:HFinner1-1})
gives 
\[
\left\langle e_{\lambda},e_{\lambda}\right\rangle _{\mathscr{H}}=\frac{1}{2}\left(1+\lambda^{2}\right)+\frac{1}{2}\left(1+1\right)=\frac{\lambda^{2}+3}{2}
\]
as in (\ref{eq:sp3}).\end{example}
\begin{cor}
\label{cor:exporg}Let $A_{\theta}\supset-iD$ be any selfadjoint
extension in $\mathscr{H}_{F}$. If $\lambda,\mu\in spect\left(A_{\theta}\right)$,
such that $\lambda\neq\mu$, then $\left\langle e_{\lambda},e_{\mu}\right\rangle _{\mathscr{H}_{F}}=0$. \end{cor}
\begin{proof}
It follows from (\ref{eq:HFinner1-1}) that 
\begin{align}
2\left\langle e_{\lambda},e_{\mu}\right\rangle _{\mathscr{H}} & =\left\langle e_{\lambda},e_{\mu}\right\rangle _{0}+\lambda\mu\left\langle e_{\lambda},e_{\mu}\right\rangle _{0}+\left(1+e^{i\left(\mu-\lambda\right)}\right)\nonumber \\
 & =\left(1+\lambda\mu\right)\left\langle e_{\lambda},e_{\mu}\right\rangle _{0}+\left(1+e^{i\left(\mu-\lambda\right)}\right)\nonumber \\
 & =\left(1+\lambda\mu\right)\frac{e^{i\left(\mu-\lambda\right)}-1}{i\left(\mu-\lambda\right)}+\left(1+e^{i\left(\mu-\lambda\right)}\right)\label{eq:HFinner-1-2}
\end{align}
By \corref{spext}, eq. (\ref{eq:dev8-1}), we have 
\[
e^{i\lambda}=\frac{1-i\lambda}{1+i\lambda}e^{i\theta},\quad e^{i\mu}=\frac{1-i\mu}{1+i\mu}e^{i\theta}
\]
and so 
\[
e^{i\left(\mu-\lambda\right)}=\frac{\left(1-i\mu\right)\left(1+i\lambda\right)}{\left(1+i\mu\right)\left(1-i\lambda\right)}.
\]
Substitute this into (\ref{eq:HFinner-1-2}) yields 
\[
2\left\langle e_{\lambda},e_{\mu}\right\rangle _{\mathscr{H}}=\frac{-2\left(1+\lambda\mu\right)}{\left(1+i\mu\right)\left(1-i\lambda\right)}+\frac{2\left(1+\lambda\mu\right)}{\left(1+i\mu\right)\left(1-i\lambda\right)}=0.
\]
\end{proof}
\begin{cor}
\label{cor:Fext}Let $F\left(x\right)=e^{-\left|x\right|}$, $\left|x\right|<1$.
Let $D_{F}\left(F_{\varphi}\right)=F_{\varphi'}$, $\forall\varphi\in C_{c}^{\infty}\left(0,1\right)$,
and $A_{\theta}\supset-iD_{F}$ be a selfadjoint extension  in $\mathscr{H}_{F}$.
Set $e_{\lambda}\left(x\right)=e^{i\lambda x}$, and 
\begin{equation}
\Lambda_{\theta}:=spect\left(A_{\theta}\right)\left(=\mbox{discrete subset in }\mathbb{R}\mbox{ by Cor. }\ref{cor:spdiscrete}\right)\label{eq:sp5}
\end{equation}
Then 
\begin{equation}
\widetilde{F}_{\theta}\left(x\right)=\sum_{\lambda\in\Lambda_{\theta}}\frac{2}{\lambda^{2}+3}e_{\lambda}\left(x\right),\;\forall x\in\mathbb{R}\label{eq:sp6}
\end{equation}
is a continuous p.d. extension of $F$ to the real line. Note that
both sides in eq. (\ref{eq:sp6}) depend on the choice of $\theta$.
\end{cor}
The type 1 extensions are indexed by $\theta\in[0,2\pi)$ where $\Lambda_{\theta}$
is given in (\ref{eq:sp5}), see also (\ref{eq:dev9}) in \corref{spext}.
\begin{cor}[Sampling property of the set $\Lambda_{\theta}$ ]
Let $F\left(x\right)=e^{-\left|x\right|}$ in $\left|x\right|<1$,
$\mathscr{H}_{F}$, $\theta$, and $\Lambda_{\theta}$ be as above.
Let $T_{F}$ be the corresponding Mercer operator. Then for all $\varphi\in L^{2}\left(0,1\right)$,
we have 
\[
\left(T_{F}\varphi\right)\left(x\right)=2\sum_{\lambda\in\Lambda_{\theta}}\frac{\widehat{\varphi}\left(\lambda\right)}{\lambda^{2}+3}e^{i\lambda x},\;\mbox{for all }x\in\left(0,1\right).
\]
\end{cor}
\begin{proof}
This is immediate from \corref{Fext}.\end{proof}
\begin{rem}
Note that the system $\left\{ e_{\lambda}\:|\:\lambda\in\Lambda_{\theta}\right\} $
is orthogonal in $\mathscr{H}_{F}$, but \emph{not} in $L^{2}\left(0,1\right)$.\index{Theorem!Mercer's-}\end{rem}
\begin{proof}
We saw that $A_{\theta}$ has pure atomic spectrum. By (\ref{cor:elambda1}),
the set 
\[
\left\{ \sqrt{\frac{2}{\lambda^{2}+3}}e_{\lambda}:\lambda\in\Lambda_{\theta}\right\} 
\]
is an ONB in $\mathscr{H}_{F}$. Hence, for $F=F_{0}=e^{-\left|x\right|}$,
we have the corresponding p.d. extension: 
\begin{align}
F_{\theta}\left(x\right) & =\sum_{\lambda\in\Lambda_{\theta}}\frac{1}{\left\Vert e_{\lambda}\right\Vert _{\mathscr{H}_{F}}^{2}}\left\langle e_{\lambda},F\right\rangle _{\mathscr{H}_{F}}e_{\lambda}\left(x\right)\nonumber \\
 & =\sum_{\lambda\in\Lambda_{\theta}}\frac{2}{\lambda^{2}+3}e_{\lambda}\left(x\right),\;\forall x\in\left[0,1\right].\label{eq:sp4-2}
\end{align}
where $\left\langle e_{\lambda},F\right\rangle _{\mathscr{H}_{F}}=\overline{e_{\lambda}\left(0\right)}=1$
by the reproducing property. But the RHS of (\ref{eq:sp4-2}) extends
to $\mathbb{R}$. See \figref{expExt}.
\end{proof}
\begin{figure}
\includegraphics[scale=0.6]{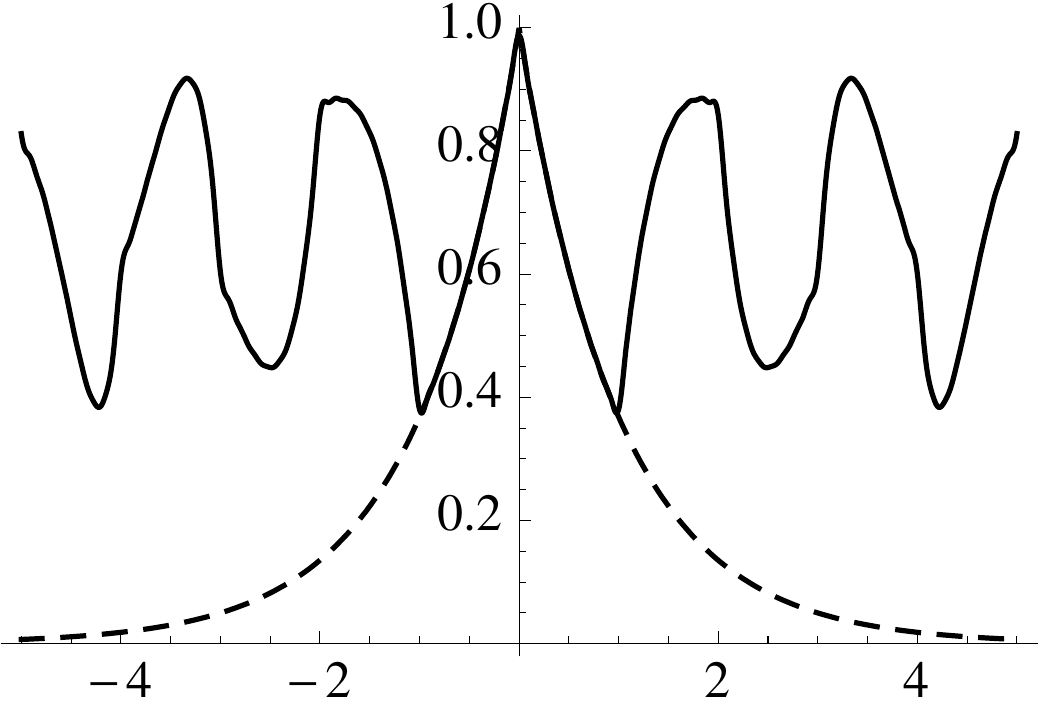}

\protect\caption{\label{fig:expExt}$\theta=0$. A type 1 continuous p.d. extension
of $F\left(x\right)=e^{-\left|x\right|}\big|_{\left[-1,1\right]}$
in $\mathscr{H}_{F}$.}
\end{figure}

\begin{cor}
\label{cor:exponb}Let $F\left(x\right)=e^{-\left|x\right|}$ in $\left(-1,1\right)$,
and let $\mathscr{H}_{F}$ be the RKHS. Let $\theta\in[0,2\pi)$,
and let $\Lambda_{\theta}$ be as above; then $\left\{ \sqrt{\frac{2}{\lambda^{2}+3}}e_{\lambda}\:|\:\lambda\in\Lambda_{\theta}\right\} $
is an ONB in $\mathscr{H}_{F}$.
\end{cor}

\section*{A summary of relevant numbers from the Reference List}

For readers wishing to follow up sources, or to go in more depth with
topics above, we suggest: 

The pioneering paper here is \cite{Aro50} and the intervening decades
have witnessed a host of applications. And by now there are books
dealing with various aspects of reproducing kernel Hilbert spaces
(RKHS). A more comprehensive citation list is: \cite{AD86,JPT14,Nus75,Ru63,MR1839648,MR2354721,AJSV13,Aro50,Nel59a,Sch64,SZ07,SZ09}.

\begin{appendix}

\part{Appendix}

\chapter{An overview of Functional Analysis books (cast of characters) \label{chap:books}}

\chaptermark{An overview of FA books}
\begin{quotation}
If people do not believe that mathematics is simple, it is only because
they do not realize how complicated life is. 

--- John von Neumann\sindex[nam]{von Neumann, J., (1903-1957)}\vspace{2em}
\end{quotation}
Below we offer a list of related Functional Analysis books; they cover
a host of diverse areas of functional analysis and applications, some
different from what we select here: Our comments are telegraphic-review
form (by P.J.):\\

Akhiezer and Glazman\emph{, ``Theory of linear operators in Hilbert
space''} \cite{AG93} 

-- a classic book set covering the detailed structure of unbounded
operators and their applications; and now in a lovely Dover edition.\\

\href{http://senate.universityofcalifornia.edu/news/images/arveson.jpg}{Arveson}\emph{,
``An invitation to $C^{*}$-algebras''} \cite{Arv76} 

-- an introduction to $C^{*}$-algebras and their representations
on Hilbert spaces. -- covers the most basic ideas, as simply and concretely
as we could. -- Hilbert spaces are separable and $C^{*}$-algebras
are GCR. Representations are given a concrete parametric description,
including the irreducible representations of any $C^{*}$-algebra,
even if not GCR. For someone interested in Borel structures, see Chapter
3. Chapter 1 is a bare-bones introduction to $C^{*}$-algebras.\\

Bachman and Narici, \emph{``Functional analysis''} \cite{MR1819613}

The book by Bachman and Narici's is a systematic introduction to the
fundamentals of functional analysis. It is easier to follow than say
Rudin\textquoteright s Functional Analysis book, but it doesn\textquoteright t
go as far either. Rather it helps readers reinforcing topics from
real analysis and other masters level courses. It serves to bridge
the gap between more difficult treatments of functional analysis.
(Dover reprints classics in a cheap paper back format.)\\

Bratteli and Robinson, \emph{``Operator algebras and quantum statistical
mechanics''} \cite{BR79,bratteli1981operator}

This is a widely cited two volume book-set, covering the theory of
operator algebras, and its applications to quantum statistical mechanics.
It is one of the more authoritative treatments of the many exciting
applications of functional analysis to physics. Both books are self-contained;
with complete proofs; -- a useful text for students with a prior exposure
to basic functional analysis. One of the main themes in v1 is decomposition
theory, and the use of Choquet simplices. Example: the set of KMS-states
for a $C^{*}$-algebraic dynamical system typically forms a Choquet
simplex\index{simplex}. An introductory chapter covers algebraic
techniques and their use in statistical physics; this is followed
up in v2. Indeed, a host of applications are covered in v2. The new
edition has a more comprehensive discussion of dissipative operators
and analytic elements; -- and it includes a positive resolution of
the question of whether maximal orthogonal probability measure on
the state space of algebra is automatically maximal among all the
probability measures on the space.\\

Conway, \emph{``A course in functional analysis''} \cite{Con90} 

-- a comprehensive introduction to functional analysis. The style
is formal and rigorous. -- is designed to be used in grad courses.
Through its eleven chapters, J. Conway masterfully wrote a beautiful
exposition of this core subject.\\

Dunford and Schwartz, \emph{``Linear operators''} \cite{DS88a,DS88b,DS88c}

This classic three-volume book set, the first Functional analysis,
and the second the theory of linear operators. And for the theory
of unbounded operators it is unsurpassed. -- written by two notable
mathematicians, it constitutes a comprehensive survey of the general
theory of linear operations, and their diverse applications. Dunford
and Schwartz are influenced by von Neumann, and they emphasize the
significance of the relationships between the abstract theory and
its applications. The two first volumes are for the students. -- treatment
is relatively self-contained. Now a paperback edition of the original
work, unabridged, in three volumes.\\

Kolmogorov and Fomin, \emph{``Introductory real analysis''} \cite{MR0377445}

This book is two books bound as one; and in the lovely format from
Dover. Part 1: metric spaces, and normed linear spaces. Part 2: Lebesgue
integration and basic functional analysis. Numerous examples are sprinkled
through the text. To get the most out of this book, it helps if you
have already seen many of the results presented elsewhere. History:
The book came from original notes from Andrei Kolmogorov's lectures
given at Moscow's Lomonosov University in the 1940's, and it still
stands as timely introduction to real and functional analysis. Strengths:
step by step presentation of all the key concepts needed in the subject;
proceeding all the way from set theory to Fredholm integral equations.
Offers a wonderful and refreshing insight. Contents (sample): Elements
of Set Theory; Metric and Topological Spaces; Normed and Topological
Linear Spaces; Linear Functionals and Linear Operators; Elements of
Differential Calculus in Linear Spaces; Measure, Measurable Functions,
Integral; Indefinite Lebesgue Integral, Differentiation Theory; Spaces
of Summable Functions; Trigonometric Series, Fourier Transformation;
Linear Integral Equations.\\

Kadison and Ringrose, \emph{``Fundamentals of the theory of operator
algebras''} \cite{MR1468229,MR1468230}

Here we cite the first two volumes in a 4-volume book-set. It begins
with the fundamentals in functional analysis, and it aims at a systematic
presentations of the main areas in the theory of operator algebras,
both $C^{*}$-algebras, von Neumann algebras, and their applications,
so including subfactors, Tomita-Takesaki theory, spectral theory,
decomposition theory, and applications to ergodic theory, representations
of groups, and to mathematical physics.\\

Lax, \emph{``Functional analysis''} \cite{MR1892228} 

The subject of functional analysis, while fundamental and central
in the landscape of mathematics, really started with seminal theorems
due to Banach, Hilbert, von Neumann, Herglotz, Hausdorff, Friedrichs,
Steinhouse,...and many other of, the perhaps less well known, founding
fathers, in Central Europe (at the time), in the period between the
two World Wars. It gained from there because of its many success stories,
-- in proving new theorems, in unifying old ones, in offering a framework
for quantum theory, for dynamical systems, and for partial differential
equations. The Journal of Functional Analysis, starting in the 1960ties,
broadened the subject, reaching almost all branches of science, and
finding functional analytic flavor in theories surprisingly far from
the original roots of the subject. Peter Lax has himself, -- alone
and with others, shaped some of greatest successes of the period,
right up to the present. That is in the book!! And it offers an upbeat
outlook for the future. It has been tested in the class room, -- it
is really user-friendly. At the end of each chapter P. Lax offers
personal recollections; -- little known stories of how several of
the pioneers in the subject have been victims, -- in the 30ties and
the 40ties, of Nazi atrocities. The writing is crisp and engaged.\\

MacCluer, \emph{``Elementary functional analysis''} \cite{MR2462971} 

I received extremely positive student-feedback on MacCluer\textquoteright s
very nice book. It covers elementary functional analysis, is great
for self-study, and easy to follow. It conveys the author\textquoteright s
enthusiasm for her subject. It includes apposite quotes, anecdotes,
and historical asides, all making for a wonderful personal touch and
drawing the reader into dialogue in a palpable way. Contents: six
chapters, each introduced by a well-chosen quote, often hinting in
a very useful manner at the material that is to follow. I particularly
like MacCluer\textquoteright s choice of Dunford and Schwartz to start
off her third chapter: ``In linear spaces with a suitable topology
one encounters three far-reaching principles concerning continuous
linear transformations\dots '' We find out quickly that these ``Big
Three'' (as the chapter is titled) are uniform boundedness, the open
mapping theorem, and Hahn-Banach. MacCluer quickly goes on to cover
these three gems in a most effective and elegant manner, as well as
a number of their corollaries or, in her words, ``close cousins,''
such as the closed graph theorem and Banach-Steinhaus. The book takes
the reader from Hilbert space preliminaries to Banach- and $C^{*}$-algebras
and, to the spectral theorem.\\

Nelson, \emph{``Topics in Dynamics I: Flows'' }\cite{Ne69}

This is a book in the Princeton Math Lecture Notes series, appearing
first in 1972, but since Prof Nelson kindly made it available on his
website. In our opinion, it is the best account of general multiplicity
for normal operators, bounded and unbounded, and for abelian {*}-algebras.
In addition it contains a number of applications of functional analysis
to geometry and to physics.\\

Riesz et al., \emph{``Functional analysis''} \cite{MR1068530} 

A pioneering book in F.A., first published in the early 50s, and now
in a Dover edition, very readable. The book starts with an example
of a continuous function which is not differentiable and then proves
Lebesgue's theorem which tells you when a function does have a derivative.
The 2nd part of the book is about integral equations which again starts
with some examples of problems from the 19th century mathematicians.
The presentation of Fredholm's method is a gem.\\

Rudin, \emph{``Functional analysis''} \cite{Rud73}

``Modern analysis'' used to be a popular name for the subject of
this lovely book. It is as important as ever, but perhaps less ``modern''.
The subject of functional analysis, while fundamental and central
in the landscape of mathematics, really started with seminal theorems
due to Banach, Hilbert, von Neumann, Herglotz, Hausdorff, Friedrichs,
Steinhouse,...and many other of, the perhaps less well known, founding
fathers, in Central Europe (at the time), in the period between the
two World Wars. In the beginning it generated awe in its ability to
provide elegant proofs of classical theorems that otherwise were thought
to be both technical and difficult. The beautiful idea that makes
it all clear as daylight: Wiener's theorem on absolutely convergent
(AC) Fourier series of $1/f$ if you can divide, and if $f$ has AC
Fourier series, is a case in point. The new subject gained from there
because of its many success stories, -- in proving new theorems, in
unifying old ones, in offering a framework for quantum theory, for
dynamical systems, and for partial differential equations. And offering
a language that facilitated interdisciplinary work in science! The
topics in Rudin's book are inspired by harmonic\index{harmonic} analysis.
The later part offers one of the most elegant compact treatment of
the theory of operators in Hilbert space, I can think of. Its approach
to unbounded operators is lovely.\\

Sakai, \emph{``$C^{*}$-algebras and $W^{*}$-algebras''} \cite{MR0442701} 

The presentation is succinct, theorem, proof, ... qed; but this lovely
book had a profound influence on the subject. It\textquoteright s
scope cover nearly all major results in the subject up until that
time. In order to accomplish this goal (without expanding into multiple
volumes), the author omits examples, motivation,\dots{} . It is for
students who already have an interest in operator theory. As a student,
myself (PJ), I learned a lot from this wonderful book. \\

Shilov, \emph{``Elementary functional analysis''} \cite{MR1375236} 

Elementary Functional Analysis by Georgi E. Shilov is suitable for
a beginning course in functional analysis and some of its applications,
e.g., to Fourier series, to harmonic analysis, to partial differential
equations (PDEs), to Sobolev spaces, and it is a good supplement and
complement to two other popular books in the subject, one by Rudin,
and another by Edwards. Rudin's book is entitled ``Functional Analysis''
includes new material on unbounded operators in Hilbert space. Edwards'
book ``Functional Analysis: Theory and Applications;'' is in the
Dover series, and it is twice as thick as Shilov's book. Topics covered
in Shilov: Function spaces, $L^{p}$-spaces, Hilbert spaces, and linear
operators; the standard Banach, and Hahn-Banach theorems. It includes
many exercises and examples. Well motivated with applications. Book
Comparison: Shilov book is gentler on students, and it is probably
easier to get started with: It stresses motivation a bit more, the
exercises are easier, and finally Shilov includes a few applications;
fashionable these days. \\

Stein et al., \emph{``Functional analysis''} \cite{MR2827930} 

This book is the fourth book in a series: Elias Stein's and Rami Shakarchi's
Princeton lectures in analysis. Elias Stein is a world authority on
harmonic analysis. The book is of more recent vintage than the others
from our present list. The book on functional analysis is actually
quite different from other texts in functional analysis. For instance
Rudin's textbook on functional analysis has quite a different emphasis
from Stein's. Stein devotes a whole chapter to applications of the
Baire category theory while Rudin devotes a page. Stein does this
because it provides some insights into establishing the existence
of a continuous but nowhere differentiable function as well as the
existence of a continuous function with Fourier series diverging a
point. A special touch in Stein: Inclusion of Brownian motion, and
of process with independent increments, a la Doob's. Stein's approach
to the construction of Brownian motion is different and closer to
the approaches taken in books on financial math. Stein et al develop
Brownian motion in the context of solving Dirichlet's problem.\\

Stone, \emph{``Linear Transformations in Hilbert Space and Their
Applications to Analysis'' }\cite{Sto90}

Stone's book is a classic, came out in 1932, and was the unique source
on spectral multiplicity, and a host of applications of the theory
of unbounded operators to analysis, to approximation theory, and to
special functions. The last two chapters illustrate the theory with
a systematic study of (infinite $\times$ infinite) Jacobi matrices;
i.e., tri-diagonal infinite matrices; assumed formally selfadjoint
(i.e., Hermitian). Sample results: A dichotomy: Their von Neumann
indices must be $(0,0)$ or $(1,1)$. Some of the first known criteria
for when they are one or the other are given; plus a number of applications
to classical analysis.\index{matrix!Jacobi} \\

Takesaki, \emph{``Theory of operator algebras''} \cite{MR548728} 

-- written by one of the most prominent researchers of the area, provides
an introduction to this rapidly developing theory. ... These books
are recommended to every graduate student interested in this exciting
branch of mathematics. Furthermore, they should be on the bookshelf
of every researcher of the area.\\

Trèves, \emph{``Topological vector spaces, distributions and kernels''}
\cite{MR2296978}

Covers topological vector spaces and their applications, and it is
a pioneering book. It is antidote for those who mistakenly believe
that functional analysis is about Banach and Hilbert spaces. It's
also about Fréchet spaces, LF spaces, Schwartz distributions (generalized
functions), nuclear spaces, tensor products, and the Schwartz Kernel
Theorem (proved by Grothendieck). Trèves's book provides the perfect
background for advanced work in linear differential, pseudodifferential,
or Fourier integral operators.\\

Yosida \emph{``Functional analysis''} \cite{Yos95}

Yosida\textquoteright s book is based on lectures given decades ago
at the University of Tokyo. It is intended as a textbook to be studied
by students on their own or to be used in a course on Functional Analysis,
i.e., the general theory of linear operators in function spaces together
with salient features of its application to diverse fields of modern
and classical analysis. Necessary prerequisites for the reading of
this book are summarized, with or without proof, in Chapter 0 under
titles: Set Theory, Topological Spaces, Measure Spaces and Linear
Spaces. Then, starting with the chapter on Semi-norms, a general theory
of Banach and Hilbert spaces is presented in connection with the theory
of generalized functions of S.L. Sobolev and L. Schwartz. The reader
may pass, e.g., from Chapter IX (Analytical Theory of Semi-groups)
directly to Chapter XIII (Ergodic Theory and Diffusion Theory) and
to Chapter XIV (Integration of the Equation of Evolution). Such materials
as ``Weak Topologies and Duality in Locally Convex Spaces'' and
``Nuclear Spaces'' are presented in the form of the appendices to
Chapter V and Chapter X, respectively.\\
\\

\noindent \textbf{Some relevant books: Classics, and in the Dover
series:}\\

Banach, \emph{``Theory of Linear Operations'' }\cite{MR1357166}

Georgi, \emph{``Weak Interactions and Modern Particle Theory'' }\cite{georgi2009weak}

Prenter, \emph{``Splines and Variational Methods'' }\cite{MR1013116}\\
\\

\index{Schrödinger equation}

In chapters \ref{chap:GNS} and \ref{chap:KS} above we have cited
pioneers in quantum physics, the foundations of quantum mechanics.
The most central here are Heisenberg (matrix mechanics), Schrödinger
(wave mechanics, the Schrödinger equation), and Dirac (Dirac\textquoteright s
equation is a relativistic wave equation, describes all spin-\textonehalf{}
massive particles free form, as well as electromagnetic interactions).
We further sketched von Neumann\textquoteright s discovery of the
equivalence of the answers given by Heisenberg and Schrödinger, and
the Stone-von Neumann uniqueness theorem. The relevant papers and
books are as follows: \cite{Hei69,Sch32,vN31,HN28,Dir35,Dir47}.

\begin{table}
\[ 
\xymatrix{
& \txt{Functional Analysis}\ar[ld]\ar[rd]\ar@{<->}[ddd] \\ 
\txt{Linear\\Operators (A)}\ar@{<->}[d]\ar@{<->}[ddr] &  & 
\txt{Mathematical\\Physics (C)}\ar@{<->}[d]\ar@{<->}[ddl]\\
\txt{Harmonic\\Analysis (B)}\ar@{<->}[dr] &  & 
\txt{Representation\\Theory (D)}\ar@{<->}[dl]\\  & 
\txt{Probability Theory / Statistics (E)} } \]

\protect\caption{}
\end{table}

\begin{table}[H]
\begin{tabular}{cl|l}
\multicolumn{3}{c}{(A)}\tabularnewline
\hline 
\multirow{2}{*}{$\Big\updownarrow$} & bounded & differential operators, ODE/PDE\tabularnewline
 & unbounded & generators of diffusion\tabularnewline
 &  & \tabularnewline
\multirow{2}{*}{$\Big\updownarrow$} & geometry & \tabularnewline
 & spectral theory & Schrödinger operators\tabularnewline
 & spectral representation  & wave operators\tabularnewline
 & single operators & scattering operators\tabularnewline
 & system of operators & \tabularnewline
 & operator commutation relations & \tabularnewline
 &  & \tabularnewline
\end{tabular}

\protect\caption{}
\end{table}

\begin{table}[H]
\begin{tabular}{cc|c}
\multicolumn{3}{c}{(B)}\tabularnewline
\hline 
 &  & analysis / synthesis\tabularnewline
 &  & Fourier analysis, wavelet analysis\tabularnewline
\multirow{2}{*}{$\Big\updownarrow$} & commutative & \tabularnewline
\cline{3-3} 
 & non-commutative & Applications:\tabularnewline
 &  & signal processing\tabularnewline
 &  & physics\tabularnewline
 &  & statistics\tabularnewline
 &  & analysis on fractals\tabularnewline
\end{tabular}

\protect\caption{}
\end{table}

\begin{table}[H]
\begin{tabular}{cc|c}
\multicolumn{3}{c}{(C)}\tabularnewline
\hline 
 &  & \tabularnewline
\multirow{2}{*}{$\Big\updownarrow$} & quantum physics & \tabularnewline
 & classical mechanics & \tabularnewline
 &  & quantum information\tabularnewline
\multirow{2}{*}{$\Big\updownarrow$} & statistical physics & states and decomposition\tabularnewline
 & quantum field theory & equilibrium: Gibbs, KMS, ...\tabularnewline
 &  & \tabularnewline
\multirow{2}{*}{$\Big\updownarrow$} & relativistic & \tabularnewline
 & non-relativistic & \tabularnewline
\end{tabular}

\protect\caption{}
\end{table}

\begin{table}[H]
\begin{tabular}{cc|c}
\multicolumn{3}{c}{(D)}\tabularnewline
\hline 
 &  & \tabularnewline
\multirow{3}{*}{$\Bigg\updownarrow$} & groups (abelian, non-abelian) & locally compact, non-locally compact\tabularnewline
 & algebras & \tabularnewline
 & generators and relations & Lie groups $\longleftrightarrow$ Lie algebras\tabularnewline
\multirow{2}{*}{} &  & induced representations\tabularnewline
 &  & decomposition of representations\tabularnewline
 &  & groups over $\mathbb{R}$, $\mathbb{C}$, or other local fields\tabularnewline
\end{tabular}

\protect\caption{}
\end{table}

\begin{table}[H]
\begin{tabular}{c|>{\centering}p{0.6\columnwidth}}
\multicolumn{2}{c}{(E)}\tabularnewline
\hline 
 & discrete\tabularnewline
 & continuous\tabularnewline
 & Gaussian \tabularnewline
 & Brownian motion\tabularnewline
stochastic processes & non-Gaussian\tabularnewline
 & Lévy\tabularnewline
\cline{2-2} 
 & solutions of diffusion equations with the use of functional integrals
(i.e., probability measure on infinite-dimensional spaces such as
$C\left(\mathbb{R}\right)$ or Schwartz space $\mathcal{S}$)\tabularnewline
\end{tabular}

\protect\caption{}
\end{table}

\chapter{Terminology from neighboring areas\label{chap:term}}

\textbf{Classical Wiener measure/space.} Classical Wiener space (named
Norbert Wiener) is the sample-space part of a probability space, a
triple (sample-space, sigma-algebra, and probability measure). The
sample space may be taken to be the collection of all continuous functions
on a given domain (usually an interval for time). So sample paths
are continuous functions. The sigma-algebra is generated by cylinder-sets,
and the probability measure is called the Wiener measure; (its construction
is subtle, see \chapref{bm} above.) It has the property that the
stochastic processes which samples the paths in the model is a Gaussian
process with independent increments, the so called Brownian motion.
It is also called the Wiener process. And it should perhaps be named
after L. Bachelier, whose work predates that of Einstein.

It, and the related process \textquotedblleft white noise\textquotedblright ,
are important in pure and applied mathematics. It is a core ingredient
in stochastic analysis: the study of stochastic calculus, diffusion
processes, and potential theory. Applications include engineering
and physics: models of noise in electronics engineering, in instrument
errors in filtering theory, and in control theory. In atomic physics,
it is used in the study of diffusion, the Fokker\textendash Planck
and Langevin equations; and in path-space integrals; the Feynman\textendash Kac
formula, in the solution of the Schrödinger equation. In finance,
it is used in the solution to the Black\textendash Scholes equation
for option prices. References include \cite{MR2966130, MR3050315, AAR13, CW14, MR887102, MR0265548, Hi80, Ito06, MR2254502, Jor14, MR0161189, Nel67, Sch32, MR0107312, MR2883397};
and we refer to Sections \ref{sec:Hilbert}, \ref{sec:pspace}, \ref{sec:dbm},
and \ref{sec:stoch}.\textbf{}\\
\textbf{}\\
\textbf{Hilbert's sixth problem.} This is not a \textquotedblleft yes/no
problem\textquotedblright ; rather the 6th asks for a mathematical
axiomatization of physics. In a common English translation, it reads:
6. Give a Mathematical Treatment of the Axioms of Physics. A parallel
is drawn to the foundations of geometry: \emph{To treat in the same
manner, by means of axioms, those physical sciences where mathematics
plays an important part; in the first rank are the theory of probabilities
and mechanics.}\index{Hilbert's six problem}\index{axioms}

Hilbert: \textquotedbl{}As to the axioms of the theory of probabilities,
it seems to me desirable that their logical investigation should be
accompanied by a rigorous and satisfactory development of the method
of mean values in mathematical physics, and in particular in the kinetic
theory of gases. ... Boltzmann's work on the principles of mechanics
suggests the problem of developing mathematically the limiting processes,
there merely indicated, which lead from the atomistic view to the
laws of motion of continua.\textquotedbl{}

In the 1930s, probability theory was put on an axiomatic and sound
foundation by Andrey Kolmogorov. In the 1960s, we have the work of
A. Wightman, R. Haag, J. Glimm, and A. Jaffe on quantum field theory.
This was followed by the Standard Model in particle physics and general
relativity. Still unresolved is the theory of quantum gravity.\index{quantum field theory}\index{Standard Model}\textbf{}\\
Ref. Sections \ref{sec:app} (pg. \pageref{sec:app}), \ref{sub:QM}
(pg. \pageref{sub:QM}), \ref{sec:pspace} (pg. \pageref{sec:pspace});
see also \cite{Wig76,MR3246256}.\textbf{ }\\
\textbf{}\\
\textbf{Quantum mechanics (QM); (quantum physics, or quantum theory)}
is a branch of physics which describes physical phenomena at \textquotedblleft small\textquotedblright{}
scales, atomic and subatomic length scales. The action is on the order
of the Planck constant. QM deals with observation of physical quantities
that can only change and interact, by discrete amounts or \textquotedblleft steps\textquotedblright{}
(hence \textquotedblleft quantum\textquotedblright ), and behave probabilistically
rather than deterministically. The \textquotedbl{}steps\textquotedbl{}
are too tiny even for microscopes. Any description must be given in
terms of a wave function, as opposed to particles. (For details, see
\cite{Dir47,bratteli1981operator}.)\\
Ref. Sections \ref{sub:QM} (pg. \pageref{sub:QM}), \ref{sec:dga}
(pg. \pageref{sec:dga}), \ref{sec:qmview} (pg. \pageref{sec:qmview}),
and Chapter \ref{chap:KS} (pg. \pageref{chap:KS}).\\
\\
\textbf{Quantum field theory (QFT)} is a mathematical framework used
in physics for constructing models of subatomic particles (quantum
mechanical). It covers such areas as particle physics and condensed
matter physics. A QFT treats particle-wave duality as excited states
of an underlying physical field, called field quanta. Of interest
are quantum mechanical interactions between particles and the corresponding
underlying fields. (For details, see \cite{MR887102}.)\\
Ref. \secref{app} (pg. \pageref{sec:app}).\\
\\
\textbf{Signal processing (SP)} is an engineering discipline dealing
with transmission of signals (information, speech or images, over
wires, or wireless. An important tool in the area involves subdivision
of time signals into frequency bands, and it involves effective algorithms
for implementations of processing or transferring information contained
in a variety of different symbolic, or abstract formats broadly designated
as signals. SP uses mathematical, statistical, computational tools.
(For details, see \cite{MR0042667,BJ02}.)\\
Ref. Section \ref{sec:Comments} (pg. \pageref{sec:Comments}), Chapters
\ref{chap:cp} (pg. \pageref{chap:cp}), \ref{chap:KS} (pg. \pageref{chap:KS}).\index{signal processing}\\
\\
\textbf{Stochastic processes (SP), or random process}, are part of
probability theory. They are used when deterministic quantities are
not feasible: random variables are measures in their respective probability
distributions (also called \textquotedblleft laws.\textquotedblright )
A SP is an indexed family of random variables (representing measurements,
or samples), for example, if a SP is indexed by time, it represents
the evolution or dynamics of some system. A SP is the probabilistic
counterpart to a deterministic process (or a deterministic system).
An example is Brownian motion (BM), the random motion of particles
(e.g., pollen) suspended in a fluid (a liquid or a gas). BM results
from their collision of the pollen with atoms or molecules making
up the gas or liquid. BM also refer to the mathematical model used
to describe such random movements. (For details, see \cite{MR0184066,MR2053326,Ito06,MR0107312,MR2883397}.)\\
Ref. Sections \ref{sub:dmeas} (pg. \pageref{sub:dmeas}), \ref{sec:stoch}
(pg. \pageref{sec:stoch}); and Chapters \ref{chap:bm} (pg. \pageref{chap:bm}),
\ref{chap:groups} (pg. \pageref{chap:groups}).\\
\\
\textbf{Unitary representations (UR)} of a groups $G$ are homomorphisms
from the group $G$ in question into the group of all unitary operators
in some Hilbert space; the Hilbert space depending on the UR. The
theory is best understood in in the case of strongly continuous URs
of locally compact (topological) groups. Applications include quantum
mechanics. Books by Hermann Weyl, Pontryagin, and George Mackey have
influenced our presentation. The theory of UR is closely connected
with harmonic analysis; -- for non-commutative groups, non-Abelian
harmonic analysis. Important groups in physics are non-commutative,
so this case is extremely important, although it is also rather technical.
There is a vast literature, though. Important papers are cited in
the books by George Mackey. In fact, versions of the Plancherel theorem
exist for some non-commutative Lie groups (of direct relevance to
physics), but they are subtle, and the non-commutative analysis is
carried out on a case-by-case basis. The best known special case it
that of compact groups where we have the Peter-Weyl theorem. But the
important symmetry groups in relativistic physics are non-compact,
see \cite{MR887102}\textbf{. }\\
Ref. Sections \ref{sec:Hilbert} (pg. \pageref{sec:Hilbert}), \ref{sec:dga}
(pg. \pageref{sec:dga}), and Chapters \ref{chap:GNS} (pg. \pageref{chap:GNS}),
\ref{chap:cp} (pg. \pageref{chap:cp}), \ref{chap:groups} (pg. \pageref{chap:groups}).\\
\textbf{}\\
\textbf{Wavelets} are wave-like functions; they can typically be visualized
as \textquotedbl{}brief oscillations\textquotedbl{} as one might see
recorded in seismographs, or in heart monitors. What is special about
wavelet functions is that they allow for effective algorithmic construction
of bases in a variety of function spaces. The algorithms in turn are
based on a notion of resolution and scale-similarity. The last two
features make wavelet decompositions more powerful than comparable
Fourier analyses. Wavelets can be localized, while Fourier bases cannot.
Wavelets are designed to have specific properties that make them of
practical use in signal processing. (For details, see \cite{BJ02}.)\\
Ref. Sections \ref{sub:ONB} (pg. \pageref{sub:ONB}), \ref{sec:KMil}
(pg. \pageref{sec:KMil}), \ref{sec:Comments} (pg. \pageref{sec:Comments}),
and Chapter \ref{chap:cp} (pg. \pageref{chap:cp}).

\chapter{Often cited above\label{chap:bios}}
\begin{quotation}
\textbf{mathematical ideas originate in empirics.} But, once they
are conceived, the subject begins to live a peculiar life of its own
and is \dots{} governed by almost entirely aesthetical motivations.
In other words, at a great distance from its empirical source, or
after much ``abstract'' inbreeding, a mathematical subject is in
danger of degeneration. Whenever this stage is reached the only remedy
seems to me to be the rejuvenating return to the source: the reinjection
of more or less directly empirical ideas. 

--- von Neumann\sindex[nam]{von Neumann, J., (1903-1957)}\vspace{2em}
\end{quotation}
Inside the book, the following authors are cited frequently, W. Arveson,
M. Atiyah, L. Bachelier, S. Banach, H. Bohr, M. Born, N. Bohr, P.
Dirac, W. Döblin, F. Dyson, K. Friedrichs, I. Gelfand, L. Gårding,
I. Gelfand, W. Heisenberg, D. Hilbert, K. It\={o}, R. Kadison, S.
Kakutani, M. Krein, P. Lax, G. Mackey, E. Nelson, R. Phillips, F.
Riesz, M. Riesz, E. Schrödinger, H.A. Schwarz, L. Schwartz, J. Schwartz,
I. Segal, I. Singer, M. Stone, J. von Neumann, N. Wiener. Below a
short bio. \medskip{}

\textbf{\href{http://senate.universityofcalifornia.edu/news/images/arveson.jpg}{William Arveson}}
(1934 -- 2011) \cite{Arv72,Arv98}. Cited in connection with $C^{*}$-algebras
and their states and representations.

W. Arveson, known for his work on completely positive maps, and their
extensions; powerful generalizations of the ideas of Banach, Krein,
and Stinespring . An early results in this area is an extension theorem
for completely positive maps with values in the algebra of all bounded
operators. This theorem led to injectivity of von-Neumann algebras
in general, and work by Alain Connes relating injectivity to hyperfiniteness.
In a series of papers in the 60's and 70's, Arveson introduced non-commutative
analogues of several concepts from classical harmonic analysis including
the Shilov and. Choquet boundaries. \medskip{}

\textbf{\href{http://upload.wikimedia.org/wikipedia/commons/thumb/a/af/Michael_Francis_Atiyah.jpg/220px-Michael_Francis_Atiyah.jpg}{Sir Michael Francis Atiyah}}
(1929 -- ). Of the Atiyah-Singer Index Theorem. The Atiyah-Singer
index of a partial differential operator (PDO) is related to the Fredholm
index; -- it equates an index, i.e., the difference of the number
of independent solutions of two geometric, homogeneous PDEs (one for
the operator and the other for its adjoint) to an associated list
of invariants in differential geometry. It applies to many problems
in mathematics after they are translated into the problem of finding
the number of independent solutions of some PDE. The Atiyah\textendash Singer
index theorem gives a formula for the index of certain differential
operators, in terms of geometric and topological invariants.\index{index!Fredholm-}

The Hirzebruch-Riemann-Roch theorem is a special cases of the Atiyah\textendash Singer
index theorem. In fact the index theorem gave a more powerful result,
because its proof applied to all compact complex manifolds, while
Hirzebruch's proof only worked for projective manifolds.

Related: In 1959 by Gelfand noticed homotopy invariance via an index,
and he asked for more general formulas for topological invariants.
For spin manifolds, Atiyah suggested that integrality could be explained
as an index of a Dirac operator (Atiyah and Singer, 1961).\medskip{}

\textbf{\href{https://en.wikipedia.org/wiki/Louis_Bachelier\#/media/File:LouisBachelier.jpg}{Louis Bachelier}}
(1870 -- 1946), a French probabilist, is credited with being the inventor
of the stochastic process, now called Brownian motion; it was part
of his PhD thesis, The Theory of Speculation, (1900). It discusses
use of random walks, and Brownian motion, to evaluate stock options,
and it is considered the first paper in mathematical finance. Even
though Bachelier's work was more mathematical, and predates Einstein's
Brownian motion paper by five years, it didn\textquoteright t receive
much attention at the time, and it was only \textquotedblleft discovered\textquotedblright{}
much later by the MIT economists Paul Samuelson, in the 1960ties.

\medskip{}

\textbf{\href{http://kielich.amu.edu.pl/Stefan_Banach/jpg/00xx.jpg}{Stefan Banach}}
(1892 -- 1945) \cite{MR1357166}. The Banach of \textquotedblleft Banach
space.\textquotedblright{} Banach called them \textquotedblleft B-spaces\textquotedblright{}
in his book. They were also formalized by Norbert Wiener (who traveled
in Europe in the 1920ties.) But the name \textquotedblleft Banach
space\textquotedblright{} stuck.

S. Banach, one of the founders of modern functional analysis and one
of the original members of the Lwów School of Mathematics, in Poland
between the two World Wars. His 1932 book, \emph{Théorie des opérations
linéaires} (Theory of Linear Operations), is the first monograph on
the general theory of functional analysis.\medskip{}

\textbf{\href{http://www-history.mcs.st-andrews.ac.uk/BigPictures/Bohr_Harald_2.jpeg}{Harald August Bohr}
}(1887 -- 1951) was a Danish mathematician and soccer player. Best
known for his theory of almost periodic functions. -- In modern language
it became the Bohr-compactification. (Different from the alternative
compactifications we discussed above.) He is the brother of the physicist
Niels Bohr.\medskip{}

\textbf{\href{http://www.nobelprize.org/nobel_prizes/physics/laureates/1922/bohr.jpg}{Niels Henrik David Bohr}}
(1885 -- 1962) was a Danish physicist who made foundational contributions
to understanding atomic structure and quantum theory, the \textquotedblleft Bohr-atom\textquotedblright ,
justifying the Balmer series for the visible spectral lines of the
hydrogen atom; received the Nobel Prize in Physics in 1922; -- ``for
his services in the investigation of the structure of atoms, and of
the radiation emanating from them''. Based on his liquid drop model
of the nucleus, Bohr concluded that it was the uranium-235 isotope,
and not the more abundant uranium-238, that was primarily responsible
for fission. 

In September 1941, at the start of WWII, Heisenberg, who had become
head of the German nuclear energy project, visited Bohr in Copenhagen.
During this meeting the two had discussions about possible plans by
the two sides in the War, for a fission bomb, the content of the discussions
have caused much speculation. Michael Frayn's 1998 play \textquotedblleft Copenhagen\textquotedblright{}
explores what might have happened at the 1941 meeting between Heisenberg
and Bohr.\medskip{}

\textbf{\href{https://commons.wikimedia.org/wiki/File:Max_Born.jpg}{Max Born}}
(1882 -- 1970) \cite{BP44}, a German physicist and mathematician,
a pioneer in the early development of quantum mechanics; also in solid-state
physics, and optics. Won the 1954 Nobel Prize in Physics for his ``fundamental
research in Quantum Mechanics, especially in the statistical interpretation
of the wave function.'' His assistants at Göttingen, between the
two World Wars, included Enrico Fermi, Werner Heisenberg, and Eugene
Wigner, among others. His early education was at Breslau, where his
fellow students included Otto Toeplitz and Ernst Hellinger. In 1926,
he formulated the now-standard interpretation of the probability density
function for states (represented as equivalence classes of solutions
to the Schrödinger equation.) After the Nazi Party came to power in
Germany in 1933, Born was suspended. Subsequently he held positions
at Johns Hopkins University, at Princeton University, and he settled
down at St John's College, Cambridge (UK). A quote: \textquotedblleft I
believe that ideas such as absolute certitude, absolute exactness,
final truth, etc. are figments of the imagination which should not
be admissible in any field of science. On the other hand, any assertion
of probability is either right or wrong from the standpoint of the
theory on which it is based.\textquotedblright{} Max Born (1954.)

\medskip{}

\textbf{\href{http://upload.wikimedia.org/wikipedia/commons/thumb/7/7d/Dirac_3.jpg/220px-Dirac_3.jpg}{Paul Adrien Maurice Dirac}}
(1902 -- 1984) \cite{Dir35,Dir47}. Cited in connection with the \textquotedblleft Dirac
equation\textquotedblright{} and especially our notation for vectors
and operators in Hilbert space, as well as the axioms of observables,
states and measurements.

P. Dirac, an English theoretical physicist; fundamental contributions
to the early development of both quantum mechanics and quantum electrodynamics.
He was the Lucasian Professor of Mathematics at the University of
Cambridge. Notable discoveries, the Dirac equation, which describes
the behavior of fermions and predicted the existence of antimatter.
Dirac shared the Nobel Prize in Physics for 1933 with Erwin Schrödinger,
``for the discovery of new productive forms of atomic theory.''
A rare interview with Dirac; see the link: \url{http://www.math.rutgers.edu/~greenfie/mill_courses/math421/int.html}\medskip{}

\textbf{\href{http://www.apprendre-en-ligne.net/ephemerides/matheux/m13-18.jpg}{Wolfgang Döblin}}
(1915 -- 40), French-German mathematician, and probabilist. Studied
probability theory in Paris, under Fréchet. Served in the French army
in the Ardennes when World War II broke out in 1939. There, he wrote
down his work on the Chapman-Kolmogorov equation. And he sent it in
a sealed envelope to the French Academy of Sciences. In 1940, after
burning his mathematical notes, he took his own life as the German
troops came in sight. In 2000, the sealed envelope was opened, revealing
that, at the time, Döblin had anticipated the theory of Markov processes,
It\={o}'s lemma (now the It\={o}\textendash Döblin lemma), and parts
of stochastic calculus.

\medskip{}

\textbf{\href{http://upload.wikimedia.org/wikipedia/commons/thumb/3/3d/Freeman_Dyson.jpg/400px-Freeman_Dyson.jpg}{Freeman John Dyson}
}(1923 -- ) theoretical physicist and mathematician; -- known for
his contributions to quantum electrodynamics, solid-state physics,
astronomy, and to nuclear engineering. Within mathematics, he is known
for his work on random matrices; his discovery of a perturbation expansion
(the Dyson expansion.) He is a regular contributor to The New York
Review of Books. Awards: the Lorentz Medal, the Max Planck Medal,
and the Enrico Fermi Award.\medskip{}

\textbf{\href{http://www5.in.tum.de/wiki/uploads/3/31/Friedrichs.jpeg}{Kurt Otto Friedrichs}
}(1901 -- 1982) \cite{Fri80,FL28}. Is the Friedrichs of the Friedrichs
extension; referring to the following Theorem: Every semibounded operator
$S$ with dense domain in Hilbert space has a selfadjoint extension
having the same lower bound as $S$. There are other semibounded and
selfadjoint extensions of $S$; -- they were found later by M. Krein.

K. Friedrichs, a noted German-American mathematician; a co-founder
of The Courant Institute at New York University and recipient of the
National Medal of Science.

\emph{A story}: Selfadjoint operators, and the gulf between the lingo
and culture of mathematics and of physics:

Peter Lax relates the following conversation in German between K.O.
Friedrichs and W. Heisenberg, to have been taken place in the late
1950ties, in New York, when Heisenberg visited The Courant Institute
at NYU. (The two had been colleagues in Germany before the war.) As
a gracious host, Friedrichs praised Heisenberg for having created
quantum mechanics. -- After an awkward silence, Friedrich went on:
\textquotedblleft ..and we owe to von Neumann our understanding of
the crucial difference between a selfadjoint operator and one that
is merely symmetric.\textquotedblright{} Another silence, and then
-- Heisenberg: \textquotedblleft What is the difference?\textquotedblright \medskip{}

\textbf{\href{http://www.maa.org/sites/default/files/images/upload_library/46/0_Halmos_Photos/lot_0091/e_ph_0093_01.jpg}{Lars Gårding}
}(1919 -- 2014). The \textquotedblleft G\textquotedblright{} in Gårding
vectors (representations of Lie groups), and in Gårding-Wightman quantum
fields.\medskip{}

\textbf{\href{http://upload.wikimedia.org/wikipedia/en/thumb/5/5e/IM_Gelfand.jpg/220px-IM_Gelfand.jpg}{Israel Moiseevich Gelfand}
}(1913 -- 2009)\textbf{ }\cite{MR0112604,MR0111991,MR0126719}. Is
the \textquotedblleft G\textquotedblright{} in GNS (Gelfand-Naimark-Segal),
the correspondence between states and cyclic representations. 

I. Gelfand, also written Israïl Moyseyovich Gel'fand, or Izrail M.
Gelfand, a Russian-American mathematician; major contributions to
many branches of mathematics: representation theory and functional
analysis. The recipient of numerous awards and honors, including the
Order of Lenin and the Wolf Prize, -- a lifelong academic, serving
decades as a professor at Moscow State University and, after immigrating
to the United States shortly before his 76th birthday, at Rutgers
University.\medskip{}

\textbf{\href{http://upload.wikimedia.org/wikipedia/commons/b/b0/Heisenberg_10.jpg}{Werner Karl Heisenberg}
}(1901 -- 1976) \cite{Hei69}. Is the Heisenberg of the Heisenberg
uncertainty principle for the operators $P$ (momentum) and $Q$ (position),
and of matrix mechanics, as the first mathematical formulation of
quantum observables. In Heisenberg\textquoteright s picture, the dynamics,
the observables are studied as function of time; by contrast to Schrödinger\textquoteright s
model which have the states (wave-functions) functions of time, and
satisfying a PDE wave equation, now called the Schrödinger equation.
In the late 1920ties, the two pictures, that of Heisenberg and of
Schrödinger were thought to be irreconcilable. Work of von Neumann
in 1932 demonstrated that they in fact are equivalent.

W. Heisenberg; one of the key creators of quantum mechanics. A 1925
paper was a breakthrough. In the subsequent series of papers with
Max Born and Pascual Jordan, this matrix formulation of quantum mechanics
took a mathematical rigorous formulation. In 1927 he published his
uncertainty principle. Heisenberg was awarded the Nobel Prize in Physics
for 1932 ``for the creation of quantum mechanics.'' He made important
contributions to the theories of the hydrodynamics of turbulent flows,
the atomic nucleus.\medskip{}

\textbf{\href{http://www.nndb.com/people/735/000087474/david-hilbert-1.jpg}{David Hilbert}
}(1862 -- 1943)\textbf{ }\cite{MR1512197,MR1512123,MR1557926}. Cited
in connection with the early formulations of the theory of operators
in (what is now called) Hilbert space. The name Hilbert space was
suggested by von Neumann who studied with Hilbert in the early 1930ties,
before he moved to the USA. (The early papers by von Neumann are in
German.) 

D. Hilbert is recognized as one of the most influential and universal
mathematicians of the 19th and early 20th centuries. Discovered and
developed invariant theory and axiomatization of geometry. In his
1900 presentation of a collection of research problems, he set the
course for much of the mathematical research of the 20th century.\medskip{}

\textbf{\href{http://estebanmoro.org/wp-content/uploads/2008/11/ito.jpeg}{Kiyoshi It\={o}}
}(1915 -- 2008) \cite{MR2397781,MR2053326}. Cited in connection with
Brownian motion, It\={o}-calculus, and stochastic processes. Making
connection to functional analysis via the theory of semigroups of
operators (Hille and Phillips.) \medskip{}

\textbf{\href{http://math.ecnu.edu.cn/RCOA/events/2012SpringOperatorAlgebrasProgram/Richard\%20Kadison/Kadison-Richard-photo.jpg}{Richard V. Kadison}
}\cite{MR0123922,MR1468229} (1925 -- ) \dots{} known for his contributions
to the study of operator algebras. Is the \textquotedblleft K\textquotedblright{}
in the Kadison-Singer problem (see \cite{MSS15}); and the \textquotedblleft K\textquotedblright{}
in the Fuglede-Kadison determinant. He is a Gustave C. Kuemmerle Professor
in the Department of Mathematics of the University of Pennsylvania;
was awarded the Leroy P. Steele Prize for Lifetime Achievement, in
1999.\medskip{}

\textbf{\href{http://www.log24.com/log/pix04A/040818-Kakutani.jpg}{Shizuo Kakutani}
}(1911 -- 2004). Of his theorems in functional analysis, there is
the Kakutani fixed-point theorem (a generalization of Brouwer's fixed-point
theorem); -- with such applications as to the Nash equilibrium in
game theory. Also notable is his solution of the Poisson equation
using the methods of stochastic analysis; as well as his pioneering
advances in our understanding of two-dimensional Brownian motion;
and its applications to PDE, and to potential theory.\medskip{}

\textbf{\href{http://www.emomi.com/history/mechanics_odessa/university/images/krein.jpg}{Mark Grigorievich Krein}}
\cite{Kre46,Kr55} (1907--1989). Is the Krein of Krein-Milman on convex
weak $*$-compact sets. Soviet mathematician; known for pioneering
works in operator theory, mathematical physics, the problem of moments,
functional and classical analysis, and representation theory. Winner
of the Wolf Prize, 1982. His list of former students includes David
Milman, Mark Naimark, Izrail Glazman, Moshe Livshits. 

\medskip{}

\textbf{\href{http://upload.wikimedia.org/wikipedia/commons/thumb/c/c3/Peter_Lax_in_Tokyo.jpg/220px-Peter_Lax_in_Tokyo.jpg}{Peter David Lax}
}(1926 -- ). The \textquotedblleft L\textquotedblright{} in Lax-Phillips
scattering theory, and the Lax-Milgram lemma. A pioneer in PDE, and
in many areas of applied mathematics; -- especially as they connect
to functional analysis.\medskip{}

\textbf{\href{http://news.harvard.edu/gazette/2006/04.06/photos/20-mackey.jpg}{George Whitelaw Mackey}
}(1916 -- 2006). The first \textquotedblleft M\textquotedblright{}
in the \textquotedblleft Mackey-machine,\textquotedblright{} a systematic
tool for constructing unitary representations of Lie groups, as induced
representations. A pioneer in non-commutative harmonic analysis, and
its applications to physics, to number theory, and to ergodic theory.\medskip{}

\textbf{\href{https://web.math.princeton.edu/~nelson/100_0028s.jpg}{Edward (Ed) Nelson}
}(1932 -- 2014) \cite{Ne69,Nel59}. Cited in connection with spectral
representation, and Brownian motion. \medskip{}

\textbf{\href{http://www.maa.org/sites/default/files/images/upload_library/46/0_Halmos_Photos/lot_0235/e_ph_0238_01.jpg}{Ralph Saul Phillips}
}(1913 -- 1998) \cite{LP89}. Cited in connection with the foundations
of functional analysis, especially the theory of semigroups of bounded
operators acting on Banach space. \medskip{}

\textbf{\href{https://en.wikipedia.org/wiki/File:Frigyes_Riesz.jpeg}{Frigyes Riesz}}
(1880 -- 1956) made fundamental contributions to functional analysis,
and to the theory of operators in Hilbert space. We frequently use
his Riesz representation theorem. He also did some of the fundamental
work, developing functional analysis for applications, especially
to spectral theory, and ergodic theory; both important in physics.
And with his brother, Marcel Riesz, work in harmonic analysis.

\medskip{}

\textbf{\href{http://matematica.unibocconi.it/sites/default/files/Riesz_Marcel_2.jpg}{Marcel Riesz}
}(1886 -- 1969), born in Hungary, was the younger brother of the mathematician
Frigyes Riesz (the two are known for the F. and M. Riesz theorem).
Both are pioneers in Functional Analysis. M. Riesz moved to Sweden
in 1911 where he taught at Stockholm University and at Lund University.
His former students include Harald Cramér, Einar Hille (of Hille-Phillips),
Otto Frostman (potential theory), Lars Hörmander (PDE), and Olaf Thorin
(harmonic analysis).\medskip{}

\textbf{\href{http://upload.wikimedia.org/wikipedia/commons/2/26/Erwin_Schr\%C3\%B6dinger.jpg}{Erwin Rudolf Josef Alexander Schrödinger}
}(1887 -- 1961) \cite{MR1782215,Sch40,Sch32}. Is the Schrödinger
of the Schrödinger equation; the PDE which governs the dynamics of
quantum states (as wave-functions). 

E. Schrödinger, a Nobel Prize in physics. -- quantum theory forming
the basis of wave mechanics: he formulated the wave equation (stationary
and time-dependent Schrödinger equation) , and he Schrödinger proposed
an original interpretation of the physical meaning of the wave function;
formalized the notion of entanglement. He was critical the conventional
Copenhagen interpretation of quantum mechanics (using e.g. the paradox
of Schrödinger's cat).\medskip{}

\textbf{\href{http://upload.wikimedia.org/wikipedia/commons/thumb/d/d3/HermannSchwarz.jpeg/220px-HermannSchwarz.jpeg}{Karl Hermann Amandus Schwarz}}
(1843 -- 1921) \cite{Schwarz197001}. Is the Schwarz of the Cauchy-Schwarz
inequality. H.A. Schwarz is German and is a contemporary of K. Weierstrass. 

H.A. Schwarz, a German mathematician, known for his work in complex
analysis. At Göttingen, he pioneered of function theory, differential
geometry and the calculus of variations.\medskip{}

\textbf{\href{http://www-history.mcs.st-and.ac.uk/BigPictures/Schwartz_5.jpeg}{Laurent-Moïse Schwartz}}
(1915 -- 2002) \cite{MR1609224,MR0107312,MR0107812}. Is the Schwartz
(French) of the theory of distributions (dating the 1950ties), also
now named \textquotedblleft generalized functions\textquotedblright{}
in the books by Gelfand et al. Parts of this theory were developed
independently on the two sides of the Iron-Curtain;-- in the time
of the Cold War. \medskip{}

\textbf{\href{http://www.crpc.rice.edu/newsletters/win96/photos/schwartz.jpg}{Jacob Theodore "Jack" Schwartz}
}(1930 -- 2009)\textbf{ }\cite{DS88a,DS88b,DS88c}. Is the Schwartz
of the book set ``linear operators'' by Dunford and Schwartz. Vol
II \cite{DS88b} is one of the best presentation of the theory of
unbounded operators. \medskip{}

\textbf{\href{http://66.90.101.64/arkblog/Irving\%20Segal.jpg}{Irving Ezra Segal}
}(1918 -- 1998) \cite{Seg50}. Cited in connection with the foundations
of functional analysis, and pioneering research in mathematical physics.
Is the \textquotedblleft S\textquotedblright{} in GNS (Gelfand-Naimark-Segal).
Segal proved the Plancherel theorem in a very general framework: locally
compact unimodular groups. For any locally compact unimodular group,
Segal established a Plancherel formula; see \cite{Seg50}. Segal showed
that there is a Plancherel formula, despite the fact that it may not
be feasible, for all locally compact unimodular groups, to \textquotedblleft write
down\textquotedblright{} all the irreducible unitary representations.
\medskip{}

\textbf{\href{http://upload.wikimedia.org/wikipedia/commons/f/f2/Isadore_Singer_1977.jpeg}{Isadore Manuel Singer}
}(1924 -- ) is an Institute Professor at the Massachusetts Institute
of Technology; He is the \textquotedblleft S\textquotedblright{} in
the Atiyah--Singer index theorem (1962), Michael Atiyah is the \textquotedblleft A.\textquotedblright{}
Also of note: The Atiyah--Hitchin--Singer theorem, and The Atiyah--Patodi--Singer
eta-invariant.

\medskip{}

\textbf{\href{http://www.math.bme.hu/~schmidt/Stone.jpg}{Marshall Harvey Stone}
}(1903 -- 1989)\textbf{ }\cite{Sto51,Sto90}. Is the \textquotedblleft S\textquotedblright{}
in the Stone-Weierstrass theorem; and in the Stone-von Neumann uniqueness
theorem (see \cite{vN32c,vN31}); the latter to the effect that any
two representations of Heisenberg\textquoteright s commutation relations
in the same (finite!) number of degrees of freedom are unitarily equivalent.
Stone was the son of Harlan Fiske Stone, Chief Justice of the United
States in 1941-1946. Marshall Stone completed a Harvard Ph.D. in 1926,
with a thesis supervised by George David Birkhoff. He taught at Harvard,
Yale, and Columbia University. And he was promoted to a full Professor
at Harvard in 1937. In 1946, he became the chairman of the Mathematics
Department at the University of Chicago. His 1932 monograph titled
\textquotedblleft Linear transformations in Hilbert space and their
applications to analysis\textquotedblright{} develops the theory of
selfadjoint operators, turning it into a form which is now a central
part of functional analysis. Theorems that carry his name: The Banach-Stone
theorem, The Glivenko-Stone theorem, Stone duality, The Stone-Weierstrass
theorem, Stone's representation theorem for Boolean algebras, Stone's
theorem for one-parameter unitary groups, Stone-\v{C}ech compactification,
and The Stone-von Neumann uniqueness theorem. \medskip{}

\textbf{\href{http://upload.wikimedia.org/wikipedia/commons/5/5e/JohnvonNeumann-LosAlamos.gif}{John von Neumann}
}(1903 -- 1957) \cite{vN31,vN32a}. Cited in connection with the Stone-von
Neumann uniqueness theorem, the deficiency indices which determine
parameters for possible selfadjoint extensions of given Hermitian
(formally selfadjoint, symmetric) with dense domain in Hilbert space. 

J. von Neumann, Hungarian-American; inventor and polymath. He made
major contributions to: foundations of mathematics, functional analysis,
ergodic theory, numerical analysis, physics (quantum mechanics, hydrodynamics,
and economics (game theory), computing (von Neumann architecture,
linear programming, self-replicating machines, stochastic computing
(Monte-Carlo\footnote{``Monte-Carlo'' means ``simulation'' with computer generated random
number.})), -- was a pioneer of the application of operator theory to quantum
mechanics, a principal member of the Manhattan Project and the Institute
for Advanced Study in Princeton. -- A key figure in the development
of game theory, cellular automata, and the digital computer.\medskip{}

\textbf{\href{http://image2.findagrave.com/photos250/photos/2003/124/7304697_1052135619.jpg}{Norbert Wiener}
}(1894 -- 1964)\textbf{ }\cite{MR0057817,MR0057763}. Cited in connection
with Brownian motion, Wiener measure, and stochastic processes. And
more directly, the \textquotedblleft Wiener\textquotedblright{} of
Paley-Wiener spaces; -- at the crossroads of harmonic analysis and
functional analysis. Also the Wiener of filters in signal processing;
high-pass/low-pass etc.\index{signal processing}

\vspace{2cm}

\begin{quote}
Der skal et par dumheder

med i en bog ....

for at også de dumme

skal syns, den er klog. 

--- Piet Hein.\sindex[nam]{Hein, P., (1905-1996)}\emph{\medskip{}
}

\emph{Translation: }

\emph{Your book should include a few stupidities}

\emph{mixing them in, -- this is art. }

\emph{so that also the stupid will think it is smart.}
\end{quote}

\chapter{Prizes and Fame}

\paragraph*{Nobel Prize in Physics:}
\begin{itemize}
\item N. Bohr, M. Born, P. Dirac, A. Einstein, W. Heisenberg, E. Wigner. 
\end{itemize}

\paragraph*{Fields Medal (Math):}
\begin{itemize}
\item Sir Michael Atiyah, A. Connes, D. Mumford.
\end{itemize}

\paragraph*{The National Medal of Science: }
\begin{itemize}
\item K. Friedrichs, P. Lax, I.M. Singer, N. Wiener, E. Wigner.
\end{itemize}

\paragraph*{The Wolf Prize: }
\begin{itemize}
\item I. Gelfand, K. It\={o}, A.N. Kolmogorov, M.G. Krein, P. Lax. 
\end{itemize}

\paragraph*{The Abel Prize: }
\begin{itemize}
\item Sir Michael Atiyah, P. Lax, I.M. Singer.
\end{itemize}

\paragraph*{Presidential Medal of Freedom: }
\begin{itemize}
\item J. von Neumann.
\end{itemize}
\phantomsection
\cleardoublepage
%\addcontentsline{toc}{chapter}{Quotes: index of credits}

\printindex[nam]{}

\end{appendix}

\backmatter

\phantomsection
\cleardoublepage
\addcontentsline{toc}{chapter}{\listexercisename}
\listofexercises

\phantomsection
\cleardoublepage
\addcontentsline{toc}{chapter}{\listfigurename}
\listoffigures

\cleardoublepage

\bibliographystyle{amsalpha}
\phantomsection\addcontentsline{toc}{chapter}{\bibname}\bibliography{ref_FA}

\phantomsection
\cleardoublepage
\addcontentsline{toc}{chapter}{Index}

\printindex[idx]{}
\end{document}